\documentclass[12pt]{article}

\usepackage[utf8]{inputenc}
\usepackage{amsmath}
\usepackage{amsfonts}

\usepackage{amssymb}

\usepackage{float}

\def\ba{\begin{eqnarray}}
\def\ea{\end{eqnarray}}

\def\P{\hbox{\bf P}}
\def\R{\hat{R}}
\def\K{\hbox{\bf K}}
\def\Sha{{\amalg\hspace{-.16cm} \amalg}}

\def\beps{\bar{\epsilon}}

\def\lb{\label}
\def\be{\begin{equation}}
\def\ee{\end{equation}}

\newtheorem{proposition}{Proposition}

\def\qed{\rule{5pt}{5pt}}

\begin{document}

\title{\bf LECTURES ON QUANTUM GROUPS AND YANG-BAXTER EQUATIONS}
\author{\Large {\bf A. P. Isaev} \\[0.2cm]
  Bogoliubov Lab. of Theoretical Physics, \\
 JINR, 141 980 Dubna, Moscow Reg. \\ and \\ 
 Lomonosov Moscow State University, \\
 Physics Faculty}

\date{}
\maketitle


\begin{center}
{\bf \large Abstract.}
\end{center}
The principles of the theory of quantum groups are reviewed from the point of view of the
possibility of their use for deformations of symmetries in physical models. The R-matrix approach
to the theory of quantum groups is discussed in detail and is taken as the basis of the quantization of
classical Lie groups, as well as some Lie supergroups.
We start by laying out the foundations of non-commutative and non-cocommutative Hopf algebras. Much attention has been paid
to Hecke and Birman-Murakami-Wenzl (BMW)
 R-matrices and related quantum
matrix algebras.
Trigonometric solutions of the Yang-Baxter
equation associated with the quantum groups $GL_q(N)$, $SO_q(N)$, $Sp_q(2n)$ and supergroups
$GL_q(N|M)$, $Osp_q(N|2m)$, as well as their rational (Yangian) limits, are presented. Rational $R$-matrices for exceptional
Lie algebras and elliptic
solutions of the Yang-Baxter equation are also considered. The basic concepts of the group algebra of
the braid group and its finite dimensional quotients (such as Hecke
and BMW algebras)
are outlined. A sketch of the representation theories of the Hecke and BMW algebras is given (including methods for finding
idempotents and their quantum dimensions). Applications of the theory of quantum groups and Yang-Baxter equations in various
 areas of theoretical physics are briefly discussed.

This is a modified version of the review paper published in 2004 as
a preprint \cite{IsBonn} of the Max-Planck-Institut für Mathematik in Bonn.

\vspace{1.5cm}



\newpage

{\footnotesize
\tableofcontents
}

\newpage

\section{INTRODUCTION}

In modern theoretical and mathematical physics, the ideas of symmetry and invariance play a very important role. Sets of symmetry transformations form groups and therefore the most natural language for
describing symmetries is the group theory language.

About 40 years ago, in the study of quantum integrable systems \cite{2,3,4}, in particular in the
framework of the quantum inverse scattering method \cite{1}, new algebraic structures arose, the
generalizations of which were later called quantum groups \cite{13}\footnote{In pure mathematics
analogous structures appeared as nontrivial
examples of ``ring-groups'' introduced by
G.I. Kac; see e.g. \cite{Enock} and references therein.}. The Yang-Baxter equations became a unifying basis
of all these investigations. The most important nontrivial
examples of quantum groups are quantizations (or deformations)
of ordinary classical Lie groups and algebras (more precisely,
one considers the deformations of the algebra
of functions of a Lie group and the universal enveloping
of a Lie algebra). The quantization is accompanied by the introduction of an additional parameter
$q$ (the deformation parameter), which plays a role analogous to the role of Planck's constant in
quantum mechanics. In the limit $q
\to 1$, the quantum groups and algebras go over into the classical ones. Although quantum
groups are deformations of the usual groups;
 they nevertheless still possess several properties
that make it possible to speak of them as "symmetry groups". Moreover, one can claim that the
quantum groups serve as symmetries and provide integrability in exactly solvable quantum systems
(see, for example, Refs. \cite{5,5a,5'})\footnote{The Yangian symmetries are
the symmetries of the same type.}. In this connection, the idea naturally arises of looking for
and constructing other physical models possessing such quantum symmetries. Some of the realizations
of this idea use the similarity of the representation theories of quantum and classical Lie groups
and algebras (for
$q$ is not equal to the root of unity).
 Thus, 
 there were attempts to apply
 quantum groups and algebras 
 in the classification and phenomenology of elementary particles and in nuclear spectroscopy investigations. Further, it is
natural to investigate already existing field-theoretical models
 (especially gauge quantum field theories)
  from the point of view of 
  relations (see, e.g.,
Ref. \cite{6}) to the noncommutative geometry \cite{6a}\footnote{After the quantization of any
Poisson manifold
\cite{Konts} and the appearance of papers \cite{6'}
the subject of field theories on noncommutative
spaces became very popular from the point of view
of string theories.}  and, in
particular, the possibility of their invariance with respect to quantum-group transformations. An
attractive idea is that of relating the deformation parameters of quantum groups to the mixing
angles that occur in the Standard Model as free parameters. One of the possible realizations of
this idea was proposed in Ref. \cite{IsPop7}
(see also \cite{7'}). Of course,
it is necessary to mention here numerous
attempts to deform the Lorentz and Poincar\'{e} groups
and construct a covariant quantum Minkowski
space-time corresponding to these deformations \cite{8,9Og,9K,9}.

It is clear that the approaches listed above (associated with deformations of symmetries in
physics) present only a small fraction of all the
 applications of the theory of quantum groups.
Quantum groups and Yang-Baxter equations naturally
 arise in many problems of theoretical physics,
and this makes it possible to speak of them
and their theories as an important paradigm in
mathematical physics. Unfortunately, the strict limits of these Lectures make it impossible to
discuss in detail all applications of quantum groups and Yang-Baxter equations. I have therefore
restricted myself to a brief listing of certain areas in theoretical and mathematical physics in
which quantum groups and Yang-Baxter equations play an important role. The incomplete list is given
in  Section {\bf \ref{appli}}. In Section {\bf \ref{hopf}}, the mathematical foundations of the theory of quantum groups are
outlined. A significant part of Sec. {\bf \ref{ybeqg}},
is a detailed exposition of the results of the famous work
by Faddeev, Reshetikhin and Takhtajan
\cite{10} who formulated the $R$-matrix approach to the theory of quantum groups. In this
Section, we also consider problems of invariant Baxterization of $R$-matrices, multi-parameter
deformations of Lie groups, and the quantization of some Lie supergroups.
The rational solutions of the Yang-Baxter equation for exceptional
Lie algebras are also considered in this Section.
At the end of Sec. {\bf \ref{ybeqg}}, we
present the basic notions of the theory of quantum Knizhnik - Zamolodchikov equations and discuss
elliptic solutions of the Yang-Baxter equation for which the algebraic basis (the type of quantum
universal enveloping Lie algebras $U_q(g)$ in the case of trigonometric solutions) has not yet been
completely clarified
(see, however \cite{19}, \cite{SkOd}).
In Sec. {\bf \ref{gabg}}, we briefly discuss the
 group algebra of a braid group and its finite dimensional
 quotients such as the Hecke
and  Birman-Murakami-Wenzl (BMW) algebras. A sketch of the representation theories of the Hecke and BMW algebras is given (including methods for finding
idempotents and their quantum dimensions). The part of
 the content of Sec. {\bf \ref{gabg}}
 can be considered as a different
presentation of some facts from Sec. {\bf \ref{ybeqg}}.
In Sec. {\bf \ref{appli}}, as we have already mentioned, some applications
of quantum groups and the Yang-Baxter equations are outlined.

\vspace{0.3cm}
These Lectures are an introductory review
 based on the paper \cite{IsFirst} and the
 MPIM (Bonn) preprint \cite{IsBonn} published in 1995 and 2004.
 Comparing to the previous
 version \cite{IsBonn}, this text has been
 considerably changed only in Subsections {\bf \ref{NuRm}},
{\bf \ref{qmatalg}}, {\bf \ref{qalur}}, {\bf \ref{baxtel}},
{\bf \ref{qdlink}} and {\bf \ref{ybfd}}.
 The structure of the review has also been significantly modified.
Three new subsections {\bf \ref{difcal}},
{\bf \ref{solYB}} and {\bf \ref{repbmw}}
 have been added, and subsection {\bf \ref{bmwalg}} has been extended.
 In the course of the presentation we have also tried briefly mention
some of the new results. According to this, we have refreshed
 the list of references.

I would like to thank G.~Arutyunov, L.~Castellani,
S.~Derkachov, A.N.~Kirillov,
P.P.~Kulish, A.I.~Molev, O.V.~Ogievetsky, Z.~Popowicz,
P.N.~Pyatov and  V.O.~Tarasov for valuable discussions and comments.
I am especially grateful to O.V. Ogievetsky for
discussing the content of Sect. {\bf \ref{gabg}}.

\vspace{0.5cm}

\noindent
{\sf  In memory of Petr Petrovich Kulish (1944 -- 2016)}

\section{HOPF ALGEBRAS\label{hopf}}


\subsection{\bf \em Coalgebras\label{hopf1}}
\setcounter{equation}0

We consider an associative unital algebra ${\cal A}$ (over the field of complex numbers
$\mathbb{C}$;
in what follows, all algebras that are introduced will also be understood to be over the
field of complex numbers). Each element of ${\cal A}$
can be expressed as a linear combination of
basis elements $e_{i} \in {\cal A}$, where $i = 1,2,3, \dots$, and
the identity element $I$ is given by the formula
$$
I = E^{i}\, e_{i} \;\;\;\;\;\; (E^{i} \in \mathbb{C}) \; ,
$$
(we imply summation over repeated indices).
Then for any two elements
$e_{i}$ and $e_{j}$ we define their multiplication in the form
\be
\lb{2.1}
{\cal A} \otimes {\cal A} \stackrel{m}{\longrightarrow} {\cal A}
\;\;\;\; \Rightarrow \;\;\;\;
e_{i} \cdot e_{j} = m^{k}_{ij} e_{k} \, ,
\ee
where $m^{k}_{ij}$ is a certain set of complex numbers that satisfy the condition
\be
\lb{2.2}
E^{i}m^{k}_{ij} = m^{k}_{ji}E^{i} =  \delta^{k}_{j}
\ee
for the identity element, and also the condition
\be
\lb{2.3}
m^{l}_{ij} m^{n}_{lk} = m^{n}_{il} m^{l}_{jk} \equiv m^{n}_{ijk} \; ,
\ee
that is equivalent to the condition of associativity
for the algebra ${\cal A}$:
\be
\lb{2.4}
(e_{i} e_{j}) e_{k} = e_{i} (e_{j} e_{k}).
\ee
The condition of associativity (\ref{2.4}) for the multiplication (\ref{2.1}) can obviously be
represented in the form of the commutativity of the diagram:

\unitlength=0.7cm
\begin{picture}(15,4.5)(-3,-1)
\put(3.3,2.5){${\cal A} \otimes {\cal A} \otimes {\cal A}$}
\put(6.5,2.7){\vector(1,0){2.2}}
\put(6.8,3){\footnotesize $id \otimes m$}
\put(4.2,0.5){${\cal A} \otimes {\cal A}$}
\put(5,2.2){\vector(0,-1){1.2}}
\put(3,1.5){\footnotesize $m \otimes id$}
\put(6,0.7){\vector(1,0){3.2}}
\put(7.3,1){\footnotesize $m$}
\put(9,2.5){${\cal A} \otimes {\cal A}$}
\put(9.8,2.2){\vector(0,-1){1.2}}
\put(10,1.5){\footnotesize $m$}
\put(9.5,0.5){${\cal A}$}

\put(4,-0.7){Fig. 1. Associativity axiom.}

\end{picture}

In Fig.1 the map $m$ represents multiplication:
${\cal A} \otimes {\cal A} \stackrel{m}{\longrightarrow} {\cal A}$, and $id$ denotes the
identity mapping. The existence of the unit element $I$ means that one can define a mapping i:
$\mathbb{C} \to {\cal A}$ (embedding of $\mathbb{C}$ in ${\cal A}$)
\be
\lb{2.5}
k \stackrel{\bf i}{\longrightarrow} k \cdot I \; , \;\;
k \in \mathbb{C}
\ee
For $I$ we have the condition (\ref{2.2}), which is visualized
 as the diagram in Fig.2:

\unitlength=0.8cm
\begin{picture}(15,5)(-3,0)
\put(4.2,2.5){$\mathbb{C} \otimes {\cal A}$}
\put(5,2.2){\vector(2,-1){1.5}}
\put(5,2.2){\vector(-2,1){0.1}}
\put(5,3){\vector(1,1){1}}
\put(4.3,3.5){{\bf i}$ \otimes id$}
\put(6.5,4.2){${\cal A} \otimes {\cal A}$}
\put(9,3){\vector(-1,1){1}}
\put(8.7,3.5){$id \otimes ${\bf i}}
\put(8.5,2.5){${\cal A} \otimes \mathbb{C}$}
\put(9,2.2){\vector(-2,-1){1.5}}
\put(9,2.2){\vector(2,1){0.1}}
\put(6.8,1.2){${\cal A}$}
\put(7,4){\vector(0,-1){2.2}}
\put(7.2,3){$m$}
\put(3.5,0.5){Fig. 2. \it Axioms for the identity.}
\end{picture}

\noindent
where the mappings
\be
\lb{natis}
\mathbb{C} \otimes {\cal A} \leftrightarrow {\cal A} \;\;
{\rm  and} \;\;
 {\cal A} \otimes \mathbb{C} \leftrightarrow {\cal A}
 \ee
 are natural isomorphisms.
One of the advantages of the diagrammatic language used here is that it leads
directly to the definition of a new fundamental object -- the coalgebra -- if we
reverse all the arrows in the diagrams of Fig.1 and Fig.2.


\newtheorem{def1}{Definition}
\begin{def1} \label{def1}
{\it A coalgebra ${\cal C}$ is a vector space
(with the basis $\{ e_{i} \}$) equipped with the
mapping
$\Delta : {\cal C} \rightarrow {\cal C} \otimes {\cal C}$
\be
\lb{2.6}
\Delta (e_{i}) = \Delta^{kj}_{i} e_{k} \otimes e_{j} \; ,
\ee
which is called comultiplication, and also equipped with the mapping
$\epsilon : {\cal C} \rightarrow \mathbb{C}$,
which is called the coidentity. The coalgebra ${\cal C}$ is called coassociative if the mapping
$\Delta$ satisfies the condition of coassociativity
(cf. the diagram in Fig.1 with the arrows reversed and
 the symbol $m$ changed to $\Delta$)
\be
\lb{2.7}
(id \otimes \Delta) \Delta = (\Delta \otimes id) \Delta
\;\;\;\;  \Rightarrow  \;\;\;\;
\Delta^{nl}_{i}\Delta^{kj}_{l} =
\Delta^{lj}_{i}\Delta^{nk}_{l} \equiv \Delta^{nkj}_{i} \; .
\ee
The coidentity $\epsilon$ must satisfy the following conditions (cf. the diagram in Fig.2 with arrows reversed and
 symbols $m,{\bf i}$ changed to $\Delta,\epsilon$)
\be
\lb{2.8}
m \bigl((\epsilon \otimes id) \Delta ({\cal C})\bigr) =
m \bigl((id \otimes \epsilon) \Delta ({\cal C})\bigr) = {\cal C}
\;\;\;\; \Rightarrow \;\;\;\;
\epsilon_{i}\Delta^{ij}_{k} = \Delta^{ji}_{k} \epsilon_{i} =
\delta^{j}_{k} \; .
\ee
Here $m$ realizes the natural isomorphisms (\ref{natis})
as a multiplication map:
 $m(c \otimes e_{i}) =  m(e_{i} \otimes c) = c \cdot e_{i}$
($\forall c \in \mathbb{C}$),
and the complex numbers $\epsilon_{i}$ are determined from the relations $\epsilon(e_{i}) = \epsilon_{i}$.}
\end{def1}

For algebras and coalgebras, the concepts of modules and comodules can be introduced.
Thus, if ${\cal A}$ is an algebra, the left ${\cal A}$-module
can be defined as a vector space $N$
and a mapping $\psi: \;\; {\cal A} \otimes N \rightarrow N$
(action of ${\cal A}$ on $N$) such that the diagrams on Fig.3
are commutative.

\unitlength=0.7cm
\begin{picture}(15,5)(-4,-1)
\put(0.4,2.5){${\cal A} \otimes {\cal A} \otimes N$}
\put(3.5,2.7){\vector(1,0){2.2}}
\put(4,3){\footnotesize  $id \otimes \psi$}
\put(1.1,0.5){${\cal A} \otimes  N$}
\put(2,2.2){\vector(0,-1){1.2}}
\put(0.3,1.5){\footnotesize $m \otimes id$}
\put(3,0.7){\vector(1,0){3.2}}
\put(4.3,1){\footnotesize $\psi$}
\put(6,2.5){${\cal A} \otimes N$}
\put(6.8,2.2){\vector(0,-1){1.2}}
\put(7,1.5){\footnotesize $\psi$}
\put(6.5,0.5){$N$}
\put(9,1.5){$\mathbb{C} \otimes N$}
\put(10,1.2){\vector(2,-1){1.5}}
\put(10,1.2){\vector(-2,1){0.1}}
\put(10,2){\vector(1,1){1}}
\put(9.3,2.5){\footnotesize {\bf i}$ \otimes id$}
\put(11.1,3.2){${\cal A} \otimes N$}
\put(11.7,0.2){$N$}
\put(12,3){\vector(0,-1){2.2}}
\put(12.2,2){\footnotesize $\psi$}
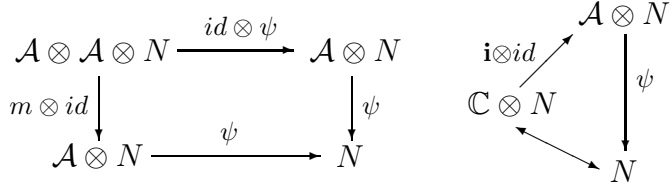
\put(1.2,-0.8){Fig. 3. Axioms for left ${\cal A}$-module.}
\end{picture}

 \noindent
In other words, the space $N$ is the representation
space of the algebra ${\cal A}$.

If $N$ is a (co)algebra and the mapping $\psi$ preserves the (co)algebraic structure
of $N$ (see below), then $N$ is called the
 {\it left ${\cal A}$-module (co)algebra}. The concepts of
{\it right module (co)algebra} is introduced similarly. If $N$ is
simultaneously the left and the right ${\cal A}$-module,
then $N$ is called the
{\it two-sided ${\cal A}$-module}.
It is obvious that the algebra ${\cal A}$ itself
is a two-sided ${\cal A}$-module for which the
left and right actions are given by the left and
right multiplications in the algebra.

Now suppose that ${\cal C}$ is a coalgebra; then a left
${\cal C}$-comodule can be defined as a
space $M$ together with a mapping
 $\Delta_{L}$: $M \rightarrow {\cal C} \otimes M$
(coaction of ${\cal C}$ on $M$)
satisfying the axioms in Fig.4 (in the diagrams of Fig.3, where the modules
were defined, it is necessary to reverse all the arrows):

\unitlength=0.7cm
\begin{picture}(15,4.5)(-4,-1)
\put(0.5,2.5){${\cal C} \otimes {\cal C} \otimes M$}
\put(5.5,2.7){\vector(-1,0){2}}
\put(4,3){\footnotesize $id \otimes \Delta_{L}$}
\put(1.3,0.5){${\cal C} \otimes  M$}
\put(2,1){\vector(0,1){1.2}}
\put(0.3,1.5){\footnotesize $\Delta \otimes id$}
\put(6.2,0.7){\vector(-1,0){3}}
\put(4.3,0.9){\footnotesize $\Delta_{L}$}
\put(6,2.5){${\cal C} \otimes M$}
\put(6.8,1){\vector(0,1){1.2}}
\put(7,1.5){\footnotesize $\Delta_{L}$}
\put(6.5,0.5){$M$}

\put(10,1.5){$\mathbb{C} \otimes M$}
\put(11,1.2){\vector(2,-1){1.5}}
\put(11,1.2){\vector(-2,1){0.1}}
\put(12,3){\vector(-1,-1){1}}
\put(10,2.5){\footnotesize $\epsilon \otimes id$}
\put(12,3.2){${\cal C} \otimes M$}
\put(12.8,0.2){$M$}
\put(13,0.8){\vector(0,1){2.2}}
\put(13.2,1.9){\footnotesize $\Delta_{L}$}

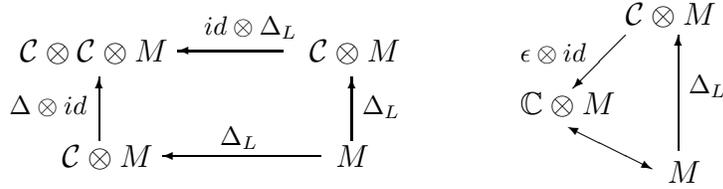
\put(3,-0.7){Fig. 4. Axioms for left ${\cal A}$-comodule.}
\end{picture}

If $M$ is a (co)algebra and the mapping $\Delta_L$ preserves the (co)algebraic
structure (for example, is a homomorphism; see below), then $M$ is called
a left ${\cal C}$-comodule (co)algebra. Right comodules are introduced
similarly, after which two-sided comodules are defined in the natural manner.
It is obvious that the coalgebra ${\cal C}$ is a two-sided ${\cal C}$-comodule.

Let ${\cal V}, \; \tilde{{\cal V}}$ be two vector spaces with bases
$\{ e_{i} \}, \; \{ \tilde{e}_{i} \}$.
We denote by ${\cal V}^{*}, \; \tilde{{\cal V}}^{*}$
the corresponding dual linear spaces whose
basis elements are linear functionals
$\{ e^{i} \}: \; {\cal V} \rightarrow \mathbb{C}$,
$\{ \tilde{e}^{i} \}: \tilde{{\cal V}} \rightarrow \mathbb{C}$.
For the values of these functionals, we use the expressions
$ \langle e^{i}\,|e_{j}\rangle $ and
$\langle \tilde{e}^{i}\,|\tilde{e}_{j}\rangle $.
For every mapping $L: \;\; {\cal V} \rightarrow \tilde{\cal V}$ it is possible to define a
unique mapping $L^{*}: \;\; \tilde{\cal V}^{*} \rightarrow {\cal V}^{*}$
induced by the equations
\be
\lb{2.9}
\langle \tilde{e}^{i} \, | L(e_{j}) \rangle  =
\langle  L^{*}(\tilde{e}^{i}) \, | e_{j} \rangle  ,
\ee
if the matrix $\langle e^{i}\, |e_{j}\rangle $ is invertible.
In addition, for the dual objects there exists the linear injection
$$
\rho: \;\; {\cal V}^{*} \otimes \tilde{\cal V}^{*} \rightarrow
( {\cal V} \otimes \tilde{\cal V} )^{*} \; ,
$$
which is given by the equations
$$
\langle \rho(e^{i} \otimes \tilde{e}^{j}) \,
| e_{k} \otimes \tilde{e}_{l} \rangle  =
\langle e^{i} \,| e_{k}\rangle \; \langle \tilde{e}^{j}\,
| \tilde{e}_{l} \rangle  \; .
$$
A consequence of these facts is that for every coalgebra
$({\cal C}, \; \Delta, \; \epsilon)$
it is possible to define an algebra ${\cal C}^{*}=
{\cal A}$ (as dual object to ${\cal C}$) with
multiplication $m =  \Delta^{*} \cdot \rho$ and
the unit element $I$ that satisfy the relations
($\forall a,a' \in {\cal A}$, $\forall c \in {\cal C}$):
$$
\langle a|c_{(1)}\rangle \langle a'|c_{(2)}\rangle  =
\langle \rho( a \otimes a')| \Delta(c)\rangle
= \langle \Delta^* \cdot \rho( a \otimes a')| c \rangle  =
\langle a \cdot a' \, | \, c \rangle \; ,
\;\;\;\;\;\;\; \langle  I | c \rangle  = \epsilon(c) \; .
$$
Here we denote $a \cdot a' := \Delta^* \cdot \rho( a \otimes a')$
 and use the convenient Sweedler notation of Ref. [11] for comultiplication
in ${\cal C}$ (cf. eq. (\ref{2.6})):
 \be
 \lb{sweed}
 \Delta(c) = \sum_{c} c_{(1)} \otimes c_{(2)} \; .
 \ee
The summation symbol $\sum_{c}$ will usually be
omitted in the equations. We also use the Sweedler notation
 for left and right coactions
$\Delta_L(v) = \bar{v}^{(-1)} \otimes v^{(0)}$ and
$\Delta_R(v) = v^{(0)} \otimes \bar{v}^{(1)}$,
 where index $(0)$ is
reserved for the comodule elements and summation symbols $\sum_{v}$ are also omitted.

Thus, duality in the diagrammatic definitions of the algebras and coalgebras
(reversal of the arrows) has in particular the consequence that the algebras
and coalgebras are indeed duals to each other.

It is natural to expect that an analogous duality can also be traced for
modules and comodules. Let ${\cal V}$ be a left comodule
for ${\cal C}$. Then the left
coaction of ${\cal C}$ on ${\cal V}$:
$v \mapsto \sum_{v} \; \bar{v}^{(-1)} \otimes v^{(0)} \;\;
(\bar{v}^{(-1)}\in {\cal C}, \;\; v^{(0)} \in {\cal V})$
induces a right action of ${\cal C}^{*}={\cal A}$ on
${\cal V}$:
$$
(v,a) \;\;\; \mapsto \;\;\;
v \triangleleft a =
\langle a\, |\bar{v}^{(-1)}\rangle  \; v^{(0)} \; ,
\;\;\;\;\;\; a \in {\cal A} \; ,
$$
and therefore ${\cal V}$ is a right module for ${\cal A}$.
Conversely, the right coaction of
${\cal C}$ on ${\cal V}$:
$v \mapsto v^{(0)} \otimes \bar{v}^{(1)}$ induces the left action of
${\cal A} = {\cal C}^{*}$ on ${\cal V}$:
$$
(a,v) \;\;\; \mapsto \;\;\; a \triangleright v
= v^{(0)} \langle a|\bar{v}^{(1)}\rangle   \; .
$$
From this we immediately conclude that the coassociative
coalgebra ${\cal C}$ (which coacts on itself by the coproduct) is
a natural module for its dual algebra ${\cal A}={\cal C}^{*}$. Indeed, the right
action ${\cal C} \otimes {\cal A} \rightarrow {\cal C}$ is determined by the equations
\be
\lb{2.10}
(c,a) \;\;\; \mapsto \;\;\; c \triangleleft a = \langle a|c_{(1)}\rangle  c_{(2)}
\ee
whereas for the left action ${\cal A} \otimes {\cal C} \rightarrow {\cal C}$ we have
\be
\lb{2.11}
(a,c) \;\;\; \mapsto \;\;\; a \triangleright c = c_{(1)} \langle a|c_{(2)}\rangle  \; .
\ee
Here $a \in {\cal A} \;\; c \in {\cal C}$.
The module axioms (shown as the diagrams in Fig. 3)
hold by virtue of the coassociativity of ${\cal C}$.

Finally, we note that the action of a certain algebra
$H$ on ${\cal C}$ from the
left (from the right) induces an action of
$H$ on ${\cal A} = {\cal C}^{*}$ from the right
(from the left). This obviously follows from relations of the type (\ref{2.9}).

\subsection{\bf \em Bialgebras\label{hopf2}}
\setcounter{equation}0

So-called bialgebras are the next important objects that are used in the theory of quantum groups.


\newtheorem{def2}[def1]{Definition}
\begin{def2} \label{def2}
{\it An associative algebra
${\cal A}$ with identity that is
simultaneously a coassociative coalgebra with coidentity is called
a bialgebra if the algebraic and coalgebraic structures are self-consistent.
Namely, the comultiplication and coidentity must be homomorphisms of the algebras:
\be
\lb{2.12}
\Delta(e_{i})\, \Delta(e_{j}) = m_{ij}^{k} \Delta(e_{k}) \Rightarrow
\Delta^{i'i''}_{i} \Delta^{j'j''}_{j} m^{k'}_{i'j'} m^{k''}_{i''j''} =
m_{ij}^{k}\, \Delta^{k'k''}_{k} \; , \;\;
\ee
$$
\Delta(I) = I \otimes I \; , \;\;
\epsilon(e_{i}\cdot e_{j}) = \epsilon(e_{i})\,
\epsilon(e_{j}) \; , \;\;
\epsilon(I) = E^i \, \epsilon_i = 1 \; .
$$
}
\end{def2}
Note that for every bialgebra we have a certain freedom in
the definition of the multiplication (\ref{2.1}) and the comultiplication (\ref{2.6}).
Indeed, all the axioms (\ref{2.3}), (\ref{2.7}), and
(\ref{2.12}) are satisfied if instead of (\ref{2.1})  we take
$$
e_{i} \cdot e_{j} = m^{k}_{ji} \; e_{k},
$$
or instead of (\ref{2.6}) choose
\be
\lb{2.13}
\Delta^{\sf op}(e_{i}) =
\Delta^{jk}_{i} \, e_{k} \otimes e_{j} \; ,
\ee
(such algebras are denoted as ${\cal A}^{op}$ and ${\cal A}^{cop}$, respectively).
Then the algebra ${\cal A}$ is called noncommutative if
$m^{k}_{ij} \neq m^{k}_{ji}$, and noncocommutative if $\Delta^{ij}_{k} \neq
\Delta^{ji}_{k}$.

In quantum physics, it is usually assumed that all algebras
of observables are bialgebras. Indeed, a coalgebraic structure is
needed to define the action of the algebra ${\cal A}$ of observables on the
state $|\psi_{1}\rangle \otimes |\psi_{2}\rangle$ of the system that is the
composite system formed
from two independent systems with wave functions
$|\psi_{1} \rangle$ and $|\psi_{2} \rangle$
\be
\lb{phys}
a \triangleright (|\psi_{1}\rangle \otimes |\psi_{2}\rangle) =
\Delta(a) \, (|\psi_{1}\rangle \otimes |\psi_{2}\rangle) =
a_{(1)} \, |\psi_{1}\rangle \otimes a_{(2)} \, |\psi_{2}\rangle \;\;
(\forall a \in {\cal A}) \; .
\ee
In other words, for bialgebras it is possible to
formulate a theory of representations in which new representations
can be constructed by direct multiplication of old ones.

A classical example of a bialgebra is the universal enveloping algebra of a Lie algebra $\mathfrak{g}$,
in particular, the spin algebra $\mathfrak{su}(2)$
 in three-dimensional space.
To demonstrate this, we consider the Lie algebra
 $\mathfrak{g}$
with generators $J_{\alpha} \;\;
(\alpha = 1,2,3, \dots)$, that satisfy the antisymmetric multiplication rule
(defining relations)
\be
\lb{lie}
[J_{\alpha}, \, J_{\beta} ] =
t_{\alpha\beta}^{\gamma} J_{\gamma} .
\ee
Here $t^{\gamma}_{\alpha\beta} = - t^{\gamma}_{\beta\alpha}$ are structure
constants which satisfy Jacoby identity. The enveloping algebra of this
algebra is the algebra $U(\mathfrak{g})$ with basis elements consisting of the
identity $I$ and the elements
$e_{i} = J_{\alpha_{1}} \cdots J_{\alpha_{n}}
\; \forall n \geq 1$, where the products of the
generators $J$ are ordered lexicographically, i.e.,
$\alpha_1 \leq \alpha_2 \leq \dots \leq \alpha_n$.
The coalgebraic structure for the algebra $U(\mathfrak{g})$ is specified by means of the mappings
\be
\lb{2.14}
\Delta(J_{\alpha})
=J_{\alpha} \otimes I  + I \otimes J_{\alpha} \; , \;\;
\epsilon(J_{\alpha}) = 0 \; , \;\;
\epsilon(I) = 1 \; ,
\ee
which satisfy all the axioms of a bialgebra. The mapping $\Delta$ in (\ref{2.14})
is none other than the rule for addition of spins. In fact one can
quantize the coalgebraic structure (\ref{2.14}) for universal enveloping algebra $U(\mathfrak{g})$
and consider the noncocommutative comultiplications $\Delta$. Such quantizations
will be considered below in Section {\bf \ref{semicl}}
and leads to the definition of
Lie bialgebras.

Considering exponentials of elements of a Lie algebra, one can arrive
at the definition of a group bialgebra of the group $G$ with structure mappings
\be
\lb{2.15}
\Delta(h) = h \otimes h \; , \;\; \epsilon(h) =1 \;\;\;\;
(\forall \; h \in G),
\ee
which obviously follow from (\ref{2.14}). The next important example
of a bialgebra is the algebra ${\cal A}(G)$ of functions
$f$ on a group
($f: G \rightarrow \mathbb{C}$).
This algebra is dual to the group algebra of the group $G$, and its structure
mappings have the form ($f,f' \in {\cal A}(G); \;\; h,h' \in G$):
\be
\lb{2.15b}
(f \cdot f')(h) = f(h)f'(h) \, , \;\;\;\;
f(h \cdot h') = (\Delta(f)) (h,h') = f_{(1)}(h)\; f_{(2)}(h') \, , \;\;\;\;
\epsilon(f) = f(I) \; ,
\ee
where $I_G$ is the identity element in the group $G$.
In particular, if the functions $T^{i}_{j}$ realize a matrix representation of the
group $G$, then we have
\be
\lb{2.15a}
T^{i}_{j}(hh') = T^{i}_{k}(h)T^{k}_{j}(h') \;\; \Rightarrow
 \;\; \Delta(T^{i}_{j}) =  T^{i}_{k}  \otimes T^{k}_{j} \; ,
\ee
(the functions $T^{i}_{j}$ can be regarded as generators
of a subalgebra in the algebra ${\cal A}(G)$).
Note that if $\mathfrak{g}$ is non-Abelian, then
 $U(\mathfrak{g})$ and $G$ are noncommutative but
cocommutative bialgebras, whereas ${\cal A}(G)$ is a commutative but noncocommutative bialgebra.
Anticipating, we mention that the most
 interesting quantum groups are
associated with noncommutative and noncocommutative bialgebras.

It is obvious that for a bialgebra ${\cal H}$ it is also possible to introduce
the concepts of left (co)modules and (co)module (co)algebras
 (right (co)modules and (co)module (co)algebras are
 introduced in exactly
the same way). Moreover, for the bialgebra ${\cal H}$ it is possible to introduce the
concept of a left (right) bimodule $B$, i.e., a left (right) ${\cal H}$-module that is
simultaneously a left (right) ${\cal H}$-comodule; at the same time, the module and
comodule structures must be self-consistent:
$$
\Delta_{L} ({\cal H} \triangleright B) =
\Delta({\cal H}) \triangleright \Delta_{L}(B) \; ,
$$
$$
(\epsilon \otimes id) \Delta_{L} (b) = b \; , \;\; b \in B \; ,
$$
where $\Delta_{L} (b) = \bar{b}^{(-1)} \otimes b^{(0)}$
and $\bar{b}^{(-1)} \in {\cal H}$, $b^{(0)} \in B$.
On the other hand, in the case of bialgebras the conditions of conserving
of the (co)algebraic structure of (co)modules can be represented in a
more explicit form. For example, for the left ${\cal H}$-module algebra
${\cal A}$ we have
$(a,b \in {\cal A}; \;\; h \in {\cal H})$:
$$
h \triangleright (ab) =
(h_{(1)} \triangleright a)
(h_{(2)} \triangleright b) \; , \;\; h \triangleright I_{A} =
\epsilon(h) I_{A} \; .
$$
In addition, for the left ${\cal H}$-module coalgebra ${\cal A}$ we must have
$$
\Delta(h \triangleright a) = \Delta(h) \triangleright \Delta(a) =
(h_{(1)} \triangleright a_{(1)}) \otimes
(h_{(2)} \triangleright a_{(2)}) \; , \;\;
\epsilon(h \triangleright a) = \epsilon(h)\epsilon(a) \; .
$$
Similarly, the algebra ${\cal A}$ is a left ${\cal H}$-comodule algebra if
$$
\Delta_{L}(ab) = \Delta_{L}(a)\, \Delta_{L}(b) \; , \;\;
\Delta_{L}(I_{A}) = I_{\cal H} \otimes I_{A} \; ,
$$
and, finally, the coalgebra ${\cal A}$ is a left  ${\cal H}$-comodule coalgebra if
\be
\lb{2.16}
(id \otimes \Delta)\Delta_{L}(a) =
m_{\cal H}(\Delta_{L} \otimes \Delta_{L})\Delta(a) \; , \;\;
(id \otimes \epsilon_{A})\Delta_{L}(a) =
I_{\cal H}\epsilon_{A}(a) \; ,
\ee
where
$$
m_{\cal H}(\Delta_{L}\otimes\Delta_{L})(a \otimes b) =
\bar{a}^{(-1)} \bar{b}^{(-1)}\otimes a^{(0)} \otimes b^{(0)} \; .
$$

We now consider the bialgebra ${\cal H}$, which acts on a certain module algebra ${\cal A}$. One
further important property of bialgebras is that we can define a new associative algebra ${\cal A}
\sharp {\cal H}$ as the cross product (smash product) of ${\cal A}$ and ${\cal H}$. Namely:


\newtheorem{def3}[def1]{Definition}
\begin{def3} \label{def3}
{\it The left smash product ${\cal A} \sharp {\cal H}$ of the bialgebra ${\cal H}$ and its left
module algebra ${\cal A}$ is an associative
algebra such that: \\
1) As a vector space,
${\cal A} \sharp {\cal H}$ is identical to ${\cal A} \otimes {\cal H}$ \\
2) The product is defined in the sense ($h,g \in {\cal H}$; $a,b \in {\cal A}$)
\be
\lb{2.17}
(a \sharp g)\, (b \sharp h) = \sum_{g} a(g_{(1)} \triangleright b) \sharp
(g_{(2)} h) \equiv (a \sharp I)\;
(\Delta(g) \triangleright ( b \sharp h)) \; ;
\ee
3) The identity element is $I \sharp I$.}
\end{def3}
If the algebra ${\cal A}$ is the bialgebra dual to the bialgebra ${\cal H}$, then the relations
(\ref{2.17}) and (\ref{2.11}) define the rules for interchanging the elements $(I \sharp g)$ and
$(a \sharp I)$:
\be
\lb{2.18}
(I \sharp g) \, (a \sharp I) =
(a_{(1)} \sharp I) \, \langle g_{(1)}|a_{(2)}\rangle \,
 (I \sharp g_{(2)}) \; .
\ee
Thus, the subalgebras ${\cal A}$ and ${\cal H}$ in
${\cal A} \, \sharp \, {\cal H}$ do not commute with each other.
The smash product depends on which action
(left or right) of the algebra ${\cal H}$ on ${\cal A}$ we choose.
In addition, the smash product generalizes the concept of the semidirect product.
In particular, if we take as bialgebra ${\cal H}$
 the Lorentz group algebra
 (see (\ref{2.15})),
and as module ${\cal A}$ the group of translations in Minkowski space,
then the smash product ${\cal A} \, \sharp \, {\cal H}$ defines the structure of the Poincare group.

The coanalog of the smash product, the smash coproduct
${\cal A} \, \underline{\sharp} \, {\cal H}$, can also be defined.
For this, we consider the bialgebra
${\cal H}$ and its comodule coalgebra ${\cal A}$.
Then on the space ${\cal A} \otimes {\cal H}$ it is
possible to define the structure of a coassociative coalgebra:
\be
\lb{2.19}
\Delta (a \, \underline{\sharp} \, h) =
(a_{(1)} \, \underline{\sharp} \,
 { \bar{a}_{(2)} }^{\;\;(-1)} h_{(1)}) \otimes
(a_{(2)}^{\;\;(0)} \, \underline{\sharp} \,  h_{(2)}) \; , \;\;
\epsilon(a \, \underline{\sharp} \,  h) =
\epsilon(a)\epsilon(h) \; .
\ee
The proof of the coassociativity reduces to verification of the identity
$$
(m_{\cal H}( \Delta_{L}\otimes \Delta_{\cal H}) \otimes id)
(id \otimes \Delta_{L})\Delta_{\cal A}(a) =
(id \otimes id \otimes \Delta_{L})
(id \otimes \Delta_{\cal A} )\Delta_{L}(a) \; ,
$$
which is satisfied if we take into account the axiom (\ref{2.16}) and the comodule axiom
\be
\lb{2.20}
(id \otimes \Delta_{L})\Delta_{L}(a) =
(\Delta_{\cal H} \otimes id)\Delta_{L}(a) \; .
\ee

Note that from the two bialgebras ${\cal A}$ and ${\cal H}$, which act and coact
on each other in a special manner, it is possible to organize a
new bialgebra that is simultaneously the smash product and smash
coproduct of ${\cal A}$ and ${\cal H}$ (bicross product; see Ref. \cite{14}).

\subsection{\bf \em Hopf algebras. Universal ${\cal R}$-matrices\label{hopf3}}
\setcounter{equation}0

We can now introduce the main concept in the theory of quantum groups,
namely, the concept of the Hopf algebra.


\newtheorem{def4}[def1]{Definition}
\begin{def4} \label{def4}
{\it A bialgebra ${\cal A}$ equipped with an additional mapping $S: \;\; {\cal A} \rightarrow {\cal
A}$ such that
\be
\lb{2.21}
\begin{array}{c}
m(S \otimes id) \Delta = m(id \otimes S) \Delta =
\hbox{\bf i} \cdot \epsilon  \Rightarrow  \\
S(a_{(1)}) \, a_{(2)}  =
a_{(1)} \, S(a_{(2)}) = \epsilon(a) \cdot I \;\; (\forall a \in {\cal A})
\end{array}
\ee
is called a Hopf algebra. The mapping $S$ is called the antipode and is an
antihomomorphism with respect to both multiplication and comultiplication:
\be
\lb{2.22}
S(ab) = S(b)S(a) \; , \;\; (S \otimes S)\Delta(a) = \sigma \cdot
\Delta(S(a)) \; ,
\ee
where $a,b \in {\cal A}$ and $\sigma$ denotes the operator of transposition, $\sigma (a \otimes b)
= (b \otimes a)$.}
\end{def4}
If we set
\be
\lb{2.23}
S(e_{i}) = S^{j}_{i}e_{j} \; ,
\ee
then the axiom (\ref{2.21}) can be rewritten in the form
\be
\lb{2.24}
\Delta^{ij}_{k}S^{n}_{i}m^{l}_{nj} =
\Delta^{ij}_{k}S^{n}_{j}m^{l}_{in} = \epsilon_{k} E^{l} \; .
\ee
From the axioms for the structure mappings of a Hopf algebra,
it is possible to obtain the useful equations
\be
\lb{2.25}
\begin{array}{c}
S^{i}_{j}\epsilon_{i} = \epsilon_{j} \; , \;\;
S^{i}_{j}E^{j} = E^{i} \; , \\  \\
\Delta^{ji}_{k}(S^{-1})^{n}_{i}m^{l}_{nj} =
\Delta^{ji}_{k}(S^{-1})^{n}_{j}m^{l}_{in} = \epsilon_{k} E^{l} \; ,
\end{array}
\ee
which we shall use in what follows. Note that, in general, the antipode $S$
is not necessarily invertible. An invertible antipode is called bijective.

In quantum physics the existence of the antipode $S$ is needed to define a space of
contragredient states $\langle \psi |$ (contragredient module of ${\cal A}$)
with pairing $\langle \psi | \phi \rangle$:
$\langle \psi | \otimes | \phi \rangle \to \mathbb{C}$.
Left actions of the Hopf algebra ${\cal A}$
of observables to the contragredient states are (cf. the actions (\ref{phys}) of ${\cal A}$ to the
states $| \psi_1 \rangle \otimes | \psi_2 \rangle$):
\be
\lb{contg}
a \triangleright \langle \psi | : =
\langle \psi | \, S(a) \;\; (a \in {\cal A}) \; ,
\ee
$$
a \triangleright (\langle \psi_1 | \otimes \langle \psi_2 |) : =
(\langle \psi_1 | \otimes \langle \psi_2 |) \Delta(S(a)) =
\langle \psi_1 | \, S(a_{(2)}) \otimes \langle \psi_2 | \, S(a_{(1)}) \; .
$$
The states $\langle \psi |$ are called left dual to the states
$| \phi \rangle$; the right dual ones are introduced with the help of the inverse
antipode $S^{-1}$ (see e.g. \cite{17}).
Then, the covariance of the pairing $\langle \psi | \phi \rangle$
under the left action of ${\cal A}$ can be established:
$$
\begin{array}{c}
a \triangleright \langle \psi | \phi \rangle \equiv
(a_{(1)} \triangleright \langle \psi |)\; (a_{(2)}
 \triangleright |\phi \rangle) =
\langle \psi | S(a_{(1)}) \, a_{(2)} | \phi \rangle = \epsilon(a) \langle \psi | \phi \rangle
\; , \\ [0.2cm]
a \triangleright \langle \psi_1 | \phi_1 \rangle \;
\langle \psi_2 | \phi_2 \rangle = a \triangleright
(\langle \psi_1 | \otimes \langle \psi_2 |)
(| \phi_1 \rangle \otimes | \phi_2 \rangle) = \\ [0.2cm]
= \langle \psi_1 | S(a_{(2)})a_{(3)}| \phi_1 \rangle
\langle \psi_1 | S(a_{(1)})a_{(4)}| \phi_1 \rangle =
\epsilon(a) \, \langle \psi_1 | \phi_1 \rangle \;
\langle \psi_2 | \phi_2 \rangle \; .
\end{array}
$$

The universal enveloping algebra $U(\mathfrak{g})$
 and the group bialgebra of the group $G$
that we considered above can again serve as examples of cocommutative
Hopf algebras. An example of a commutative
 (but noncocommutative)
Hopf algebra is the bialgebra ${\cal A}(G)$,
which we also considered above. The antipodes for these algebras have the form
\begin{equation}
\label{2.25a}
\begin{array}{c}
U(\mathfrak{g}): \;\; S(J_{\alpha}) = - J_{\alpha} \; , \;\;
S(I) = I \; ,
\\ \\
G: \;\; S(h) = h^{-1} \; ,
\\ \\
{\cal A}(G): \;\; S(f)(h) = f(h^{-1}) \; ,
\end{array}
\end{equation}
and satisfy the relation $S^{2} =id$,
which holds for all commutative or cocommutative Hopf algebras.

From the point of view of the axiom (\ref{2.21}), $S(a)$ looks like the
inverse of the element $a$, although in the general case $S^{2} \neq id$.
We recall that if a set ${\cal G}$
of elements with associative
multiplication ${\cal G} \otimes {\cal G} \rightarrow {\cal G}$
and with identity (semigroup) also contains all the
inverse elements, then such a set ${\cal G}$ becomes a group. Thus, from the point of
view of the presence of the mapping $S$, a Hopf algebra  generalizes the notion
of the group algebra (for which $S(h) = h^{-1}$), although by itself it obviously need
not be a group algebra. In accordance with Drinfeld's definition [13] the concepts
of a Hopf algebra and a quantum group are more or less equivalent. Of course, the most
interesting examples of quantum groups arise when one considers
noncommutative and noncocommutative Hopf algebras.
\vspace{0.2cm}

Consider a noncommutative Hopf algebra ${\cal A}$ which is also noncocommutative
$\Delta \neq \Delta^{\sf op} \equiv \sigma \Delta$,
where $\sigma$ is the transposition operator
$\sigma(a \otimes b) = b \otimes a$ ($\forall a,b \in {\cal A}$).
\newtheorem{def5}[def1]{Definition}
\begin{def5} \label{def5}
{\it A Hopf algebra ${\cal A}$ for which there exists an invertible element
${\cal R} \in {\cal A} \otimes
{\cal A}$ such that
\be
\lb{2.26}
\Delta^{\sf op}(a) = {\cal R} \Delta(a) {\cal R}^{-1} \; ,
\;\;\;\;\; \forall a \in {\cal A} \; ,
\ee
\be
\lb{2.27}
(\Delta \otimes id) ({\cal R}) = {\cal R}_{13} {\cal R}_{23} \; , \; \;
(id \otimes \Delta) ({\cal R}) = {\cal R}_{13} {\cal R}_{12}
\ee
is called quasitriangular. Here the element
\be
\lb{2.28}
{\cal R}=\sum_{ij} R^{(ij)} e_{i} \otimes e_{j}
\ee
is called the universal ${\cal R}$ matrix,
$R^{(ij)} \in \mathbb{C}$ are
the constants and the symbols ${\cal R}_{12}, {\cal R}_{13},
{\cal R}_{23}$ have the meaning
\be
\lb{2.30a}
{\cal R}_{12} = \sum_{ij} R^{(ij)} e_{i} \otimes e_{j}
\otimes I \; , \;\;
{\cal R}_{13} = \sum_{ij} R^{(ij)} e_{i} \otimes I \otimes e_{j}  \; , \;\;
{\cal R}_{23} =
\sum_{ij} R^{(ij)} I \otimes e_{i} \otimes e_{j} \; .
\ee}
\end{def5}
The relation (\ref{2.26}) shows that the noncocommutativity in a quasitriangular
Hopf algebra is kept "under control." It can be shown \cite{13'} that for such
a Hopf algebra the universal ${\cal R}$ matrix (\ref{2.28}) satisfies
the Yang-Baxter equation
\be
\lb{2.30}
{\cal R}_{12}{\cal R}_{13}{\cal R}_{23} = {\cal R}_{23}{\cal R}_{13}{\cal R}_{12} \; ,
\ee
(to which a considerable part of the review will be devoted)
and the relations
\be
\lb{2.29}
(id \otimes \epsilon){\cal R} = (\epsilon \otimes id){\cal R} = I \; ,\;\;
\ee
\be
\lb{2.29a}
\begin{array}{c}
(S \otimes id){\cal R} = {\cal R}^{-1} \;\; \Leftrightarrow \;\;
(S^{-1} \otimes id){\cal R}^{-1} = {\cal R} \; , \\ [0.2cm]
(id \otimes S){\cal R}^{-1} = {\cal R}
\;\; \Leftrightarrow \;\;
(id \otimes S^{-1}){\cal R} = {\cal R}^{-1}  \; .
\end{array}
\ee
The Yang-Baxter eq. (\ref{2.30}) follows from
(\ref{2.26}) and (\ref{2.27}):
\be
\lb{2.31}
 {\cal R}_{12} \, {\cal R}_{13} \, {\cal R}_{23} =
 {\cal R}_{12} \; (\Delta \otimes id)({\cal R}) =
   (\Delta^{\sf op} \otimes id)({\cal R})\; {\cal R}_{12} =
   {\cal R}_{23} \, {\cal R}_{13} \, {\cal R}_{12} \; .
\ee
It is easy to derive the relations (\ref{2.29}) by applying $(\epsilon \otimes id \otimes id)$
and
$(id \otimes id \otimes \epsilon)$ respectively to
the first and second relation in (\ref{2.27}),
and then taking into account (\ref{2.8}).
Next, we prove the equalities in (\ref{2.29a}).
 We consider expressions
${\cal R} \cdot (S \otimes id) \, {\cal R}$ and
${\cal R} \cdot (id \otimes S^{-1}) \, {\cal R}$
  and make use of the
Hopf algebra axioms (\ref{2.21})
and equations (\ref{2.27}), (\ref{2.29}):
$$
\begin{array}{c}
{\cal R}_{23} \cdot (id \otimes S \otimes id) \, {\cal R}_{23} =
(m_{12} \otimes id_3) \bigl( {\cal R}_{13} \,
(id \otimes S \otimes id) {\cal R}_{23} \bigr) = \\ [0.1cm]
= (m_{12} \otimes id_3)\,
(id \otimes S \otimes id) {\cal R}_{13} \, {\cal R}_{23} =
(m_{12} \otimes id_3)\,
\bigl( (id \otimes S)\Delta \otimes id \bigr){\cal R} =
 (\hbox{\bf i} \cdot \epsilon \otimes id) \, {\cal R} = I \; ,
 \end{array}
$$
$$
\begin{array}{c}
{\cal R}_{12} \, (id \otimes S^{-1}) \, {\cal R}_{12} =
 (id_1  \otimes   m_{23})\,
(id \otimes id \otimes S^{-1} ) {\cal R}_{12} \, {\cal R}_{13} =
\\ [0.2cm]
= (id_1  \otimes   m_{23})\,
\bigl(id \otimes  (id \otimes S^{-1})\Delta^{\sf op} \bigr){\cal R} =
 (id \otimes \hbox{\bf i} \cdot \epsilon) \, {\cal R} = I \; ,
 \end{array}
$$
where the ultimate equality follows from (\ref{2.21})
which is written in the form
$a_{(2)} S^{-1}(a_{(1)}) = \epsilon(a) I$.

The next important concept that we shall need in what follows
is the concept of the Hopf algebra ${\cal A}^{*}$ that is the dual of the
Hopf algebra ${\cal A}$.
We choose in ${\cal A}^{*}$ basis elements $\{ e^{i} \}$ and define multiplication, the identity,
comultiplication, the coidentity, and the antipode for ${\cal A}^{*}$ in the form
\be
\lb{2.33}
 e^{i} e^{j}  = m^{ij}_{k} e^{k} \; , \;\; I = \bar{E}_{i}e^{i} \; , \;\;
 \Delta(e^{i}) = \Delta^{i}_{jk} e^{j} \otimes e^{k} \; , \;\;
\epsilon(e^{i}) = \bar{\epsilon}^{i} \; , \;\;
S(e^{i}) = \bar{S}^{i}_{j} e^{j} \; .
\ee


\newtheorem{def6}[def1]{Definition}
\begin{def6} \label{def6}
{\it Two Hopf algebras ${\cal A}$ and ${\cal A}^{*}$ with corresponding bases
 $\{e_{i} \}$ and $\{e^{i} \}$ are
 said to be dual to each other if there exists a
 nondegenerate pairing $\langle  . | . \rangle $:
${\cal A}^{*} \otimes {\cal A} \rightarrow \mathbb{C}$ such that
\be
\lb{2.34}
\begin{array}{c}
\langle  e^{i}e^{j} | e_{k}\rangle  \equiv \langle  e^{i} \otimes e^{j} | \Delta(e_{k})\rangle  =
\langle  e^{i} | e_{k'}\rangle  \Delta_{k}^{k'k''} \langle  e^{j} | e_{k''}\rangle  \\ \\
\langle  e^{i} | e_{j} e_{k}\rangle  \equiv \langle  \Delta(e^{i})| e_{j} \otimes e_{k}\rangle  =
\langle  e^{i'} | e_{j}\rangle  \Delta^{i}_{i'i''} \langle  e^{i''} | e_{k}\rangle  \\ \\
\langle S(e^{i}) | e_{j} \rangle  = \langle e^{i} | S(e_{j}) \rangle  \; , \;\;
\langle  e^{i} | I \rangle  = \epsilon (e^{i}) \; , \;\;
\langle  I | e_{i} \rangle  = \epsilon (e_{i}) \; .
\end{array}
\ee}
\end{def6}
Since the pairing $\langle .|.\rangle $ (\ref{2.34}) is nondegenerate,
we can always choose basis elements $\{ e^{i} \}$ such that
\be
\lb{2.35}
\langle  e^{i} | e_{j} \rangle  = \delta^{i}_{j} \; .
\ee
Then from the axioms for the pairing (\ref{2.34}) and from the
definitions of the structure maps (\ref{2.1}), (\ref{2.23}), and (\ref{2.33}) in
${\cal A}$ and ${\cal A}^{*}$ we readily deduce
\be
\lb{2.36}
m^{ij}_{k} = \Delta^{ij}_{k} \; , \;\;
m_{ij}^{k} = \Delta_{ij}^{k} \; , \;\;
\bar{S}^{i}_{j} = S^{i}_{j} \; , \;\;
\bar{\epsilon}^{i} = E^{i} \; , \;\;
\bar{E}_{i} = \epsilon_{i} \; .
\ee
Thus, the multiplication, identity, comultiplication, coidentity,
and antipode in a Hopf algebra define, respectively, comultiplication,
coidentity, multiplication, identity, and antipode in the dual Hopf algebra.

\vspace{0.3cm}


\noindent
{\bf Remark.} In \cite{Pontr} L.S. Pontryagin showed that the set of characters of an abelian locally compact
group $G$ is an abelian group, called the dual group $G^*$ of $G$. The group $G^*$ is also locally
compact. Moreover, the dual group of $G^*$ is isomorphic to $G$. This beautiful theory becomes
wrong if $G$ is a noncommutative group, even if it is finite. To restore the duality principle one
can replace the set of characters for a finite noncommutative group $G$ by the category of its
irreducible representations (irreducible representations for the commutative groups are exactly
characters). Indeed, T. Tannaka and M. Krein showed that the compact group $G$ can be recovered
from the set of its irreducible unitary representations. They proved a duality theorem for compact
groups, involving irreducible representations of $G$ (although no group-like structure is to be put
on that class, since the tensor product of two irreducible representation may no longer be
irreducible). However, the tensor product of two irreducible representations of the compact group
$G$ can be expanded as a
sum of irreducible representations and, thus, the dual object has the structure of an algebra. Recall (see (\ref{2.15a})) that
matrix representations of group $G$ are realized by the sets
of special functions $T^i_{\; k}$.
One can consider the group algebra ${\cal G}$ of finite group $G$
and the algebra ${\cal A}(G) \equiv {\cal G}^\star$ of functions on the
group $G$ as simplest
examples of the Hopf algebras. The structure mappings for these algebras
have been defined in (\ref{2.15}), (\ref{2.15b}) and (\ref{2.25a}).
Note that the algebras ${\cal G}$ and ${\cal G}^\star$
are Hopf dual to each other.
The detailed structure of ${\cal G}^\star$ follows
from the representation theory of finite groups
(see e.g. \cite{Sierre}).

\subsection{\bf \em Heisenberg and Quantum doubles.
 Yetter--Drinfeld modules\label{hopf4}}
\setcounter{equation}0

In Subsection 2.2 we have defined (see Definition \ref{def3}) the notion of the smash
(cross) product of the bialgebra and its module algebra.
Since the Hopf dual algebra ${\cal A}^*$ is the natural right and left module
algebra for the Hopf algebra ${\cal A}$ (\ref{2.10}), (\ref{2.11}),
one can immediately define the right
${\cal A}^* \sharp {\cal A}$ and the left ${\cal A} \sharp {\cal A}^*$
cross products of the algebra ${\cal A}$ on ${\cal A}^*$.
These cross-product algebras are called
 Heisenberg doubles of ${\cal A}$
and they are the associative algebras with
nontrivial cross-multiplication rules (cf. eq. (\ref{2.18})):
\begin{equation}
\label{25d}
a \, \bar{a} = (a_{(1)} \triangleright \bar{a}) \,
a_{(2)} = \bar{a}_{(1)} \, \langle a_{(1)} \, | \, \bar{a}_{(2)} \rangle \, a_{(2)}  \; ,
\end{equation}
\begin{equation}
\label{25dd}
\bar{a} \, a = a_{(1)} (\bar{a} \triangleleft a_{(2)}) \,
 = a_{(1)} \, \langle \bar{a}_{(1)} \, | \, a_{(2)} \rangle \, \bar{a}_{(2)}  \; ,
\end{equation}
where $a \in {\cal A}$ and $\bar{a} \in {\cal A}^*$.
Here we discuss only the left cross product algebra ${\cal A} \sharp {\cal A}^*$
(\ref{25d}) (the other one (\ref{25dd})
is considered analogously).

As in the previous subsection, we denote $\{ e^{i} \}$ and $\{ e_{i} \}$
the dual basis elements of ${\cal A}^*$ and ${\cal A}$, respectively.
In terms of this basis we rewrite (\ref{25d}) in the form
\begin{equation}
\label{25f}
e_r \, e^n =  e^i \, \Delta^n_{if} \langle e_j | e^f \rangle  \, \Delta^{jk}_r  \, e_k
= m^n_{ij} \,  e^i \, e_k \, \Delta^{jk}_r  \; .
\end{equation}
 Let us define a right ${\cal A}^*$ - coaction and
a left ${\cal A}$ - coaction
on the algebra ${\cal A} \sharp {\cal A}^*$, such that these coactions
respect the algebra structure of ${\cal A} \sharp {\cal A}^*$:
\begin{equation}
\label{26} \Delta_R( z) = C \, (z \otimes 1) \, C^{-1} \; , \;\;\;
\Delta_L( z) = C^{-1} \, (1 \otimes z) \, C \; , \;\;\;
C \equiv e_{i} \otimes e^{i} \; .
\end{equation}
The inverse of the canonical element $C$ is
$$
C^{-1} = S(e_{i}) \otimes e^{i} = e_{i} \otimes S(e^{i}) \; ,
$$
and $\Delta_R$, $\Delta_L$ (\ref{26}) are represented in the form
\begin{equation}
\label{27}
\Delta_R( z) =
(e_{k(1)} \, z \, S(e_{k(2)})) \otimes e^{k} \; , \;\;\;
\Delta_L( z) =
e_{k} \otimes  S(e^k_{(1)}) \, z \, e^k_{(2)}
\; .
\end{equation}
Note that $\Delta_R(\bar{z}) = \Delta(\bar{z})$ $\forall \bar{z} \in {\cal A}^*$
and $\Delta_L(z) = \Delta(z)$ $\forall z \in {\cal A}$
(here ${\cal A}$ and ${\cal A}^*$ are understood as the Hopf subalgebras
in ${\cal A} \sharp {\cal A}^*$ and $\Delta$ are corresponding
comultiplications). Indeed, for $z \in {\cal A}$ we have
$$
\Delta_L(z) = e_{k} \otimes  S(e^k_{(1)}) \, z \, e^k_{(2)} =
e_{k} \otimes  S(e^k_{(1)}) \, e^k_{(2)} \, \langle z_{(1)} \, | \, e^k_{(3)} \rangle
z_{(2)} =
$$
$$
= e_{k} \,  \langle z_{(1)} \, | \, e^k \rangle \otimes
z_{(2)} = z_{(1)} \otimes  z_{(2)} \; ,
$$
(the proof of $\Delta_R(\bar{z}) = \Delta(\bar{z})$ is similar).
The axioms
$$
(id \otimes \Delta)\Delta_R = (\Delta_R \otimes id)\Delta_R
\;\; , \;\;\;
(id \otimes \Delta_L)\Delta_L = (\Delta \otimes id)\Delta_L \; ,
$$
$$
(id \otimes \Delta_R)\Delta_L(z) = C^{-1}_{13} \,
(\Delta_L \otimes id)\Delta_R(z) \, C_{13} \; ,
$$
can be verified directly by using relations (cf. (\ref{2.27}))
$$
(id \otimes \Delta) C_{12} = C_{13} \, C_{23} \; , \;\;\;
(\Delta \otimes id) C_{12} = C_{13} \, C_{23} \; ,
$$
and the pentagon identity \cite{BSk} for $C$
\be
\lb{c5}
C_{12} \, C_{13} \, C_{23} = C_{23} \, C_{12} \; .
\ee
The proof of (\ref{c5}) is straightforward (see (\ref{25f})):
$$
C_{12} \, C_{13} \, C_{23} = e_i \, e_j \otimes e^i \, e_k \otimes e^j \, e^k =
e_n \otimes  m^n_{ij} \,  e^i \, e_k \, \Delta^{jk}_r \otimes e^r =
$$
$$
= e_n \otimes  e_r \, e^n  \otimes e^r = C_{23} \, C_{12} \; .
$$
The pentagon identity (\ref{c5}) is used for the construction
of the explicit solutions of the tetrahedron equations
(3D generalizations of the Yang-Baxter equations).

Although ${\cal A}$ and ${\cal A}^*$ are Hopf algebras, their Heisenberg
doubles  ${\cal A} \sharp {\cal A}^*$, ${\cal A}^* \sharp {\cal A}$ are not Hopf algebras. But
as we have seen just before the algebra ${\cal A} \sharp {\cal A}^*$
(as well as ${\cal A}^* \sharp {\cal A}$) still possesses some
covariance properties, since the coactions (\ref{26}) are
covariant transformations (homomorphisms) of the algebra
${\cal A} \sharp {\cal A}^*$.

The natural question is the following: is it possible to invent such
a cross-product of the Hopf algebra and its dual Hopf algebra to
obtain a new Hopf algebra?
Drinfeld \cite{13} showed that there exists a quasitriangular
Hopf algebra ${\cal D}({\cal A})$ that is a special smash product of the Hopf
algebras ${\cal A}$ and ${\cal A}^{o}$:
${\cal D}({\cal A}) = {\cal A} \Join {\cal A}^{o}$, which is called the quantum double.
Here we denote by ${\cal A}^{o}$ the algebra ${\cal A}^{*}$ with opposite
comultiplication: $\Delta(e^{i}) = m^{i}_{kj} e^{j} \otimes e^{k}$,
${\cal A}^o = ({\cal A}^*)^{cop}$.
It follows
from (\ref{2.25}) that the antipode for ${\cal A}^{o}$ will be not $S$ but
the skew antipode $S^{-1}$. Thus, the structure mappings for ${\cal A}^{o}$ have the form
\be
\lb{2.37}
e^{i}e^{j} = \Delta^{ij}_{k} e^{k} \; , \; \;
\Delta(e^{i}) = m^{i}_{kj} e^{j} \otimes e^{k} \; , \; \;
S(e^{i}) = (S^{-1})^{i}_{j} e^{j} \; .
\ee
The algebras ${\cal A}$ and ${\cal A}^{o}$ are said to be antidual,
and for them we can introduce the antidual
pairing $\langle \langle .|.\rangle \rangle $:
${\cal A}^{o} \otimes {\cal A} \rightarrow \mathbb{C}$, which satisfies the conditions
\be
\lb{2.38}
\begin{array}{c}
\langle \langle  e^{i}e^{j} | e_{k}\rangle \rangle  \equiv \langle \langle  e^{i} \otimes e^{j} | \Delta(e_{k})\rangle \rangle  =
 \Delta_{k}^{ij}  \; , \\ \\
\langle \langle  e^{i} | e_{k} e_{j}\rangle \rangle  \equiv \langle \langle  \Delta(e^{i})| e_{j} \otimes e_{k}\rangle \rangle  =
m^{i}_{kj}  \; , \\ \\
\langle \langle S(e^{i}) | e_{j} \rangle \rangle  = \langle \langle e^{i} | S^{-1}(e_{j}) \rangle \rangle
= (S^{-1})^{i}_{j} \; , \;\; \\ \\
\langle \langle e^{i} | S(e_{j}) \rangle \rangle  = \langle \langle S^{-1}(e^{i}) | e_{j}\rangle \rangle  = S^{i}_{j}
\; , \;\; \\ \\
\langle \langle  e^{i} | I \rangle \rangle  = E^{i} \; , \;\;
\langle \langle  I | e_{i} \rangle \rangle  = \epsilon_{i} \; .
\end{array}
\ee
The universal $R$-matrix can be expressed in the form
of the canonical element
\be
\lb{2.39}
{\cal R} = (e_{i} \Join I) \otimes (I \Join e^{i}) \; ,
\ee
and the multiplication in ${\cal D}({\cal A})$
is defined in accordance with (the summation signs are omitted)
\be
\lb{2.40}
(a \Join \alpha)(b \Join \beta) =
a \left( (\alpha_{(3)} \triangleright b) \triangleleft
S(\alpha_{(1)}) \right) \Join \alpha_{(2)} \beta \; ,
\ee
where $\alpha, \beta \in {\cal A}^{o}$,
$a,b \in {\cal A}$, $\Delta^{2}(\alpha) =
\alpha_{(1)} \otimes \alpha_{(2)} \otimes \alpha_{(3)}$ and
\be
\lb{2.41}
\alpha \triangleright b = b_{(1)} \langle \langle \alpha | b_{(2)} \rangle \rangle  \; , \;\;
b \triangleleft \alpha = \langle \langle \alpha | b_{(1)} \rangle \rangle  b_{(2)}\; .
\ee

The coalgebraic structure on the quantum double is defined by
the direct product of the coalgebraic structures on the Hopf algebras
${\cal A}$ and ${\cal A}^{o}$:
\be
\lb{2.42}
\Delta(e_{i} \Join e^{j}) =
\Delta(e_{i} \Join I) \Delta(I \Join e^{j}) =
\Delta^{nk}_{i} m^{j}_{lp}(e_{n} \Join e^{p})
\otimes (e_{k} \Join e^{l}) \; .
\ee
Finally, the antipode and coidentity for ${\cal D}({\cal A})$ have the form
\be
\lb{2.43}
S(a \Join \alpha) = S(a) \Join S(\alpha) \; , \;\;
\epsilon(a \Join \alpha) = \epsilon(a) \epsilon(\alpha)  \; .
\ee

All the axioms of a Hopf algebra can be verified for ${\cal D}({\cal A})$
by direct calculation. A simple proof of the associativity of the
multiplication (\ref{2.40}) and the coassociativity of the
comultiplication (\ref{2.42}) can be found in Ref. \cite{15}.

Taking into account (\ref{2.41}), we can rewrite (\ref{2.40}) as the
commutator for the elements $(I \Join \alpha)$ and $(b \Join I)$:
\be
\lb{2.44a}
(I \Join \alpha)(b \Join I) = \langle \langle  S( \alpha_{(1)}) | b_{(1)} \rangle \rangle
(b_{(2)} \Join I) (I \Join \alpha_{(2)})
\langle \langle   \alpha_{(3)} | b_{(3)} \rangle \rangle
\ee
or, in terms of the basis elements
$\alpha = e^{t}$ and
$b = e_{s}$ we have \cite{13}
\be
\lb{2.44}
\begin{array}{c}
(I \Join e^{t})(e_{s} \Join I) =  m^{t}_{klp} \Delta_{s}^{njk}
(S^{-1})^{p}_{n} (e_{j} \Join I)(I \Join e^{l}) \equiv \\ \\
\left( m^{t}_{ip}(S^{-1})^{p}_{n} \Delta^{nr}_{s} \right)
\left( m^{i}_{kl} \Delta^{jk}_{r} \right)
(e_{j} \Join I)(I \Join e^{l}) \; ,
\end{array}
\ee
where $m^{t}_{klp}$ and $\Delta^{njk}_{s}$ are defined in (\ref{2.3}) and (\ref{2.7}),
and $(S^{-1})^{p}_{n}$  is the matrix of the skew antipode.

The consistence of definitions of left and right bimodules over the quantum double ${\cal D}({\cal A})$
should be clarified in view of the nontrivial
structure of the cross-multiplication rule (\ref{2.44a}), (\ref{2.44})
for subalgebras ${\cal A}$ and ${\cal A}^o$. It can be done (see, e.g. \cite{Rosso2}) if one considers
left or right coinvariant bimodules (Hopf modules): $M^L =  \{ m: \;\; \Delta_L(m) = 1 \otimes m  \}$ or
$M^R =  \{ m: \;\; \Delta_R(m) = m \otimes 1 \}$. E.g., for $M^R$ one can define the left ${\cal A}$
and left ${\cal A}^o$-module actions as
\be
\lb{2.44b}
a \triangleright m = a_{(1)} \, m \, S(a_{(2)}) \; ,
\ee
\be
\lb{2.44bb}
\alpha \triangleright m = \langle \langle S(\alpha) , m_{(-1)} \rangle \rangle \, m_{(0)} \; ,
\ee
where  $\Delta_L(m) = m_{(-1)} \otimes m_{(0)}$ is the left  ${\cal A}$-coaction on $M^R$
and $a \in {\cal A}$, $\alpha \in {\cal A}^o$.
Note that  left ${\cal A}$-module action (\ref{2.44b}) respects the right coinvariance of $M^R$.
The compatibility condition for the left ${\cal A}$-action (\ref{2.44b})
and left ${\cal A}$-coaction  $\Delta_L$
is written in the form (we represent $\Delta_L(a \triangleright m)$ in two different ways):
\be
\lb{2.44c}
(a \triangleright m)_{(-1)} \otimes (a \triangleright m)_{(0)} =
a_{(1)} \, m_{(-1)} \, S(a_{(3)})
\otimes a_{(2)} \triangleright m_{(0)} \; .
\ee
A module with the property (\ref{2.44c}) is called Yetter-Drinfeld module. Then,
using (\ref{2.44b}), (\ref{2.44c}) and opposite coproduct for ${\cal A}^o$, we obtain
\be
\lb{2.44d}
\begin{array}{c}
\alpha \triangleright (a \triangleright m) = \alpha \triangleright (a_{(1)} \, m \, S(a_{(2)})) =
 \langle \langle S(\alpha) , a_{(1)} \, m_{(-1)} \, S(a_{(3)})  \rangle \rangle \,
a_{(2)} \triangleright  m_{(0)} = \\ \\
=  \langle \langle S(\alpha_{(1)}) , \, a_{(1)} \rangle \rangle
\,\langle \langle \alpha_{(3)}, \,   a_{(3)}  \rangle \rangle \,
a_{(2)} \triangleright  (\alpha_{(2)} \triangleright m ) \; ,
\end{array}
\ee
and one can recognize in eq. (\ref{2.44d}) the quantum double multiplication
formula (\ref{2.44a}).

It follows from Eqs. (\ref{2.3}), (\ref{2.7}) and from the identities for the skew antipode
(\ref{2.25}) that
\be
\lb{2.45}
\left( m^{q}_{tk} \Delta^{ks}_{m} \right)
\left( m^{t}_{ip}(S^{-1})^{p}_{n} \Delta^{nr}_{s} \right) =
\delta^{q}_{i}\delta^{r}_{m} \; ,
\ee
and this enables us to rewrite (\ref{2.44}) in the form
$$
\left( m^{q}_{tk} \Delta^{ks}_{m} \right)
(I \Join e^{t})(e_{s} \Join I) =  ( m^{q}_{kl} \Delta_{m}^{jk} )
(e_{j} \Join I)(I \Join e^{l}) \; .
$$
This equation is equivalent to the axiom (\ref{2.26}) for the universal matrix ${\cal R}$ (\ref{2.39}).
The relations (\ref{2.27}) for ${\cal R}$ (\ref{2.39}) are readily verified.
Thus, ${\cal D}({\cal A})$ is indeed a quasitriangular Hopf algebra with universal
${\cal R}$ matrix represented by (\ref{2.39}).

In conclusion, we note that many relations for the structure constants
of Hopf algebras [for example, the relation (2.45)] can be obtained
and represented in a transparent form by means of the following diagrammatic technique: \\
\unitlength=1cm
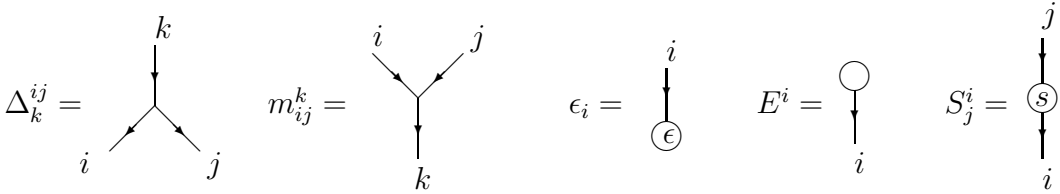
\begin{picture}(15,4)
\put(0,2){$\Delta_{k}^{ij}=$}
\put(2,2.1){\line(-1,-1){0.6}}
\put(1.7,1.8){\vector(-1,-1){0.1}}
\put(1,1.2){$i$}
\put(2,2.1){\line(0,1){0.8}}
\put(2,2.5){\vector(0,-1){0.1}}
\put(2,3){$k$}
\put(2,2.1){\line(1,-1){0.6}}
\put(2.3,1.8){\vector(1,-1){0.1}}
\put(2.7,1.2){$j$}
\put(3.5,2){$m^{k}_{ij}=$}
\put(5.5,2.2){\line(1,1){0.6}}
\put(5.8,2.5){\vector(-1,-1){0.1}}
\put(6.2,2.9){$j$}
\put(5.5,2.2){\line(0,-1){0.8}}
\put(5.5,1.8){\vector(0,-1){0.1}}
\put(5.45,1){$k$}
\put(5.5,2.2){\line(-1,1){0.6}}
\put(5.2,2.5){\vector(1,-1){0.1}}
\put(4.9,2.9){$i$}
\put(7.5,2){$\epsilon_{i} =$}
\put(8.8,1.9){\line(0,1){0.7}}
\put(8.8,2.3){\vector(0,-1){0.1}}
\put(8.8,1.7){\circle{0.4}}
\put(8.75,1.6){$\epsilon$}
\put(8.8,2.7){$i$}
\put(10,2){$E^{i} =$}
\put(11.3,2.3){\line(0,-1){0.7}}
\put(11.3,2){\vector(0,-1){0.1}}
\put(11.3,2.5){\circle{0.4}}
\put(11.3,1.2){$i$}
\put(12.5,2){$S^{i}_{j} =$}
\put(13.8,2.4){\line(0,1){0.6}}
\put(13.8,2.75){\vector(0,-1){0.1}}
\put(13.8,3.2){$j$}
\put(13.8,2.2){\circle{0.4}}
\put(13.8,2){\line(0,-1){0.6}}
\put(13.8,1.75){\vector(0,-1){0.1}}
\put(13.7,2.1){$s$}
\put(13.8,1){$i$}
\end{picture}

For example, the axioms of associativity (2.3) and coassociativity
(\ref{2.7}) and the axioms for the antipode (\ref{2.24}) can be
represented in the form

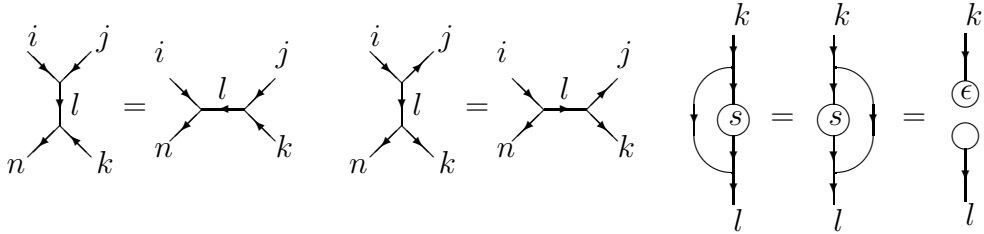
\unitlength=0.7cm
\begin{picture}(20,5)
\put(1,2.1){\line(-1,-1){0.6}}
\put(0.7,1.8){\vector(-1,-1){0.1}}
\put(0,1.2){$n$}
\put(1,2.1){\line(0,1){0.8}}
\put(1,2.5){\vector(0,-1){0.1}}
\put(1.2,2.3){$l$}
\put(1,2.1){\line(1,-1){0.6}}
\put(1.3,1.8){\vector(-1,1){0.1}}
\put(1.7,1.2){$k$}
\put(1,2.9){\line(1,1){0.6}}
\put(1.3,3.2){\vector(-1,-1){0.1}}
\put(1.7,3.6){$j$}
\put(1,2.9){\line(-1,1){0.6}}
\put(0.7,3.2){\vector(1,-1){0.1}}
\put(0.4,3.6){$i$}
\put(2.2,2.3){$=$}
\put(3.7,2.4){\line(-1,-1){0.6}}
\put(3.4,2.1){\vector(-1,-1){0.1}}
\put(2.8,1.5){$n$}
\put(3.7,2.4){\line(1,0){0.8}}
\put(4.1,2.4){\vector(-1,0){0.1}}
\put(4,2.6){$l$}
\put(4.5,2.4){\line(1,-1){0.6}}
\put(4.8,2.1){\vector(-1,1){0.1}}
\put(5.1,1.5){$k$}
\put(4.5,2.4){\line(1,1){0.6}}
\put(4.8,2.7){\vector(-1,-1){0.1}}
\put(5.1,3.3){$j$}
\put(3.7,2.4){\line(-1,1){0.6}}
\put(3.4,2.7){\vector(1,-1){0.1}}
\put(2.8,3.3){$i$}

\put(7.5,2.1){\line(-1,-1){0.6}}
\put(7.2,1.8){\vector(-1,-1){0.1}}
\put(6.5,1.2){$n$}
\put(7.5,2.1){\line(0,1){0.8}}
\put(7.5,2.5){\vector(0,-1){0.1}}
\put(7.7,2.3){$l$}
\put(7.5,2.1){\line(1,-1){0.6}}
\put(7.8,1.8){\vector(1,-1){0.1}}
\put(8.2,1.2){$k$}
\put(7.5,2.9){\line(1,1){0.6}}
\put(7.8,3.2){\vector(1,1){0.1}}
\put(8.2,3.6){$j$}
\put(7.5,2.9){\line(-1,1){0.6}}
\put(7.2,3.2){\vector(1,-1){0.1}}
\put(6.9,3.6){$i$}
\put(8.7,2.3){$=$}
\put(10.2,2.4){\line(-1,-1){0.6}}
\put(9.9,2.1){\vector(-1,-1){0.1}}
\put(9.3,1.5){$n$}
\put(10.2,2.4){\line(1,0){0.8}}
\put(10.6,2.4){\vector(1,0){0.1}}
\put(10.5,2.6){$l$}
\put(11,2.4){\line(1,-1){0.6}}
\put(11.3,2.1){\vector(1,-1){0.1}}
\put(11.6,1.5){$k$}
\put(11,2.4){\line(1,1){0.6}}
\put(11.3,2.7){\vector(1,1){0.1}}
\put(11.6,3.3){$j$}
\put(10.2,2.4){\line(-1,1){0.6}}
\put(9.9,2.7){\vector(1,-1){0.1}}
\put(9.3,3.3){$i$}

\put(13.8,2.5){\line(0,1){1.3}}
\put(13.8,2.8){\vector(0,-1){0.1}}
\put(13.8,3.5){\vector(0,-1){0.1}}
\put(13.8,4){$k$}
\put(13.8,2.2){\circle{0.6}}
\put(13.8,1.9){\line(0,-1){1.3}}
\put(13.7,2.1){$s$}
\put(13.8,0.1){$l$}
\put(13.8,2.2){\oval(1.5,2)[l]}
\put(13.05,2.2){\vector(0,-1){0.1}}
\put(13.8,1.5){\vector(0,-1){0.1}}
\put(13.8,0.9){\vector(0,-1){0.1}}

\put(14.5,2.1){$=$}
\put(15.7,2.5){\line(0,1){1.3}}
\put(15.7,2.8){\vector(0,-1){0.1}}
\put(15.7,3.5){\vector(0,-1){0.1}}
\put(15.7,4){$k$}
\put(15.7,2.2){\circle{0.6}}
\put(15.7,1.9){\line(0,-1){1.3}}
\put(15.6,2.1){$s$}
\put(15.7,0.1){$l$}
\put(15.7,2.2){\oval(1.5,2)[r]}
\put(16.45,2.2){\vector(0,-1){0.1}}
\put(15.7,1.5){\vector(0,-1){0.1}}
\put(15.7,0.9){\vector(0,-1){0.1}}

\put(17,2.1){$=$}
\put(18.2,2.95){\line(0,1){0.9}}
\put(18.2,3.4){\vector(0,-1){0.1}}
\put(18.2,2.7){\circle{0.5}}
\put(18.1,2.6){$\epsilon$}
\put(18.2,4){$k$}

\put(18.2,1.65){\line(0,-1){1}}
\put(18.2,1.15){\vector(0,-1){0.1}}
\put(18.2,1.9){\circle{0.5}}
\put(18.2,0.2){$l$}
\end{picture}

Now we make three important remarks relating to the further
development of the theory of Hopf algebras.

\subsection{\bf \em Twisted, ribbon and quasi-Hopf
 algebras\label{trqH}}
\setcounter{equation}0

\noindent
{\bf Remark 1.}{\it Twisted Hopf algebras.} \\
Consider a Hopf algebra
${\cal A}$ $(\Delta, \, \epsilon, \, S)$.
Let ${\cal F}$
be an invertible element of ${\cal A} \otimes {\cal A}$ such that:
\be
\lb{ef}
(\epsilon \otimes id) {\cal F} = 1 = (id \otimes \epsilon ) {\cal F} \;  ,
\ee
and we denote ${\cal F} = \sum_i \alpha_i \otimes \beta_i$,
${\cal F}^{-1} = \sum_i \, \gamma_i \otimes \delta_i$,
$I \equiv 1$.
Following the twisting procedure \cite{17}
one can define a new Hopf algebra
${\cal A}^{(F)}$
$(\Delta^{(F)}, \, \epsilon^{(F)}, \, S^{(F)})$ (twisted Hopf algebra)
with the new structure mappings
\be
\lb{2}
\Delta^{(F)}(a) = {\cal F} \, \Delta(a) \, {\cal F}^{-1} \; ,
\ee
\be
\lb{222}
\epsilon^{(F)}(a) = \epsilon (a) \; , \;\;\;  S^{(F)}(a) = U \, S(a) \, U^{-1}  \;\;\;
(\forall a \in {\cal A}) \; ,
\ee
where the twisting element ${\cal F}$ satisfies the cocycle equation
\be
\lb{cocycl}
{\cal F}_{12} \, (\Delta \otimes id) {\cal F} =
{\cal F}_{23} \, (id \otimes \Delta) {\cal F} \; ,
\ee
and the element $U =  \alpha_i \, S(\beta_i)$ is invertible and obeys
\be
\lb{U}
U^{-1} =  S(\gamma_i) \, \delta_i
 \; , \;\;\;
 S(\alpha_i) \, U^{-1} \, \beta_i = 1 \; ,
\ee
(the summation over $i$ is assumed). First of all we show
that the algebra ${\cal A}^{(F)}$
$(\Delta^{(F)}, \, \epsilon)$ is a bialgebra. Indeed,
the cocycle equation (\ref{cocycl})
guarantees the coassociativity condition (\ref{2.7})
for the new coproduct $\Delta^{(F)}$ (\ref{2}).
Then the axioms for counit $\epsilon$ (\ref{2.8}) are easily deduced from
(\ref{ef}). Considering the identity
$$
m(id \otimes S \otimes id) \left( {\cal F}^{-1}_{23} \, {\cal F}_{12} \,
(\Delta \otimes id) {\cal F} \right) = m(id \otimes S \otimes id)
(id \otimes \Delta ) {\cal F}
$$
we obtain the form for $U^{-1}$ (\ref{U}).
The second relation in (\ref{U}) is obtained from the identity:
$m(S \otimes id) {\cal F}^{-1}{\cal F} = 1$.

Now the new antipode $S^{(F)}$ (\ref{222}) follows from equation
$$
m (id \otimes S) \, (\Delta^{(F)}(a) \, {\cal F}) =
m (id \otimes S) \, ({\cal F} \, \Delta(a) ) \; ,
$$
which is rewritten in the form $\tilde{a}_{(1)} U S(\tilde{a}_{(2)}) = \epsilon(a) \, U$,
where $\Delta^{(F)}(a) = \tilde{a}_{(1)} \otimes \tilde{a}_{(2)}$.

If the algebra ${\cal A}$ is a quasitriangular Hopf algebra with the universal
${\cal R}$ matrix
(\ref{2.26}), then the new Hopf algebra ${\cal A}^{(F)}$ is also quasitriangular
and a new universal ${\cal R}$-matrix is
\be
\lb{2a}
{\cal R}^{(F)} = {\cal F}_{21} \, {\cal R} \, {\cal F}^{-1} \; ,
\ee
since we have
$$
{\Delta^{(F)}}' = {\cal F}_{21} \, \Delta^{\sf op} \, {\cal F}_{21}^{-1} =
{\cal F}_{21} \, {\cal R} \, \Delta \, {\cal R}^{-1} \, {\cal F}_{21}^{-1} =
\left( {\cal F}_{21} \, {\cal R} \, {\cal F}^{-1} \right) \Delta^{(F)} \left( {\cal F} \,
{\cal R}^{-1} \, {\cal F}_{21}^{-1} \right) \; .
$$
The Yang - Baxter eq. (\ref{2.30}) for $R$-matrix (\ref{2a}) can be
directly checked with the help of (\ref{2.26}) and (\ref{cocycl}).

Impose additional relations on ${\cal F}$
\be
\lb{4}
(\Delta \otimes id) {\cal F} =
{\cal F}_{13} \, {\cal F}_{23} \;\; , \;\;\;
(id \otimes \Delta) {\cal F} =
{\cal F}_{13} \, {\cal F}_{12} \; ,
\ee
which, together with (\ref{cocycl}), implies the Yang-Baxter equation for
${\cal F}$.
Using (\ref{2.26}) one deduces from (\ref{4}) equations
\be
\lb{5aa}
{\cal R}_{12} \, {\cal F}_{13} \, {\cal F}_{23} =
{\cal F}_{23} \, {\cal F}_{13} \, {\cal R}_{12} \;\; , \;\;\;
{\cal F}_{12} \, {\cal F}_{13} \, {\cal R}_{23} =
{\cal R}_{23} \, {\cal F}_{13} \, {\cal F}_{12} \; .
\ee
Eqs. (\ref{5aa}) and the Yang-Baxter relations for universal elements ${\cal R}$, ${\cal F}$
define the twist which is proposed in
\cite{16} (the additional condition
${\cal F}^{21}{\cal F} = 1 \otimes 1$ is assumed in \cite{16}).

Note that if ${\cal A}$ is the Hopf algebra
of functions on the group algebra of group $G$ (\ref{2.15b}),
then eq. (\ref{cocycl})
can be written in the form of 2-cocycle equation:
$$
{\cal F}(a,b) \, {\cal F}(ab,c) = {\cal F}(b,c){\cal F}(a,bc) \; ,
\;\;\; (\forall a,b,c \in G) \; ,
$$
for the projective representation $\rho$ of $G$:
$\rho(a) \rho(b) = {\cal F}(a,b) \, \rho(a \, b)$.
That is why  eq. (\ref{cocycl}) is called cocycle equation.

Many explicit solutions
of the cocycle equation (\ref{cocycl}) are known (see e.g.
\cite{tw1}--\cite{tw3} and references therein).

\noindent
{\bf Remark 2.}{\it Ribbon Hopf algebras.} \\
 Here we explain the notion of the ribbon Hopf algebras \cite{16'}.
Consider quasitriangular Hopf algebra ${\cal A}$ and represent the universal
${\cal R}$-matrix in the form
\be
\lb{unrab}
{\cal R}= \sum_{\mu} \alpha_{\mu} \otimes \beta_{\mu} \; ,
  \;\;\;\; {\cal R}^{-1} =
 \sum_{\mu} \, \gamma_{\mu} \otimes \delta_{\mu} \; ,
\ee
where $\alpha_{\mu}, \, \beta_{\mu}, \, \gamma_{\mu}, \, \delta_{\mu} \in {\cal A}$.
By using the right equalities in
(\ref{2.29a}) we represent the identities
$(id \otimes S)\bigl( {\cal R}{\cal R}^{-1} \bigr)= I
= (id \otimes S)\bigl({\cal R}^{-1}{\cal R}\bigr)$ as
\be
\lb{close}
 \alpha_{\mu} \, \alpha_\nu \otimes \beta_{\nu} \, S(\beta_\mu)
  =I = \alpha_{\mu} \, \alpha_\nu \otimes
 S(\beta_{\nu}) \, \beta_\mu \; .
\ee
(the summation over repeated
indices $\mu$ and $\nu$ is assumed
and we write $I$ instead of $(I \otimes I)$)
while for $(S \otimes id){\cal R}{\cal R}^{-1} = I =
(S \otimes id){\cal R}^{-1}{\cal R}$
we have
\be
\lb{close1}
 S(\gamma_{\mu}) \, \gamma_\nu \otimes \delta_{\nu} \, \delta_\mu
  = I =
 \gamma_{\mu} \, S(\gamma_\nu) \otimes \delta_{\nu} \, \delta_\mu \; .
\ee
We use identities (\ref{close}) and (\ref{close1})
below in Subsection {\bf \ref{qtrace}}
(Remark 1).

Consider the element
$u=\sum_{\mu} \, S(\beta_{\mu}) \, \alpha_{\mu}$ for which
the following proposition holds

\begin{proposition}\label{prop1}
{\it (see \cite{13'}). \\
1.) For any $a \in {\cal A}$ we have
\be
\lb{0.2}
S^{2}(a) \, u = u \, a \; .
\ee
2.) the element $u$ is invertible, with
 \be
 \lb{u-inv}
 u^{-1} =  S^{-1}(\delta_{\mu}) \, \gamma_{\mu} \; .
 \ee
}
\end{proposition}
{\bf Proof.}
1.) From the relation (\ref{2.26}) it follows that
$\forall a \in {\cal A}$ (the summation signs are omitted):
$$
\alpha_{\mu} \, a_{(1)} \otimes \beta_{\mu} \, a_{(2)} \otimes a_{(3)} =
a_{(2)} \, \alpha_{\mu} \otimes a_{(1)} \, \beta_{\mu} \otimes a_{(3)} \; ,
$$
where $a_{(1)} \otimes  a_{(2)} \otimes a_{(3)} =
(\Delta \otimes id)\Delta(a)$.
From this we obtain
$$
S^{2}( a_{(3)} ) \, S( \beta_{\mu} \, a_{(2)}) \, \alpha_{\mu} a_{(1)} =
S^{2}( a_{(3)} ) \, S( a_{(1)} \, \beta_{\mu} ) \, a_{(2)} \alpha_{\mu} \; ,
$$
or
$$
S^{2}( a_{(3)} ) \, S( a_{(2)}) \, u \, a_{(1)} =
S^{2}( a_{(3)} ) \, S(\beta_{\mu}) \, S( a_{(1)} ) \, a_{(2)} \alpha_{\mu} \; .
$$
Applying to this equation the axioms (\ref{2.21}),
we obtain (\ref{0.2}). \\
2.) Putting $w =  S^{-1}(\delta_{\mu}) \, \gamma_{\mu}$, we have
$$
u \, w =  u \, S^{-1}(\delta_{\mu}) \, \gamma_{\mu} =
 S(\delta_{\mu})\, u  \, \gamma_{\mu} =
 S(\beta_{\nu} \, \delta_{\mu})\, \alpha_{\nu}
 \, \gamma_{\mu} \; .
$$
Since ${\cal R} \cdot {\cal R}^{-1}
= \alpha_{\nu} \gamma_{\mu} \otimes
\beta_{\nu} \delta_{\mu} = I$, we have $u \, w = I$.
It follows from last equation and from (\ref{0.2}) that $S^{2}(w) \, u =1$,
and therefore the element $u$ has both a right and left inverse
(\ref{u-inv}). \hfill \qed

\noindent
Thus, the element $u$ is invertible and we can rewrite (\ref{0.2}) in the form
\be
\lb{0.3}
S^{2}(a) = u \, a \, u^{-1} \; .
\ee
This relation shows, in particular, that the operation of taking
the antipode is not involutive.

\begin{proposition}\label{prop1b}
{\it (see \cite{13'}). \\
Define the following elements:
\be
\lb{u1234}
u_{1} \equiv u =  S(\beta_{\mu}) \, \alpha_{\mu}
\; , \;\;\;
u_{2} =  S(\gamma_{\mu}) \, \delta_{\mu} \; , \;\;\;
u_{3} =  \beta_{\mu} \, S^{-1}(\alpha_{\mu})  \; , \;\;\;
u_{4} = \gamma_{\mu} \, S^{-1}(\delta_{\mu})  \; .
\ee
The relations (\ref{0.3})
are satisfied if we take any of the elements
$u_i$ from (\ref{u1234}):
\be
\lb{0.3b}
S^{2}(a) = u_i \, a \, u_i^{-1} \; , \;\;\;\;\;
\forall a \in {\cal A} \; .
\ee
In addition we have
$S(u_1)^{-1} = u_{2}$, $S(u_{3})^{-1} = u_{4}$,
and it turns out
that all $u_{i}$ commute with each other,
 while the elements
$u_{i}u_{j}^{-1} = u_{j}^{-1}u_{i}$ are central in ${\cal A}$.
Consequently, the element $u \, S(u)= u_{1} \, u^{-1}_{2}$ is also central.} \end{proposition}
{\bf Proof.} In view of relation (\ref{u-inv}) we have
$S(u_1)^{-1} = S(u^{-1}) = S(\gamma_\mu) \delta_\mu = u_2$
and $u_2^{-1} = S(u) = S^{-1}(u)=
S^{-1}(\alpha_{\mu}) \, \beta_{\mu}$, where we use
the identity $S^{2}(u)=u$ which follows from (\ref{0.3}).
Applying the map $S$ to both parts of (\ref{0.3}), we deduce
$S^3(a) = u_2 S(a) u_2^{-1}$ which is equivalent to
(\ref{0.3b}) for $i=2$.
Note that from (\ref{2.29a}) we have
${\cal R}^{\pm} = (S^{-1} \otimes S^{-1}){\cal R}^{\pm}$.
Thus,
one can make in all formulas above the
substitution $\alpha_\mu \to S^{-1}(\alpha_\mu )$,
$\beta_\mu \to S^{-1}(\beta_\mu )$ and
$\gamma_\mu \to S^{-1}(\gamma_\mu )$,
$\delta_\mu \to S^{-1}(\delta_\mu )$
to exchange the elements $u_1$ and $u_2$
respectively to the elements $u_3$ and $u_4$. It means that
equations (\ref{0.3b}) are valid for $i=3,4$ and we have
$u_4 = S(u_3)^{-1}$. Relations (\ref{0.3b})
 yields $S^2(u_j)=u_j$ $(\forall j)$ and substitution
 $a=u_j$ to (\ref{0.3b})
gives $u_i u_j = u_j u_i$ $(\forall i,j=1,...,4)$.
Finally, for any $a \in {\cal A}$, we have
$u_{j}^{-1}u_{i} \, a \, u_{i}^{-1}u_{j} =
u_{j}^{-1}\, S^2(a) \, u_{j} = a$, which means that
elements $u_{j}^{-1}u_{i} = u_{i} u_{j}^{-1}$
are central. \hfill \qed

\vspace{0.2cm}
\noindent
In Ref. \cite{13'} it was noted that
$$
\Delta(u) = ({\cal R}_{21}{\cal R}_{12})^{-1}(u \otimes u) =
(u \otimes u) ({\cal R}_{21}{\cal R}_{12})^{-1} \; .
$$
\vspace{0.2cm}

On the basis of all these propositions, we introduce the important concept of a ribbon Hopf algebra
(see \cite{16'}):


\newtheorem{def7}[def1]{Definition}
\begin{def7} \label{def7}
{\it Consider a quasitriangular Hopf algebra $({\cal A}, \; {\cal R})$. Then the triplet
$({\cal A}, \; {\cal R}, \; v)$ is called a ribbon Hopf algebra if
$v$ is a central element in ${\cal A}$ and
$$
v^{2} = u \, S(u) \; , \;\; S(v) = v \; , \;\; \epsilon(v) =1 \; ,
$$
$$
\Delta(v) = ({\cal R}_{21} \, {\cal R}_{12})^{-1} \, (v \otimes v) \; .
$$}
\end{def7}
For each quasitriangular Hopf algebra ${\cal A}$, we can define
${\cal A}$-colored ribbon graphs \cite{16'}. If, moreover, ${\cal A}$
is a ribbon Hopf algebra, then for each
${\cal A}$-colored ribbon graph we can associate the central element of
${\cal A}$ that generalizes the
Jones polynomial being an invariant of a knot in
$\mathbb{R}^3$
(see \cite{16'}, \cite{18}).

\noindent
{\bf Remark 3.}{\it Quasi-Hopf algebras.} \\
One can introduce a generalization of a Hopf algebra, called
a quasi-Hopf algebra, \cite{17} which is defined as an associative unital algebra ${\cal A}$
with homomorphism $\Delta: \;\;
{\cal A} \rightarrow {\cal A} \otimes {\cal A}$, homomorphism
$\epsilon: \;\; {\cal A} \rightarrow \mathbb{C}$,
antiautomorphism $S: \;\; {\cal A} \rightarrow {\cal A}$
and invertible element $\Phi \in
{\cal A}\otimes {\cal A}\otimes {\cal A}$. At the same time $\Delta$, $\epsilon$, $\Phi$
and $S$ satisfy the axioms
\be
\lb{2.47}
(id \otimes \Delta) \Delta(a) =
\Phi \cdot (\Delta \otimes id) \Delta(a) \cdot \Phi^{-1} , \;\;
a \in {\cal A},
\ee
\be
\lb{2.48}
(id \otimes id \otimes \Delta)(\Phi )
\cdot (\Delta \otimes id\otimes id) (\Phi) =
(I \otimes \Phi) \cdot (id \otimes \Delta \otimes id)(\Phi ) \cdot
(\Phi \otimes I),
\ee
\be
\lb{2.48a}
(\epsilon \otimes id) \Delta = id = (id \otimes \epsilon) \Delta \; ,\;\;
(id \otimes \epsilon \otimes id)\Phi = I \otimes I
\ee
\be
\lb{2.49}
S(a_{(1)}) \, \alpha \,  a_{(2)} = \epsilon(a) \, \alpha \; , \;\;
a_{(1)} \, \beta \, S(a_{(2)}) = \epsilon(a) \, \beta \; ,
\ee
$$
\phi_i \, \beta \,  S(\phi'_i) \, \alpha \, \phi''_i = I \; , \;\;
S(\bar{\phi}_i) \, \alpha \,  \bar{\phi}'_i \, \beta \, S(\bar{\phi}''_i) = I \; ,
$$
where $\alpha$   and  $\beta$   are   certain   fixed   elements   of
${\cal A}$, $\Delta(a) = a_{(1)} \otimes a_{(2)}$, and
$$
\Phi : = \phi_i \otimes \phi'_i \otimes \phi''_i \; , \;\;
\Phi^{-1} : = \bar{\phi}_i \otimes \bar{\phi}'_i
\otimes \bar{\phi}''_i \; ,
$$
(summation over $i$ is assumed).
Thus, a quasi-Hopf algebra differs from an ordinary Hopf algebra in
that the axiom of coassociativity is replaced by the weaker condition (\ref{2.47}).
In other words, a quasi-Hopf algebra is noncoassociative, but this
noncoassociativity is kept under control by means of the element $\Phi$.
The axioms (\ref{2.49}) (which looks like different definitions
of the left and right antipodes) generalize the axioms (\ref{2.21}) for usual Hopf algebras
and consequently the elements $\alpha$ and $\beta$ involve into the play with
the contragredient representations of the quasi-Hopf algebras.

To make the pentagonal condition (\ref{2.48}) more transparent, let us consider
(following Ref. \cite{17}) the algebra ${\cal A}$ as the algebra of functions on
a "noncommutative" space $X$ equipped with a $*$ product: $X \times X \rightarrow X$.
Then, elements $a \in {\cal A}$, $b \in {\cal A} \otimes {\cal A}, \dots$
are written in the form $a(x)$, $b(x,y)$ $\dots$ and $\Delta(a)$ is represented as
$a(x * y)$. The homomorphism $\epsilon$ defines the point in $X$, which we denote $1$ and
instead of $\epsilon(a)$ we write $a(1)$. Then, equations (\ref{2.47}) -- (\ref{2.48a})
are represented in the form \cite{17}:
$$
a(x *(y*z)) = \Phi(x,y,z) \, a((x *y )*z) \, \Phi(x,y,z)^{-1} \; ,
$$
\be
\label{5phi}
\Phi(x,y,z*u) \, \Phi(x*y,z,u) = \Phi(y,z,u)  \, \Phi(x,y*z,u) \, \Phi(x,y,z) \; ,
\ee
$$
a(1 * x) = a(x) = a(x *1) \; , \;\;\; \Phi(x,1,z) = 1 \; .
$$
Now it is clear that (\ref{5phi}) (and respectively (\ref{2.48}))
is the sufficient condition for the commutativity of the diagram:

\unitlength=0.7cm
\begin{picture}(20,5)
\put(1,3){$a\left(((x *y)*z)*u \right)$}
\put(6.5,3){$\longrightarrow$}
\put(8,3){$a((x *y)*(z*u))$}
\put(13.5,3){$\longrightarrow$}
\put(15,3){$a(x *(y*(z*u)))$}
\put(4,2.5){\vector(0,-1){1}}
\put(1,0.8){$a((x *(y*z))*u)$}
\put(9,0.955){\line(1,0){3}}
\put(11.5,0.8){$\rightarrow$}
\put(17,2.5){\vector(0,-1){1}}
\put(15,0.8){$a(x *((y*z)*u))$}

\end{picture}

\noindent
{\bf Remark.} Applications of the theory of quasi-Hopf algebras to the solutions of the
Knizhnik-Zamolodchikov equations are discussed in Ref. \cite{17}. On the other hand, one can
suppose that, by virtue of the occurrence of the pentagonal relation (\ref{2.48}) for the element
$\Phi$, quasi-Hopf algebras will be associated with multidimensional generalizations of Yang-Baxter
equations.

\section{THE YANG--BAXTER EQUATION AND \\ QUANTIZATION OF LIE GROUPS\label{ybeqg}}
\setcounter{equation}0

In this section, we discuss the $R$-matrix approach
to the theory of quantum groups
\cite{10}, on the basis of which we perform a
quantization of classical Lie groups and also some
Lie supergroups. We present trigonometric solutions of
the Yang--Baxter equation invariant under the adjoint action of the quantum
groups $GL_q(N)$, $SO_q(N)$, $Sp_q(2n)$ and supergroups
$GL_q(N|M)$, $Osp_q(N|2m)$. We briefly discuss the corresponding
Yangian (rational) solutions,
and $Z_{N} \otimes Z_{N}$ symmetric elliptic solutions
of the Yang--Baxter equation. We also show that for every (trigonometric)
solution $R(x)$ of the Yang--Baxter equation one can construct
 the set of difference equations which are called
 quantum Knizhnik--Zamolod\-chikov equations.

\subsection{\bf \em Numerical $R$-matrices\label{NuRm}}
\setcounter{equation}0

This Subsection is based on the results
 presented in \cite{18}, \cite{18a}.

 \subsubsection{Invertible and skew-invertible $R$-matrices}

Let ${\cal A}$ be a quasitriangular Hopf algebra.
Consider representations $T^{(\nu)}$
of ${\cal A}$
 in $N_\nu$-dimensional vector spaces $V_\nu$ (the index $\nu$
 enumerates representations). In view of (\ref{2.15a})
and (\ref{2.30}), the matrix
$(R_{(\nu,\mu)})^{i_\nu j_\mu}_{\;\; k_\nu l_\mu} =
(T^{(\nu)i_\nu}_{\;\;\;\; k_\nu} \otimes
T^{(\mu)j_\mu}_{\;\;\;\; l_\mu}) {\cal R}$,
where ${\cal R} \in {\cal A} \otimes {\cal A}$ is
the universal element (\ref{2.28}),
satisfies the generalized matrix Yang-Baxter equation
\be
\lb{ybenu}
(R_{(\nu,\mu)})^{i_{\nu}i_{\mu}}_{\;\; j_{\nu}j_{\mu}} \,
(R_{(\nu,\lambda)})^{j_{\nu}i_{\lambda}}_{\;\; k_{\nu}j_{\lambda}} \,
(R_{(\mu,\lambda)})^{j_{\mu}j_{\lambda}}_{\;\; k_{\mu}k_{\lambda}} =
(R_{(\mu,\lambda)})^{i_{\mu}i_{\lambda}}_{\;\; j_{\mu}j_{\lambda}} \,
(R_{(\nu,\lambda)})^{i_{\nu}j_{\lambda}}_{\;\; j_{\nu}k_{\lambda}} \,
(R_{(\nu,\mu)})^{j_{\nu}j_{\mu}}_{\;\; k_{\nu}k_{\mu}} \; .
\ee
Here the summation over repeated indices
$j_{\nu},j_{\mu},j_{\lambda}$ is assumed.
Let the representations $T^{(\nu)}$, $T^{(\mu)}$, $T^{(\lambda)}$
 be equivalent to a representation $T$ which acts
 in $N$-dimensional vector space $V$.
 In this case, according to (\ref{ybenu}),
  the image $R^{i j}_{\;\; k \ell} =
  (T^i_{\;\; k} \otimes T^j_{\;\; \ell}) {\cal R}$
 of the universal
element ${\cal R} \in {\cal A} \otimes {\cal A}$ satisfies the
standard matrix Yang-Baxter equation
\be
\lb{ybe}
R^{i_{1}i_{2}}_{\;\; j_{1}j_{2}} \,
R^{j_{1}i_{3}}_{\;\; k_{1}j_{3}} \,
R^{j_{2}j_{3}}_{\;\; k_{2}k_{3}} =
R^{i_{2}i_{3}}_{\;\; j_{2}j_{3}} \,
R^{i_{1}j_{3}}_{\;\; j_{1}k_{3}} \,
R^{j_{1}j_{2}}_{\;\; k_{1}k_{2}} \; .
\ee

 A lot of numerical
 solutions of the Yang--Baxter equations
 (\ref{ybenu}), (\ref{ybe}) can be constructed as
 representations of the universal
 ${\cal R}$-matrices. However,
not all numerical solutions $R$ of eqs. (\ref{ybenu})
 and (\ref{ybe})  are images
  $(T^{(\nu)} \otimes T^{(\mu)}) {\cal R}$ and
$(T \otimes T) {\cal R}$ of the universal
element ${\cal R}$ for some
quasitriangular Hopf algebra ${\cal A}$.
Below we consider solutions  $R \in {\rm End}(V \otimes V)$
 of the standard matrix Yang-Baxter equation (\ref{ybe})
that are not necessarily the universal
  ${\cal R}$-matrix representations.

 First we assume that a solution
  $R$ of eqs. (\ref{ybe}) is the invertible matrix
\be
\lb{inv}
R^{i_1 i_2}_{\;\; k_1 \ell_2} \, (R^{-1}) ^{k_1 \ell_2}_{\;\; j_1 j_2} =
\delta^{i_1}_{j_1} \, \delta^{i_2}_{j_2} =
(R^{-1})^{i_1 i_2}_{\;\; k_1 \ell_2} \,  R^{k_1 \ell_2}_{\;\; j_1 j_2}  \; .
\ee
Note that for all images $R = (T \otimes T) {\cal R}$,
 such invertibility follows
from the invertibility (\ref{2.29a})
of the universal element ${\cal R}$.
In terms of the concise matrix notation \cite{10},
we write relations (\ref{inv}) and (\ref{ybe})
in the following equivalent forms
 \be
 \lb{inv1}
 R_{12} R_{12}^{-1} = I_{12} = R_{12}^{-1} R_{12}
 \;\; \Longleftrightarrow \;\;
 \R_{12} \R_{12}^{-1} = I_{12} = \R_{12}^{-1} \R_{12} \; ,
 \ee
\be
R_{12} \, R_{13} \, R_{23}  =  R_{23} \, R_{13} \, R_{12}
\;\; \Longrightarrow
\lb{3.1.2}
\ee
\be
\lb{3.1.2i}
\R_{12} \, \R_{23} \, \R_{12}  =  \R_{23} \, \R_{12} \, \R_{23}
\;\; \Longrightarrow
\ee
\be
\lb{3.1.2is}
\R_{23} \, \R^{-1}_{12} \, \R^{-1}_{23}  =
\R^{-1}_{12} \, \R^{-1}_{23} \, \R_{12} \; , \;\;\;
\R_{12} \, \R^{-1}_{23} \, \R^{-1}_{12}  =
\R^{-1}_{23} \, \R^{-1}_{12} \, \R_{23} \; .
\ee
Here
$\R := P \, R$, the matrix $P$ is the permutation:
\be
\lb{perM}
P^{i_{1}i_{2}}_{j_{1}j_{2}} = \delta^{i_{1}}_{j_{2}} \delta^{i_{2}}_{j_{1}} \; , \;\;\;
\R^{i_{1}i_{2}}_{j_{1}j_{2}} = (P \, R)^{i_{1}i_{2}}_{j_{1}j_{2}} =
R^{i_{2}i_{1}}_{j_{1}j_{2}} \; ,
\ee
$I_{12} := I \otimes I$
($I \in {\rm Mat}(N)$ is the unit matrix in $V$)
and indices $1,2,3$ label the vector spaces $V$ in $V^{\otimes 3}$
where the corresponding
matrices $R_{12},R_{23},...$ act non-trivially, e.g.
 $R_{12} = R \otimes I$, $R_{23} = I \otimes R$, etc.
 We also note, that if matrix $R_{12}$ satisfies
 the Yang-Baxter equation (\ref{ybe}), then the
 matrix $R_{12}'=R_{21}$ also satisfies
 the Yang-Baxter equation
 \be
 \lb{ybe01}
 R_{21} \, R_{31} \, R_{32}  =  R_{32} \, R_{31} \, R_{21}
\;\;\;\; \stackrel{1 \leftrightarrow 3}{\Longleftrightarrow}
\;\;\;\;
R_{23} \, R_{13} \, R_{12}  = R_{12} \, R_{13} \, R_{23}  \; .
 \ee

In what follows, we introduce matrices
\be
\lb{3.1.3i}
\R_a := \R_{a, a+1} = I^{\otimes (a-1)} \otimes
\R \otimes I^{\otimes (M-a)} \; , \;\;\;
(a=1,\dots,M) \; ,
\ee
 which act in the space $V^{\otimes (M+1)}$ and, according to
 the Yang-Baxter equations
 (\ref{3.1.2i}), we have braid relations
\be
\R_{a} \, \R_{a+1} \, \R_{a}  =  \R_{a+1} \, \R_{a} \, \R_{a+1}
\;\;\;\;\; (a=1,...,M) \; .
\lb{3.1.3}
\ee
 In view of these relations and locality relations
\be
\lb{3.1.3is}
[\R_a , \, \R_b] =0 \;\; ({\rm for} \;\; |a-b| >1) \; ,
\ee
the invertible
 matrices $\R_a$ define a representation of generators of the
 braid group ${\cal B}_{M+1}$ (see definition
 in Subsec. {\bf \ref{gabg1}}).
 The name "braid group" is
justified since relation (\ref{3.1.3}) admits
the graphic visualization

\vspace{0.5cm}

\unitlength=4.5mm
\begin{picture}(17,5.5)(0,-0.5)

\put(0,1.9){$\R_{a+1} \cdot \R_a \cdot \R_{a+1} \;    =$}

\put(10.15,5.2){\line(1,-1){1}}
\put(11.0,4){\vector(-1,-1){0.8}}
\put(11.15,4.2){\vector(1,-1){1}}
\put(11.4,4.4){\line(1,1){0.8}}

\put(8.15,5.2){\vector(0,-1){2}}

\put(8,5){$\bullet$}
\put(10,5){$\bullet$}
\put(12,5){$\bullet$}

\put(8,5.5){\tiny $a$}
\put(9.7,5.5){\tiny $a+1$}
\put(11.7,5.5){\tiny $a+2$}


\put(8,3){$\bullet$}
\put(10,3){$\bullet$}
\put(12,3){$\bullet$}

\put(8.15,3.2){\line(1,-1){1}}
\put(9.0,2){\vector(-1,-1){0.8}}
\put(9.15,2.2){\vector(1,-1){1}}
\put(9.4,2.4){\line(1,1){0.8}}
\put(12.15,3.2){\vector(0,-1){2}}

\put(8,1){$\bullet$}
\put(10,1){$\bullet$}
\put(12,1){$\bullet$}


\put(10.15,1.2){\line(1,-1){1}}
\put(11.0,0){\vector(-1,-1){0.8}}
\put(11.15,0.2){\vector(1,-1){1}}
\put(11.4,0.4){\line(1,1){0.8}}

\put(8.15,1.2){\vector(0,-1){2}}

\put(8,-1){$\bullet$}
\put(10,-1){$\bullet$}
\put(12,-1){$\bullet$}


\put(14,1.9){$=$}

\put(16,5.5){\tiny $a$}
\put(17.7,5.5){\tiny $a+1$}
\put(19.7,5.5){\tiny $a+2$}

\put(16,5){$\bullet$}
\put(18,5){$\bullet$}
\put(20,5){$\bullet$}

\put(16.15,5.2){\line(1,-1){1}}
\put(17.0,4){\vector(-1,-1){0.8}}
\put(17.15,4.2){\vector(1,-1){1}}
\put(17.4,4.4){\line(1,1){0.8}}
\put(20.15,5.2){\vector(0,-1){2}}

\put(16,3){$\bullet$}
\put(18,3){$\bullet$}
\put(20,3){$\bullet$}


\put(18.15,3.2){\line(1,-1){1}}
\put(19.0,2){\vector(-1,-1){0.8}}
\put(19.15,2.2){\vector(1,-1){1}}
\put(19.4,2.4){\line(1,1){0.8}}

\put(16.15,3.2){\vector(0,-1){2}}

\put(16,1){$\bullet$}
\put(18,1){$\bullet$}
\put(20,1){$\bullet$}


\put(16,-1){$\bullet$}
\put(18,-1){$\bullet$}
\put(20,-1){$\bullet$}

\put(16.15,1.2){\line(1,-1){1}}
\put(17.0,0){\vector(-1,-1){0.8}}
\put(17.15,0.2){\vector(1,-1){1}}
\put(17.4,0.4){\line(1,1){0.8}}
\put(20.15,1.2){\vector(0,-1){2}}

\put(21.5,1.9){$  = \;\; \R_a \cdot \R_{a+1} \cdot \R_a$ \; ,}

\end{picture}
\vspace{-1.6cm}
\be
\lb{figR3}
{}
\ee

\vspace{0.8cm}

\noindent
that is an identity of two braids
with three strands (the third Reidemeister move).
Here we make use of the graphic representation

\unitlength=5mm
\begin{picture}(17,4.2)

\put(4.5,1.9){$\R_a \;\;  =$}

\put(7,3){$\bullet$}
\put(9,3){$\bullet$}

\put(9.7,3){$\dots$}
\put(11,3){$\bullet$}
\put(13,3){$\bullet$}
\put(14,3){$\dots$}
\put(15.5,3){$\bullet$}

\put(7,3.5){\tiny $1$}
\put(9,3.5){\tiny $2$}
\put(11,3.5){\tiny $a$}
\put(12.7,3.5){\tiny $a+1$}
\put(15.5,3.5){\tiny $M+1$}

\put(7,0.7){\tiny $1$}
\put(9,0.7){\tiny $2$}
\put(11,0.7){\tiny $a$}
\put(12.7,0.7){\tiny $a+1$}
\put(15.5,0.7){\tiny $M+1$}

\put(11.15,3.2){\line(1,-1){1}}
\put(12.0,2){\vector(-1,-1){0.75}}
\put(12.15,2.2){\vector(1,-1){0.9}}
\put(12.4,2.4){\line(1,1){0.8}}
\put(7.15,3.2){\vector(0,-1){1.9}}
\put(9.15,3.2){\vector(0,-1){1.9}}
\put(15.65,3.2){\vector(0,-1){1.9}}

\put(7,1){$\bullet$}
\put(9,1){$\bullet$}

\put(9.7,1){$\dots$}
\put(11,1){$\bullet$}
\put(13,1){$\bullet$}
\put(14,1){$\dots$}
\put(15.5,1){$\bullet$}
\put(17,1.9){.}
\end{picture}
\vspace{-1.6cm}
\be
\lb{figR}
{}
\ee

\vspace{0.3cm}
\noindent
We discuss in detail the braid group ${\cal B}_{M+1}$,
 its group algebra
$\mathbb{C}[{\cal B}_{M+1}]$ and finite dimensional
quotients of $\mathbb{C}[{\cal B}_{M+1}]$ in Section  {\bf \ref{gabg}} below.

Let $X(\R_a)$ be
a formal series in $\R_a^{\pm 1}$.
 The direct consequences of (\ref{3.1.3}) are equations
\be
X(\R_a) \, \R_{a+1} \, \R_a  =  \R_{a+1} \, \R_a \, X(\R_{a+1}) \; , \;\;
\R_a \, \R_{a+1} \, X(\R_a)  =  X(\R_{a+1}) \, \R_a \, \R_{a+1} \; ,
\lb{3.1.4}
\ee
which make it possible to carry functions $X(\R_{a})$,
 $X(\R_{a+1})$
through the operators $\R_{a+1} \R_a$ and $\R_a \R_{a+1}$.


\newtheorem{def8}[def1]{Definition}
\begin{def8} \label{def8}
{\it The matrix $R \in {\rm End}(V^{\otimes 2})$
is called skew-invertible if there
exists matrix $\Psi \in {\rm End}(V^{\otimes 2})$
such that (cf. (\ref{inv}))
\be
\lb{skewb}
R^{mk}_{\;\; jn} \, \Psi^{in}_{\;\; ml} =\delta^i_j \delta^k_l =
\Psi^{ni}_{\;\; lm} \, R^{km}_{\;\; nj} \; .
\ee
The index-free forms of these relations
are\footnote{The form (\ref{skew}) is very convenient for
calculations (see below) and was proposed in \cite{18c}.
Equations (\ref{skewa}) are equivalently written as
$\Psi^{t_2}_{12} \, R^{t_2}_{12} = I_{12}=
R^{t_2}_{12} \, \Psi^{t_2}_{12}$.}
\be
\lb{skewa}
R^{t_1}_{12} \, \Psi^{t_1}_{12}  = I_{12} \, , \;\;\;\;
\Psi^{t_1}_{12} \, R^{t_1}_{12} = I_{12} \; ,
\ee
\be
\lb{skew}
Tr_2 (\R_{12} \, \hat{\Psi}_{23} ) = P_{13} = Tr_2 \left( \hat{\Psi}_{12} \, \R_{23} \right) \; ,
\ee
where $\hat{\Psi} = P \, \Psi$. We say that the
invertible and skew-invertible $R$-matrix is
completely invertible if the
inverse matrix $R^{-1}$ is also skew-invertible, i.e. there exists
a matrix $\Phi \in {\rm End}(V^{\otimes 2})$ such that
$$
\Phi_{12}^{t_2} \, (R^{-1})_{12}^{t_2} = I_{12} =
(R^{-1})_{12}^{t_2} \, \Phi_{12}^{t_2}  \;\;  \Leftrightarrow \;\;
\Phi^{i_1 i_2}_{\; k_2j_1} \, (R^{-1})^{k_2 i_3}_{\; j_3 i_2} =
\delta^{i_1}_{j_3} \, \delta^{i_3}_{j_1} =
(R^{-1})^{i_1 i_2}_{\; k_2 j_1} \, \Phi^{k_2 i_3}_{\; j_3 i_2}
\;\; \Rightarrow
$$
\be
\lb{skew2}
Tr_2 \left( \hat{\Phi}_{12} \, \R^{-1}_{23} \right)
 = P_{13} = Tr_2 (\R^{-1}_{12} \, \hat{\Phi}_{23} ) \; ,
\ee
where $\R^{-1} = R^{-1} \, P$ and $\hat{\Phi} = \Phi \, P$.}
\end{def8}
 The skew-invertible
 $R$-matrices were considered in \cite{18},
where operator $\Psi_{12}$ was denoted as $((R_{12}^{t_1})^{-1})^{t_1}$
(cf. (\ref{skewa})).

\vspace{0.1cm}

\subsubsection{Quantum traces\label{qtrace}}

Now we define four matrices
\be
\lb{dmatr}
D_1 = Tr_2(\hat{\Psi}_{12}) \; , \;\;\; Q_2 = Tr_1 (\hat{\Psi}_{12}) \; ,
\ee
\be
\lb{dmatr1}
\overline{D}_1 = Tr_2(\hat{\Phi}_{12}) \; , \;\;\; \overline{Q}_2 = Tr_1 (\hat{\Phi}_{12}) \; ,
\ee
which are important for our consideration below.
\begin{proposition}\label{prop2}
{\it Let the Yang-Baxter matrix $R$ be invertible
and skew-invertible, then the following identities hold
\be
\lb{qtrs}
 Tr_2 ( \R_{12} D_2 ) = I_1 \; , \;\;\;
Tr_1 (Q_{1} \, \R_{12} )  = I_2 \; ,
\ee
\be
\lb{qtrs1}
 Tr_2 (\R^{-1}_{12} \, \overline{D}_2 ) = I_1 \; , \;\;\;
Tr_1 ( \overline{Q}_{1} \, \R^{-1}_{12} ) = I_2 \; ,
\ee
\be
\lb{sk2}
D_0 \, P_{02} = Tr_{3} D_3
\R_{23}^{- 1} \, \R_{03} \; , \;\;\;
D_0 \, P_{02} = Tr_{3} D_3
\R_{23} \, \R_{03}^{- 1} \; ,
\ee
\be
\lb{sk3}
Q_0 \, P_{02} = Tr_{1} Q_1
\R_{12}^{- 1} \, \R_{10} \; , \;\;\;
Q_0 \, P_{02} = Tr_{1} Q_1
\R_{12} \, \R_{10}^{- 1} \; ,
\ee
\be
\lb{sk4}
D_2 \, \R_{12}^{-1}  = \hat{\Psi}_{21} \, D_1 \;\; , \;\;\;
\R_{12}^{-1} \, D_2 = D_1 \, \hat{\Psi}_{21} \; ,
\ee
\be
\lb{sk6}
Q_1 \, \R_{12}^{-1}  = \hat{\Psi}_{21} \, Q_2 \;\; , \;\;\;
\R_{12}^{-1} \, Q_1 = Q_2 \, \hat{\Psi}_{21}  \; ,
\ee
where the matrices $D$ and $Q$ commute and satisfy
\be
\lb{sk9}
D_2 \, Q_2 = Q_2 \, D_2 = Tr_3(D_3 \, \R_2^{-1}) = Tr_1(Q_1 \, \R_1^{-1}) \; .
\ee
If the matrix $R$ is completely invertible, then
\be
\lb{sk5}
\R_{12} \, D_2 = D_1 \, \hat{\Phi}_{21} \;\; , \;\;\;
D_2 \, \R_{12}  = \hat{\Phi}_{21} \, D_1 \; .
\ee
\be
\lb{sk7}
\R_{12} \, Q_1 = Q_2 \, \hat{\Phi}_{21} \;\; , \;\;\;
Q_1 \, \R_{12}  = \hat{\Phi}_{21} \, Q_2 \; .
\ee
and the matrices $D$ and $Q$ are invertible such that
\be
\lb{sk8}
D^{-1} = \overline{Q} \; , \;\;\; Q^{-1} = \overline{D} \; .
\ee
Conversely, if the matrix $D$ (or $Q$) is invertible, then the
matrix $R$ is completely invertible.
For invertible matrices $D$ and $Q$ one has the relations
\be
\lb{sskk}
Tr_1(D_1^{-1} \, \R_{12}^{-1}) = I_2 =
 Tr_3(Q_3^{-1} \, \R_{23}^{-1}) \; .
\ee
}
\end{proposition}
{\bf Proof.} Identities (\ref{qtrs}) and (\ref{qtrs1})
follow from (\ref{skew}) and (\ref{skew2}). To obtain
(\ref{sk2}) and (\ref{sk3}) we multiply both sides of
equations (\ref{3.1.4}) (for $a=1$) from the left
by $\hat{\Psi}_{01}$ and from the right by $\hat{\Psi}_{34}$ and take the trace
$Tr_{13}$ ($\equiv Tr_1 \, Tr_3$). Using (\ref{skew}) we obtain
\be
\lb{ybe1}
Tr_1 \hat{\Psi}_{01} X(\R_1) P_{24} \R_1 =
Tr_3 \R_2 \, P_{02} \, X(\R_2) \, \hat{\Psi}_{34} \; ,
\ee
\be
\lb{ybe2}
Tr_1 \hat{\Psi}_{01} \R_1 \, P_{24} \, X(\R_1) =
Tr_3 X(\R_2) \, P_{02} \, \R_2 \, \hat{\Psi}_{34} \; ,
\ee
where $\R_a \equiv \R_{a\, a+1}$ (see (\ref{3.1.3i})).
We put $X(\R) = \R^{-1}$ in (\ref{ybe1}), (\ref{ybe2})
and take the traces $Tr_4$ or $Tr_0$.
Using (\ref{skew}) and (\ref{dmatr}) we obtain four identities
\be
\lb{4a}
D_0 \, I_2 = Tr_{3} D_3
\R_2^{\mp 1} \, P_{02} \, \R_2^{\pm 1} \; , \;\;\;
Q_0 \, I_2 = Tr_{1} Q_1
\R_1^{\pm 1} \, P_{02} \, \R_1^{\mp 1} \; ,
\ee
which can be easily written as
(\ref{sk2}) and (\ref{sk3}).
Applying to the both sides of the first relation in (\ref{sk2})
the operation $Tr_{0} (\hat{\Psi}_{10} ...)$
and to the both sides of the second relation in (\ref{sk2})
 the operation $Tr_{2} (\hat{\Psi}_{12} ...)$
we obtain identities
 $$
Tr_0 (\hat{\Psi}_{10} D_0 P_{02}) =
Tr_{03} (\hat{\Psi}_{10} D_3 \R_{23}^{-1} \R_{03}) \, , \;\;\;
Tr_{2} (D_0 \hat{\Psi}_{12} P_{02}) =
Tr_{23} (\hat{\Psi}_{12} D_3 \R_{23} \R_{03}^{-1}) \; ,
 $$
which, by means of (\ref{skew}), give (\ref{sk4}).
Similarly, applying to the both sides of
the first relation in (\ref{sk3}) the
operation $Tr_{0} (\dots \hat{\Psi}_{03})$
and to the both sides of the second relation in (\ref{sk3}) the
operation $Tr_{2} (\dots \hat{\Psi}_{23})$
we obtain (\ref{sk6}). Taking the traces $Tr_2(\dots)$ and $Tr_1(\dots)$
of (\ref{sk4}) and (\ref{sk6}), respectively, we deduce (\ref{sk9}).

If the matrix $R$ is completely invertible, then
acting to the first relation (\ref{sk2}) by
$Tr_{2} (\hat{\Phi}_{12} ...)$
and to the second relation (\ref{sk2}) by $Tr_{0} (\hat{\Phi}_{10} ...)$
we obtain (\ref{sk5}).
Analogously, acting to the first relation (\ref{sk3}) by $Tr_{2} (\dots \hat{\Phi}_{23})$
and to the second relation (\ref{sk3}) by
$Tr_{0} (\dots \hat{\Phi}_{03})$
we find (\ref{sk7}). Eqs. (\ref{sk8}) are obtained
by taking traces $Tr_{2}(\dots)$
and $Tr_{1}(\dots)$ of (\ref{sk5}) and
(\ref{sk7}), respectively, and applying (\ref{dmatr1}),
(\ref{qtrs}). Thus, for the completely invertible $R$ the matrices
$D$ and $Q$ are invertible.

Conversely, if the matrix $D$ is invertible, then $D_1 \hat{R}_{21} D_2^{-1}$
(cf. (\ref{sk5})) is the skew-inverse matrix for $\hat{R}^{-1}$. Indeed,
$$
Tr_2 \left( \hat{R}^{-1}_{12} \, D_2 \, \hat{R}_{32} \, D_3^{-1} \, \right) =
 Tr_2 \left( \hat{R}^{-1}_{12} \, D_2 \, \hat{R}_{32} \right) \, D_3^{-1}  =
$$
$$
=D_1 \, Tr_2 \left( \hat{\Psi}_{21} \, \hat{R}_{32} \right) \, D_3^{-1}  =
D_1 \, P_{13} \, D_3^{-1}  = P_{13} \; ,
$$
where in the second equality we apply second relation in (\ref{sk4}).
Thus, the $R$-matrix is completely invertible.
For the invertible matrix $Q$,
 the proof of the fact, that the Yang-Baxter $R$-matrix is
completely invertible, is similar. For
 invertible matrices $D$ and $Q$ we have
(\ref{sk8}) and one can rewrite relations (\ref{qtrs1}) as
 (\ref{sskk}). \hfill \qed

\vspace{0.3cm}
\noindent
{\bf Corollary 1.} Let $\R$ be skew-invertible
and the matrix $A_{12}$ be one of the matrices
$\{ \R_{12}, \, \R_{12}^{-1}, \, \hat{\Psi}_{12} \}$.
Then, from (\ref{sk4}), (\ref{sk6}) we obtain
\be
\lb{sk10}
[ A_{12}, \, D_1 \, D_2]  = 0 =
[ A_{12}, \, Q_1 \, Q_2  ] \; ,
\ee
\be
\lb{sk10a}
A_{12} \, (D \, Q)_1 = (D \, Q)_2 \, A_{12} \; .
\ee
If $\R$ is completely invertible, then
matrices $\hat{\Psi}_{12}$, $\hat{\Phi}_{12}$ are invertible
 $$
 \hat{\Psi}_{12}^{-1} =   D_1^{-1} \R_{21} D_2 =
Q_2^{-1} \R_{21}  Q_1  \, , \;\;\;
\hat{\Phi}_{12}^{-1} = D_2 \R_{21}^{-1} D_1^{-1}=
Q_1 \R_{21}^{-1}  Q_2^{-1}  \; .
 $$
In this case, by using
(\ref{sk5}) and (\ref{sk7}), we prove eqs. (\ref{sk10}),
(\ref{sk10a}) for $A_{12} = \hat{\Phi}_{12}$ and
deduce the relation on the matrices $\hat{\Phi}$ and $\hat{\Psi}$:
$$
\hat{\Phi}^{-1}_{12} = D^2_2 \, \hat{\Psi}_{12} \, D^{-2}_1  =
Q^2_1 \, \hat{\Psi}_{12} \, Q^{-2}_2 \; .
$$

\vspace{0.3cm}
\noindent
{\bf Corollary 2.} For any quantum $(N \times N)$ matrix $E$
(with noncommutative entries $E^i_j$) one can find the following identities
\be
\lb{RER}
Tr ( D \, E ) \, I_1  = Tr_{2} \left( D_2
\R_1^{\mp 1} \, E_1 \, \R_1^{\pm 1} \right) \, , \;\;
Tr (Q  \, E) \, I_2 = Tr_{1} \left( Q_1
\R_1^{\pm 1} \, E_2 \, \R_1^{\mp 1} \right),
\ee
that demonstrate the invariance properties of the quantum traces
\be
\lb{qtrDQ}
{\rm Tr} ( D \, E ) =: {\rm Tr}_{\cal D}(E) \; , \;\;\;\;
{\rm Tr} ( Q \, E ) =: {\rm Tr}_{\cal Q}(E) \; .
\ee
To prove identities (\ref{RER})
we multiply eqs. (\ref{4a}) by the matrix $E_0$
and take the trace $Tr_0(\dots)$.
Note that in view of (\ref{sk10}) the multiple quantum traces
satisfy cyclic property:
\be
\lb{cycl}
\begin{array}{c}
Tr_{{\cal D}(1 \dots m)} \left( X(\R) \cdot Y \right) =
Tr_{{\cal D}(1 \dots m)} \left(Y \cdot  X(\R) \right) \; , \\
Tr_{{\cal Q}(1 \dots m)} \left(X(\R) \cdot Y \right) =
Tr_{{\cal Q}(1 \dots m)} \left(Y \cdot X(\R) \right) \; ,
\end{array}
\ee
where $X(\R) \in {\rm End}(V^{\otimes m})$ denotes arbitrary element of the group algebra of the
braid group ${\cal B}_m$ in $R$-matrix representation (\ref{3.1.3}),
(\ref{3.1.3is})
and $Y \in {\rm End}(V^{\otimes m})$ are arbitrary
(quantum) operators.

\vspace{0.3cm}
\noindent
{\bf Corollary 3.} Let $R$ be completely
invertible matrix. We multiply the first and the second Yang-Baxter eqs. in (\ref{3.1.2is})
respectively from the right and left by the matrix $D_3$:
$$
\R_{2} \, \R^{-1}_{1} \, \R^{-1}_{2} \, D_3 =
\R^{-1}_{1} \, \R^{-1}_{2} \, \R_{1} \, D_3 \; , \;\;\;
D_3 \, \R_{1} \, \R^{-1}_{2} \, \R^{-1}_{1}   =
D_3 \, \R^{-1}_{2} \, \R^{-1}_{1} \, \R_{2}  \; ,
$$
and use relations (\ref{sk4}).
As a result we deduce
\be
\lb{ppr}
\R_{23} \, \hat{\Psi}_{21} \, \hat{\Psi}_{32} =
\hat{\Psi}_{21} \, \hat{\Psi}_{32} \, \hat{R}_{12} \; , \;\;\;
\R_{12} \, \hat{\Psi}_{32} \, \hat{\Psi}_{21} =
\hat{\Psi}_{32} \, \hat{\Psi}_{21} \, \hat{R}_{23} \; .
\ee
Analogously, if we multiply Yang-Baxter eqs. (\ref{3.1.4})
(for $X(\R_a) = \R_a^{-1}$ and $a=1$)
from the left and right  by the matrix $Q_1$
and use relations (\ref{sk7}),
we respectively deduce
\be
\lb{pph}
\hat{\Phi}_{21} \, \hat{\Phi}_{32} \, \hat{R}_{12} =
\R_{23} \, \hat{\Phi}_{21} \, \hat{\Phi}_{32} \; , \;\;\;
\hat{\Phi}_{32} \, \hat{\Phi}_{21} \, \hat{R}_{23} =
\R_{12} \, \hat{\Phi}_{32} \, \hat{\Phi}_{21} \; .
\ee

\vspace{0.2cm}
\noindent
{\bf Corollary 4.} The trace $Tr_{04}(\dots)$ of eq. (\ref{ybe1})
(or (\ref{ybe2})) gives
\be
\lb{yyy1}
Tr_1 Q_1 X(\R_1)  = Tr_3 D_3 X(\R_2) \equiv Y_2(X) \; ,
\ee
where we redefined the arbitrary function
 $X$: $X(\R) \R \rightarrow X(\R)$.
In particular, for $X=1$, we obtain $Tr(D) = Tr(Q)$.
 Eq. (\ref{yyy1}) leads to the following identity
\be
\lb{yis1}
 Tr_{12} \left( Q_1 Q_2 X(\R_1) \right) =
 Tr_{23} \left( D_3 Q_2 X(\R_2) \right) =
 Tr_{34} \left( D_3 D_4 X(\R_3) \right) \; .
\ee

\begin{proposition}\label{prop3a}
{\it For any polynomial
$X \in \mathbb{C}[\R_{1},\R_{1}^{-1}]$
the matrix $Y(X)$ defined in (\ref{yyy1})
satisfies $[D, \, Y] = 0 = [Y, \, Q]$  and
\be
\lb{yyy2}
Y_2(X) \; \R_1^{\pm 1} = \R_1^{\pm 1} \; Y_1(X) \; .
\ee
Matrices
\be
\lb{yyy2a}
Y_2^{(n)} := Y(\R^{n}) =
Tr_3 \left( D_3 \, \R_2^{n} \right) = Tr_1 \left( Q_1 \, \R_1^{n} \right)\; , \;\;\;\;\;\; \forall \; n \in \mathbb{Z} \; ,
\ee
generate a commutative set.}
\end{proposition}
{\bf Proof.} From (\ref{yyy1}) and (\ref{sk10}) we have
$$
 \begin{array}{c}
D_2 \, Y_2 = Tr_3 \bigl( D_2 \, D_3 \, X(\R_2) \bigr) =
Tr_3 \bigl( X(\R_2) \, D_2 \, D_3 \bigr) = Y_2 \, D_2 \; , \\ [0.2cm]
Q_2 \, Y_2 = Tr_1 \bigl( Q_1 \, Q_2 \, X(\R_1) \bigr) =
Tr_3 \bigl( X(\R_1) \, Q_1 \, Q_2 \bigr) = Y_2 \, Q_2 \; .
\end{array}
$$
The left hand side of (\ref{yyy2}) is transformed as following
$$
 \begin{array}{c}
Y_2(X) \; \R_1^{\pm 1} =
{\rm Tr}_3 \bigl( D_3 X(\R_2) \R_1^{\pm 1} \R_2^{\pm 1}
\R_2^{\mp 1} \bigr) = \R_1^{\pm 1} \; {\rm Tr}_3
 \bigl( D_3 \R_2^{\pm 1} X(\R_1) \R_2^{\mp 1} \bigr) = \\ [0.2cm]
=  \R_1^{\pm 1} \; Tr_2 \bigl( D_2 \, X(\R_1) \bigr) =
\R_1^{\pm 1} \; Y_1(X) \; ,
\end{array}
$$
where we used (\ref{3.1.4}) and first relation in (\ref{RER}).

The commutativity of the matrices $Y_2^{(n)}$ follows from (\ref{yyy2}), since
for $n$ is even and odd we have, respectively
$$
Y_2(X) \, Y_2^{(2k)} = Tr_3 \left( D_3 \, Y_2 \, \R_2^{2k} \right) =
Tr_3 \left( D_3  \, \R_2^{2k} \, Y_2 \right) =  Y_2^{(2k)} \, Y_2(X) \; ,
$$
$$
Y_2(X) \, Y_2^{(2k+1)} = Tr_1 \left( Q_1 \, Y_2  \, \R_1^{2k+1} \right)
= Tr_1 \left( Q_1  \, \R_1^{2k+1} \, Y_1 \right) =
Tr_1 \left(Y_1 \, Q_1  \, \R_1^{2k+1}  \right) =
$$
$$
 = Tr_1 \left( Q_1  \,  Y_1 \,  \R_1^{2k+1}  \right) =
Tr_1 \left( Q_1  \, \R_1^{2k+1} \, Y_2 \right) =
Y_2^{(2k+1)} \, Y_2(X) \; .
$$
For $X(\R)=\R^m$ $(m \in \mathbb{Z})$ we obtain commutativity of
matrices (\ref{yyy2a}). \hfill \qed


\begin{proposition}\label{prop3}
{\it The identity (\ref{yis1}) is generalized as
\be
\lb{yis2}
Tr_{1 \dots n} ( Q_1 \cdots Q_k D_{k+1} \cdots D_n X_{1 \to n} ) =
Tr_{1 \dots n} ( D_1 \cdots D_n X_{1 \to n} ) \;\;\;
(\forall n \geq 2 , \; k=1,...,n) \, ,
\ee
where $X_{1 \to n} := X(\R_1, \dots , \R_{n-1}) \in
\mathbb{C}[\R_1^{\pm 1},...,\R_{n-1}^{\pm 1}]$ is an element
 of the group algebra of the braid group ${\cal B}_{n}$
 in the $R$-matrix representation\footnote{In view of the graphical
 representation (\ref{figR}) any monomial in $X_{1 \to n}$
 is interpreted as a braid with $n$ strands.}
  (\ref{3.1.3}), (\ref{3.1.3is}).
}
\end{proposition}

\noindent
{\bf Proof.} Indeed, from (\ref{3.1.3}), (\ref{3.1.3is})
we have
$
\R_1 \cdots \R_{n} \, X_{1 \to n} =
X_{2 \to n+1} \, \R_1 \cdots \R_{n}$. Multiplying both
sides of this eq. by the matrices $Q_{1}$ and
$D_{n+1}$ from the left and right and taking the trace $Tr_{1}Tr_{n+1} \left( \dots \right)$ we
deduce (by means of (\ref{qtrs}))
$$
Tr_1 \left( Q_1 \, \R_1  \cdots \R_{n-1} \, X_{1 \to n} \right) =
Tr_{n+1} \left( D_{n+1} \, X_{2 \to n+1} \, \R_2 \cdots \R_{n} \right) \; ,
$$
which is written, after the redefinition
 $X_{1 \to n} \to (\R_1  \cdots \R_{n-1})^{-1} \, X_{1 \to n}$, in the form
\be
\lb{yis3}
Tr_1 \left( Q_1 \,  X_{1 \to n} \right) =
Tr_{n+1} \left( D_{n+1} \, (\R_2  \cdots \R_{n})^{-1} \,
X_{2 \to n+1} \, \R_2 \cdots \R_{n} \right) \; .
\ee
Then, applying the trace $Tr_2 (Q_2 \dots)$ to (\ref{yis3})
(and again using (\ref{yis3})) we obtain
\be
\lb{yis4}
\begin{array}{c}
Tr_{12} \left( Q_1 Q_2 \,  X_{1 \to n} \right) =
Tr_{n+1} D_{n+1} Tr_2 \left(  Q_2 \, (\R_2  \cdots \R_{n})^{-1} \,
X_{2 \to n+1} \, \R_2 \cdots \R_{n} \right)= \\[0.1cm]
= Tr_{n+1,n+2} \left( D_{n+1} D_{n+2} \, (\R_3  \cdots \R_{n+1})^{-2} \,
X_{3 \to n+2} \, (\R_3 \cdots \R_{n+1})^2 \right) \; .
\end{array}
\ee
Applying the trace $Tr_3 (Q_3 \dots)$ to (\ref{yis4}) etc. we obtain
\be
\lb{yis44}
\begin{array}{c}
Tr_{1 \dots k} \left( Q_1 \cdots Q_k \,  X_{1 \to n} \right) = \\ [0.1cm]
= Tr_{n+1 \dots n+k} \left( D_{n+1} \cdots D_{n+k} \, (\R_{k+1}  \cdots \R_{n+k-1})^{-k} \,
X_{k+1 \to n+k} \, (\R_{k+1} \cdots \R_{n+k-1})^k \right) \; ,
\end{array}
\ee
and finally multiplying the both sides of (\ref{yis44})
by $D_{k+1} \cdots D_n$ from the left and taking the trace
$Tr_{k+1 \dots n}$
 (applying $Tr_{k+1 \dots n}(D_{k+1} \cdots D_n \, \dots)$ to
 both sides of (\ref{yis44}))
we deduce (\ref{yis2})
\be
\lb{yis5}
\begin{array}{c}
Tr_{1 \dots n} \left( Q_1 \cdots Q_k D_{k+1} \cdots D_n \,  X_{1 \to n} \right) = \\[0.3cm]
= Tr_{k+1 \dots n+k} \left( D_{k+1} \cdots D_{n+k} \, (\R_{k+1}  \cdots \R_{n+k-1})^{-k} \,
X_{k+1 \to n+k} \, (\R_{k+1} \cdots \R_{n+k-1})^k \right) = \\[0.3cm]
= Tr_{k+1 \dots n+k}  \left( D_{k+1} \cdots D_{n+k} \,
X_{k+1 \to n+k}  \right) \; ,
\end{array}
\ee
where we have used the cyclic property (\ref{cycl}). \hfill \qed

\vspace{0.1cm}

\noindent
{\bf Remark 1.} A numerical $R$-matrix which is the
 image $(T \otimes T){\cal R}$ of the universal ${\cal R}$ matrix
 (\ref{unrab}) for the
quasitriangular Hopf algebra is obliged to be skew-invertible. Indeed, relations (\ref{close}) are written in the matrix form
$$
\begin{array}{c}
\delta^i_j \delta^k_\ell =
T^i_j(\alpha_{\mu}\alpha_{\nu})\;
T^k_\ell(\beta_\nu S(\beta_\mu)) =
R^{m k}_{j\, n}\;\;
T^i_m(\alpha_{\mu}) \, T^n_\ell( S(\beta_\mu)) \; , \\ [0.2cm]
\delta^k_\ell \delta^i_j  =
T^k_\ell(\alpha_{\mu}\alpha_{\nu})\;
T^i_j(S(\beta_\nu) \beta_\mu) =
T^n_\ell(\alpha_{\nu})\, T^i_m(S(\beta_\nu)) \;\;
R^{km}_{nj} \; ,
\end{array}
$$
and, thus, relations (\ref{close}) are the algebraic
 counterparts of (\ref{skewb}),
where the matrix $\Psi$ is given by the equation
\be
\lb{psi}
\Psi^{in}_{ml} = T^i_m( \alpha_{\mu}) \, T^n_l( S(\beta_\mu) ) = \hat{\Psi}^{ni}_{ml} \; .
\ee
Moreover, in view of (\ref{contg}), the transposed
matrix $\Psi^{t_2}$ of (\ref{psi}) is interpreted as the image
$(T \otimes \overline{T}){\cal R}$, where $\overline{T}$  denotes a contragredient
representation to $T$, i.e. $\overline{T}(a)=T^t(S(a))$
$(\forall a \in {\cal A})$. Then, the second equation in (\ref{ppr}) is
nothing but the image of the universal
Yang-Baxter equation (\ref{2.30}) in the representation
$(T\otimes T \otimes \overline{T})$.

The image $(T \otimes T){\cal R}^{-1}= R^{-1}$ is also skew-invertible. The matrix
$\Phi_{12}$ in (\ref{skew2}) is given by
\be
\lb{phi}
\Phi^{in}_{ml} = T^i_m(S(\gamma_\mu)) \, T^n_l(\delta_{\mu}) = \hat{\Phi}^{in}_{lm} \; .
\ee
and the algebraic counterpart of (\ref{skew2}) is (\ref{close1}).
The second equation in (\ref{pph}) is the image of the universal
equation
${\cal R}_{23} {\cal R}_{12}^{-1} {\cal R}_{13}^{-1}=
{\cal R}_{13}^{-1} {\cal R}_{12}^{-1}{\cal R}_{23}$ (see (\ref{2.30})) in the representation
$(\overline{T} \otimes T \otimes T)$.
From eqs. (\ref{psi}), (\ref{phi})
 we also have the universal formulas for
matrices (\ref{dmatr}), (\ref{dmatr1}) of quantum traces
\be
\lb{dmatr2}
\begin{array}{c}
D = T(S(\beta_\mu)\alpha_\mu) = T(u) \; , \;\;\;
\overline{D} = T(S(\gamma_\mu)\delta_\mu) = T(u_2) \; , \\ \\
Q = T(\alpha_\mu S(\beta_\mu)) = T(S(u_3)) = T(u_4^{-1}) \; , \;\;\;
\overline{Q} = T(\delta_\mu S(\gamma_\mu)) = T(S(u_4)) = T(u_3^{-1})\; ,
\end{array}
\ee
where elements $u_1=u,u_2,u_3,u_4$ were introduced in
(\ref{u1234}) in Subsection {\bf \ref{trqH}}. Then, in view of
 Proposition {\bf \ref{prop1b}}, all matrices
 (\ref{dmatr2}) commute with each other and the products
 $DQ$, $D\overline{Q}=I$, $\overline{D}Q=I$ and
 $\overline{D}\, \overline{Q}=(DQ)^{-1}$
  are images of central elements $(u_i u_j^{-1}) \in {\cal A}$
  in the representation $T$.


\subsubsection{$R$-matrix formulation of
 link and knot invariants}

The $R$-matrix formulation of
link and knot invariants was developed in
\cite{16'} , \cite{TuRo}, \cite{18} (see also references therein).
Taking into account the fact that $R$-matrices satisfy (by definition)
 the third Reidemeister move (\ref{figR3}) we see that
Propositions {\bf \ref{prop2}} and {\bf \ref{prop3}}
are important for
constructing of link and knot invariants.
 Indeed, using graphic
representation (\ref{figR}), one can visualize relations
(\ref{qtrs}) and (\ref{sskk}) from Proposition {\bf \ref{prop2}}
 as the first Reidemeister moves:

\unitlength=5.5mm
\begin{picture}(17,3.5)
\put(0,1.9){$ Tr_2 \left( \R_{12} D_2 \right) =$}
\put(5.3,3){\vector(1,-1){2}}
\put(6.5,2.2){\line(1,1){0.8}}
\put(6.1,1.8){\vector(-1,-1){0.8}}

\put(7.3,2){\oval(1,2)[r]}
\put(7.8,1.7){\vector(0,1){0.4}}

\put(8,1.9){$D$}
\put(9.5,1.9){$=$}
\put(11,3){\vector(0,-2){2}}
\put(11.2,1.9){$I_{1}\;\;$,}


\put(14,1.9){$ Tr_2 \left(\R_{12}^{-1} Q_2^{-1} \right) =$}
\put(20,3){\line(1,-1){0.8}}
\put(21,2){\vector(-1,-1){1}}
\put(21.2,1.8){\vector(1,-1){0.8}}
\put(21,2){\line(1,1){1}}
\put(22,2){\oval(1,2)[r]}
\put(22.5,1.7){\vector(0,1){0.4}}

\put(22.7,1.9){$Q^{-1}$}

\put(24.5,1.9){$=$}
\put(25.8,3){\vector(0,-2){2}}
\put(26,1.9){$I_{1}\;\;$,}

\end{picture}
\vspace{-1cm}
\be
\lb{figa5}
{}
\ee

\unitlength=5.5mm
\begin{picture}(17,3.5)
\put(-0.2,1.9){$ Tr_1 \left( Q_1 \R_{12} \right) =$}
\put(6.5,3){\vector(1,-1){2}}
\put(7.7,2.2){\line(1,1){0.8}}
\put(7.3,1.8){\vector(-1,-1){0.8}}

\put(6.5,2){\oval(1,2)[l]}
\put(6,1.7){\vector(0,1){0.4}}
\put(5.2,1.9){$Q$}

\put(9.2,1.9){$=$}
\put(10.7,3){\vector(0,-2){2}}
\put(11,1.9){$I_{2}\;\;$,}


\put(14,1.9){$ Tr_2 \left(D_1^{-1} \R_{12}^{-1} \right) =$}
\put(22,3){\line(1,-1){0.8}}
\put(23,2){\vector(-1,-1){1}}
\put(23.2,1.8){\vector(1,-1){0.8}}
\put(23,2){\line(1,1){1}}
\put(22,2){\oval(1,2)[l]}
\put(21.5,1.7){\vector(0,1){0.4}}

\put(20,1.9){$D^{-1}$}

\put(24.5,1.9){$=$}
\put(25.8,3){\vector(0,-2){2}}
\put(26,1.9){$I_{2}\;\;$.}

\end{picture}
\vspace{-1cm}
\be
\lb{figa6}
{}
\ee
\noindent
These pictures show that
the elementary braids $\R$ and $\R^{-1}$ are closed
by matrices $D$, $Q^{-1} = \overline{D}$ on the right, and
 by matrices $D^{-1}$, $Q=\overline{D}^{\, -1}$
  on the left, to obtain trivial braids. We note that in general
 $D \neq Q^{-1}$. We stress however that for many
explicit numerical $R$-matrices we
 have\footnote{For $R = (T \otimes T){\cal R}$
 the matrix $DQ= T(u_1 u_4^{-1})$ is the
 image of central element and for irreducible representation
 $T$ we have $DQ \sim I$; see also Examples 1 and 2 below.}
 $Q^{-1} \sim D$ and therefore, after the special normalization
 of $R$-matrices, we deal with the standard first
 Reidemeister moves.
Finally, for the case of the skew-invertible $R$-matrices,
Proposition {\bf \ref{prop3}} demonstrates the equivalence of the
complete closuring of braids\footnote{Here the braids
$X_{1 \to n}$ are elements of the
braid group ${\cal B}_n$ in the $R$-matrix representations.}
 $X_{1 \to n}$ from the left and from the right
 by means of the quantum traces respectively
 with matrices $Q$ and $D$. Thus, for any braid $X_{1\to n}$
 with $n$ strands ($X_{1\to n}$ is a monomial constructed
  as a product of any number of $R$-matrices
 $\{\R_1,..., \R_{n-1}\}$) the characteristic (\ref{yis2})
 \be
 \lb{qxinv}
 {\sf Q}(X_{1 \to n}) :=
 Tr_{1 \dots n} ( Q_1 \cdots Q_n X_{1 \to n} ) \equiv
Tr_{1 \dots n} ( D_1 \cdots D_n X_{1 \to n} )  \; ,
 \ee
  gives (by closing
  of the braid $X_{1\to n}$) the invariant for link/knot.

 \vspace{0.2cm}

 \noindent
 {\bf Remark 2.} Let $T$ be the
  representation
  of the quasitriangular Hopf algebra
  ${\cal A}$ in the space $V$. Consider a special matrix representation
 of the universal ${\cal R}$-matrix
  \be
  \lb{TnR1}
  R_{(k,m)} = \sum_\mu \; T^{\otimes k}(\alpha_{\mu}) \otimes
  T^{\otimes m}(\beta_\mu) \equiv
  (T^{\otimes k} \otimes T^{\otimes m})\, {\cal R} \; ,
  \ee
  where $T^{\otimes k}$ acts to the first factor in ${\cal R}$,
  $T^{\otimes m}$ acts to the second factor in ${\cal R}$
  and we have used the notation (\ref{unrab}).
  Then applying (\ref{2.27}) we deduce
  \be
  \lb{TnR2}
  R_{(k,m)} = R_{1\to k; k+m} \cdots R_{1\to k; k+2}
  \cdot R_{1\to k; k+1}
  = (P_{1\to k; k+m} \cdots  P_{1\to k; k+1})
  \; \hat{R}_{(k,m)} \; ,
  \ee
    \be
  \lb{TnR3}
  \hat{R}_{(k,m)} = \R_{(m \to k+m-1)} \cdots \R_{(2 \to k+1)} \, \R_{(1 \to k)}  \; ,
\ee
where
  $$
  \begin{array}{c}
  R_{1\to k; k+ \ell} :=
  R_{1, k+ \ell} \, R_{2, k+ \ell}
  \cdots R_{k, k+ \ell} =
  P_{1\to k; k+ \ell} \cdot \R_{(1 \to k-1)}
  \R_{k, k+ \ell} \; ,   \\ [0.2cm]
     P_{1\to k; k+ \ell} := P_{1, k+ \ell} \, P_{2, k+ \ell}
   \cdots P_{k, k+ \ell} \, , \;\;\;
   \R_{(k \to \ell)} := \hat{R}_{k} \, \hat{R}_{k+1}
\cdots \hat{R}_{\ell}  \; ,
  \end{array}
  $$
   $R_{ij} := ((T \otimes T){\cal R})_{ij}$
  and the braid $\hat{R}_{(k,m)}$ is obtained from matrix
 $R_{(k,m)}$ by substitution $R_{ij} = P_{ij} \, \hat{R}_{ij}$
 and shifting all permutation matrices $P_{ij}$ to the left.
  The braid $\hat{R}_{(k,m)}$ defined in (\ref{TnR3})
  can be visualized, by means of (\ref{figR}), as the intersection
  of two cables (or two ribbons) with $m$ and $k$ strands

  \unitlength=10mm
\begin{picture}(17,4)

\put(0.5,1.9){$\hat{R}_{(k,m)} \;\;  =$}

\put(3,3){$\bullet$}
\put(5.2,3){$\bullet$}

\put(4.1,3){$\bullet$}
\put(4.15,3.15){\vector(2,-1){4}}
\put(3.5,3.1){$\dots$}
\put(4.5,3.1){$\dots$}
\put(5.7,1.8){$\ddots$}
\put(6.2,2.1){$\ddots$}

\put(7.7,3){$\dots$}
\put(7,3){$\bullet$}
\put(9.2,3){$\bullet$}

\put(3,3.4){$_1$}

\put(5.3,3.4){$_m$}
\put(6.7,3.4){$_{m+1}$}
\put(8.7,3.4){$_{m+k}$}

\put(7.15,3.15){\line(-2,-1){0.7}}
\put(9.2,3.1){\line(-2,-1){1.7}}


\put(6.15,2.65){\line(-2,-1){0.4}}
\put(7.2,2.1){\line(-2,-1){0.4}}

\put(4.8,2){\vector(-2,-1){1.6}}
\put(6.13,1.6){\vector(-2,-1){0.85}}

\put(3.2,3.15){\vector(2,-1){4}}
\put(5.3,3.15){\vector(2,-1){4}}

\put(3,1){$\bullet$}
\put(5.2,1){$\bullet$}
\put(3.9,1){$\dots$}
\put(7,1){$\bullet$}
 \put(9.2,1){$\bullet$}

 \put(8,1){$\bullet$}
 \put(7.4,1.1){$\dots$}
 \put(8.5,1.1){$\dots$}

 \put(10.5,1.9){\large \bf .}




\end{picture}
\vspace{-2cm}
\be
\lb{figas1b}
{}
\ee

\vspace{0.3cm}
 \noindent
 This pictorial presentation demonstrates
the fact that in general the $R$-matrix approach could describe
invariants not only for ordinary links and knots,
but also for ribbon (cable) links and knots.
In this case the right (or left) closuring for braids with
 matrices $D_1 \cdots D_n \sim Q_1^{-1} \cdots Q_n^{-1}$
 (or $Q^{\otimes n} \sim (D^{-1})^{\otimes n}$) are
 also different for the cable (ribbon) braids
 $\R_{(n,n)}$ and $\R_{(n,n)}^{-1}$ (cf. (\ref{figa5}), (\ref{figa6})).
 For example, for right closuring
  it follows from the visualization of the moves
 which are shown on the pictures:

\unitlength=5mm
\begin{picture}(17,6)(-3,-0.5)

\qbezier(-0.8,0.9)(2.1,6.6)(3,2.5)
\qbezier(-0.3,0.9)(2.1,5.6)(2.6,2.5)

\qbezier(0.5,2)(3,0)(3,2.5)
\qbezier(0.8,2.3)(2.5,1)(2.6,2.5)

\qbezier(-1,4.5)(-1,3.5)(-0.1,2.6)
\qbezier(-0.6,4.6)(-0.6,3.6)(0.1,2.9)

\put(-0.5,1.2){\vector(0,-1){1}}
\put(-0.8,4.1){\vector(0,1){1}}

\put(3.2,2.5){\scriptsize $(Q^{-1})^{\otimes n}$}

\put(5.9,2.5){$=$}

\qbezier(7.2,0.3)(7.2,1.6)(7.5,1.9)
\qbezier(7.5,1.9)(8,2.7)(7.55,3.45)
\qbezier(7.3,3.8)(7,4.2)(7.1,5.1)

\qbezier(7.9,0.3)(7.9,1.4)(7.7,1.7)
\qbezier(7.4,2.1)(6.95,2.8)(7.4,3.6)
\qbezier(7.4,3.6)(7.8,4.1)(7.7,5)


\put(8.5,2.5){$\equiv \; {\sf D}$}






\put(10.5,2.5){\Large \bf ,}

\qbezier(12.2,4.1)(15.1,-1.6)(16,2.5)
\qbezier(12.7,4.1)(15.1,-0.6)(15.6,2.5)

\qbezier(13.5,3)(16,5)(16,2.5)
\qbezier(13.8,2.7)(15.5,4)(15.6,2.5)

\qbezier(12,0.5)(12,1.5)(12.9,2.4)
\qbezier(12.4,0.4)(12.4,1.4)(13.1,2.1)

\put(12.2,0.9){\vector(0,-1){1}}
\put(12.5,3.9){\vector(0,1){1}}

\put(16.2,2.5){\scriptsize $D^{\otimes n}$}

\put(17.9,2.5){$=$}

\qbezier(19.2,5.1)(19.2,3.8)(19.5,3.5)
\qbezier(19.5,3.5)(20,2.7)(19.55,1.95)
\qbezier(19.3,1.6)(19,1.2)(19.1,0.3)

\qbezier(19.9,5.1)(19.9,4)(19.7,3.7)
\qbezier(19.4,3.3)(18.95,2.6)(19.4,1.8)
\qbezier(19.4,1.8)(19.8,1.3)(19.7,0.4)

\put(20.5,2.5){$\equiv \; \overline{\sf D}= {\sf D}^{-1}$}

\put(25,2.5){\Large \bf ,}

\end{picture}
\vspace{-1.5cm}
\be
\lb{figa7}
{}
\ee

\noindent
where we pull the ribbons along the arrows
on the left hand side of the equalities and obtain
two differently twisted ribbons (as spirals)
in the right hand side of the equalities. Thus, for
ribbon (cable) links/knots, to obtain the
first Reidemeister moves, we need to insert
matrices ${\sf D}$ and $\overline{\sf D}$
in the closuring of braids $D^{\otimes n}\cdot {\sf D}$ and
 $(Q^{-1})^{\otimes n}\cdot \overline{\sf D}$
 (here the "ribbon" matrices ${\sf D}$ and
 $\overline{\sf D}$  are defined in (\ref{figa7})):

 \unitlength=5mm
\begin{picture}(17,6)(-3,-0.5)


\qbezier(1.2,0.9)(4.3,6.6)(5.3,2.5)
\qbezier(1.7,0.9)(4.1,5.6)(4.6,2.5)

\qbezier(2.5,2)(5.2,0)(5.3,2.5)
\qbezier(2.8,2.3)(4.5,1)(4.6,2.5)

\qbezier(1,4.5)(1,3.5)(1.9,2.6)
\qbezier(1.4,4.6)(1.4,3.6)(2.1,2.9)

\put(1.5,1.2){\vector(0,-1){1}}
\put(1.2,4.1){\vector(0,1){1}}

\put(4.7,2.3){\scriptsize $\overline{\sf D}$}
\put(4.4,3.3){\line(1,0){0.6}}
\put(4.4,1.9){\line(1,0){0.75}}

\put(4.6,4.1){\scriptsize $(Q^{-1})^{\otimes n}$}

\put(6.4,2.5){$=$}

\qbezier(7.7,5.1)(7.7,3.8)(8,3.5)
\qbezier(8,3.5)(8.5,2.7)(7.8,1.6)
\qbezier(7.8,1.6)(7.5,1.2)(7.6,0.3)

\qbezier(8.2,5.1)(8.2,3.8)(8.5,3.5)
\qbezier(8.5,3.5)(9,2.7)(8.3,1.6)
\qbezier(8.3,1.6)(8,1.2)(8.1,0.3)


\put(9.5,2.5){\Large \bf ,}

\qbezier(13.2,4.1)(16.2,-1.6)(17.3,2.5)
\qbezier(13.7,4.1)(16.1,-0.6)(16.6,2.5)

\qbezier(14.6,3)(17.1,5)(17.3,2.5)
\qbezier(14.8,2.7)(16.5,4)(16.6,2.5)

\qbezier(13,0.5)(13,1.5)(13.9,2.4)
\qbezier(13.4,0.4)(13.4,1.4)(14.1,2.1)

\put(13.2,0.9){\vector(0,-1){1}}
\put(13.5,3.9){\vector(0,1){1}}

\put(16.7,2.3){\scriptsize ${\sf D}$}
\put(16.4,1.8){\line(1,0){0.6}}
\put(16.45,3.1){\line(1,0){0.7}}

\put(17,3.8){\scriptsize $D^{\otimes n}$}

\put(17.9,2.5){$=$}

\qbezier(19.2,5.1)(19.2,3.8)(19.5,3.5)
\qbezier(19.5,3.5)(20,2.7)(19.3,1.6)
\qbezier(19.3,1.6)(19,1.2)(19.1,0.3)

\qbezier(19.7,5.1)(19.7,3.8)(20,3.5)
\qbezier(20,3.5)(20.5,2.7)(19.8,1.6)
\qbezier(19.8,1.6)(19.5,1.2)(19.6,0.3)

\put(21,2.5){\Large \bf .}

\end{picture}
\vspace{-1.2cm}
\be
\lb{figa8}
{}
\ee

\noindent
Thus, in the right hand side of the relations (\ref{figa8}), we
 obtain the unit operators in $V^{\otimes n}$.

\subsubsection{Spectral decomposition of $R$-matrices
and examples of knot/link invariants}
We now assume that the invertible Yang-Baxter
 $\R$ matrix obeys the characteristic equation
\begin{equation}
\lb{3.1.27}
(\R - \mu_{1})(\R - \mu_{2}) \cdots
(\R - \mu_{M}) = 0 \; ,
\end{equation}
where $\mu_{i} \in \mathbb{C}$,
$\mu_{i} \neq \mu_{j}$ if $i \neq j$ and
$\mu_i \neq 0$ $\forall i$. This equation can be represented in the form
\be
\lb{chars}
\R^M - \sigma_1(\mu) \, \R^{M-1} + \dots
+ (-1)^{M-1} \, \sigma_{M-1}(\mu) \, \R +
(-1)^M \, \sigma_M(\mu) \; {\bf 1} = 0 \; ,
\ee
where  ${\bf 1}$ is a unit matrix in $V^{\otimes 2}$ and
$$
\sigma_k(\mu) = \sum_{i_1 < i_2 < \dots < i_k} \mu_{i_1} \dots \mu_{i_k} \; ,
$$
are elementary symmetric polynomials of $\mu_i$ $(i=1, \dots , M)$.
For $R$-matrices satisfying (\ref{3.1.27}), one can introduce a
complete set of $M$ projectors:
\begin{equation}
\lb{3.1.28}
\hbox{\bf P}_{k} = \prod\limits_{j \neq k}
\frac{(\R - \mu_{j})}{(\mu_{k} -\mu_{j})} \; , \;\;\;\;\;\;
\sum_k \, \hbox{\bf P}_{k} = {\bf 1} \; ,
\end{equation}
which project the $\R$ matrix to its eigenvalues
$\hbox{\bf P}_{k} \, \R = \R \, \hbox{\bf P}_{k} = \mu_k \, \hbox{\bf P}_{k}$
and can be used for the spectral decomposition
of an arbitrary function $X$ of $R$:
\be
\lb{3.1.29}
X(\R) = \sum_{k=1}^{M} X(\mu_{k}) \hbox{\bf P}_{k} \; .
\ee
In particular, for
$X={\bf 1}$ we obtain the completeness condition (see second
equation in (\ref{3.1.28})). The derivation of formulas
(\ref{3.1.28}) can be found, for example, in \cite{IsRub2},
\cite{Cvit}.

In the calculations
of the knot/link invariants (\ref{qxinv}), the characteristic
equations (\ref{chars}) play the role of the skein relations.
We also note, that for many known explicit examples of the
completely invertible Yang-Baxter $\R$ matrices,
which satisfy the characteristic equation
(\ref{3.1.27}), all matrices $Y(X)$, defined in (\ref{yyy1}),
 are proportional
to the identity matrix (see Proposition 4 in \cite{IsBonn}).

\vspace{0.2cm}

\noindent
{\bf Examples.} Here we consider two special cases $M=2,3$ for the
characteristic equation (\ref{3.1.27}).
By renormalizing the matrix $\R$, it is
always possible to fix first two eigenvalues
in (\ref{3.1.27}) so that $\mu_1 = q \neq 0$ and
$\mu_2 = - q^{-1} \neq 0$, where $q \in \mathbb{C}$. \\
{\bf 1.} For $M=2$ Eqs. (\ref{3.1.27}),
(\ref{chars}) are represented in the form of the Hecke condition
\begin{equation}
\lb{3.3.7}
\begin{array}{c}
(\R - q)(\R + q^{-1}) = 0 \;\; \Rightarrow \;\;
\R^2 = \lambda \, \R + {\bf 1}
\;\; \Rightarrow \;\; \\ [0.2cm]
\R - \lambda \, {\bf 1}  - \R^{-1} = 0 \; , \;\;\;
\lambda := (q-q^{-1}) \; .
\end{array}
\end{equation}
In this case, for all $n \in \mathbb{Z}$ we obtain
 \be
 \lb{Rn}
 \R^n = \alpha_n \R + \alpha_{n-1} {\bf 1} \; , \;\;\;\;\;\;\;
 \alpha_n := \frac{q^{n}- (-q)^{-n}}{q+q^{-1}} \; ,
 \ee
and according to (\ref{qtrs})
all matrices $Y(\R^n) \equiv Y^{(n)}$ in (\ref{yyy2a})
 are proportional to the
identity matrix
\be
 \lb{trRn}
Tr_2 (D_2 \R^n_{12}) = (\alpha_n  + \alpha_{n-1} Tr(D)) \, I_1 \; .
\ee
In particular one
can immediately find (see (\ref{qtrs}), (\ref{sk9}))
\be
 \lb{QD}
Tr (Q) = Tr (D) \, , \;\;\;\;\;\;
Y_1(\R^{-1}) = Q_1 \, D_1 = Tr_2 (D_2 \R_1^{-1}) =
\left(1 - \lambda \, Tr(D) \right) \, I_1 =
q^{-2d} \,  I_1 \; ,
 \ee
 where we introduce useful parametrization $q^{-2d} = (1 - \lambda \, Tr(D))$.
Equation (\ref{QD}) means that for the skew invertible Hecke $R$-matrix, in
the case $\lambda Tr(D) \neq 1$, the matrices $D$ and $Q$
are always invertible and
$Q^{-1} = q^{2d} \, D $.

\noindent
{\bf 2.} For $M=3$, equations (\ref{3.1.27}),
(\ref{chars}) are the
Birman-Murakami-Wenzl cubic relations
(cf. eq. (\ref{3.7.1}) below)
 \be
 \lb{3bmw}
 \begin{array}{c}
(\R - q)(\R + q^{-1})(\R -\nu) = 0 \;\; \Rightarrow \;\;
\R^3 - (\lambda +\nu) \, \R^2  + (\lambda \nu -1) \R +
 \nu {\bf 1} = 0  \;\; \Rightarrow  \\ [0.2cm]
 \hat{K}\R = \R \hat{K} = \nu \hat{K} \; , \;\;\;\;
 \hat{K} := \frac{1}{\lambda \nu} (q-\R)(q^{-1} + \R) =
 \frac{1}{\lambda \nu}\Bigl( {\bf 1} + \lambda \R
  -  \R^2 \Bigr) \; .
 \end{array}
 \ee
 where $\lambda = (q- q^{-1})$.
In this case we have
  \be
 \lb{3bmwk2}
 \hat{K}^2 = \mu \; \hat{K} \; , \;\;\;\;\; \mu :=
 \frac{1}{\lambda}(\nu^{-1} + \lambda - \nu) \; ,
 \ee
and for all $n \in \mathbb{Z}$ we obtain
 \be
 \lb{Rn2}
 \R^n = \alpha_n \R + \alpha_{n-1} {\bf 1} + \beta_n \hat{K} \; ,
 \ee
 $$
 \alpha_n := \frac{q^{n}- (-q)^{-n}}{q+q^{-1}} \; , \;\;
 \beta_n := \frac{\lambda \nu}{q+q^{-1}}
  \left(\frac{(\nu^n -  (-q)^{-n})}{(\nu + q^{-1})} -
 \frac{(\nu^n - q^{n})}{(\nu -q)} \right) \; .
 $$
 Let the matrix $\hat{K}$ be a one-dimensional projector
 in $V^{\otimes 2}$, i.e.
 $\hat{K}^{i_1i_2}_{k_1 k_2}= \bar{C}^{i_1i_2} C_{k_1 k_2}$.
 In this case one can define
 the quantum trace (\ref{yyy1}) as following
 (see eq. (\ref{3.7.29}) in Section {\bf \ref{qBCD}} below)
 $$
 \hat{K}_{23} \, X(\R_{12}) \, \hat{K}_{23} =
 \nu^{-1} Tr_2\bigl( X(\R_{12}) D_2 \bigr) \; \hat{K}_{23} \; ,
 \;\;\;\;\; D^i_{\; j} : = \nu \, \bar{C}^{ik} C_{jk} \; ,
 $$
 and we deduce
 $$
 \begin{array}{c}
 \hat{K}_{23} \, \R_{12} \, \hat{K}_{23} = \nu^{-1} \hat{K}_{23} \; ,
 \;\;\; \hat{K}_{23} \, \hat{K}_{12} \, \hat{K}_{23} = \hat{K}_{23}
 \; , \\ [0.2cm]
 \hat{K}_{23} \, {\bf 1} \, \hat{K}_{23}  =
 \frac{1}{\lambda}(\nu^{-1} + \lambda - \nu) \, \hat{K}_{23} =
 \nu^{-1} Tr(D) \hat{K}_{23} \;\;\;  \Rightarrow \;\;\;
 Tr(D) = \frac{(q-\nu)(q^{-1} + \nu)}{\lambda} \; .
 \end{array}
 $$
 Using these relations and (\ref{Rn2}) we obtain
 \be
 \lb{trRn1}
 Tr_2 (\R^n_{12} D_2) =
 \Bigl( \alpha_n  +
 \frac{(q-\nu)(q^{-1} + \nu)}{\lambda} \, \alpha_{n-1}
 + \nu \beta_n \Bigr) \, I_1 \; ,
 \ee
 where $\alpha_n$ and $\beta_n$ were introduced in (\ref{Rn2}).
Thus, for the cubic characteristic equation (\ref{3bmw}),
all matrices $Y(\R^n)$ (\ref{yyy2a})
 are also proportional to the identity matrix and for $n=-1$
 we find $Q^{-1} = \nu^{-2} \, D$.

\begin{figure}[h]
 \unitlength=4.4mm
\begin{picture}(17,13)(0,-4.3)

\put(8,3){$n$}
\put(9.2,3.2){\oval(1,8.5)[l]}

\put(14.3,-1.1){\line(0,1){8.3}}

 \put(10.2,-3){\line(0,1){2}}
 \put(10.2,7.4){\line(0,1){2}}

\put(9.7,8.8){\scriptsize $1$}
\put(11.7,7.5){\scriptsize $2$}

\put(9.7,-2.8){\scriptsize $1$}
\put(11.7,-2){\scriptsize $2$}

\put(10,7){$\bullet$}
\put(12,7){$\bullet$}

\put(10.15,7.2){\line(1,-1){1}}
\put(11.0,6){\vector(-1,-1){0.8}}
\put(11.15,6.2){\vector(1,-1){1}}
\put(11.4,6.4){\line(1,1){0.8}}

\put(10,5){$\bullet$}
\put(12,5){$\bullet$}


\put(10.15,5.2){\line(1,-1){1}}
\put(11.0,4){\vector(-1,-1){0.8}}
\put(11.15,4.2){\vector(1,-1){1}}
\put(11.4,4.4){\line(1,1){0.8}}

\put(10,3){$\bullet$}
\put(12,3){$\bullet$}


\put(11,1.8){$\vdots$}


\put(10,1){$\bullet$}
\put(12,1){$\bullet$}

\put(10.15,1.2){\line(1,-1){1}}
\put(11.0,0){\vector(-1,-1){0.8}}
\put(11.15,0.2){\vector(1,-1){1}}
\put(11.4,0.4){\line(1,1){0.8}}

\put(10,-1){$\bullet$}
\put(12,-1){$\bullet$}


\put(13.3,7.1){\oval(2,2)[t]}

\put(13.3,-0.9){\oval(2,2)[b]}

\put(16,-2.3){\bf A}


\put(23,3){$n$}
\put(24.2,3.2){\oval(1,8.5)[l]}

\put(29.3,-1.1){\line(0,1){8.3}}
 \put(31.3,-1.1){\line(0,1){8.3}}

\put(24.8,7.5){\tiny $1$}
\put(26.7,7.5){\tiny $2$}

\put(24.8,-1.5){\tiny $1$}
\put(26.7,-1.5){\tiny $2$}

\put(25,7){$\bullet$}
\put(27,7){$\bullet$}

\put(25.15,7.2){\line(1,-1){1}}
\put(26.0,6){\vector(-1,-1){0.8}}
\put(26.15,6.2){\vector(1,-1){1}}
\put(26.4,6.4){\line(1,1){0.8}}

\put(25,5){$\bullet$}
\put(27,5){$\bullet$}


\put(25.15,5.2){\line(1,-1){1}}
\put(26.0,4){\vector(-1,-1){0.8}}
\put(26.15,4.2){\vector(1,-1){1}}
\put(26.4,4.4){\line(1,1){0.8}}

\put(25,3){$\bullet$}
\put(27,3){$\bullet$}


\put(26,1.8){$\vdots$}


\put(25,1){$\bullet$}
\put(27,1){$\bullet$}

\put(25.15,1.2){\line(1,-1){1}}
\put(26.0,0){\vector(-1,-1){0.8}}
\put(26.15,0.2){\vector(1,-1){1}}
\put(26.4,0.4){\line(1,1){0.8}}

\put(25,-1){$\bullet$}
\put(27,-1){$\bullet$}


\put(28.3,7.1){\oval(6,4)[t]}
\put(28.3,7.1){\oval(2,2)[t]}

\put(28.3,-0.9){\oval(6,4)[b]}
\put(28.3,-0.9){\oval(2,2)[b]}

\put(33.5,-2.3){\bf B}

\end{picture}

 \vspace{-1cm}

\caption{\label{zamyk}{\em Closure of the
 braid $\R^n$ (the right picture {\bf B})
 gives toroidal knots for odd $n$ and
 links for even $n$.} }
\end{figure}
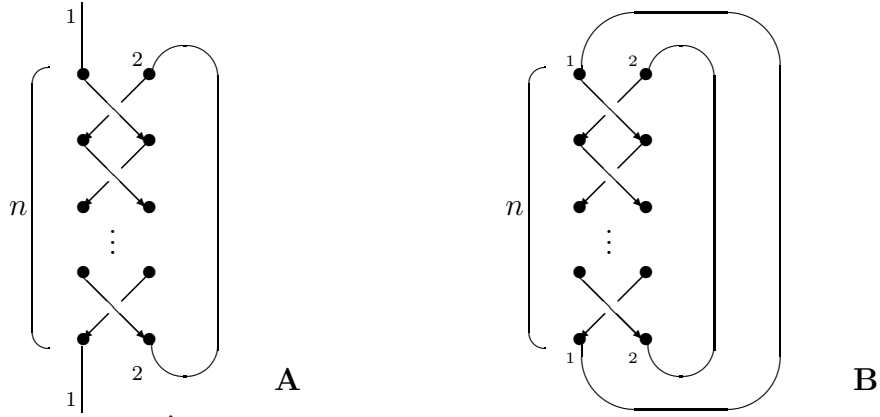

 \noindent
{\bf Remark 3.} Equations (\ref{trRn}) and (\ref{trRn1}) are
visualized  in {\bf Fig.\ref{zamyk}A} and give
(for the cases $M=2,3$) invariants of links and
knots (\ref{qxinv})
 $$
 \begin{array}{l}
 M=2: \;\;\;  {\sf Q}(\R^n_{12}) =
 Tr_{12} (\R^n_{12} D_1 D_2) =
  \bigl(\alpha_n  + \alpha_{n-1} Tr(D) \bigr) \, Tr(D) \; ,
 \\ [0.3cm]
 M=3: \;\;\; {\sf Q}(\R^n_{12}) =
 Tr_{12} (\R^n_{12} D_1 D_2) =
 \Bigl( \alpha_n  +
 \frac{(q-\nu)(q^{-1} + \nu)}{\lambda} \, \alpha_{n-1}
 + \nu \beta_n \Bigr) \frac{(q-\nu)(q^{-1} + \nu)}{\lambda} \; .
 \end{array}
 $$
 which are presented in {\bf Fig.\ref{zamyk}B}

\noindent
The explicit examples of $R$-matrices subject
to (\ref{3.3.7}) and
(\ref{3bmw}) with fixed values $Tr(D)$
and $\nu$ are given in Sections  {\bf \ref{qgrsl}}, {\bf \ref{qsuper}}
and {\bf \ref{qBCD1}}, {\bf \ref{mpsosp2}}, {\bf \ref{ospsup}}
 below.


\subsection{\bf \em Quantum matrix algebras\label{qmatalg}}
\setcounter{equation}0

\subsubsection{$RTT$ algebras}

We consider an algebra ${\cal A}^*$ of functions on a quasitriangular Hopf algebra ${\cal A}$ and
assume that generators of ${\cal A}^*$ are the identity element
$1$ and elements of $N \times N$ matrix $T=||T^{i}_{j}||$ $(i,j=1, \dots ,N)$, which
define $N$- dimensional matrix representation of ${\cal A}$.
We will use the following notation: $f(a) = \langle a , \, f \rangle$
for the functions $f \in {\cal A}^*$ of elements $a \in {\cal A}$. For the image
$R_{12} = \langle  {\cal R}, \, T_1 \otimes T_2 \rangle$ of the universal matrix
${\cal R} \in {\cal A} \otimes {\cal A}$
we deduce $\forall a \in {\cal A}$ the identity
(by using (\ref{sweed}) and (\ref{2.15a}))
$$
R_{12} \, \langle a , \, T_1 \, T_2 \rangle =
R_{12} \, \langle a_{(1)}, \, T_1 \rangle \, \langle a_{(2)}, \, T_2 \rangle =
\langle \,  {\cal R} \, \Delta(a)  \, , \,  T_1 \otimes T_2 \, \rangle =
$$
$$
= \langle \,   \Delta^{\sf op}(a) \, {\cal R} \, , \, T_1 \otimes T_2 \, \rangle
= \langle \,   \Delta^{\sf op}(a)  \, , \, T_1 \otimes T_2 \, \rangle \, R_{12} =
 \langle  a, \, T_2 \, T_1 \rangle \, R_{12} \; .
$$
Since the element $a \in {\cal A}$ is not fixed here, one can conclude
(for the nondegenerate pairing) that the elements $T^{i}_{j}$  satisfy the
following quadratic relations ($RTT$ relations):
\be
\lb{3.1.1}
R^{i_{1}i_{2}}_{j_{1}j_{2}} \,
T^{j_{1}}_{k_{1}} \, T^{j_{2}}_{k_{2}} =
T^{i_{2}}_{j_{2}} \, T^{i_{1}}_{j_{1}} \,
R^{j_{1}j_{2}}_{k_{1}k_{2}}        \Leftrightarrow
R_{12}T_{1}T_{2} = T_{2}T_{1}R_{12} \Leftrightarrow
\R \, T_1 \, T_2 = T_1 \, T_2 \, \R \; ,
\ee
where the indices $1$ and $2$ label the matrix spaces
and the matrix $R_{12}$ satisfies Yang-Baxter equation (\ref{ybe}), (\ref{3.1.2}).

In the case of nontrivial $R$-matrices satisfying (\ref{ybe})
the relations (\ref{3.1.1})
define a noncommutative quadratic algebra (as the algebra of functions
with the generators $\{1, \, T^i_j \}$),
which is called $RTT$ algebra.
We stress that
one can consider $RTT$ algebra (\ref{3.1.1}) with the Yang-Baxter $R$-matrix
which is not in general the image of any universal ${\cal R}$ matrix.
The Yang-Baxter equation for $R$ is necessary to ensure that on monomials of third degree
in $T$ no relations additional to (\ref{3.1.1}) arise.
We shall assume that $R_{12}$ is a skew-invertible matrix.
In this case matrices $D$ and $Q$ (\ref{dmatr}) define 1-dimensional representations
$\rho_D(T^i_j) = D^i_j$ and $\rho_Q(T^i_j) = Q^i_j$ for RTT algebra (\ref{3.1.1}) (see (\ref{sk10})).
In some cases below we also
assume that $R_{12}$ is a lower triangular block matrix and
its elements depend on the numerical parameter $q=exp(h)$,
which is called the deformation parameter.

Suppose that the $RTT$ algebra can be extended in such a way
that it also contains all elements $(T^{-1})^{i}_{j}$:
\be
\lb{TT-1}
(T^{-1})^{i}_{k} \, T^{k}_{j} =
T^{i}_{k} \, (T^{-1})^{k}_{j} = \delta^{i}_{j} \cdot 1 \; .
\ee
Then this algebra becomes a Hopf algebra with structure mappings
\be
\Delta (T^{i}_{k}) = T^{i}_{j} \otimes T^{j}_{k} \ , \;\;
\epsilon(T^{i}_{j}) = \delta^{i}_{j} \  , \;\;
 S (T^{i}_{j}) = (T^{-1})^{i}_{j} \ ,
\lb{3.1.5}
\ee
which, as is readily verified, satisfy the standard axioms
(see sections {\bf \ref{hopf2}, \ref{hopf3}}):
\be
\lb{3.1.6}
\begin{array}{c}
(id \otimes \Delta) \Delta (T^{i}_{j}) =
(\Delta \otimes id) \Delta (T^{i}_{j}) \; ,  \cr
(\epsilon \otimes id) \Delta(T^{i}_{j}) =
(id \otimes \epsilon) \Delta(T^{i}_{j}) = T^{i}_{j} \; ,  \cr
m( S \otimes id) \Delta(T^{i}_{j}) =
m(id \otimes  S) \Delta(T^{i}_{j}) =\epsilon(T^{i}_{j}) 1 \; .
\end{array}
\ee
The antipode $S$ is not an involution, since instead of $S^{2} = id$ we have
equation
\be
\lb{3.1.7}
S^{2}(T^{i}_{j}) \, D^{j}_{l} = D^{i}_{k} \, T^{k}_{l}   \; ,
\ee
which can be rewritten in the form
\be
\lb{3.1.8}
D^{j}_{l} T^{l}_{k} S(T^{i}_{j}) = D^{i}_{k}  \; ,
\ee
and the matrix $D$ has been defined in (\ref{dmatr}).
The relations (\ref{3.1.7}) and (\ref{3.1.8}) can be interpreted as the rules of permutation
of the operations of taking the inverse matrix and the transposition ($t$):
\be
\lb{tmin}
D^{t} (T^{-1})^{t} = (T^{t})^{-1} D^{t}
\ee

To prove (\ref{3.1.7}) we note that the $RTT$ relations
(\ref{3.1.1}) can be represented in the form
$$
T_1^{-1} \, \R_{12} \, T_1 = T_2 \, \R_{12} \, T_2^{-1} \; .
$$
We multiply this relation by $\hat{\Psi}_{01}$ from the left and by $\hat{\Psi}_{23}$
from the right and take the traces $Tr_{12}(\dots)$.
Then, taking into account eqs.
(\ref{skew}), we arrive at the relations:
$$
Tr_1 \left( \hat{\Psi}_{01} \, T_1^{-1} \, P_{13} \, T_1 \right) =
Tr_2 \left( T_2 \, P_{02} \, T_2^{-1} \,  \hat{\Psi}_{23} \right) \; .
$$
Acting to this relations by $Tr_3(\dots)$ and $Tr_0(\dots)$
we obtain, respectively
\be
\lb{3.1.10a}
D_0 = Tr_2 T_2 \, P_{02} \, T_2^{-1} \, D_2 \; , \;\;\;
Tr_1 Q_1 \, T_1^{-1} \, P_{13} \, T_1 \, = Q_3  \; .
\ee
The first equation in (\ref{3.1.10a}) is identical to (\ref{3.1.7}) and (\ref{3.1.8}),
while the second one gives
\be
\lb{second}
 S(T^j_i) \, T^k_l \, Q^l_j =  Q^k_i \; .
\ee

As it was shown in subsection {\bf \ref{qtrace}},
the matrices $D^{i}_{j}$ and $Q^i_j$ (\ref{dmatr}), entering the conditions (\ref{3.1.10a}),
define {\it the quantum traces} \cite{10},\cite{18}.
To explain the features
of the quantum trace, we consider the $N^{2}$-dimensional adjoint
${\cal A}^*$-comodule $E$
(in what follows we continue
to use the concise notation ${\cal A}^*$ for the $RTT$ algebra).
We represent its basis elements in the form of
an $N \times N$ matrix $E=||E^{i}_{j}||, \; i,j=1,...,N$. The adjoint coaction is
\be
\lb{3.1.11}
E^{i}_{j} \rightarrow T^{i}_{i^{'}}
S(T^{j^{'}}_{j}) \otimes E^{i^{'}}_{j^{'}}
=: (TET^{-1})^{i}_{j} \; ,
\ee
where 
in the right hand
side of (\ref{3.1.11}) we have introduced
abbreviations that we shall use in what follows
(we omit the sign of the tensor product
and should only remember that the elements $E^i_j$ commute with the
elements $T^k_m$). We stress that
there is a different form of the adjoint coaction:
\be
\lb{3.1.12}
E^{i}_{j} \rightarrow
E^{i^{'}}_{j^{'}} \otimes
S(T_{i^{'}}^{i}) T^{j^{'}}_{j} =:
(T^{-1}ET)^{i}_{j} \; .
\ee

One can check that in (\ref{3.1.11}) and (\ref{3.1.12})
the elements $E^i_j$ form, respectively,
left and right comodules. The matrix $||T^i_j||$
is interpreted now as the matrix of linear noncommutative
adjoint transformations. Both left and right comodules $E$ are reducible,
and irreducible subspaces in $E$ can be distinguished by
means of the quantum traces. For the case (\ref{3.1.11}),
the quantum trace has the form (cf. (\ref{qtrDQ}))
\be
Tr_{D}E :=
Tr(DE) \equiv \sum_{i,j=1}^{N} D^{i}_{j} E^{j}_{i}
\label{3.1.13}
\ee
and satisfies the following invariance property,
which follows from Eqs. (\ref{3.1.7}), (\ref{3.1.8}) and first
relation in (\ref{3.1.10a}):
\be
\lb{3.1.14}
Tr_{D}(TET^{-1}) = Tr_{D}(E) \;,
\ee
For the case (\ref{3.1.12}), the definition of the quantum trace must be changed to
\be
\label{3.1.15}
Tr_{Q}E :=
Tr(Q \, E) \equiv \sum_{i,j=1}^{N} \, Q^{i}_{j} \, E^{j}_{i} \; , \;\;
Tr_{Q}(T^{-1}ET) = Tr_{Q}(E) \; ,
\ee
this follows from the second relation in (\ref{3.1.10a}).
Thus, $Tr_{D}(E)$ and $Tr_{Q}(E)$ are, respectively,
the scalar parts of the comodules
$E$ (\ref{3.1.11}) and (\ref{3.1.12}), whereas the $q$-traceless part of
$E$ generates $(N^{2}-1)$-dimensional (reducible in the general case and irreducible
in the case of linear quantum groups) ${\cal A}^*$-adjoint comodules.
Note that, if the matrix $D$ is invertible, one can substitute
$Q \rightarrow const \cdot D^{-1}$ in (\ref{3.1.15}) since eq. (\ref{3.1.8})
is rewritten in the form $(D^{-1})^k_i = S(T^j_i) \, T^k_l \, (D^{-1})^l_j$
(cf. (\ref{second})). We also note that formulas
(\ref{3.1.14}) and (\ref{3.1.15}) of the adjoint
invariancy of the quantum traces can be considered as universal
analogs of (\ref{RER}).

An important consequence of the definition of the quantum trace
(\ref{qtrDQ}),
(\ref{3.1.14}), (\ref{3.1.15}) and the $RTT$ relations (\ref{3.1.1}) is the fact that
\be
\lb{yyy4ab}
\begin{array}{c}
T^{-1}_{1} \, X(\R) \, T_{1} = T_2 \, X(\R) \, T_{2}^{-1}
 \;\;\;\; \Rightarrow \\ [0.2cm]
T^{-1}_{1} Tr_{D (2)}(X(\R)) T_{1} = Tr_{D (2)}(X(\R))
\, , \;\;\;\;
T_{2} \, Tr_{Q(1)}(X(\R)) \, T^{-1}_{2} = Tr_{Q(1)}(X(\R))\, ,
\end{array}
\ee
where $X(.)$ is an arbitrary function, while $Tr_{Q(1)}$
and $Tr_{D(2)}$ are
the quantum traces over the first and second space, respectively.
Eq. (\ref{yyy4ab}) indicates that the matrices
$Y_2 = Tr_{D(3)}(X(\R_2))= Tr_{Q(1)}(X(\R_1))$ (see (\ref{yyy1})) must be
proportional to the identity matrix if $T^i_j$ are
functions which define an irreducible
representation of the quantum group ${\cal A}$.
 In particular, we must have
\be
\lb{3.1.16}
Tr_{D(3)}(\R^{k}_{23}) =
Tr_{Q(1)}(\R^{k}_{12}) = c_{k} I_{2} \; ,
\ee
where $c_{k}$ are certain constants, e.g. $c_1 =1$ (\ref{qtrs}), (\ref{qtrs1}), and $I_{k}$
is the identity matrix in the $k$-th space.
Note that a direct consequence of (\ref{sk9}) is
\be
\lb{3.1.17}
Tr_{D(2)}(\R_{12}^{-1}) = c_{-1} \cdot I_{1} = D_1 \, Q_1 \; ,
\ee
and for $c_{-1} \neq 0$ matrices $D$, $Q$ are invertible.
As we will see below, for the quantum groups of the classical series,
the fact (\ref{3.1.16}) does indeed hold. In what follows, we shall
attempt to restrict consideration to either left or right adjoint comodules
with quantum traces (\ref{3.1.13}) or (\ref{3.1.15}). The analogous relations for right (or left) comodules
 can be considered exactly in the same way.

\vspace{0.2cm}

\subsubsection{Faddeev-Reshetikhin-Takhtajan $L^{\pm}$ algebras}

It can be seen from comparison of the relations (\ref{3.1.1}) and (\ref{ybe}), (\ref{3.1.2})
that for the generators $T^{i}_{j}$ it is possible to choose the
following finite-dimensional matrix representations:
\be
\lb{3.1.18}
(T^{i}_{j})^{k}_{l} = R^{ik}_{jl} \; , \;\;\;
(T^{i}_{j})^{k}_{l} = (R^{-1})^{ki}_{lj}  \; .
\ee
In these representations the images of invariance relations
 (\ref{3.1.14}), (\ref{3.1.15}) coincide with (\ref{RER}).
Since the $R$-matrix satisfies the Yang-Baxter equation, there exist linear functionals
$(L^{\pm})^{i}_{j} \in {\cal A}$ that realize the homomorphisms (\ref{3.1.18}), i.e., we have
\be
\lb{3.1.19}
\langle L^{+}_{2},T_{1}\rangle  = R_{12} := R_{12}^{(+)} \; , \;\;
\langle L^{-}_{2},T_{1}\rangle  = R^{-1}_{21} := R_{12}^{(-)} \; ,
\ee
{}For the case $R_{12}  = \langle{\cal R}, T_1 \otimes T_2 \rangle$
we immediately construct
the mapping from ${\cal A}^*$ to ${\cal A}$ (see, for example, \cite{18,18''})
\be
\lb{0.5}
\begin{array}{c}
\langle  {\cal R} \, , \, id  \otimes T^{i}_{j} \rangle  = (L^{+})^{i}_{j} \; , \;\;\;
\langle  {\cal R} \, , \, S(T^{i}_{j}) \otimes id \rangle  = (L^{-})^{i}_{j} \; , \\ \\
\langle  {\cal R} \, , \, T^{i}_{j} \otimes id \rangle  = S((L^{-})^{i}_{j}) \; .
\end{array}
\ee

Eqs. (\ref{3.1.19}) are generalized in the following form
$$
\langle  L_f^{\pm}, T_{1} \, T_{2} \dots T_{k} \rangle  =
R^{(\pm)}_{1f} \, R^{(\pm)}_{2f} \dots R^{(\pm)}_{kf} \; .
$$
The Yang-Baxter equation (\ref{3.1.2}) can now be reproduced from the $RTT$ relations
(\ref{3.1.1}) by averaging them with the $L^{\pm}$ operators.

From the requirement that elements
 $(L^{\pm})^{i}_{j} \in {\cal A}$ generate
  the algebra that is the dual to
the algebra ${\cal A}^*$ (the definition of the dual algebra is given in Def. \ref{def6}, Subsec.
{\bf \ref{hopf3}}),
we obtain the following commutation relations for the generators $L^{(\pm)}$:
\be
\lb{3.1.20b}
\R_{12} \, L^{\pm}_2 \, L^{\pm}_1 = L^{\pm}_2 \, L^{\pm}_1 \, \R_{12} \; , \;\;\;
\ee
\be
\lb{3.1.20}
\R_{12} \, L^{+}_2 \, L^{-}_1 = L^{-}_2 \, L^{+}_1 \, \R_{12} \; .
\ee
The same equations are obtained from the universal Yang-Baxter equation
(\ref{2.30}) by the averaging it with $(T_1 \otimes T_2 \otimes id)$,
$(id \otimes T_1 \otimes T_2 )$, $(T_1  \otimes id \otimes T_2)$
and using (\ref{0.5}).
The algebra (\ref{3.1.20b}), (\ref{3.1.20})
is obviously a Hopf algebra with comultiplication, antipode, and coidentity:
\be
\lb{3.1.20a}
\triangle (L^{\pm})^{i}_{j} =
(L^{\pm})^{i}_{k} \otimes (L^{\pm})^{k}_{j} \; , \;\;
S(L^{\pm}) = (L^{\pm})^{-1} \; , \;\;
\ee
\be
\lb{3.1.20c}
\epsilon((L^{\pm})^{i}_{j}) = \langle (L^{\pm})^{i}_{j},1\rangle =
\delta^{i}_{j} \; ,
\ee
where we have assumed that the matrices
$L^{\pm}$ are invertible.

We call the Hopf algebras with generators $\{ (L^{\pm})^{i}_{j} \}$, defining relations
(\ref{3.1.20b}), (\ref{3.1.20}) and structure mappings (\ref{3.1.20a}), (\ref{3.1.20c})
as {\it Faddeev-Reshetikhin-Takhtajan (FRT) algebras}.
 As was shown in Ref. \cite{10}, for the
$R$-matrices of the
quantum groups of the classical series
$A_{n}$, $B_{n}$, $C_{n}$, $D_{n}$ (respectively, $SL_{q}(n+1)$, $SO_{q}(2n+1)$, $Sp_{q}(2n)$,
$SO_{q}(2n)$),
the relations (\ref{3.1.20b}), (\ref{3.1.20}) define quantum universal enveloping Lie algebras in
which the elements $(L^{\pm})^{i}_{j}$ play the role of the quantum analog of the Cartan-Weyl
generators. We will investigate the case of $SL_{q}(n)$
 below in Subsec. {\bf \ref{qgrsl}}.

One can construct (see, e.g., \cite{Vlad}) the FRT algebra
(\ref{3.1.20b}) -- (\ref{3.1.20c}) as a Drinfeld
double of two dual Hopf subalgebras ${\cal B}^{+}$ and ${\cal B}^{-}$ with generators
$(L^{+})^i_j$ and $(L^{-})^i_j$, defining relations (\ref{3.1.20b}) and structure mappings
(\ref{3.1.20c}), and (cf. (\ref{3.1.20a}))
\be
\lb{3.1.20aa}
\begin{array}{c}
\triangle (L^{+})^{i}_{j} =
(L^{+})^{i}_{k} \otimes (L^{+})^{k}_{j} \; , \;\;\;
\triangle_{\rm op} \, (L^{-})^{i}_{j} =
(L^{-})^{k}_{j} \otimes (L^{-})^{i}_{k} \; , \\ \\
S(L^{+}) = (L^{+})^{-1} \; , \;\;
S^{-1}_{\rm op}  (L^{-}) = (L^{-})^{-1} \; , \;\;
\end{array}
\ee
In this case, the Hopf algebras ${\cal B}^{+}$ and ${\cal B}^{-}$
are dual to each other
with respect to the pairing \cite{Vlad}
\be
\lb{pairv}
\langle \! \langle L_1^{-} , L_2^{+} \rangle \! \rangle  =
 R^{-1}_{12} \; .
\ee
We denote by ${\cal B}^{-O}$ the Hopf algebra with
generators $(L^{-})^i_j$, and with comultiplication and antipode
(\ref{3.1.20a}) opposite to that of (\ref{3.1.20aa}). The algebras
${\cal B}^{+}$ and ${\cal B}^{-O}$ are antidual with respect to the pairing
(\ref{pairv}). As it was shown in Subsec. {\bf \ref{hopf4}},
 from the antidual Hopf algebras
${\cal B}^{+}$ and ${\cal B}^{-O}$ it is possible to construct a
Drinfeld quantum double ${\cal B}^{+} \Join {\cal B}^{-O}$, for which the cross commutation
relations have the form (\ref{3.1.20}). Thus, for the algebras ${\cal B}^{\pm}$ in (\ref{3.1.20b}),
one can propose a special cross-product
(quantum Drinfeld double),
given by (\ref{3.1.20}), which is again a Hopf algebra
(with structure mappings (\ref{3.1.20a}), (\ref{3.1.20c})), and which was used in Ref. \cite{10}
for $R$-matrix formulation of quantum deformations of the universal enveloping Lie algebras.

Note that the $FRT$ algebra (\ref{3.1.20b}), (\ref{3.1.20}) is a covariant algebra (comodule algebra)
with respect to
the left and right cotransformations
\be
\lb{3.1.21}
\begin{array}{c}
(L^{\pm})^{i}_{j} \rightarrow
(T^{-1})^{k}_{j} \otimes (L^{\pm})^{i}_{k}
\equiv (L^{\pm} \, T^{-1})^{i}_{j}
\; , \;\;       \\ \\
(L^{\pm})^{i}_{j} \rightarrow
(L^{\pm})^{k}_{j} \otimes (T^{-1})^{i}_{k}
\equiv (T^{-1} \, L^{\pm})^{i}_{j} \; ,
\end{array}
\ee
(we forget here for a moment
that the matrices $T$ and $L^{\pm}$ could have the different triangular properties).
Thus, the matrices
\be
\lb{3.1.21a}
L^{i}_{j} = (S(L^{-})L^{+})^{i}_{j} \; , \;\;\;
\bar{L}^{i}_{j} = (L^{+}S(L^{-}))^{i}_{j} \; ,
\ee
realize, respectively,
the left and right adjoint comodules (\ref{3.1.11}) and (\ref{3.1.12}).
It is readily verified that the coinvariants
\be
\lb{3.1.22}
p_{M}
= Tr_{D} \left( L^{M} \right) \; , \;\;\;
\overline{p}_M = Tr_{Q} \left( \bar{L}^{M} \right)
\ee
are central elements (see also \cite{10}) for the $FRT$ algebra (\ref{3.1.20b}),
(\ref{3.1.20}) (below we show that $p_M =\overline{p}_M$ for the realizations
(\ref{3.1.21a})).
Indeed, one can obtain from (\ref{3.1.20b}), (\ref{3.1.20}) the relations
\be
\label{iiii}
L_2^{M} \, L_1^{\pm} =
L_1^{\pm} \, \R^{\pm 1} \, L_1^{M} \, \R^{\mp 1} \; , \;\;
L_2^{\mp} \, \bar{L}_1^{M} =
\R^{\pm 1} \, \bar{L}_2^{M} \, \R^{\mp 1} \, L^{\mp}_2  \; ,
\ee
where $\R := \R_{12}$.
Then, by taking the traces $Tr_{D(2)}$ and $Tr_{Q(1)}$
respectively of the first and second relations (\ref{iiii}),
 and using (\ref{RER}) we prove $[p_M, L^{\pm}] =0=
 [\bar{p}_M, L^{\pm}] $ and therefore
demonstrate the centrality of the elements (\ref{3.1.22})
 for the algebra (\ref{3.1.20b}), (\ref{3.1.20}).

The equality $p_M =\overline{p}_M$ for the elements (\ref{3.1.22}) (where
$L$ and $\overline{L}$ are composed from $L^{\pm}$
(\ref{3.1.21a})) is proved as follows:
$$
Tr_{D} \left( L^{M} \right)
= Tr_{D(2)} \left( S(L^{-}_2) \, \bar{L}^{M}_2 \, L^{-}_2 \right)
= Tr_{Q(1)}Tr_{D(2)}
\left( S(L^{-}_2) \, \R \, \bar{L}_2^{M} \, L^{-}_2 \right) =
$$
$$
= Tr_{Q(1)}Tr_{D(2)}
\left( \bar{L}_1^{M} \, S(L^{-}_2) \, \R \, L^{-}_2   \right) =
Tr_{Q(1)}Tr_{D(2)}
\left( \bar{L}_1^{M} \, L^{-}_1 \, \R \, S(L^{-}_1)   \right) =
Tr_{Q} \left( \bar{L}^{M}  \right)  \; ,
$$
where we have used Eqs. (\ref{qtrs}), (\ref{qtrs1}), (\ref{3.1.20b}),
(\ref{3.1.20}), (\ref{iiii}).

\subsubsection{Reflection equation algebras}

Note also that the generators $L^{i}_{j}$ and $\overline{L}^{i}_{j}$ (\ref{3.1.21a}) satisfy equations
\be
\lb{3.1.23a}
\R_{12} \, L_1 \, \R_{12} \, L_1 = L_1 \, \R_{12} \, L_1 \, \R_{12}    \; ,
\ee
\be
\lb{3.1.23b}
\R_{12} \, \overline{L}_2 \, \R_{12} \, \overline{L}_2 = \overline{L}_2 \, \R_{12} \, \overline{L}_2 \, \R_{12}  \; ,
\ee
In Subsection {\bf \ref{factsc}} below,
 we will
see that (\ref{3.1.23a}) and (\ref{3.1.23b})
are the special limits of the reflection equations with
spectral parameters. In view of this, algebras with generators $L^i_j$ and
$\overline{L}^i_j$ and defining
 relations (\ref{3.1.23a}), (\ref{3.1.23b})
are called {\it left and right reflection equation algebras}, since (\ref{3.1.23a}) and (\ref{3.1.23b})
 are covariant under the left and right co-actions (co-transformations)
(\ref{3.1.11}), (\ref{3.1.12}). A set (which is incomplete in general; see below)
of central elements for these algebras are
represented by the same formulas as in (\ref{3.1.22}). Indeed, one can deduce
from (\ref{3.1.23a}), (\ref{3.1.23b}) the relations
\be
\lb{3.1.23d}
 L_1 \, \R_{12} \, L^M_1 \, \R_{12}^{-1} = \R_{12}^{-1} \, L^M_1 \, \R_{12} \, L_1    \; , \;\;
 \overline{L}_2 \, \R_{12} \, \overline{L}_2^M \, \R_{12}^{-1} =
 \R_{12}^{-1} \, \overline{L}^M_2 \, \R_{12} \, \overline{L}_2   \; .
 \ee
Then, taking the quantum traces $Tr_{D(2)}(\dots)$ and $Tr_{Q(1)}(\dots)$
of the first and second relations and using (\ref{RER}) we prove the centrality of the elements (\ref{3.1.22})
for the algebras (\ref{3.1.23a}), (\ref{3.1.23b})
\be
\lb{3.1.23i}
[L^i_j , \; Tr_{D} \left( L^{M} \right)] = 0  \; , \;\;\;\;\;\;
[\bar{L}^i_j , \;  Tr_{Q} \left( \bar{L}^{M} \right)] = 0 \; .
\ee

The algebra (\ref{3.1.23a}) (and similarly the second algebra (\ref{3.1.23b})) decomposes into the
direct sum of two subalgebras, namely, into the abelian algebra with generator $p_{1} := Tr_{_D}(L)$
and the algebra with $(N^{2} -1)$ traceless generators
$\tilde{L}^{i}_{j}$ (we assume that $Tr_{_D}(I) \neq 0$)
\be
\lb{3.1.24}
L^{i}_{j} =   p_{1}' \, \delta^{i}_{j} +
\lambda \, \tilde{L}^{i}_{j} \;\;\;\; \Rightarrow \;\;\;\;
\tilde{L}^{i}_{j} = \frac{1}{\lambda}
\left(  L^i_j - p_{1}' \, \delta^i_j \right) \; , \;\;\;\;\;
p_1' : = \frac{p_{1}}{Tr_{_D}(I)} \; ,
\ee
where the factor $\lambda := q-q^{-1}$ is introduced to ensure that the operators
$\tilde{L}$ have the correct classical limit
 for $q \rightarrow 1$.
For the latest algebra, it is easy to obtain the commutation relations
\be
\lb{3.1.25d}
\R_{12} \, \tilde{L}_1 \, \R_{12} \, \tilde{L}_1 -
\tilde{L}_1 \, \R_{12} \, \tilde{L}_1 \, \R_{12}  =
\frac{p_1'}{\lambda} ( \tilde{L}_1 \, \R_{12}^{2} -
\R_{12}^{2} \, \tilde{L}_1 ) \; ,
\ee
which after normalization
$\tilde{L}_1 \to - p_1' \, \tilde{L}_1$
(for $p_1' \neq 0$) gives
\be
\lb{3.1.25}
\R_{12} \, \tilde{L}_1 \, \R_{12} \, \tilde{L}_1 -
\tilde{L}_1 \, \R_{12} \, \tilde{L}_1 \, \R_{12}  =
\frac{1}{\lambda} (\R_{12}^{2} \, \tilde{L}_1 -
 \tilde{L}_1 \, \R_{12}^{2} ) \; .
\ee
These relations
can be regarded (for an arbitrary
Yang-Baxter $R$-matrix)
as a deformation of the commutation relations for Lie algebras. For the Hecke type $R$-matrix
(\ref{3.3.7}), the relations (\ref{3.1.25}) are equivalent to
\be
\lb{3.1.25a}
\R_{12} \, \tilde{L}_1 \, \R_{12} \, \tilde{L}_1 -
\tilde{L}_1 \, \R_{12} \, \tilde{L}_1 \, \R_{12}  =
\R_{12} \, \tilde{L}_1  - \tilde{L}_1 \, \R_{12}   \; ,
\ee
and corresponding algebra has a projector type
representation $\varrho$:
$(\tilde{L}^i_j)^\alpha_\beta = A^{i\alpha} B_{j\beta}$,
where numerical rectangular matrices $A$ and $B$ are such that
${\rm Tr}_{\varrho}(\tilde{L}^i_j)=
 B_{j\alpha} A^{i\alpha} = Q^i_j$ (for any matrix
$Q$ that satisfies $Tr_{1}Q_1 \R_{12} = I_2$;
see (\ref{qtrs})).

The relations (\ref{3.1.23a}), (\ref{3.1.23b}),
 (\ref{3.1.25}) and (\ref{3.1.25a}) are extremely important
and arise, for example, in the construction of a differential calculus on quantum
groups as the commutation relations for invariant vector fields
(see \cite{Wor}
 -- \cite{18''}
and refs. therein; see also Section {\bf \ref{difcal}} below).

Note that, instead of (\ref{0.5}), one can use a somewhat different linear
mapping from ${\cal A}^{*}$ to ${\cal A}$ \cite{13',18,18'',18''',18b}
(which is completely determined by (\ref{0.5})):
\be
\lb{0.5a}
\langle  \sigma({\cal R}) \, {\cal R} \, , \, id \otimes a \rangle  =
 \alpha \;\;\; (a \in {\cal A}^{*} \; , \alpha \in {\cal A}) \; ,
\ee
where $\sigma (a \otimes b) = (b \otimes a)$, $\forall a,b \in {\cal A}$. The explicit calculations
give
\be
\lb{RRTL}
\langle  \sigma({\cal R}) \, {\cal R} \, , \, id \otimes T^i_j \rangle  = L^i_j \; ,
\ee
\be
\lb{0.5aaa}
\begin{array}{c}
\langle  \sigma({\cal R}) \, {\cal R} \, , \, id \otimes T_1 \, T_2 \rangle  =
S(L^{-}_1) \, L_2 \, L^+_1 =
L_1 \, \R_1 \, L_1 \, \R_1^{-1} \; , \\ \\
\!\!\!\!\!\!\!\!\!
\langle  \sigma({\cal R}) \, {\cal R} \, , \, id \otimes T_1 \, T_2 \, T_3 \rangle  =
S(L^{-}_1) \, S(L^{-}_2) \, L_3 \, L^+_2 \, L^+_1 =
L_{\underline{1}} \, L_{\underline{2}} \, L_{\underline{3}}
\equiv L_{\overline{3}} \, L_{\overline{2}} \, L_{\overline{1}}
\; , \\
. \; . \; . \; . \; . \; . \; . \; . \; . \; . \; . \; . \; , \\
\langle  \sigma({\cal R}) \, {\cal R} \, , \, id \otimes T_1 \dots \, T_k \rangle  =
L_{\underline{1}} \, L_{\underline{2}} \dots L_{\underline{k}} \equiv
L_{\overline{k}} \dots L_{\overline{2}} \, L_{\overline{1}} \; ,
 \end{array}
\ee
where
\be
\lb{0.5ab}
L_{\underline{k+1}} = \R_{k} \, L_{\underline{k}}
\, \R_{k}^{-1} \; ,  \;\;\;\; L_{\overline{k+1}} = \R_{k}^{-1} \, L_{\overline{k}} \, \R_{k} \; , \;\;\;\;\;
L_{\underline{1}} = L_{\overline{1}}  = L_1 \; ,
\ee
and we have used eqs. (\ref{2.27}), (\ref{0.5}) and (\ref{iiii}).
If we confine ourselves to the fairly general case of quasitriangular
Hopf algebras ${\cal A}$, for which the mapping
 (\ref{0.5a}) is invertible (such Hopf algebras
 are called factorizable \cite{18'''}), one can map
the identities for the $RTT$ algebra into the identities for the reflection equation algebra and
vice versa. For this we need to use relations (\ref{0.5a}) (for more details see
\cite{18b}).

In view of (\ref{0.5aaa}) one can represent the reflection equation algebra (\ref{3.1.23a})
in the ``universal'' form
$$
{\cal R}_{32} \, ({\cal R}_{31} {\cal R}_{13}) \, {\cal R}_{23} \,
({\cal R}_{21} \, {\cal R}_{12}) = ({\cal R}_{21} \, {\cal R}_{12}) \,
{\cal R}_{32} \, ({\cal R}_{31} {\cal R}_{13}) \, {\cal R}_{23}  \; ,
$$
where the notation ${\cal R}_{ij}$ has been introduced in (\ref{2.30a}). The pairing of this
relation with $(id \otimes T \otimes T)$ gives (\ref{3.1.23a}). The algebra (\ref{3.1.23b}) has an
analogous representation if we start with
\be
\lb{RRTbL}
\langle  \sigma({\cal R}) \, {\cal R} \, , \, T^i_j \otimes id  \rangle  = \overline{L}^i_j \; .
\ee

We note that the identity (which has been obtained in (\ref{0.5aaa}))
\be
\lb{L123}
L_{\underline{1}} \, L_{\underline{2}} \dots L_{\underline{k}} =
L_{\overline{k}} \dots L_{\overline{2}} \, L_{\overline{1}} \; ,
\ee
valid in more general case of any reflection equation algebra (\ref{3.1.23a})
(even not realized in the form (\ref{3.1.21a})). Below we also use the
following identity (which can be proved by induction)
\be
\lb{L124}
L_{\underline{k+1}} \, L_{\underline{k+2}} \dots L_{\underline{k+n}} =
U_{(k,n)} \, L_{\underline{1}} \, L_{\underline{2}} \dots L_{\underline{n}} \; U_{(k,n)}^{-1} \; .
\ee
where the operator $U_{(k,n)}$ is represented as a
product of $k$ or $n$ factors (cf. (\ref{TnR3})):
\be
\lb{04.4a}
U_{(k,n)} = \R_{(k \to n+k-1)} \dots \R_{(2 \to n+1)} \, \R_{(1 \to n)} \equiv
\R_{(k \leftarrow 1)} \, \R_{(k+1 \leftarrow 2)} \dots
\R_{(n+k-1 \leftarrow n)} \; ,
\ee
\be
\lb{04.4b}
\R_{(k \leftarrow m)} := \hat{R}_{k} \, \hat{R}_{k-1}
\cdots \hat{R}_{m} \; , \;\;\;
\R_{(m \to k)} := \hat{R}_{m} \, \hat{R}_{m+1}
\cdots \hat{R}_{k} \; .
\ee


\subsubsection{Central and commuting subalgebras
for reflection equation and $RTT$ algebras\label{cencom}}

As we prove in previous Subsection the elements
(\ref{3.1.22}) are central for the $RLRL$
(reflection equation) algebras (\ref{3.1.23a}), (\ref{3.1.23b}).
Now the description of a more general set of central elements for  reflection equation algebra is
in order.


\begin{proposition}\label{prop5}
{\it Let $X_{(1 \to m)}$ be an arbitrary element of the group algebra of the
braid group ${\cal B}_{m}$ generated by skew-invertible $R$-matrices $\hat{R}_a$ $(a=1,\dots,m-1)$
with defining  relations (\ref{3.1.3}), (\ref{3.1.3i}), (\ref{3.1.3is}). Then, the elements
\be
\lb{center}
z_m(X) = Tr_{{\cal D}(1 \dots m)} \left(X_{(1 \to m)} \, L_{\underline{1}} \, L_{\underline{2}}
\dots L_{\underline{m}}  \right) \; \;\;\; (m=1,2,\dots) \; ,
\ee
belong to the center $Z(L)$ of the reflection equation algebra (\ref{3.1.23a}), where we recall $\hat{R}_{12} \equiv \hat{R}_1$.}
\end{proposition}

\noindent
{\bf Proof.} First of all we note that $z_m(X)$ (\ref{center}) satisfies
\be
\lb{cent1}
\begin{array}{c}
z_m(X) \, I_1 = Tr_{{\cal D}(2 \dots m+1)} \left(X_{(2 \to m+1)} \,
L_{\underline{2}} \, L_{\underline{3}} \dots L_{\underline{m+1}}  \right) = \\
= Tr_{{\cal D}(2 \dots m+1)} \left(X_{(2 \to m+1)} \,
L_{\overline{m+1}} \dots  L_{\overline{3}} \, L_{\overline{2}}  \right) \; ,
\end{array}
\ee
where $X_{(2 \to m+1)} \in {\cal B}_{m+1}$
 is obtained from $X_{(1 \to m)}$ by the shift $R_a \to
R_{a+1}$ $(\forall a)$. The first equality follows from the chain of relations
$$
\begin{array}{c}
Tr_{{\cal D}(2 \dots m+1)} \left(X_{(2 \to m+1)} \,
L_{\underline{2}} \dots L_{\underline{m+1}}  \right) = \\
= Tr_{{\cal D}(2 \dots m+1)} \left(X_{(2 \to m+1)} \, \R_1 \cdots \R_m \,
L_{\underline{1}} \dots L_{\underline{m}} \, \R_m^{-1} \cdots \R_1^{-1} \right)  = \\
= Tr_{{\cal D}(2 \dots m+1)} \left(  \R_1 \cdots \R_m \, (X_{(1 \to m)} \,
L_{\underline{1}} \dots L_{\underline{m}}) \, \R_m^{-1} \cdots \R_1^{-1} \right)  = \\
= Tr_{{\cal D}(2 \dots m)} \left(  \R_1 \cdots \R_{m-1} \,  [Tr_{{\cal D}(m)} (X_{(1 \to m)} \,
L_{\underline{1}} \dots L_{\underline{m}})] \, \R_{m-1}^{-1} \cdots \R_1^{-1} \right) = \\
= \dots =  I_1 \, Tr_{{\cal D}(1 \dots m)} \left(X_{(1 \to m)} \,
L_{\underline{1}} \, L_{\underline{2}} \dots L_{\underline{m}}  \right) \; ,
\end{array}
$$
where we have applied (\ref{RER}) many times. The second equality in
(\ref{cent1}) is proved in the same way, or by using the generalization
of the identity (\ref{L123})
$$
L_{\underline{m}} \, L_{\underline{m+1}} \dots L_{\underline{k}} =
L_{\overline{k}} \dots L_{\overline{m+1}} \, L_{\overline{m}} \;\;\; (m < k) \; .
$$
Then, the proof of the commutativity of the arbitrary generator of the reflection equation (RE) algebra (\ref{3.1.23a})
with elements $z_m(X)$ is straightforward
$$
L_1 \, z_m(X) = Tr_{{\cal D}(2 \dots m+1)} \left(X_{(2 \to m+1)} \, L_{\underline{1}}
L_{\underline{2}} \, L_{\underline{3}} \dots L_{\underline{m+1}}  \right) =
$$
$$
 = Tr_{{\cal D}(2 \dots m+1)} \left(X_{(2 \to m+1)} \,
L_{\overline{m+1}} \dots L_{\overline{2}} \, L_{\overline{1}} \right) = z_m(X) \, L_1 \; .
$$
\hfill ~~\rule{2.5mm}{2.5mm}\par

{\bf Remark 1.} If, in the definition of central generators (\ref{center}), we take
the set of elements $X=X_\alpha$, $\alpha =1,2,\dots$,
which are all
primitive idempotents for any finite dimensional
 quotient ${\cal B}_m'$ of the group algebra of ${\cal B}_m$,
then, the set of central elements $z_m(X_\alpha)$ form a basis in the subspace of $Z(L)$ generated by elements
(\ref{center}) for any matrices $X \in {\cal B}_m'$.

{\bf Remark 2.} The ``power sums'' (\ref{3.1.22}) belong to the space $Z(L)$. Indeed, the
substitution of
$X = \R_{(m-1 \leftarrow 1)}:= \R_{m-1} \dots \R_1$  in (\ref{center}) gives
$$
z_m(X) = Tr_{{\cal D}(1 \dots m)} \left( L_{\underline{1}} \dots L_{\underline{m-1}} \, (\R_{(m-1
\leftarrow 1)} L_{1}
\R_{(m-1 \leftarrow 1)}^{-1}) \, \R_{(m-1 \leftarrow 1)} \right) =
$$
\be
\lb{sigms}
= Tr_{{\cal D}(1 \dots m)} \left( L_{\underline{1}} \dots L_{\underline{m-2}} \,
(\R_{(m-2 \leftarrow 1)} L_{1}
\R_{(m-2 \leftarrow 1)}^{-1})
 \R_{m-1} \R_{(m-2 \leftarrow 1)} L_{1}
 \right) =
\ee
$$
= Tr_{{\cal D}(1 \dots m-1)} \left( L_{\underline{1}} \dots L_{\underline{m-2}} \,
 \R_{(m-2 \leftarrow 1)} L_{1}^2 \right) = \dots =  Tr_{D(1)} (L_1^m) = p_m \; ,
$$
where in the first line we used the cyclic property
of the quantum trace (\ref{cycl}) and in
the second line we applied (\ref{qtrs}).

\vspace{0.2cm}

Now we discuss the set of commuting elements in
the $RTT$ algebra (\ref{3.1.1}). For this algebra
one can construct \cite{18'} the following
elements:
\be
\lb{0.4}
Q_{k} =  Tr_{Y(1 \dots k)} (\R_{(k-1 \leftarrow 1)}
\, T_{1} \, T_{2} \cdots T_{k} )
= Tr_{Y(1 \dots k)} (\R_{(1 \to k-1)} \,
 T_{1} \, T_{2} \cdots T_{k} )   \; ,
\ee
where
$$
Tr_{Y(1 \dots k)}(X_{1...k})  :=
Tr_{1} \dots Tr_{k} (Y_1 \dots Y_k \, X_{1...k}) \; ,
$$
and the matrices $Y$ are such that $Y_1 Y_2 \R_1 = \R_1 Y_1 Y_2$
(e.g., $Y = D$ or $Y = Q$, see (\ref{sk10})).
The second equality in (\ref{0.4}) is obtained as follows:
\be
\lb{04.4}
\begin{array}{c}
Tr_{Y(1 \dots k)} (\hat{R}_{1} \cdots \hat{R}_{k-1} \,
T_{1} \cdots T_{k} )  =
 Tr_{Y(1 \dots k)} (\hat{R}_{k-1} \, T_{1} \cdots T_{k} \,
\hat{R}_{1} \cdots \hat{R}_{k-2}) =  \\
= Tr_{Y(1 \dots k)} (\hat{R}_{1} \cdots \hat{R}_{k-3} \,
\hat{R}_{k-1} \, \hat{R}_{k-2} \, T_{1}  \cdots T_{k} \, ) = \\
= Tr_{Y(1 \dots k)} (
\hat{R}_{k-1} \, \hat{R}_{k-2} \, T_{1}  \cdots T_{k} \, \hat{R}_{1} \cdots \hat{R}_{k-3} )
= \\ =  \dots
= Tr_{Y(1 \dots k)} (\hat{R}_{k-1} \, \hat{R}_{k-2}
\cdots \hat{R}_{1} \, T_{1}  \cdots T_{k} ) \; .
\end{array}
\ee
Note that by means of (\ref{0.5a}) we map the elements $Q_k$ (\ref{0.4}) (for $Y = D$) to the
central elements $p_k$ (\ref{3.1.22}) of the reflection equation algebra.
\begin{proposition}\label{comset}
The elements (\ref{0.4}) generate
a commutative subalgebra in the $RTT$ algebra (\ref{3.1.1}).
\end{proposition}
{\bf Proof.} Our proof of the commutativity
of the elements $Q_{k}$ is based (see \cite{18c})
on the fact that there exists the operator $U_{(k,n)}$ (\ref{04.4a}) which satisfies
$$
U_{(k,n)} \, \R_{i} \, U_{(k,n)}^{-1} = \hat{R}_{i+k} \; , \;\;
i = 1, ... , n-1 \, , \;\;\;\;\;\;
U_{(k,n)} \, \R_{n+j} \, U_{(k,n)}^{-1} = \hat{R}_{j} \, , \;\;
j = 1, ... , k-1 \; .
$$
Using the operator $U_{(k,n)}$ we obtain the commutativity of $Q_k$:
\be
\lb{qq}
\begin{array}{c}
Q_k \, Q_n = Tr_{Y(1 \dots k)} (\hat{R}_{(1 \to k-1)} \,
T_{1} \cdots T_{k} ) \, Tr_{Y(1 \dots n)} (\hat{R}_{(1 \to n-1)} \,
T_{1} \cdots T_{n} ) = \\ [0.2cm]
= Tr_{Y(1 \dots k+n)} (\hat{R}_{(1 \to k-1)} \, \hat{R}_{(k+1 \to k+n-1)} \,
T_{1} \cdots  T_{k+n} ) = \\ [0.2cm]
= Tr_{Y(1 \dots k+n)} ( U_{(k,n)} \, \hat{R}_{(n+1 \to n+k-1)} \, \hat{R}_{(1 \to n-1)} \,
U_{(k,n)}^{-1} \, T_{1} \cdots  T_{k+n} ) = \\ [0.2cm]
= Tr_{Y(1 \dots k+n)} (\hat{R}_{(1 \to n-1)} \, \hat{R}_{(n+1 \to n+k-1)} \,
U_{(k,n)}^{-1} \, T_{1} \cdots  T_{k+n} \,  U_{(k,n)} ) = Q_n \, Q_k \; .
\end{array}
\ee
\hfill \qed

In fact, applying the same method as in (\ref{qq}), one can prove \cite{18c} that the set of
commuting elements in the $RTT$ algebra is wider then the set (\ref{0.4}) and consist of all
elements of the form
\be
\lb{charac}
Q_{k}(X) =
Tr_{Y(1 \dots k)} \left( X(\R_1, \dots \R_{k-1})  \, T_{1} \, T_{2} \cdots T_{k} \right)
\ee
where $X(\dots)$ run over basis elements of the braid group algebra
with generators $\{ \R_i \}$
$(i=1, \dots, k-1)$.

Our conjecture is that, for the Hecke
type $R$-matrices (\ref{3.3.7}), the set of elements (\ref{0.4})
$$
Q_k := Q_k(\R_{1 \to k-1}) \equiv Q_k(\R_{k-1 \leftarrow 1}) \; ,
$$
 is complete and all $Q_k(X)$ (for any braid $X$ with $k$ strands)
 are expressed as polynomials of the
commuting variables $\{ Q_1, ... , Q_k \}$ and deformation parameter $q$. These polynomials,
 if we add some additional constraints dictated by
 Markov (Reidemeister) moves
for the braids $X$ (see Sect. 1 in \cite{Jon1}), could be related
 to link polynomials. On the
other hand, eq. (\ref{charac}) defines $q$-analogs of characters for representations of the algebra
${\cal A}$
(\ref{3.1.20b}) -- (\ref{3.1.20c}) and for $RTT$ algebra
${\cal A}^*$. These representations
 are characterized by special choices of the elements $X(\dots)$ being
central idempotents in the Hecke algebra generated by
matrices $\{ \R_i \}$ $(i=1, \dots, k-1)$. We will discuss these ideas in detail
in Sect. {\bf \ref{qdlink}} below.


\subsubsection[Heisenberg double (HD) for the $RTT$
and reflection equation (RE) algebras]{Heisenberg 
double for the $RTT$
and reflection equation algebras}

Since the $RTT$ algebra ${\cal A}^*$ (\ref{3.1.1}), (\ref{3.1.5}) and the
quantum algebra ${\cal A}$ (\ref{3.1.20b}), (\ref{3.1.20}),
(\ref{3.1.20a}), (\ref{3.1.20c}) are
Hopf dual to each other
(with respect to the pairing (\ref{3.1.19})) one can define the
left and right Heisenberg doubles (HD) of these
algebras (about HD see Subsec. {\bf \ref{hopf4}}). Their cross-multiplication rules  (\ref{25d}), (\ref{25dd}) are
written for the left HD  in the form
\be
\lb{lhd}
L^{+}_1 \, T_2 = T_2 R_{21} L^{+}_1 \; , \;\;\;
L^{-}_1 \, T_2 = T_2 R^{-1}_{12} L^{-}_1 \; ,
\ee
and for the right one we have
\be
\lb{rhd}
T_1 L^{+}_2 =  L^{+}_2 R_{12} T_1 \; , \;\;\;
T_1 L^{-}_2 =  L^{-}_2 R^{-1}_{21} T_1 \; .
\ee
The corresponding cross products of the $RTT$ algebra  and
the  reflection equation
algebras (\ref{3.1.21a}), (\ref{3.1.23a}), (\ref{3.1.23b})
are described by the cross-multiplication rules
\be
\lb{tangb}
\overline{L}_1 \, T_2 = T_2 \R_{12} \overline{L}_2 \R_{12} \; , \;\;\;
T_1 L_2 = \R_{12} L_1 \R_{12} T_1
\ee
in the case of the left (\ref{lhd}) and right
(\ref{rhd}) HD, respectively. A remarkable property
\cite{alfad}
of these cross products is the existence of automorphisms of the HD algebras
\be
\lb{discr}
\{ T \, , \, \overline{L} \} \; \stackrel{\overline{m}_n}{\rightarrow} \;
\{ T  \, \overline{L}^n \, , \, \overline{L} \} \; , \;\;\;
\{ T \, , \, L \} \;  \stackrel{m_n}{\rightarrow} \; \{ L^n \, T \, , \, L \} \; ,
\ee
i.e. we have (the same is valid
for the automorphisms $\overline{m}_n$)
 \be
 \lb{discr2}
\R_{12} \, (L^n \, T)_1 \, (L^n \, T)_2  =
(L^n \, T)_1 \, (L^n \, T)_2 \, \R_{12} \; , \;\;\;\;\;\;
(L^n \, T)_1 \, L_2 =
\R_{12} \,L_1 \,\R_{12} \, (L^n \, T)_1 \;\; \Rightarrow
 \ee
 $$
(L^n \, T)_1 \, L_2^k =
(\R_{12} \,L_1 \,\R_{12})^k \, (L^n \, T)_1 \; ,
\;\;\;\;\;\; \forall n,k \in \mathbb{Z}_{\geq 0} \; .
$$
One can check these properties by induction using equations (\ref{3.1.1}), (\ref{3.1.23a}), (\ref{3.1.23b}) and
(\ref{tangb}). The maps $m_n, \overline{m}_n$ define discrete
time evolutions on the $RTT$ algebra. For the Hecke type $R$-matrices
(\ref{3.3.7}) the automorphisms (\ref{discr}) can be generalized
in the form
\be
\lb{discr1}
\{ T \, , \, \overline{L} \} \; \stackrel{\overline{m}'_n}{\rightarrow} \;
\{ T  \, (\sum_{m=0}^n \bar{x}_m \overline{L}^m)\, , \, \overline{L} \} \; , \;\;\;
\{ T \, , \, L \} \;  \stackrel{m'_n}{\rightarrow} \; \{ (\sum_{m=0}^n x_m L^m)
\, T \, , \, L \} \; ,
\ee
for any parameters $x_m,\bar{x}_m \in \mathbb{C}$. This generalization
follows from the fact that any symmetric function of two variables
$L_1$ and $\R_1 L_1 \R_1$ commutes with $\R_1$.

For the left  and right Heisenberg doubles (\ref{lhd}), (\ref{rhd}),
(\ref{tangb}) one can define new reflection equation algebras, generated by the elements of matrices $L$ and $\overline{L}$
transformed by the adjoint action of the
$RTT$ algebra
$$
\overline{Y} = T \, \overline{L}^{-1} \, T^{-1} \; , \;\;\; Y = T^{-1} \, L^{-1} \, T \; ,
$$
for which we have \cite{Isaev1}
(cf. (\ref{3.1.23a}), (\ref{3.1.23b}), (\ref{tangb})):
$$
\R_{12} \, \overline{Y}_1 \, \R_{12} \, \overline{Y}_1 =
\overline{Y}_1 \, \R_{12} \, \overline{Y}_1 \, \R_{12} \; , \;\;\;
\R_{12} \, Y_2 \, \R_{12} \, Y_2 =
Y_2 \, \R_{12} \, Y_2 \, \R_{12} \; ,
$$
$$
T_1 \, \overline{Y}_2 = \R_{12} \, \overline{Y}_1 \, \R_{12} \, T_1 \; , \;\;\;
Y_1 \, T_2 = T_2 \, \R_{12} \, Y_2 \, \R_{12} \; .
$$
The elements of these matrices satisfy:
$[\overline{Y}_2 , \, \overline{L}_1] = 0 = [Y_1 , \, L_2]$.
In the differential calculus on quantum groups
matrices $L$ and $M := Y^{-1}$ are interpreted
(see \cite{Isaev1} and \cite{IsPyaE})
as invariant vector fields
on the $RTT$ algebras (see Proposition
{\bf \ref{propis1}} below).

The cross-multiplication rules (\ref{tangb}) for the HD
of the $RTT$ and reflection equation algebras
were extensively exploited in the context of the $R$-matrix approach
to the differential calculus on quantum groups
\cite{AsCast} -- \cite{18''}
 (see also Section {\bf \ref{difcal}}
 below). Another cross-multiplications
(of the $RTT$ and reflection equation matrix algebras) which characterized by the relations
\be
\lb{crtL}
\overline{L}_1 \, T_2 = T_2 \R_{12} \overline{L}_2 \R_{12}^{-1}
 \; , \;\;\;\;\;\;
T_1 L_2 = \R_{12} L_1 \R_{12}^{-1} T_1  \; ,
\ee
were also considered in various investigations
\cite{Isaev1}, \cite{22}, \cite{Isaev2} of a
noncommutative differential geometry on quantum groups.
\begin{proposition}\label{propis1}
1.) For cross-multiplication
of the $RTT$ and reflection equation algebras (REA)
with generators $T^i_j$ and $L^i_j$ subject defining
relations (see (\ref{3.1.1}), (\ref{3.1.23a})
and (\ref{crtL}))
\be
\lb{crtL2}
\R_{12} \, T_1 T_2  =
T_1 T_2 \, \R_{12} \; , \;\;\;
\R_{12} L_1 \, \R_{12} L_1 = L_1 \, \R_{12} L_1 \, \R_{12}
 \; , \;\;\;
T_1 L_2 = \R_{12} L_1 \R_{12}^{-1} T_1 \; ,
\ee
we have the following equations\cite{Isaev2}
\be
\lb{crtL2a}
\R_{12} \, (L \, T)_1 \, (L \, T)_2  =
(L \, T)_1 \, (L \, T)_2 \, \R_{12} \; ,
\ee
(we however stress that it is impossible to define
 the whole discrete evolution (\ref{discr2})
for the double algebra (\ref{crtL2})). \\
2.) Let $L^i_j$,
$\tilde{L}^i_j$ be generators of the
REA (\ref{3.1.23a}) and $\tilde{L}^i_j$ subject the
following cross-commutation relations \cite{Isaev2}
with generators of (\ref{crtL2})
\be
\lb{crtL3}
 T_1 \tilde{L}_2 = \R_{12} \tilde{L}_1 \R_{12}^{-1} T_1
 \; , \;\;\; \R_{12}^{-1} \tilde{L}_1 \, \R_{12} L_1 =
L_1 \, \R_{12}^{-1} \tilde{L}_1 \, \R_{12} \; ,
\ee
Then we have \cite{Isaev2}, \cite{14}
\be
\lb{crtL4}
\begin{array}{c}
\R_{12} \, (\tilde{L} \, T)_1 \, L_2  =
\R_{12} L_1 \, \R_{12}^{-1} (\tilde{L} \, T)_1
\, , \;\;\; \;
\R_{12} \, (\tilde{L} \, T)_1 \, (\tilde{L} \, T)_2  =
(\tilde{L} \, T)_1 \, (\tilde{L} \, T)_2 \, \R_{12} \; ,
\\ [0.2cm]
\R_{12} (L \, \tilde{L})_1 \, \R_{12} (L \, \tilde{L})_1 =
(L \, \tilde{L})_1 \, \R_{12} (L \, \tilde{L})_1 \, \R_{12}
 \; .
 \end{array}
\ee
\end{proposition}

\subsubsection{Quantum matrix algebras in general setting}


\newtheorem{def9}[def1]{Definition}
\begin{def9} \label{def9}
{\it Let ${\cal A}(\Delta,\epsilon,S)$ be a Hopf algebra and functions
$T^i_j \in {\cal A}^*$
($i,j=1,\dots,N$) define $N$- dimensional representation of ${\cal A}$.
Consider an element $m \in {\cal A} \otimes {\cal A}$. The unital subalgebra of
${\cal A}$ with $N^2$ generators: $M^i_j := \langle m , \, id \otimes T^i_j \rangle
\in {\cal A}$
(or  $M^i_j :=  \langle m , T^i_j \otimes id \rangle$) is called quantum matrix algebra
${\cal M}(m) \subset {\cal A}$.

A unital algebra ${\cal M}$ is called right (or left) ${\cal A}$-coideal algebra
(${\cal A}$-covariant algebra) if it is ${\cal
A}$-comodule algebra with respect to any coactions (homomorphisms of ${\cal M}$):
${\cal M} \rightarrow {\cal M} \otimes {\cal A}$
(or ${\cal M} \rightarrow {\cal A} \otimes {\cal M}$).}
\end{def9}
For example, consider a quasi-triangular Hopf algebra ${\cal A}(\Delta,\epsilon,S,{\cal R})$ and an
element $m \in {\cal A} \otimes {\cal A}$ such that $(\Delta \otimes id)m = m_{13}m_{23}$
(cf. (\ref{2.27})). Then elements  $M^i_j := \langle m , \, T^i_j \otimes id \rangle$ generate
the matrix algebra ${\cal M}(m)$ which is nothing but
$RTT$ algebra (\ref{3.1.1}). Reflection equation algebras (\ref{3.1.23a}), (\ref{3.1.23b}) are
also quantum matrix algebras in view of the relations (\ref{RRTL}),
(\ref{RRTbL}). The reflection
equation algebras are $RTT$-coideal algebras since they are covariant under coactions
(\ref{3.1.11}), (\ref{3.1.12}).

Now we present a definition of a more general quantum matrix
algebra
${\cal M}(\R, \hat{F})$
generated by $(N \times N)$ matrix components $M^i_j$ subject to the relation
\be
\lb{qmam}
\R_{12} \, M_1 \, \hat{F}_{12} \, M_1 \, \hat{F}_{12} =
M_1 \, \hat{F}_{12} \, M_1 \, \hat{F}_{12} \, \R_{12} \; ,
\ee
where the pair of the Yang-Baxter operators $\{ \R , \hat{F} \} \in {\rm End}(V_N^{\otimes 2})$
satisfy the conditions
\be
\lb{compat}
\R_{12} \, \hat{F}_{23} \, \hat{F}_{12}  = \hat{F}_{23} \, \hat{F}_{12} \, \R_{23} \; , \;\;\;
\R_{23} \, \hat{F}_{12} \, \hat{F}_{23}  = \hat{F}_{12} \, \hat{F}_{23} \, \R_{12} \; .
\ee
The algebra ${\cal M}(\R, \hat{F})$ is a quantum matrix algebra ${\cal M}( {\cal \sigma(R) F})$
 since we can reproduce (\ref{qmam}) (for details see  \cite{18ddd}) by means of identifications
$$
M^i_j :=  \langle {\cal \sigma(R) F} , id  \otimes T^i_j\rangle \; , \;\;\;
\hat{F}_{12} := P_{12} \langle {\cal F}, \, T_1 \otimes T_2 \rangle \; ,
$$
where ${\cal F}$ is a twisting matrix (\ref{cocycl}), (\ref{4}) and
$P_{12}$ is the permutation matrix (\ref{perM}). Note that eqs.
(\ref{compat}) are the images of equations (\ref{5aa}). It means that,
for the pair of the Yang-Baxter operators $\{ \R , \hat{F} \}$ (\ref{compat}), the matrix
\be
\lb{twistP}
\hat{R}^{F}_{21} = \hat{F}_{12} \hat{R}_{12} \hat{F}^{-1}_{12} =
\langle {\cal R}^{\cal F} , T_1 \otimes T_2 \rangle P_{12} \; ,
\ee
is the Yang-Baxter matrix as well. Specializing to $\hat{F} = P$ or $\hat{F} = \R$ one reproduces
from (\ref{qmam}) the $RTT$ or reflection equation algebras, respectively. The algebras
${\cal M}(\R, \hat{F})$
(\ref{qmam}) and their modifications were discussed in \cite{18c}, \cite{18ddd}, \cite{18d}, \cite{GPS05}.

At the end of this subsection we introduce a notion
of a coideal subalgebra of
the quantum algebra (\ref{3.1.20b}), (\ref{3.1.20}).
Let $R_{12}$ be a Yang-Baxter $R$-matrix and
there are numerical matrices $G^i_j$,
$\overline{G}^i_j$ which satisfy the conditions
\be
\lb{coid}
\begin{array}{c}
R_{12} \, G_2 \, R_{12}^{t_2} \, G_1 =
G_1 \, R_{21}^{t_1} \, G_2 \, R_{21}^{t_1 t_2} \; , \\ [0.3cm]
R_{12} \, \overline{G}_1 \, R_{12}^{t_1} \, \overline{G}_2 =
\overline{G}_2 \, R_{21}^{t_2} \,
\overline{G}_1 \, R_{21}^{t_1 t_2}\; .
\end{array}
\ee
Using relations (\ref{3.1.20b}), (\ref{3.1.20})
and conditions (\ref{coid}) it can be shown directly that
the elements of quantum matrices
$$
K = L^{-} \, G \, (L^{+})^t \; , \;\;\;
\overline{K} = S(L^{+}) \, \overline{G} \, (S(L^{-}))^t \; ,
$$
obey the following commutation relations
\be
\lb{coid1}
\begin{array}{c}
R_{12} \, K_2 \, R_{12}^{t_2} \, K_1 =
K_1 \, R_{21}^{t_1} \, K_2 \, R_{21}^{t_1 t_2} \; , \\ [0.3cm]
R_{12} \, \overline{K}_1 \, R_{12}^{t_1} \, \overline{K}_2 =
\overline{K}_2 \, R_{21}^{t_2} \,
\overline{K}_1 \, R_{21}^{t_1 t_2}\; ,
\end{array}
\ee
which we consider as the defining relations for a new type of quantum
matrix algebras ${\cal K}(\R,G)$ and $\overline{{\cal K}}(\R, \overline{G})$.
The defining relations (\ref{coid1})
are covariant\footnote{Here the notion {\it covariant}
is equivalent to the statement that (\ref{blin})
are homomorphisms for the algebras defined by (\ref{coid1}).}
under the left and right ${\cal A}$-coactions:
\be
\lb{blin}
K^i_j \longrightarrow   (L^{-})^i_k \, (L^{+})^j_n \otimes K^k_n \; , \;\;\;
\overline{K}^i_j \longrightarrow
\overline{K}^k_n \otimes S (L^{+})^i_k \, S (L^{-})^j_n  \; .
\ee
Thus, the unital algebras ${\cal K}$
and $\overline{{\cal K}}$ (with generators $K^i_j$ and
$\overline{K}^i_j$, respectively) are
left and right ${\cal A}$-comodule algebras and these algebras
are called coideal subalgebras
of ${\cal A}$.

One can consider two more such algebras with generators
$K' = L^+ G' (L^-)^t$ and $\overline{K}' =
S(L^-) \overline{G}' S(L^+)^t$ which obey the following defining relations
$$
\begin{array}{c}
R^{-1}_{12} \, K'_1 \, (R_{12}^{-1})^{t_1} \, K'_2 =
K'_2 \, (R_{21}^{-1})^{t_2} \, K'_1 \, (R_{21}^{-1})^{t_1 t_2} \; ,
\\ [0.3cm]
R^{-1}_{12} \, \overline{K}'_2 \, (R_{12}^{-1})^{t_2} \, \overline{K}'_1 =
\overline{K}'_1 \, (R_{21}^{-1})^{t_1} \,
\overline{K}'_2 \, (R_{21}^{-1})^{t_1 t_2}\; .
\end{array}
$$
Note, that these relations can be obtained from (\ref{coid1})
by the substitution $R_{12} \rightarrow R_{12}^{-1}$.

For the special case of $GL_q(N)$ $R$-matrices
(see Subsec. {\bf \ref{qgrsl}}) the algebras
(\ref{coid1})
have been considered in \cite{18mol}, \cite{18nou} (see also references therein).
In this case the coideal subalgebras coincide with quantized
enveloping algebras introduced earlier by Gavrilik and Klimyk \cite{18gav}.

Representation theory for compact quantum groups
has been considered in \cite{soib}.
In \cite{DKM}
a universal solution to the reflection
equation has been introduced and
general problems of the representation theory for the reflection equation algebra
were discussed
(representations and characters for some special reflection equation
algebras were considered in \cite{Sap}).

\subsection{\bf \em The semiclassical limit (Sklyanin
 brackets and Lie bialgebras)\label{semicl}}
\setcounter{equation}0

We assume that the $R$-matrix introduced in (\ref{3.1.1}) has the following expansion
in the limit $h \rightarrow 0$ $(q = e^{h} \rightarrow 1)$:
\be
\lb{3.2.1}
R_{12} = \hbox{\bf 1} + h \, r_{12} + O(h^{2}) \; .
\ee
Here ${\hbox{\bf 1}}= I \otimes I$ denotes the $(N^{2} \times N^{2})$
unit matrix. One says that such
$R$-matrices have semiclassical behavior, and $r_{12}$ is called a classical
$r$ matrix. It is readily found from the quantum Yang-Baxter equation
(\ref{3.1.3}) that $r_{12}$ satisfies the so-called classical Yang-Baxter equation
\be
\lb{3.2.2}
[r_{12}, \; r_{13} + r_{23}] + [r_{13}, \; r_{23}] = 0 \; .
\ee
Substituting the expansion (\ref{3.2.1}) in the $RTT$ relations (\ref{3.1.1}), we obtain
\be
\lb{3.2.3}
[T_{1},T_{2}] = h [ T_{1}T_{2},r_{12} ] + O(h^{2}) \; .
\ee
This equation demonstrates the fact that the
$RTT$ relations (\ref{3.1.1}) can be interpreted as a quantization (deformation) of the
classical Poisson bracket (Sklyanin bracket \cite{19}):
\be
\lb{3.2.4}
\{ T_{1},T_{2} \} =  [ T_{1}T_{2},r_{12} ]  \; ,
\ee
(here the elements $T^i_j$ are commutative coordinates of some Poisson manifold).
The classical Yang-Baxter equation (\ref{3.2.2}) guarantees fulfillment of the
Jacobi identity for the bracket (\ref{3.2.4}). From the requirement of
antisymmetry of the Poisson bracket (\ref{3.2.4}), we obtain
\be
\lb{3.2.5}
\{ T_{1},T_{2} \} =  [ T_{1}T_{2},-r_{21} ]  \; .
\ee
Thus, the classical $r$ matrix $r^{(-)}_{12} = -r_{21}$ corresponding to the
representation $R^{(-)}$
(\ref{3.1.18}) must also be a solution of Eq. (\ref{3.2.2}), as is readily shown
by making the substitution  $3 \leftrightarrow 1$ in (\ref{3.2.2}). On the other hand, comparing
(\ref{3.2.4}) and (\ref{3.2.5}), we obtain
\be
\lb{3.2.6}
T_{1}T_{2} (r_{12} + r_{21}) = (r_{12} + r_{21}) T_{1}T_{2}  \; .
\ee
Thus,
\be
\lb{3.2.7}
t_{12} = \frac{1}{2} (r_{12} + r_{21})
\ee
is an invariant with respect to the adjoint action of
the matrix $T_{1}T_{2}$ (it is an ad-invariant). We introduce the new classical
$r$ matrix
\be
\lb{3.2.8}
\tilde{r}_{12} = \frac{1}{2} (r_{12} - r_{21})  \; .
\ee
Then the Sklyanin bracket can be represented in the manifestly antisymmetric form
\be
\lb{3.2.9}
\{ T_{1},T_{2} \} =  [ T_{1}T_{2},\tilde{r}_{12} ]  \; ,
\ee
and the matrix $\tilde{r}$ (\ref{3.2.8}) satisfies the modified classical
Yang-Baxter equation
\be
\lb{3.2.10}
[\tilde{r}_{12}, \; \tilde{r}_{13} + \tilde{r}_{23}]
+ [\tilde{r}_{13}, \; \tilde{r}_{23}] =
{1 \over 4} \, [r_{23}+r_{32}, \; r_{13} + r_{31}] = [t_{23}, \; t_{13}] \; .
\ee

Note that the reflection equation
algebras (\ref{3.1.23a}), (\ref{3.1.23b}) can also be regarded as the result of
quantization of a certain Poisson structure. For example,
 for these algebras, after substitution of (\ref{3.2.1}),
 we have \cite{STS} (see also \cite{Isaev2})
$$
\{ L_{2}, \, L_{1} \} = [ L_{1}, \,  [ L_{2}, \, \tilde{r}_{12} ] ] +
L_{1} t_{12} L_{2} - L_{2} t_{12} L_{1} \; ,
$$
$$
\{ \overline{L}_{2}, \, \overline{L}_{1} \} =
- [ \overline{L}_{1}, \,  [ \overline{L}_{2}, \, \tilde{r}_{12} ] ] +
\overline{L}_{1} t_{12} \overline{L}_{2} - \overline{L}_{2} t_{12} \overline{L}_{1} \; ,
$$
where again we must assume that $[L_{1}L_{2}, \, t_{12}] = 0 =
[\overline{L}_{1}\overline{L}_{2}, \, t_{12}]$
[cf. (\ref{3.2.6})]. On the other hand, the relations (\ref{3.1.25})
in the zeroth order in $h$ give the equations
$$
[ \tilde{L}_{1}, \, \tilde{L}_{2} ] =
[ t_{12} , \, \tilde{L}_{1} ] \;\; , \;\;\;\;
\left( [t_{12} , \, \tilde{L}_{1} + \tilde{L}_{2}]=0 \right) \; ,
$$
and this enables us to regard (\ref{3.1.25}) as a deformation of the
defining relations of a Lie algebra.

Now we consider the universal enveloping $U(\mathfrak{g})$ of a
Lie algebra $\mathfrak{g}$ with defining relations (\ref{lie})
as a bialgebra (\ref{2.14}) and assume that the
cocommutative comultiplication $\Delta$ (\ref{2.14}) is quantized
$\Delta \rightarrow \Delta_h$ in such a way that $\Delta_h$ is noncocommutative.
The semiclassical expansion of $\Delta_h$ is\footnote{The terms
$h \phi_\alpha^{\beta \gamma} J_\beta^1 J_\gamma^1$ and
$h \phi_\alpha^{\beta \gamma} J_\beta^2 J_\gamma^2$ are gauged out
by triviality transformation from this expansion (see, e.g., \cite{18a}).}
\be
\lb{colie}
\Delta_h (J_\alpha)  = J^1_\alpha + J^2_\alpha + h \,
\mu^{\beta \gamma}_{\alpha} \, J^1_\beta \, J^2_\gamma + h^2 \,
\left(  \mu^{\beta_1 \beta_2, \gamma}_{\alpha} \, J^1_{\beta_1} \, J^1_{\beta_2} \, J^2_\gamma
+  \mu^{\beta , \gamma_1 \gamma_2}_{\alpha} \,
J^1_{\beta} \, J^2_{\gamma_1} \, J^2_{\gamma_2} \right) + h^3 \dots
\ee
where $J^1_\alpha = J_\alpha \otimes 1$, $J^2_\alpha = 1 \otimes J_\alpha$,
the term of zeroth order in $h$ in (\ref{colie}) is the classical comultiplication (\ref{2.14})
and $\mu^{\beta \gamma}_{\alpha}, \mu^{\beta_1 \beta_2, \gamma}_{\alpha}, \dots$
are some constants. The comultiplication map
(\ref{colie}) (as well as the opposite comultiplication $\Delta^{\sf op}_h$; see (\ref{2.13}))
should be a homomorphic map for the Lie algebra (\ref{lie}):
\be
\lb{cl2}
[ \Delta_h(J_{\alpha}) , \,  \Delta_h(J_{\beta}) ] =
t_{\alpha\beta}^{\gamma} \, \Delta_h(J_{\gamma}) \; , \;\;\;
[ \Delta^{\sf op}_h(J_{\alpha}) , \,  \Delta^{\sf op}_h(J_{\beta}) ] =
t_{\alpha\beta}^{\gamma} \, \Delta^{\sf op}_h(J_{\gamma}) \; .
\ee
Then, the subtraction of the second relation
of (\ref{cl2}) from the first one gives the following equation
$$
[ \Delta^-_h(J_{\alpha}) , \,  \Delta^+_h(J_{\beta}) ] +
[ \Delta^+_h(J_{\alpha}) , \,  \Delta^-_h(J_{\beta}) ] =
t_{\alpha\beta}^{\gamma} \, \Delta^-_h(J_{\gamma}) \; ,
$$
(here we define $\Delta^-_h : = \Delta_h - \Delta^{\sf op}_h$ and
$\Delta^+_h : = \frac{1}{2} \, ( \Delta_h + \Delta^{\sf op}_h)$)
which is rewritten (in the first order of $h$) as
\be
\lb{bilie1}
[ \delta(J_{\alpha}) , \,  J^1_{\beta} +  J^2_{\beta} ] +
[ J^1_{\alpha} +  J^2_{\alpha} , \,  \delta(J_{\beta}) ] =
t_{\alpha\beta}^{\gamma} \, \delta(J_{\gamma}) \; ,
\ee
where the map $\delta$: $g \rightarrow g \wedge g$ is
\be
\lb{bilie2}
\delta(J_{\alpha}) = \delta^{\beta \gamma}_{\alpha} \, J_\beta \otimes J_\gamma \; , \;\;\;
\delta^{\beta \gamma}_{\alpha} := \mu^{\beta \gamma}_{\alpha} -
\mu^{\gamma \beta }_{\alpha} \; .
\ee
Eq. (\ref{bilie1}) is nothing but the cocycle condition for
$\delta^{\beta \gamma}_{\alpha}$:
$$
\left( \delta^{\rho \mu}_\alpha \, t^\kappa_{\beta \rho} -
\delta^{\rho \kappa}_\alpha \, t^\mu_{\beta \rho} \right) -
\left( \delta^{\rho \mu}_\beta \, t^\kappa_{\alpha \rho} -
\delta^{\rho \kappa}_\beta \, t^\mu_{\alpha \rho} \right) =
t^\gamma_{\alpha \beta} \, \delta^{\mu \kappa}_\gamma \; .
$$
On the other hand the structure constants $(\Delta^-)^{ij}_k = \Delta^{ij}_k
- \Delta^{ji}_k$ satisfy the co-Jacobi identity
$$
(\Delta^-)^{jk}_i \, (\Delta^-)^{nm}_j + (\Delta^-)^{jn}_i \, (\Delta^-)^{mk}_j +
(\Delta^-)^{jm}_i \, (\Delta^-)^{kn}_j = 0 \; ,
$$
as it is evident from the co-associativity condition (\ref{2.7}). This identity
for the comultiplication (\ref{colie}) in the order $h^2$ reduces to the
co-Jacoby identity for the
structure constants $\delta^{\beta \gamma}_{\alpha}$ (\ref{bilie2}):
\be
\lb{bilie3}
\delta^{\beta \gamma}_{\alpha} \, \delta^{\rho \xi}_{\beta} + ({\rm cycle}
\;\; \gamma, \rho, \xi) = 0 \; .
\ee
Thus, we have arrived to the following definition \cite{13}:


\newtheorem{def10}[def1]{Definition}
\begin{def10} \label{def10}
{\it The vector space $\mathfrak{g}$ with the basis $\{ J_\alpha \}$ equipped with a linear map $\delta$: $g
\rightarrow g
\wedge g$ (\ref{bilie2}) satisfying the co-Jacobi identity (\ref{bilie3}) is called Lie coalgebra.
A Lie bialgebra is a Lie algebra (\ref{lie}) which is in the same time is a Lie coalgebra with the
map $\delta$: $g \rightarrow g \wedge g$ (\ref{bilie2}), (\ref{bilie3}) satisfying the cocycle
condition (\ref{bilie1}). }
\end{def10}

Let $\mathfrak{g}$ be a Lie bialgebra.
If there exists an element $r \in g \otimes g$ such that the map $\delta$
has the form
$$
\delta(J) =  [ J \otimes 1 + 1 \otimes J, \, r] \;\;\; \forall J \in g \; ,
$$
then the Lie bialgebra $\mathfrak{g}$ is called co-boundary or $r$-matrix bialgebra.

\subsection[\bf \em The quantum groups $GL_q(N)$, $SL_q(N)$ and
their quantum algebras and hyperplanes]{\bf \em The quantum groups $GL_q(N)$, $SL_q(N)$ and
their quantum \\ algebras and hyperplanes\label{qgrsl}}
\setcounter{equation}0

\subsubsection{$GL_q(N)$ quantum hyperplanes and $R$ matrices}

In this subsection, we discuss the simplest nontrivial quantum
groups, which are the quantizations (deformations) of the linear Lie
groups $GL(N)$ and $SL(N)$. We begin with the definition of a quantum hyperplane.
We recall that the Lie group $GL(N)$ is the set of nondegenerate
$N \times N$ matrices $T^{i}_{j}$
that act on an $N$-dimensional vector space, whose coordinates we
denote by $x^{i}$, ($i=1, \dots N$). Thus, we have the transformations
\be
\lb{3.3.1}
x^{i} \rightarrow \tilde{x^{i}} = T^{i}_{j}x^{j} ,
\ee
which we can regard from a different point of view.
Namely, let $\{ T^{i}_{j} \}$ and $\{ x^{i} \}$ $(i,j=1, \dots ,N)$
be the generators of two Abelian (commuting) algebras
\be
\lb{3.3.2}
[ x^{i}, x^{j}] = [T^{i}_{j},T^{k}_{l}] = [ T^{i}_{j}, x^{k}] =0 \; .
\ee
Then the transformation (\ref{3.3.1}) can be regarded as an action
(more precisely, it is a coaction) of the algebra
$\{ T \}$ on the algebra $\{ x \}$
\be
\lb{3.3.1bb}
x^{i} \rightarrow \delta_T(x^{i}) \equiv \tilde{x^{i}} = T^{i}_{j} \otimes x^{j} ,
\ee
that preserves the Abelian structure of the latter, i.e.,
we have $[ \tilde{x}^{i}, \, \tilde{x}^{j} ] = 0$.
We introduce a deformed $N$-dimensional "vector space" whose
coordinates $\{ x^{i} \}$ commute as follows:
\be
\lb{3.3.3}
x^{i}x^{j} = q x^{j} x^{i}, \; i < j
\ee
where $q$ is some number (the deformation parameter). In other words, we now
have a noncommutative associative algebra with $N$ generators $\{ x^{i} \}$.
In accordance with (\ref{3.3.3}), any element of this algebra,
which is a monomial of arbitrary degree
\be
\lb{3.3.4}
x^{i_{1}}x^{i_{2}} \cdots  x^{i_{K}},
\ee
can be uniquely ordered lexicographically, i.e., in such a way
that $i_{1} \leq i_{2} \leq \dots \leq i_{K}$. Of such algebras, one says that they possess the
Poincare-Birkhoff-Witt (PBW) property. An
algebra with $N$ generators satisfying (\ref{3.3.3}) is called
an $N$-dimensional quantum hyperplane \cite{20},\cite{21}.
The relations (\ref{3.3.3}) can be written in the matrix form
\be
\lb{3.3.5}
R^{i_{1}i_{2}}_{j_{1}j_{2}}x^{j_{1}}x^{j_{2}} =
q x^{i_{2}} x^{i_{1}} \Leftrightarrow R_{12}x_{1}x_{2} = q x_{2}x_{1}
\Leftrightarrow \R \, x_1 \, x_2 = q \, x_1 \, x_2 \; .
\ee
Here the indices $1$ and $2$ label the vector spaces on which the
$R$-matrix, realized in the tensor square
$Mat(N) \otimes Mat(N) =: Mat(N)_{1} \; Mat(N)_{2}$, acts.
Thus, the indices $1$ and $2$ of the $R$-matrix show how the
$R$-matrix acts on the direct product of the first and second vector spaces.
We emphasize that the $R$-matrix depends on the parameter $q$ and, generally
speaking, its explicit form is recovered nonuniquely from the relations (\ref{3.3.3}).
However, if we require that the $R$-matrix (\ref{3.3.5}) be constructed by means of two
$GL(N)$-invariant tensors $\hbox{\bf 1}_{12}$ and $P_{12}$,
 i.e.,\footnote{The form of $R$-matrix (\ref{anzdym})
 proves to be very fruitful for
the construction of solutions for the dynamical Yang-Baxter equations; see \cite{21'}, \cite{21''} and references therein (see also
subsection {\bf \ref{baxtel}} below).}
\be
\lb{anzdym}
R^{i_{1}i_{2}}_{j_{1}j_{2}} =
(\delta^{i_{1}}_{j_{1}} \delta^{i_{2}}_{j_{2}} )
\cdot a_{i_{1}i_{2}} +
(\delta^{i_{1}}_{j_{2}} \delta^{i_{2}}_{j_{1}} )
\cdot b_{i_{1}i_{2}} \; ,
\ee
and also satisfy the Yang-Baxter equation (\ref{ybe})
 and have lower-triangular
block form ($R^{i_{1}i_{2}}_{j_{1}j_{2}} = 0$, $i_{1} < j_{1}$),
then we obtain the explicit expression  \cite{Jimb1}, \cite{10}
\begin{equation}
\lb{3.3.6a}
R_{12}=  q \sum_{i} e_{ii} \otimes e_{ii} +
\sum_{i \neq j} e_{ii} \otimes e_{jj} +
\lambda \sum_{i > j} e_{ij} \otimes e_{ji} \; ,
\end{equation}
\begin{equation}
\lb{3.3.6b}
\R_{12}= P_{12} \, R_{12}= q \sum_{i} e_{ii} \otimes e_{ii} +
\sum_{i \neq j} e_{ij} \otimes e_{ji} +
\lambda \sum_{i > j} e_{jj} \otimes e_{ii} \; ,
\end{equation}
where  $e_{ij}|_{_{i,j = 1, \dots ,N}}$
are matrix units: $(e_{ij})^{k}_{l} =
\delta^{ik} \delta_{jl}$ ,
 $P_{12} := \sum_{k,\ell} e_{k\ell} \otimes e_{\ell k}$
 is a permutation matrix
 and here and below we often use notation
 \fbox{$\; \lambda : = q-q^{-1} \,$} .
Eq. (\ref{3.3.6a}) is represented in the components
in the form
\begin{equation}
\lb{3.3.6}
\begin{array}{c}
R^{i_{1}i_{2}}_{j_{1}j_{2}}=
\delta^{i_{1}}_{j_{1}} \delta^{i_{2}}_{j_{2}}(1+(q-1)\delta^{i_{1}i_{2}}) +
 \lambda \delta^{i_{1}}_{j_{2}} \delta^{i_{2}}_{j_{1}}
\Theta_{i_{1}i_{2}} \; ,  \\ [0.3cm]
\hat{R}^{i_{1}i_{2}}_{j_{1}j_{2}}=
\delta^{i_{1}}_{j_{2}} \delta^{i_{2}}_{j_{1}} \, q^{\delta_{i_{1}i_{2}}} +
 \lambda \delta^{i_{1}}_{j_{1}} \delta^{i_{2}}_{j_{2}}
\Theta_{i_{2}i_{1}} \; ,  \\  [0.3cm]
\Theta_{ij} = \{ 1 \;\; {\rm if} \;\; i>j , \;\; 0 \;\; {\rm if} \;\; i \leq j \}.
\end{array}
\end{equation}
 It can be verified (by using, e.g.,
the diagrammatic technique of Sec. {\bf \ref{multpar}}) that this
$R$-matrix satisfies the Hecke relation (\ref{3.3.7})
(a special case of (\ref{3.1.27}))
\be
\lb{3.3.7aa}
R_{12} - \lambda \, P_{12} - R^{-1}_{21} = 0
\;\; \Rightarrow \;\;
\R - \lambda I - \R^{-1} = 0 \;\;\; \Leftrightarrow \;\;\;
\R^2 = \lambda \R + I = 0 \; ,
\ee
where $I^{i_{1}i_2}_{j_{1}j_2} = \delta^{i_{1}}_{j_{1}} \delta^{i_{2}}_{j_{2}}$ is a unit operator.
The following helpful relations also follow from the explicit form
(\ref{3.3.6a}), (\ref{3.3.6}) for the $GL_q(N)$ $R$-matrix:
\be
\lb{3.3.7a1q}
R_{12}[\frac{1}{q}] = R_{12}^{-1}[q]
\;\;\;\; \Leftrightarrow \;\;\;\;
\R_{12}[\frac{1}{q}] = \R_{21}^{-1}[q] \; ,
\ee
\be
\lb{3.3.7ad}
R_{12}^{t_{1}t_{2}} = R_{21} \; ,  \;\;\; R_{12}^{t_{1}} \, R_{12} = R_{12} \, R^{t_{1}}_{12} \; .
\ee
The $R$-matrix (\ref{anzdym}) (where without loss of generality one can fix $b_{ii}=0$) is
skew-invertible iff $a_{ij} \neq 0$
$(\forall i,j)$ and $\det(||b_{ij}+a_{ii} \delta_{ij}||) \neq 0$. Then the skew-inverse matrix
$\Psi_{12}$ (\ref{skew}) is represented in the form:
\be
\lb{skewdym}
\hat{\Psi}^{i_1 i_2}_{j_1 j_2} = \frac{1}{a_{_{i_2 i_1}}} \,
\delta^{i_1}_{j_2} \delta^{i_2}_{j_1}  -
d_{i_2 i_1} \; \delta^{i_1}_{j_1} \delta^{i_2}_{j_2}  \; ,
\ee
where coefficients $d_{i j}$ are defined by the matrix equation
$$
d = A^{-1} \, B \, (A + B)^{-1} \; , \;\;\; A_{ij} : = a_{ii} \delta_{ij}
 \; , \;\;\; B  = || b_{ij} || \; .
$$
For the given $R$-matrix
(\ref{3.3.6}), the matrix $\hat{\Psi}_{12}$ (\ref{skewdym}) is calculated in the form
\cite{18a}
\be
\lb{skewdym1}
\begin{array}{c}
\hat{\Psi}^{i_1 i_2}_{j_1 j_2} = q^{- \delta_{i_1 i_2}}
\delta^{i_1}_{j_2} \delta^{i_2}_{j_1}  - \lambda \Theta_{i_2 i_1} q^{2(i_1 - i_2)}
\delta^{i_1}_{j_1} \delta^{i_2}_{j_2}  \; , \\ [0.2cm]
\hat{\Psi}_{12} = q^{-1} \sum_{i} e_{ii} \otimes e_{ii} +
\sum_{i \neq j} e_{ij} \otimes e_{ji} -
\lambda \sum_{i < j} q^{2(i-j)} e_{ii} \otimes e_{jj} \; ,
\end{array}
\ee
Then the quantum trace matrices $D$, $Q$ (\ref{dmatr}) and the related quantum traces
(\ref{3.1.13}) are
\be
\lb{3.3.14}
\begin{array}{c}
D_1 \equiv  Tr_{2}\left( \hat{\Psi}_{12} \right)
= diag \{ q^{-2N+1},q^{-2N+3}, \dots , q^{-1} \} \; , \;\;\;
D^i_j = q^{2(i-N)-1} \delta^i_j \; , \\ [0.3cm]
Q_2 \equiv  Tr_{1}\left( \hat{\Psi}_{12} \right)
= diag \{ q^{-1},   \dots ,q^{-2N+3}, q^{-2N+1}\} \; , \;\;\;
Q^i_j = q^{1-2i} \delta^i_j \; ,
\end{array}
\ee
$$
Tr_{D}A := Tr(DA)
\equiv \sum_{i=1}^{N} q^{2(i-N)-1} A^{i}_{i} \; , \;\;
Tr_{Q}A := Tr(QA)
\equiv \sum_{i=1}^{N} q^{1-2i} A^{i}_{i} \; .
$$
We also note the useful relations [cf. (\ref{qtrs}), (\ref{3.1.17})]
\be
\lb{3.3.15}
\begin{array}{c}
Tr_{D}(I) = Tr(D) = q^{-N} \, [N]_{q}  = Tr(Q) = Tr_{Q}(I) \; , \\ [0.3cm]
q^N \, Tr_{D(3)}\R_{23}^{\pm 1} = q^{\pm N} \cdot I_{(2)} =
q^N \, Tr_{Q(1)}\R^{\pm 1} \; , \;\;\;
\end{array}
\ee
where $[N]_{q} = {q^{N} -q^{-N} \over q-q^{-1}}$. One can readily prove the cyclic property of the
quantum traces
\be
\lb{3.3.16}
Tr_{D(12)} ( \R E_{12} ) = Tr_{1 2} (D_1 D_2  \R E_{12} ) = Tr_{1 2} (\R D_1 D_2 \, E_{12} ) =
Tr_{D(12)} ( E_{12} \R ) \; ,
\ee
where $E_{12} \in Mat(N) \otimes Mat(N)$ is a matrix with noncommutative entries. In
(\ref{3.3.16}) we have used  the fact that the matrix
$D$,  by definition, obeys eq. $[\R, \, D_1 D_2] = 0$ (\ref{sk10})
(note that for $R$-matrices of the type (\ref{anzdym}), all diagonal matrices $D$ satisfy this equation).
The same cyclic property $Tr_{Q(12)} ( \R E_{12} ) =Tr_{Q(12)} ( E_{12} \R )$ is also valid for the
traces $Tr_{Q}$.

In semiclassical limit (\ref{3.2.1}), relation (\ref{3.3.7aa})
can be written in the form
\be
\lb{3.3.8}
r_{12}  +  r_{21} = 2P_{12} \; .
\ee
Thus, for the Lie-Poisson structure on the group
$GL(N)$ the transposition matrix $P_{12}$ is taken as the $ad$-invariant
tensor $t_{12}$. For the $\tilde{r}$ matrix (\ref{3.2.8}) determining the
Sklyanin bracket, we obtain from
(\ref{3.3.6}) the expression
\be
\lb{3.3.9}
\tilde{r}_{12} = \sum_{i > j}
[ e_{ij} \otimes e_{ji} - e_{ji}  \otimes e_{ij} ]
\in gl(N) \wedge gl(N) \; .
\ee

In accordance with (\ref{3.3.7}), (\ref{3.1.28}), and (\ref{3.1.29}) for
$q^{2} \neq -1$ the matrix $\R$ has the spectral decomposition
\be
\R  = q \P^{+} - q^{-1} \P^{-} \; ,
\lb{3.3.10}
\ee
with projectors
\be
\lb{3.3.11}
\P^{\pm}  = (q + q^{-1})^{-1} \{ q^{\mp 1} {\bf 1} \pm \R \} \; ,
\ee
which are the quantum analogs of the symmetrizer
($\P^{+}$) and antisymmetrizer ($\P^{-}$),
as can be seen by setting $q=1$ in (\ref{3.3.11}). Using the projector $\P^{-}$,
we can represent the definition (\ref{3.3.3}) of the quantum hyperplane in the form
\be
\lb{3.3.12}
\P^{-} \, x_1 \, x_2 = 0.
\ee
Note that the relations
\be
\lb{3.3.13}
\P^{+} \, x_1 \, x_2 = 0 \Leftrightarrow (x^{i})^{2} = 0 \; , \;\;
x^{i} x^{j} = - q^{-1} x^{j} x^{i} \;\; (i < j)
\ee
define a fermionic $N$-dimensional quantum hyperplane that is a deformation
of the algebra of $N$ fermions: $x^{i}x^{j}= - x^{j}x^{i}$.

\subsubsection{Quantum groups $Fun(GL_q(N))$, $Fun(SL_q(N))$ and
$q$-determinants}

A natural question now is about of the properties
of the $N \times N$ matrix elements $T^{i}_{j}$ that determine the transformations
(\ref{3.3.1bb}) of the quantum bosonic (\ref{3.3.3}), (\ref{3.3.12}) and fermionic
(\ref{3.3.13}) hyperplanes. These properties should be such that the transformed coordinates
$\tilde{x}^i$ form the same quantum algebras ($q$- hyperplanes)
(\ref{3.3.12}) and (\ref{3.3.13}). It is readily seen that the elements of the
$N \times N$ matrix $T^{i}_{j}$ must satisfy the both conditions
\be
\lb{3.3.17}
\P^{+} \, T_1 \,  T_2 \, \P^{-} = 0 \; , \;\;\;\;\;
\P^{-} \, T_1 \,  T_2 \, \P^{+} = 0 \; .
\ee
Indeed, we have for bosonic $x^-$ and fermionic $x^+$ hyperplanes
(we omit the symbol $\otimes$ in (\ref{3.3.1bb}))
$$
0 = \P^{\pm} \, \tilde{x}^{\pm}_1 \, \tilde{x}^{\pm}_2 =
\P^{\pm} \, T_1 \, T_2 \, x^{\pm}_1 \, x^{\pm}_2 =
$$
$$
= \P^{\pm} \, T_1 \, T_2 \, (\P^+ + \P^-) \, x^{\pm}_1 \, x^{\pm}_2 =
 \P^{\pm} \, T_1 \, T_2 \, \P^{\mp} \, x^{\pm}_1 \, x^{\pm}_2 \; ,
$$
and we deduce (\ref{3.3.17}) (otherwise new quadratic relations on the coordinates $x^{\pm}$ should be
imposed). Eqs. (\ref{3.3.17}) are equivalent to the $RTT$ relations (\ref{3.1.1}) for the
elements of the $N \times N$ quantum matrix $||T^i_j||$:
 \be
\lb{3.3.17i}
\R \, T_1 \,  T_2 - T_1 \,  T_2 \, \R  =
(q+q^{-1}) \bigl( \P^{+} \, T_1 \,  T_2 \, \P^{-} -
\P^{-} \, T_1 \,  T_2 \, \P^{+}  \bigr) = 0 \; .
 \ee
We note that one can define the quantum matrix
algebra when only one of two relations in (\ref{3.3.17})
is fulfilled. In this case the quantum matrix algebras, generated
by $T^i_j$, are called {\it half-quantum} or {\it Manin} matrix
algebras \cite{IsOg5}, \cite{ChFR}.

For the $R$-matrix (\ref{3.3.6}) the $RTT$ relations (\ref{3.1.1}) and (\ref{3.3.17i}) can be
written in the component form
\be
\lb{glrtt}
\begin{array}{c}
T^i_k \, T^j_k = q \, T^j_k \, T^i_k \; , \;\;\;
T_i^k \, T_j^k = q \, T_j^k \, T_i^k \; , \;\; (i < j, \, k=1, \dots,N) \; , \\ [0.3cm]
\; [T^{i_1}_{j_1}, \, T^{i_2}_{j_2}] =
(q-q^{-1}) \, T^{i_1}_{j_2} \, T^{i_2}_{j_1} \; , \;\;\;
[ T^{i_1}_{j_2}, \, T^{i_2}_{j_1}] = 0 \; , \;\;
(i_1 < i_2, \; j_1 < j_2)  \; .
\end{array}
\ee
The $RTT$ algebra with defining relations (\ref{glrtt}) is a bialgebra with the structure mappings
$\Delta,\epsilon$ presented in (\ref{3.1.5}).
 The simplest special case
$(N=2)$ of this algebra is defined by
\be
\lb{glrtt2}
\begin{array}{c}
T^1_{\; k} \, T^2_{\; k} = q \, T^2_{\; k} \, T^1_{\; k}
\; , \;\;\;
T_{\; 1}^k \, T_{\; 2}^k = q \, T_{\; 2}^k \, T_{\; 1}^k
\; , \;\; (k=1,2) \; , \\ [0.3cm]
\; [T^{1}_{\; 1}, \, T^{2}_{\; 2}] =
(q-q^{-1}) \, T^{1}_{\; 2} \, T^{2}_{\; 1} \; , \;\;\;
[ T^{1}_{\; 2}, \, T^{2}_{\; 1}] = 0  \; .
\end{array}
\ee
One can directly check that
${\det}_q(T) := T^{1}_{\; 1} \, T^{2}_{\; 2} -
q  T^{1}_{\; 2} \, T^{2}_{\; 1} \equiv
T^{2}_{\; 2} \, T^{1}_{\; 1} -
q^{-1}  T^{1}_{\; 2} \, T^{2}_{\; 1}$ is a central element for the algebra (\ref{glrtt2}). This element is
called quantum determinant for $(2\times 2)$ quantum matrix $||T^i_j||$, since for $q=1$ the
element ${\det}_q(T)$ coincides with the usual determinant.
Let ${\det}_q(T)$ be invertible element. Then
the inverse matrix $T^{-1}$ is
$$
T^{-1} = \left(
\begin{array}{cc}
T^{2}_{\; 2}& - q^{-1} T^{1}_{\; 2} \\
- q T^{2}_{\; 1} & T^{1}_{\; 1}
\end{array}
\right) \frac{1}{{\det}_q(T)} \;\;\;\;\; \Rightarrow \;\;\;\;\;
{\det}_{q^{-1}} (T^{-1}) = {\det}_q^{-1}(T) \; .
$$

Now we generalize the definition of the quantum determinant
for the case of $(N \times N)$ quantum
matrices
$||T^i_j||$. We introduce the quantum determinant ${\det}_q(T)$, which is a deformation of the ordinary
determinant and also is a central element for the $RTT$-algebra (\ref{glrtt}). For this aim we
introduce the
$q$-deformed  antisymmetric tensors
${\cal E}_{j_{1}j_{2} \dots j_{N}}$ and
${\cal E}^{j_{1}j_{2} \dots j_{N}}$ $(\forall j_k =1, \dots, N)$ as follows:
$$
\sum_{j_1 \dots j_N =1}^N \,
{\cal E}_{j_1 j_2 \dots j_N} \, {\cal E}^{j_1 j_2 \dots j_N}  =
{\cal E}_{\langle 12 \dots N} {\cal E}^{12 \dots N \rangle}  = 1 \; ,
$$
\be
\lb{3.3.19}
\begin{array}{c}
{\cal E}_{\langle 12 \dots N} \P^{+}_{k,k+1} =
{\cal E}_{\langle 12 \dots N} (\R_{k,k+1} +q^{-1}) = 0 \;, \;\; 1 \leq k < N \; , \\ [0.3cm]
 \P^{+}_{k,k+1} {\cal E}^{12 \dots N\rangle } =
 (\R_{k,k+1} +q^{-1}) {\cal E}^{12 \dots N\rangle } = 0 \;, \;\; 1 \leq k < N \; ,
\end{array}
\ee
where we have used concise matrix notations. Namely, we denoted by
$12 \dots N\rangle$ and
$\langle 12 \dots N$ -- sets of incoming and outgoing indices, where
$1,2, \dots N$ are numbers of
the $N$-dimensional vector spaces $V_N$, and
$\P^{+}_{k,k+1}= I^{\otimes (k-1)}\otimes
\P^{+} \otimes I^{\otimes (N-k-1)}$ are the
symmetrizers (\ref{3.3.10}) acting in the vector spaces $V_N$ labeled by numbers $k$ and $k+1$.
Note that, in view of the $RTT$ relations (\ref{3.1.1}),
(\ref{glrtt}), the tensors
${\cal E}_{\langle 12 \dots N}
(T_{1} T_{2} \cdots T_{N})$ and
$\;(T_{1} T_{2} \cdots T_{N}){\cal E}^{12 \dots N\rangle}$
possess the same symmetry\footnote{It is not true for
 the half-matrix algebras (see definition after eq. (\ref{3.3.17i})).}
(\ref{3.3.19}) as tensors
${\cal E}_{\langle 12 \dots N}$ and ${\cal E}^{12 \dots N\rangle }$,
respectively. Supposing that the ${\cal E}$-tensors are unique
(up to a normalization), one can write
\be
\lb{3.3.18}
\begin{array}{c}
{\det}_{q}(T) \, {\cal E}_{j_{1}j_{2} \dots j_{N}}  =
{\cal E}_{i_{1}i_{2} \dots i_{N}}
T^{i_{1}}_{j_{1}} \cdot T^{i_{2}}_{j_{2}} \cdots
T^{i_{N}}_{j_{N}}  \; , \\ [0.3cm]
 {\cal E}^{ i_{1}i_{2} \dots i_{N}} \, {\det}_{q}(T)  =
T^{i_{1}}_{j_{1}} \cdot T^{i_{2}}_{j_{2}} \cdots
T^{i_{N}}_{j_{N}} {\cal E}^{ j_{1}j_{2} \dots j_{N}} \; ,
\end{array}
\ee
or in concise matrix notations we have
\be
\lb{3.3.20}
{\det}_{q}(T) \, {\cal E}_{\langle 12 \dots N}
= {\cal E}_{\langle 12 \dots N} \,
T_{1} \cdot T_{2} \cdots T_{N} \; , \;\;\;\;\;\;
 {\cal E}^{12 \dots N\rangle } \, {\det}_{q}(T)
= T_{1} \cdot T_{2} \cdots T_{N} \,
 {\cal E}^{12 \dots N\rangle } \; ,
\ee
where $T_{m} := I^{\otimes (m-1)} \otimes T \otimes I^{\otimes (N-m)}$. The scalar coefficient
${\det}_{q}(T)$:
\be
\lb{3.3.20a}
{\det}_{q}(T) = {\cal E}_{\langle 12 \dots N}
\left( T_{1} T_{2} \cdots T_{N} \right) {\cal E}^{12 \dots N\rangle } =
 Tr_{12 \dots N} (A_{1 \to N} \, T_1 T_2 \cdots T_N) \; ,
\ee
is called the quantum determinant for the $(N\times N)$ quantum matrix $||T^i_j||$. In
(\ref{3.3.20a}) we introduced the rank 1 projector
\be
\lb{Aeps}
\begin{array}{c}
A_{1 \to N} := {\cal E}^{12 \dots N\rangle } {\cal E}_{\langle 12 \dots N} \; , \;\;\;
A_{1 \to N} \, A_{1 \to N} = A_{1 \to N} \; , \\[0.2cm]
A_{1 \to N} \, \P^{+}_{k,k+1} = \P^{+}_{k,k+1} \, A_{1 \to N} = 0
\; , \;\; 1 \leq k < N \; ,
\end{array}
\ee
which acts as an $q$-antisymmetrizer in the tensor product $V_N^{\otimes N}$ of $N$ copies of
vector spaces $V_N$. It is worth noting that the $q$-antisymmetrizers
$A_{1 \to 2} := \P^{-}_{1,2}$ (\ref{3.3.11}) and $A_{1 \to N}$ are two special
representatives of the set of antisymmetrizers $\{ A_{1 \to m}\}$ $(m = 2,3, \dots ,N)$ which act
in the tensor product of $m$ vector spaces $V_N$ and satisfy
\be
\lb{Aeps1}
\begin{array}{c}
A_{1 \to m} \, A_{1 \to m} = A_{1 \to m} \; , \\[0.2cm]
A_{1 \to m} \, \P^{+}_{k,k+1} = \P^{+}_{k,k+1} \, A_{1 \to m} = 0
\; , \;\; 1 \leq k < m \; .
\end{array}
\ee
All of them can be explicitly constructed in terms of the $R$-matrices (\ref{3.3.6a}),
(\ref{3.3.6}) (see, e.g., \cite{Gur}, \cite{Jimb1}
 and Subsection {\bf \ref{hRMa}} below).

The fact that ${\det}_{q}(T)$ is indeed a central element in the $RTT$ algebra
(\ref{glrtt}) can be obtained as follows
\be
\lb{3.3.21a}
 {\cal E}_{\langle 12 \dots N} \, {\det}_{q}(T) \, T_{N+1}
= {\cal E}_{\langle 12 \dots N} \,
T_{1} \, T_{2} \cdots T_{N} \, T_{N+1} =
\ee
$$
= {\cal E}_{\langle 12 \dots N} \, (R_{1,N+1}  \cdots R_{N,N+1})^{-1}
 T_{N+1} \, T_{1} \, T_{2} \cdots T_{N} \, (R_{1,N+1}  \cdots R_{N,N+1}) =
$$
$$
=  q^{-1} \, T_{N+1}  {\cal E}_{\langle 12 \dots N} \, T_{1} \, T_{2} \cdots T_{N}
\, (R_{1,N+1}  \cdots R_{N,N+1}) =
T_{N+1} \, {\det}_q(T) \,  {\cal E}_{\langle 12 \dots N} \; ,
$$
where we have used the definition
(\ref{3.3.20}), the $RTT$ relations presented in the form
$T_{m} T_{m+1}= R_{m,m+1}^{-1} T_{m+1} T_{m} R_{m,m+1}$,
and the equations
\be
\lb{3.3.21}
\begin{array}{c}
q I_{N+1} \, {\cal E}_{\langle 12 \dots N}
= {\cal E}_{\langle 12 \dots N} \,
R_{1,N+1} \cdot R_{2,N+1} \cdots R_{N,N+1} \\ [0.2cm]
q^{-1} I_{N+1} \, {\cal E}_{\langle 12 \dots N}
= {\cal E}_{\langle 12 \dots N} \,
R_{N+1,1}^{-1} \cdot R_{N+1,2}^{-1} \cdots R_{N+1,N}^{-1} \; .
\end{array}
\ee
In fact we have only used the first eq. in (\ref{3.3.21}). The second one is needed if we apply
$RTT$ relations in different manner:
$T_{m} T_{m+1}= R_{m+1,m} T_{m+1} T_{m} R_{m+1,m}^{-1}$.

The relations (\ref{3.3.21}) are deduced from the expressions
(\ref{3.3.20a}) for quantum determinants. Indeed, we have
\be
\lb{3.3.22}
{\det}_{q}(R^{(\pm)}_{N+1})=  {\cal E}_{\langle 12 \dots N} \,
R^{(\pm)}_{1,N+1} \cdots R^{(\pm)}_{N,N+1}  {\cal E}^{12 \dots N\rangle } =
q^{\pm 1} \, I_{N+1} \; ,
\ee
where matrices $R^{(\pm)}$ are representations for elements $T^{i}_{j}$ which were defined in
(\ref{3.1.18}), (\ref{3.1.19}). The last equality in
(\ref{3.3.22}) follow from the fact
that $R^{(+)}$ and $R^{(-)}$ are, respectively, upper and lower triangular block matrices
with diagonal blocks of the form
$$
({R^{(\pm)}}^{i}_{i})^{k}_{l} = \delta^{k}_{l} q^{\pm \delta_{ik}} \; .
$$

Assume that the quantum determinant (\ref{3.3.20a}) is invertible central element. Consider an
extension of the $RTT$ algebra (\ref{glrtt}) by the central element
${\det}_q^{-1}(T)$ which is inverse element for the quantum determinant (\ref{3.3.20a}).
 Then, one can use the ${\cal E}$- tensor
(\ref{3.3.19}), the identity
${\cal E}_{j j_2 \dots j_N} \,  {\cal E}^{i j_2 \dots j_N}
 = \frac{q^N}{[N]_q} \, D^i_j$ (see eq. (\ref{trNA}) below;
 matrix $D$ is defined in (\ref{3.3.14}))
 and inverse element
${\det}_q^{-1}(T)$ to find
an explicit form for the inverse matrix
$T^{-1}$:
\be
\lb{explT-1}
(T^{-1})^i_j = M^{i}_{k} \,  (D^{-1})^k_{j} \, {\det}_q^{-1}(T)
\;\;\;\; \Rightarrow \;\;\;\;
  T^\ell_i \; (T^{-1})^i_j = \delta^\ell_j \; ,
\ee
where $M^{i_1}_{j_1} := q^{-N}[N]_q \,
{\cal E}_{j_1 j_2 \dots j_N} \, T^{j_2}_{i_2} \cdots
 T^{j_N}_{i_N} \, {\cal E}^{i_1 i_2 \dots i_N}$
  are quantum minors of the elements $T^{i_1}_{j_1}$. So, the existence of the inverse
matrix $||(T^{-1})^i_j||$ for the RTT-algebra with R-matrix (\ref{3.3.6}) is equivalent to the
invertibility of the central element ${\det}_q(T)$.
We note that eq. (\ref{3.3.20}) can be written as
\be
\lb{3.3.20i}
{\cal E}_{\langle 12 \dots N}\, T_{N}^{-1} \cdots T_{1}^{-1}
= {\det}_{q}^{-1}(T) \, {\cal E}_{\langle 12 \dots N} \,
 \; , \;\;\;\;\;\;
 T_{N}^{-1} \cdots T_{1}^{-1}  \, {\cal E}^{12 \dots N\rangle } \,
= {\cal E}^{12 \dots N\rangle }  {\det}_{q}^{-1}(T) \; .
\ee


\newtheorem{def11}[def1]{Definition}
\begin{def11} \label{def11}
{\it A Hopf algebra generated by unit element $1$, $N^2$ elements
$T^{i}_{j}$ $(i,j = 1, \dots, N)$ which satisfy relations
(\ref{3.1.1}) with $R$-matrix (\ref{3.3.6}) and element ${\det}_q^{-1}(T)$ is called the algebra
of functions on the linear quantum group $GL_q(N)$
and denoted by $Fun(GL_q(N))$.}
\end{def11}
The structure mappings for the algebra $Fun(GL_q(N))$ are presented in (\ref{3.1.5}), where
elements $(T^{-1})^i_j$ are defined in (\ref{explT-1}).

The algebra $Fun(SL_{q}(N))$ can be obtained from the algebra $Fun(GL_q(N))$ by imposing the
additional condition $det_{q}(T)=1$ and, in accordance with
(\ref{3.3.22}), the matrix representations (\ref{3.1.19}) for
$T^{i}_{j} \in Fun(SL_{q}(N))$ are given by formulas
\be
\lb{LTsl}
\langle L^{+}_{2},T_{1} \rangle = \frac{1}{q^{1/N}} R_{12} \; , \;\;
\langle L^{-}_{2},T_{1} \rangle = q^{1/N} R^{-1}_{21} .
\ee
Conversely, formulas (\ref{3.1.19}), (\ref{LTsl})
can be interpreted as matrix representations of
elements $(L^{\pm})^i_j$ which are generators
(see next subsection {\bf \ref{qalur}})
of the universal enveloping algebras
 $U_q(g\ell(N))$, $U_q(s\ell(N))$.

\vspace{0.2cm}

\noindent
{\bf Remark 1.} The complexification of
the linear quantum groups can be introduced as follows.
We first consider the case of the group
$GL_q(N)$ and assume that $q$ {\it is a real number}.
We have to define an involution $*$-operation,
or simply $*$-involution, (which is the antihomomorphism) on the algebra
$Fun(GL_q(N))$ or, in other words, we must define the conjugated
algebra $Fun(\widetilde{GL_{q}}(N))$ with generators\footnote{We recall that
$(T^{-1})^t \neq (T^{t})^{-1}$ in the case of the quantum matrices (see (\ref{tmin})).}
\be
\lb{Tdag}
\tilde{T} = (T^{\dagger})^{-1} \; , \;\;\;\;\;\;\;\;
 T^{\dagger} := (T^{*})^{t} \;\; \Leftrightarrow \;\;
(T^\dagger)^i_j := (T^j_i)^*  \; ,
\ee
and defining
relations identical to (\ref{3.1.1}):
\be
\lb{3.3.24}
R_{12} \, \tilde{T}_{1} \, \tilde{T}_{2} =
\tilde{T}_{2} \, \tilde{T}_{1} \, R_{12} \;\;\; \Rightarrow \;\;\;
\R_{12} \, \tilde{T}_{1} \, \tilde{T}_{2} =
\tilde{T}_{1} \, \tilde{T}_{2} \, \R_{12}  \; .
\ee
Then we introduce the extended algebra with generators $\{ T^i_j, \; \tilde{T}^k_l \}$
that is the cross (smash) product of the algebras (\ref{3.1.1}) and
(\ref{3.3.24}) with additional cross commutation relations
(see, for example, Refs. \cite{9K}, \cite{9} and \cite{10})
\begin{equation}
\R \, T_1 \, \tilde{T}_2 = \tilde{T}_1 \, T_2 \, \R \; .
\label{3.3.25}
\end{equation}
It is natural to relate this extended double algebra to
$Fun(GL_q(N,\mathbb{C}))$.

The case of $SL_{q}(N,\mathbb{C})$ can be obtained from
$GL_{q}(N,\mathbb{C})$ by imposing two subsidiary conditions on the central elements:
\be
\lb{3.3.26}
det_{q}(T)=1, \;\;\;\;\; det_{q}(\tilde{T})=1 \; .
\ee
The real form $U_{q}(N)$ is extracted
  from $GL_{q}(N,\mathbb{C})$ if we require
\be
\lb{3.3.27}
T = \tilde{T} = (T^{\dagger})^{-1}
\ee
and if, in addition to this,
 we impose the conditions (\ref{3.3.26}), then the group
$SU_q(N)$ is distinguished.

In the case $|q| = 1$, the definition of $*$-involutions
on the linear quantum groups $GL_{q}(N)$ and $SL_{q}(N)$ is
a nontrivial problem that can be solved
only after an imbedding of these quantum groups
into the algebra of functions on their cotangent bundles
(see Remark 2 in next subsection).

\subsubsection{Quantum algebras $U_q(gl(N))$ and $U_q(sl(N))$.
Universal ${\cal R}$ matrix
for $U_q(\mathfrak{g})$\label{qalur}}

The quantum universal enveloping algebras $U_q(gl(N))$ and $U_q(sl(N))$ appear in the $R$-matrix
approach \cite{10} as the algebras with defining relations
(\ref{3.1.20b}), (\ref{3.1.20}). To show this we
consider the upper and lower triangular matrices  $L^+$, $L^-$ in the form (cf. \cite{10}, \cite{DiFre})
\be
\lb{lplu}
L^{+} =
\left(
\begin{array}{cccc}
q^{H_1} & 0  & \dots & 0 \\
0 & q^{H_2} &  \dots & 0 \\
\vdots & \vdots & \ddots &  \vdots \\
0 & 0 & \dots  & q^{H_N}
\end{array}
\right)
\,
\left(
\begin{array}{ccccc}
1 & \lambda f_1 & \lambda f_{1 3} & \dots & *  \\
0 & 1 & \lambda f_2  & \dots  & \dots \\
\vdots & \ddots & \ddots & \ddots & \vdots \\
0 & \dots & 0 &  1 & \lambda f_{N-1} \\
0 & 0 & \dots  & 0 & 1
\end{array}
\right) \; , \;\;\;
\ee
\be
\lb{lmin}
L^{-} =
\left(
\begin{array}{ccccc}
1  & 0 & \dots & \dots & 0 \\
-\lambda e_1 & 1 & 0  & \dots & 0 \\
-\lambda e_{3 1} & -\lambda e_2 & 1 & \dots & 0 \\
\vdots & \vdots & \ddots & \ddots & \vdots \\
\, * & \dots & \dots  & -\lambda e_{N-1} & 1
\end{array}
\right) \,
\left(
\begin{array}{cccc}
q^{-\tilde{H}_1} & 0  & \dots & 0 \\
0 & q^{-\tilde{H}_2} &  \dots & 0 \\
\vdots & \vdots & \ddots &  \vdots \\
0 & 0 & \dots  & q^{-\tilde{H}_N}
\end{array}
\right)
 \; ,
\ee
where $e_\alpha$ and $f_\alpha$ denote respectively
positive and negative root generators of $U_q(sl(N))$.
Here we took into account
 definitions (\ref{0.5}) of matrices $L^{\pm}$
and the convention that the universal $R$-matrix has the form
${\cal R}^{\sf op} \sim \sum_{\alpha}  r_\alpha(H_i)
(f_\alpha \otimes e_\alpha)$
that is in the agreement to the low-triangular
expression (\ref{3.3.6a}) of the $GL_q(N)$ $R$-matrix.
In particular, from (\ref{lplu})
and (\ref{lmin}), we have
\be
\lb{lll}
(L^{+})^i_i = q^{H_i} \; , \;\;\;
(L^{-})^i_i = q^{- \tilde{H}_i} \; , \;\;\;
 (L^{+})^i_{i+1} = \lambda \, q^{H_i} \, f_i
\; , \;\;\; (L^{-})_i^{i+1} = - \lambda \, e_i
\, q^{-\tilde{H}_i} \; .
\ee
For $R$-matrices (\ref{3.3.6a}), (\ref{3.3.6}) the relations (\ref{3.1.20b}),
(\ref{3.1.20}) are
represented in the component form as
\be
\lb{gl1}
(L^{\pm})^i_k \, (L^{\pm})^j_k = q \, (L^{\pm})^j_k \, (L^{\pm})^i_k \; , \;\;\;
(L^{\pm})_i^k \, (L^{\pm})_j^k = q \, (L^{\pm})_j^k \, (L^{\pm})_i^k \; , \;\; (i > j) \; ,
\ee
\be
\lb{gl2}
[(L^{\pm})^{i_1}_{j_1}, \, (L^{\pm})^{i_2}_{j_2}] =
\lambda \, (L^{\pm})^{i_1}_{j_2} \, (L^{\pm})^{i_2}_{j_1}  , \;\;
[ (L^{\pm})^{i_1}_{j_2}, \, (L^{\pm})^{i_2}_{j_1}] = 0, \;\;
(i_1 > i_2, \, j_1 > j_2),
\ee
\be
\lb{gl3}
(L^{+})^i_k \, (L^{-})^j_k = q \, (L^{-})^j_k \, (L^{+})^i_k \; , \;\;\;
(L^{-})_i^k \, (L^{+})_j^k = q \, (L^{+})_j^k \, (L^{-})_i^k \; , \;\; (i < j) \; ,
\ee
\be
\lb{gl4}
\; [(L^{\mp})^{i_1}_{j_1}, \, (L^{\pm})^{i_2}_{j_2}] = 0 \; , \;\; (i_1 > i_2, \; j_1 > j_2)
\; , \;\;\;
[ (L^{+})^{i}_{i}, \, (L^{-})^{i}_{i}] = 0 \; ,
\ee
\be
\lb{gl5}
 [(L^{-})^{i_1}_{j_1}, \, (L^{+})^{i_2}_{j_2}] =
\lambda \, \left( (L^{+})^{i_1}_{j_2} \, (L^{-})^{i_2}_{j_1} -
(L^{-})^{i_1}_{j_2} \, (L^{+})^{i_2}_{j_1} \right) \; , \;\;
(i_1 > i_2, \; j_1 < j_2)  \; ,
\ee
(there are no summation over repeated indices).
We have written only the terms and relations
which survive under the condition that
$(L^+)^i_{\; j} = 0 = (L^{-})^j_{\; i}$, $i > j$.

The substitution of (\ref{lll}) into eqs. (\ref{gl1}) -
(\ref{gl5}) gives the Drinfeld-Jimbo \cite{Jimb1} formulation of  $U_q(gl(N))$. Indeed, from eqs.
(\ref{gl1}),
(\ref{gl3}) and (\ref{gl4}) one can obtain that $q^{H_i - \tilde{H}_i}$ are the
central elements. Thus, the matrices $L^{\pm}$ can be renormalized
(by multiplying them with diagonal matrices) in such a way that
elements $q^{H_i-\tilde{H}_i}$ are fixed as units, i.e.
$H_i=\tilde{H}_i$. Then, from eq. (\ref{gl3}) we find
\be
\lb{dj1}
f_i q^{H_j} = q^{\delta_{j,i} - \delta_{j,i+1}} \,  q^{H_j} f_i \; , \;\;\;
e_i q^{H_j} = q^{\delta_{j,i+1} - \delta_{j,i}} \,  q^{H_j} e_i \; .
\ee
The first eq. in (\ref{gl4}) gives $e_i f_j = f_j e_i$
for $i \neq j$
and taking into account (\ref{gl5}) we derive
\be
\lb{dj2}
e_i f_j - f_j e_i = \delta_{i,j} \frac{q^{H_i - H_{i+1}} - q^{H_{i+1} - H_{i}}}{\lambda} \; .
\ee
The first eq. in (\ref{gl2}) yields a part
 of the Serre relations
\be
\lb{dj3}
e_i e_j = e_j e_i \; , \;\;\; f_i f_j = f_j f_i \; , \;\; (|i-j| \geq 2) \; .
\ee
and gives the expressions
 of the composite roots via the simple roots $\{ e_i, f_j \}$:
\be
\lb{dj4}
\begin{array}{c}
f_{i-1,i+1} = (f_i f_{i-1} - q^{-1} f_{i-1} f_i)  =
 \lambda^{-1} \, q^{-H_{i-1}} \, (L^+)^{i-1}_{i+1} \; , \\ [0.3cm]
e_{i+1,i-1} = (e_{i-1} e_i - q e_i e_{i-1}) =
 - \lambda^{-1} \, (L^-)^{i+1}_{i-1} \, q^{H_{i-1}} \; .
\end{array}
\ee
Using these definitions and eqs. (\ref{gl1}) we deduce another
part of Serre relations
\be
\lb{dj5}
\begin{array}{c}
e_i^2 e_{i\pm 1} - (q + q^{-1}) e_i  e_{i\pm 1} e_i +  e_{i\pm 1} e_i^2 = 0 \;\;
(1 \leq i,i\pm 1 \leq N) \; , \\ [0.3cm]
f_i^2 f_{i\pm 1} - (q + q^{-1}) f_i  f_{i\pm 1} f_i +  f_{i\pm 1} f_i^2 = 0 \;\;
(1 \leq i,i\pm 1 \leq N) \; .
\end{array}
\ee
So, we see that equations (\ref{3.1.20b}), (\ref{3.1.20}) (with the form of
$L^{\pm}$ given in (\ref{lplu}), (\ref{lmin})) not only
yield the commutation relations
(\ref{dj1}), (\ref{dj2}) for the elements of the Chevalley basis, but also
present the Serre relations (\ref{dj3}), (\ref{dj5}) and define the composite root elements (\ref{dj4}) as
the $q$- commutators of the simple root elements.
In this sense the generators $(L^{\pm})^i_j$
(\ref{lplu}), (\ref{lmin}) play the role of a
quantum analog of elements of the Cartan-Weyl basis for
$U_q(gl(N))$, where
$q^{H_k}$, $(L^{+})^i_j$ and $(L^{-})^j_i$ ($i <j$) are
respectively analogs of Cartan elements, negative and positive
root generators. The quantum Casimir operators are given by eqs.
(\ref{3.1.21a}) and (\ref{3.1.22}).
The co-multiplication, antipode and coidentity
in terms of the generators $\{H_i, \, e_i , \, f_i \}$
can be deduced from (\ref{3.1.20a}), (\ref{3.1.20c})
$$
\begin{array}{c}
\Delta(q^{H_i}) = q^{H_i} \otimes q^{H_i} \; , \;\;\;
\Delta(f_i) = 1 \otimes f_i + f_i \otimes q^{H_{i+1} - H_i} \; , \;\;\; \Delta(e_i) =
e_i \otimes 1 + q^{H_{i} - H_{i+1}} \otimes e_i \; , \\ [0.2cm]
S(H_i) = - H_i , \;\;\; S(e_i) = - q^{H_{i+1}-H_i} e_i , \;\;\;
S(f_i) = -  f_i q^{H_i-H_{i+1}} , \;\;\;
\varepsilon(H_{i},e_i,f_i)=0 \; .
\end{array}
$$

Note that $\sum_i H_i$ is a central element in the algebra
$U_q(gl(N))$ and the condition $\sum_i H_i = 0$ reduces
$U_q(gl(N))$ to the algebra $U_q(sl(N))$ with generators
$\{ h_i := H_i - H_{i+1}, \; e_i, \; f_i \}$ subject to the relations
\be
\lb{DriJim}
[q^{h_i} , \, q^{h_j}] = 0 \; , \;\;\;
q^{h_j} \, f_i = q^{- a_{ij}} \, f_i \, q^{h_j}  \; , \;\;\;
q^{h_j} \, e_i = q^{a_{ij}} \, e_i \, q^{h_j}  \; ,
\ee
\be
\lb{DriJim1}
e_i f_j - f_j e_i = \delta_{ij} \frac{q^{d_i h_i} - q^{- d_i h_{i}}}{q^{d_i} - q^{-d_i}}
\; ,
\ee
and Serre relations
\be
\lb{DriJim2}
\sum_{k=0}^{1-a_{ij}} (-1)^k \left[
\begin{array}{c}
1 - a_{ij} \\
k
\end{array}
\right]_{q^{d_i}} (e_i)^k e_j (e_i)^{1 - a_{ij} -k} = 0 \; , \;\;\;
(e_i \rightarrow f_i) \; ,
\ee
where
\be
\lb{binom}
\left[
\begin{array}{c}
n \\
k
\end{array}
\right]_q =
\frac{[n]_q!}{[k]_q! [n-k]_q!} \; , \;\; [k]_q = \frac{q^k - q^{-k}}{q-q^{-1}} \; ,
\;\; [k]_q! := [1]_q [2]_q \cdots [k]_q \; , \;\; [0]_q! := 1 \; ,
\ee
$a_{ij}= 2\delta_{ij} - \delta_{j i+1} - \delta_{i j+1}$ is Cartan matrix
for $sl(N)$, $d_i$ are smallest positive integers (from the set 1,2,3)
such that $d_i a_{ij} \equiv a_{ij}^{\sf sym}$ is symmetric
Cartan matrix (for $sl(N)$ case $d_i =1$).
For the quantum
algebra (\ref{DriJim}) -- (\ref{DriJim2}),
the structure mappings are
 \be
 \lb{strg}
\begin{array}{c}
\Delta(q^{h_i}) = q^{h_i} \otimes q^{h_i} \; , \;\;\;
\Delta(f_i) = 1 \otimes f_i + f_i \otimes q^{- d_i h_i} \; , \;\;\; \Delta(e_i) =
e_i \otimes 1 + q^{d_i h_{i}} \otimes e_i \; , \\ [0.2cm]
S(h_i) = - h_i , \;\;\; S(e_i) = - q^{-d_i h_i} e_i , \;\;\;
S(f_i) = -  f_i q^{d_i  h_i} , \;\;\;
\varepsilon(h_{i},e_i,f_i)=0 \; .
\end{array}
 \ee

\vspace{0.3cm}
\noindent
{\bf Remark 1.} The relations (\ref{DriJim}) - (\ref{DriJim2})
are used for the Drinfeld - Jimbo
\cite{13}, \cite{Jimbo2} formulation of the quantum universal enveloping algebra
$U_q(\mathfrak{g})$ for any simple Lie algebra
 $\mathfrak{g}$ (the elements $e_i,f_i,q^{h_i}$
are related to the simple root $\alpha_i|_{i=1,...,r}$
of the Lie algebra $\mathfrak{g}$ of the rank $r$).
By using this formulation of the quasitriangular
Hopf algebra $U_q(\mathfrak{g})$, one can explicitly construct
the corresponding universal ${\cal R}$ matrix
(the definition via canonical element is given in
(\ref{2.39})). In the case of algebra $U_q(sl_N)$, the explicit multiplicative formula for the
universal ${\cal R}$ matrix has been invented in
\cite{Rosso}. This result was generalized in \cite{KiRe} for the case of
$U_q(\mathfrak{g})$, where $\mathfrak{g}$ is
any semisimple Lie algebra. For the case of quantum Lie superalgebras the universal ${\cal R}$ matrix has
been found in \cite{37}. Finite dimensional representations for the quantum simple Lie algebras
$U_q(\mathfrak{g})$
(\ref{DriJim}) - (\ref{DriJim2}) were considered, e.g.,
 in \cite{Rosso1}.

Here we give (without proof) the explicit
multiplicative formula for the
$U_q(\mathfrak{g})$ universal ${\cal R}$ matrix,
where $\mathfrak{g}$ is
a simple Lie algebra. We give this formula in the form proposed
by Khoroshkin and Tolstoy \cite{37} (see also \cite{KiRe},
\cite{KlSch}). For this we need the notion of the normal ordering
 \cite{37} of the system $\Delta_+$ of
positive roots of Lie algebra $\mathfrak{g}$. We say that the system
$\Delta_+$ is in the normal ordering $\Delta_+^{(n)}$,
 if
 each composite root $\alpha+\beta \in \Delta_+$,
where $\alpha,\beta \in \Delta_+$, has to be placed
  in the ordering between $\alpha$ and $\beta$. It is
 clear that there is an arbitrariness in
 such a normal ordering $\Delta_+^{(n)}$ of positive roots.
 \begin{proposition}\label{the15}
{\it For any quantized Lie algebra $U_q(\mathfrak{g})$ with
  defining relations (\ref{DriJim}) -- (\ref{DriJim2}) and
   for any normal ordering $\Delta_+^{(n)}$ of the positive root
   system $\Delta_+$ of $\mathfrak{g}$,
   the universal ${\cal R}$ matrix such that
   ${\cal R}^{-1} \Delta {\cal R} = \Delta^{\sf op}$,
   where $\Delta$ is the comultiplication
   (\ref{strg}),
   is given by the formula
   \be
   \lb{khtoR}
   {\cal R} = K \cdot
   \overrightarrow{\prod_{\beta \in \Delta_+^{(n)}}} \;
   \exp_{q_\beta}\bigl((q-q^{-1})(e_\beta \otimes f_{\beta}) \bigr) \; ,
   \ee
   $$
   K  := q^{\sum_{ij} d_{ij} h_i \otimes h_j} \; , \;\;\;\;
   \exp_{q}\bigl(x\bigr) :=\sum_{n \geq 0} x^n/(n)_q! \; , \;\;\;\;
   (n)_q := (q^n-1)/(q-1) \; ,
   $$
   where $q_\beta = q^{(\beta,\beta)}$, $d_{ij}$
   is an inverse matrix for the symmetrized
   Cartan matrix $a_{ij}^{\sf sym}=d_i a_{ij}$
   (see definition of $d_i$ after (\ref{binom})) and the
 ordered product runs over the normal ordering $\Delta_+^{(n)}$
    of the positive roots.}
 \end{proposition}
 For the case of $U_q(s\ell(2))$ algebra
 $$
 [e,f] = \frac{q^h-q^{-h}}{q-q^{-1}} \; , \;\;\;\;
 q^{h} \, f = q^{- 2} \, f \, q^{h}\; , \;\;\;\;
 q^{h} \, e = q^{2} \, e \, q^{h} \; ,
 $$
 the formula (\ref{khtoR}) is simplified
 $$
{\cal R} = q^{\frac{1}{2} h \otimes h} \cdot
\exp_{q^2}\Bigl((q-q^{-1})(e \otimes f) \Bigr) \; .
$$
 Finally we note that in the paper \cite{37} the authors
 used another comultiplication $\Delta'$
 for $U_q(\mathfrak{g})$
  \be
 \lb{strg2}
\Delta'(q^{h_i}) = q^{h_i} \otimes q^{h_i} \; , \;\;\;
\Delta'(f_i) = 1 \otimes f_i + f_i \otimes q^{d_i h_i} \; , \;\;\; \Delta'(e_i) =
e_i \otimes 1 + q^{-d_i h_{i}} \otimes e_i \; ,
 \ee
  which is
 related to the comultiplication (\ref{strg})
  by twisting
 $\Delta' = K^{-1} \, P_{12} \, \Delta \, P_{12} \, K
 \equiv K^{-1} \, \Delta^{\sf op} \, K$.
 This explains why our formula (\ref{khtoR}) differs
 by twisting from the formula for ${\cal R}$
 given in \cite{37}.

\vspace{0.3cm}
\noindent
{\bf Remark 2.}
The $*$-involution on the algebra $U_q(sl(N))$ (\ref{DriJim}) -- (\ref{DriJim2})
for real $q$ is defined if we note that the algebra with generators $T,\widetilde{T}$
(\ref{3.3.1}), (\ref{3.3.24}), (\ref{3.3.25}) coincides to the $L^{\pm}$-algebra
(\ref{3.1.20b}), (\ref{3.1.20}) after an
identification: $L^{-} = T^{-1}$, $L^{+} = \widetilde{T}^{-1}$.
Then, according to
(\ref{Tdag}) we require $(L^{+})^\dagger = (L^-)^{-1}$. In terms of the Chevalley generators \ref{lplu}, (\ref{lmin})
it means that
\be
\lb{starU}
h_i^* = h_i \; , \;\;\; f_i^* = q \, q^{-h_i} \, e_i \; , \;\;\;
e_i^* = q^{-1} \, f_i \, q^{-h_i} \; .
\ee
One can directly check that eqs. (\ref{DriJim}) -- (\ref{DriJim2}) respect
the antihomomorphism (\ref{starU}).

In the case $|q| = 1$, the definition of $*$-involutions
on the linear quantum groups and algebras is
a nontrivial problem that can be solved \cite{23}
 only after extension of the algebra
of functions on the quantum groups to the algebra of functions on their
cotangent bundles, i.e., to the algebra which is a Heisenberg
double of $Fun(GL_q(N))$ and $U_q(gl(N))$
with cross-multiplication rules (\ref{lhd}) - (\ref{tangb}).


\subsection{\bf \em Hecke type $R$-matrices.
Related quantum matrix algebras\label{hRMa}}
\setcounter{equation}0

The material in this Section is based in part on
the results of papers \cite{IsBonn}, \cite{Gur},  \cite{HIOPT}.

\subsubsection{Definitions. (Anti)symmetrizers for
Hecke type $R$-matrices}

\newtheorem{def11d}[def1]{Definition}
\begin{def11d} \label{def11d}
{\it The Yang-Baxter $R$-matrices which obey
(\ref{3.1.2i}) and Hecke condition (\ref{3.3.7}),
 (\ref{3.3.7aa}) are called {\it Hecke type $R$-matrices}. }
\end{def11d}
First of all we note that
the $GL_q(N)$ matrices (\ref{3.3.6a}), (\ref{3.3.6}) are
examples of Hecke type $R$-matrices
since they satisfy the
 Hecke condition (\ref{3.3.7}), (\ref{3.3.7aa}). We also
note that if $\R[q]$ satisfies Hecke condition (\ref{3.3.7aa}), then
 $\R[-q^{-1}]$ and  $-\R[q^{-1}]$ are also satisfy (\ref{3.3.7aa}). In this subsection we present
some general facts about Hecke type $R$-matrices and related quantum algebras.

The antisymmetrizers $A_{1 \to m}$ (\ref{Aeps1}) can be explicitly constructed in terms of the
Hecke type $R$-matrices by using the following inductive procedure \cite{Jimb1}
(the same procedure was used in \cite{Gur}; see also Subsec. {\bf \ref{jmel}} below):
\be
\lb{antirr}
A_{1 \to k} =    A_{2 \to k}
\left( \frac{\R_{1}( q^{k-1} )}{[k]_q} \right) A_{2 \to k} =
A_{1 \to k-1}
\left( \frac{\R_{k-1}( q^{k-1} )}{[k]_q} \right) A_{1 \to k-1} =
\ee
$$
= \frac{1}{[k]_q!} \, A_{1 \to k-1} \,
\R_{k-1}\left( q^{k-1} \right) \,
\R_{k-2}\left( q^{k-2} \right)\dots
 \R_{2}\left( q^{2} \right) \, \R_{1}\left( q \right) \; , \;\;\; (k = 2,3, \dots N) \; ,
$$
where $A_{1 \to 1}=1$, $\R(x) = (x^{-1}\R - x\R^{-1})/\lambda$ -- Baxterized
$R$-matrix (see below Subsect. {\bf \ref{baxtel}}),
$\R$ is a Hecke type $R$-matrix,
$[k]_q = (q^k - q^{-k})/\lambda$ and as usual
\be
\lb{3.3.23}
\R_{k} =I^{\otimes (k-1)} \otimes \R \otimes I^{\otimes (N-k)}
\in Mat(N)^{\otimes (N+1)} \; .
\ee


\newtheorem{def12}[def1]{Definition}
\begin{def12} \label{def12}
{\it We say that the Hecke type $R$-matrix is of the height $N$, if $A_{1 \to M} = 0$ $\forall M
>N$ and ${\sf rank}(A_{1 \to N}) = 1$.}
\end{def12}

\noindent
Note that, for the $GL_q(N)$ type $R$-matrix (\ref{3.3.6a}), (\ref{3.3.6}), the operator
$A_{1 \to N+1} = 0$ and
$A_{1 \to N}$ is the highest $q$-antisymmetrizer in the sequence of the
antisymmetrizers (\ref{antirr}). Moreover we have
${\sf rank}(A_{1 \to N}) = 1$ in this case.
The latter can easily be
understood by considering the fermionic quantum hyperplane (\ref{3.3.13}).
Since the operators $A_{1 \to k}$ (\ref{antirr}) satisfy
(cf. (\ref{3.3.19}))
\be
\lb{raar}
\R_{j} \, A_{1 \to k} = A_{1 \to k} \, \R_{j}  =
- q^{-1} \, A_{1 \to k} \;\; (j = 1, \dots, k-1) \; ,
\ee
they are symmetry operators for
$k$-th order monomials $x^{i_1} \dots x^{i_k}$ in the
q-fermionic algebra (\ref{3.3.13}).
In view of the explicit relations (\ref{3.3.13}) one can conclude that
there is only one independent monomial of the order $N$ and all monomials
$x^{i_1} \dots x^{i_k}$, for $k >  N$, are equal to zero. This statement is equivalent to
the conditions ${\sf rank}(A_{1 \to N}) = 1$
 and $A_{1 \to N+1} = 0$.

In view of the definition (\ref{antirr}), the condition $A_{1 \to N+1} = 0$
leads to (for arbitrary Hecke $\R$-matrix):
\be
\lb{3.3.21ab}
A_{1 \to N} \, \R^{\pm 1}_N \, A_{1 \to N} = \frac{q^{\pm N}}{[N]_q}
\, A_{1 \to N} \, I_{N+1} \, , \;\;\;\;
A_{2 \to N+1} \, \R^{\pm 1}_1 \, A_{2 \to N+1} = \frac{q^{\pm N}}{[N]_q}
\, A_{2 \to N+1} \, I_{1} \, .
\ee
In the case of skew-invertible Hecke $\R$-matrices,
by applying (\ref{skew}) to eqs. (\ref{3.3.21ab}), we obtain
\be
\lb{gov1}
A_{1 \to N} \, P_{_N\, _{N+1}} \, A_{1 \to N} =
\frac{q^{N}}{[N]_q} \, A_{1 \to N} \, Q_{_{N+1}} \; , \;\;\;
A_{2 \to N+1} \, P_{_1 \, _2} \, A_{2 \to N+1} =
\frac{q^N}{[N]_q} D_1 \, A_{2 \to N+1} \; ,
\ee
and for completely invertible $\R$-matrices we have in addition
$\overline{Q} = q^{2N} \, Q$, $\overline{D} =  q^{2N} \, D$.
Acting respectively by $Tr_{_{N+1}}$ and $Tr_1$
to the first and second eq. in (\ref{gov1}),
we deduce (cf. (\ref{3.3.15}))
\be
\lb{trDtrQ}
Tr(Q) = Tr(D) = q^{-N} [N]_q \;\; \Rightarrow \;\;
Tr(\overline{Q}) = Tr(\overline{D}) = q^{N} [N]_q \; ,
\ee
while applying $Tr_{(1 \dots N)}$ and $Tr_{(2 \dots N+1)}$ to the
same eqs. we obtain \cite{IsBonn}
 \be
 \lb{trNA}
Tr_{(1 \dots N-1)} A_{1 \to N} =  \frac{{\sf rank}(A_{1 \to N})}{Tr(Q)} Q_N \; , \;\;\;\;\;\;
Tr_{(2 \dots N)} A_{1 \to N} =  \frac{{\sf rank}(A_{1 \to N})}{Tr(D)} D_1 \; .
 \ee

On the other hand, applying quantum traces $Tr_{D(N-k+1 \dots N)}$ and
$Tr_{Q(1 \dots k)}$ to the antisymmetrizers $A_{1 \dots N}$
we deduce \cite{IsBonn} ($0 \leq k \leq N$)
\be
\lb{detD3}
\left[ \!\!
\begin{array}{c}
N \\[-0.1cm]
k
\end{array} \!\! \right]_q
\, Tr_{D(k+1 \dots N)} \left( A_{1 \dots N}  \right)  = q^{(k-N) \, N} \, A_{1 \dots k} \; , \;\;\;\;\;\;\;\;
\left. A_{1 \dots k} \right|_{k=0} \; :=1  \; ,
\ee
\be
\lb{detQ3}
\left[ \!\!
\begin{array}{c}
N \\[-0.1cm]
k
\end{array} \!\! \right]_q
Tr_{Q(1 \dots k)} \left( A_{1 \dots N}  \right)  = q^{-k \, N} \, A_{k+1 \dots N}  \;  , \;\;\;\;\;\;\;\;
\left. A_{k+1 \dots N} \right|_{k=N} \; :=1 \; ,
\ee
where q-binomial coefficients are defined in
(\ref{binom}) and we have used eqs. (\ref{antirr}) and identities
\be
\lb{trqRx}
Tr_{_{D(k+1)}} \R_k(x) = Tr_{_{Q(k-1)}} \R_{k-1}(x) =
\left( \frac{x^{-1} - x}{\lambda} + x \, Tr(D) \right)  I_k
= \frac{x^{-1} - x q^{-2 \, N}}{\lambda} \, I_k  \; ,
\ee
which follow from (\ref{qtrs}), (\ref{trDtrQ}). In view of eqs. (\ref{sk10}) matrices $D$ and $Q$
can be considered as one dimensional representations of $RTT$ algebra (\ref{3.1.1}): $\rho_D(T^i_j)
= D^i_j$,
$\rho_Q(T^i_j) = Q^i_j$. Thus, we have
\be
\lb{detD}
A_{1 \dots N} D_1 D_2 \dots D_N = det_q(D) \, A_{1 \dots N} \; , \;\;\;
A_{1 \dots N} Q_1 Q_2 \dots Q_N = det_q(Q) \, A_{1 \dots N} \; .
\ee
and taking $k=0$ in (\ref{detD3}) and $k=N$ in (\ref{detQ3}) we obtain
\be
\lb{detD1}
{\det}_q(D) = q^{-N^2} \; , \;\;\;
{\det}_q(Q) = q^{-N^2} \; .
\ee

For the Hecke type $R$-matrix
one can construct (in addition to the q- antisymmetrizer $A_{1 \to k}$ (\ref{antirr}))
the q-symmetrizer $S_{1 \to k}$:
\be
\lb{symrr}
S_{1 \to k} =    S_{2 \to k}
\left( \frac{\R_{1}( q^{1-k} )}{[k]_q} \right) S_{2 \to k} =
S_{1 \to k-1}
\left( \frac{\R_{k-1}( q^{1-k} )}{[k]_q} \right) S_{1 \to k-1} \; ,
\ee
(see also Sec. {\bf \ref{gabg}} below). Using identities (\ref{antirr}), (\ref{symrr})
and  (\ref{trqRx}) one can calculate
q-ranks for the projectors $A_{1 \to k}$ and $S_{1 \to k}$:
$$
{\sf rank}_q (A_{1 \to k}) : =
Tr_{D(1 \dots k)} A_{1 \to k} = \frac{(q^{k-1} Tr(D) - [k-1]_q)}{[k]_q} \,
Tr_{D(1 \dots k-1)} A_{1 \to k-1} = \dots =
$$
\be
\lb{trA}
= \frac{1}{[k]_q !} \, \prod_{m=1}^{k} \left( q^{m-1} Tr(D) - [m-1]_q \right) \; ,
\ee
and analogously
\be
\lb{trS}
{\sf rank}_q (S_{1 \to k}) : =
Tr_{D(1 \dots k)} S_{1 \to k}
= \frac{1}{[k]_q !} \, \prod_{m=1}^{k} \left( q^{1-m} Tr(D) + [m-1]_q \right) \; .
\ee

By substituting (\ref{trDtrQ}) into (\ref{trA}), (\ref{trS})
we deduce for the Hecke type $R$-matrix (of the height $N$)
the following "q-dimensions" of the
antisymmetrizer and symmetrizer:
$$
Tr_{D(1 \dots k)} A_{1 \to k} = q^{-k N} \,
\left[ \!\!
\begin{array}{c}
N \\[-0.1cm]
k
\end{array} \!\! \right]_q \;\; (k \leq N) \; , \;\;\;
Tr_{D(1 \dots k)} A_{1 \to k} =  0 \;\; (k > N) \; ,
$$
$$
Tr_{D(1 \dots k)} S_{1 \to k} = q^{-k N} \,
\left[ \!\!
\begin{array}{c}
N +k -1 \\[-0.1cm]
k
\end{array} \!\! \right]_q \; .
$$

The general formula for q-dimensions
of any Young q-symmetrizer (related to any Young diagram),
which is rational function
 of the Hecke type $R$-matrices\footnote{These symmetrizers
 are images of the idempotents of the
 Hecke algebra; see Subsections {\bf \ref{jmel2}}
 and {\bf \ref{idbaxt}} below.}, is known and
can be found in \cite{W1}, \cite{IsOg3} (see
Subsect. {\bf \ref{qdlink}} below).

\vspace{0.3cm}


Sometimes it is convenient to have Eqs.
(\ref{3.3.21}) not only for $GL_q(N)$ type $R$-matrices,
but in a more general form
which is valid for any Hecke $R$-matrix such that
$A_{1 \to N+1}=0$.
For this we consider identity (see, e.g., \cite{HIOPT}):
\be
\lb{haha}
A_{2 \to N+1} \,
\R_{1}^{\pm 1} \cdot \R_{2}^{\pm 1} \cdots \R_{N}^{\pm 1}
  = (-1)^{N-1} \, q^{\pm 1} \, [N]_q \, A_{2 \to N+1} \, A_{1 \to N} \; .
\ee
(we demonstrate a connection of (\ref{3.3.21}) with (\ref{haha}) below).
The mirror counterpart of the relations (\ref{haha}) are also valid
\be
\lb{hahad}
\R_{N}^{\pm 1} \cdots \R_{2}^{\pm 1} \cdot \R_{1}^{\pm 1} \, A_{2 \to N+1}
  = (-1)^{N-1} \, q^{\pm 1} \, [N]_q \, A_{1 \to N} \, A_{2 \to N+1} \; .
\ee
Eqs. (\ref{haha}) and (\ref{hahad}) can be readily
deduced from equation
\be
\lb{hahac}
A_{2 \to N+1} \, \R^{\pm 1}_1 \dots  \R^{\pm 1}_N =
 \R^{\pm 1}_1 \dots  \R^{\pm 1}_N  \, A_{1 \to N} \; ,
\ee
which is obtained from the fact that antisymmetrizers
are expressed in terms of $R$-matrices
(\ref{antirr}) and by using identities
$$
\R_k \; (\R^{\pm 1}_1 \dots  \R^{\pm 1}_N) =
 (\R^{\pm 1}_1 \dots  \R^{\pm 1}_N)  \; \R_{k-1} \; ,
 \;\;\;\;\;\; \forall \; k =2,...,N \; ,
$$
followed from the braid relations (\ref{3.1.3}), (\ref{3.1.4}). Acting on (\ref{hahac}) by $A_{1 \to N}$  from the left, and
making use of eqs. (\ref{raar}), (\ref{3.3.21ab})
 and again eq. (\ref{hahac}),
 we deduce (\ref{hahad}). Eq. (\ref{haha}) can be proved in the same way. Multiplying identities (\ref{haha}) and
 (\ref{hahad}) by $A_{2 \to N}$ respectively
 from the right and left, we find
\be
\lb{hahaf}
A_{1 \to N} \, A_{2 \to N} \, A_{1 \to N} = [N]_q^{-2} \, A_{1 \to N} \; , \;\;\;
A_{2 \to N} \, A_{1 \to N} \, A_{2 \to N} = [N]_q^{-2} \, A_{2 \to N} \; ,
\ee
and then multiplying (\ref{hahad}) by
 (\ref{haha}) from the left and (\ref{haha}) by
 (\ref{hahad}) from the right we obtain equations:
$$
A_{1 \to N} ( \R_N \dots \R_2 \R^{2}_1 \R_2 \dots  \R_N - q^2) = 0 \; , \;\;\; A_{2 \to N+1} ( \R_1
\R_2 \dots \R^{2}_N \dots \R_2 \R_1 - q^2) = 0 \; ,
$$
which followed from (\ref{3.3.21ab})
and are equivalent to $A_{1 \to N+1} =0$
(see (\ref{sant03}) below).

The identity (\ref{haha})
is valid for any Hecke $R$-matrix of the height $N$
and can be transformed into eq. (\ref{3.3.21}).
Indeed, for the case when ${\sf rank}(A_{1 \to N})=1$ and, thus,
$A_{1 \to N}$ is given by the first eq. in
(\ref{Aeps}), one can
act on (\ref{haha}) by ${\cal E}_{\langle 2 \dots N+1}$
from the left and as a result the counterpart of (\ref{3.3.21}) is obtained
\be
\lb{haha'}
 {\cal E}_{\langle 23 \dots N+1} \,
\R_{1}^{\pm 1} \cdot \R_{2}^{\pm 1} \cdots \R_{N}^{\pm 1}
  = q^{\pm 1} \, {\sf N}_{\langle N+1}^{ \; 1\rangle} \, {\cal E}_{\langle 12 \dots N} \; .
\ee
Here we have introduced the matrix:
\be
\lb{matN}
{\sf N}_{\langle N+1}^{ \; 1 \rangle} : = (-1)^{N-1} \, [N]_q  \,
 {\cal E}_{\langle 2 \dots N+1} \, {\cal E}^{1 \dots N\rangle} \; ,
\ee
which, for the case of
$GL_q(N)$ $R$-matrix (\ref{3.3.6a}), is equal to the unit matrix ${\sf N}^i_j = \delta^i_j$
(cf. (\ref{3.3.21})).
Analogously by acting of ${\cal E}^{2 \dots N+1 \rangle}$ on (\ref{hahad})
from the right we deduce
\be
\lb{haha''}
\R_{N}^{\pm 1} \cdots \R_{2}^{\pm 1} \cdot \R_{1}^{\pm 1} \,
 {\cal E}^{23 \dots N+1 \rangle}
  = q^{\pm 1}\, {\cal E}^{ 12 \dots N \rangle}
\, ({\sf N}^{-1})_{\langle 1}^{ \; N+ 1\rangle}  \; ,
\ee
where matrix
\be
\lb{matN1}
({\sf N}^{-1})_{\langle 1}^{ \; N+1 \rangle}
: = (-1)^{N-1} \, [N]_q  \,
 {\cal E}_{\langle 1 \dots N} \,
 {\cal E}^{2 \dots N+1 \rangle} \; ,
\ee
is inverse to the matrix (\ref{matN}) in view of (\ref{hahaf}).

\subsubsection{Quantum determinants for $RTT$
 and $RLRL$ algebras}

For the $RTT$ algebra defined by
the Hecke $R$-matrix of the height $N$
(Definition {\bf \em \ref{def12}}), one can
introduce a generalization of the $GL_q(N)$ $q$-determinant
by using the same formulas
(\ref{3.3.18}), (\ref{3.3.20a}):
 \be
\lb{detheck}
\begin{array}{c}
{\det}_{q}(T) \, {\cal E}_{\langle 12 \dots N}
= {\cal E}_{\langle 12 \dots N} \,
T_{1} \cdot T_{2} \cdots T_{N} \; , \;\;\;\;\;\;
 {\cal E}^{12 \dots N\rangle } \, {\det}_{q}(T)
= T_{1} \cdot T_{2} \cdots T_{N} \,
 {\cal E}^{12 \dots N\rangle } \; , \\ [0.2cm]
 {\det}_{q}(T) = {\cal E}_{\langle 12 \dots N}
\left( T_{1} T_{2} \cdots T_{N} \right)
{\cal E}^{12 \dots N\rangle } =
 Tr_{12 \dots N} (A_{1 \to N} \, T_1 T_2 \cdots T_N) \; .
 \end{array}
\ee
In the case when matrix ${\sf N}^i_j$ is not proportional
 to the unit matrix the chain of relations
(\ref{3.3.21a}) gives \cite{Gur}, \cite{HIOPT}:
\be
\lb{nonc}
\begin{array}{c}
{\cal E}_{\langle 1 \dots N} \, {\det}_{q}(T) \, T_{N+1} =
{\cal E}_{\langle 1 \dots N} \, T_1 \dots T_N \, T_{N+1} = \\
= {\cal E}_{\langle 1 \dots N} \, \R_N^{-1} \dots \R_1^{-1} \, T_1 \dots T_N \, T_{N+1}
\, \R_1 \dots \R_N  = \\
= q^{-1} (N^{-1})^{N+1 \rangle}_{\langle 1} \, T_1 \,
{\cal E}_{\langle 2 \dots N+1} \, T_2 \dots T_{N+1}
\, \R_1 \dots \R_N  = \\
= q^{-1} ({\sf N}^{-1})^{N+1 \rangle}_{\langle 1} \, T_1 \, {\det}_q(T) \,
{\cal E}_{\langle 2 \dots N+1} \, \R_1 \dots \R_N  =
({\sf N}^{-1} \, T {\sf N})_{N+1} \,
 {\det}_q(T) {\cal E}_{\langle 1 \dots N}  \; ,
\end{array}
\ee
where an explicit form of ${\sf N}^{-1}$ can be extracted from (\ref{matN1}).
It means, that for $RTT$ algebras defined by
general Hecke type $R$-matrix of the limited height,
the element ${\det}_q(T)$ is not necessary central.
However, let the noncentral
element ${\det}_q(T)$ be invertible. In this case
one can also define
the inverse matrix $T^{-1}$ (cf. (\ref{explT-1}); see also
eq. (6.16) in \cite{HIOPT})
$$
(T^{-1})^{1 \rangle}_{\langle N+1} = \frac{(-1)^{N-1}}{[N]_q} \,
 {\sf N}_1^{-1} \, {\det}_q^{-1}(T) \,  {\cal E}_{\langle 2 \dots N+1} \,
 T_2 \cdots T_N \, {\cal E}^{1 \dots N\rangle} \;\;
 \Rightarrow \;\;  (T^{-1})^{1 \rangle}_{\langle N+1} \, T_{N+1} =
 I^{1 \rangle}_{\langle N+1} \; ,
$$
where the matrices ${\sf N}$ and ${\sf N}^{-1}$ are
 given in (\ref{matN}) and (\ref{haha''}).

\vspace{0.2cm}

In the case of the
Hecke type $R$-matrix of the height $N$,
 the structures (\ref{tangb}) and (\ref{crtL})
  of the cross-products
(doubles) for the $RTT$ and reflection equation algebras
 \begin{eqnarray}
\R_{12} T_1 T_2 = T_1  T_2 \R_{12}  \; , \;\;\;\;
 T_1 L_2 = \R_{12} L_1 \R_{12}^{\, \pm 1} T_1 \; , \;\;\;\; \R_{12}  L_1 \R_{12}  L_2 =
L_1 \R_{12} L_1 \R_{12} \; , \lb{tangb0} \\ [0.2cm]
\R_{12} T_1 T_2 = T_1  T_2 \R_{12} \; , \;\;\;\;
\overline{L}_1 \, T_2 = T_2 \R_{12}
 \overline{L}_2 \R_{12}^{\, \pm 1}  \; , \;\;\;\;
 \R_{12}  \overline{L}_2 \, \R_{12}  \overline{L}_2 =
 \overline{L}_2 \R_{12} \overline{L}_2 \R_{12} \; ,
 \lb{tangb2}
 \end{eqnarray}
 help us \cite{Isaev1}, \cite{IsBonn},
 \cite{18''} to introduce
the notion of the quantum determinants $Det_q(L)$, $Det_q(\overline{L})$ for the
corresponding reflection equation algebras (\ref{3.1.23a}) and (\ref{3.1.23b}) with generators $L$ and $\overline{L}$.
It can be done by using
the definition (\ref{detheck}) of
${\det}_q(T)$ for the $RTT$ algebras
with the Hecke type $R$-matrix of the height $N$.
In view of the automorphism (\ref{discr}), the
quantum matrix $(L\, T)$
satisfies the same $RTT$ relation (\ref{3.1.1}) and, thus, one can
consider the same quantum determinant ${\det}_q(.)$
for the quantum matrix $(L\, T)$ as for the matrix $T$.
This determinant is divisible from the right
by ${\det}_q(T)$ and the quotient
depends on the matrix $L$ only. This quotient is called the quantum
determinant for the reflection equation algebra (\ref{3.1.23a}).
 We consider only the case of the double with
  structure (\ref{tangb}) with upper signs in
  (\ref{tangb0}), (\ref{tangb2})
  (the case with lower signs is considered in the same way;
  see (\ref{crtL2a}))
\be
\lb{Ldet}
\begin{array}{c}
\displaystyle
Det_q(L) := {\det}_q(L\, T) \frac{1}{{\det}_q(T)} =
{\cal E}_{\langle 1 \dots N} ( L_{1} T_1 \, L_{2} T_2
\dots L_{N} T_N )
 {\cal E}^{1 \dots N \rangle} \, \frac{1}{{\det}_q(T)} = \\ [0.3cm]
 \displaystyle
= {\cal E}_{\langle 1 \dots N} ( L_{\widetilde{1}} \, L_{\widetilde{2}} \dots L_{\widetilde{N}})
T_1 \dots T_N \, {\cal E}^{1 \dots N \rangle} \, \frac{1}{{\det}_q(T)}
= {\cal E}_{\langle 1 \dots N} ( L_{\widetilde{1}} \, L_{\widetilde{2}} \dots L_{\widetilde{N}}) \,
 {\cal E}^{1 \dots N \rangle}  = \\ [0.3cm]
 =  {\cal E}_{\langle 1 \dots N}
( L_{\widetilde{N}} \dots  L_{\widetilde{2}} \, L_{\widetilde{1}})
 {\cal E}^{1 \dots N \rangle} =
Tr_{1 \dots N} \left( A_{1 \dots N} \, L_{\widetilde{1}} \, L_{\widetilde{2}}
\dots L_{\widetilde{N}} \right) \; ,
\end{array}
\ee
where
\be
\lb{Ldett}
 L_{\widetilde{k+1}} = \R_k L_{\widetilde{k}} \R_k \; , \;\;\;\;\;
L_{\widetilde{1}} := L_1
\ee
 are operators that form
 a commutative set
  $[L_{\widetilde{k}}, \, L_{\widetilde{\ell}}]=0$.
  The definition (\ref{Ldet}) is generalized as following
  \be
\lb{LdetE}
 L_{\widetilde{1}} \, L_{\widetilde{2}}
\dots L_{\widetilde{N}} \, A_{1 \dots N} =
A_{1 \dots N} \, Det_q(L) \; .
\ee

For the second type algebra (\ref{3.1.23b}),
(\ref{tangb}) (the algebra (\ref{tangb2})
with upper sign) the definition is
analogous:
$$
\begin{array}{c}
Det_q(\overline{L}) = det^{-1}_q(T) \, {\det}_q(T \overline{L}) =
det^{-1}_q(T) \, {\cal E}_{\langle 1 \dots N} T_1 \dots T_N \,
( \overline{L}_{\widetilde{1}} \, \overline{L}_{\widetilde{2}} \dots \overline{L}_{\widetilde{N}})
{\cal E}^{1 \dots N \rangle} = \\ [0.3cm]
= {\cal E}_{\langle 1 \dots N}
( \overline{L}_{\widetilde{N}} \dots  \overline{L}_{\widetilde{2}} \, \overline{L}_{\widetilde{1}})
 {\cal E}^{1 \dots N \rangle} =
Tr_{1 \dots N} \left( A_{1 \dots N} \, \overline{L}_{\widetilde{1}} \, \overline{L}_{\widetilde{2}}
\dots \overline{L}_{\widetilde{N}} \right) \; ,
 \end{array}
$$
where $\overline{L}_{\widetilde{k}}
= \R_k \, \overline{L}_{\widetilde{k+1}} \, \R_k$ ,
$\;\;\overline{L}_{\widetilde{N}} := \overline{L}_{N}$
and $[\overline{L}_{\widetilde{k}}, \,
\overline{L}_{\widetilde{\ell}}]=0$.

Below, we restrict ourselves to considering only the case
of the left reflection equation algebra
(\ref{3.1.23a}) (and the double (\ref{tangb0})
with upper sign), since the case of the
right algebra (\ref{3.1.23b}) is
investigated analogously.
An interesting property of the determinant $Det_q(L)$
(followed from discrete evolutions
(\ref{discr}), (\ref{discr2}))
is of its multiplicativity \cite{IsBonn}
\be
\lb{Ldet1}
 \begin{array}{c}
 \displaystyle
\!\!\! Det'_q(L^{n+1})  :=
{\det}_q(L^{n+1}\, T) \frac{1}{{\det}_q(T)}
= {\det}_q(L^{n+1}\, T) \frac{1}{{\det}_q(L \, T)}
{\det}_q(L \, T) \frac{1}{{\det}_q(T)} = \\ [0.4cm]
= \bigl( {\cal E}_{\langle 1 \dots N}  (L_{\widetilde{1}})^n \,
(L_{\widetilde{2}})^n \dots (L_{\widetilde{N}})^n \,
 {\cal E}^{1 \dots N \rangle}  \bigr) \; Det_q(L) =
Det'_q(L^{n}) \, Det_q(L)  =\left( Det_q(L)  \right)^n \; ,
\end{array}
\ee
where, for $n=1$, we have $Det'_q(L^n) = Det_q(L)$.

In view of the automorphism (\ref{discr1}) for $n=1$,
such that $T \to (L+x)\, T$,
one can define (in the same way as in (\ref{Ldet}))
a quantum determinant $Det_q(L;x)$ \cite{Isaev1}:
\be
\lb{Ldet2}
Det_q(L;x) := {\det}_q( (L+x)\, T) \frac{1}{{\det}_q(T)}
= Tr_{1 \dots N} \left( A_{1 \dots N} (L_{\widetilde{1}} +x) \dots (L_{\widetilde{N}}+x)
\right) \; ,
\ee
where $x\in \mathbb{C}$ is a parameter and
$L_{\widetilde{k}}$ are given in (\ref{Ldett}). Thus,
 we introduce
the characteristic polynomial for the quantum matrix $L$.
Here we prefer to use the notation $Det_q(L;x)$ instead of $Det_q(L+x)$,
since the dependence on $(L+x)$ seems to be broken in view
of the last expression of (\ref{Ldet2}).
Taking into account (\ref{3.3.20}), the determinant
(\ref{Ldet2}) can be also derived as follows
$$
 {\cal E}^{1 \dots N \rangle} \, Det_q(L;x) =  (L_{1} +x)T_1\,
(L_{2}+x)T_2 \dots (L_{N}+x)T_N \, {\cal E}^{1 \dots N \rangle} \; \frac{1}{{\det}_q(T)} =
$$
\be
\lb{Ldet55}
= \left( (L_{\widetilde{1}} +x)\,
(L_{\widetilde{2}}+x) \dots (L_{\widetilde{N}}+x) \right)
 {\cal E}^{1 \dots N \rangle} \; .
\ee
The expansion of (\ref{Ldet2}) over the parameter $x$ gives
\be
\lb{Ldet3}
Det_q(L;x) = \sum_{k=0}^N \, x^k \, a_{N-k}(L) \; .
\ee
Here $a_0(L) =1$, $a_{N}(L) = Det_q(L)$,
\be
\lb{Ldets}
a_{m}(L) =  Tr_{1 \dots N}
 \Bigl(  A_{1 \dots N} \!\!\!\!
\sum_{_{1 \leq k_1 < \dots < k_m \leq N}} \!\!\!\!
L_{\widetilde{k_1}} L_{\widetilde{k_2}} \dots L_{\widetilde{k}_m} \Bigr)
 =   \alpha^{(m)}_N \; Tr_{1 \dots N}
 \left(  A_{1 \dots N} L_{\widetilde{1}} \, L_{\widetilde{2}}
\cdots L_{\widetilde{m}} \right)
 \; ,
\ee
where in the second equality we applied identities
(for all $k_1 < \dots < k_m$)
$$
Tr_{1 \dots N} \Bigl(  A_{1 \dots N}
L_{\widetilde{k_1}} L_{\widetilde{k_2}} \dots L_{\widetilde{k}_m} \Bigr) =
q^{- 2(k_1 + k_2 + \dots k_m)+m(m+1)} Tr_{1 \dots N}
 \left(  A_{1 \dots N} L_{\widetilde{1}} \, L_{\widetilde{2}}
\cdots L_{\widetilde{m}} \right) \; ,
$$
and \cite{IsBonn}
\be
\lb{Ldeta}
\alpha^{(m)}_N = q^{m(m+1)} \,
 \sum_{_{1 \leq k_1 < \dots < k_m \leq N}} \!\!\!\!
\, q^{- 2(k_1 + k_2 + \dots k_m)}
= q^{m(m-N)} \, \left[ \!\!
\begin{array}{c}
N \\[-0.2cm]
m
\end{array} \!\! \right]_q \; ,
\ee
(q-binomial coefficients
$\left[ \!\!
\begin{array}{c}
N \\[-0.2cm]
m
\end{array} \!\! \right]_q$
were introduced in (\ref{binom})).
The sums in (\ref{Ldeta}) are readily calculated by means
 of their generating function
$$
\alpha(t) = \sum_{m=0}^N t^{N-m} q^{-m(m+1)} \alpha^{(m)}_N =
\prod_{m=1}^{N} (t + q^{-2m}) \; ,
$$
 which leads to the equation
 $\alpha^{(m)}_{N+1}=\alpha^{(m)}_N +
 q^{2(m-N-1)}\alpha^{(m-1)}_N$
 solved by (\ref{Ldeta}).
We note that by using (\ref{detD}), (\ref{detD1})
and then evaluating the trace $Tr_{D(m+1 \dots N)}$
by means of (\ref{detD3}),
the elements $a_m(L)$ can also be
written in the form \cite{IsBonn}
\be
\lb{sigm}
a_m(L) = \frac{\alpha^{(m)}_N}{{\det}_q(D)} \; Tr_{D(1 \dots N)}
 \left(  A_{1 \dots N} L_{\widetilde{1}} \cdots L_{\widetilde{m}} \right)
 = q^{m^2} \; Tr_{D(1 \dots m)}
 \left(  A_{1 \dots m} L_{\widetilde{1}} \cdots L_{\widetilde{m}} \right)
 \; .
\ee
Then we have
\be
\lb{Ldeti}
\begin{array}{c}
L_{\widetilde{1}} L_{\widetilde{2}} \cdots L_{\widetilde{m}} =
[ L_1 (\R_1 L_1) (\R_2\R_1L_1) \dots (\R_{m-1} \dots \R_1 L_1)]
\R_{(1 \to m-1)} \dots \R_{(1 \to 2)} \R_1 = \\ [0.2cm]
= L_{\underline{1}} L_{\underline{2}} \cdots L_{\underline{m}}
(\R_1 \dots \R_{(m-1 \leftarrow 1)}) (\R_{(1 \to m-1)} \dots \R_1)
= L_{\underline{1}} L_{\underline{2}} \cdots L_{\underline{m}} \, y_2 \, y_3 \dots y_m \; ,
\end{array}
\ee
and analogously
\be
\lb{Ldeti0}
L_{\widetilde{1}} L_{\widetilde{2}} \cdots L_{\widetilde{m}} =
L_{\widetilde{m}} \cdots L_{\widetilde{2}} L_{\widetilde{1}} =
y_2 \, y_3 \dots y_m \, L_{\overline{m}} \dots L_{\overline{2}} L_{\overline{1}}  \; .
\ee
where we used notation $\R_{(1 \to m)}$
and $\R_{(m \leftarrow 1)}$ given in (\ref{04.4b}),
matrices $y_1 =1$, $y_2 = \R_1^2,...$,
$y_{k+1} = \R_k y_k \R_k$
define a commutative set
 $[y_k, \; y_\ell]=0$ and elements
$$
L_{\underline{k+1}}=\R_k L_{\underline{k}} \R_k^{-1} \; ,
\;\;\;\;\;\;\;\;
L_{\overline{k+1}}=\R_k^{-1} L_{\overline{k}} \R_k \; ,
$$
were introduced in (\ref{0.5ab}).
According to the identities (\ref{Ldeti}), (\ref{Ldeti0})
and taking into account (\ref{L123}) one can write (\ref{sigm})
as
\be
\lb{sigm1}
a_m(L) = q^{m} \, Tr_{D(1 \dots m)}  \left(
 \, A_{1 \to m} L_{\underline{1}} L_{\underline{2}} \cdots L_{\underline{m}}
 \right) = q^{m} \, Tr_{D(1 \dots m)}  \left(
 \, A_{1 \to m} L_{\overline{m}} \cdots L_{\overline{2}} L_{\overline{1}}
 \right) \; .
\ee

The elements $a_m(L)$
(\ref{Ldets}), (\ref{sigm}), (\ref{sigm1}) are central elements for the
reflection equation algebra (\ref{3.1.23a}).
Indeed, these elements are obtained from the general center elements (\ref{center})
by substitution $X = A_{1 \dots m}$.

\subsubsection[Differential calculus on the $RTT$
algebras. QG
covariant connections and curvatures]{Differential calculus on the $RTT$
algebras. Quantum group
covariant \protect\\ connections and curvatures\label{difcal}}

\noindent
{\bf a. Bicovariant differential algebras and
quantum $BRST$ operator} \\
 For the Hecke type $R$ matrix (\ref{3.3.7}),
(\ref{3.3.7aa})
one can define  \cite{Sudb} --
\cite{Isaev2} (see also references therein)
 the bicovariant differential complex on the $RTT$ algebra
 \be
 \lb{difcom}
 \R \, T_1 \, T_2 = T_1 \, T_2 \, \R \; , \;\;\;\;
 \R \, dT_1 \, T_2 = T_1 \, dT_2 \, \R^{-1} \; , \;\;\;\;
 \R \, dT_1 \, dT_2 = - dT_1 \, dT_2 \, \R^{-1} \; ,
 \ee
 where $\R := \R_{12}$ and generators $dT^i_j$ $(i,j=1,...,N)$
 are interpreted as differentials
 of the elements $T^i_j$. The algebra (\ref{difcom})
  is a graded (exterior) Hopf algebra \cite{Brze}
  with structure mappings \cite{IsPop}, \cite{Isaev2}
$$
\begin{array}{c}
  \Delta(T) =
  T \otimes T \, , \;\;\;\;
 \epsilon(T) = I \, , \;\;\;\;  S(T) = T^{-1} \, , \\ [0.2cm]
\Delta(dT^i_{\;\; j}) =
dT^i_{\;\; k} \otimes T^k_{\;\; j} +
T^i_{\;\; k} \otimes dT^k_{\;\; j}  \, , \;\;\;\;
 \epsilon(dT^i_{\;\; j}) = 0 \, , \;\;\;\;
 S(dT^i_{\;\; j}) = -
 (T^{-1} \, dT \, T^{-1})^i_{\;\; j} \, ,
 \end{array}
$$
where in the first line we
use the index free matrix notation, the grading is
$\mathfrak{gr}(T^i_{\; k})=0$,
$\mathfrak{gr}(d T^i_{\; k})=1$ and
 $\otimes$ is the graded tensor product.
   We extend
  the graded algebra (\ref{difcom}) by adding \cite{SWZum},
  \cite{ZumB}, \cite{Isaev1} -- \cite{22''} new generators
 $\partial^i_{\, j} := \partial/\partial T^j_{\, i}$
 (quantum derivatives) and ${\bf i}^i_{\; j}$
 (quantum inner derivatives) with
 commutation relations (cf. (\ref{difcom}))
  \be
 \lb{difcom1}
 \R \, \partial_2 \, \partial_1 =
 \partial_2 \, \partial_1 \, \R \, , \;\;\;
  \R \, {\bf i}_2 \, \partial_1  =
  \partial_2 \, {\bf i}_1 \, \R^{-1} \, , \;\;\;
 \R \, {\bf i}_2 \, {\bf i}_1 =
 -{\bf i}_2 \, {\bf i}_1 \, \R^{-1} \, .
 \ee
 Assume that the matrix $\partial$
 is formally invertible. Then
 the algebra (\ref{difcom1})
  is a graded Hopf algebra with structure mappings
  \cite{Isaev1}
$$
\begin{array}{c}
  \Delta(\partial) =
  \partial \otimes \partial \, , \;\;\;\;
 \epsilon(\partial) = I \, , \;\;\;\;
 S(\partial) = \partial^{-1} \, , \\ [0.2cm]
\Delta({\bf i}) =
{\bf i} \otimes \partial
+ \partial \otimes {\bf i}  \, , \;\;\;\;
 \epsilon({\bf i}) = 0 \, , \;\;\;\;
 S({\bf i}) = - \partial^{-1} \, {\bf i} \, \partial^{-1} \, ,
 \end{array}
$$
where
$\mathfrak{gr}(\partial^i_{\; k})=0$,
$\mathfrak{gr}({\bf i}^i_{\; k})=1$ and
 $\otimes$ is the graded tensor product.
 Now we define the cross-product
  (Heisenberg double)
 of the algebras (\ref{difcom}) and (\ref{difcom1})
 with cross-commutation
 relations
 \be
 \lb{difcom2}
  \partial_1 \, \R \,  dT_1 =
 dT_2 \, \R^{-1} \, \partial_2 \, ,  \;\;\;
  T_2 \, \R \,  {\bf i}_2 =
 {\bf i}_1  \, \R^{-1} \, T_1 \, ,
 \ee
 \be
 \lb{difcom2b}
   \partial_1 \, \R^{\bar{\epsilon}} \,  T_1 =
 T_2 \, \R^{-\bar{\epsilon}} \,
 \partial_2 + I_1 I_2 \, ,
   \;\;\; {\bf i}_1 \, \R^{\epsilon} \,  dT_1 +
 dT_2 \, \R^{\epsilon} \, {\bf i}_2 = I_1I_2 \; .
 \ee
 where $\epsilon = \mp1$ and
 $\bar{\epsilon}= \pm 1$.
 These relations are postulated in such a way that 
 they are
 bi-covariant with respect to the left and
  right coaction of the $RTT$ algebra:
  $$
 T \to T_{L} \otimes T \otimes T_{R} \; , \;\;\;\;
 dT \to T_{L} \otimes dT \otimes T_R
 \; , \;\;\;\; \partial \to
 T_R^{-1}\otimes \partial \otimes T_L^{-1}
 \; , \;\;\;\; {\bf i} \to
 T_R^{-1}\otimes {\bf i} \otimes T_L^{-1} \; .
  $$
  Here elements of matrices
  $T_{L}$ and $T_{R}$ generate two $RTT$ algebras
  (the first relation in (\ref{difcom})).
 We note that the
 bi-covariance do not fix  uniquely
 the relations (\ref{difcom2}),
 (\ref{difcom2b}). Actually
 we have four graded bi-covariant
 algebras ${\cal A}^{\epsilon,\bar{\epsilon}}$
 with generators $\{T,dT,\partial,{\bf i} \}$
 and with different choice of signs
 $\epsilon$ and $\bar{\epsilon}$ in (\ref{difcom2b}).
 It was shown in \cite{22''}, \cite{Isaev1}, that there are explicit inner automorphisms in
 ${\cal A}^{\epsilon,\bar{\epsilon}}$
 which relates all these algebras.

  Now we introduce
 the left $\overline{L}$ and the right $L$
 invariant vector fields in the
 algebra ${\cal A^{\epsilon,\bar{\epsilon}}}$, corresponding
 left $\overline{\Omega}
 \in {\cal A^{\epsilon,\bar{\epsilon}}}$
 and the right $\Omega\in
 {\cal A^{\epsilon,\bar{\epsilon}}}$
 invariant  differential
 1-forms,
 inner derivatives $\overline{\cal I}$,
  ${\cal I} \in  {\cal A^{\epsilon,\bar{\epsilon}}}$
 and special invariant operators $\overline{W}$,
 $W\in  {\cal A^{\epsilon,\bar{\epsilon}}}$:
 \be
 \lb{covops}
 \begin{array}{rl}
 \overline{L} \, := \,
 I - \bar{\epsilon}
 \lambda \, \partial \, T \, , & \;\;\;
 L \, := \, I - \bar{\epsilon}
 \lambda \, T \, \partial
 \; \equiv \; T \, \overline{L} \, T^{-1}
  \; , \\ [0.1cm]
 \overline{\Omega} \, := \, T^{-1} \, dT
 \, , \;\;\;\;\;\; \overline{\cal I} : ={\bf i}\, T
  \, , & \;\;\;
 \Omega \, :=  \, dT \, T^{-1} \; \equiv \;
 T \, \overline{\Omega} \, T^{-1} \, , \;\;\;\;\;\;
 {\cal I}: =T\,{\bf i}  \; , \\ [0.1cm]
 \overline{W}
  \, := \, 1 - \epsilon \lambda\; {\bf i} \, dT
 = 1 - \epsilon \lambda \; \overline{\cal I} \, \overline{\Omega} \, , & \;\;\;
 W \, := \,
 1 - \epsilon \lambda \, dT \, {\bf i}
 = 1 - \epsilon \lambda \; \Omega \, {\cal I} \; ,
 \end{array}
  \ee
  where as usual $\lambda = q - q^{-1}$.
  Last relations in (\ref{difcom}), (\ref{difcom1})
  and (\ref{difcom2}) lead to
  (see \cite{Isaev1})
 \be
 \lb{difcom7}
 \overline{W}_{\!1} \, dT_2 =
 dT_2 \, \R^\epsilon \, \overline{W}_{\!2}\,
 \R^\epsilon \, , \;\;\;
 dT_1 \, W_2 = \R^\epsilon \, W_1 \, \R^\epsilon \,
  dT_1 \, , 
  \ee
  \be
 \lb{difcom7b}
 \begin{array}{c}
 {\bf i}_2 \, \overline{W}_{\!1} =
   \R^\epsilon \, \overline{W}_{\!2}\,
    \R^\epsilon \, {\bf i}_2 \, ,
 \;\;\; W_2 \, {\bf i}_1 =
  {\bf i}_1 \, \R^\epsilon \, W_1 \, \R^\epsilon \, ,
  \end{array}
 \ee
  and from these identities we immediately
  obtain \cite{Isaev1}
  \be
  \lb{difcom8}
  \begin{array}{c}
  \overline{W}_{\!1}\, W_2 =
  W_2 \, \overline{W}_{\!1} \, , \;\;\;
  \overline{W}_{\!2}\,
  \R^\epsilon \, \overline{W}_{\!2}\, \R^\epsilon =
  \R^\epsilon \, \overline{W}_{\!2} \,
 \R^\epsilon  \, \overline{W}_{\!2} \, , \;\;\;
   W_1 \,\R^\epsilon \, W_1 \, \R^\epsilon
   = \R^\epsilon \, W_1 \, \R^\epsilon \, W_1 \; .
 \end{array}
  \ee
  Operators (\ref{covops}) also
  obey (see e.g. \cite{SWZum},
  \cite{Isaev1}, \cite{22''})  the following
  relations (cf. (\ref{3.1.23a}),
  (\ref{3.1.23b}), (\ref{tangb}))
  \be
 \lb{difcom3}
 \begin{array}{c}
 \left\{
 \begin{array}{c}
 \overline{L}_1 \, T_2 = T_2 \R^{-\bar{\epsilon}}
 \overline{L}_2 \R^{-\beps} \; , \;\;\;
T_1 L_2 = \R^{\beps} L_1 \R^{\beps}
T_1 \; ,  \\ 
\partial_2 \, \overline{L}_1 = \R^{-\beps}
\, \overline{L}_2 \, \R^{-\beps} \,
\partial_2  \; , \;\;\;
L_2 \, \partial_1 = \partial_1 \,
\R^{\beps} \, L_1 \, \R^{\beps} \; ,
\end{array} \Rightarrow \right.
\\ [0.4cm]
 \overline{L}_1 \, L_2 =
 L_2 \,  \overline{L}_1 \, , \;\;\;
\R^{-\beps} \, \overline{L}_2 \,
\R^{-\beps} \, \overline{L}_2 =
\overline{L}_2 \, \R^{-\beps}
\, \overline{L}_2 \, \R^{-\beps} \, , \;\;\;
 \R^{\beps} \, L_1 \, \R^{\beps} \, L_1 =
 L_1 \, \R^{\beps} \, L_1 \,
 \R^{\beps} \; ,
\end{array}
\ee
and in addition we have
  \be
 \lb{difcom9}
 \begin{array}{c}
 \left\{
 \begin{array}{c}
 \overline{L}_1 \, dT_2 = dT_2 \R^{-1}
 \overline{L}_2 \R \; , \;\;\;
dT_1 L_2 = \R^{-1} L_1 \R dT_1 \; ,  \\ 
{\bf i}_2 \, \overline{L}_1 = \R^{-1}
\, \overline{L}_2 \, \R \, {\bf i}_2  \; , \;\;\;
L_2 \, {\bf i}_1 = {\bf i}_1 \,
\R^{-1} \, L_1 \, \R \; ,
\end{array} \Rightarrow \right.
\\ [0.4cm]
\R^{-1} \, \overline{L}_2 \,\R \, \overline{W}_2 =
\overline{W}_2 \, \R^{-1}
\, \overline{L}_2 \, \R \, , \;\;\;
 W_1 \, \R^{-1} \, L_1 \, \R  =
  \R^{-1}\, L_1 \, \R \, W_1  \; , \\ [0.2cm]
 \R^{-1} \, \overline{L}_2 \,\R \,
 \overline{\cal I}_2 =
\overline{\cal I}_2 \, \R^{-\beps}
\, \overline{L}_2 \, \R^{-\beps} \, , \;\;\;
 W_1 \, \R^{-1} \, L_1 \, \R  =
  \R^{-1}\, L_1 \, \R \, W_1  \; .
\end{array}
\ee
\begin{proposition}\label{propn}
(See \cite{Isaev1}) The defining relations (\ref{difcom}),
(\ref{difcom1}) and (\ref{difcom2}),
(\ref{difcom2b})  are in an
agreement with formulas
 \be
 \lb{TxdTy}
\begin{array}{c}
\R T_1(x) T_2(x) = T_1(x) T_2(x) \R \; , \;\;\;
\R dT_1(y) T_2(x) = T_1(x) dT_2(y) \R^{-1}
 \; , \\ [0.2cm]
\R dT_1(y) \, dT_2(y) =
-dT_1(y) \,dT_2(y) \R^{-1} \; \Leftrightarrow \;
\R^\epsilon dT_1(y) \, dT_2(y) =
-dT_1(y) \,dT_2(y) \R^{-\epsilon}  \; .
\end{array}
 \ee
where $x,y$ are parameters, $\R := \R_{12}$ and
$$
T(x) = \frac{1}{x + T \partial} \, T \equiv
\frac{1}{T^{-1} \, x + \partial}  \; , \;\;\;\;\;
dT(y) = \frac{1}{y + dT \, {\bf i} } \, dT \equiv
 dT \, \frac{1}{y + {\bf i} \, dT }  \; .
$$
\end{proposition}
{\bf Proof.} We write the first relation in
(\ref{TxdTy}) as
$$
(T^{-1}_2 \, x + \partial_2)
(T^{-1}_1 \, x + \partial_1)\, \R =
\R \, (T^{-1}_2 \, x + \partial_2)
(T^{-1}_1 \, x + \partial_1) \; .
$$
The terms of order $x^2$ and $x^0$ give respectively
first relations in (\ref{difcom}) and (\ref{difcom1}).
The terms of order $x^1$ yield relation compatible
to the first relation in (\ref{difcom2b}). For
 the second relation in (\ref{TxdTy}) we obtain
 $$
 (y + dT_2 {\bf i}_2) (T^{-1}_1 \, x + \partial_1)
 \, \R \, dT_1 =dT_2 \, \R^{-1} \,
 (T^{-1}_2 \, x + \partial_2)
 (y + {\bf i}_1 dT_1 ) \; ,
 $$
 and in the orders $(xy)^0$, $x^1$, $(xy)^1$, $y^1$
 we respectively deduce second relations in
 (\ref{difcom1}), (\ref{difcom2}),
 (\ref{difcom}) and first relation
 in (\ref{difcom2}). Finally, the third relation
 in (\ref{TxdTy}) is represented as
 $$
  (y + dT_1 \, {\bf i}_1 )\R dT_1
  \frac{1}{(y + {\bf i}_1 \,dT_1 )} \, dT_2 =
  - dT_1 \,
 \frac{1}{(y + dT_2 {\bf i}_2  )} dT_2 \R^{-1}
 (y + {\bf i}_2 dT_2 ) \; ,
 $$
 and, after the change of parameter $y \to \frac{y-1}{\epsilon \lambda}$, introducing  operators $W,\overline{W}$ (\ref{covops}),
 and applying
 eqs. (\ref{difcom7}), we write it
 in the form
$$
 (y + \R^\epsilon W_1 \R^\epsilon)
 (y + W_1 )\R^\epsilon dT_1 \, dT_2 = - dT_1 \, dT_2
 \R^{-\epsilon} (y + \overline{W}_2 )
 (y + \R^\epsilon \overline{W}_2 \R^\epsilon) \; .
 $$
 The terms of order $y^2$ give the last relation in
 (\ref{difcom}). The terms of order $y^0$ and $y^1$
 are identities in view of the relations
 (\ref{difcom7}) and (\ref{difcom8}) which encode last relations in
 (\ref{difcom1}) and (\ref{difcom2b}). \hfill \qed

\vspace{0.2cm}

\noindent
{\bf Corollary}. The formulas (\ref{TxdTy}) (incorporating
 the whole differential algebra
 (\ref{difcom}), (\ref{difcom1}) and (\ref{difcom2}))
 have the structure of eqs. (\ref{difcom})
 for which one can easily establish
  the PBW property. This indicates that
   the whole differential algebra (\ref{difcom}), (\ref{difcom1}) and (\ref{difcom2}), (\ref{difcom2b})) is also of the PBW type (the flat deformation of the
   differential algebra obtained in the classical limit
   $q \to 1$, $\R \to P$).
 We also note that the signs of powers of
 $\R$-matrices are flashing
 $\R^{\beps} \to \R^{-\beps}$
 in the first relation
 in (\ref{difcom2b}) and not flashing
 $\R^{\epsilon} \to \R^{\epsilon}$ in the
 second relation in (\ref{difcom2b}). This is important,
 otherwise relations (\ref{TxdTy}) are not
 fulfilled.

 \vspace{0.2cm}

Here we present additional commutation relations
 for the invariant operators (\ref{covops}). These
 relations are useful from a technical point of view.
  \be
 \lb{difcom4}
 \begin{array}{rl}
 \overline{\Omega}_1 \, T_2 = T_2 \, \R^{-1} \,
 \overline{\Omega}_2 \, \R^{-1}  \, , & \;\;\;
 T_1 \, \Omega_2 = \R_{12} \, \Omega_1 \, \R_{12}\, T_1
\; ,  \\ [0.2cm]
 \overline{W}_1 \, T_2 = T_2 \, \R \,
 \overline{W}_2 \, \R^{-1} \, , & \;\;\;
 T_1 \; W_2 = \R \, W_1 \R^{-1} \, T_1 \, ,
 \\ [0.2cm]
 \overline{\cal I}_1 \, T_2 = \R \, T_2 \, \R \,
 \overline{\cal I}_2 \, , & \;\;\;
  T_1 \; {\cal I}_2 = \R^{-1} \, {\cal I}_1
  \R^{-1} \, T_1 \, ,  \\ [0.2cm]
 \overline{L}_1 \, dT_2 = dT_2 \, \R_{12}^{-1} \,
 \overline{L}_2 \, \R_{12} \; , & \;\;\;\;
dT_1 \, L_2 = \R_{12}^{-1} \, L_1 \, \R_{12}\, dT_1
\; ,  \\ [0.2cm]
\R^{-1} \, \overline{\Omega}_2 \,
\R^{-1} \, \overline{\Omega}_2 =
- \overline{\Omega}_2 \, \R^{-1}
\, \overline{\Omega}_2 \, \R
\; , & \;\;\;\;
\R_{12} \, \Omega_1 \, \R \, \Omega_1 =
- \Omega_1 \, \R \, \Omega_1
\, \R^{-1} \; ,  \\ [0.2cm]
\R^{-\beps} \, \overline{L}_2 \,
\R^{-\beps} \, \overline{\Omega}_2 =
\overline{\Omega}_2 \, \R^{-1} \,
 \overline{L}_2 \, \R \; , & \;\;\;\;
\Omega_1 \, \R^{\beps} \, L_1 \,
\R^{\beps} =
\R^{-1} \, L_1 \, \R \, \Omega_1 \; .
\end{array}
 \ee
 \be
 \lb{difcom5}
 \begin{array}{rl}
 \R_{12} \, \overline{\cal I}_2 \,
\R_{12} \, \overline{\cal I}_2 =
- \overline{\cal I}_2 \, \R_{12}
\, \overline{\cal I}_2 \, \R_{12}^{-1}
\; , & \;\;\;\;  \R_{12}^{-1} \, {\cal I}_1 \,
\R_{12}^{-1} \, {\cal I}_1 =
- {\cal I}_1 \, \R_{12}^{-1} \, {\cal I}_1
\, \R_{12} \; , \\ [0.2cm]
 \R_{12} \overline{\cal I}_2 \, \R_{12}^{\epsilon} \,
\overline{\Omega}_2 \,  +
\overline{\Omega}_2 \, \R_{12}^{\epsilon}
\, \overline{\cal I}_2  \, \R_{12} = \R_{12}
\, , & \;
{\cal I}_1 \, \R_{12}^{\epsilon} \, \Omega_1
\,\R_{12}  +
 \R_{12}  \, \Omega_1 \, \R_{12}^{\epsilon}
  \, {\cal I}_1 =  \R_{12} \; ,  \\ [0.2cm]
  \overline{\Omega}_2 \, \R_{12}^{\epsilon}
\, \overline{W}_2  \, \R_{12}^{\epsilon} =
 \R_{12}\, \overline{W}_2 \, \R_{12}^{-1} \,
\overline{\Omega}_2  \,  , & \;
\R_{12}^{\epsilon} \, W_1 \,
\R_{12}^{\epsilon} \, \Omega_1 =
 \Omega_1 \, \R_{12} \, W_1 \, \R_{12}^{-1} \; ,
 \\ [0.2cm]
 \R_{12}^{\epsilon} \, \overline{W}_2 \,
 \R_{12}^{\epsilon} \, \overline{\cal I}_2   =
\overline{\cal I}_2 \, \R_{12}
\, \overline{W}_2  \, \R_{12}^{-1} \,  , & \;
{\cal I}_1  \, \R_{12}^{\epsilon} \, W_1
\, \R_{12}^{\epsilon} =
\R_{12}\, W_1 \, \R_{12}^{-1} \, {\cal I}_1   \; ,
 \end{array}
 \ee
where we denote $\overline{\cal I}={\bf i}\, T$  and
  ${\cal I}=T\,{\bf i}$.

  Now instead of the left $\overline{L}$ and
  the right $L$ invariant vector fields we introduce
  new left and right invariant operators \cite{Isaev1}
  $$
  \overline{\cal L} := \overline{W}\, \overline{L} =
 (1 - \epsilon \lambda\; {\bf i} \, dT)
  (1 - \bar{\epsilon} \lambda\; \partial \, T)  \; , \;\;\;\;\;
    {\cal L} := W \, L =
 (1 - \epsilon \lambda\; dT \, {\bf i})
  (1 - \bar{\epsilon} \lambda\; T \, \partial) \; .
  $$
  As we mentioned above, all algebras
  ${\cal A}^{\epsilon,\beps}$ are equivalent
  for different choice of the signs $\epsilon,\beps$.
  For simplicity, further we consider only the
  left invariant operators and fix $\epsilon=1$ and
  $\beps = -1$. In accordance with formulas (\ref{difcom4}), (\ref{difcom5}) for $\epsilon =1$ and $\beps = -1$, we have the following statement.
  \begin{proposition}\label{extalg}
  The complete set of commutation relations
  for the exterior differential algebra
  $\Gamma^{\wedge} \subset
  {\cal A}^{\epsilon,\beps}|_{\epsilon =-\beps=1}$ with generators
  $T$, $\overline{\cal L}$, $\overline{\cal I}$,
  $\overline{\Omega}$
  is \cite{SWZum},
  \cite{18''}, \cite{Isaev1}, \cite{IsOgd2}:
  \be
  \lb{difcom10a}
 \R T_1 T_2 = T_1 T_2 \R \, , \;\;\;
 \overline{\Omega}_1 \, T_2 = T_2 \, \R^{-1} \,
 \overline{\Omega}_2 \, \R^{-1}  \, , \;\;\;
 \R^{-1} \, \overline{\Omega}_2 \,
\R^{-1} \, \overline{\Omega}_2 =
- \overline{\Omega}_2 \, \R^{-1}
\, \overline{\Omega}_2 \, \R \; ,
\ee
  \be
  \lb{difcom10}
  \begin{array}{c}
  \overline{\cal L}_2 \R \overline{\cal L}_2 \R =\R
 \overline{\cal L}_2 \R \overline{\cal L}_2 \, , \;\;\;
 \overline{\Omega}_2 \R \overline{\cal L}_2 \R =\R
 \overline{\cal L}_2 \R \overline{\Omega}_2 \, , \;\;\;
 \overline{\cal I}_2 \R \overline{\cal L}_2 \R =\R
 \overline{\cal L}_2 \R \overline{\cal I}_2 \, ,
 \\ [0.2cm]
 \overline{\cal L}_1 T_2 =
 T_2 \R \overline{\cal L}_2 \R \, , \;\;\;
 \overline{\cal I}_1 \, T_2 = \R \, T_2 \, \R \,
 \overline{\cal I}_2 \, , \;\;\;
  \R_{12} \, \overline{\cal I}_2 \,
\R_{12} \, \overline{\cal I}_2 =
- \overline{\cal I}_2 \, \R_{12}
\, \overline{\cal I}_2 \, \R_{12}^{-1} \, ,
 \end{array}
  \ee
\be
  \lb{difcom10b}
  \R_{12} \overline{\cal I}_2 \, \R_{12} \,
\overline{\Omega}_2 \,  +
\overline{\Omega}_2 \, \R_{12}
\, \overline{\cal I}_2  \, \R_{12} = \R_{12} \, ,
  \ee
  where $\R:=\R_{12}$ is the Hecke type $R$-matrix.
  \end{proposition}
 We note that, for the Hecke type $\R$-matrix, the differential algebra (\ref{difcom10a}) -- (\ref{difcom10b})
 is identical to the differential algebra $\Gamma^\wedge$, proposed in the papers
 \cite{SWZum}, \cite{18''}, \cite{IsOgd2},
 up to the relation (\ref{difcom10b})
 which is written
  in that papers as
 $\R_{12} \overline{\cal I}_2 \, \R_{12} \,
\overline{\Omega}_2 \,  +
\overline{\Omega}_2 \, \R_{12}
\, \overline{\cal I}_2  \, \R_{12} = - \R_{12}$.
The change of the sign in the right hand side
of the relation (\ref{difcom10b})
 can be achieved by the transformation
 $\overline{\cal I} \to -\overline{\cal I}$.

 By using the general construction \cite{IsOgd}
 of the BRST charge for an arbitrary quantum Lie algebra, we have constructed in \cite{IsOgd2}
 a BRST operator ${\sf Q}$ for the differential algebra
(\ref{difcom10}) in the following
form\footnote{In all formulas in \cite{IsOgd2} we
should make the change of notation:
$\omega \to \overline{\Omega}$,
$J \to - \overline{\cal I}$,
$L \to \overline{\cal L}$.}:
\be
\lb{Qbrst}
{\sf Q} =  Tr_{_{\!\!\cal Q}}
\left( \overline{\Omega} \,
\frac{(\overline{\cal L}- 1 )}{\lambda}
+ \overline{\Omega}  \, \overline{\cal L}\,
\frac{\overline{\Omega}  \, \overline{\cal I}}
{(1  - \lambda \overline{\Omega} \,
 \overline{\cal I} )} \right) = -\frac{1}{\lambda}
\, Tr_{_{\!\!\cal Q}} (\overline{\Omega}) +
\frac{1}{\lambda} \,
Tr_{_{\!\!\cal Q}} \left( \Theta  \right) \, ,
\ee
where $\Theta := \overline{\Omega} \,
\overline{\cal L} \,
(1  - \lambda \overline{\Omega} \,
 \overline{\cal I} )^{-1}$ and
 $Tr_{_{\!\!\cal Q}}(X) := q^{2d} \, Tr(Q X)$ is the second
 quantum trace in (\ref{qtrDQ}). The
 normalization factor $q^{2d}$ is introduced to have $Tr_{_{\!\!\cal Q}1}(\R^{-1})=I_2$
 (see (\ref{sk9}), (\ref{QD})). We note that
 the left-invariant operator $\widetilde{W}:=(1  - \lambda
 \overline{\Omega} \, \overline{\cal I})$, appeared
 in (\ref{Qbrst}), differs from the operator
 $\overline{W}=(1  - \lambda \overline{\cal I} \,
 \overline{\Omega})$ defined in (\ref{covops}).
 For the operator $\widetilde{W}$ we have
 \be
 \lb{difcom11}
 \begin{array}{c}
 \widetilde{W}_2 \, \R \, \overline{\cal L}_2 \, \R =
  \R \, \overline{\cal L}_2 \, \R \,  \widetilde{W}_2  \; , \;\;\;\;
 \widetilde{W}_2 \, \R \, \overline{\Omega}_2 \, \R^{-1} =
 \R^{-1} \,  \overline{\Omega}_2 \, \R^{-1}  \,  \widetilde{W}_2  \; , \\ [0.2cm]
 \widetilde{W}_2 \, \R  \, \overline{\cal I}_2  \R^{-1} =
 \R  \, \widetilde{\cal I}_2  \R  \,  \widetilde{W}_2   ,  \;\;\;
 \widetilde{W}_2 \, \R \, \widetilde{W}_2 \, \R =
 \R \,  \widetilde{W}_2 \, \R  \,  \widetilde{W}_2   ,  \;\;\;
 \Theta_2  \R^{-1}  \widetilde{W}_2  \R =
 \R   \widetilde{W}_2  \R    \Theta_2  .
 \end{array}
 \ee
 In the definition (\ref{Qbrst}),
 the differential 1-forms
 $\overline{\Omega}^i_{\; k}$ and the inner derivatives
 $\overline{\cal I}^i_{\; k}$ play the role of
 the ghost and anti-ghost variables.
 One can check directly \cite{IsOgd2}
 that the BRST operator ${\sf Q}$
given by (\ref{Qbrst}) satisfies:
\begin{equation}
\label{diffq}
 {\sf Q}^2 = 0 \; ,  \;\;\;\;\;\;\;
[ {\sf Q} , \, \overline{\cal L} ] = 0 \; ,
\end{equation}
\begin{equation}
\label{diffT}
[{\sf Q} , \, T] = T \, \overline{\Omega}
\equiv d \, T ,
\;\;\;\;\;\;  [{\sf Q} , \, \overline{\Omega}]_+ =
- \overline{\Omega}^2
\equiv d \, \overline{\Omega} ,
\end{equation}
\begin{equation}
\label{cartan}
[ {\sf Q} , \, \overline{\cal I} ]_+ =
\frac{1}{\lambda} \, (\overline{\cal L}-1)  \; .
\end{equation}
The (anti)commutator with ${\sf Q}$
(relations (\ref{diffT})) defines the exterior differential operator over the differential
algebra (\ref{difcom});
it provides the structure of the de Rham complex on the
subalgebra with generators $T^i_j$ and
$\overline{\Omega}^i_j$
(the de Rham complex over the quantum
group $GL_q(N)$ has been firstly considered
by Yu.I.Manin, G.Maltsiniotis and B.Tsygan \cite{ManMT}).

To obtain relations (\ref{diffq}) -- (\ref{cartan})
one has to use the invariance property (\ref{RER})
of the quantum trace ${\rm Tr}_{_{\cal Q}}$
and relations
\begin{eqnarray}
\label{th}
\hat{R} \, \Theta_2 \, \hat{R}^{-1} \,
 \overline{\Omega}_2
= - \overline{\Omega}_2 \, \hat{R}^{-1} \,
\Theta_2 \, \hat{R} \; , & \;\;\;\;\;\;
\hat{R} \, \Theta_2 \, \hat{R}^{-1} \, \Theta_2 =
- \Theta_2 \, \hat{R}^{-1} \, \Theta_2 \, \hat{R}^{-1}  \; , \\ [0.2cm]
\label{TL}
\hat{R}^{-1} \, \Theta_2 \, \hat{R} \,
\overline{\cal L}_2 =
\overline{\cal L}_2 \, \hat{R} \, \Theta_2 \, \hat{R}^{-1}  \; ,  & \;\;\;\;\;\;
 \Theta_1 \, T_2 =
T_2 \, \hat{R}^{-1} \, \Theta_2 \, \hat{R}  \; ,
\;\;\;\;\;\; \quad \quad
  \end{eqnarray}
 \be
 \label{JT}
 \overline{\cal I}_2 \,
 \hat{R} \, \Theta_2 \, \hat{R}^{-1} +
\hat{R}^{-1} \, \Theta_2 \, \hat{R} \,
\overline{\cal I}_2
= \overline{\cal L}_2 \,
\widetilde{W}^{-1}_2
\, \hat{R}^{-1} \,
 \widetilde{W}_2  \; ,
\ee
which follow from  (\ref{difcom10a}) --
(\ref{difcom10b}) and (\ref{difcom11}).
In particular, the condition ${\sf Q}^2 =0$
follows from
the last equation in (\ref{difcom10a}) and equations
(\ref{th}) which lead to identities
\begin{equation}
\left( Tr_{_{\!\!\cal Q}} (\overline{\Omega})
\right)^2 = 0\,  ,\ \
\left( Tr_{_{\!\!\cal Q}} ( \Theta ) \right)^2 = 0
\,  , \ \
[ Tr_{_{\!\!\cal Q}} ( \Theta ) , \,
Tr_{_{\!\!\cal Q}} (\overline{\Omega}) ]_{+} = 0 \, .
\end{equation}
Here we take into account that $Tr_{_{\!\!\cal Q}}(\Theta^2)=0=
Tr_{_{\!\!\cal Q}} (\overline{\Omega}^2)$
(see Sect. 4 in \cite{Isaev1}).
Finally we note (for details see \cite{IsOgd2})
that the operator ${\sf Q}$ given by (\ref{Qbrst}) has the
correct classical limit
for $q \rightarrow 1$, $\lambda = q - q^{-1}
\rightarrow 0$, $\R_{12} \to P_{12}$ and
$\overline{\cal L} \rightarrow {\bf 1} + \lambda X$, where
elements $X^i_{\; k}$ are interpreted as Lie
algebra generators.

\vspace{0.3cm}

\noindent
{\bf b. Quantum group
covariant connections and curvatures} \\
To proceed further we introduce the
$Z_{2}$-graded algebra (denoted by ${\cal E}$) of quantum
 hyperplane with generators
$\{ e_{i},(de)_{i} \}$   $\ (i=1,2,...N)$
 satisfying commutation relations:
\be
\lb{2.5d}
\begin{array}{c}
R_{12} e_{1\rangle } e_{2\rangle } =
ce_{2\rangle }e_{1\rangle } \, , \;\;\;\;
(\pm)cR_{12}(de)_{1\rangle }e_{2\rangle } =
e_{2\rangle } (de)_{1\rangle } \, , \;\;\;\;
R_{12}(de)_{1\rangle }(de)_{2\rangle }
= - {1\over c} (de)_{2\rangle }(de)_{1\rangle } ,
\end{array}
\ee
One can recognize in these
relations (for $(\pm)=+1$) the Wess-Zumino formulas of the covariant
differential calculus on  the bosonic ($c=q$)
and fermionic ($c=-1/q$) quantum hyperplanes \cite{10},
\cite{20,21,WZ,SWZum}, where
$ e_{i} $ are the coordinates of the quantum hyperplanes and
$(de)_{i}$ are the associated differentials
(differential $1$-forms).
 The $Z_{2}$-graded
  algebra ${\cal E} = \sum_{k \geq 0}
 \Omega^{k}({\cal E})$ is the sum of subspaces
 $\Omega^{k}({\cal E})$ of differential $k$-forms.

The left-coaction
$\Delta_{l}$   of the $Z_{2}$-graded
Hopf algebra (\ref{difcom})
to the generators of the algebra ${\cal E}$ is given by
 the following homomorphism:
\be
\lb{2.8aa}
e_{i} \stackrel{\Delta_{l}}{\longrightarrow} \
\widetilde{e}_{i} = T_{ij} \otimes e_{j}, \;\;\;\;\;\;\;
(de)_{i} \stackrel{\Delta_{l}}{\longrightarrow}
(\widetilde{de})_{i} = (dT)_{ij} \otimes
e_{j} + T_{ij} \otimes (de)_{j} .
\ee
The algebra ${\cal E}$ with generators $\{ e, de \}$ becomes now a
left-comodule algebra with respect
to the coaction (\ref{2.8aa}),  since
the all axioms for
the comodule algebras are fulfilled \cite{IsPop}.

Now we assume that the algebra ${\cal E}$ can be extended to
$\bar{\cal E}$ by adding new elements $A_{ij}|_{i,j=1,...,N}$.
We also assume that the differential $d$
can be extended onto the
whole algebra $\bar{\cal E}$ and hence again
this algebra is decomposed as
$\bar{\cal E} = \sum_{k \geq 0} \Omega^{k} (\bar{\cal E})$.
Then, we postulate, first, that the
elements $A_{ij}$ belong to the subspace
$\Omega^{1} (\bar{\cal E})$ and, second, the elements
$(\nabla e)_{i} \in \Omega^{1}(\bar{\cal E})$ defined as
\be
\lb{2.13d}
(\nabla e)_{i} = (de)_{i} - A_{ij} e_{j} ,
\ee
are transformed homogeneously under (\ref{2.8aa})
as the left-comodule
\be
\lb{2.14d}
(\nabla e)_{i} \stackrel{\Delta_{l}}{\longrightarrow}
T_{ij} \otimes (\nabla e)_{j} =
T_{ij} \otimes \left( (de)_{j} - A_{jk}e_{k}\right) .
\ee
According to the classical case we interpret
the operator $A_{ij}$ satisfying (\ref{2.14d})
as a quantum deformation of a gauge potential 1-form
and the operator $\nabla$ as the quantum
version of the covariant derivative. The second action
of $\nabla$ on the both sides of (\ref{2.13d}) gives
\be
\lb{2.22d}
(\nabla(\nabla e))_i = - \left( d(A) - A^{2}\right)_{ij} e_j
 = - F_{ij}\, e_j ,
\ee
where we define the noncommutative analog of the
field strength (curvature) 2-form $F$.
The next action of the covariant
derivative to the formula (\ref{2.22d}) yields the
Bianchi identity which is written in the
standard form $d(F)=[A , \; F]$.
Using (\ref{2.8aa}), (\ref{2.14d})
and (\ref{2.22d}) one can deduce
the noncommutative analog
of the  gauge transformation for the noncommutative
connection 1-form and curvature 2-form as
\be
\lb{2.15d}
A_{ik} \stackrel{\Delta_{l}}{\longrightarrow}
\widetilde{A}_{ik} =
T_{ij} T^{-1}_{lk} \otimes A_{jl}
+ dT_{ij} T^{-1}_{jk} \otimes 1 \; , \;\;\;\;\;\;
F_{ij}{ \stackrel{\Delta_{l}}{\longrightarrow}}
\ \widetilde{F}_{ij}
= \left( T_{ik} T^{-1}_{lj} \right) \otimes F_{kl} .
\ee

As it was argued in \cite{IsPop}, \cite{Isaev7},
the possible choice of the covariant algebra
of the connection 1-form $A$ and curvature 2-form $F$
is given by the defining relations
\be
\lb{isapo1}
F_1 \R_{12} A_1 \R_{12} = \R_{12} A_1 \R_{12} F_1 \; , \;\;\;\;\;
F_1 \R_{12} F_1 \R_{12} = \R_{12} F_1 \R_{12} F_1 \; ,
\ee
\be
\lb{isapo2}
\R_{12} A_1 \R_{12} A_1 + A_1 \R_{12} A_1 \R_{12}^{-1}
 =  \lambda \, g \; ( \R_{12} F_1 + F_1 \R_{12}^{-1} ) \; ,
\ee
where $\lambda = q-q^{-1}$ and
$g$ is an arbitrary parameter. In particular, to
 check the commutation relations (\ref{isapo2})
 for the elements $A_{ij}$ we remark
that there is a representation
for the generators $A_{ij}$,
namely $A = dTT^{-1} \otimes 1$,
 which is related to the flat connection $F_{ij}=0$. Using this
representation and formulas (\ref{difcom}) we
conclude that the generators $A_{ij}$
have to satisfy relation (\ref{isapo2}) with the r.h.s.
equals to zero. In what follows we consider only
the case $g=0$. Note, that the algebra
(\ref{isapo1}), (\ref{isapo2}) (for $g=0$) is covariant
not only under coaction (\ref{2.15d})
of the $RTT$ algebra, but also
is a braided comodule algebra with respect to the
braided coaction of the $RLRL$ algebra \cite{Isaev2}.

Let $\R$ be skew-invertible $R$-matrix for which
we define the quantum traces (\ref{qtrDQ})
  with properties (\ref{RER}) (see also (\ref{3.1.13}),
 (\ref{3.1.14}) and (\ref{3.1.15})).
By analogy with the
classical case, we can consider the noncommutative
version of the invariant Chern characters
\cite{IsPop}, \cite{Isaev7},
\cite{Isaev2}
\be
\lb{2.31d}
C^{(k)} =Tr_{q}(F^{k}) = Tr(DF^{k})=
D_{ij}F_{jj_{1}}F_{j_{1}j_{2}}
\dots F_{j_{k-1}i}\; ,
\ee
where we have used the quantum trace (\ref{3.1.13}),
 with matrix $D$. Chern characters (\ref{2.31d})
 are central elements for the algebra (\ref{isapo1})
 (the proof is the same as proof of (\ref{3.1.23i})).
Applying (\ref{3.1.14}) we immediately obtain
that 2k-forms $C^{(k)}$ (\ref{2.31d})
 are co-invariants under the adjoint co-transformation
 of $F$ given in (\ref{2.15d}). Moreover,
 $C^{(k)}$ are the
closed 2k-forms. Indeed, from the Bianchi identities
$dF= [A, F]$ we deduce
\be
\lb{2.30d}
d C^{(k)} =Tr_{q}(A \, F^{k} - F^{k} \, A)=0 \; ,
\ee
where we have taken into account
(see Eqs. (\ref{qtrs}), (\ref{sk10}), (\ref{isapo1}))
\be
\lb{2.25d}
\begin{array}{c}
Tr_{q}(A \, F^{k})=
Tr_{q1}Tr_{q2}(\R^{-1}_{12} \, \R_{12} \,
A_1 \, \R_{12} \, F^{k}_1)
=Tr_{q1}Tr_{q2}(\R^{-1}_{12} \, F^{k}_1 \, \R_{12} \,
A_1 \, \R_{12} ) = \\ [0.2cm]
=Tr_{q1}Tr_{q2}(F^{k}_1 \, \R_{12} \, A_1)
=Tr_{q}(F^{k} \,  A) \; .
\end{array}
\ee
We note that $Tr_{q}(A^{2k})=0$ for the algebra
(\ref{isapo2}), when $g=0$
(see Proposition 4 in \cite{Isaev1}).
In view of this, the natural conjecture is
that $C^{(k)}$ have to be presented, for $g=0$,
as the exact form
$C^{(k)}= dL_{CS}^{(k)}$, where the noncommutative
Chern-Simons
$(2k - 1)$-forms $L_{CS}^{(k)}$ are
\be
\lb{2.26d}
L^{(k)}_{CS}=Tr_{q} \Bigl( A(dA)^{k-1} + \frac{1}{h^{(k)}_{1}} A^{3}(dA)^{k-2} +
\dots + \frac{1}{h^{(k)}_{k-1}} A^{2k-1} \Bigr) \; ,
\ee
and unknown coefficients $h_j^{(k)}$
depend on the choice of the Hecke matrix $\R$ (in the classical case
$\R_{12} = P_{12}$ and $q=1$, all these coefficients are known
 \cite{ZumZee}). We checked this
conjecture in the case $k=2$, for
$GL_q(N)$ $R$-matrix (\ref{3.3.6a}) and
the special algebra (\ref{isapo2}), when $g=0$. In this case
 we obtained \cite{IsPop}, \cite{Isaev7}
 a noncommutative analog of the three-dimensional
 Chern-Simons term in the form:
 \be
\lb{2.32d}
L^{(2)}_{CS}=Tr_{q} \Bigl( A \, dA + \frac{1}{h^{(2)}_{1}} A^{3} \Bigr) \; , \;\;\;\;\;
h^{(2)}_{1} = -\; \frac{q^{2} + 1 + q^{-2}}{q^{2} + q^{-2}} \; .
\ee
{\bf Remark.} The elements of
the differential calculus on the $RLRL$
(reflection equation) algebra were considered in papers
\cite{9Og}, \cite{9K}, \cite{Isaev2},
\cite{GPS0} (see also references therein).

\subsubsection[$\alpha$-Deformation of the
HD of $RTT$ and $RLRL$ algebras. 
 Quantum CHN identities]{$\alpha$-Deformation
 of the Heisenberg double of $RTT$ and $RLRL$
 algebras. Quantum
 Cayley-Hamilton-Newton identities}

Now, for the right HD (\ref{tangb}) (the algebra
 (\ref{tangb0}) with upper sign),
 we calculate the commutation relations of the
 elements $a_m(L)$ with generators $T^i_j$ of the
$RTT$ algebra defined by the Hecke type $R$-matrix. Note that in the case of the Heisenberg double
of $Fun(SL_q(N))$ and $U_q(sl(N))$  we need to renormalize the
Hecke $R$-matrix: $\R \rightarrow q^{-1/N}\R$
according to (\ref{LTsl}). This leads to the following generalization of the
cross-multiplication rules (\ref{tangb})
(we consider only the right HD)
\be
\lb{tangb1}
T_1 L_2 = \alpha \, \R_{12} L_1 \R_{12} T_1 \; ,
\ee
where $\R$ is a Hecke
$R$-matrix (\ref{3.3.7aa}) of the height $N$ and
 for the special case of the $SL_q(N)$-type HD
we have to fix $\alpha = q^{-2/N}$ (but
generally the constant $\alpha \neq 0$ is arbitrary). So,
the commutation relations (\ref{tangb1}) define
 the one-parameter deformation of the Heisenberg double
 of $RTT$ and $RLRL$ algebras for the Hecke type
 $R$-matrix. Note that
the automorphism (\ref{discr1}) is only correct
for the choice $\alpha = 1$
in (\ref{tangb1}). For example, in view of (\ref{tangb1})
the quantum matrices $(L+x)T$ start to obey the modified
$RTT$ relations
\be
\lb{mrtt}
\R_1 \, (L_1 +x)T_1 \; (\alpha^{-1} L_2  + x)T_2 =
(L_1 +x)T_1 \; (\alpha^{-1} L_2 + x)T_2 \, \R_1 \; .
\ee
However, for the general choice of $\alpha$ in
(\ref{tangb1}), the definition of the
characteristic polynomial (\ref{Ldet3}) is not changed, since instead of (\ref{Ldet55})
we can take
\be
\lb{Ldet5}
\begin{array}{c}
 (L_{1} +x)T_1\,
(\alpha^{-1} L_{2}+x)T_2 \dots (\alpha^{1-N} L_{N}+x)T_N \,
{\cal E}^{1 \dots N \rangle} \, det^{-1}_q(T) =
\\[0.3cm]
= \left( (L_{\widetilde{1}} +x)\,
(L_{\widetilde{2}}+x) \dots (L_{\widetilde{N}}+x) \right)
 {\cal E}^{1 \dots N \rangle} =  {\cal E}^{1 \dots N \rangle} \, Det_q(L;x) \; .
\end{array}
\ee
(according to (\ref{mrtt}) we modify the first line in (\ref{Ldet55})
but it does not affect the final expressions for $Det_q(L;x)$
and coefficients $a_m(L)$). To calculate the commutation relations of $a_m(L)$ with $T^i_j$ we find (by using (\ref{Ldet5}),
(\ref{mrtt}), (\ref{haha''}) and (\ref{nonc}))
\be
\lb{Ldet51}
\begin{array}{l}
 (L_{1} +x)T_1\,  Det_q(L;\alpha x) \, {\cal E}^{2 \dots N+1 \rangle} = \\[0.2cm]
= (L_{1} +x)T_1\,
(\alpha^{-1} L_{2}+x)T_2 \dots (\alpha^{-N} L_{N+1}+x)T_{N+1} \,
{\cal E}^{2 \dots N+1 \rangle} \, \frac{\alpha^N}{{\det}_q(T)} = \\[0.2cm]
= \R_1 \dots \R_N (L_{1} +x)T_1  \dots (\alpha^{-N} L_{N+1}+x)T_{N+1} \,
\R_N^{-1} \dots \R_1^{-1} \, {\cal E}^{2 \dots N+1 \rangle} \,
\frac{\alpha^N}{{\det}_q(T)} =  \\[0.2cm]
= {\cal E}^{2 \dots N+1 \rangle} N^{1 \rangle}_{\langle N+1}
\, Det_q(L;x) \,{\det}_q(T) \bigl(\alpha^{-N} L_{N+1}+x \bigr)T_{N+1} \,
 (N^{-1})^{N+1 \rangle}_{\langle 1} \,
\frac{\alpha^N}{{\det}_q(T)} =  \\[0.2cm]
=  \alpha^N  {\cal E}^{2 \dots N+1 \rangle} N^{1 \rangle}_{\langle N+1}
\, Det_q(L;x) \, \bigl( q^2 (N^{-1} L N)_{_{N+1}}+ x \bigr)
(N^{-1}TN)_{_{N+1}}  \, (N^{-1})^{N+1 \rangle}_{\langle 1}  = \\[0.2cm]
= \alpha^N \, {\cal E}^{2 \dots N+1 \rangle} \, Det_q(L;x)  \,
( q^2 \, L_1 + x) \, T_1 \; ,
\end{array}
\ee
where we have taken into account
the commutation relations of ${\det}_q(T)$ and $L^i_j$
deduced by the standard method
\be
\lb{dqTL}
\begin{array}{c}
{\cal E}_{\langle 1 \dots N} \; {\det}_q(T) \; L_{N+1} =
{\cal E}_{\langle 1 \dots N} \, T_1 \dots T_N \, L_{N+1} = \\[0.2cm]
= \alpha^N \, {\cal E}_{\langle 1 \dots N} \,
\R_N \dots \R_1 L_1 \R_1 \dots \R_N
T_1 \dots T_N  = q^2 \, \alpha^N \; (N^{-1} \, L\, N)_{_{N+1}}
\; {\det}_q(T) \; {\cal E}_{\langle 1 \dots N} \; ,
\end{array}
\ee
(eqs. (\ref{tangb1}) and (\ref{haha'}) were applied). Thus, we have the
following relations (see (\ref{Ldet51}))
\be
\lb{dqTL1}
 (L_{1} +x)T_1 \,  Det_q(L;\alpha x) = \alpha^N  \, Det_q(L;x)  \,
( q^2 \, L_1 + x) \, T_1 \; .
\ee
The expansion of (\ref{dqTL1}) over $x$ gives the recurrent
equation for desired commutation relations of
$a_k(L)$ with $T^i_j$ ($k \geq 0$):
$$
\alpha^{-k} \, L T \, a_k + \alpha^{-1-k} \, T \, a_{k+1} =
q^2 \, a_k \, L T + a_{k+1} T \; , \;\;\; T \, a_0 = a_0 \, T \; .
$$
 These equations are easy to solve by iteration, and the solution is
\be
\lb{dqTL2}
\alpha^{-k} \, T a_{k} = a_{k} T -
(q^2 -1) \sum_{m =1}^{k} (-1)^m a_{k-m} L^{m} \, T \; .
\ee
Since the matrix $T$ is invertible,
we write this equation in the form
\be
\lb{dqTL3}
\alpha^{-k} \, T a_{k} \, T^{-1} = a_{k}  -
(q^2 -1) \sum_{m =1}^{k} (-1)^m a_{k-m} L^{m}  \; .
\ee
For the left hand side of (\ref{dqTL3}),
by using the definition (\ref{sigm1}) of $a_k$, we deduce
$$
\begin{array}{c}
\alpha^{-k} \, q^{-k} \, T_1 a_{k} \, T_1^{-1} =
\alpha^{-k} \, T_1  \, Tr_{D(2 \dots k+1)}
\Bigl( A_{2 \dots k+1} \,  L_{2} \R_2 L_2 \R_2^{-1}
 \dots \R_{k \leftarrow 2}
 L_2 \R_{k \leftarrow 2}^{-1} \Bigr) \, T_1^{-1}  = \\ [0.2cm]
=  Tr_{D(2 \dots k+1)} \left(A_{2 \dots k+1} \, \R_1 \dots \R_k \,
 L_{\underline{1}} \dots L_{\underline{k}}\, \R_k \dots \R_1 \, A_{2 \dots k+1} \right) = \\ [0.2cm]
=  Tr_{D(2 \dots k+1)} \Bigl( \R_{(1 \to k)} \, A_{1 \dots k} \,
 L_{\underline{1}} \dots L_{\underline{k}}  \, A_{1 \dots k} \,
  \R_{(k \leftarrow 1)} \Bigr) = \\ [0.2cm]
=  Tr_{D(2 \dots k+1)} \left(\R_1 \dots \R_{k} \,
 A_{1 \dots k} \,  L_{\underline{1}} \dots L_{\underline{k}}
 \, A_{1 \dots k}  (\R_{k}^{-1} +\lambda) \,
 \R_{(k-1 \leftarrow 1)}\right) =
  \end{array}
$$
\be
\lb{idei}
\begin{array}{c}
=  Tr_{D(2 \dots k)} \left(\R_{(1 \to k-1)}
\left[ Tr_{D(k)} ( A_{1 \dots k}  L_{\underline{1}}
\dots L_{\underline{k}} \,
  A_{1 \dots k} ) \right] \R_{(k-1 \leftarrow 1)} \right) +
\\ [0.3cm]
+ \lambda  \, q^{2(1-k)} \,  Tr_{D(2 \dots k)} \left(
A_{1 \dots k}  L_{\underline{1}} \dots L_{\underline{k}} \right)
 \; ,
\end{array}
\ee
where in the last transformation we apply the first identities
in (\ref{qtrs}) and (\ref{RER}). By repeating this
 transformation for (\ref{idei}) many times we obtain
$$
\alpha^{-k} \, q^{-k} \, T_1 a_{k} \, T_1^{-1}
= q^{-k} a_k + \lambda (1 + q^{-2} + \dots q^{2(1-k)} )
Tr_{D(2 \dots k)} \left(
A_{1 \dots k}  L_{\underline{1}} \dots L_{\underline{k}} \right) =
$$
\be
\lb{left}
 =
 q^{-k} a_k + q (1 - q^{-2k})  Tr_{D(2 \dots k)} \left(
A_{1 \dots k}  L_{\underline{1}}
\dots L_{\underline{k}} \right) \; .
\ee
Comparing (\ref{dqTL3}) and (\ref{left})  we obtain the
remarkable identities for quantum RE matrices $L$
(the so-called Cayley-Hamilton-Newton identities \cite{IOPa})
\be
\lb{chn}
[k]_q \, Tr_{D(2 \dots k)}
(A_{1 \to k} L_{\underline{1}} \dots L_{\underline{k}}) =
- \sum_{m=1}^{k} (-1)^m a_{k-m} L_1^m \; .
\ee
It follows from (\ref{chn}) (apply $Tr_{D(1)}$ to the both sides) that the two basic sets
(\ref{3.1.22}), (\ref{sigm1}) of central elements for RE algebra (defined by Hecke type $R$-matrix) are related by
the
$q$-analogue of the Newton relations:
\be
\lb{Ldet7}
\frac{[k]_q}{q^k} \, a_k +
\sum_{m =1}^{k} (-1)^m a_{k-m} \, p_m = 0 \; , \;\;\;\;\;\;
k=1,...,N \; ,
\ee
where we introduce
power sums $p_m =  Tr_D( L^{m})$, $m = 1, \dots ,N$, and we imply $a_0 =1$.
Note that in view of (\ref{LdetE}), (\ref{Ldeti}) and
(\ref{detD3}) we have
$$
[N]_q \, Tr_{D(2 \dots N)}
(A_{1 \to N} L_{\underline{1}} \dots L_{\underline{N}}) =
[N]_q \, Tr_{D(2 \dots N)}
\bigl(A_{1 \to N} L_{\widetilde{1}} \dots L_{\widetilde{N}}
(y_2 \cdots y_N)^{-1} \bigr) =
$$
$$
 = [N]_q \, Det_q(L) \, q^{(N-1)N} \,  Tr_{D(2 \dots N)}
(A_{1 \to N}) = Det_q(L)\, I_1 \equiv a_N \, I_1 \; .
$$
Thus, for $k=N$ the relation
(\ref{chn}) provides the characteristic identity for the quantum matrix $L$
($q$-analogue of the Cayley-Hamilton theorem):
\be
\lb{Ldet4}
\sum_{k=0}^N \, (-L)^k \, a_{N-k}(L) = 0 \; .
\ee
This identity can formally be obtained by the substitution of $x = -L$ in the characteristic
polynomial (\ref{Ldet3}). Therefore, in view of (\ref{Ldet7}) and (\ref{Ldet4}), the elements $a_m(L)$
can be interpreted as noncommutative analogs of elementary symmetric functions for eigenvalues of
the quantum matrix $L$ (see details in \cite{IsPyaE}).

Introduce generating functions $a(t),p(t)$ for elementary symmetric functions and power sums
$$
a(t) = \sum_{k \geq 0} a_k \, t^k \; , \;\;\; p(t) = \sum_{k \geq 1} p_k \, t^k \; .
$$
Then, it is worth noting \cite{IsBonn}
 that quantum Newton relations (\ref{Ldet7})
 can be written as a finite
difference equation for $a(t)$:
\be
\lb{fdifNR}
a(t) \, p(-t) = \frac{a(q^{-2} t) - a(t)}{q-q^{-1}} \; .
\ee
This equation shows that the power sums $p_k$ can always be expressed as polynomials of the
elementary symmetric functions $a_k$.

The Cayley-Hamilton-Newton identities (\ref{chn}) for the $GL(n)$-type quantum matrix algebra were
invented in \cite{IOPa}. It seems that these matrix identities were unknown even for the case of
usual commutative
matrices ($q=1$). For the reflection equation algebra, in the case
$N=2$, the identity
(\ref{Ldet4}) were considered in
\cite{31} and in \cite{48}. For general $N$ these
identities were proved in \cite{NT}, \cite{PySa}. Newton's relations (\ref{Ldet7}) have been
obtained in \cite{PySa}. Identities (\ref{Ldet51}),
(\ref{dqTL}), (\ref{dqTL1}) and (\ref{dqTL2}) and their
special cases were essentially used in \cite{IsPyaE}
 (see Propositions {\sf 3.21} and {\sf 3.24} there).
Moreover, the spectral properties of the reflection
equation matrices $L^i_j$ were investigated in \cite{IsPyaE}.
In fact, all algebraic relations and identities of this Section
were important for the investigations \cite{IsPyaE}
of the theory of the $q$-deformed isotropic top \cite{alfad},
\cite{23}.

For $GL(m|n)$-type quantum super-matrix algebras Cayley-Hamilton identities  were obtained in
\cite{GPS05}. For  orthogonal and symplectic types quantum matrix algebras the
Cayley-Hamilton identities and Newton relations were derived in  \cite{OgPyaSO}.

\subsection{\bf \em Multi-parameter deformations of linear
 groups\label{multpar}}
\setcounter{equation}0

In this subsection, we consider a multi-parameter deformation of the linear group $GL(N)$ (Refs.
\cite{16,21}, and \cite{24}--\cite{27}). A multi-parameter quantum hyperplane is defined by the
relations
\be
\lb{3.4.1}
x^{i}x^{j} = r_{ij} x^{j}x^{i} \; , \;\; i<j  \; ,
\ee
which can be written in the $R$ -matrix form (\ref{3.3.5}) if we introduce an
additional parameter $q$. Thus, we have $N(N - 1)/2+1$ deformation parameters:
$r_{ij}, \; i<j$ and $q$. The corresponding $R$-matrix is
(see e.g. \cite{16})
\begin{equation}
\label{rmult}
R_{12} =
q \, \sum_i e_{i,i} \otimes e_{i,i} + \sum_{i \neq j}
\,  ( e_{i,i} \otimes e_{j,j} ) a_{ji} + (q - q^{-1})
\sum_{i > j} e_{i,j} \otimes e_{j,i} \; ,
\end{equation}
where $a_{ij} = 1/a_{ji} = r_{ij}/q$ (for $i >j$), and
 it can be represented in components as
\be
\lb{3.4.2}
R^{i_{1},i_{2}}_{j_{1},j_{2}}
= \delta^{i_{1}}_{j_{1}} \delta^{i_{2}}_{j_{2}}
\left( q \delta^{i_{1}i_{2}} +
 \Theta_{i_{2}i_{1}} \frac{q}{r_{i_{1}i_{2}}}    +
 \Theta_{i_{1}i_{2}} \frac{r_{i_{2}i_{1}}}{q} \right) +
 (q-q^{-1}) \delta^{i_{1}}_{j_{2}} \delta^{i_{2}}_{j_{1}}
\Theta_{i_{1}i_{2}} \; ,
\ee
where $\Theta_{ij}$ is defined in (\ref{3.3.6}). The $R$-matrix (\ref{rmult}) is obtained  by the
twisting of the standard one-parameter
$R$-matrix (\ref{3.3.6a}) (see Subsection {\bf \ref{trqH}}
 and eqs. (\ref{2a}), (\ref{twistP}))
\be
\lb{twistgl}
R_{12} \rightarrow F_{21} R_{12} F^{-1}_{12} \Leftrightarrow
\hat{R}_{12} \rightarrow F_{12} \hat{R}_{12} F^{-1}_{12} \; , \;\;\;
F_{12} = \sum_{i,j} ( e_{i,i} \otimes e_{j,j} ) \, f_{ij} \; ,
\ee
where $a_{ij} = f_{ij}/f_{ji}$ and $\hat{F} = P F$
satisfies the twisting matrix conditions (\ref{compat}).
Thus, the multi-parametric $R$-matrix (\ref{3.4.2})
is reduced to the one-parameter $R$-matrix
with the help of the appropriate twisting
(see also Refs. \cite{16} and \cite{27}).

By the construction, in view of the twisting procedure
 (\ref{twistgl}), (\ref{compat}), the $R$-matrix (\ref{3.4.2}) satisfies
 the Yang-Baxter equation (\ref{3.1.3}) and
the same Hecke condition (\ref{3.3.7aa}) as in the one-parameter case.

Now, to justify expression (\ref{3.4.2}),
 we try to find the most general Yang-Baxter solution
 $R_{12}$ of the form (\ref{anzdym}). We only require that the $R$-matrix (\ref{anzdym}) has the lower-triangular
block form: $b_{ij} = 0$ for $i \geq j$ (as it was shown in \cite{18a}
this condition is not restrictive). When we check the
fulfillment of the Yang-Baxter equation, it is convenient to use the diagrammatic technique \cite{IsFirst}
\be
\lb{3.4.3}
\R=\hat{R}^{i_{1},i_{2}}_{j_{1},j_{2}}
= \delta^{i_{1}}_{j_{2}} \delta^{i_{2}}_{j_{1}}
\left( a^{0}_{i_{1}} \delta^{i_{1}i_{2}} +
 \Theta_{i_{2}i_{1}} a^{-}_{i_{1}i_{2}} +
 \Theta_{i_{1}i_{2}} a^{+}_{i_{1}i_{2}} \right) +
 b_{i_{1}i_{2}} \delta^{i_{1}}_{j_{1}} \delta^{i_{2}}_{j_{2}}
\Theta_{i_{2}i_{1}} =
\ee

\unitlength=8mm
\begin{picture}(17,4)
\put(2,3.2){$i_{1}$}
\put(2,3){\line(1,-1){1}}
\put(2,1){\line(1,1){1}}
\put(2,0.5){$j_{1}$}
\put(2.7,2.8){$a^{0}_{i_{1}}$}
\put(2.5,2.5){\line(1,0){1}}
\put(2.6,2.4){\line(1,0){0.8}}
\put(4,0.5){$j_{2}$}
\put(3,2){\line(1,-1){1}}
\put(3,2){\line(1,1){1}}
\put(4,3.2){$i_{2}$}
\put(4.5,1.9){$+$}

\put(6,3.2){$i_{1}$}
\put(6,3){\line(1,-1){1}}
\put(6,1){\line(1,1){1}}
\put(6,0.5){$j_{1}$}
\put(6.6,2.85){$a^{-}_{i_{1}i_{2}}$}
 {\thicklines
\put(6.5,2.5){\line(1,0){0.4}}
\put(7.5,2.5){\vector(-1,0){0.7}}
 }
\put(8,0.5){$j_{2}$}
\put(7,2){\line(1,-1){1}}
\put(7,2){\line(1,1){1}}
\put(8,3.2){$i_{2}$}
\put(8.5,1.9){$+$}

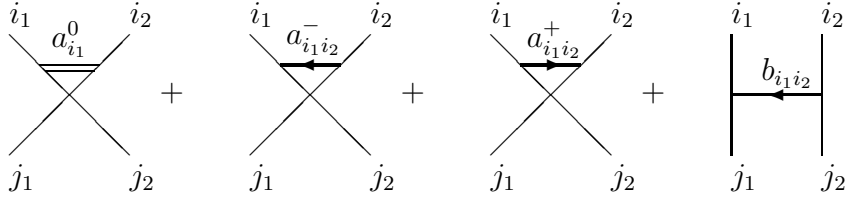
\put(10,3.2){$i_{1}$}
\put(10,3){\line(1,-1){1}}
\put(10,1){\line(1,1){1}}
\put(10,0.5){$j_{1}$}
\put(10.6,2.8){$a^{+}_{i_{1}i_{2}}$}
{\thicklines
\put(10.6,2.5){\vector(1,0){0.6}}
\put(10.5,2.5){\line(1,0){1}}
}
\put(12,0.5){$j_{2}$}
\put(11,2){\line(1,-1){1}}
\put(11,2){\line(1,1){1}}
\put(12,3.2){$i_{2}$}

\put(12.5,1.9){$+$}
 {\thicklines
\put(15.5,2){\vector(-1,0){0.9}}
\put(14,2){\line(1,0){0.7}}
 }
\put(14.5,2.25){$b_{i_{1}i_{2}}$}
\put(14,3.2){$i_{1}$}
\put(14,0.5){$j_{1}$}
\put(14,3){\line(0,-2){2}}
\put(15.5,1){\line(0,2){2}}
\put(15.5,0.5){$j_{2}$}
\put(15.5,3.2){$i_{2}$}
\end{picture}

It turns out that not all solutions of the Yang-Baxter equation
(\ref{3.1.3}) that can be represented in the form (\ref{3.4.3}) are exhausted
by the multi-parameter $R$-matrices (\ref{3.4.2}). Indeed, if we substitute the matrix
(\ref{3.4.3}) into the Yang-Baxter equation (\ref{3.1.3}),
 we obtain the following general
conditions on the coefficients $a^{0}_{i}, \; a^{\pm}_{ij}, \; b_{ij}$:
\be
\lb{3.4.4}
b_{ij} = b \; , \;\; a^{+}_{ij}a^{-}_{ji} = c \; , \;\;
(a^{0}_{i})^{2} - b a^{0}_{i} - c = 0 \;\; (\forall i,j) \; .
\ee
We normalize (\ref{3.4.3}) in such a way
that $c=1$ and choose for convenience, instead of the parameter $b$, a different parameter
$q$, setting $b = q - q^{-1}$. Then $a^{0}_{i}$ can take two
values $\pm q^{\pm 1}$.
For such a normalization, the solution of the Yang-Baxter
equation of the form (\ref{3.4.3}) automatically satisfies the Hecke relation
(\ref{3.3.7aa}). If we set $a^{0}_{i}= q$ (or
$a^{0}_{i}= -q^{-1}$) for all $i$, then we arrive at the
many-parametric case $GL_{q,r_{ij}}(N)$ (\ref{3.4.2}) (up to
exchange $q \rightarrow - q^{-1}$ in the case $a^{0}_{i}= -q^{-1}$). If, however, we set
\be
\lb{3.4.5}
a^{0}_{i}= q \; (1 \leq i \leq M) \; , \;\;\;\;\;\;
a^{0}_{i}= -q^{-1} \;\; ( M+1 \leq i \leq N) \; ,
\ee
then the $R$-matrix (\ref{3.4.3}) does not reduce to (\ref{3.4.2}) and will correspond to a
multi-parameter deformation of the supergroup  $GL(M|N-M)$:
\begin{equation}
\label{srmult}
\R_{12} =
\sum_i (-1)^{[i]} \, q^{1 - 2[i]} \, e_{ii} \otimes e_{ii} +
 \sum_{i \neq j} a_{ij}^{+} \, e_{ij} \otimes e_{ji} +
\lambda \,  \sum_{j > i} \, e_{ii} \otimes e_{jj} \; ,
\end{equation}
where $i,j = 1, \dots , N+M$, $[i] = 0,1$ (mod$(2)$),
we take into account (\ref{3.4.5})
and $a_{ij}^{+} = 1/a_{ji}^{+}$ for $i>j$.
We consider this case (for a special choice of $a_{ij}^{+}$)
 below in Sec. {\bf \ref{qsuper}}.

By virtue of the fulfillment of the Hecke identity (\ref{3.3.7aa}) for the multi-parameter case, we
can introduce the same projectors $\P^{-}$ and $\P^{+}$ as in the one-parameter case
(\ref{3.3.10}); the first of them defines the bosonic quantum hyperplane
(\ref{3.4.1}) (the relations (\ref{3.3.5}) with R-matrix (\ref{3.4.2})), and the second one
defines the fermionic quantum hyperplane
\be
\lb{3.4.6}
\hbox{\bf P}^{+} x_{1}x_{2} = 0 \;\;\;\;\;\; \Leftrightarrow
\;\;\;\;\;\;
(x^{i})^{2} = 0 \; , \;\;\; q^{2} x^{i}x^{j} = - r_{ij}x^{j}x^{i} \;
(i>j) \; .
\ee
Regarding (\ref{3.4.1}) and (\ref{3.4.6}) as comodules for the multi-parameter quantum group
$GL_{q,r_{ij}}(N)$, we find that the generators
$T^{i}_{j}$ of the algebra $Fun(GL_{q,r_{ij}}(N))$ satisfy the same
$RTT$ relations (\ref{3.1.1}) but with $R$-matrix (\ref{3.4.2}). Note, however, that the
quantum determinant $det_{q}(T)$ (\ref{3.3.18}) is not central in the multi-parameter case
\cite{26}. This is due to the fact that in general for the multi-parameter
$R$-matrix we have $N \neq const \cdot I$
in the equations (\ref{matN}), (\ref{nonc}). Therefore, reduction
to the $SL$ case by means of the condition $det_{q}(T)=1$ is possible
only under certain restrictions on the parameters $q,r_{ij}$.
A detailed discussion of these facts can be found in Refs. \cite{26} and \cite{27}.

The algebra (\ref{3.1.20b}), (\ref{3.1.20}) (with the multi-parameter $R$-matrix (\ref{3.5.3}))
which is dual to the algebra $Fun(GL_{q,r_{ij}}(N))$ can also be considered. It appears that this
algebra is isomorphic to the one-parameter deformation of $gl(N)$ (\ref{dj1}) -- (\ref{dj5}). One
can find details about the dual algebras for the special case of $Fun(GL_{q,p}(2))$ in papers \cite{31}, \cite{Dobr}.

\subsection{\bf \em The quantum supergroups $GL_q(N|M)$
and $SL_q(N|M)$\label{qsuper}}
\setcounter{equation}0

We choose the Hecke type $R$-matrix (\ref{3.4.3}),
 (\ref{srmult})
 and write it in the form
 (cf. Refs. \cite{29}, \cite{Wall})
\be
\lb{3.6.1aa}
\R=
\sum_i (-1)^{[i]} \, q^{1 - 2[i]} \, e_{ii} \otimes e_{ii} +
 \sum_{i \neq j} (-1)^{[i][j]}  \, e_{ij} \otimes e_{ji} +
\lambda \,  \sum_{j > i} \, e_{ii} \otimes e_{jj} \; ,
\ee
where
we have set (see (\ref{srmult}))
$$
a^{0}_{i} = (-1)^{[i]}q^{1-2[i]} \; , \;\; a^{+}_{ij} = (a^{-}_{ij})^{-1} =
(-1)^{[i][j]} \; , \;\;
b = q-q^{-1} = \lambda \; .
$$
We stress here that the matrix units $e_{ij}$ and tensor products
in (\ref{3.6.1aa}) {\it are not graded}, as follows form
the previous Section {\bf \ref{multpar}}.
The component presentation of (\ref{3.6.1aa}) is
\be
\lb{3.6.1}
\begin{array}{c}
\hat{R}^{i_{1},i_{2}}_{j_{1},j_{2}}
=
\delta^{i_{1}}_{j_{2}} \delta^{i_{2}}_{j_{1}} \, (-1)^{[i_1][i_2]} \,
q^{\delta_{_{i_{1}i_{2}}} \, (1-2[i_{_{1}}])}
+ \delta^{i_{1}}_{j_{1}} \delta^{i_{2}}_{j_{2}} \, \lambda \,
\Theta_{i_{2}i_{1}} \; .
\end{array}
\ee
Thus, the parameters $a^{0}_{i}$ take the two values $\pm q^{\pm 1}$
and, as we assumed it in Sec. {\bf \ref{multpar}},
the $R$-matrix (\ref{3.6.1aa}), (\ref{3.6.1}) must correspond to some
supergroup. Indeed, in the limit $q \rightarrow 1$,
we find that $\R$ tends to the supertransposition operator:
\be
\lb{stran}
\R^{i_{1}i_{2}}_{j_{1}j_{2}} \rightarrow (-1)^{[i_1][i_2]}
\delta^{i_{1}}_{j_{2}} \delta^{i_{2}}_{j_{1}} \equiv {\cal P}_{12} \; .
\ee
Suppose that the $R$-matrix acts in the space of
the direct product $x \underline{\otimes} \, y$
of two supervectors
 $x$ and $y$ with coordinates
 $x^{j_{1}}$ and $y^{j_{2}}$,
 and $[i] = 0,1$ denotes the parity (grading)
 of the components\footnote{There are two equivalent descriptions of super vector spaces $V$. The first one is to consider the graded basis vectors $e_i$ while coordinates $x^i$ of super vectors $e_i x^i\in V$
  are ordinary numbers. Another (dual) approach is that vectors $e_i$ form a bases of an ordinary vector space, but coordinates $x^i$  are graded in such a way that $e_i x^i$ belongs to the superspace $V$.
  Here we use the second approach.}
  $x^{i}$ and $y^{i}$. According to
(\ref{stran}), we write the condition
for the graded tensor product $\underline{\otimes}$ as
$$
x^{j_{1}} \underline{\otimes} \, y^{j_{2}} =
{\cal P}^{j_{1}j_{2}}_{k_{1}k_{2}}
(1 \underline{\otimes} \, y^{k_{1}}) (x^{k_{2}} \underline{\otimes} \, 1) \;\;\; \Rightarrow \;\;\;
x^{j_{1}} \underline{\otimes} \, y^{j_{2}} = (-1)^{[j_{1}][j_{2}]}
(1 \underline{\otimes} \, y^{j_{2}}) (x^{j_{1}} \underline{\otimes} \, 1) \; .
$$
For definiteness, we will assume that
\be
\lb{pari}
[i] =0 \;\;\; (1 \leq i \leq N) \; , \;\;\;
[i] =1 \;\;\; (N+1 \leq i \leq N+M) \; .
\ee

As we noted in Sec. {\bf \ref{multpar}},
the $R$-matrix (\ref{3.6.1}) satisfies
 Yang-Baxter equation
(\ref{3.1.3}) (in the braid group form)
and the Hecke relation (\ref{3.3.7}).
In addition to the matrix
$\R$, we introduce the new $R$-matrix
$$
R_{12} = {\cal P}_{12} \R_{12} = (-)^{(1)(2)} P_{12} \R_{12} =
$$
\be
\lb{newSR}
= \sum_i  q^{1 - 2[i]} \, e_{ii} \otimes e_{ii} +
 \sum_{i \neq j}  e_{ii} \otimes e_{jj} +
\lambda \,  \sum_{i > j} \, (-1)^{[i][j]}  \, e_{ij} \otimes e_{ji} \; .
\ee
with the semiclassical behavior (\ref{3.2.1}). Here and
below we use notation
\be
\lb{matnot}
\left( (-)^{(1)(2)} \right)^{i_1 i_2}_{j_1 j_2} : =
(-1)^{[i_1][i_2]} \, \delta^{i_1}_{j_1} \, \delta^{i_2}_{j_2} \; .
\ee
Then, for the new $R$-matrix (\ref{newSR}),
 we obtain from equation (\ref{3.1.3}) the graded form
\cite{32} of the Yang-Baxter equation:
\be
\lb{3.6.2}
R_{12}(-)^{(2)(3)}R_{13}(-)^{(2)(3)}R_{23}  =
R_{23}(-)^{(2)(3)}R_{13}
(-)^{(2)(3)}R_{12} \; ,
\ee
where we have taken into account the fact that
$R_{12}$ is an even $R$-matrix, i.e.,
$$
R^{i_{1}i_{2}}_{j_{1}j_{2}} \neq 0 \;\; \hbox{if} \;\; [i_1] + [j_1] +[i_2] +[j_2] = 0 \;
({\rm mod}(2)) \; \Rightarrow
$$
$$
 (-1)^{[i_3]([i_1] +[i_2])} R^{i_{1}i_{2}}_{j_{1}j_{2}} \delta^{i_3}_{j_3} =
R^{i_{1}i_{2}}_{j_{1}j_{2}} \delta^{i_3}_{j_3} \,
(-1)^{[j_3]([j_1] +[j_2])} \; \Leftrightarrow
$$
$$
(-)^{(3)((1) +(2))} \, R_{12} \, I_3 = R_{12} \, I_3 \, (-)^{(3)((1) + (2))} \; .
$$
In the last relation we set
$$
\left( (-)^{(3)((1) +(2))} \right)^{i_{1}i_{2}i_3}_{j_{1}j_{2}j_3} =
(-1)^{[i_3] ([i_1]+[i_2])} \,
\delta^{i_1}_{j_1} \, \delta^{i_2}_{j_2} \, \delta^{i_3}_{j_3} \; .
$$
For the $R$-matrix (\ref{3.6.1aa}), (\ref{3.6.1}) we will also use the properties:
\be
\lb{3.6.5ra}
\hat{R}_{12}  =  (-)^{(1)(2)}\hat{R}_{12} (-)^{(1)(2)} \; , \;\;
(-)^{(1)+(2)} \hat{R}_{12}  =  \hat{R}_{12} (-)^{(1)+(2)} \; .
\ee
 Since $GL_q(N|M)$ $R$-matrix (\ref{3.6.1}) satisfies the
Hecke condition, we find
$$
(\hat{R}^{-1})^{i_{1},i_{2}}_{j_{1},j_{2}}
= \hat{R}^{i_{1},i_{2}}_{j_{1},j_{2}} - \lambda \, \delta^{i_{1}}_{j_{1}} \delta^{i_{2}}_{j_{2}}
 = \delta^{i_{1}}_{j_{2}} \delta^{i_{2}}_{j_{1}} \, (-1)^{[i_1][i_2]} \,
q^{\delta_{_{i_{1}i_{2}}} \, (2[i_{_{1}}]-1)} - \lambda \,
 \delta^{i_{1}}_{j_{1}} \delta^{i_{2}}_{j_{2}} \, \Theta_{i_{1}i_{2}} \; ,
$$
and we have the identities (cf. (\ref{3.3.7a1q}))
\be
\lb{3.6.5rr}
\hat{R}_{12}^{-1}[q^{-1}] = \hat{R}_{21}[q] \; .
\ee
 Finally, the skew-inverse matrix $\Psi_{12}$ (\ref{skew})
  for the
$GL_q(N|M)$ R-matrix, defined in (\ref{3.6.1aa}), (\ref{3.6.1}),
 has the form
\be
\lb{skglnm}
\begin{array}{c}
\hat{\Psi}_{12} = \sum_{i} e_{ii} \otimes e_{ii} (-1)^{[i]} q^{2[i]-1} +
\sum_{i \neq j} (-1)^{[i][j]} e_{ij} \otimes e_{ji} - \\ [0.2cm]
- \lambda \sum_{i < j} e_{ii} \otimes e_{jj} (-1)^{[i]+[j]}
 q^{(-1)^{^{[i]}} (2 i - 2N -1) -(-1)^{^{[j]}} (2 j - 2N -1)}  \; , \\ \\
\hat{\Psi}^{i_1 i_2}_{j_1 j_2} = (-1)^{^{[i_1][i_2]}} \,
q^{\delta_{_{i_{1}i_{2}}} \, (2[i_{_{1}}] -1)} \,
\delta^{i_1}_{j_2} \delta^{i_2}_{j_1}  - \\ [0.2cm]
 - (-1)^{^{[i_1] +[i_2]}} \lambda  \,
 q^{(-1)^{^{[i_1]}} (2 i_1 - 2N -1)} \, q^{(-1)^{^{[i_2]}}(1+ 2N -2 i_2) )} \,
\Theta_{i_2 i_1} \, \delta^{i_1}_{j_1} \delta^{i_2}_{j_2}  \; ,
\end{array}
\ee
which follows from the general formula (\ref{skewdym})
(for the case $[i_1]=[i_2]=0$ we reproduce matrix (\ref{skewdym1})).
The corresponding
matrices of quantum supertraces are
\be
\lb{dqsup}
\begin{array}{c}
D_1 \equiv  Tr_{2}\left( \hat{\Psi}_{12} \right) \;\; \Rightarrow  \;\;\;
D^i_j = (-1)^{^{[i]}} \, q^{2M + (-1)^{^{[i]}}(2i-2N-1)} \, \delta^i_j \; , \\ \\
Q_2 \equiv  Tr_{1}\left( \hat{\Psi}_{12} \right)
\;\; \Rightarrow  \;\;\;
Q^i_j = (-1)^{^{[i]}} \, q^{-2M + (-1)^{^{[i]}}(2N+1-2i)} \delta^i_j \; , \\ \\
Tr(D) = Tr(Q) = (1 - q^{2(M-N)})/\lambda  = q^{(M-N)} \, [N-M]_q \; .
\end{array}
\ee
Note that the quantum supertrace $Tr_D$, which is constructed by means of the matrix $D$ (see the first line in
(\ref{dqsup})), coincides up to the factor $q^{(3M-N)/2}$ with the quantum trace presented in
\cite{31}. For $q \rightarrow 1$ quantum supertraces $Tr_{\!_D}$ and
 $Tr_{\!_Q}$ tend to the usual supertraces.

The quantum multidimensional superplanes for the Hecke type $R$-matrix (\ref{3.6.1}) are defined as
algebras ${\cal V}_{\pm}$ with generators $x_i$ $(i=1,\dots,N+M)$ and defining relations
(see, for example, Refs. \cite{21}, \cite{30} and \cite{31}):
\be
\lb{3.6.3}
\begin{array}{l}
\!\!\! {\cal V}_{-}: \;  (\R - q)x^{1\rangle}x^{2\rangle} = 0
\;\;\; \Leftrightarrow \;\;\; x^{i}x^{j} = (-1)^{[i][j]} q
x^{j}x^{i}   \;\; (i<j) \; ,
\;  (x^{i})^{2} = 0 \; \hbox{if} \; [i]=1 \; , \\ \\
\!\!\! {\cal V}_{+}: \; (\R + q^{-1})x^{1\rangle}x^{2\rangle} = 0
\;\; \Leftrightarrow \;\; q x^{i}x^{j} =
-(-1)^{[i][j]} x^{j}x^{i}   \;\; (i<j)  \; ,
\;  (x^{i})^{2} = 0 \; \hbox{if} \; [i]=0 \; .
\end{array}
\ee
The super-hyperplane ${\cal V}_{+}$ can be interpreted (see, e.g., \cite{Isaev1}) as an exterior algebra
of differentials $dx^{i}$ of the coordinates $x^{i}$ for the first hyperplane ${\cal V}_{-}$.

We take the left coaction (\ref{3.3.1bb}) of the quantum supergroup
with generators $T^i_j$ to the quantum
superspaces ${\cal V}_{\pm}$, defined in (\ref{3.6.3}),
 and consider this coaction to the
spaces ${\cal V}_{\pm}\, \underline{\otimes} \, {\cal V}_{\pm}$:
\be
\lb{3.6.4de}
(T^{i_1}_{j_1} \underline{\otimes} \, x^{j_1}) \, (T^{i_2}_{j_2} \underline{\otimes} \, x^{j_2}) =
(-1)^{[j_1]([i_2] + [j_2])} \,
(T^{i_1}_{j_1} \, T^{i_2}_{j_2}) \underline{\otimes} \, (x^{j_1} \, x^{j_2}) \; ,
\ee
where $\underline{\otimes} \,$ is understood as a graded direct product.
We postulate the gradings of the elements $T^{i}_{j}$ and $x^i$ as
$[T^{i}_{j}] = [i]+[j]$ and $[x^i]=[i]$. In \cite{32}, the right
coaction of the quantum supergroup
was considered with another signs in the formulas,
but it can be shown that this difference is not essential.

{}From the condition of covariance of the relations
(\ref{3.6.3}) under coaction (\ref{3.6.4de}),
we deduce the graded form of the $RTT$ equations:
\be
\lb{3.6.4d}
\R^{i_1 i_2}_{k_1 k_2} \,
T^{k_1}_{j_1}  \, (-1)^{[j_1][k_2]} \, T^{k_2}_{j_2} \, (-1)^{[j_1][j_2]} =
T^{i_1}_{k_1}  \, (-1)^{[k_1][i_2]} \, T^{i_2}_{k_2} \, (-1)^{[k_1][k_2]} \,
\R^{k_1 k_2}_{j_1 j_2} \; ,
\ee
written, with the help of the
 concise matrix notation (\ref{matnot}), as (cf. (\ref{qmam}))
\be
\lb{3.6.4}
\begin{array}{c}
\R T_{1} (-)^{(1)(2)} T_{2} (-)^{(1)(2)} =
T_{1} (-)^{(1)(2)} T_{2} (-)^{(1)(2)} \R \Leftrightarrow \\ \\
R_{12} T_{1} (-)^{(1)(2)} T_{2} (-)^{(1)(2)} =
(-)^{(1)(2)} T_{2} (-)^{(1)(2)} T_{1} R_{12} \; ,
\end{array}
\ee
and in the component form
 (we use the one-parametric
 $R$-matrix (\ref{3.6.1}); the multi-parametric case was considered in \cite{21})
  we have
\be
\lb{3.6.4f}
\begin{array}{c}
T^{i_1}_{j_1} \, T^{i_2}_{j_2} - (-1)^{([i_1] + [j_1])([i_2] + [j_2])}
\,  T^{i_2}_{j_2} \, T^{i_1}_{j_1}
= (q-q^{-1}) \, (-1)^{([j_2][i_2] + [j_1]([i_2] + [j_2])}
\, T^{i_1}_{j_2} \, T^{i_2}_{j_1} \\ [0.2cm] (i_1 < i_2 \, , \;\; j_1 < j_2) \; , \\ [0.2cm]
T^{i_1}_{j_1} \, T^{i_2}_{j_2} = (-1)^{([i_1] + [j_1])([i_2] + [j_2])}
\,  T^{i_2}_{j_2} \, T^{i_1}_{j_1}  \;\;\; (i_1 < i_2 \, , \;\; j_1 > j_2) \; , \\ [0.2cm]
T^{i_1}_{j_1} \, T^{i_2}_{j_1} = (-1)^{[i_1][i_2]} \, q \,
\,  T^{i_2}_{j_1} \, T^{i_1}_{j_1}  \;\; ([j_1]=0 \, , \; i_1 < i_2) \; ,  \\ [0.2cm]
T^{i_1}_{j_1} \, T^{i_2}_{j_1} = (-1)^{([i_1]+1)([i_2]+1)} \, q^{-1} \,
\,  T^{i_2}_{j_1} \, T^{i_1}_{j_1}  \;\; ([j_1]=1 \, , \; i_1 < i_2) \; ,  \\ [0.2cm]
T^{i_1}_{j_1} \, T^{i_1}_{j_2} = (-1)^{[j_1][j_2]} \, q \,
\,  T^{i_1}_{j_2} \, T^{i_1}_{j_1}  \;\; ([i_1]=0 \, , \; j_1 < j_2) \; , \\ [0.2cm]
T^{i_1}_{j_1} \, T^{i_1}_{j_2} = (-1)^{([j_1]+1)([j_2]+1)} \, q^{-1} \,
\,  T^{i_1}_{j_2} \, T^{i_1}_{j_1}  \; ([i_1]=1 \, , \; j_1 < j_2) \; ,  \\ [0.2cm]
(T^{i_1}_{j_1})^2=0 \;\; ([i_1] \neq [j_1]) \; .
\end{array}
\ee
Relations (\ref{3.6.4d}) -- (\ref{3.6.4f}) are the defining relations for the generators
$T^{i}_{j}$ of the graded quantum algebra $Fun(GL_q(N|M))$.

By using (\ref{3.2.1}), the semiclassical analog of (\ref{3.6.4}) can readily be deduced
$$
\left\{T_1 \, , \; (-)^{(1)(2)} T_2 (-)^{(1)(2)} \right\} =
[ T_1 (-)^{(1)(2)} T_2 (-)^{(1)(2)} , \, r_{12} ] \; ,
$$
or in the component form we have
$$
\begin{array}{c}
(-1)^{[j_1]([i_2] +[j_2])} \,
\left\{ T^{i_1}_{j_1} , \,  T^{i_2}_{j_2} \right\}  = \\ \\
= T^{i_1}_{k_1}  \, (-1)^{[k_1][i_2]} \, T^{i_2}_{k_2} \, (-1)^{[k_1][k_2]} \,
r^{k_1 k_2}_{j_1 j_2}  -   r^{i_1 i_2}_{k_1 k_2} \,
T^{k_1}_{j_1}  \, (-1)^{[j_1][k_2]} \, T^{k_2}_{j_2} \, (-1)^{[j_1][j_2]}
\; ,
\end{array}
$$
where $\{.,.\}$ denotes the Poisson super-brackets
$$
\{ T^{i_1}_{j_1} , \, T^{i_2}_{j_2} \} = - (-1)^{([i_1] + [j_1])([i_2] + [j_2])} \,
\{ T^{i_2}_{j_2} , \,  T^{i_1}_{j_1} \} \; .
$$

The matrix $|| T^{i}_{j} ||$ is represented in the block form
\begin{equation}
\lb{3.6.5}
T^{i}_{j} =
 \left(
\begin{tabular}{c|c}
$A $        &     $ B $ \\  \hline \\ [-0.4cm]
$C $ &     $ D  $    \\
\end{tabular}
\right)
\end{equation}
where the elements of the $N \times N$ matrix
$||A^r_s||$ and of the $M \times M$ matrix
$||D^\alpha_\beta||$ form the algebras $Fun(GL_{q}(N))$ and $Fun(GL_{q^{-1}}(M))$, respectively.
Indeed, from (\ref{3.6.4f}) we have
\be
\lb{rAArDD}
\hat{R}^{i_1 i_2}_{k_1 k_2}[q] A^{k_1}_{\; j_1} \, A^{k_2}_{\; j_2} =
A_{\; k_1}^{i_1} \, A_{\; k_2}^{i_2} \hat{R}_{j_1 j_2}^{k_1 k_2}[q] \; , \;\;
\hat{R}^{\alpha_1 \alpha_2}_{\gamma_1 \gamma_2}[q^{-1}] D^{\gamma_1}_{\; \beta_1} \, D^{\gamma_2}_{\; \beta_2} =
D_{\; \gamma_1}^{\alpha_1} \, D_{\; \gamma_2}^{\alpha_2} \hat{R}_{\beta_1 \beta_2}^{\gamma_1
\gamma_2}[q^{-1}] \, ,
\ee
where $\hat{R}^{i_1 i_2}_{k_1 k_2}[q]$ and $\hat{R}^{\alpha_1 \alpha_2}_{\gamma_1
\gamma_2}[q^{-1}]$ are standard $Fun(GL_{q}(N))$ and $Fun(GL_{q^{-1}}(M))$
$R$-matrices defined in (\ref{3.3.6}). We assume that the quantum matrices
$||A^r_s||$ and $||D^\alpha_\beta||$ are invertible. It means that the algebra
$Fun(GL_{q}(N|M))$ should be extended by the elements ${\det}_q^{-1}(A)$ and
${\det}_{q^{-1}}^{-1}(D)$ (see (\ref{explT-1})
and Definition {\bf \em \ref{def11}}). In this case, from
(\ref{rAArDD}) and (\ref{3.3.7a1q}) we obtain
$$
\begin{array}{c}
\hat{R}^{\alpha_1 \alpha_2}_{\gamma_1 \gamma_2}[q] (D^{-1})^{\gamma_1}_{\; \beta_1}
\, (D^{-1})^{\gamma_2}_{\; \beta_2} =
(D^{-1})_{\; \gamma_1}^{\alpha_1} \, (D^{-1})_{\; \gamma_2}^{\alpha_2}
\hat{R}_{\beta_1 \beta_2}^{\gamma_1 \gamma_2}[q] \; , \\ [0.2cm]
\hat{R}^{i_1 i_2}_{k_1 k_2}[q^{-1}] (A^{-1})^{k_1}_{\; j_1} \, (A^{-1})^{k_2}_{\; j_2} =
(A^{-1})_{\; k_1}^{i_1} \, (A^{-1})_{\; k_2}^{i_2} \hat{R}_{j_1 j_2}^{k_1 k_2}[q^{-1}]  \; .
\end{array}
$$

 For the elements of the rectangular matrices
$||B^r_\beta||$ and $||C^\alpha_s||$, we obtain from (\ref{3.6.4f})
the commutation relations:
$$
\begin{array}{c}
\hat{R}^{i_1 i_2}_{k_1 k_2}[q] B^{k_1}_{\; \alpha_1} \, B^{k_2}_{\; \alpha_2} =
- B_{\; \beta_1}^{i_1} \, B_{\; \beta_2}^{i_2} \hat{R}^{\beta_1 \beta_2}_{\alpha_1
\alpha_2}[q^{-1}] \; , \;\; B^{i_1}_{\; \alpha_1} C_{\; i_2}^{\alpha_2} = -
 C_{\; i_2}^{\alpha_2}  B^{i_1}_{\; \alpha_1} \; , \\ [0.2cm]
 -\hat{R}^{\alpha_1 \alpha_2}_{\beta_1 \beta_2}[q^{-1}] C_{\; j_1}^{\beta_1}
\, C_{\; j_2}^{\beta_2}
= C_{\; i_1}^{\alpha_1}  C_{\; i_2}^{\alpha_2} \hat{R}^{i_1 i_2}_{j_1 j_2}[q] \, ,
\\ [0.2cm]
 \hat{R}^{\alpha_1 \alpha_2}_{\beta_1 \beta_2}[q] C_{\; j_2}^{\beta_2}
\, D_{\; \gamma_1}^{\beta_1}
= P^{\alpha_1 \alpha_2}_{\beta_1 \beta_2} D^{\beta_1}_{\;  \gamma_1}  C_{\; j_2}^{\beta_2}  \, ,
\;\;
 B^{i_2}_{\; \beta_2} \, D^{\alpha_1}_{\;  \beta_1} P_{\gamma_1 \gamma_2}^{\beta_1 \beta_2}
=  D^{\alpha_1}_{\;  \beta_1}  B^{i_2}_{\; \beta_2}
\hat{R}_{\gamma_1 \gamma_2}^{\beta_1 \beta_2}[q^{-1}] \, , \\ [0.2cm]
 A_{\; i_2}^{k_2}
\, C_{\; i_1}^{\alpha_1}  \hat{R}^{i_1 i_2}_{j_1 j_2}[q^{-1}]
=  C^{\alpha_1}_{\;  i_1}  A_{\; i_2}^{k_2} P^{i_1 i_2}_{j_1 j_2} \, ,
\;\;
\hat{R}^{i_1 i_2}_{j_1 j_2}[q^{-1}]  A^{j_2}_{\; k_2} B^{j_1}_{\; \beta_1}
=  P^{i_1 i_2}_{j_1 j_2} B^{j_1}_{\; \beta_1}  A^{j_2}_{\; k_2} \, , \\ [0.2cm]
 A_{\; j}^{i} D^\alpha_\beta - D^\alpha_\beta  A_{\; j}^{i} =
 (q-q^{-1}) C^\alpha_{\; j} B^i_{\; \beta} \; .
\end{array}
$$
By using these relations, one can prove that the elements of the matrix $X=(A-BD^{-1}C)$ satisfy the $RTT$
commutation relations
\be
\lb{3.6.6X}
\hat{R}^{i_1 i_2}_{k_1 k_2}[q] X^{k_1}_{\; j_1} \, X^{k_2}_{\; j_2} =
X_{\; k_1}^{i_1} \, X_{\; k_2}^{i_2} \hat{R}_{j_1 j_2}^{k_1 k_2}[q] \; .
\ee
It means that elements $X^i_j$ $(i,j=1,\dots,N)$ generate a subalgebra
$Fun(GL_q(N))$ in $Fun(GL_{q}(N|M))$. Assume that the quantum matrix
$X=A-BD^{-1}C$ is also invertible. Then the same is valid for the matrix
$||T^{i}_{j}||$, as it follows from the Gauss decomposition
\begin{equation}
\lb{3.6.6}
T = \left(
\begin{tabular}{c|c}
$ A $        &     $ B $   \\  \hline \\ [-0.4cm]
$ C $        &     $ D $    \\
\end{tabular}
\right) =
 \left(
\begin{tabular}{c|c}
$1  $        &     $BD^{-1}\!\!$ \\  \hline \\ [-0.4cm]
$0 $ &     $ 1  $    \\
\end{tabular}
\right)
 \left(
\begin{tabular}{c|c}
$X$        &     $ 0 $ \\  \hline \\ [-0.4cm]
$0 $ &     $ D $    \\
\end{tabular}
\right)
 \left(
\begin{tabular}{c|c}
$1  $        &     $ 0$ \\  \hline \\ [-0.4cm]
$D^{-1}C $ &     $ 1  $    \\
\end{tabular}
\right) \; ,
\end{equation}
and we have
\begin{equation}
\lb{3.6.6T1}
T^{-1} =  \left(
\begin{tabular}{c|c}
$X^{-1}$        &     $- X^{-1} \, B D^{-1}$   \\  \hline \\ [-0.4cm]
$-D^{-1} C \, X^{-1} $        &     $ Y $    \\
\end{tabular}
\right) =
\left(
\begin{tabular}{c|c}
$X^{-1}$        &     $ - A^{-1} B Y $   \\  \hline \\ [-0.4cm]
$- Y \, C  A^{-1}  $        &     $ Y $    \\
\end{tabular}
\right)
 \; ,
\end{equation}
where  $Y=D^{-1}CX^{-1}BD^{-1} + D^{-1}=(D-CA^{-1}B)^{-1}$.
By inverting relations (\ref{3.6.4}), we obtain that
elements of the inverse matrix $T^{-1}$ satisfy
\be
\lb{sT1T1}
\R_{12} \, (-)^{(1)(2)} T_{2}^{-1} (-)^{(1)(2)} T_{1}^{-1} =
(-)^{(1)(2)} T_{2}^{-1} (-)^{(1)(2)} T_{1}^{-1} \, \R_{12} \; .
\ee
 Now we note that eq.
(\ref{3.6.6X}) is simply obtained from (\ref{rAArDD}) and (\ref{sT1T1}) if
we take into account
(\ref{3.6.5rr}) and (\ref{3.6.6T1}).

In view of the existing of the inverse element (\ref{3.6.6T1}),
the algebra $Fun(GL_{q}(N|M))$ with defining relations
(\ref{3.6.4}) is a Hopf algebra with usual structure mappings (\ref{3.1.5}):
$$
\Delta (T^{i}_{k}) = T^{i}_{j} \underline{\otimes} \, T^{j}_{k} \ , \;\;
\epsilon(T^{i}_{j}) = \delta^{i}_{j} \  , \;\;
 S (T^{i}_{j}) = (T^{-1})^{i}_{j} \ ,
$$
where the product $\underline{\otimes}$
in the definition of $\Delta$ is understood as the graded
direct product.

We define dual quantum multidimensional superplanes as a
$Fun(GL_q(N|M))$-comodule algebra  ${\cal V}^*$ with generators $y_i$
$(i=1,2,\dots,N+M)$ and left coaction (cf. (\ref{3.3.1bb}))
\be
\lb{3.6.4df}
y_i \to \delta_T(y_i) =
(1 \underline{\otimes} \,  y_{j}) \,
((T^{-1})^{j}_{i} \, \underline{\otimes} \, 1) \equiv
(-1)^{([i]+[j])[j]} \;  ((T^{-1})^{j}_{i}
\, \underline{\otimes} \,  y_{j}) \; .
\ee
This coaction is such that the pairing
\be
\lb{brst}
Q = (y_i \, x^i) \; ,
\ee
 is a co-invariant element $\delta_T(Q) = 1 \underline{\otimes} \, Q$
 if the
  generators $x^i$ of the algebras ${\cal V}_{\pm}$ (\ref{3.6.3}) are
transformed according to (\ref{3.3.1bb}). Assume that the grading of the coordinate $y_i$ is opposite to the grading of $x^i$, i.e.,
$[y_i] = [i]+1$.
Then the dual algebras ${\cal V}^*$, which are covariant under the transformations (\ref{3.6.4df}),
 has the following defining relations (cf. (\ref{3.6.3}))
\be
\lb{dualss}
{\cal V}^*_{-}:
\;\; y_{\langle 2}y_{\langle 1} \left( \R_{12}^{\; '}  - q \right)=  0 \; , \;\;\;
{\cal V}^*_{+}:
\;\; y_{\langle 2}y_{\langle 1} \left( \R_{12}^{\; '}  + q^{- 1} \right)=  0 \; ,
\ee
where we have used new Hecke type Yang-Baxter $R$-matrix:
$\R_{12}^{\; '} = (-)^{(1)} \R_{12}(-)^{(1)}$. We check directly the covariance of relations (\ref{dualss})
under coaction (\ref{3.6.4df}):
$$
y_{\langle 2}y_{\langle 1} \left( \R_{12}^{\; '}  \pm q^{\mp 1} \right) \to
 y_{\langle 2} T_2^{-1} \, y_{\langle 1} \, T_1^{-1} \, \left( \R_{12}^{\; '}  \pm q^{\mp 1}
 \right) =
 $$
 $$
 = y_{\langle 2} \, y_{\langle 1} (-)^{((1)+1)(2)} T_2^{-1}
 (-)^{((1)+1)(2)} T_1^{-1} \, \left( \R_{12}^{\; '}  \pm q^{\mp 1}
 \right) =
$$
 $$
 = y_{\langle 2} \, y_{\langle 1} (-)^{(1)(2)+(2)} T_2^{-1}
 (-)^{(1)(2)} T_1^{-1} \, \left( \R_{12}  \pm q^{\mp 1}
 \right)(-)^{(2)} =
$$
 $$
 = y_{\langle 2} \, y_{\langle 1} (-)^{(1)(2)+(2)}  \left( \R_{12}  \pm q^{\mp 1}
 \right) T_2^{-1}  (-)^{(1)(2)} T_1^{-1} \, (-)^{(2)} =
$$
 $$
 = y_{\langle 2} \, y_{\langle 1} \left( (-)^{(1)}\R_{12} (-)^{(1)}  \pm q^{\mp 1}
 \right) (-)^{((1)+1)(2)}   T_2^{-1}  (-)^{((1)+1)(2)} T_1^{-1}  \; .
$$
Here we have used concise notation
$(1 \underline{\otimes} \,  y_{j}) \, ((T^{-1})^{j}_{i} \
 underline{\otimes} \, 1) \equiv y_{j}(T^{-1})^{j}_{i}$.
In the component form, eqs. (\ref{dualss}) are
$$
\begin{array}{l}
{\cal V}^*_{-}: \;\;
 q \, y_{i} \, y_{j} =
-(-1)^{([i]+1)([j]+1)} y_{j}\, y_{i}   \;\; (i<j) \; ,
\;\;\;\;   (y_{i})^{2} = 0 \;\; \hbox{if} \;\; [i]=1 \; ,   \\ [0.2cm]
{\cal V}^*_{+}: \;\;
y_{i}\, y_{j} = (-1)^{([i]+1)([j]+1)} q \, y_{j} \, y_{i}   \;\; (i<j)  \; ,
\;\;\;\;  (y_{i})^{2} = 0 \;\; \hbox{if} \;\; [i]=0 \; .
\end{array}
$$
Covariant (with respect to coactions
(\ref{3.3.1bb}),
(\ref{3.6.4df})) cross-commutation relations for generators
 $x^i \in {\cal V}_{\pm}$ and
$y_j \in {\cal V}_{\pm}^*$ are:
$$
x^{2 \rangle} \, y_{\langle 2} = (-)^{(2)} \, y_{\langle 1} \,  \R_{12} \, x^{1 \rangle} \; .
$$
Using these relations we define covariant algebras ${\cal V}_{\pm} \, \sharp \, {\cal V}_{\pm}^{*}$
and ${\cal V}_{\pm} \, \sharp \, {\cal V}_{\mp}^{*}$ which are the cross products of algebras
${\cal V}_{\pm}$ and ${\cal V}_{\pm}^{*}$. For $q^2 \neq -1$
 one can easily check  that the element
$Q \in {\cal V}_{\mp} \, \sharp \, {\cal V}_{\pm}^*$
(defined in (\ref{brst}) and having the grading $[Q]=1$)
 satisfies $Q^2=0$. Let $d$:
 ${\cal V}_{-} \, \sharp \, {\cal V}_{+}^* \to
{\cal V}_{-} \, \sharp \, {\cal V}_{+}^*$ be the linear map
$d (f) = f \cdot Q$, where
$f \in {\cal V}_{-} \, \sharp \, {\cal V}_{+}^*$, and
we have $d^2 (f) =0$. Put
$H({\cal V}_{-} \, \sharp \, {\cal V}_{+}^*)={\rm Ker}(d)/{\rm Im}(d)$.
The map $d$ defines the structure of the Koszul complex on
 ${\cal V}_{-} \, \sharp \, {\cal V}_{+}^*$.


\begin{proposition}\label{prop6}
(see \cite{21}). {\it $H({\cal V}_{-} \, \sharp \, {\cal V}_{+}^*)$ is an
one-dimensional subspace generated by
 \be
 \lb{sdetf}
\prod_{[i]=0} y_i \, \prod_{[j]=1} x^j \;\;\;
{\rm mod} \; ({\rm Im}(d)) \; ,
 \ee
and
\be
\lb{sdetb}
\delta_T (\prod_{[i]=0} y_i \, \prod_{[j]=1} x^j) =
{\rm sdet}^{-1}_q(T) \, \underline{\otimes} \,
\prod_{[i]=0} y_i \, \prod_{[j]=1} x^j \;\;\; {\rm mod}
\; ({\rm Im}(d)) \; ,
\ee
where $\Delta ({\rm sdet}_q(T)) =
{\rm sdet}_q(T) \, \underline{\otimes} \, {\rm sdet}_q(T)$.
 The element
${\rm sdet}_q(T)$ is called the quantum Berezinian (or quantum super-determinant).}
\end{proposition}

We now compare the relations (\ref{3.6.4}) with the graded Yang-Baxter equation
(\ref{3.6.2}). From this comparison we readily see that the finite-dimensional
matrix representations for the generators $T^{i}_{j}$ of the quantum algebra
$Fun(GL_{q}(N|M))$ (the super-analogs of the representations (\ref{3.1.18}))
can be chosen in the form
\be
\lb{3.6.7}
(T_{1})_3 = (-)^{(1)(3)} R_{13} (-)^{(1)(3)}
 \equiv R^{(+)}_{13} \; , \;\;
(T_{1})_3 = (R^{-1})_{31} \equiv R^{(-)}_{13} \; .
\ee
From this, in an obvious manner, we obtain definitions of the quantum
superalgebras which are dual to the algebras $Fun(GL_{q}(N|M))$
 (cf. Eqs. (\ref{3.1.19})):
\be
\lb{3.6.8}
\langle L^{+}_{2},T_{1} \rangle = (-)^{(1)(2)} R_{12} (-)^{(1)(2)} = R_{12} \; , \;\;
\langle L^{-}_{2},T_{1} \rangle =
R^{-1}_{21} ,
\ee
were operator-valued matrices $L^{\pm}$ satisfy
\be
\lb{frtsup}
\begin{array}{c}
\hat{R}_{12} \,  L^{\pm}_2 \, (-)^{(1)(2)} \, L^{\pm}_1  \, (-)^{(1)(2)} =
L^{\pm}_2 \, (-)^{(1)(2)} \, L^{\pm}_1   \, (-)^{(1)(2)} \, \hat{R}_{12}  \; , \\[0.2cm]
\hat{R}_{12} \,  L^{+}_2 \,(-)^{(1)(2)} \, L^{-}_1  \, (-)^{(1)(2)} =
L^{-}_2 \, (-)^{(1)(2)} \, L^{+}_1  \, (-)^{(1)(2)} \,  \hat{R}_{12} \; .
\end{array}
\ee
By using the identity $\R_{12} (-)^{(1)(2)} = (-)^{(1)(2)} \R_{12}$
(see (\ref{3.6.5ra})) for the $R$-matrix
(\ref{3.6.1aa}), one can deduce from (\ref{frtsup}) the standard reflection
equation (\ref{3.1.23a}) for the matrix $L = S(L^{-}) L^+$.

Recall that the $R$-matrix (\ref{3.6.1aa}) for $GL_{q}(N|M)$
 is such that its diagonal blocks
$\hat{R}^{i_1 i_2}_{k_1 k_2}[q]$ and $\hat{R}^{\alpha_1 \alpha_2}_{\gamma_1
\gamma_2}[q]$ are standard $Fun(GL_{q}(N))$ and $Fun(GL_{q}(M))$
$R$-matrices of the Hecke type and we have commutation relations
 (\ref{rAArDD}), (\ref{3.6.6X}) for matrices $D$ and $X = A-BD^{-1}C$.
Then, one can write the quantum superdeterminant
for $GL_{q}(N|M)$ by means of the definition (\ref{sdetb}) in the
form \cite{21}, \cite{33'}, \cite{33}
\be
\lb{sdeta}
{\rm sdet}_{q}^{-1}(T) =
{\det}_{q^{-1}} (A - B D^{-1} C)^{-1} \, {\det}_{q^{-1}}(D) \; .
\ee
where ${\det}_q$ is defined in (\ref{3.3.20a}). Then, the algebra $Fun(SL_{q}(N|M))$ is
distinguished by the relation $sdet_{q}(T)=1$.

\vspace{0.1cm}

\noindent
{\bf Remark.} The standard formula for
superdeterminant of supermatrix $T$ is deduced
 from the integral representation
\be
\lb{ber}
{\rm sdet}(T^{M}_{N}) \sim \int \, \prod_{M} (dE^{M} \, dF_{M}) \,
\exp [ i E^{N} \, (T^{M}_{N}) \, F_{M} ] \; ,
\;\;\;\;\; T^{M}_{N} = \left(
\begin{array}{cc}
A^{\xi}_{\eta} & B_{\eta}^{m} \\
C_{n}^{\xi} & D^{m}_{n}
\end{array} \right) \; ,
\ee
where supermatrix $T$ is given in the block
form  and we respectively divide the supervectors
$E^{N} = (b^{\eta}, \, \beta^{n})$ and $F_{M} = (c_{\xi}, \, \gamma_{m})$
on even  $\{\beta$, $\gamma\}$ and odd $\{b$, $c\}$ parts.
Then, we transform the quadratic expression
$E^{N} \, (T^{M}_{N}) \, F_{M}$
to the ''diagonal'' form
$$
E \, (T) \, F = b \, (A - B \, D^{-1} \, C) \, c +
\widetilde{\beta} \, D \, \widetilde{\gamma} \; ,
$$
by making the linear change of even variables
$\beta = \widetilde{\beta} - b \, B \, D^{-1}$,
$\gamma = \widetilde{\gamma} -  D^{-1} \, C \, c$ .
The Jacobian for such change of variables is equal to $1$.
After integration over $b,c$ and
$\widetilde{\beta},\widetilde{\gamma}$ in (\ref{ber}), we obtain
 \be
 \lb{sdet0}
{\rm sdet}(T^{M}_{N}) =
\det (A - B \, D^{-1} \, C) \, \frac{1}{\det D}  \; .
 \ee
The element ${\rm sdet}_{q}^{-1}(T)$,
which appeared in
eqs. (\ref{sdetb}) and (\ref{sdeta}), is denoted as an
inverse of the super-determinant ${\rm sdet}_{q}(T)$
since the element ${\rm sdet}^{-1}_{q}(T)$ tends
to ${\rm sdet}^{-1}(T)$ for $q \to 1$ in
view of the standard formula (\ref{sdet0}).
We also note that the meaning of
Proposition {\bf \em \ref{prop6}} is to find the
Jacobian of the super-coordinate transformation for the measure (\ref{sdetf}) of
an integration over multi-dimensional quantum superplane.

 \vspace{0.1cm}

The quantum supergroup $GL_{q}(N|M)$ was studied in detail from somewhat
different positions in Ref. \cite{33}. The simplest example of a quantum supergroup,
$GL_{q}(1|1)$, has been investigated in many studies (see, for example, Refs.
\cite{31} and \cite{34}).
The $R$-matrices (\ref{3.4.3}) can be used
(see next Sect. {\bf \ref{baxtel}}) to construct the supersymmetric
Baxterized solutions of the Yang-Baxter equation (\ref{3.5.5})
obtained in Ref. \cite{35}. The Yangian limits of these
solutions\footnote{The corresponding RTT algebra
 defines the Yangian
of the Lie superalgebra $gl(n|m)$ \cite{33'}.} were used to formulate integrable
supersymmetric spin chains (see e.g. \cite{36}).
The universal R-matrices for the linear quantum
supergroups (and more generally for
quantum deformations of finite dimensional
contragredient Lie (super)algebras) were constructed in Ref. \cite{37}.

\subsection{\bf \em $GL_q(N)$- and $GL_q(N|M)$- invariant Baxterized $R$-matrices. Dynamical $R$-matrices\label{baxtel}}
\setcounter{equation}0

By Baxterization, we mean the construction of an $R$-matrix that
depends not only on a deformation parameter $q$ but also on an additional
complex spectral parameter $x$. We wish to find a solution $\R(x)$ of the
Yang-Baxter equation with spectral parameter $x$
 (see eq. (\ref{3.5.2}) below) satisfying the
 quantum invariance condition
$$
T_{1} T_{2} \; \R(x) \; (T_{1} T_{2})^{-1} = \R(x) \; , \;\;\;\;\;\;
\left( T^{i}_{j} \in Fun(GL_{q}(N)) \right) \; .
$$
Then we must seek it in the form \cite{Isaev1}
\be
\lb{3.5.1}
\R(x) = b(x)(\hbox{\bf 1} + a(x) \R) \; ,
\ee
(here $a(x)$ and $b(x)$ are certain functions of $x$),
since by virtue of the Hecke relation (\ref{3.3.7aa}) there exist only
two basis matrices ${\bf 1}$ and $\R$
 that are invariants in the sense of the relations
$T_{1} T_{2} \; \R \; (T_{1} T_{2})^{-1} = \R$
followed from (\ref{3.1.1}). The Yang-Baxter equation with dependence on the spectral parameter
is chosen in the form
\be
\lb{3.5.2}
\R_{12}(x) \, \R_{23}(xy) \, \R_{12}(y)  =  \R_{23}(y) \, \R_{12}(xy) \, \R_{23}(x) \; .
\ee
Only the function a(x) is fixed by this equation. Indeed, we substitute here
(\ref{3.5.1}) and take into account (\ref{3.1.3}) and the Hecke condition
(\ref{3.3.7}). As a result we obtain the equation \cite{Isaev1}
\be
\lb{3.5.3}
a(x) + a(y) + \lambda a(x)a(y) = a(xy),
\ee
which is readily solved by the change of variables
$a(x) = (1/\lambda)(\tilde{a}(x) -1)$. After this, we obtain for
the function $a$ the general solution
\be
\lb{3.5.4}
a(x) = (1/\lambda) (x^{\xi}-1),
\ee
where for simplicity the arbitrary parameter $\xi$ can be set equal to $-2$.
For convenience, we choose the normalizing function $b(x)= x$.
Then the Baxterized $R$-matrix satisfying the Yang-Baxter equation (\ref{3.5.2})
will have the form \cite{Jimb1}, \cite{Jon3}, \cite{CGX},
 \cite{Isaev1}
\be
\lb{3.5.5}
\R(x)= b(x) \left( \hbox{\bf 1} +
 (1/\lambda) (x^{-2}-1) \R \right) =
  {1 \over \lambda} \, (x^{-1}\R - x \R^{-1}) \; .
\ee
Remarkably, this matrix is written as the rational function
 of $\R$
\be
\lb{3.5.5b}
\R(x) =  \frac{(a^{-1} x - a x^{-1})}{\lambda\, x^2} \;
\frac{\R -a \, x^2}{\R -a \, x^{-2}} \; ,  \;\;\;\;\;\;\;
a = \mp q^{\pm 1} \; .
\ee
Below we call this $R$-matrix as the {\it
Hecke type Baxterized $R$-matrix}.
For the normalization adopted in (\ref{3.5.5})
we obtain
\be
\lb{normgl}
\R(1) = 1 \; , \;\;\; \P^{\pm} = \frac{1}{[2]_q} \, \R(q^{\mp 1}) \; ,
\ee
and the unitarity condition holds\footnote{Strictly speaking we have to renormalize
the $R$-matrix (\ref{3.5.5}):
$\hat{R}(x) \rightarrow \lambda (x^{-1} q - x q^{-1})^{-1} \hat{R}(x)$,
to obtain the unitarity condition
with the unit matrix in
the right hand side of (\ref{conuni}).}
\be
\lb{conuni}
\R(x) \, \R(x^{-1}) = \left( 1 - \frac{(x - x^{-1})^2}{\lambda^2} \right) \; .
\ee
This unitarity follows from rational representation (\ref{3.5.5b})
and can be readily deduced from the
spectral decomposition
$$
\R(x) = \frac{(x^{-1} q - x q^{-1})}{\lambda} \P^{+} +
 \frac{(x q - (x q)^{-1})}{\lambda} \P^{-} \; ,
$$
where projectors $\P^{\pm}$
were defined in (\ref{3.3.10}), (\ref{3.3.11}).
Note that we have obtained the Baxterized solution (\ref{3.5.5}) of the Yang-Baxter equation
(\ref{3.5.2}) only using the braiding relations (\ref{3.1.2i})
and Hecke condition (\ref{3.3.7aa}) for the constant matrix $\R$. Thus, any
constant Hecke solution of (\ref{3.1.2i}) (e.g. the
multi-parametric solution (\ref{rmult})) can be used for the construction
of the Baxterized $R$-matrices (\ref{3.5.5}).

For the Baxterized $R$-matrix (\ref{3.5.5}), constructed via
skew-invertible Hecke type $R$-matrix, one can deduce
 the cross-unitarity conditions
\be
\lb{cross01}
\begin{array}{c}
Tr_{D(2)}\bigl( \R_{1}(x) P_{01} \R_1(z) \bigr) =
\eta(x,z) \, D_0 \, I_1 \; , \;\;\;\;\;\; 
Tr_{Q(1)}\bigl( \R_{1}(x) P_{23} \R_1(z) \bigr) =
\eta(x,z) \, Q_{3} \, I_2 \; ,
\end{array}
\ee
where
$$
\eta(x,z)=  \frac{(x-x^{-1})(z-z^{-1})}{\lambda^2} \; ,
$$
$P_{01},P_{23}$ are permutations, matrices $D,Q$ were
  defined in (\ref{dmatr})
and spectral parameters $x,z$ are constrained by the condition
$$
(x\, z)^2 = \frac{1}{1-\lambda Tr(D)} =:b^2 \; .
$$
We stress that for the $GL(N|M)$-type $R$-matrix
we have $b^2 = q^{2(N-M)}$.

Let $\hat{\Psi}_{12}$
be a skew-inverse matrix (\ref{skew}) for the Hecke $R$-matrix (\ref{3.3.7}). Then, for the
Baxterized $R$-matrix (\ref{3.5.5}), one can define the skew-inverse
Baxterized matrix $\hat{\Psi}(x)$:
\be
\lb{skinh}
\hat{\Psi}_{12}(x) = \frac{\lambda}{x^{-1}-x} \, \Bigl( \hat{\Psi}_{12}  +
\frac{\lambda}{b^{-2} - x^{-2}} \, D_1 Q_2 \Bigr) \; ,
\ee
such that
\be
\lb{skinh1}
Tr_2 \bigl( \hat{\Psi}_{12}(x) \, \R_{23}(x) \bigr) \; = \; P_{13}
\; = \;
Tr_2 \bigl( \R_{12}(x)  \, \hat{\Psi}_{23}(x)  \bigr) \; .
\ee

\vspace{0.2cm}

Let $x_i$ and $p_{j}$ $(i,j=1,...,N)$
be generators of the Heisenberg algebra
\be
\lb{10a}
[x_{i}, \, p_{j} ] = i \, \hbar \, \delta_{ij}
 \;\;\;\;\;\; (i,j \leq N-1) \; ,
\ee
where $\hbar$ is a Planck constant. The dynamical
Yang-Baxter equation is defined as follows
\cite{GerNev}, \cite{Feld}, \cite{21'}
(see also \cite{BytFad}, \cite{EtVar}, \cite{HIOPT})
\be
\lb{14a}
(Q^{-1}_{3} \; \hat{R}_{12}(p) \; Q_{3})
\; \hat{R}_{23}(p) \;
(Q^{-1}_{3} \; \hat{R}_{12}(p) \; Q_{3})
 = \hat{R}_{23}(p) \;
(Q^{-1}_{3} \; \hat{R}_{12}(p) \; Q_{3})
 \; \hat{R}_{23}(p) \; ,
\ee
where $Q:= {\rm diag}(e^{ix_1}, e^{ix_2},...,e^{ix_N})$.
We seek the solution of (\ref{14a}) in the form
(cf. (\ref{anzdym}))
\be
\lb{15ii}
\hat{R}_{12} =
\hat{R}^{i_{1}i_{2}}_{j_{1}j_{2}}(p) =
\delta^{i_{1}}_{j_{2}}\delta^{i_{2}}_{j_{1}} \, a_{i_{1}i_{2}}(p) +
\delta^{i_{1}}_{j_{1}}\delta^{i_{2}}_{j_{2}} \, b_{i_{1}i_{2}} (p) \; ,
\ee
and require that this $R$-matrix satisfies Hecke condition (\ref{3.3.7}).
Without limitation of generality one can put $b_{ii}(p) = 0$.
Now we substitute
(\ref{15ii}) to the dynamical Yang-Baxter
 equation (\ref{14a}) and obtain
 the following constraints
\be
\lb{19ii}
a_{ij}(p_{1}, \dots , p_{N}) = a_{ij}(p_{i},p_{j}) \; , \;\;
b_{ij}(p_{1}, \dots , p_{N}) = b_{ij}(p_{i},p_{j}) \; , \;\;
\ee
 and equations \cite{21'}
 \be
\lb{19heck}
a_{i}^{2} - \lambda a_{i} - 1 = 0 \; , \;\;\;\;
b_{ij}(p_{i},p_{j}) + b_{ji}(p_{j},p_{i}) = \lambda \; ,
\;\;\;\;\;\; i \neq j  \; ,
\ee
 \be
\lb{20heck}
a_{ij}(p_{i},p_{j}) \, a_{ji}(p_{j},p_{i}) -
b_{ij}(p_{i},p_{j}) \, b_{ji}(p_{j},p_{i}) = 1 \; ,
\;\;\;\;\;\; i \neq j \; ,
\ee
\be
\lb{20ii}
b_{ij} \, b_{jk} \, b_{ki} +
b_{ik} \, b_{kj} \, b_{ji} =0 \; ,
\;\;\;\;\;\; i \neq j \neq k \neq i \; ,
\ee
\be
\lb{21ii}
b_{ij}(p_{i} + \hbar, p_{j}) = \frac{ b_{ij}(p_{i},p_{j}) \, a_{i}}{1/a_{i} +
b_{ij}(p_{i},p_{j})} \; , \;\;\;\;\;
b_{ij}(p_{i} , p_{j}+ \hbar) = \frac{  b_{ij}(p_{i},p_{j}) / a_{j}}{
a_{j} - b_{ij}(p_{i},p_{j})} \; ,
\ee
where $a_i := a_{ii}$, equations (\ref{19heck}),
(\ref{20heck}) are consequences
 of the Hecke condition (\ref{3.3.7}), while equations
 (\ref{20ii}), (\ref{21ii}) follow from (\ref{14a}).
 The general solution of these equations
 for coefficients $b_{ij}(p)$ are \cite{21'}
\be
\lb{24ii}
b_{ij}(p_{i},p_{j}) = \frac{ \lambda \, a_{i}^{p_{i}/\hbar} \, a_{j}^{-p_{j}/\hbar}
\, b_{ij}^{0} }{
a_{i}^{p_{i}/\hbar}\, a_{j}^{-p_{j}/\hbar}\, b_{ij}^{0} +
a_{i}^{-p_{i}/\hbar} \, a_{j}^{p_{j}/\hbar} \, b_{ji}^{0} } \; ,
\ee
where constants $b_{ij}^{0} := b_{ij}(0,0)$ have to
obey the algebraic relations:
\be
\lb{25ii}
b^{0}_{ii}=0 \; , \;\;\;\; b^{0}_{ij} + b^{0}_{ji} = \lambda \; ,  \;\;\;\;
b^{0}_{ij} \, b^{0}_{jk} \, b^{0}_{ki} +
b^{0}_{ik} \, b^{0}_{kj} \, b^{0}_{ji} =0 \; .
\ee
The first equation in (\ref{19heck}) has two solutions
$a_i = \pm q^{\pm 1}$. Recall (see Subsection {\bf \ref{qsuper}}),
that if we take $a_{i} = q \; , \; \forall i$
(or $a_{i} = -q^{-1} \; , \; \forall i$) then we will have the case
of the standard quantum group $GL_{q}(N)$ (or $GL_{-1/q}(N)$).
 But if we consider
the mixing case, $a_{i} = q$ for $1 \leq i \leq K$ and
$a_{i} = - q^{-1}$ for $ K + 1 \leq i \leq N$, then we come to the
case of supergroups $GL_{q}(K|N-K)$. By considering
 the solution (\ref{24ii}), it is clear that if
$a_{i} = a_{j}$ (indices $i$ and $j$ `have the same grading')
then $b_{ij}(p_{i},p_{j}) = b_{ij}(p_{i}-p_{j})$, but if
$a_{i} = -1/a_{j}$ (the case of supergroups when indices
$i$ and $j$ `have the opposite grading') then we deduce that
$b_{ij}(p_{i},p_{j}) = b_{ij}(p_{i}+p_{j})$.
Note that the only conditions on the parameters $a_{ij}(p)$
needed for fulfillment of the dynamical Yang-Baxter equation
are listed in (\ref{20heck}).

Now we demonstrate that every
solution $R(p)$ given in (\ref{15ii}), (\ref{20heck}),
 (\ref{24ii}) will lead to
the solution $R(p,z)$ for the dynamical
Yang-Baxter equation with spectral parameters
\be
\lb{28ii}
\hat{R}_{12}(p,\, y) \, Q_3 \, \hat{R}_{23}(p,\, y  z)
 \, Q_3^{-1} \, \hat{R}(p,\, z) =
 Q_3 \, \hat{R}_{23}(p,\, z) \, Q_3^{-1} \,
 \hat{R}_{12}(p,\, y  z) \,
 Q_3 \, \hat{R}_{23}(p,\, y) \, Q_3^{-1} \, .
\ee
Indeed, it is not difficult to check, by using (\ref{14a})
and the Hecke relation for $\hat{R}(p)$, that the following
matrices (cf. (\ref{3.5.5}))
$$
\hat{R}(p,\, y) = y^{-1} \, \hat{R}(p) - y \, \hat{R}(p)^{-1}
\; , \;\;
$$
are the solutions of (\ref{28ii}).
We note that these solutions satisfy the identity
(cf. (\ref{conuni}))
$$
\hat{R}(p, \, y) \, \hat{R}(p, \, y^{-1}) =
\left( \lambda^{2} - (y-y^{-1})^{2} \right) \; ,
$$
which is a kind of unitary condition
for $\hat{R}(p)$ (if $y^{*} = y^{-1}$).

\subsection{\bf \em Quantum matrix algebras with spectral parameters.
 Yangians $Y_q(gl_N)$ and $Y(gl_N)$}
 \setcounter{equation}0

It is a remarkable fact that the relations (\ref{3.1.20b}),
(\ref{3.1.20}), with the Hecke $\R$-matrix,
are written as follows:
\be
\lb{3.5.6}
\R_{12}(x) \, L_2(xy) \, L_1(y) = L_2(y) \, L_1(xy) \, \R_{12}(x) ,
\ee
\be
\lb{3.5.7}
L(x) := x^{-1} L^{+} - x L^{-} .
\ee
where $x$ and $y$ are arbitrary spectral parameters
and $\R(x)$ is Baxterized $R$-matrix (\ref{3.5.5}).
Moreover, if we take the pairing of the
relation (\ref{3.5.6}) with the representation matrix
$T^{i}_{j}$ acting in the third space and use (\ref{3.1.19}), we obtain the Yang-Baxter equation
(\ref{3.5.2}) for the solution (\ref{3.5.5}).
Thus, in a certain sense, eq. (\ref{3.5.6}) generalizes (\ref{3.5.2}).

Now we take $GL_q(N)$ type Baxterized $R$-matrix (\ref{3.5.5})
and consider eqs. (\ref{3.5.6}) as defining relations for
new infinite dimensional algebras with generators
$(L_{(r)})^i_j$, which appeared
in the expansion
 \be
 \lb{3.5.7b}
L^i_j(x) = \sum_{r \geq 0}
(L_{(r)})^i_j \; x^{- 2r} \; .
 \ee
This algebra is called {\em quantum}
Yangian $Y_q(gl_N)$ and it is a subalgebra in
a quantum affine algebra $U_q(\hat{gl}_N)$
(the $RTT$ definition of $U_q(\hat{gl}_N)$ is given in \cite{RSTS},
\cite{DiFre}).
Formula (\ref{3.5.7}) defines a homomorphism
$Y_q(gl_N) \rightarrow U_q(gl_N)$ which
is called {\it evaluation homomorphism}.

Let $R(x)$ be some
solution of the Yang-Baxter equation (\ref{3.5.2})
 and there is an algebra with defining relations
 (\ref{3.5.6}).
It is known that (\ref{3.5.2}) are associativity conditions of
the unique ordering of monomials of third degree
$L_{1}(x)L_{2}(y)L_{3}(z)$ for the algebra (\ref{3.5.6}).
Indeed, we have the following diagram
 (the so-called "diamond" condition):

\unitlength=6mm
\begin{picture}(17,3.7)(-2,0)
\put(-0.5,1.6){\footnotesize $L(x)L(y)L(z)$}
\put(3.2,2.2){\vector(3,2){1.2}}
\put(3.2,1.3){\vector(3,-2){1.2}}
\put(4.5,2.8){\footnotesize  $L(y)L(x)L(z)$}
\put(4.5,0.3){\footnotesize  $L(x)L(z)L(y)$}
\put(8.2,3){\vector(1,0){1.4}}
\put(8.2,0.5){\vector(1,0){1.4}}
\put(9.8,2.8){\footnotesize $L(y)L(z)L(x)$}
\put(9.8,0.3){\footnotesize  $L(z)L(x)L(y)$}
\put(13.5,2.8){\vector(3,-2){1.2}}
\put(13.5,0.5){\vector(3,2){1.2}}
\put(14.7,1.5){\footnotesize $L(z)L(y)L(x)$}
\end{picture}

\noindent
This diagram means that two different ways of reordering
$L(x)L(y)L(z) \to L(z)L(y)L(x)$ (by means of (\ref{3.5.2}))
gives the same result.

Now we stress that, for the quantum algebra (\ref{3.5.6})
with special $R$-matrix (\ref{3.5.5}), defined
by the Hecke $R$-matrix of the height $N$,
the quantum determinant (the analog of (\ref{3.3.20a}))
can also be constructed \cite{Jimb1} :
\be
\lb{hahaha0}
{\det}_{q}(L(x)) \, {\cal E}_{\langle 12 \dots N}
= {\cal E}_{\langle 12 \dots N} \,  L_N(q^{N-1} x) \cdots  L_2(q x) \, L_1(x) \;\; \Leftrightarrow
\ee
$$
{\det}_q(L(x)) = Tr_{1 \dots N} \left( A_{1 \to N} \,
 L_N(q^{N-1} x) \cdots  L_2(q x) \, L_1(x) \right) =
$$
\be
\lb{hahaha}
= Tr_{1 \dots N} \left(  L_N( x) \, L_{N-1}(q x) \cdots L_1( q^{N-1} x) \, A_{1 \to N}
 \right) \; ,
\ee
\be
\lb{hahaha1}
{\det}_{q}(L(x)) \, A_{1 \to N}
=   L_N( x) \, L_{N-1}(q x) \cdots L_1( q^{N-1} x) \, A_{1 \to N} \; ,
\ee
where the rank 1 antisymmetrizer $A_{1 \to N}$ has been introduced in (\ref{antirr}).
Eq. (\ref{hahaha0}) is self-consistent since its right hand side has the same symmetry as the
left hand side (the action on both sides of this equation by the projectors (\ref{normgl})
$\P^{+}_k \sim \R_k(q^{-1})$ gives zero).
The last form (\ref{hahaha}) of the quantum determinant
${\det}_{q}(L(x))$ is obtained with the
help of (\ref{antirr}) and (\ref{3.5.6}).


\begin{proposition}\label{prop7}
{\it The $q$-determinant
${\det}_q(L(x))$ is a generating function of central elements
for the algebra (\ref{3.5.6}) with $GL_q(N)$-type
 Baxterized $R$-matrix (\ref{3.5.5}).}
\end{proposition}

\vspace{0.2cm}
\noindent
{\bf Proof.} The centrality of ${\det}_q(L(x))$
means that $[L^i_j(x y), \, {\det}_q(L(x))]=0$
$\forall x,y$. Indeed,
\be
\lb{hahaha5}
L_{N+1}(x y) \,  Tr_{1 \dots N} \left(  L_N( x) \, L_{N-1}(q x) \cdots
L_1( q^{N-1} x) \, A_{1 \to N} \right) =
\ee
\be
\lb{hahaha3}
\begin{array}{c}
 Tr_{1 \dots N} \left( \hat{R}^{-1}_N(y)  \dots
\hat{R}^{-1}_1(q^{1-N}y)  \,  L_{N+1}( x) \cdots
L_2( q^{N-1} x) L_{1}(x y) \cdot \right. \\ [0.2cm]
\left. \cdot   \hat{R}_1(q^{1-N}y)  \dots
\hat{R}_N (y) \, A_{1 \to N} \right) \; .
\end{array}
\ee
Using the Yang-Baxter eq. (\ref{3.5.2}) and the representation of $A_{1 \to N}$
in terms of the Baxterized elements (\ref{antirr}) we deduce
$$
 \hat{R}_1(q^{1-N}y)  \dots  \hat{R}_N (y) \, A_{1 \to N} =
 A_{2 \to N+1} \,  \hat{R}_1(y)  \dots  \hat{R}_N (q^{1-N}y) \,
$$
By means of this relation and eq. (\ref{hahaha1})
one can rewrite (\ref{hahaha3}) in the form
$$
\begin{array}{c}
{\det}_{q}(L(x)) \,  Tr_{1 \dots N} \left( \hat{R}^{-1}_N(y)  \dots
\hat{R}^{-1}_1(q^{1-N}y)  \,  L_{1}(x y)
 A_{2 \to N+1} \cdot \right. \\ [0.2cm]
\left.
\cdot  \hat{R}_1(y)  \dots  \hat{R}_N (q^{1-N}y)
 A_{1 \to N}  \right) =
 \end{array}
$$
\be
\lb{hahaha4}
= {\det}_{q}(L(x)) \, \left( N^{-1}(y)
\,  L(x y) \, N(y) \right)_{N+1}  =
{\det}_{q}(L(x)) \,  L_{N+1}(x y)
\ee
where matrices $N(y)$ and $N(y)^{-1}$ are defined by
$$
(N(y))_{\langle N+1}^{1 \rangle}= {\cal E}_{\langle 2 \dots N+1} \,
\hat{R}_1(y)  \dots  \hat{R}_N (q^{1-N}y) \, {\cal E}^{1 \dots N \rangle} \; ,
$$
$$
(N(y)^{-1})_{\langle 1}^{N+1 \rangle}= {\cal E}_{\langle 1 \dots N} \,
\hat{R}^{-1}_N(y)  \dots
\hat{R}^{-1}_1(q^{1-N}y) \, {\cal E}^{2 \dots N+1 \rangle} \; .
$$
and, for $GL(N)$-type Baxterized $R$-matrices, they are proportional to the unit matrix. Comparing
(\ref{hahaha5}) and
(\ref{hahaha4}) we obtain the statement of the Proposition.
\hfill \qed

\noindent
We stress here that not for all
Hecke-type Baxterized $R$-matrices the element
${\det}_{q}(L(x))$ is central for the algebra (\ref{3.5.6}).
The example is given by
multi-parametric Hecke $R$-matrices (\ref{3.4.2}).

\vspace{0.2cm}

We now note that from the algebra (\ref{3.5.6}), disregarding the particular representation
(\ref{3.5.7}) for the $L(x)$ operator, we can obtain a realization for the
Yangian $Y(gl(N))$ \cite{Drin85}, \cite{13}
(see also review paper \cite{Mol}). Indeed, in (\ref{3.5.2}) and (\ref{3.5.6})
we make the change of spectral parameters
\be
\lb{3.5.8}
x = exp(- \frac{1}{2}\lambda (\theta - \theta')) \; , \;\;
y = exp(- \frac{1}{2}\lambda \theta') \; .
\ee
Then the relations (\ref{3.5.2}) and (\ref{3.5.6}) can be rewritten in the form
\be
\lb{3.5.9}
\R_{12}(\theta - \theta') \, \R_{23}(\theta) \, \R_{12}(\theta')
=  \R_{23}(\theta') \, \R_{12}(\theta) \, \R_{23}(\theta -\theta') \Rightarrow
\ee
\be
\lb{3.5.9a}
R_{23}(\theta - \theta') \, R_{13}(\theta) \, R_{12}(\theta')
=  R_{12}(\theta') \, R_{13}(\theta) \, R_{23}(\theta -\theta') \; ,
\ee
\be
\lb{3.5.10}
\R_{12}(\theta -\theta') \, L_2(\theta) \, L_1(\theta') =
L_2(\theta') \, L_1(\theta) \, \R_{12}(\theta - \theta') ,
\ee
where we redefine $L$-operator
$L(\theta) :=
L \left(\exp(-{\lambda \over 2}\theta) \right)$
 and $R$-matrix
\be
\lb{3.5.10a}
\R(\theta) := \R \left( e^{-{\lambda \over 2}\theta} \right) =
\cosh \left(\lambda \theta /2 \right) +
\frac{1}{\lambda} \, \sinh \left(\lambda \theta /2 \right)
(\R + \R^{-1}) \; .
\ee
Eqs. (\ref{3.5.9a}) have a beautiful graphical representation in
the form of the triangle equation \cite{3} \\
\unitlength=8mm
\begin{picture}(15.5,5)(0,-1)
\put(3.7,0){\footnotesize $3$}
\put(4.2,0.6){\tiny $\theta$}
\put(6.1,1.9){\tiny $\theta'$}
\put(4.1,2.3){\tiny $\theta\! -\! \theta'$}
\put(4,4){\vector(0,-1){4}}
\put(3.3,3.6){\vector(3,-2){3.2}}
\put(6.2,1.2){\footnotesize $2$}
\put(6.2,2.6){\footnotesize $1$}
\put(3.3,0.4){\vector(3,2){3.2}}
\put(7.5,1.9){$=$}

\put(10.5,4){\vector(0,-1){4}}
\put(8.3,2.5){\vector(3,-2){3.2}}
\put(8.3,1.5){\vector(3,2){3}}
\put(9.5,1.9){\tiny $\theta'$}
\put(10.7,2.7){\tiny $\theta$}
\put(10.6,0.2){\tiny $\theta\! -\! \theta'$}
\put(11.3,0.6){\footnotesize $2$}
\put(11.1,3.6){\footnotesize $1$}
\put(10.1,0){\footnotesize $3$}
\end{picture}

\vspace{-2cm}

\be
\lb{zzz}
{}
\ee
where the arrowed lines show trajectories of point particles
and the $R$-matrix \\

\unitlength=8mm
\begin{picture}(15,2)
\put(5,1){$R_{ij}(\theta) =$}
\put(6.8,2){$i$}
\put(8.8,2){$j$}
\put(7.7,1.3){\tiny $\theta$}
\put(7,0.2){\vector(1,1){1.5}}
\put(8.5,0.2){\vector(-1,1){1.5}}
\end{picture}
\\
describes a single act of the scattering of these particles.
We now take the limit $\lambda = q- q^{-1} \rightarrow 0$
in eq. (\ref{3.5.10}).
On the basis of (\ref{3.5.5}), (\ref{3.5.10a}) we readily find that in this limit the
matrix $\R(\theta)$
is equal to the Yang matrix:
\be
\lb{3.5.11}
\R(\theta) =
({\hbox{\bf 1}} + \theta \, P_{12}) \;\;\;\;
\Rightarrow \;\;\;\;  R_{12}(\theta) =
\theta \, \Bigl({\hbox{\bf 1}} +  \frac{P_{12}}{\theta}\Bigr)
 \; .
\ee
For the operators $L(\theta)$, we shall assume the expansion
\be
\lb{3.5.12}
L(\theta)^{i}_{j} =
\sum_{k=0}^{\infty} {T^{(k)}}^{i}_{j} \theta^{-k} \; ,
\ee
where ${T^{(0)}}^{i}_{j} = \delta^{i}_{j}$
and  ${T^{(k)}}^{i}_{j}$ $(k >0)$ are the generators of the
Yangian $Y(gl(N))$ (see Ref. \cite{13}). The defining relations for the Yangian
$Y(gl(N))$ are obtained from (\ref{3.5.10}) by substituting
(\ref{3.5.11}) and (\ref{3.5.12}) (we give
these relations in more general form of the super
Yangian $Y(gl(N|M))$; see below (\ref{yaglnm})). The comultiplication for $Y(gl(N))$ obviously has the form
\be
\lb{3.5.13}
\Delta(L(\theta)^{i}_{j} ) = L(\theta)^{i}_{k} \otimes
L(\theta)^{k}_{j} \; .
\ee
The Yangian $Y(sl(N))$ is obtained from
$Y(gl(N))$ after the imposition of an additional
  condition for the generators ${T^{(k)}}^{i}_{j}$:
$$
{\det}_{q}(L(\theta)) = 1 \; ,
$$
where the Yangian quantum determinant \cite{KuSk82}
\be
\lb{ydetq}
{\det}_{q}(L(\theta)) =  Tr_{1 \dots N} \left( A_{1 \to N}^{cl} \,
 L_N(\theta - N +1) \cdots  L_2(\theta -1) \, L_1(\theta) \right)
\ee
is obtained from (\ref{hahaha}) after substitution
$q = e^h$, $x = \exp(-{\lambda \over 2} \theta) \sim e^{- h \theta}$
and taking the limit $h \to 0$ (or $\lambda \to 0$). In (\ref{ydetq}) we denote by
$A_{1 \to N}^{cl}$ a classical antisymmetrizer:
 $$
 A_{1 \to N}^{cl}=
 \lim\limits_{q \to 1} A_{1 \to N} = \frac{1}{N!}
(1 + P_{N-1} +  ... +
P_{1} \cdots P_{N-1}) \cdots
(1 + P_2 + P_{1} P_{2})(1 + P_1) \; .
 $$

Since the $\R$ matrix (\ref{3.6.1aa}), (\ref{3.6.1}) (for the group $GL_q(N|M)$) satisfies the
Hecke condition (\ref{3.3.7aa}), the same Baxterized $R$-matrix (\ref{3.5.5}) is appropriate for
the supersymmetric case. Almost all statements of this subsection can be readily reformulated for
the supersymmetric case. In particular, the Yangian $R$-matrix for $Y(gl(N|M))$ is deduced from
(\ref{3.5.10a}) and has the form
(cf. (\ref{3.5.11})):
\be
\lb{sgly}
\R(\theta) =
({\hbox{\bf 1}} + \theta \, {\cal P}_{12}) \; ,
\ee
where ${\cal P}_{12}$ is a supertransposition operator introduced in (\ref{stran}).
The defining relations (\ref{3.5.6}) should be modified
for the super Yangian $Y(gl(N|M))$ (cf. (\ref{3.6.4})):
\be
\lb{sgly1}
\R_{12}(\theta - \theta') \, (-)^{(1)(2)} \,  L_2(\theta) \, (-)^{(1)(2)} \, L_1(\theta') =
(-)^{(1)(2)} \, L_2(\theta') \, (-)^{(1)(2)} \, L_1(\theta) \, \R_{12}(\theta - \theta') ,
\ee
while the form of the comultiplication (\ref{3.5.13})
(where $\otimes$ is the graded tensor product) is unchanged. Taking into account
(\ref{3.5.12}) and (\ref{sgly}) we obtain the component form
of the defining relations (\ref{sgly1}) for $Y(gl(N|M))$
\be
\lb{yaglnm}
\, [ {T^{(r)}}^i_j , \,  {T^{(s+1)}}^k_l \} -
[ {T^{(r+1)}}^i_j , \,  {T^{(s)}}^k_l \} = (-1)^{[k][i] + [k][j] + [i][j]}
\left( {T^{(s)}}^k_j  \,  {T^{(r)}}^i_l - {T^{(r)}}^k_j \,  {T^{(s)}}^i_l \right) \; ,
 \ee
where $r,s \geq 0$, ${T^{(0)}}^i_j = (-1)^{[i]} \, \delta^i_j$, the grading $[i]= 0,1$ mod$(2)$ is
defined in (\ref{pari}) and $[a,b \}$ denotes a supercommutator
$[a,b \} := a \, b - (-1)^{[a] [b]}b \, a$, $[a] = {\rm deg}(a)$.

The relations (\ref{3.5.10}), (\ref{sgly1})
play an important role in the quantum inverse
scattering method \cite{1}. Equations (\ref{3.5.9a}) are the
conditions of factorization of the $S$ matrices in certain exactly solvable
two-dimensional models of quantum field theory (see Ref. \cite{3}).
The matrix representations for the operators
(\ref{3.5.7}) satisfying (\ref{3.5.6}) lead to the formulation of lattice
integrable systems (see, for example, Ref. \cite{28}).
 These questions
will be discussed in more detail in the final section of the review.

Another interesting presentations of the
 quantum operators $L(x)$, which satisfy (\ref{3.5.6}),
 are given in \cite{KashR},
\cite{IsSerg}. In paper \cite{IsSerg},
to construct $L$ operator, we use
equation (\ref{3.5.6}) with Baxterized $R$-matrix which
is defined by means of the
multi-parametric $R$-matrix (\ref{3.4.2}).
These $L$ operators were applied
to the formulation of 3-dimensional integrable models.

The super Yangians $Y(gl(N|M))$ and their representations have been discussed in \cite{Zhang}, \cite{33'}. The quantum
Berezinian for the Yangian (an analog of (\ref{sdeta}) and super-analog of (\ref{ydetq})) was introduced in \cite{33'}.

\subsection{\bf \em The quantum groups $SO_q(N)$ and $Sp_q(2n)$
(B, C, and D series)\label{qBCD}}
\setcounter{equation}0

\subsubsection{Spectral decomposition for
 $SO_q(N)$ and $Sp_q(2n)$ type R-matrices\label{qBCD1}}

In the remarkable paper \cite{10}, the quantum
groups\footnote{We often use the short notation
 $SO_q(N)$ and $Sp_q(N)$ (for quantum groups) instead of more precise
notation for algebras $Fun(SO_q(N))$ and $Fun(Sp_q(N))$.}
$SO_q(N)$ and $Sp_q(N)|_{N=2n}$ were studied as Hopf algebras
with the defining $RTT$ relations
(\ref{3.1.1}). These quantum groups
 are quantum deformations of Lie groups
$SO(N)$ ($B_n$ and $D_n$ series respectively for $N=2n+1$
and $N=2n$) and $Sp(2n)$ ($C_n$ series).
It was shown in \cite{10} that $SO_q(N)$
 and $Sp_q(2n)$ type $R$-matrices (used in the $RTT$
 algebra (\ref{3.1.1})) has the form
 \be
\lb{frtosp}
R_{12} = \sum_{i,j} \, q^{(\delta_{ij} - \delta_{ij'})} \, e_{ii} \otimes e_{jj}
+ \lambda \sum_{i > j} \, e_{ij} \otimes e_{ji} - \lambda \,
 \sum_{i > j} \, q^{\rho_i - \rho_j} \, \epsilon_i \, \epsilon_j \,
 e_{ij} \otimes e_{i'j'} \; ,
\ee
 \be
 \lb{frtosph}
\hat{R}_{12} := P_{12} \, R_{12} =
  \sum_{i,j} \, q^{(\delta_{ij} - \delta_{ij'})} \,
 e_{ij} \otimes e_{ji} +
\lambda \sum_{i < j} \, e_{ii} \otimes e_{jj} - \lambda \,
 \sum_{i > j} \,
 q^{\rho_i - \rho_j} \, \epsilon_i \, \epsilon_j \,
 e_{i'j} \otimes e_{ij'} \; ,
 \ee
where
{\footnotesize
$$
\begin{array}{c}
\epsilon_{i} = +1 \; \forall i \;\; ({\rm for} \;\; SO_{q}(N)) \; , \;\;
 \\ [0.1cm]
\epsilon_{i} = +1 \; (1 \leq i \leq n), \;\;
\epsilon_{i} = -1 \; (n+1 \leq i \leq 2n) \;\; ({\rm for} \;\; Sp_q(2n)),
\end{array}
$$
\be
\lb{rris}
 (\rho_{1}, \dots , \rho_{N}) =  \left\{
\begin{array}{ll}
 (n-\frac{1}{2}, \, n-\frac{3}{2}, \dots , \frac{1}{2}, \, 0 , \,
-\frac{1}{2}, \, \dots , -n+\frac{1}{2} ) \; ,
\;\; & B: (SO_{q}(2n+1)) \; , \\
 (n, \, n-1, \, \dots , 1, \, -1,
 \dots , \, 1-n , \, -n ) \; , \;\;\;\;\;\; & C: (Sp_{q}(2n)) \; , \\
 (n-1, \, n-2, \dots , 1, \, 0, \, 0, \, -1,
 \dots , 1-n ) \; , \;\; & D: (SO_{q}(2n)) \; .
\end{array}
\right.
\ee}
 We deduce these $R$-matrices in Subsection {\bf \ref{mpsosp2}} below.
 The matrices (\ref{frtosp}) satisfy
not only the Yang - Baxter eq. (\ref{ybe}), (\ref{3.1.2i}) but also
the cubic characteristic equation (\ref{3bmw})
(see eq. (\ref{3.1.27}) for $M=3$):
\be
\lb{3.7.1}
(\R -q \hbox{\bf 1})(\R + q^{-1}\hbox{\bf 1})
(\R - \nu \hbox{\bf 1}) = 0
\ee
where $\nu = \epsilon q^{\epsilon -N}$ is a "singlet" eigenvalue of $\R$, and
the case $\epsilon = +1$ corresponds to the orthogonal groups $SO_q(N)$
($B_n$ and $D_n$ series), while the case $\epsilon = -1$ corresponds to the symplectic groups
$Sp_q(2n)$ ($C_n$ series). The projectors (\ref{3.1.28})
arising from the characteristic
equation (\ref{3.7.1}) can be written as follows \cite{10}
\be
\lb{3.7.2}
\begin{array}{c}
\displaystyle
\P^{\pm} = \frac{(\R \pm q^{\mp 1}\hbox{\bf 1})
(\R - \nu \hbox{\bf 1})}
{(q + q^{-1})(q^{\pm 1} \mp \nu )}
\equiv \frac{1}{q+q^{-1}} \left( \pm \R + q^{\mp 1}\hbox{\bf 1}
+ \mu_{\pm} \hbox{\bf K} \right) \; , \\ [0.5cm]
\displaystyle
\P^{0} = \frac{(\R - q \hbox{\bf 1}) (\R + q^{-1}\hbox{\bf 1})}
{(\nu - q)(q^{-1} + \nu )} \equiv
\mu^{-1} \hbox{\bf K}\; ,
\end{array}
\ee
where
$$
 \mu = \frac{(q-\nu)(q^{-1}+\nu)}{\lambda\nu} =
\frac{\lambda + \nu^{-1} -\nu}{\lambda} =
(1 + \epsilon [N-\epsilon]_{q}) \; ,
$$
$$
\mu_{\pm} = \pm \frac{\lambda}{(1 \mp q^{\pm 1} \nu^{-1})} =
\mp \frac{\nu \pm q^{\mp 1}}{\mu} \; , \;\;\;\;\;\;\;
\lambda := q - q^{-1} \; .
$$
We also give the relations between the
 parameters $\nu, \; \mu, \; \mu_{\pm}$
 that we introduced:
$$
q \mu_{+} - q^{-1} \mu_{-} = \nu (\mu_{+} + \mu_{-}) \; , \;\;
\mu_{+} + \mu_{-} = - \frac{q + q^{-1}}{\mu} \; ,
$$
which are helpful in various calculations with projectors
(\ref{3.7.2}). For convenience, we define
 in (\ref{3.7.2}) the operator $\hbox{\bf K}^{i_1 i_2}_{j_1 j_2}$,
which projects $\R$ onto the "singlet" eigenvalue $\nu$:
\be
\lb{3.7.3}
 \K \, \R = \R \, \K  = \nu \, \K \; , \;\;\;  (\K^{2} = \mu \K) \; .
\ee
Then, the characteristic equation
(\ref{3.7.1}) is written in another form  (cf. (\ref{3.3.7}))
\be
\lb{3.7.4}
 \R - \R^{-1} - \lambda + \lambda \, \hbox{\bf K}  = 0 \; .
\ee
The spectral decomposition (\ref{3.1.29}) for
the $SO_q(N)$ and $Sp_q(2n)$ type $R$-matrices is
$$
\R = q \; \P^{+} - q^{-1} \; \P^{-} + \nu \; \hbox{\bf K} \; .
$$

Note that in the semiclassical limit (\ref{3.2.1}), when
$q = e^{h} \to 1$,
the characteristic equation
(\ref{3.7.4})
is reduced to the relation:
\be
\lb{3.7.6}
\frac{1}{2} (r_{12} + r_{21}) =
P_{12} - \epsilon \K^{(0)}_{12} \; ,
\ee
 where $(P_{12})^{i_1 i_2}_{j_1 j_2} \equiv
 (P)^{i_1 i_2}_{j_1 j_2}
  = \delta^{i_1}_{j_2} \delta^{i_2}_{j_1}$ is
the permutation matrix,
and in the right hand side of (\ref{3.7.6}) we
 obtain split Casimir operators for $so$ and $sp$ Lie algebras
 (see \cite{IsKriv}, and refs. therein).
Thus, as in the $GL_{q}(N)$ case (3.3.8), the semiclassical limit
(\ref{3.7.6}) of the characteristic equation fixes the $ad$-invariant part of the
classical $r$ matrix. Here we have used an expansion of the matrix
$\K = \K^{(0)} + h \K^{(1)} + O(h^{2})$, where the first term is
\be
\lb{3.7.6a}
(\K^{(0)})^{i_{1}i_{2}}_{j_{1}j_{2}} =
(C_{0})^{i_{1}i_{2}}  (C_{0}^{-1})_{j_{1}j_{2}} \Rightarrow
\K^{(0)}_{12} = C_{0}^{12 \rangle} \, (C_{0}^{-1})_{\langle 12} \: .
\ee
The matrices $(C_{0})^{ij}$: $(C_{0})^{2} = \epsilon$,
$(C_{0})^{t} = \epsilon C_{0}$ are the metric (symmetric) and symplectic
(antisymmetric) matrices, respectively for the groups $SO(N)$ and $Sp(2n)$.
The semiclassical expansion for the projectors (\ref{3.7.2}) and
(\ref{3.7.5}) has the form
\be
\lb{3.7.7}
\begin{array}{c}
\P_{cl}^{\pm} =
\frac{1}{2} \left( ( \hbox{\bf 1} \pm P) \pm
h P \tilde{r} - (1 \pm \epsilon) \P_{cl}^{0} \right)
\; , \\ \\
\P_{cl}^{0} =
\frac{\epsilon}{N} \left( \K^{(0)}
+ h \K^{(1)} \right)
\end{array}
\ee
where the semiclassical matrix $\tilde{r}$ (\ref{3.2.8})
(which satisfies the modified classical Yang-Baxter equation) is
given by the formula:
$$
\tilde{r} = r_{12} - P_{12} + \epsilon K^{(0)}_{12} =
- r_{21} + P_{12} - \epsilon K^{(0)}_{12} \; .
$$
The ranks of the quantum projectors (\ref{3.7.2}) are equal (for $q$
which is not the root of unity) to the ranks of the projectors (\ref{3.7.7}),
which are readily calculated in the classical limit $h=0$. Accordingly, we have \cite{10}: \\
1) for the groups $SO_q(N)$
\be
\lb{3.7.8}
{\rm rank} (P^{(+)})= \frac{N(N+1)}{2} -1, \;
{\rm rank} (P^{(-)})= \frac{N(N-1)}{2} , \; {\rm rank} (P^{(0)})=1  \;  ;
\ee
2)  for the groups $Sp_q(2n)$
\be
\lb{3.7.9}
{\rm rank} (P^{(+)})= \frac{N(N+1)}{2} , \;
{\rm rank} (P^{(-)})= \frac{N(N-1)}{2} - 1, \; {\rm rank} (P^{(0)})=1  \;  .
\ee
Since the rank of projector $P^{(0)}$ is equal to $1$, we can write
\be
\lb{3.7.9i}
(P^{(0)})^{i_1i_2}_{\; j_1 j_2} = \alpha C^{i_1i_2} \overline{C}_{\; j_1 j_2}
\;\;\; \Rightarrow \;\;\;
\K^{i_1i_2}_{\; j_1 j_2} = C^{i_1i_2} \overline{C}_{\; j_1 j_2}
\ee
where in view of the second equation in (\ref{3.7.3})  we have
$\overline{C}_{\; i_1 i_2} C^{i_1 i_2} = \mu = \alpha^{-1}$.

\subsubsection{Quantum algebras $Fun(SO_{q}(N))$, $Fun(Sp_{q}(2n))$
and their dual algebras}

The number of generators $T^{i}_{j}$
$(i,j=1,...,N)$ for the algebras $Fun(SO_{q}(N))$ and $Fun(Sp_{q}(2n))$
($2n=N$), which satisfy $RTT$-relations (\ref{3.1.1}):
 \be
 \lb{3.1.1b}
R^{i_{1}i_{2}}_{j_{1}j_{2}} \,
T^{j_{1}}_{k_{1}} \, T^{j_{2}}_{k_{2}} =
T^{i_{2}}_{j_{2}} \, T^{i_{1}}_{j_{1}} \,
R^{j_{1}j_{2}}_{k_{1}k_{2}} \; ,
 \ee
coincides with dimensions of the groups
 $SO(N)$ and $Sp(2n)$ in the undeformed case,
since for $T^{i}_{j}$
the following subsidiary conditions are imposed:
\be
\lb{3.7.10}
TCT^{t}C^{-1} = CT^{t}C^{-1}T =1   \Rightarrow
\ee
\be
\lb{3.7.11}
T_{1}T_{2} \, C^{12 \rangle} = C^{12 \rangle} \; , \;\;
C^{-1}_{\langle 12} \, T_{1}T_{2} =  C^{-1}_{\langle 12} \; .
\ee
These relations directly generalize the classical conditions for the
elements of the groups $SO(N)$ and $Sp(2n)$. The matrices
$C^{ij}, \; C^{-1}_{kl}$, which are
understood in (\ref{3.7.11}) as elements in $V_N \otimes V_N$
($1$ and $2$ label the spaces $V_N$)
are the $q$ analogs of the metric and symplectic matrices
$C_{0}$ for $SO(N)$ and $Sp(N)$,
respectively. The explicit form of these matrices, which is given in
Ref. \cite{10} (see
also Sec. {\bf \ref{mpsosp}}), is not important for us, but we stress that the
following equation holds
\be
\lb{3.7.12}
C^{-1}= \epsilon C \; ,
\ee
where $\epsilon=+1$ and $\epsilon=-1$ respectively for
$SO_q(N)$ and $Sp_q(N)$ cases.
Substituting the $R$-matrix representations
(\ref{3.1.18}) for $T^{i}_{j}$ in the relations
(\ref{3.7.10}), we obtain the following conditions on the $R$-matrices:
\be
\lb{3.7.13}
R_{12} = C_{1} (R_{12}^{t_{1}})^{-1}C_{1}^{-1} =
C_{2} (R_{12}^{-1})^{t_{2}}C_{2}^{-1} \; ,
\ee
where, as usual, $C_{1} = C \otimes I$ and $C_{2} = I \otimes C$. As consequences of
(\ref{3.7.13}) we have the equation
\be
\lb{3.7.14}
R_{12}^{t_{1}t_{2}} = C_{1}^{-1}C_{2}^{-1} R_{12} C_{1} C_{2}  \; ,
\ee
and also subsidiary conditions on the generators of the dual algebra
(\ref{3.1.20b}), (\ref{3.1.20}):
\be
\lb{3.7.14a}
L^{\pm}_{2} \, L^{\pm}_{1} \, C^{12 \rangle} = C^{12 \rangle} \; , \;\;
C^{-1}_{\langle 12} \, L^{\pm}_{2} \, L^{\pm}_{1} = C^{-1}_{\langle 12} \; .
\ee

The semiclassical analogs of the conditions
(\ref{3.7.13}) and (\ref{3.7.14}) have the form
$$
r_{12} = -(C_{0})_{1} r_{12}^{t_{1}} (C_{0})_{1}^{-1} =
- (C_{0})_{2} r_{12}^{t_{2}} (C_{0})_{2}^{-1} =
(C_{0})_{1} (C_{0})_{2} r_{12}^{t_{1}t_{2}}
 (C_{0})_{1}^{-1} (C_{0})_{2}^{-1}  \; .
$$

It follows from Eqs. (\ref{3.7.10}) and (\ref{3.7.12}) that the antipode
$S(T) = C \, T^{t} \, C^{-1}$ for the Hopf algebras $Fun(SO_{q}(N))$
and $Fun(Sp_{q}(N))$ satisfy the relation
\be
\lb{3.7.15}
S^{2}(T) = (CC^{t}) T (CC^{t})^{-1} \; ,
\ee
which is analogous to (\ref{3.1.7}). Thus, the matrix $D$ that defines the
quantum trace for the quantum groups of the $SO$ and $Sp$ series
can be chosen in the form
\be
\lb{3.7.16}
D = \epsilon \nu \, CC^{t} \;\; \Leftrightarrow \;\; D^{i}_{j} =
\nu \, C^{ik}\, C^{-1}_{jk} \; .
\ee
 where we take into account (\ref{3.7.12}).
Here we choose the numerical factor $\epsilon \nu$ in order to relate (\ref{3.7.16})
with the general definitions of $D$-matrix (\ref{dmatr}), (\ref{qtrs}).

We now note that the matrix
$C^{12 \rangle} \, C^{-1}_{\langle 12} \in Mat(N) \otimes Mat(N)$
projects any vector $X^{12 \rangle}$ onto the vector $C^{12 \rangle}$,
i.e., the rank of the
projector $C^{12 \rangle} C^{-1}_{\langle 12}$ is 1.
In addition, from (\ref{3.7.11}) we have
$$
C^{12 \rangle} C^{-1}_{\langle 12} \, T_1 \, T_2
= T_1 \, T_2 \,  C^{12 \rangle} C^{-1}_{\langle 12} \; ,
$$
which means that the projector $C^{12 \rangle} C^{-1}_{\langle 12}$
should be a polynomial in $\R$.
Therefore $C^{12 \rangle} C^{-1}_{\langle 12} \sim \P^{0}_{12}$,
and, as it was established in Ref. \cite{10}
(cf. (\ref{3.7.9i})),
\be
\lb{3.7.17}
C^{12 \rangle} C^{-1}_{\langle 12}  \equiv \K_{12}  \; .
\ee
Using this relation,
$RTT$ relations (\ref{3.1.1b}) and equations (\ref{3.7.4})
one can deduce
\be
\lb{3.7.18}
T_1 \, T_2 \, \K_{12} = \K_{12} \, T_1 \, T_2 =
\tau(T) \; \K_{12} \; ,
\ee
where we defined the scalar element
$\tau = \mu^{-1}C^{-1}_{\langle 12} \, T_1 T_2 \, C^{12 \rangle}$.
Comparing of Eq. (\ref{3.7.18}) with Eqs.
(\ref{3.7.10}) and
(\ref{3.7.11}) we conclude that $\tau =1$. Therefore, for the correct
definition
of the quantum groups $SO_{q}(N)$ and $Sp_{q}(N)$ we should require the
 centrality of the element $\tau$ in the $RTT$ algebra (the
 centrality of the element $\tau$ is discussed below
 after eq. (\ref{cen2})).

We note that eqs. (\ref{3.7.3}), (\ref{3.7.17}) are equivalent to
the relations
\be
\lb{3.7.18a}
\R_{12} \, C^{12 \rangle} = \nu \,  C^{12 \rangle} \; , \;\;\;
C^{-1}_{\langle 12} \, \R_{12} = \nu \, C^{-1}_{\langle 12} \; ,
\ee
which give the possibility to rewrite conditions (\ref{3.7.14a}) for the
generators of the reflection
equation algebras (\ref{3.1.21a}):
$$
L_1 \, \R_{12} \, L_1 \, C^{12 \rangle} = \nu \, C^{12 \rangle} \; , \;\;\;
C^{-1}_{\langle 12} \, L_1 \, \R_{12} \, L_1  = \nu \, C^{-1}_{\langle 12} \; ,
$$
$$
\overline{L}_2 \, \R_{12} \, \overline{L}_2 \, C^{12 \rangle} = \nu \, C^{12 \rangle} \; , \;\;\;
C^{-1}_{\langle 12} \, \overline{L}_2 \, \R_{12} \, \overline{L}_2  =
\nu \, C^{-1}_{\langle 12} \; .
$$

We now present some important relations for the matrices $\R$ and $\K$;
many of them are given, in some form or other, in Ref. \cite{10}.
We note first that in accordance with (\ref{3.1.4})
we have
\be
\lb{3.7.19}
\K_{12} \, \R_{23} \, \R_{12}   =
\R_{23} \, \R_{12} \, \K_{23} \Leftrightarrow
\R_{12} \, \R_{23} \, \K_{12}   =  \K_{23} \, \R_{12} \, \R_{23} .
\ee
Further, from Eqs. (\ref{3.7.13}) and (\ref{3.7.17}) (or substituting
the matrix representations
(\ref{3.1.18}) in (\ref{3.7.18}) for $\tau =1$), we obtain
\be
\lb{3.7.20}
\begin{array}{c}
\R^{\pm 1}_{12}  \, {\R_{23}}^{\pm 1} \, \K_{12}   =
P_{12} \, P_{23} \, \K_{12}  = \K_{23} \, P_{12} \, P_{23}  \; , \\ \\
\K_{12} \, {\R_{23}}^{\pm 1} \, \R^{\pm 1}_{12}   =
\K_{12} \, P_{23} \, P_{12}  = P_{23} \, P_{12} \, \K_{23} \; .
\end{array}
\ee
A consequence of these relations is the equations
\be
\lb{3.7.21}
\begin{array}{c}
{\R_{23}}^{\pm 1} \, \K_{12} \, {\R_{23}}^{\pm 1} =
\R^{\mp 1}_{12} \,  \K_{23} \, \R^{\mp 1}_{12}
\Leftrightarrow   \R_{12} \, \R_{23} \, \K_{12}   =
\K_{23} \, \R^{-1}_{12} \, {\R_{23}}^{-1} \; , \\ \\
  \R_{23} \, \R_{12} \, \K_{23}  =  \K_{12} \, {\R_{23}}^{-1} \, \R^{-1}_{12}  \; .
\end{array}
\ee
In particular, taking into account the characteristic equation
(\ref{3.7.4}), we obtain the identity
$$
(\R_{12}  - \lambda) \, \K_{23} \, (\R_{12}  - \lambda) =
(\R_{23} - \lambda) \, \K_{12} \, (\R_{23} - \lambda)
$$
or
\be
\lb{3.7.22}
\begin{array}{c}
\R_{12} \, \K_{23} \R_{12}  = {\R_{23}}^{-1} \, \K_{12} \, {\R_{23}}^{-1} = \\ \\
\R_{23} \, \K_{12} \, \R_{23} + \lambda ( \R_{12} \,\K_{23}  -
\K_{12} \, \R_{23} - \R_{23} \, \K_{12}   +  \K_{23} \, \R_{12} ) +
\lambda^{2} (\K_{12}  - \K_{23})  \; ,
\end{array}
\ee
which will be used in Sec. {\bf \ref{ospbax}}. Equation (\ref{3.7.17})
leads to the identities
\be
\lb{3.7.23}
\K_{12} \, \K_{23}  = \K_{12} \, P_{23} \, P_{12}  =
P_{23} \, P_{12} \, \K_{23} \; , \;\;
\K_{23} \, \K_{12}   = P_{12} \, P_{23} \, \K_{12}  =  \K_{23} \, P_{12} \, P_{23}  \; ,
\ee
from which we immediately obtain
\be
\lb{3.7.24}
\K_{12} \, \K_{23} \, \K_{12}   = \K_{12}  \; , \;\;
\K_{23} \, \K_{12} \, \K_{23} =  \K_{23}  \; .
\ee
We now compare the relations (\ref{3.7.20}) and (\ref{3.7.23}).
The result of this comparison is the equations
\be
\lb{3.7.25}
\begin{array}{c}
{\R_{23}}^{\pm 1} \, \R^{\pm 1}_{12} \, \K_{23}  =
\K_{12} \, \K_{23}  =  \K_{12} \, {\R_{23}}^{\pm 1} \R^{\pm 1}_{12}  \; ,
\\ \\
\R^{\pm 1}_{12} \, {\R_{23}}^{\pm 1} \, \K_{12}   =
\K_{23} \, \K_{12}  =  \K_{23} \, \R^{\pm 1}_{12} \, {\R_{23}}^{\pm 1} ,
\end{array}
\ee
We now apply to the first of the chain of equations in
(\ref{3.7.25}) the matrix $\K_{12}$ from the right (or $\K_{23}$ from the left) and take
into account (\ref{3.7.3}) and (\ref{3.7.24}). We then obtain
\be
\lb{3.7.26}
\K_{23} \, \R^{\pm 1}_{12} \, \K_{23}  = \nu^{\mp 1} \K_{23} \; , \;\;
\K_{12} \, {\R_{23}}^{\pm 1} \, \K_{12}   = \nu^{\mp 1} \K_{12}   \; .
\ee

The braid relation (\ref{3.1.3})
and eqs. (\ref{3.7.3}), (\ref{3.7.4}), (\ref{3.7.26})
define the $R$- matrix representation of the
Birman - Murakami - Wenzl algebra \cite{W2}
(see also Sec. {\bf \ref{bmwalg}} below). Eqs. (\ref{3.7.19}),
(\ref{3.7.21}), (\ref{3.7.22}), (\ref{3.7.24}) and (\ref{3.7.25})
directly follow from this definition.
As we shall see in Sec. {\bf \ref{ospbax}}, the relations
for the Birman - Murakami - Wenzl algebra will be sufficient for the construction of
$SO_{q}(N)$ and $Sp_{q}(2n)$-symmetric Baxterized $\R(x)$ matrices. The relations
(\ref{3.7.19}), (\ref{3.7.21}), and (\ref{3.7.24})-(\ref{3.7.26}) have a natural
graphical representation in the form of relations for braids and links
if we use the diagrammatic technique (only 3 of these operators are independent
in view of (\ref{3.7.4}))

\unitlength=6mm
\begin{picture}(17,4)
\put(0,1.9){$\hbox{\bf R} =$}
\put(1.3,3){\line(1,-1){1}}
\put(2.1,1.8){\vector(-1,-1){0.8}}
\put(2.3,2){\vector(1,-1){1}}
\put(2.5,2.2){\line(1,1){0.8}}

\put(4,1.9){$\hbox{\bf R}^{-1} =$}
\put(6,3){\line(1,-1){0.8}}
\put(7,2){\vector(-1,-1){1}}
\put(7.2,1.8){\vector(1,-1){0.8}}
\put(7,2){\line(1,1){1}}

\put(8.7,1.9){$I_{1}I_{2} =$}
\put(11,3){\vector(0,-2){2}}
\put(11.8,3){\vector(0,-2){2}}

\put(13,1.9){$\hbox{\bf K} =$}
\put(15.3,3){\oval(1.5,1)[b]}
\put(14.6,2.8){\vector(0,-1){0.1}}
\put(16,2.8){\vector(0,-1){0.1}}
\put(15.3,1){\oval(1.5,1)[t]}
\put(14.55,1){\vector(0,-1){0.1}}
\put(16.05,1){\vector(0,-1){0.1}}


\end{picture}
\vspace{-1cm}
\be
\lb{figa}
{}
\ee

\subsubsection{Quantum traces and quantum hyperplanes
 for $SO_{q}(N)$ and $Sp_{q}(N)$}

We now give some important relations for the quantum trace
(\ref{3.1.13}) corresponding to the quantum groups $SO_{q}(N)$ and $Sp_{q}(N)$.
Similar relations for the $q$- trace (\ref{3.1.15}) can be derived in exactly
the same way. From the definitions of the matrix
$\K$ (\ref{3.7.17}) and the matrix $D$ (\ref{3.7.16}), we obtain
\be
\lb{3.7.27}
Tr_{q2}(\K_{12}  ) = \nu \, I_{(1)} \; .
\ee
We use the relations (\ref{3.7.14}) and the definition of the quantum trace
(\ref{3.1.13}) with matrix $D$ (\ref{3.7.16}); then, for an arbitrary quantum matrix
$E^{i}_{j}$, we obtain the relations:
\be
\lb{3.7.28}
\begin{array}{c}
\nu \, \R^{n}_{12} \, E_1 \, \K_{12}  = Tr_{q2}(\K_{12} \, E_1 \,
\R^{n}_{12}) \K_{12}  \; , \\ \\
\nu \, \K_{12} \, E_1 \, \R^{n}_{12} = \K_{12} \,
Tr_{q2}(\R^{n}_{12} \, E_1 \,  \K_{12} ) \; , \forall n \; ,
\end{array}
\ee
\be
\lb{3.7.29}
\nu \, \K_{12} \, E_1 \, \K_{12}  = Tr_{q}( E) \K_{12}  \; .
\ee
Calculating the $Tr_{q2}$ of (\ref{3.7.28}), we deduce
\be
\lb{3.7.30}
Tr_{q2}(\R^{n}_{12} \, E_1 \, \K_{12}  ) = Tr_{q2}(\K_{12} \, E_1 \, \R^{n}_{12})
\; , \;\; \forall n \; .
\ee
Further, from the first identity of (\ref{3.7.26}), averaging it by means
of $Tr_{q2}$, we readily obtain for the algebras $Fun(SO_{q}(N))$ and $Fun(Sp_{q}(N))$
analogs of (\ref{3.1.16}). These take the form
\be
\lb{3.7.31}
Tr_{q2}(\R^{\pm 1}_{12} ) \equiv \epsilon \, \nu \, Tr_{2}(CC^{t} \R^{\pm 1}_{12}) =
 \nu^{1 \mp 1} I_{(1)}.
\ee
Using this relation and Eq. (\ref{3.7.4}), we can calculate
\be
\lb{3.7.32}
Tr_{q} (I) = Tr (D) = \nu (1 + \epsilon \, [N-\epsilon]_{q}) = \nu \, \mu \; .
\ee
We now separate irreducible representations for the left adjoint comodules
(\ref{3.1.11}). For an arbitrary $N \times N$ quantum matrix $E^{i}_{j}$, we have
\be
\lb{3.7.33}
\begin{array}{c}
E_1 = \nu^{-1} Tr_{q2}(E_1 \, \K_{12} ) =
E^{(0)}_1 + E^{(+)}_1 + E^{(-)}_1 \; ,
\end{array}
\ee
where $E^{(i)}_1 = \nu^{-1} \, Tr_{q2}(\P^{i}_{12} \, E_1 \, \K_{12}) =
\nu^{-1} \, Tr_{q2}(\K_{12} \, E_1 \, \P^{i}_{12})$.
It is obvious that the tensors $E^{(i)}, \; (i=\pm, \; 0)$ are
invariant with respect to the adjoint coaction (\ref{3.1.11}) and
$Tr_{q2}(\P^{(j)}E^{(i)}\K) = 0$ (if $i \neq j$) by virtue of (\ref{3.7.28}).
Thus, (\ref{3.7.33}) is the required decomposition of the adjoint
comodule $E$ into irreducible components. It is clear that the component $E^{(0)}$
is proportional to the unit matrix, $(E^{(0)})^{i}_{j} = e^{(1)} \cdot \delta^{i}_{j}$
($e^{(1)}$ is a constant), and,
thus, applying $Tr_{q1}$ to (\ref{3.7.33}), we obtain
\be
\lb{3.7.34}
Tr_{q}(E) =  e^{(1)} Tr_{q}(I) =
 \nu \mu E^{(1)}  \; ,
\ee
where we have used the property (\ref{3.3.16}), which also holds for the case of the
quantum groups $SO_{q}(N)$ and
$Sp_{q}(2n)$.
To conclude this subsection, we note that, as in the case of the
linear quantum groups, we can define fermionic and bosonic quantum hyperplanes
covariant with respect to the coactions of the groups $SO_{q}(N)$ and
$Sp_{q}(2n)$. Taking into account the ranks of the projectors
(\ref{3.7.8}) and (\ref{3.7.9}), we can formulate definitions of the hyperplanes
for  $SO_{q}(N)$ ($\epsilon =1$) and for $Sp_{q}(N)$ ($\epsilon =-1$) in the form
\be
\lb{3.7.35}
(\P^{-} + (\epsilon-1) \K) x x'=0
\ee
for the bosonic hyperplane [number of relations $N(N-1)/2$] and
\be
\lb{3.7.36}
(\P^{+} + (\epsilon+1) \K) x x'=0
\ee
for the fermionic hyperplane [number of relations $N(N+1)/2$].
For all these algebras, the elements $\K x x'$ are central elements, and it is
obvious that for $Sp_{q}(N)$ bosons and $SO_{q}(N)$ fermions we have $\K x x' = 0$.
It is interesting that the projectors $\P^{\pm}$
(\ref{3.7.2}) can be represented as
\be
\lb{3.7.5}
\P^{\pm} =
\frac{1}{q+q^{-1}}
\left( \pm \R' + q^{\mp 1}\hbox{\bf 1} \right)
- \frac{1}{2\mu} (1 \pm \epsilon )   \hbox{\bf K} \; ,
\ee
where the matrix
$$
\hat{R}' = \R - \frac{1}{2} [ \mu_{-} (1 + \epsilon) +
\mu_{+} (\epsilon -1) ] \K \; ,
$$
satisfies the Hecke condition (\ref{3.3.7}). However,
using (\ref{3.7.19}) -- (\ref{3.7.26}) one can
directly check that $\hat{R}'$ does not obey the Yang-Baxter eq. (\ref{3.1.2i}).

Note that
the conditions (\ref{3.7.10}) and (\ref{3.7.11}) can be understood as conditions of
invariance of the quadratic forms $x_{(1)} C^{-1} x_{(2)}$ and $y_{(1)} C y_{(2)}$
with respect to
left and right transformations of the hyperplanes $x_{(k)}, \; y_{(k)}$:
$$
x_{(k)}^{i} \rightarrow
T^{i}_{j} \otimes x_{(k)}^{j} \; , \;\;
{y_{(k)}}_{i} \rightarrow
{y_{(k)}}_{j} \otimes T^{j}_{i} \; .
$$

\subsection[\bf \em The multi-parameter
deformations $SO_{q,a_{ij}}$,
$Sp_{q,a_{ij}}$ and 
q-supergroups $Osp_{q}(N|2m)$]{\bf \em The multi-parameter
deformations $SO_{q,a_{ij}}(N)$ and
$Sp_{q,a_{ij}}(2n)$ of orthogonal and symplectic
groups and the q-supergroups $Osp_{q}(N|2m)$\label{mpsosp}}
\setcounter{equation}0

\subsubsection{General multi-parametric $R$-matrices
of the $OSp$-type\label{mpsosp1}}

In this subsection we show that it is possible to define multi-parameter deformations of the quantum groups
$SO_q(N)$ and $Sp_q(2n)$ and also the quantum supergroups
$OSp_{q}(N|2m)$ (as $RTT$-algebras) if we consider for the $R$-matrix the ansatz:
\be
\lb{3.8.01}
\R= \sum_{i,j=1}^K \, a_{ij} \, e_{ij} \otimes e_{ji}  +
\sum_{i < j} \, b_{ij} \, e_{ii} \otimes e_{jj}  +
\sum_{i> j} \, d^i_{j} \, e_{i'j} \otimes e_{ij'}  \; \Rightarrow
\ee
\be
\lb{3.8.1}
\hat{R}^{i_{1},i_{2}}_{j_{1},j_{2}}=
\delta^{i_{1}}_{j_{2}} \delta^{i_{2}}_{j_{1}} \, a_{i_{1}i_{2}} +
\delta^{i_{1}}_{j_{1}} \delta^{i_{2}}_{j_{2}}  \, b_{i_{1}i_{2}} \,
\Theta_{i_{2}i_{1}} +
\delta^{i_{1}i'_{2}} \delta_{j_{1}j'_{2}} \,
d^{i_{2}}_{j_{1}} \, \Theta^{i_{2}}_{\; j_{1}} =
\ee

\unitlength=8mm
\begin{picture}(17,4)
\put(2,1.9){$ \;\;\; =$}
\put(4,3.2){$i_{1}$}
\put(4,3){\line(1,-1){1}}
\put(4,1){\line(1,1){1}}
\put(4,0.5){$j_{1}$}
\put(4.6,2.7){$a_{i_{1}i_{2}}$}
\put(4.5,2.5){\line(1,0){1}}
\put(6,0.5){$j_{2}$}
\put(5,2){\line(1,-1){1}}
\put(5,2){\line(1,1){1}}
\put(6,3.2){$i_{2}$}
\put(6.5,1.9){$+$}
\put(8,2){\line(1,0){1.5}}
{\thicklines
\put(8.9,1.99){\vector(-1,0){0.3}}
}
\put(8.5,2.3){$b_{i_{1}i_{2}}$}
\put(8,3.2){$i_{1}$}
\put(8,0.5){$j_{1}$}
\put(8,3){\line(0,-2){2}}
\put(9.5,1){\line(0,2){2}}
\put(9.5,0.5){$j_{2}$}
\put(9.5,3.2){$i_{2}$}
\put(10.5,1.9){$+$}
\put(11.5,3.2){$i_{1}$}
\put(11.5,0.5){$j_{1}$}
\put(12.5,3){\oval(1.5,1)[b]}
\put(11.9,1.4){\line(1,1){1.2}}
{\thicklines
\put(12.9,2.4){\vector(-1,-1){0.4}}
}
\put(12.5,1){\oval(1.5,1)[t]}
\put(13,1.9){$d^{i_{2}}_{j_{1}}$}
\put(13.2,3.2){$i_{2}$}
\put(13.2,0.5){$j_{2}$}
\end{picture}
\\
where $\Theta^{i}_{\; j} := \Theta_{ij}$,
$j' = K+1-j$; $K=N$ for the groups $SO(N)$,
$K=2n$ for the groups $Sp(2n)$, and $K= N+2m$ for the
groups $O\!sp(N|2m)$. The expression (\ref{3.8.1}) is a natural generalization of the expression
(\ref{3.4.3}) for the multi-parameter $R$-matrix corresponding to the linear quantum groups.
Namely, the third term in
(\ref{3.8.1})   is   constructed   from   the   $SO$-invariant   tensor
$\delta^{i_{1}i_{2}'} \delta_{j_{1}j_{2}'}$, which takes
into account the presence of the invariant metrics for the
considered groups. The functions $\Theta$
are introduced in (\ref{3.8.1}) (they
are indicated as arrows in the graphical representation)
 in order to ensure that the matrix $R_{12} = P_{12}\R$ has
lower triangular block form. This is necessary for the correct definition of the
operators $L^{(\pm)}$ by means of the expressions (\ref{3.1.19}).
We demonstrate below that the ansatz (\ref{3.8.1}) for the
solution of the Yang-Baxter equation (\ref{3.1.2i}) automatically
define the family of the Birman-Murakami-Wenzl $R$- matrices
with fixed parameter $\nu$ which corresponds to the quantum groups
$SO_{q}(N)$, $Sp_{q}(2n)$ and $Osp_{q}(N|2m)$.

We substitute the ansatz (\ref{3.8.1}) for $R$-matrix
in the Yang-Baxter equation (\ref{3.1.2i}).
It is obvious that the first two terms in (\ref{3.8.1}) make contributions to the
Yang-Baxter equation that are analogous to the contributions
of the general $R$-matrix ansatz in the case of the linear quantum
groups (see Sec. {\bf \ref{multpar}}). It is therefore clear that for the
parameters $a_{ij}$ and $b_{ij}$ we reproduce
almost the same conditions (\ref{3.4.4}),
which in the convenient normalization $c=1, \; b = q-q^{-1}$ have the form
\be
\lb{3.8.2}
b_{ij} = b = \lambda \;\; (\forall i,j), \;\;
a_{ii} = a^{0}_{i} = \pm q^{\pm 1} \;\; (i \neq i'), \;\;
a_{ij} \, a_{ji} = 1 \;\; (i \neq j , \; i \neq j') \; .
\ee
Note that the conditions in (\ref{3.8.2}) are somewhat weaker than in
(\ref{3.4.4}) (because of the restrictions $i \neq i'$, $i \neq j'$).
 This is due to the fact that the contributions to the Yang-Baxter
 equation proportional to $a_{ii'}$ begin to be canceled by the
 contributions from the third term in (\ref{3.8.1}). The corresponding condition on $a_{ii'}$
 fulfilling the Yang-Baxter equation can be expressed as follows:
\be
\lb{3.8.3}
\begin{array}{c}
a_{jj'} = \kappa^{-1}_{j}(a^{0}_{j} - b)
\; , \;\;\; a_{j'j} = \kappa_{j}(a^{0}_{j} - b)
\;\;\;\;\; (j \neq j') \; , \;\;\;
\Leftrightarrow \\ \\
a^{0}_{j} \, a_{jj'} = \kappa^{-1}_{j} \; , \;\;\;
 a^{0}_{j} \, a_{j'j} = \kappa_{j}  \;\;\;\;\;\; (j \neq j') \; ,
\end{array}
\ee
where in addition for the constants $a^{0}_{j}$ and $\kappa_{i}$ we have
\be
\lb{3.8.3a}
\kappa_{j} \kappa_{j'} = 1 \; , \;\;\; a^{0}_{j} = a^{0}_{j'} \; .
\ee
Taking into account Eqs. (\ref{3.8.2}), the relations (\ref{3.8.3}) are
equivalent to the pair of possibilities ($ j \neq j'$):
\be
\lb{3.8.4}
\begin{array}{l}
1.) \;\; a^{0}_{j} = q \;\;\; \rightarrow \;\;\;
 a_{j'j} \, \kappa_{j}^{-1} = q^{-1} = a_{jj'} \, \kappa_{j}
\; , \\ \\
2.) \;\; a^{0}_{j} = -q^{-1} \;\;\;  \rightarrow \;\;\;
 a_{j'j} \, \kappa_{j}^{-1} = -q = a_{jj'} \, \kappa_{j} \; .
\end{array}
\ee
We shall see below that if we restrict our consideration to the first
possibility for all $j$ (or only the second possibility),
then we obtain the $R$-matrices for the
quantum groups $SO_{q}(N)$ and $Sp_{q}(2n)$. If, however, we consider the
mixed case, when both possibilities
are satisfied (for different $j$), then we expect
(by analogy with the linear quantum groups;
 see Sec. {\bf \ref{multpar}}) that the
corresponding $R$-matrix will be associated with the supergroups $Osp_{q}(N|2m)$.
The case $j=j'$ is obviously realized only for groups of the
series B ($SO_{q}(2n +1)$) and for the supergroups $Osp_{q}(2n+1|2m)$, and it
follows from the Yang-Baxter equation (\ref{3.1.3}) that
\be
\lb{3.8.5}
a_{jj'} = 1 \; , \;\;\; \kappa_{j} = 1 \; , \;\;\;
{\rm for} \;\; j = j'= \frac{K+1}{2} \; .
\ee
For the groups $SO_{q}(2n)$, $Sp_{q}(2n)$ and $Osp_{q}(2n|2m)$,
the parameter $a_{jj'}\; (j=j')$
is simply absent. Further consideration of the contributions to the Yang-Baxter
equation from the third term in (\ref{3.8.1}) leads to the equations
\be
\lb{3.8.6}
a_{ij} \, a_{i'j} =  \kappa_{j} \; , \;\;\;
a_{ji} \, a_{ji'} =  \kappa^{-1}_{j} \;\;\; (\forall i \neq i') \; ,
\ee
\be
\lb{3.8.7}
\lambda \, d^{j}_{k} \, \kappa_{i} + d^{j}_{i} \, d^{i}_{k} = 0 \; , \;\;
\ee
(there is no summation over repeated indices). The general solution of
Eq. (\ref{3.8.7}) has the form
\be
\lb{3.8.8}
d^{i}_{j} = - \lambda \, \kappa_{i} \, \frac{c_{j}}{c_{i}} \; ,
\ee
where $c_{i}$ are arbitrary parameters. The remaining terms in the Yang-Baxter
equation that do not cancel under the conditions (\ref{3.8.2})-(\ref{3.8.8}) give
recursion relations for the coefficients $c_{i}$:
\be
\lb{3.8.9}
c_{j'}a_{j'j} + \lambda c_{j} \Theta_{j'j} -
\lambda c_{j} \sum_{i>j} \kappa_{i} \frac{c_{i'}}{c_{i}} =
\nu c_{j} \; .
\ee
These relations can be represented graphically in the form

\unitlength=8mm
\begin{picture}(17,5)
\put(2,4){\line(1,-1){1}}
\put(2,2){\line(1,1){1}}
\put(2,1.5){$j$}
\put(2,4.1){\circle*{0.2}}
\put(1.3,4.2){$c_{j'}$}
\put(3,4.05){\oval(2,1)[t]}
\put(2.6,3.7){$a_{j'j}$}
\put(2.5,3.5){\line(1,0){1}}
\put(4,1.5){$j'$}
\put(3,3){\line(1,-1){1}}
\put(3,3){\line(1,1){1}}
\put(4.5,2.9){$+ \;\;\;\; \lambda$}
\put(5.5,4.2){$c_{j}$}
\put(6,4.1){\circle*{0.2}}
 {\thicklines
\put(6,3){\line(1,0){0.7}}
\put(7.5,3){\vector(-1,0){0.9}}
 }
\put(6.75,4){\oval(1.5,1)[t]}
\put(6,1.5){$j$}
\put(6,4){\line(0,-2){2}}
\put(7.5,2){\line(0,2){2}}
\put(7.5,1.5){$j'$}
\put(8.5,2.9){$+$}
\put(9.5,1.5){$j$}
\put(9.8,4.1){\circle*{0.2}}
\put(9.2,4.2){$c_{i'}$}
\put(10.5,4){\oval(1.5,1)[t]}
\put(10.5,4){\oval(1.5,1)[b]}
 {\thicklines
\put(10.6,3.1){\vector(-1,-1){0.1}}
 }
\put(9.9,2.4){\line(1,1){1.2}}
\put(10.5,2){\oval(1.5,1)[t]}
\put(10.7,2.8){$d^{i}_{j}$}
\put(11.2,1.5){$j'$}
\put(12.1,2.9){$ = \; \nu$}
\put(14.5,2.7){\oval(1.5,1)[t]}
\put(13.8,2.9){\circle*{0.2}}
\put(13.5,3.2){$c_{j}$}
\put(13.6,2.1){$j$}
\put(15.1,2.1){$j'$}
\end{picture}

\noindent
Another equivalent forms of (\ref{3.8.9}) are
\be
\lb{3.8.9c}
\sum_{k >m} d^i_{k'} d^k_j = d^i_j \left( \nu - \frac{c_{m'}}{c_m} a_{m'm} -
\lambda \Theta_{m'm} \right) \; ,
\ee
\be
\lb{3.8.9b}
\sum_{k < m} d^i_{k'} d^k_j = d^i_j \left( - \nu^{-1}
+ \frac{c_{m'}}{c_m} \, a_{m \, m'}^{-1} -
\lambda \Theta_{mm'} \right) \; ,
\ee
which are related to each other
by the identity
\be
\lb{3.8.9aa}
 \sum_{k} d^i_{k'} d^k_j = - \lambda \, \mu \, d^i_j \; , \;\;
 (\mu :=  (\lambda - \nu + \nu^{-1})/ \lambda) \; ,
\ee
using below.
Now the $R$-matrix (\ref{3.8.01}) is represented in the form
\be
\lb{Rgen3}
\hat{R} = a_{ij} \, e_{ij} \otimes e_{ji}  +
\lambda \, \Theta_{ji} \, e_{ii} \otimes e_{jj} + \Theta_{ij} \, d^i_j  \,
 e_{i'j} \otimes e_{ij'} \; ,
\ee
where the summation over the indices $i,j$ is assumed and
parameters $a_{ij}$ and $d^i_j$ are fixed by the conditions
(\ref{3.8.2}) -- (\ref{3.8.6}), (\ref{3.8.8}) and (\ref{3.8.9}).
This $R$-matrix satisfies the Yang-Baxter equation
(\ref{3.1.2i}) and additional relations
(cf. (\ref{3.7.3}), (\ref{3.7.4}) and (\ref{3.7.26}))
\be
\lb{char1}
\R^2 - \lambda \R - {\bf 1} = - \lambda \nu \, \K \; , \;\;\; \K \R = \R \K = \nu \K \; ,
\ee
\be
\lb{Bmw}
\K_{12} \, \R^{\pm 1}_{23} \, \K_{12} =
\nu^{\mp 1} \, \K_{12} \; ,  \;\;\; \K^2 = \mu \, \K \; ,
\ee
where we have introduced the rank 1 matrix:
\be
\lb{3.8.10a}
\K := -\lambda^{-1} \,
\sum_{i,j} \, d^i_j \, e_{i'j} \otimes e_{i j'} = \sum_{i,j} \, \kappa_i \,
\frac{c_j}{c_i} \, e_{i'j} \otimes e_{i j'} \; \Leftrightarrow
\ee
\be
\lb{3.8.10}
\K^{i_1 i_2}_{j_1 j_2}  = C^{i_1 i_2} \, C_{j_1 j_2} \; , \;\;\;
C^{ij} = \epsilon \delta^{ij'} \frac{\kappa_j}{c_{j}} \; , \;\;
C_{ij} = \frac{1}{\epsilon} \delta_{ij'} c_{i} \; .
\ee
To prove relations (\ref{char1}), (\ref{Bmw}) we have used
the definitions of $a_{ij}$ (\ref{3.8.2}) -- (\ref{3.8.5}), $d^i_j$ (\ref{3.8.8})
and take into account the identities (\ref{3.8.9c}) -- (\ref{3.8.9aa}) and
$$
\Theta_{ki'} \Theta_{kj} =
\Theta_{ki'} \Theta_{kj} (\Theta_{i'j} + \Theta_{ji'} + \delta_{i'j}) =
\Theta_{i'j} \Theta_{ki'} + (\Theta_{ji'}  + \delta_{i'j}) \Theta_{kj} \; .
$$
Thus, the $R$-matrix (\ref{Rgen3}) with constraints
(\ref{3.8.2}) -- (\ref{3.8.6}), (\ref{3.8.8}) and (\ref{3.8.9})
automatically leads to the $R$-matrix representation of the
Birman - Murakami - Wenzl algebra (the definition of this algebra
is given below in Sec. {\bf \ref{bmwalg}}).
In (\ref{3.8.10}) we define the quantum metric,
or quantum symplectic, matrices $C$ (cf. (\ref{3.7.9i}), (\ref{3.7.17})).
The parameter $\epsilon$ (see Sec. {\bf \ref{mpsosp}}) is introduced
in (\ref{3.8.10}) in order
to match the definition of the matrices $C$ to the study of Ref. \cite{10},
where $\epsilon = \pm 1$.


Note that the conditions
(\ref{3.8.2}) - (\ref{3.8.6}) can be solved as:
\be
\lb{Rgen}
a_{ij} = (a_i^0)^{(\delta_{ij} - \delta_{ij'})} \, \frac{f_{ij}}{f_{ji}} \; , \;\;\;
\frac{f_{ij} f_{i'j}}{f_{ji} f_{j i'}} =\kappa_j \;\; (\forall i \neq i')
\;\; \Rightarrow \;\; \kappa_j = \frac{f_{j'j}}{f_{j j'}} \; ,
\ee
and, after substitution of (\ref{Rgen})
in (\ref{Rgen3}), one can observe that the $R$-matrix (\ref{Rgen3}) is
\be
\lb{Rgen11}
\hat{R} = \sum_{i, j} \, (a_i^0)^{(\delta_{ij} - \delta_{ij'})} \, \frac{f_{ij}}{f_{ji}} \,
 e_{ij} \otimes e_{ji}  +
\lambda \sum_{i < j} \, e_{ii} \otimes e_{jj} - \lambda \,
 \sum_{i > j} \frac{f_{i'i}}{f_{i i'}} \, \frac{ c_j}{ c_i} \,
 e_{i'j} \otimes e_{ij'} \; ,
\ee
and produced by the twisting (\ref{twistgl}) from the matrix
\be
\lb{Rgen1}
\hat{R} = \sum_{i, j} \, (a_i^0)^{(\delta_{ij} - \delta_{ij'})} \,
 e_{ij} \otimes e_{ji}  +
\lambda \sum_{i < j} \, e_{ii} \otimes e_{jj} - \lambda \,
 \sum_{i > j} \, \frac{ \tilde{c}_j}{ \tilde{c}_i} \,
 e_{i'j} \otimes e_{ij'} \; ,
\ee
where $\tilde{c}_i = f_{ii'} c_i$
and the parameters $a^0_j = a^0_{j'}$, $c_j$ are determined in (\ref{3.8.4}),
(\ref{3.8.5}) and (\ref{3.8.9}).
In this case the relations (\ref{compat}) lead to
additional conditions on the twisting parameters $f_{ij}$:
\be
\lb{ffff}
f_{ij} f_{i'j} = \kappa_j v_j \; , \;\;\; f_{ji} f_{ji'} = v_j \;, \;\; \forall i \; ,
\ee
which are consistent with (\ref{3.8.6}), (\ref{Rgen}). It is evident that
for $R$-matrix (\ref{Rgen1}) the analogs of matrices
(\ref{3.8.10a}), (\ref{3.8.10}) are
$$
\K^{i_1 i_2}_{j_1 j_2}  = \sum_{i,j} \,
\frac{\tilde{c}_j}{\tilde{c}_i} \,
(e_{i'j})^{i_1}_{j_1} \otimes (e_{i j'})^{i_2}_{j_2} =
\widetilde{C}^{i_1 i_2} \, \widetilde{C}_{j_1 j_2} \;\; \Rightarrow \;
$$
\be
\lb{3.8.10b}
\widetilde{C}^{ij} = \epsilon \delta^{ij'} \frac{1}{\tilde{c}_{j}} \; , \;\;
\widetilde{C}_{ij} = \frac{1}{\epsilon} \delta_{ij'} \tilde{c}_{i}
\;\; \Rightarrow \;\; \widetilde{C}^{ij} \, \widetilde{C}_{jk} = \delta^i_k \; .
\ee

Now we show that
the constant $\nu$ is fixed by the relations (\ref{3.8.9}) uniquely.
We consider the solution of Eqs. (\ref{3.8.9}), which are
written in the form
\be
\lb{3.8.11}
\gamma_{j} \ a_{j'j} \ \kappa^{-1}_{j} + \lambda \ \Theta_{j'j} -
\lambda \ \sum_{i=j+1}^{K}  \ \gamma_{i} =
\nu \; ,
\ee
where
\be
\lb{3.8.12}
\;\; \gamma_{j}  = \frac{c_{j'}}{c_{j}} \ \kappa_j
= \frac{\tilde{c}_{j'}}{\tilde{c}_{j}}  = \frac{1}{\gamma_{j'}}  \; .
\ee
Equation (\ref{3.8.11}) is readily solved by the changing of variables:
$\gamma_{i} \;\; \to \;\; X_{i}$,
$$
X_{j} := q^{2j} \sum_{i=j+1}^K  \gamma_{i} \; , \;\; (X_K = 0) \; ,
$$
where the inverse transformation is
$\gamma_{j} = q^{-2j}(q^{2}X_{j-1} - X_{j})$ and we fix $\nu$
in (\ref{3.8.11}) by
taking into account the properties (\ref{3.8.12}).


\subsubsection{The case of $SO_{q}(N)$ and $Sp_{q}(N)$
 groups\label{mpsosp2}}

 First, we consider the
possibility 1 in (\ref{3.8.4}).
The possibility 2 gives,
in view of a symmetry of equation (\ref{3.8.11}), an analogous result
except for the substitution $q \rightarrow -q^{-1}$. The corresponding
form of eq. (\ref{3.8.12}) for $j > j'$ is
$$
q \, (X_{j-1} - X_j) = q^{2j} \, \nu \; ,
$$
and  we obtain the solution:
\be
\lb{3.8.11c}
X_j = \nu q^{2K-1} \, \frac{1-q^{-2(K-j)}}{1-q^{-2}} \;\; \Leftrightarrow \;\;
\gamma_j = \nu q^{2K - 2j +1} \; , \;\;\; (j > j') \; .
\ee
For the case $K=2n+1$ the possibility
$j=j' (= K+1 -j = n+1)$ is realized and Eq. (\ref{3.8.11})
(in view of (\ref{3.8.5}), (\ref{3.8.12}), (\ref{3.8.11c})) gives
\be
\lb{sol1a}
\gamma_{n+1}  = \nu q^{K-1} =1 \;\; \Rightarrow \;\; \nu = q^{1-K} \; .
\ee
For the case $K=2n$ we take  $j= {K \over 2}$,
$(j' = {K \over 2} +1 > j)$ in Eq. (\ref{3.8.11}) and obtain
$\gamma_{K \over 2} = \nu q^{K+1} - \lambda q$. On the other hand,
eq. (\ref{3.8.11c}) gives $\gamma_{{K \over 2} +1} = \nu q^{K-1}$. Thus,
in view of the condition $\gamma_{K \over 2} = \gamma_{{K \over 2} +1}^{-1}$
(\ref{3.8.12}), we deduce the equation $1 + \lambda q^K \nu = \nu^2 q^{2K}$
with two roots:
\be
\lb{sol1b}
\nu_1 = q^{1-K} \; , \;\;\; \nu_2 = - q^{-1-K} \; .
\ee
We summarize the results (\ref{3.8.11c}) -- (\ref{sol1b}),
for the solution of (\ref{3.8.11}), in the form
\be
\lb{3.8.13}
\gamma_{j} \equiv \frac{\tilde{c}_{j'}}{\tilde{c}_{j}} = \nu q^{2(N -j) + 1} \;
(j > j') \; , \;\;\;
\gamma_{j} = 1 \;
(j = j') \; , \;\;\; \nu = \epsilon q^{\epsilon - N} \; ,
\ee
(parameters $\tilde{c}_i$ were introduced in (\ref{Rgen1}))
and relate the cases $(\epsilon = +1)$ and
$(\epsilon = -1)$ to the groups $SO_{q}(N)$ and $Sp_{q}(N)$, respectively.

In order to determine from the conditions (\ref{3.8.13}) the parameters
$\tilde{c}_{j}$ and,
thus, to fix the matrices $\widetilde{C}$ (\ref{3.8.10b}),
we require fulfillment of the relation
$\widetilde{C}^{ij} = \epsilon \widetilde{C}_{ij}$
(cf. (\ref{3.7.12})).
Substitution of (\ref{3.8.10b}) gives the
equation $\tilde{c}_{j} \tilde{c}_{j'} = \epsilon$, which
together with (\ref{3.8.13}) enables us to choose $\tilde{c}_{j}$ in the form \cite{10}
\be
\lb{3.8.15}
\tilde{c}_{j} = \epsilon \, q^{j- \frac{1}{2} (N+\epsilon +1)} \;\;  (j>j') \Rightarrow
\tilde{c}_{j} = \epsilon_{j} \, q^{-\rho_{j}} \; ,
\ee
where
\be
\lb{eeee}
\begin{array}{c}
\epsilon_{i} = +1 \; \forall i \;\; ({\rm for} \;\; SO_{q}(N)) \; , \;\;
 \\ [0.1cm]
\epsilon_{i} = +1 \; (1 \leq i \leq n), \;\;
\epsilon_{i} = -1 \; (n+1 \leq i \leq 2n) \;\; ({\rm for} \;\; Sp_q(2n)),
\end{array}
\ee
and (cf. (\ref{rris}))
\be
\lb{rrrr}
 (\rho_{1}, \dots , \rho_{N}) =  \left\{
\begin{array}{l}
 (n-\frac{1}{2}, \, n-\frac{3}{2}, \dots , \frac{1}{2}, \, 0 , \,
-\frac{1}{2}, \, \dots , -n+\frac{1}{2} ) \; B: (SO_{q}(2n+1)), \\
 (n, \, n-1, \, \dots , 1, \, -1,
 \dots , \, 1-n , \, -n ) \;\;\;\;\;  C: (Sp_{q}(2n)), \\
 (n-1, \, n-2, \dots , 1, \, 0, \, 0, \, -1,
 \dots , 1-n ) \; D: (SO_{q}(2n)).
\end{array}
\right.
\ee
We note that diagonal matrices
 $\rho = {\rm diag}(\rho_1,...,\rho_N)$
are equal to $\rho = \sum_{i=1}^n \delta_i \, H_i$,
where elements $H_i = (e_{ii}-e_{N-i,N-i})$ form a
dual basis
of the Cartan subalgebras in the Lie algebras
$so(N)$, $sp(2n)$ and $\delta = (\delta_1, ...,\delta_n)$
are Weyl vectors for root systems of $so(2n),so(2n+1)$ and $sp(2n)$
(see definitions in \cite{Burb} and
in Section {\bf \ref{solYB}} below; see also \cite{IsRub2}, Sections 3.1.1 and 3.5.2).

Thus, the final expression
 for the $R$-matrix (\ref{Rgen11}) corresponding to the multi-parameter
deformation of the groups $SO(N)$ and $Sp(2n)$ (Refs. \cite{27aa}, \cite{27}) is:
$$
\R_{12}
= \sum_{i,j} \, q^{(\delta_{ij} - \delta_{ij'})} \,
\frac{f_{ij}}{f_{ji}} \, e_{ij} \otimes e_{ji} +
\lambda \sum_{i < j} \, e_{ii} \otimes e_{jj} - \lambda \,
 \sum_{i > j} \, \frac{f_{i'i}}{f_{jj'}} \, q^{\rho_i - \rho_j} \, \epsilon_i \, \epsilon_j \,
 e_{i'j} \otimes e_{ij'} \; ,
$$
where the parameters are defined in (\ref{ffff}), (\ref{eeee}), (\ref{rrrr}).
The matrix $R= P \R$ is represented in the component form as
\be
\lb{3.8.16}
R^{i_{1},i_{2}}_{j_{1},j_{2}}= \delta^{i_{1}}_{j_{1}} \delta^{i_{2}}_{j_{2}} \,
\left[ (q \, \delta^{i_{1}i_{2}}|_{i_{1} \neq i_{2}'} +
q^{-1} \, \kappa_{i_1} \, \delta^{i_{1}i_{2}'}|_{i_{1} \neq i_{2}}
+ a_{i_{2}i_{1}}|_{i_{1} \neq i_{2} \neq i'_{1}}
 + \delta^{i_{1}i_{1}'} \delta^{i_{2}i_{2}'} \right] +
\ee
$$
+ \lambda \, \delta^{i_{1}}_{j_{2}} \delta^{i_{2}}_{j_{1}} \,
\Theta_{i_{1}i_{2}} -\lambda \, \kappa_{i_1}
\delta^{i_{1}i'_{2}} \delta_{j_{1}j'_{2}} \,
\Theta^{i_{1}}_{j_{1}} \, \epsilon_{i_{1}} \, \epsilon_{j_{1}}
q^{\rho_{i_{1}} - \rho_{j_{1}}} \; ,
$$
where
\be
\lb{cond}
\begin{array}{c}
a_{ij} = 1/a_{ji} = \frac{f_{ij}}{f_{ji}} \;\; \forall j \neq i \neq j' \; , \;\;\;
a_{jj'}\kappa_j = a_{j'j}\kappa_j^{-1} = q^{-1} \;\; \forall j \neq j' \; , \\
a_{ij} \, a_{i'j} =  \kappa_{j} \; , \;\;\;
a_{ji} \, a_{ji'} =  \kappa^{-1}_{j} \; ,  \; \;
\kappa_{i} = (\kappa_{i'})^{-1} = \frac{f_{i'i}}{f_{ii'}} \; .
\end{array}
\ee

Now we clarify the role of the parameters $\kappa_i$.
Consider the $RTT$ algebra (\ref{3.1.1}) with multi-parameter
$R$-matrix (\ref{3.8.16}).
We show that for $\kappa_i \neq \pm 1$ the element $\tau$
introduced in (\ref{3.7.18}) is not central \cite{27aa}.
Indeed, we take the identity $\K_{12} \K_{23} \K_{12} = \K_{12}\, I_3$
(which is readily deduced from the explicit
representation (\ref{3.8.10a}), (\ref{3.8.10})) and multiply it by
$T_1 T_2 T_3$ from the right. For the right-hand side
of the identity we have
\be
\lb{cen1}
\K_{12} T_1 T_2 T_3 = \tau \, \K_{12} \, T_3 \; ,
\ee
while for the left-hand side we obtain
\be
\lb{cen2}
\K_{12} \K_{23} \K_{12}T_1 T_2 T_3 =
\K_{12} T_1 \K_{23} T_2 T_3 \, \K_{12} =
\K_{12} T_1 \K_{23} \, \K_{12} \tau =
 X_3 T_3 X_3^{-1} \K_{12} \tau \; ,
\ee
where
$X^i_j = C_{jk} C^{ki} = \delta^i_j \, \kappa_{i}$,
$(X^{-1})^i_j = C_{kj} C^{ik} = \delta^i_j \, \kappa_{i'}$
and we have used the identity
$T_1 \K_{23} \, \K_{12}  = \K_{23} \, \K_{12} X_3 T_3 X_3^{-1}$
followed from the definition (\ref{3.8.10}). Comparing (\ref{cen1})
and (\ref{cen2}) we obtain
 $$
 \tau \, T = X \, T \, X^{-1} \, \tau \;\;\; \Rightarrow \;\;\;
 \tau \, T^i_j = \kappa_{i} \, T^i_j \, \kappa_{j'} \; \tau \; .
 $$
Thus, only for $\kappa_i = \pm 1$ the element $\tau$ is central
and one can relate the multi-parameter
$R$-matrices (\ref{3.8.16}) with
quantum deformations $SO_{q,a_{ij}}(N)$ and $Sp_{q,a_{ij}}(2n)$
of the groups $SO(N)$ and $Sp(2n)$ (see discussions after
eq. (\ref{3.7.18})).

The conditions (\ref{cond})
show that, for $\kappa_i = \pm 1$,
the independent parameters are $q$ and $a_{ij}$ for
 $i < j \leq j'$.
The numbers of these parameters are $n(n-1)/2+1$ and $n(n+1)/2+1$
respectively for the groups of the series $C,D$ (with $N=2n$) and $B$ (with $N=2n+1$).
Note that the last term in the square brackets in the expression
(\ref{3.8.16}) is appeared only for the groups
of the series $B$. If we set $a_{ij}=1$
$(j' \neq i\neq j)$,
$\kappa_i = 1$ then the $R$-matrices (\ref{3.8.16}) are identical to the one-parameter matrices $R= P \R$ (\ref{frtosp})
deduced from (\ref{Rgen1}) and given in Ref. \cite{10}

\subsubsection{The case of $Osp_q(N|2m)$ supergroups\label{ospsup}}

For the groups $Osp(N|2m)$ $(K=N +2m)$, we choose a grading
in accordance with the
rules
\be
\lb{supar}
\begin{array}{l}
[j] = 0 \;\; {\rm for} \;\; m+1 \leq j \leq m+N \\
\, [j] = 1 \;\; {\rm for} \;\; 1 \leq j \leq m \; , \;\;\;
m + N+ 1 \leq j \leq N + 2m \; .
\end{array}
\ee
Thus, for $[j] =0$ $(j \neq j')$ and $[j] =1$, the
possibilities 1.) and 2.) in (\ref{3.8.4}) are respectively realized
\be
\lb{ospp}
a^{0}_j = (-1)^{[j]} \, q^{1 - 2[j]} = (-1)^{[j]} \, q^{(-1)^{[j]}}
\; , \;\;\; [j] = (j') \; .
\ee
In this case eq. (\ref{3.8.11}) is written as the system
of equations
\be
\lb{osp1}
\gamma_{j} \, q  +
\lambda \ \sum_{i=j+1}^{K}  \ \gamma_{i} =
- \nu \; , \;\; (N+ m+ 1 \leq j \leq N + 2m) \; ,
\ee
\be
\lb{osp2}
\gamma_{j} \, q - \lambda +
\lambda \ \sum_{i=j+1}^{K}  \ \gamma_{i} =
- \nu \; , \;\; (1 \leq j \leq m) \; ,
\ee
\be
\lb{osp3}
\gamma_{j} \, q^{(\delta_{j'j} \, -1)} + \lambda \ \Theta_{j'j} -
\lambda \ \sum_{i=j+1}^{K}  \ \gamma_{i} =
\nu \; , \;\; (m+1 \leq j \leq m+N) \; .
\ee
In (\ref{osp3}), for the case $j=j'$, we take into account (\ref{3.8.5}).
The solution of (\ref{osp1}) is (cf. (\ref{3.8.11c})):
\be
\lb{osp4}
\gamma_j = - \nu \, q^{2(j-K) - 1} \; , \;\; (N+ m+ 1 \leq j \leq N + 2m) \; ,
\ee
and we have $\lambda \sum\limits_{i > m+N} \gamma_i = \nu(q^{-2m}-1)$. Using this fact, the
eq. (\ref{osp3}) is written in the form:
\be
\lb{osp5}
\gamma_j q^{(\delta_{j'j} \, -1)}  + \lambda \ \Theta_{j'j} -
\lambda \sum_{i=j+1}^{m+N}  \ \gamma_{i} = \nu q^{-2m} \; ,
\ee
and its solution is:
\be
\lb{osp6}
\gamma_j = \nu \, q^{2(N -j) + 1} \; , \;\; (j' < j \leq N + m) \; .
\ee
In addition, for $N=2n+1$ and $j = j' = m+n+1$, we have
\be
\lb{osp7}
\gamma_{m+n+1} = \nu q^{N-2m -1} =1 \; , \;\;\; \nu = q^{2m -N +1}  \; ,
\ee
and for $N=2n$ we obtain condition:
$\gamma_{m+{N \over 2}} = \nu q^{N-2m +1} - \lambda q =
\gamma^{-1}_{m+{N \over 2} +1}$ which is equivalent to the quadratic equation on
$\nu$:
\be
\lb{osp8}
(\nu q^{N-2m} -q)(\nu q^{N-2m} +q^{-1}) = 0 \; .
\ee
Accordingly, we summarize the results
(\ref{osp4}), (\ref{osp6}) -- (\ref{osp8}) as
\be
\lb{3.8.14}
\gamma_{j} \equiv \frac{\tilde{c}_{j'}}{\tilde{c}_j}
= (-1)^{[j]} \nu q^{(-1)^{[j]} (2N-2j + 1) - [j]4m} \;\; (j >j') \; ,
\;\; \nu = \epsilon q^{\epsilon + 2m - N} \; ,
\ee
where $\epsilon = \pm 1$ and the case $\epsilon =+1$ corresponds $Osp_{q}(N|2m)$
while the case $\epsilon =-1$ we relate to a quantum group
denoted as $Osp'_{q}(2m|2n)$.
It is obvious that for the groups $Osp_{q}(2n+1|2m)$ (as well as
for $SO_{q}(2n+1)$)
we have $\gamma_{j} = \gamma_{j'} = 1$ for $j=j'$. Note that if in
(\ref{3.8.14}) we set $m=0$, or
$N=0$ and $q \rightarrow -q^{-1}$, then we reproduce (\ref{3.8.13}).

The analog of the relation (\ref{3.7.12}) for the groups $Osp_{q}(N|2m)$ is the equation
$$
\widetilde{C}^{ij} =
(-1)^{(i)} \, \epsilon \, \widetilde{C}_{ij} \; ,
$$
which is equivalent to
$(-1)^{(i)} \, \tilde{c}_{i} \, \tilde{c}_{i'} = \epsilon$,
and taking into account (\ref{3.8.14}) we obtain
$$
c_{j} = \epsilon \, q^{(-1)^{[j]}(j-m-N -\frac{1}{2}) -{\epsilon \over 2}
 + {N \over 2}} \;\;
(j>j') \;\; \Rightarrow \;\; \tilde{c}_j = \epsilon_j q^{-\rho_j} \; ,
$$
where $[j] = 0,1$ is the grading (\ref{supar}).
The parameters $(\rho_1, \dots , \rho_K)$,
$(\epsilon_1 , \dots , \epsilon_K)$ are fixed according to the
following cases:

\vspace{0.1cm}
\noindent
1.) The case $\epsilon =+1$, $\nu = q^{1 +2m-N}$ for $Osp_{q}(N|2m)$ $(N =2n+1)$:
{\footnotesize \be
\lb{er1}
\begin{array}{c}
\rho_i = ( \underbrace{ { N \over 2} -m, \dots , {N \over 2} -1}_m;
\underbrace{{N \over 2} -1, \dots ,{1 \over 2},0,
-{1 \over 2}, \dots , 1 - {N \over 2}}_{2n+1};
\underbrace{1 - {N \over 2}, \dots , m - {N \over 2}}_m ) \\ [0.4cm]
\epsilon_i = ( \underbrace{ -1, \dots ,  -1}_m;
\underbrace{+1, \dots , +1}_{2n+1};
\underbrace{ +1, \dots , +1}_m )
\end{array}
\ee}
2.) The case $\epsilon =+1$, $\nu = q^{1 +2m-N}$ for $Osp_{q}(N|2m)$ $(N =2n)$:
{\footnotesize \be
\lb{er2}
\begin{array}{c}
\rho_i = ( \underbrace{ n -m, \dots , n -1}_m;
\underbrace{n -1, \dots ,1,0,0,-1,
 \dots , 1 - n}_{2n};
\underbrace{1 - n, \dots , m - n}_m ) \\ [0.4cm]
\epsilon_i = ( \underbrace{ -1, \dots ,  -1}_m;
\underbrace{+1, \dots , +1}_{2n};
\underbrace{ +1, \dots , +1}_m )
\end{array}
\ee}
3.) The case $\epsilon =-1$, $\nu = -q^{-1 +2m-2n}$ for $Osp'_{q}(2m|2n)$:
{\footnotesize \be
\lb{er3}
\begin{array}{c}
\rho_i = ( \underbrace{n+1 -m, \dots , n}_m;
\underbrace{n, \dots ,1,-1,
 \dots ,-n}_{2n};
\underbrace{-n, \dots , m -1- n}_m ) \\ [0.4cm]
\epsilon_i = ( \underbrace{ -1, \dots ,  -1}_m;
\underbrace{+1, \dots , +1}_{n},\underbrace{-1, \dots , -1}_{n};
\underbrace{ -1, \dots , -1}_m )
\end{array}
\ee}

To conclude this subsection, we give the final expression
 for the $R$-matrix (\ref{Rgen11}) corresponding to the
quantum supergroups  $Osp_{q}(N|2m)$ $(\epsilon =+1)$ and
 $Osp'_{q}(2m|2n)$ $(\epsilon =-1)$
\be
\lb{Rosp}
\R_{12} = \sum_{i,j} \, (-1)^{(i)[j]} \, q^{(-)^{(i)}(\delta_{ij} - \delta_{ij'})} \,
 e_{ij} \otimes e_{ji}
+ \lambda \sum_{i < j} \, e_{ii} \otimes e_{jj} - \lambda \,
 \sum_{i > j} \,  q^{\rho_i - \rho_j} \, \epsilon_i \, \epsilon_j \,
 e_{i'j} \otimes e_{ij'} \; ,
\ee
where the parameters $\epsilon_i$, $\rho_j$ are defined in (\ref{er1}) - (\ref{er3}). We stress here that the matrix units $e_{ij}$
and tensor products in (\ref{Rosp}) {\it are not graded},
as follows form the discussion of the general Yang-Baxter
solution (\ref{Rgen11})
in Section {\bf \ref{mpsosp1}}.
To obtain (\ref{Rosp})
 we used the condition (\ref{ospp}) and put
$f_{ij}/f_{ji} = (-1)^{[i][j]}$ $(\forall i \neq j \neq i')$,
$f_{i'i} = f_{ii'}=1$  in (\ref{Rgen11}).
This choice of the parameters $f_{ij}$ is such that
$\R_{12}$ tends to the supertransposition matrix $(-1)^{(1)(2)}P_{12}$
 when $q \rightarrow 1$
 (for the notation $(-1)^{(1)(2)}$ see (\ref{matnot})).
 In the supergroup case, the multi-parameter $R$-matrices
are restored directly from (\ref{Rosp}) by the same
twisting (\ref{twistgl}) if we take into account conditions (\ref{ffff}).

The quantum supergroup $Osp_q(N,2m)$
is the graded algebra generated by elements
$T^i_{\; k}$ ($i,k=1,...,N+2m$)
of $(N+2m)\times (N+2m)$ supermatrix.
As in the case of quantum supergroups $GL_q(N,M)$
and $SL_q(N,M)$,
the generators $\{ T^i_{\; k} \}$ of the $Osp_q(N,2m)$
algebra satisfy the graded $RTT$
relations (\ref{3.6.4d}), (\ref{3.6.4}) but with the
$Osp_q(N,2m)$ type $R$-matrix (\ref{Rosp}).
The Hopf structure of the quantum supergroup
 $Osp_q(N,2m)$ is introduced in the same way as
 the Hopf structure of $GL_q(N,M)$
 (see Section {\bf \ref{qsuper}}).

We note that the parameters $\nu$ for the
cases $Osp_{q}$ (\ref{er2}) and
 $Osp'_{q}$ (\ref{er3}) are related
 to each other by means of the transformation:
 $q \leftrightarrow -q^{-1}$, $n \leftrightarrow m$. However, this
 transformation does not relate the corresponding
 $R$-matrices (\ref{Rosp}). Our conjecture is that for
 the cases $Osp_{q}$ (\ref{er2}) and $Osp'_{q}$ (\ref{er3})
the $R$-matrices (\ref{Rosp}) and corresponding quantum groups
are inequivalent.

\vspace{0.1cm}

The $R$-matrices constructed in this subsection for the quantum supergroups
realize $R$-matrix representations of the Birman - Murakami - Wenzl algebra
since they are the special examples of the
general $R$-matrix (\ref{Rgen3}) which
satisfy (\ref{char1}), (\ref{Bmw}). Some of these $R$-matrices
can be obtained on the basis of the results of Ref. \cite{35}, in which Baxterized
trigonometric solutions
(see next Sec. {\bf \ref{ospbax}}) of the Yang-Baxter equation associated with
the classical supergroups $Osp(N|2m)$ were obtained. Rational solutions,
some special cases, and other questions relating to the subject of the
quantum supergroups $Osp_{q}(N|2m)$ are also discussed in Refs. \cite{36}, \cite{38}
and \cite{39rag}.

\subsection{\bf \em $SO_{q}(N)$-, $Sp_{q}(2n)$- and $Osp_{q}(N|2m)$- invariant Baxterized $R$-matrices\label{ospbax}}
\setcounter{equation}0

Arguing, as in Sec. {\bf \ref{baxtel}},
  we conclude that the $SO_{q}(N)$ and $Sp_{q}(N)$
(as well as $Osp_{q}(N|2m)$) invariant
Baxterized matrices $\R(x)$ must be sought (by virtue of the fact
that the characteristic equation (\ref{3.7.1}) is cubic) in the
form of a linear combination of the three basis matrices {\bf 1}, $\R, \; \R^{2}$.
Expressing $\R^{2}$ in terms of $\K$ and $\R$, we can represent
invariant $R(x)$-matrix in the form \cite{Isaev1}
\be
\lb{3.9.1}
\R(x) = c(x) \left( \hbox{\bf 1} + a(x)\R + b(x) \K \right) \; ,
\ee
where $a(x), \; b(x)$ and  $c(x)$ are certain functions that
depend on the spectral parameter $x$. We determine the functions
$a(x), \; b(x)$ from the Yang-Baxter equation (\ref{3.5.2}).
The normalizing function $c(x)$ is not fixed by eq. (\ref{3.5.2}).
After substitution of (\ref{3.9.1}) in (\ref{3.5.2})
and using (\ref{3.7.19}) -- (\ref{3.7.26})
the following relations arise \cite{Isaev1}:
\be
\lb{3.9.5}
\begin{array}{l}
a_{1} + a_{3} + \lambda a_{1}\, a_{3} = a_{2} \; , \\ [0.3cm]
b_{3} - b_{2} - \lambda\nu a_{1}\, a_{3} + \nu a_{1}\, b_{3} -
\lambda  a_{1}\, b_{2}\, b_{3} + \lambda^{2} a_{1}\, a_{3}\, b_{2} \; +
 \\ [0.1cm]
+ \; b_{1}(1 + \nu a_{3} - \lambda a_{3}\, b_{2}
+ \mu b_{3} + \nu^{-1} a_{2}\, b_{3} + b_{2}\, b_{3})
 = 0 \; , \\  [0.3cm]
a_{2}\, b_{1} + a_{3}\, b_{1}\, b_{2} =
 a_{1}\, b_{2} + \lambda a_{1}\, a_{3}\, b_{2} \; , \;\;\;\;
a_{2}\, b_{3} + a_{1}\, b_{2}\, b_{3} = a_{3}\, b_{2} +
\lambda a_{1}\, a_{3}\, b_{2}   \; .
\end{array}
\ee
where we denoted
$$
a_{1} = a(x), \;\;\; a_{2} = a(xy), \;\;\; a_{3} = a(y), \;\;\;
b_{1} = b(x), \;\;\; b_{2} =b(xy), \;\;\; b_{3} = b(y) \; .
$$
The four relations (\ref{3.9.5}) are equivalent to the three functional equations
\be
\lb{3.9.6}
a(x) + a(y) + \lambda a(x)a(y) = a(xy) \; ,
\ee
\be
\lb{3.9.7}
\begin{array}{c}
b(y) - b(xy)
+ a(x) [ \nu b(y) - \lambda\nu a(y)
- \lambda  b(xy)b(y) + \lambda^{2} a(y)b(xy) ] + \\ [0.1cm]
+ b(x) [ 1 + \nu a(y) - \lambda a(y)b(xy) +
 \mu b(y) + \nu^{-1} a(xy)b(y) + b(xy)b(y) ] = 0 \; ,
\end{array}
\ee
\be
\lb{3.9.8}
a(xy)b(y) + a(x)b(xy)b(y) = b(xy)(a(y) + \lambda a(x)a(y))   \; ,
\ee
since the third and fourth relations in (\ref{3.9.5}) give the same equation
(\ref{3.9.8}). As it was to be expected, Eq. (\ref{3.9.6}) is identical to Eq.
(\ref{3.5.3}) obtained in the $GL_q(N)$ case, and its general solution
is given in (\ref{3.5.4}). By means of (\ref{3.9.6}), we can transform the
right-hand side of Eq. (\ref{3.9.8}) in such a way that (\ref{3.9.8})
reduces to the equation
\be
\lb{3.9.9}
\frac{a(x)}{a(xy)} = \frac{b(xy)-b(y)}{b(xy)(b(y) + 1)} \equiv
1 - \frac{b(y)(1 + b(y))^{-1}}{b(xy)(1 + b(xy))^{-1}} \; ,
\ee
We now note that Eq. (\ref{3.9.6}) can be rewritten in the form
\be
\lb{3.9.10}
\frac{a(x)}{a(xy)} =
1- \frac{a(y)(\lambda a(y) + 1)^{-1}}{a(xy)(\lambda a(xy) + 1)^{-1}}
\ee
and, comparing (\ref{3.9.9}) and (\ref{3.9.10}), we arrive at the result
\be
\lb{3.9.11}
\frac{a(y)(b(y)+1)}{(\lambda a(y) + 1)b(y)} = const
\equiv \frac{\alpha +1}{\lambda}  \; ,
\ee
where $\alpha$ denotes an arbitrary parameter. The specific choice of the
form of the constant in the right hand side of (\ref{3.9.11})
 is made for convenience in what follows.
 Substituting the solution
(\ref{3.5.4}) in (\ref{3.9.11}), we obtain the following general expression for $b(y)$:
\be
\lb{3.9.12}
b(y) =  \frac{y^{\xi} - 1}{\alpha y^{\xi} + 1} \; .
\ee
It is a remarkable fact that Eq. (\ref{3.9.7}) is satisfied identically on the functions
(\ref{3.5.4}) and (\ref{3.9.12}) if the constant
$\alpha$ satisfies the quadratic equation
\be
\lb{3.9.13}
\alpha^{2} - \frac{\lambda}{\nu}\alpha -
\frac{1}{\nu^{2}} = 0 \; .
\ee
The two solutions of this equation are readily found:
\be
\lb{3.9.14}
\alpha_{\pm} = \pm \frac{q^{\pm 1}}{\nu} \; .
\ee
where we recall that
\be
\lb{nu}
\begin{array}{ll}
\nu = \epsilon q^{\epsilon - N} \, , & {\rm for} \; {\rm groups} \;\;
SO_q(N) \;\; (\epsilon = +1), \;\; Sp_q(N) \;\; (\epsilon = -1); \\[0.2cm]
\nu = \epsilon q^{\epsilon + 2m - N} \, , & {\rm for} \; {\rm supergroups} \;\;
Osp_q(N|2m) \;\; (\epsilon = +1) \; {\rm and} \\[0.2cm]
 & \;\;\;\;\;\;\;\;\;\;\;\;\;\;\;\;\;\;\; Osp'_q(2m|N) \;\; (N=2n, \; \epsilon = -1).
\end{array}
\ee

Thus, the
 solutions of the Yang-Baxter equation (\ref{3.5.2}) can be represented in the form\footnote{The Baxterized trigonometric $R$-matrices
(\ref{3.9.15}), corresponding to the one of the parameter choice in
(\ref{3.9.14}), were first found by V.Bazhanov in 1984 and were published in
\cite{40}. The same R-matrices were independently constructed in Ref. \cite{41}.} \cite{CGX}, \cite{Isaev1}
\be
\lb{3.9.15}
\R(x) = c(x) \left( \hbox{\bf 1} +
\frac{1}{\lambda}(x^{\xi}-1) \R +
 \frac{x^{\xi} - 1}{\alpha x^{\xi} + 1} \K \right) \; .
\ee
and we have the two possibilities $\alpha = \alpha_{\pm}$ (\ref{3.9.14}),
 which are inequivalent (both for all the cases
 $SO_q(N)$, $Sp_g(N)$ and $Osp_q(N|2m)$),
since these solutions cannot be reduced to each other by any functional
transformations of the spectral parameter $x$.
 However, these solutions
are related by the transformation $q \to -q^{-1}$.
 For convenience we
choose $c(x) = x$ and $\xi = -2$ in (\ref{3.9.15}); then for the R-matrices
(\ref{3.9.15}) we can propose four equivalent expressions:
\be
\lb{3.9.15a}
\R^{\pm}(x) := \frac{1}{\lambda} \left(
 x^{-1} \R - x \R^{-1} \right) +
 \frac{\alpha_{\pm} + 1}{\alpha_{\pm} x^{-1} + x} \K  =
\ee
\be
\lb{3.9.15d}
 = \frac{(q^{\pm 2} x^{-1} - x)}{(q^{\pm 2} -1)\, x^2} \, \frac{
(\R   \pm q^{\pm 1} x^{2})}{(\R \pm q^{\pm 1} x^{-2})} =
 \ee
$$
= \frac{x-x^{-1}}{\lambda(x+ \alpha_{\pm} x^{-1})} \left( - x \R^{-1}
 - \alpha_{\pm} x^{-1} \R +
\frac{\lambda (\alpha_{\pm} +1)}{x - x^{-1}}  \right) =
$$
 \be
\lb{3.9.16}
= \frac{(x^{-1}q - x q^{-1})}{\lambda}\P^{+} +
\frac{(xq -(xq)^{-1})}{\lambda} \P^{-} +
 \frac{(q^{\pm 2} x^{-1} - x)}{(q^{\pm 2} -1)} \,
\frac{(x^{-1} + x \, \alpha_{\pm} )}{(x + x^{-1} \, \alpha_{\pm})} \P^{0}\; ,
\ee
where projectors $\P^{\pm}$ and $\P^{0}$ are defined in (\ref{3.7.2}).
The last expression is the spectral decomposition of $\R(x)$,
from which, for example, we can readily obtain identities
\be
\lb{3.9.17}
\R^{+}(\pm q) = \pm (q + q^{-1}) \, \P^{-} \; , \;\;\;\;\;
\R^{-}(\pm q^{-1}) = \pm (q + q^{-1}) \, \P^{+} \; ,
\ee
\be
\lb{3.9.17b}
\lim_{x^2 \to -\alpha_{\pm}} R^{\pm}(x) \sim  \P^{0} \; , \;\;\;\;\;
\R^{\pm}(1) = \hbox{\bf 1} \; , \;\;\;\;\;
\R^{\pm}(i) = \pm \frac{i (q + q^{-1})}{\lambda}
(\hbox{\bf 1} -2\P^{\pm})  \; .
 \ee
 From rational representations (\ref{3.9.15d}) of $R$-matrix,
 one can immediately deduce the identity
\be
\lb{3.9.17a}
\R^{\pm}(x) \, \R^{\pm}(x^{-1}) = \left(1 - \frac{(x-x^{-1})^{2}}{\lambda^{2}}
\right)\cdot  {\bf 1} \; .
\ee
Note that the relations (\ref{3.9.17}) -- (\ref{3.9.17a}) agree with the Yang-Baxter equation (\ref{3.5.2}).

The cross-unitarity condition
for the BMW type $R$-matrix (\ref{3.9.15a}) can be written in the
matrix form as (cf. (\ref{cross01}))
\be
\lb{cross02}
\begin{array}{c}
Tr_{D(2)}\left( \R^{\pm}_{1}(x) P_{01} \R^{\pm}_1(z) \right) =
\eta^{\pm}(x) \, \eta^{\pm}(z) \, D_0 \, I_1 \; , \\ \\
Tr_{Q(1)}\left( \R^{\pm}_{1}(x) P_{23} \R^{\pm}_1(z) \right) =
\eta^{\pm}(x) \, \eta^{\pm}(z) \, Q_{3} \, I_2 \; ,
\end{array}
\ee
where the matrices $D,Q$ are defined
in (\ref{dmatr}) and
$$
(x\, z)^2 = \alpha_{\pm}^2  \; , \;\;\;\;\;
\eta^{\pm}(x)=  \frac{1}{\lambda} (x-x^{-1})
\frac{(\alpha_{\pm} \nu x^2 + \nu^{-1})}{(x^2 +\alpha_{\pm})}
\; , \;\;\;\;\; \alpha_{\pm} := \pm \frac{q^{\pm 1}}{\nu} \; .
$$

The Baxterized $\R$ matrices (\ref{3.9.15}) and
(\ref{3.9.15a}) -- (\ref{3.9.16}) determine
algebras with the defining relations (\ref{3.5.6}). However, a realization of
the operators $L(x)$ in terms of the generators $L^{(\pm)}$
of the quantum algebras $U_{q}(so(N))$ and $U_{q}(sp(N))$
(analogous to (\ref{3.5.7})) is unfortunately missing
(see, however, \cite{18mol}, \cite{IsOgM1} and
discussion of the case $q \to 1$ in \cite{28}, \cite{IsKK}).

To conclude this subsection, we present the expressions for the
rational $R$-matrices of
the Yangians $Y(so(N))$, $Y(sp(N))$ and $Y(osp(N|2m))$. We give
the definition of Yangians in Section {\bf \ref{solYB}} below.
We make the ansatz $x= \exp(-\lambda \theta/ 2)$
 for the spectral parameter
 in (\ref{3.9.15a}) and rewrite the $R$-matrix
in the form (cf. (\ref{3.5.10a})):
\be
\lb{3.9.15b}
\begin{array}{c}
\R(\theta) := \R \left( e^{-{\lambda \over 2}\theta} \right)
= \cosh (\lambda \theta / 2)\, [ {\bf 1} - \K ] + {1 \over \lambda}
\sinh (\lambda \theta / 2) \, [ \R + \R^{-1} ] + \\ \\
+  [\cosh (\lambda \theta / 2) + \beta_{\pm} \sinh (\lambda \theta / 2)]^{-1} \, \K
\end{array}
\ee
where $\displaystyle \beta_{\pm} = \frac{\alpha_{\pm} -1}{\alpha_{\pm} +1}$.
The Yangian $R$-matrices can be obtained from
(\ref{3.9.15b}) after the passage to the limit
$h \rightarrow 0$ $(q=exp(h) \rightarrow 1)$. Further, it is easy
to see that the cases $\alpha = \alpha_{+}$, $\epsilon = 1$ $(SO_{q}(N))$ and
$\alpha = \alpha_{-}$, $\epsilon = -1$ $(Sp_{q}(2n))$
is reduced to the $GL(N)$-symmetric
Yang's $R$-matrix (\ref{3.5.11}). The nontrivial $SO(N)$- and $Sp(N)$-symmetric
Yangian $R$-matrices for $Y(so(N))$ and $Y(sp(N))$ correspond to
the choice
\be
\lb{3.9.18}
\alpha = \alpha_{-}, \;\;\;
\epsilon = 1 \;\;\; (SO_{q}(N)) \; ; \;\;\;\;\;\;
\alpha = \alpha_{+}, \;\;\;
\epsilon = -1 \;\;\; (Sp_{q}(N)) \; ,
\ee
and have the form
\be
\lb{3.9.19}
\R(\theta) =  \left( \hbox{\bf 1} +
\theta  \, P_{12} \right) +
\frac{2 \, \theta}{\left( 2 \epsilon  - (N + 2 \theta) \right)} \, K^{(0)}_{12} \; .
\ee
The matrix $K^{(0)}_{12}$ is defined in (\ref{3.7.6a}).
Nontrivial rational $R$-matrices for super Yangians $Y(osp(N|2m))$ and $Y'(osp(2m|2n))$
can be obtained from (\ref{3.9.15b}) in the cases
$$
\alpha = \alpha_{-}, \;\;\;
\epsilon = 1 \;\;\; (Osp_{q}(N|2m)) \; ; \;\;\;\;\;\;
\alpha = \alpha_{+}, \;\;\;
\epsilon = -1 \; , \;\; N= 2n \;\;\; (Osp'_{q}(2m|2n)) \; .
$$
The form of these supersymmetric $R$-matrices is
\be
\lb{3.9.19s}
\R(\theta) = (\hbox{\bf 1} + \theta \, {\cal P}_{12} ) +
\frac{2 \, \theta}{2 \epsilon + 2 m - (N + 2 \, \theta )} \, {\cal K}^{(0)}_{12}
\ee
where ${\cal P}^{i_1 i_2}_{j_1 j_2} =
(-1)^{[i_1][i_2]} \delta^{i_1}_{j_2} \delta^{i_2}_{j_1}$ is the
supertransposition operator
(the parity $[j]$ is defined in (\ref{supar})).
The matrix $({\cal K}^{(0)})^{i_1 i_2}_{j_1 j_2} =
{\cal C}^{i_1 i_2} {\cal C}_{j_1 j_2}$ is a classical limit ($q \rightarrow 1$)
of the rank 1 matrix $\K$ in the supersymmetric case and
the ortho-symplectic matrices
${\cal C}^{i j} = \epsilon_j \delta^{i j'}$,
${\cal C}_{ij} = \epsilon_i \delta_{i j'}$
are determined by their parameters $\epsilon_i$ (\ref{er1}) - (\ref{er3}).
Then, the defining relations for the generators (\ref{3.5.12}) of the
Yangians $Y(so(N))$, $Y(sp(N))$ and
$Y(osp(N|2m))$, $Y'(osp(2m|2n))$ are
identical to (\ref{3.5.10}) and (\ref{sgly1}), respectively,
while the comultiplication is given by (\ref{3.5.13}).

The Yangian $R$-matrix (\ref{3.9.19}) for the $SO(N)$ case was found in Ref. \cite{3},
and that for the $Sp(2n)$ case in Ref. \cite{39}. These $R$-matrices were used in
Ref. \cite{28} to construct and investigate exactly solvable
$SO(N)$- and $Sp(2n)$-symmetric magnets.
Twisted Yangians for the $SO(N)$ and $Sp(2n)$ cases have been
considered in \cite{Mol}.
The super Yangians
of the type $Y(osp(N|2m))$ and corresponding
spin chain models where discussed in \cite{39rag}.

\subsection[\bf \em Split Casimir operators and
rational solutions of Yang-Baxter equations. 
Yangians]{\bf \em Split Casimir operators and
rational solutions of Yang-Baxter equations for
simple Lie algebras. Yangians $Y(\mathfrak{g})$\label{solYB}}
\setcounter{equation}0

The material of this subsection is based on the papers
\cite{Drin85}, \cite{MaKa}; see also \cite{MaKa2}, \cite{ChaPr},
\cite{IsKriv}.

\subsubsection{Invariant
  $R$-matrices for simple Lie algebras $\mathfrak{g}$}

Let $\mathfrak{g}$ be a simple Lie algebra with  the basis
elements $X_a|_{a=1,...,\dim(\mathfrak{g})}$
 and defining relations
\be\lb{lialg}
[X_a, \; X_b] = X^d_{ab} \; X_d \; ,
\ee
 where $X^d_{ab}$ are the structure constants.
 We denote an enveloping algebra of the Lie algebra  $\mathfrak{g}$ as ${\cal U}(\mathfrak{g})$.
 Let ${\sf g}^{df}$ be the inverse matrix to the
 Cartan-Killing  metric
\be\lb{li04}
   {\sf g}_{ab} \equiv X^{d}_{ac} \, X^{c}_{bd} =
   {\rm Tr}({\rm ad}(X_a)\cdot {\rm ad}(X_b)) \; ,
\ee
where ${\rm ad}$ denotes adjoint representation. Introduce operator
 \be
 \lb{kaz-01}
\hat{C}  = {\sf g}^{ab} X_a \, \otimes \,
  X_b \equiv X_a \, \otimes \,  X^a
  \;\; \in \;\; \mathfrak{g} \, \otimes \,  \mathfrak{g}
   \;\; \subset \;\; {\cal U}(\mathfrak{g})\, \otimes \, {\cal U}(\mathfrak{g}) \; ,
 \ee
 which is called the {\it split (or polarized) Casimir operator} of the Lie algebra $\mathfrak{g}$.  The operator (\ref{kaz-01})
is related to the usual quadratic Casimir operator
 \be
 \lb{kaz-c2}
 C_{(2)} = {\sf g}^{ab} \; X_a \cdot X_b \;\; \in \;\;
 {\cal U}(\mathfrak{g}) \; ,
  \ee
 by means of the formula
 \be
 \lb{adCC1}
 \Delta(C_{(2)}) = C_{(2)} \otimes I +
 I \otimes C_{(2)} + 2 \, \hat{C}  \; ,
 \ee
 where $\Delta$: ${\cal U}(\mathfrak{g}) \to
 {\cal U}(\mathfrak{g}) \otimes {\cal U}(\mathfrak{g})$
 is the standard comultiplication
 defined by $\Delta(X_a) = X_a \otimes I + I\otimes X_a $.
 Let $u$ be a spectral
 parameter. One can check that the
 operator function
 \be
 \lb{rmatr}
 r(u) = \frac{\hat{C}}{u} = \frac{X_a \, \otimes \,  X^a}{u}
 \equiv r_{21}(u) \; ,
 \ee
obeys the semiclassical
Yang-Baxter equation (cf. (\ref{3.2.2})):
 \be
 \lb{cYBE}
 [r_{12}(u), \, r_{13}(u+v)] + [r_{13}(u+v), \, r_{23}(v)] +
 [r_{12}(u), \, r_{23}(v)] = 0 \; .
 \ee

 The aim of this subsection is to find rational (as a function
 in the spectral parameter $u$) solutions $R(u)$ of the
 Yang-Baxter equations (\ref{3.5.9a})
 ($\theta'=u$, $\theta=u+v$):
 \be
 \lb{YBEu}
 R_{12}(u) \, R_{13}(u+v) \, R_{23}(v)
 =  R_{23}(v) \, R_{13}(u+v) \, R_{12}(u) \; ,
 \ee
  that are unitary
 \be
 \lb{uniR}
 R_{12}(u) R_{21}(-u)={\bf 1} \equiv I \otimes I  \; ,
 \ee
  and possess semiclassical
 behavior $R_{12}(u) \to {\bf 1}$ as $u \to \infty$.
 We then write the expansion
 \be
 \lb{Rmatr}
 R_{12}(u) = {\bf 1} + \frac{\hat{C}}{u} + \frac{X}{u^2}
 + O\Bigl( \frac{1}{u^3} \Bigr) \; .
 \ee
 The second term here is justified by (\ref{cYBE}).
 As we will see below
 the solutions of this kind are given (up to a renormalization)
 by (\ref{3.9.19})
 for Lie algebras
 $\mathfrak{g} = \mathfrak{so}(N)$ and
 $\mathfrak{sp}(N)|_{N=2r}$ in defining representations.

 First, we use the unitarity condition (\ref{uniR}) to find $X$:
 $$
 \begin{array}{c}
 \displaystyle
 {\bf 1} = R_{12}(u) R_{21}(-u) =
 {\bf 1} - \frac{1}{u^2} \hat{C}^2 +
 \frac{1}{u^2}(X_{12}+X_{21}) + ... \; .
 \end{array}
 $$
 We search for the symmetric solutions
 $X_{12}= X_{21}$, so we have
 \be
 \lb{XX}
 X_{12} = \frac{1}{2} \hat{C}^2 \; .
 \ee
 Then we examine the limit $v \to \infty$ of (\ref{YBEu})
 \be
 \lb{3.5.9f}
 \begin{array}{c}
 \displaystyle
 R_{12}(u) \, \Bigl({\bf 1} + \frac{\hat{C}_{13}}{u+v}
 + \frac{X_{13}}{(u+v)^2} +...\Bigr)
  \,  \Bigl({\bf 1} + \frac{\hat{C}_{23}}{v}
 + \frac{X_{23}}{v^2} +...\Bigr) = \\ [0.3cm]
  \displaystyle
 =  \Bigl({\bf 1} + \frac{\hat{C}_{23}}{v}
 + \frac{X_{23}}{v^2} +...\Bigr) \, \Bigl({\bf 1} + \frac{\hat{C}_{13}}{u+v}
 + \frac{X_{13}}{(u+v)^2} +...\Bigr) \, R_{12}(u) \; .
 \end{array}
\ee
 We expand out the brackets in (\ref{3.5.9f})
 and multiply both sides by $v(u+v)$. As a result we obtain
 \be
 \lb{Rmat01}
 \begin{array}{c}
 R_{12}(u) \, \Bigl( v(\hat{C}_{23} + \hat{C}_{13})
 + u \hat{C}_{23} + \hat{C}_{13} \hat{C}_{23} +
 \frac{u+v}{v} X_{23} + \frac{v}{u+v} X_{13} + ...  \Bigr)
 = \\ [0.3cm]
 = \Bigl( v(\hat{C}_{23} + \hat{C}_{13})
 + u \hat{C}_{23} + \hat{C}_{23} \hat{C}_{13}  +
 \frac{u+v}{v} X_{23} + \frac{v}{u+v} X_{13} + ...  \Bigr)
 \, R_{12}(u) \; .
 \end{array}
 \ee
 We see that the terms of order $v$ give
 \be
 \lb{Rmat04}
 [R_{12}(u), \, (\hat{C}_{23} + \hat{C}_{13})] = 0
 \;\;\;\;\; \Rightarrow \;\;\;\;\;
 [R_{12}(u), \, I \otimes X_a + X_a\otimes I] = 0  \; ,
 \ee
 which is the condition of the invariance of $R(u)$
 under the action of $\mathfrak{g}$. Thus, by Schur's lemma, one can express the image $R_{(\mu\nu)}(u) = (T_\mu \otimes T_\nu)R(u)$
 of the operator $R(u)$ in the representation
 $(T_\mu \otimes T_\nu)$ of $\mathfrak{g}$ as following
 \be
 \lb{Rmat02}
 R_{(\mu\nu)}(u) =  \sum_{T_\lambda \subset T_\mu \otimes T_\nu}
 \tau_\lambda(u) \; P_{\lambda} \; ,
 \ee
 where $P_{\lambda}$ is the projector onto the irreducible
 subrepresentation  $T_\lambda \subset T_\mu \otimes T_\nu$ and
 $\tau_\lambda(u)$ are some rational functions of $u$, which
 is yet undetermined. At this stage, we require that the set
 of projectors $P_{\lambda}$ form the complete system of
 mutually orthogonal projectors
 \be
 \lb{Rmat09}
 \sum_\lambda P_{\lambda} = I_\mu \otimes I_\nu \;\; ,
 \;\;\;\;\;\;
 P_{\lambda} \, P_{\lambda'} =
 P_{\lambda} \delta_{\lambda \lambda'} \; .
 \ee
 We also require that the
 decomposition $T_\mu \otimes T_\nu =
 \sum_\lambda T_\lambda$ is without of
 multiplicities, otherwise, $R(u)$ acts on the isomorphic
 components $T_{\lambda_1},..., T_{\lambda_r}$ as matrix
 $||M_{ij}(u)||_{i,j=1,...,r}$ which is not, in general, be diagonalizable.

 The terms of order $v^0=1$ in (\ref{Rmat01}) give
 \be
 \lb{Rmat03}
 \begin{array}{c}
 R_{12}(u) \, \Bigl(  u \hat{C}_{23} + \hat{C}_{13} \hat{C}_{23} +
  X_{23} +  X_{13} \Bigr)
 = \Bigl( u \hat{C}_{23} + \hat{C}_{23} \hat{C}_{13}  +
  X_{23} +  X_{13} \Bigr)
 \, R_{12}(u) \; .
 \end{array}
 \ee
 We rewrite it by using (\ref{XX}) and applying identities
 $$
 \begin{array}{c}
 \hat{C}_{13}\, \hat{C}_{23} +
 \frac{1}{2}(\hat{C}_{13}^2 +\hat{C}_{23}^2) =
 \frac{1}{2}[\hat{C}_{13}, \; \hat{C}_{23}] +
 \frac{1}{2}(\hat{C}_{13} +\hat{C}_{23})^2  \; , \\ [0.2cm]
 [\hat{C}_{13}, \; \hat{C}_{23}] =
 X^a_{bc} \; X^b \otimes X^c \otimes X_a
 = - \frac{1}{2}[\Delta{C_{(2)}}, \; (I \otimes X_a)]
 \otimes X^a  \; ,
 \end{array}
 $$
 so that, because of  (\ref{Rmat04})
 we simplify (\ref{Rmat03}) as
 $$
 R_{12}(u) \, \Bigl(  u \hat{C}_{23} +
 \frac{1}{2}[\hat{C}_{13}, \; \hat{C}_{23}] \Bigr)
 = \Bigl( u \hat{C}_{23} +
 \frac{1}{2}[\hat{C}_{23}, \; \hat{C}_{13}] \Bigr)
 \, R_{12}(u) \;\; \Rightarrow
 $$
  \be
 \lb{Rmat05}
  \begin{array}{c}
 R_{12}(u) \, \Bigl(  u (I \otimes X_a) -
 \frac{1}{4}[\Delta{C_{(2)}}, \, I \otimes X_a] \Bigr)
 = \Bigl( u (I \otimes X_a)  +
 \frac{1}{4}[\Delta{C_{(2)}}, \, I \otimes X_a] \Bigr)
 \, R_{12}(u) \; .
 \end{array}
 \ee
 Now we consider the image of (\ref{Rmat05}) in the representation
  $T_\mu \otimes T_\nu$, substitute (\ref{Rmat02}) and
  act by projectors $P_\kappa$ and $P_\lambda$ from the
  right and left, respectively. As a result we deduce
  relation between coefficients $\tau_\lambda(u)$:
   \be
 \lb{Rmat06}
  \begin{array}{c}
\tau_\lambda(u) \,  \Bigl(  u  -
 \frac{1}{4} \bigl( C_{(2)}(\lambda) - C_{(2)}(\kappa) \bigr) \Bigr)
 P_\lambda \, (I \otimes X_a) P_\kappa = \\ [0.2cm]
 = \tau_\kappa(u) \, \Bigl( u  +
 \frac{1}{4}\bigl( C_{(2)}(\lambda) - C_{(2)}(\kappa) \bigr) \Bigr)
 \, P_\lambda \, (I \otimes X_a) P_\kappa \; ,
 \end{array}
 \ee
 where $C_{(2)}(\lambda)$ is the value of the
 quadratic Casimir operator (\ref{kaz-c2}) in the representation
 $T_\lambda$ and for brevity we write $(I \otimes X_a)$
 instead of $(T_\mu(I) \otimes T_\nu(X_a))$.
 We enumerate the representations $T_\lambda$
 of the Lie algebra $\mathfrak{g}$
 by their highest weights $\lambda$. In this case
 the value of $C_{(2)}(\lambda)$ is given by the formula
 \be
 \lb{kvkaz}
 C_{(2)}(\lambda) = (\lambda,\lambda + 2 \, \delta) \; ,
 \ee
 where $\delta$ is the Weil vector of the algebra $\mathfrak{g}$.
 Finally, from the equation (\ref{Rmat06}),
 in the case when
 $P_\lambda \, (I \otimes X_a) P_\kappa \neq 0$, we have
 \be
 \lb{Rmat07}
 \frac{\tau_\lambda(u)}{\tau_\kappa(u)} =
 \frac{ u  + \frac{1}{4} \bigl( C_{(2)}(\lambda)
 - C_{(2)}(\kappa) \bigr) )}{ u  -
 \frac{1}{4} \bigl( C_{(2)}(\lambda) - C_{(2)}(\kappa)\bigr)} \; .
 \ee

  We consider the condition
  $P_\lambda \, (I \otimes X_a) P_\kappa \neq 0$
  in more detail. Let $V_\lambda$
  be the space of the representation $T_\lambda$. We note that
  $(I \otimes X_a + X_a \otimes I) \cdot V_\lambda \subset V_\lambda$,
  where $V_\lambda \subset V_\mu \otimes V_\nu$,
  and, for orthogonal projectors $P_\kappa$ and $P_\lambda$,
  we deduce
  \be
 \lb{Rmat08}
  P_\lambda \, (I \otimes X_a) P_\kappa =
  \frac{1}{2}P_\lambda \, (I \otimes X_a - X_a \otimes I) P_\kappa \; .
  \ee
  One can interpret $(I \otimes X_a - X_a \otimes I)$
  as the tensor operator in the adjoint representation
  and, in according to the Wigner-Eckart theorem,  the matrix
  (\ref{Rmat08}) should be
  proportional to Clebsch-Gordan coefficients which
  transform the basis of $V_\lambda$ into the basis of
  $V_{\rm ad} \otimes V_\kappa$.
We note that
for existence of the $R$-matrix, it is necessary that the
system of equations (\ref{Rmat07}) have a solution. However,
in general, the system (\ref{Rmat07}) is overdetermined
and not always has a solution.

Further we consider the
 equivalent representations $T_\nu = T_\mu$ and require that
 the $R$-matrix is symmetric $R_{12} = R_{21}$. Then, the space
 $V_\mu \otimes V_\mu$ is splitted into symmetric
 $P_{12}^{(+)} \; (V_\mu \otimes V_\mu)$
 and antisymmetric $P_{12}^{(-)} \;
 (V_\mu \otimes V_\mu)$ parts, where
 $P_{12}^{(\pm)} := \frac{1}{2}(I \pm P_{12}) $.
 It means that the whole set
 of projectors (\ref{Rmat09}) also is divided onto subsets
 of symmetric and antisymmetric projectors
 \be
 \lb{Rmat10}
 P_\kappa^{(+)} : = P_{12}^{(+)} P_\kappa \; , \;\;\;
   P_\sigma^{(-)} : = P_{12}^{(-)} P_\sigma
   \;\;\;\;\; \Rightarrow \;\;\;\;\;
 \sum_\kappa P_\kappa^{(+)}  = P^{(+)}  \; , \;\;\;\;\;
 \sum_\sigma P_\sigma^{(-)}  = P^{(-)}  \; ,
 \ee
 and for matrices (\ref{Rmat08}) we have
 $P_\sigma^{(\pm)} \, (I \otimes X_a -
 X_a \otimes I) P_\kappa^{(\pm)}=0$.
 So, the nonzero contributions to (\ref{Rmat08}) are
 $P_\sigma^{(\pm)} \, (I \otimes X_a -
 X_a \otimes I) P_\kappa^{(\mp)}$.
  Thus, the representations
 $T_\lambda$ and $T_\kappa$ in equations
 (\ref{Rmat07}) should have the different symmetry,
 i.e.
 $$
 V_\lambda \subset P^{(+)} \, V_\mu^{\otimes 2} \; , \;\;\;\;
 V_\kappa\subset P^{(-)} \, V_\mu^{\otimes 2} \; , \;\;\;
 {\rm or} \;\;\;
 V_\lambda \subset P^{(-)} \, V_\mu^{\otimes 2} \; , \;\;\;\;
 V_\kappa\subset P^{(+)} \, V_\mu^{\otimes 2} \; ,
 $$
 and satisfy conditions
 \be
 \lb{Rmat11}
 T_\lambda \subset {\rm ad} \otimes T_\kappa \; , \;\;\;\;
  T_\kappa \subset {\rm ad} \otimes T_\lambda \; .
 \ee

 \vspace{0.2cm}

 Below, by solving of the system of
  equations (\ref{Rmat07}),
  we present $\mathfrak{g}$-invariant
  $R$-matrices  for all simple Lie algebras $\mathfrak{g}$
  (except for $\mathfrak{g} = \mathfrak{sl}_N,\mathfrak{e}_8$)
  in the defining representation $T_\mu =\square \equiv [1]$
  \cite{OgWieg} (see also \cite{IsKriv}).

  {\bf 1. The $\mathfrak{so}$ and $\mathfrak{sp}$
  algebras} \\
  The expansion of the
  tensor product of defining representations of
  algebras $\mathfrak{so}(N)$ and $\mathfrak{sp}(N)|_{N=2r}$
  is
   $[1]^{\otimes 2} = [2] + [1^2] + [\emptyset]$,
and the highest weights of subrepresentations,
Weyl vectors $\delta$
 and corresponding eigenvalues $C_{(2)}(\lambda)$,
defined in (\ref{kvkaz}), are
$$
\lambda_{[2]} = (2,0,...,0) \; , \;\;\;\;
\lambda_{[1^2]} = (1,1,0,...,0) \; , \;\;\;\;
\lambda_{[\emptyset]} = (0,...,0) \; ,
$$
for $\mathfrak{so}(N)$:
 \be
 \lb{Rmat12}
 \delta = \sum_{i=1}^{[\frac{N}{2}]} (\frac{N}{2}-i)e^{(i)}
 \; , \;\;\;\;
C_{(2)}(\lambda_{[2]}) = 2 N\; , \;\;\;\;
C_{(2)}(\lambda_{[1^2]}) = 2 N -4 \; , \;\;\;\;
 C_{(2)}(\lambda_{[\emptyset]}) = 0  \; ,
 \ee
 and for $\mathfrak{sp}(N)$:
  \be
 \lb{Rmat12b}
 \delta = \sum_{i=1}^{\frac{N}{2}} (\frac{N}{2}+1-i)e^{(i)}
 \; , \;\;\;\;
C_{(2)}(\lambda_{[2]}) = 2 N+4\; , \;\;\;\;
C_{(2)}(\lambda_{[1^2]}) = 2 N \; , \;\;\;\;
 C_{(2)}(\lambda_{[\emptyset]}) = 0 \; .
\ee
Below we also need the expansion
\be
\lb{Rmat11b}
{\sf [1^2]} \otimes {\sf [2]} = {\sf [3,1]} +
{\sf [2,1^2]} + {\sf [1^2]} + {\sf [2]} \; .
\ee

For the $\mathfrak{so}(N)$ case,
the representations ${\sf [2]}$  and  $[\emptyset]$
belong to the symmetric part of ${\sf [1]}^{\otimes 2}$,
while the adjoint representation
${\sf [1^2]}={\rm ad}$ belongs to the antisymmetric part of
${\sf [1]}^{\otimes 2}$, and we have (cf. (\ref{Rmat11}), (\ref{Rmat11b}))
${\sf [1^2]} \subset {\sf [1^2]} \otimes [\emptyset]$,
${\sf [1^2]} \subset {\sf [1^2]} \otimes {\sf [2]}$. Therefore, in view
 of (\ref{Rmat12}), the system
 of equations (\ref{Rmat07}) is written as
 \be
 \lb{Rmat13}
 \frac{\tau_{[\emptyset]}(u)}{\tau_{[1^2]}(u)} =
 \frac{u  - \frac{N}{2}+1}{ u  +\frac{N}{2}-1}
 \; , \;\;\;\;\quad
  \frac{\tau_{[2]}(u)}{\tau_{[1^2]}(u)} =
 \frac{ u +1}{ u - 1} \; ,
 \ee
 and after a renormalization the $\mathfrak{so}$-invariant
  $R$-matrix (\ref{Rmat02}) is
  \be
 \lb{Rmat14}
 R(u) = P_{[1^2]} +
 \frac{\tau_{[2]}(u)}{\tau_{[1^2]}(u)} P_{[2]}+
 \frac{\tau_{[\emptyset]}(u)}{\tau_{[1^2]}(u)} P_{[\emptyset]}
 =  P_{[1^2]} +
 \frac{ u +1}{ u - 1} P_{[2]}+
  \frac{u  - \frac{N}{2}+1}{ u  +\frac{N}{2}-1}  P_{[\emptyset]}
 \ee
 This $R$-matrix is called Zamolodchikov's solution \cite{3}
 of the Yang-Baxter equation.

\vspace{0.1cm}

For $\mathfrak{sp}(N)$ algebras the representations $\sf [1^2]$  and  ${\sf [\emptyset]}$
belong to the antisymmetric part of ${\sf [1]}^{\otimes 2}$,
while the adjoint representation
${\sf [2]}={\rm ad}$ belongs to the symmetric part of
${\sf [1]}^{\otimes 2}$, and we have (cf. (\ref{Rmat11}), (\ref{Rmat11b}))
${\sf [2]} \subset {\sf [2]} \otimes {\sf [\emptyset]}$,
${\sf [2]} \subset {\sf [2]} \otimes {\sf [1^2]}$. Therefore, in view
 of (\ref{Rmat12b}), the system
 of equations (\ref{Rmat07}) is written as
 \be
 \lb{Rmat13b}
 \frac{\tau_{[\emptyset]}(u)}{\tau_{[2]}(u)} =
 \frac{u  - \frac{N}{2}-1}{ u  +\frac{N}{2}+1}
 \; , \;\;\;\;\quad
  \frac{\tau_{[1^2]}(u)}{\tau_{[2]}(u)} =
 \frac{ u -1}{ u + 1} \; ,
 \ee
 and after a renormalization the $\mathfrak{sp}$-invariant
  $R$-matrix (\ref{Rmat02}) is
  \be
 \lb{Rmat14b}
 R(u) = P_{[2]} +
 \frac{\tau_{[1^2]}(u)}{\tau_{[2]}(u)} P_{[1^2]}+
 \frac{\tau_{[\emptyset]}(u)}{\tau_{[2]}(u)} P_{[\emptyset]}
 =  P_{[2]} +
 \frac{ u -1}{ u + 1} P_{[1^2]}+
  \frac{u  - \frac{N}{2}-1}{ u  +\frac{N}{2}+1}  P_{[\emptyset]}
 \ee

 \vspace{0.2cm}

  {\bf 2. The $\mathfrak{g}_2$ algebra} \\
  Further we denote by $[[{\sf n}]]$ the ${\sf n}$-dimensional
 representation of the Lie algebra $\mathfrak{g}$.
  The $R$-matrix operator (\ref{Rmat02}) of the
algebra $\mathfrak{g}_2$  in the minimal
fundamental representation $[[{\sf 7}]]$ acts in the reducible
49-dimensional space $[[{\sf 7}]]^{\otimes 2}$ which can be
expanded in the irreducible components as follows:
 \be
 \lb{g2jj}
 [[7]] \otimes [[7]] = \mathbb{S}([[7]]^{\otimes 2}) +
 \mathbb{A}([[7]]^{\otimes 2})=
 ([[1]] + [[27]]) + ([[7]] + [[14]]) \; .
 \ee
 Here the fundamental $[[7]]$ and adjoint $[[14]]$ representations
 of $\mathfrak{g}_2$ embed into the antisymmetric $\mathbb{A}$
 part of $[[7]]^{\otimes 2}$, while
 the representations $[[1]]$ and  $[[27]]$ compose
 the symmetric $\mathbb{S}$ part of $[[7]]^{\otimes 2}$.
 The highest weight vectors $\mu_{[n]}$ of the
 representations $[n]$
 in the right hand side of (\ref{g2jj}) and
 Weyl vector $\delta$ are \cite{Burb}
 (see also \cite{IsRub2} and references therein)
  \be
 \lb{g2wv}
 \begin{array}{c}
 \mu_{[[1]]} = (0,0,0) \; , \;\;\;\;
  \mu_{[[7]]} = \lambda_{(1)} = (0,-1,1) \; , \;\;\;\;
    \mu_{[[14]]} = \lambda_{(2)} = (-1,-1,2) \; , \\ [0.2cm]
      \mu_{[[27]]} = 2 \lambda_{(1)} = (0,-2,2) \; , \;\;\;\;
    \delta = \sum\limits_{i=1}^2 \lambda_{(i)} = (-1,-2,3) \; ,
      \end{array}
 \ee
  where $\lambda_{(1)}$ and $\lambda_{(2)}$ are fundamental
  weights of $\mathfrak{g}_2$, and  we describe
  the root space of the rank-$2$ Lie algebra
 $\mathfrak{g}_2$ as a plane ${\sf P}_u$ in 3-dimensional
 Euclidean space $\mathbb{R}^3$, normal to vector $u=(1,1,1)$.
 The values (\ref{kvkaz}) of the quadratic Casimir operators
 of the representations of the
 highest weights (\ref{g2wv}) are written as
   \be
 \lb{g2c2}
 \begin{array}{c}
 C_2^{[[1]]} = 0 \; , \;\;\;\;C_2^{[[7]]} = 12 \; , \;\;\;\;
 C_2^{[[14]]} = 24 \; , \;\;\;\; C_2^{[[27]]} = 28 \; ,
 \end{array}
 \ee
 For the case of Lie algebra $\mathfrak{g}_2$
 the conditions (\ref{Rmat11}) are (see, e.g., \cite{Yamat} about
 tensor product of $\mathfrak{g}_2$ representations)
 $$
 [[1]] \subset {\rm ad} \otimes [[14]] \; , \;\;\;\;
 [[27]] \subset {\rm ad} \otimes [[7]] \; , \;\;\;\;
 [[27]] \subset {\rm ad} \otimes [[14]] \; ,
 $$
 where ${\rm ad} = [[14]]$, and we write the system of equations
 (\ref{Rmat07}) as
  \be
 \lb{g207}
 \frac{\tau_{[[1]]}(u)}{\tau_{[[14]]}(u)} = \frac{u-6}{u+6}
 \; , \;\;\;\;\;\;
  \frac{\tau_{[[7]]}(u)}{\tau_{[[27]]}(u)} =
  \frac{u-4}{u+4} \; , \;\;\;\;\;\;
  \frac{\tau_{[[14]]}(u)}{\tau_{[[27]]}(u)} =
 \frac{ u  -1}{ u  +1} \; .
 \ee
 Finally, after a normalization, the $\mathfrak{g}_2$-invariant
$R$-matrix (\ref{Rmat02}) in the fundamental representation $[[7]]$
acquires the form \cite{OgWieg}
(see also \cite{IsKriv} with $u \to -u$)
\be\label{g2R}
 \begin{array}{c}
 \displaystyle
 R(u) =  \frac{\tau_{[[1]]}(u)}{\tau_{[[27]]}(u)} P_{[[1]]}
 + \frac{\tau_{[[7]]}(u)}{\tau_{[[27]]}(u)} P_{[[7]]} +\frac{\tau_{[14]]}(u)}{\tau_{[[27]]}(u)} P_{[[14]]} +P_{[[27]]} =
 \\ [0.4cm]
 \displaystyle = \frac{(u-6)(u-1)}{(u+6)(u+1)} P_{[[1]]}
 +\frac{u-4}{u+4} P_{[[7]]} +\frac{u-1}{u+1}P_{[[14]]} +P_{[[27]]}  \; .
 \end{array}
\ee
\vspace{0.2cm}
  {\bf 3. The $\mathfrak{f}_4$ algebra} \\
  The $\mathfrak{f}_4$-invariant $R$-matrix
  (\ref{Rmat02}) in the minimal
fundamental representation $[[{\sf 26}]]$ acts in the reducible
676-dimensional space $[[{\sf 26}]]^{\otimes 2}$ which is
expanded in the irreducible components as follows:
\be
 \lb{f4jj}
 [[26]]^{\otimes 2} =
 \mathbb{S}([[26]]^{\otimes 2}) +
 \mathbb{A}([[26]]^{\otimes 2})
 = ([[1]] + [[26]] + [[324]]) + ([[52]] + [[273]]) \, ,
 \ee
 where representations $[[{\sf 1}]],
  [[{\sf 26}]]$ and  $[[{\sf 324}]]$
 belong to the symmetric part of $[[26]]^{\otimes 2}$,
 while the adjoint
 representation $[[{\sf 52}]]$ and representation $[[{\sf 273}]]$  belong to the antisymmetric
 part of $[[26]]^{\otimes 2}$.

  The highest weight vectors of the representations
 in (\ref{f4jj}) and $\mathfrak{f}_4$ Weyl vector
 $\delta$ are \cite{Burb} (see also \cite{IsRub2} and references therein)
  \be
 \lb{f4wv}
 \begin{array}{c}
 \mu_{[[1]]} = (0,0,0,0) \; , \;\;\;\;
  \mu_{[[26]]} = \lambda_{(4)} = (0,0,0,1) \; , \;\;\;\;
    \mu_{[[52]]} = \lambda_{(1)} = (1,0,0,1) \; , \\ [0.2cm]
      \mu_{[[273]]} = \lambda_{(3)} =
      \frac{1}{2}(1,1,1,3) \; , \;\;\;\;
      \mu_{[[324]]} = 2 \lambda_{(4)} = (0,0,0,2) \; , \;\;\;\;
    \delta = \sum\limits_{i=1}^4 \lambda_{(i)} =
    \frac{1}{2}\,(5,3,1,11) \; ,
      \end{array}
 \ee
  where $\lambda_{(i)}|_{i=1,...,4}$  are fundamental
  weights of $\mathfrak{f}_4$. Now we deduce
 the values (\ref{kvkaz}) of the quadratic Casimir operators,
 which correspond to the highest weights (\ref{g2wv})
   \be
 \lb{f4c2}
 \begin{array}{c}
 C_2^{[[1]]} = 0 \; , \;\;\;\;C_2^{[[26]]} = 12 \; , \;\;\;\;
 C_2^{[[52]]} = 18 \; , \;\;\;\; C_2^{[[273]]} = 24
 \; , \;\;\;\; C_2^{[[324]]} = 26 \; ,
 \end{array}
 \ee
 The analogs of the conditions (\ref{Rmat11})
  for the Lie algebra $\mathfrak{f}_4$ have the form
 $$
 [[52]] \subset {\rm ad} \otimes [[1]] \; , \;\;\;\;
 [[324]] \subset {\rm ad} \otimes [[52]] \; , \;\;\;\;
 [[324]] \subset {\rm ad} \otimes [[273]] \; , \;\;\;\;
 [[273]] \subset {\rm ad} \otimes [[26]] \; ,
 $$
 where ${\rm ad} \equiv [[52]]$, and the system of equations
 (\ref{Rmat07}) is represented as
  \be
 \lb{f407}
 \frac{\tau_{_{[[1]]}}(u)}{\tau_{_{[[52]]}}(u)} =
 \frac{u-\frac{9}{2}}{u+\frac{9}{2}}
 \, , \;\;\;\;
  \frac{\tau_{[_{[[52]]}]}(u)}{\tau_{_{[[324]]}}(u)} =
  \frac{u-2}{u+2} \, , \;\;\;\;
  \frac{\tau_{_{[[273]]}}(u)}{\tau_{_{[[324]]}}(u)} =
 \frac{ u  -\frac{1}{2}}{ u  +\frac{1}{2}} \, , \;\;\;\;
  \frac{\tau_{_{[[26]]}}(u)}{\tau_{_{[[273]]}}(u)} =
 \frac{ u  -3}{ u  +3} \, .
 \ee
  Finally, after a normalization, the $\mathfrak{f}_4$-invariant
$R$-matrix in the fundamental representation $[[26]]$
has the form \cite{OgWieg}
(see also \cite{IsKriv} with $u \to -2u$)
\be\label{f4R}
 \begin{array}{c}
 R(u) =  \frac{\tau_{_{[[1]]}}(u)}{\tau_{_{[[324]]}}(u)} P_{[[1]]}
 + \frac{\tau_{_{[[26]]}}(u)}{\tau_{_{[[324]]}}(u)} P_{[[26]]} +\frac{\tau_{_{[[52]]}}(u)}{\tau_{_{[[324]]}}(u)} P_{[[52]]} +
 \frac{\tau_{_{[[273]]}}(u)}{\tau_{_{[[324]]}}(u)} P_{[[273]]} +
 P_{[[324]]} =
 \\ [0.4cm]
  = \frac{(u-9/2)(u-2)}{(u+9/2)(u+2)} P_{[[1]]} + \frac{(u-3)(u-1/2)}{(u+3)(u+1/2)} P_{[[26]]}  +
\frac{u-2}{u+2} P_{[[52]]} - \frac{u-1/2}{u+1/2}  P_{[[273]]}
 +  P_{[[324]]} .
 \end{array}
\ee
\vspace{0.2cm}
  {\bf 3. The $\mathfrak{e}_6$ algebra} \\
 The 78-dimensional algebra $\mathfrak{e}_6$
  has two inequivalent minimal fundamental representations
$[[{\sf 27}]]$ and $[[\overline{\sf 27}]]$. Here
we consider $\mathfrak{e}_6$-invariant $R$-matrix (\ref{Rmat02})
which acts in the space of reducible representation
 \be
 \lb{e6jj}
 [[{\sf 27}]]^{\otimes 2} =
 \mathbb{S}([[{\sf 27}]]^{\otimes 2}) +
 \mathbb{A}([[{\sf 27}]]^{\otimes 2})
 = ([[\overline{\sf 27}]] + [[{\sf 351}]]_1) \, +
 \, ([[{\sf 351}]]_2) \, ,
 \;\;\;\;\;  {\rm ad} \equiv [[{\sf 78}]] \, .
 \ee
   The highest weight vectors $\mu_{[[n]]}$ of the representations
 in (\ref{e6jj}) and  Weyl vector
 $\delta$ for $\mathfrak{e}_6$ are \cite{Burb}, \cite{IsRub2}
  \be
 \lb{e6wv}
 \begin{array}{c}
 \mu_{[[{\sf 27}]]} = \lambda_{(1)} =
 (-\frac{1}{3},-\frac{1}{3},1,0,0,0,0,\frac{1}{3}) \; , \;\;\;\;
  \mu_{[[\overline{\sf 27}]]} = \lambda_{(6)} =
  \frac{2}{3}(-1,-1,0,0,0,0,0,1) \; , \\ [0.2cm]
    \mu_{[[{\sf 351}]]_1} = 2 \lambda_{(1)} =
    (-\frac{2}{3},-\frac{2}{3},2,0,0,0,0,\frac{2}{3}) \; , \;\;\;\;
      \mu_{[[{\sf 351}]]_2} = \lambda_{(2)} =
     (-\frac{2}{3},-\frac{2}{3},1,1,0,0,0,\frac{2}{3}) \; , \\ [0.2cm]
    \delta = \sum\limits_{i=1}^6 \lambda_{(i)} =
   (-4,-4,4,3,2,1,0,4) \; ,
      \end{array}
 \ee
  where $\lambda_{(i)}|_{i=1,...,6}$  are fundamental
  weights of $\mathfrak{e}_6$, and we numerate nodes of
  $\mathfrak{e}_6$ Dynkin diagram as follows

  \unitlength=4mm
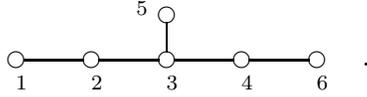
\begin{picture}(17,4)(-5,0)


\put(5.5,2){\circle{0.5}}
\put(5.75,2){\line(1,0){2}}
\put(8,2){\circle{0.5}}
\put(8.25,2){\line(1,0){2}}
\put(10.5,2){\circle{0.5}}
\put(10.75,2){\line(1,0){2}}
\put(13,2){\circle{0.5}}
\put(13.25,2){\line(1,0){2}}
\put(15.5,2){\circle{0.5}}

\put(10.5,3.5){\circle{0.5}}
\put(10.5,2.25){\line(0,1){1}}

\put(5.5,1){\scriptsize  1}
\put(8,1){\scriptsize 2}
\put(10.5,1){\scriptsize 3}
\put(13,1){\scriptsize 4}
\put(15.5,1){\scriptsize 6}
\put(9.5,3.5){\scriptsize 5}

\put(17,1.8){.}

\end{picture}

\noindent
The values (\ref{kvkaz}) of the quadratic Casimir operators,
 which correspond to the representations with
  highest weights (\ref{e6wv}), are
   \be
 \lb{e6c2}
 \begin{array}{c}
 C_2^{[[\overline{27}]]} = \frac{52}{3} \; , \;\;\;\;
 C_2^{[[351]]_1} = \frac{112}{3} \; , \;\;\;\;
 C_2^{[[351]]_2} = \frac{100}{3} \; ,
 \end{array}
 \ee
 We note, that in view of the symmetry of the
 $\mathfrak{e}_6$ Dynkin diagram, the
 values (\ref{e6c2}) are invariant under the change of  the
 fundamental weights in  (\ref{e6wv}):
 $\lambda_{(1)} \leftrightarrow \lambda_{(6)}$
 and $\lambda_{(2)} \leftrightarrow \lambda_{(4)}$.
 This symmetry also means that the $R$-matrix in the representation
 $[[\overline{27}]]^{\otimes 2}$ has the same form
 as the $R$-matrix in the representation
 $[[27]]^{\otimes 2}$.
 The analogs of the conditions (\ref{Rmat11})
  for the Lie algebra $\mathfrak{e}_6$ are
  (see e.g. \cite{Yamat})
   $$
 [[351]]_2 \subset {\rm ad} \otimes [[\overline{27}]] \; , \;\;\;\;
 [[351]]_1 \subset {\rm ad} \otimes [[351]]_2 \; ,
 $$
 where ${\rm ad} \equiv [[78]]$, and the system of equations
 (\ref{Rmat07}) is written as
  \be
 \lb{e607}
  \frac{\tau_{_{[[\overline{27}]]}}(u)}{\tau_{_{[[351]]_2}}(u)} =
  \frac{u-4}{u+4} \, , \;\;\;\;\;\;\;\;
  \frac{\tau_{_{[[351]]_1}}(u)}{\tau_{_{[[351]]_2}}(u)} =
 \frac{ u  +1}{ u  -1}  \, .
 \ee
 Finally, the $\mathfrak{e}_6$-invariant $R$-matrix
 in the representation $[[27]]$ acquires the form
 \cite{OgWieg}, \cite{IsKriv}
\be\label{e6R}
\begin{array}{c}
R(u) = \frac{\tau_{_{[[\overline{27}]]}}(u)}{\tau_{_{[[351]]_2}}(u)} P_{[[\overline{27}]]}
+\frac{\tau_{_{[[351]]_1}}(u)}{\tau_{_{[[351]]_2}}(u)} P_{[[351]]_1}
+ P_{[[351]]_2} = \\ [0.3cm]
=\frac{u-4}{u+4} P_{[[\overline{27}]]}
 +\frac{u+1}{u-1}  P_{[[351]]_1}
+  P_{[[351]]_2} \; .
 \end{array}
\ee
 {\bf 4. The $\mathfrak{e}_7$ algebra} \\
 Here
we consider $\mathfrak{e}_7$-invariant $R$-matrix (\ref{Rmat02})
which acts in the space of the representation
 \be
 \lb{e7jj}
[[{\sf 56}]]^{\otimes 2} = \mathbb{S}([[{\sf 56}]]^{\otimes 2})+
\mathbb{A}([[{\sf 56}]]^{\otimes 2}) =
\left( [[{\sf 133}]] + [[{\sf 1463}]]\right)  +
\left( [[{\sf 1}]] + [[{\sf 1539}]]\right) \, ,
 \ee
 where $[[{\sf 133}]] \equiv {\rm ad}$ is the adjoint
 representation of $\mathfrak{e}_7$.
   The highest weight vectors $\mu_{[[n]]}$ of the representations
 in (\ref{e7jj}) and  Weyl vector
 $\delta$ for the algebra
 $\mathfrak{e}_7$ are \cite{Burb}, \cite{IsRub2}
  \be
 \lb{e7wv}
 \begin{array}{c}
 \mu_{[[1]]} = (0,0,0,0,0,0,0,0) \; , \;\;\;\;
  \mu_{[[{\sf 133}]]} = \lambda_{(7)} =
  (-1,0,0,0,0,0,0,1) \; , \\ [0.2cm]
    \mu_{[[{\sf 1463}]]} = 2 \lambda_{(1)} =
    (-1,2,0,0,0,0,0,1) \; , \;\;\;\;
      \mu_{[[{\sf 1539}]]} = \lambda_{(2)} =
     (-1,1,1,0,0,0,0,1) \; , \\ [0.2cm]
     \mu_{[[56]]} = \lambda_{(1)} =
     (-\frac{1}{2},1,0,0,0,0,0,\frac{1}{2}) \; , \;\;\;\;
    \delta = \sum\limits_{i=1}^7 \lambda_{(i)} =
   (-\frac{17}{2},5,4,3,2,1,0,\frac{17}{2}) \; ,
      \end{array}
 \ee
  where $\lambda_{(i)}|_{i=1,...,7}$  are fundamental
  weights of $\mathfrak{e}_7$ and we numerate nodes
  in Dynkin diagram as following

  \unitlength=4mm
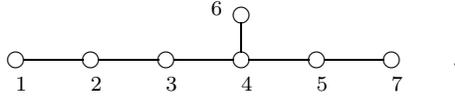
\begin{picture}(17,3.5)


\put(5.5,2){\circle{0.5}}
\put(5.75,2){\line(1,0){2}}
\put(8,2){\circle{0.5}}
\put(8.25,2){\line(1,0){2}}
\put(10.5,2){\circle{0.5}}
\put(10.75,2){\line(1,0){2}}

\put(13,2){\circle{0.5}}
\put(13.25,2){\line(1,0){2}}
\put(15.5,2){\circle{0.5}}
\put(15.75,2){\line(1,0){2}}
\put(18,2){\circle{0.5}}

\put(13,3.5){\circle{0.5}}
\put(13,2.25){\line(0,1){1}}

\put(5.5,1){\scriptsize  1}
\put(8,1){\scriptsize 2}
\put(10.5,1){\scriptsize 3}
\put(13,1){\scriptsize 4}
\put(15.5,1){\scriptsize 5}
\put(18,1){\scriptsize 7}
\put(12,3.5){\scriptsize 6}

\put(20,1.8){.}

\end{picture}

   The values (\ref{kvkaz}) of the quadratic Casimir operators,
 which correspond to the representations with
  highest weights (\ref{e7wv}), are
   \be
 \lb{e7c2}
 \begin{array}{c}
 C_2^{[[56]]} = \frac{57}{2} \; , \;\;\;\;
  C_2^{[[133]]} = 36 \; , \;\;\;\;
 C_2^{[[1463]]} = 60 \; , \;\;\;\;
 C_2^{[[1539]]} = 56 \; ,
 \end{array}
 \ee
 The analogs of the conditions (\ref{Rmat11})
  for the Lie algebra $\mathfrak{e}_7$ have the form
  (see e.g. \cite{Yamat})
   $$
 [[133]] \subset {\rm ad} \otimes [[1]] \; , \;\;\;\;
 [[1539]] \subset {\rm ad} \otimes [[1463]] \; , \;\;\;\;
 [[1539]] \subset {\rm ad} \otimes [[133]] \; ,
 $$
 where ${\rm ad} \equiv [[133]]$, and the system of equations
 (\ref{Rmat07}) is written as
  \be
 \lb{e707}
 \frac{\tau_{[[1]]}(u)}{\tau_{[[133]]}(u)} =
  \frac{u-9}{u+9} \, , \;\;\;\;\;\;\;\;
  \frac{\tau_{[[133]]}(u)}{\tau_{[[1539]]}(u)} =
  \frac{u-5}{u+5} \, , \;\;\;\;\;\;\;\;
  \frac{\tau_{[[1463]]}(u)}{\tau_{[[1539]]}(u)} =
 \frac{ u  +1}{ u  -1}  \, .
 \ee
 Thus, the $\mathfrak{e}_7$-invariant solution
(\ref{Rmat02}) of the Yang-Baxter equation
 in the defining representation $[[56]]$ has  the form
\be\lb{Re7}
R(u)   = \frac{(u-9)(u-5)}{(u+9)(u+5)} P_{[[1]]}  +
\frac{u-5}{u+5} P_{[[133]]} +
 \frac{u+1}{u-1}P_{[[1463]]} + P_{[[1539]]} .
\ee
{\bf 5. The $\mathfrak{e}_8$ algebra} \\
For the exceptional Lie algebra $\mathfrak{e}_8$
the adjoint and minimal fundamental representations
coincide and have dimension 248.
The tensor product of its two 248-dimensional representations has the following decomposition into irreducible representations \cite{Cvit,Yamat}:
\be\label{e8jj}
\begin{array}{c}
[[{\sf 248}]]^{\otimes 2} =
\mathbb{S}([[{\sf 248}]]^{\otimes 2})+
\mathbb{A}([[{\sf 248}]]^{\otimes 2}) 
= \left([[{\sf 1}]]+ [[{\sf 3875}]] + [[{\sf 27000}]]\right)  +
\left( [[{\sf 248}]] + [[{\sf 30380}]]\right) \, .
\end{array}
\ee
   The highest weight vectors $\mu_{[[n]]}$ of the representations
 in (\ref{e8jj}) and  Weyl vector
 $\delta$ for the algebra
 $\mathfrak{e}_8$ are \cite{Burb}, \cite{IsRub2}
  \be
 \lb{e8wv}
 \begin{array}{c}
 \mu_{[[248]]} = \lambda_{(1)}= (1,0,0,0,0,0,0,1) \; , \;\;\;\;
  \mu_{[[{\sf 3875}]]} = \lambda_{(8)} =
  (0,0,0,0,0,0,0,2) \; , \\ [0.2cm]
    \mu_{[[{\sf 27000}]]} = 2 \lambda_{(1)} =
    (2,0,0,0,0,0,0,2) \; , \;\;\;\;
      \mu_{[[{\sf 30380}]]} = \lambda_{(2)} =
     (1,1,0,0,0,0,0,2) \; , \\ [0.2cm]
    \delta = \sum\limits_{i=1}^8 \lambda_{(i)} =
   (6,5,4,3,2,1,0,23) \; ,
      \end{array}
 \ee
  where $\lambda_{(i)}|_{i=1,...,8}$  are fundamental
  weights of $\mathfrak{e}_8$ and we numerate nodes
  in Dynkin diagram as following

  \unitlength=4mm
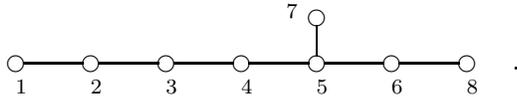
\begin{picture}(22,3.5)


\put(5.5,1.5){\circle{0.5}}
\put(5.5,0.5){\scriptsize  1}
\put(5.75,1.5){\line(1,0){2}}
\put(8,1.5){\circle{0.5}}
\put(8,0.5){\scriptsize 2}
\put(8.25,1.5){\line(1,0){2}}
\put(10.5,1.5){\circle{0.5}}
\put(10.5,0.5){\scriptsize 3}
\put(10.75,1.5){\line(1,0){2}}
\put(13,1.5){\circle{0.5}}
\put(13,0.5){\scriptsize 4}
\put(13.25,1.5){\line(1,0){2}}

\put(15.5,1.5){\circle{0.5}}
\put(15.5,0.5){\scriptsize 5}
\put(15.75,1.5){\line(1,0){2}}
\put(18,1.5){\circle{0.5}}
\put(18,0.5){\scriptsize 6}
\put(20.5,1.5){\circle{0.5}}
\put(20.5,0.5){\scriptsize 8}
\put(18.25,1.5){\line(1,0){2}}

\put(15.5,3){\circle{0.5}}
\put(14.5,3){\scriptsize 7}
\put(15.5,1.75){\line(0,1){1}}

\put(22,1.3){.}

\end{picture}

   The values (\ref{kvkaz}) of the quadratic Casimir operators,
 which correspond to the representations with
  highest weights (\ref{e8wv}), are
   \be
 \lb{e8c2}
 \begin{array}{c}
 C_2^{[[248]]} = 60 \; , \;\;\;\;
  C_2^{[[3875]]} = 96 \; , \;\;\;\;
 C_2^{[[27000]]} = 124 \; , \;\;\;\;
 C_2^{[[30380]]} = 120 \; ,
 \end{array}
 \ee
 The analogs of the conditions (\ref{Rmat11})
  for the Lie algebra $\mathfrak{e}_8$ have the form
  (see e.g. \cite{Yamat})
   $$
   \begin{array}{c}
 [[248]] \subset {\rm ad} \otimes [[1]] \; , \;\;\;\;
 [[3875]] \subset {\rm ad} \otimes [[248]] \; , \;\;\;\;
 [[27000]] \subset {\rm ad} \otimes [[248]] \; , \\ [0.2cm]
 [[30380]] \subset {\rm ad} \otimes [[27000]] \; ,\;\;\;\;
 [[30380]] \subset {\rm ad} \otimes [[3875]]\; ,
 \end{array}
 $$
 where ${\rm ad} \equiv [[248]]$, and the part of
 the system of equations
 (\ref{Rmat07}) is written as
  \be
 \lb{e807}
 \frac{\tau_{[[248]]}(u)}{\tau_{[[3875]]}(u)} =
  \frac{u-9}{u+9} \, , \;\;\;\;
  \frac{\tau_{[[248]]}(u)}{\tau_{[[27000]]}(u)} =
  \frac{u-16}{u+16} \, , \;\;\;\;
  \frac{\tau_{[[27000]]}(u)}{\tau_{[[30380]]}(u)} =
 \frac{ u  +1}{ u  -1}  \, , \;\;\;\;
  \frac{\tau_{[[30380]]}(u)}{\tau_{[[3875]]}(u)} =
 \frac{ u  +6}{ u  -6} \, .
 \ee
 It is clear that equations (\ref{e807}) are inconsistent.
 Thus, the $\mathfrak{e}_8$-invariant $R$-matrix
 in the minimal fundamental (adjoint)
 representation $[[248]]$ does not exist. This fact is in agreement with the general statement  \cite{Drin85}, \cite{ChaPr}
 that the adjoint representation of $\mathfrak{e}_8$
 can not be lifted to the representation of the Yangian
 $Y(\mathfrak{e}_8)$.

  \vspace{0.2cm}

In paper \cite{IsKriv} the $R$-matrix
solutions (\ref{Rmat14}), (\ref{Rmat14b}),
(\ref{g2R}), (\ref{f4R}), (\ref{e6R}) and (\ref{Re7})
are written as rational functions of the split Casimir
operators (\ref{kaz-01}) in the defining representations. In particular,
for Lie algebras
$\mathfrak{g} = \mathfrak{so},\mathfrak{sp},\mathfrak{e}_6$
 we have \cite{IsKriv}
$$
R(u) = \frac{u+ \frac{1}{2}(\hat{C} + a_{\mathfrak{g}})}{u-
 \frac{1}{2}(\hat{C} + a_{\mathfrak{g}})}  \; ,
$$
where $a_{\mathfrak{so}}=1$,
$a_{\mathfrak{sp}}=-1$, $a_{\mathfrak{e}_6}=2/3$
 and the split Casimir operators $\hat{C}$ are
normalized such that $\hat{C}^{ki}_{jk} = C_2(\lambda) \delta^i_j$
(here $\lambda$ is a highest weight
of the defining representation of $\mathfrak{g}$
and the value $C_2(\lambda)$ is given by (\ref{kvkaz})).

\subsubsection{Yangians $Y(\mathfrak{g})$}

The Yangians $Y(\mathfrak{g})$ can be defined by means of the
 rational $R$-matrix solutions (\ref{Rmatr}), (\ref{Rmat02})
 of the Yang-Baxter equations (\ref{YBEu}). First, we consider equations
 (cf. (\ref{3.5.10})
 \be
\lb{yang01}
R_{12}(u -v) \, L_1(u) \, L_2(v) =
L_2(v) \, L_1(u) \, R_{12}(u - v) \; ,
\ee
where $u,v$ are spectral parameters, indices $1,2$
numerate the spaces $V$ of the defining representation $T$
 in the product $V \otimes V$, and
$R_{12}(u)$ is the $\mathfrak{g}$-invariant
$R$-matrix in the representation $T\otimes T$.
We search the elements $L^i_{\; j}(u)$ of the quadratic algebra
(\ref{yang01}) in the form (cf. (\ref{3.5.12}))
\be
\lb{yang02}
L^{i}_{\; j}(u) = \delta^{i}_{j} +
\sum_{k=1}^{\infty} {T^{(k)}}^{i}_{j} u^{-k} = \delta^{i}_{j}
+ \frac{1}{u} {\bf I}^{i}_{j}
+ \frac{1}{u^2} \Bigl(\frac{1}{2}({\bf I}^2)^{i}_{\; j}
+ {\bf J}^{i}_{\; j}\Bigr) + \frac{1}{u^3} \dots \; ,
\ee
where ${T^{(k)}}^{i}_{j}$ $(k >1)$ are the generators of the
Yangian $Y(\mathfrak{g})$ \cite{13}, \cite{Drin85} and
we introduce notation
$$
{T^{(1)}}^{i}_{j} \equiv
{\bf I}^{i}_{j} = I^a (T_a)^{i}_{\; j} \; , \;\;\;\;
{T^{(2)}}^{i}_{j} \equiv \frac{1}{2}({\bf I}^2)^{i}_{\; j}
+ {\bf J}^{i}_{\; j}
\; , \;\;\;\; {\bf J}^{i}_{\; j} := J^a (T_a)^{i}_{\; j} \; .
$$
 Here $(T_a)^{i}_{\; j} = T^{i}_{\; j}(X_a)$  are generators
 of $\mathfrak{g}$ in the representation $T$.
 Now we substitute expansions (\ref{Rmatr}),
 (\ref{XX}) and (\ref{yang02})
 into (\ref{yang01}), multiply both sides by $(u-v)^2$
 and consider eq. (\ref{yang01})
 in the limit $u,v \to \infty$. We take into account identities
 (cf. (\ref{Rmat04}))
 $$
 ({\bf I}_1 + {\bf I}_2) \hat{C}_{12} =
 \hat{C}_{12} ({\bf I}_1 + {\bf I}_2) \; , \;\;\;
 ({\bf J}_1 + {\bf J}_2) \hat{C}_{12} =
 \hat{C}_{12} ({\bf J}_1 + {\bf J}_2) \; ,
 $$
 where $\hat{C}_{12}$
 is the split Casimir operator in the defining representation $T$.
 Then the terms of zero order in $u,v$ of eqs. (\ref{yang01})
 give relations (\ref{lialg}):
 \be
 \lb{yang03}
 [{\bf I}_1, \; {\bf I}_2] =
 \hat{C}_{12} {\bf I}_1 - {\bf I}_1 \hat{C}_{12}
 \;\;\;\;\;\; \Rightarrow  \;\;\;\;\;\;
 [I_a, \; I_b] = X^d_{ab} \; I_d \; ,
 \ee
 which means that coefficients $I_a$ are the basis elements of the
 Lie algebra $\mathfrak{g}$. The terms of order $u^{-2}v$
 and $uv^{-2}$ of eqs. (\ref{yang01})
 give commutation relations
 \be
 \lb{yang04}
 [{\bf I}_1, \; {\bf J}_2] =
 \hat{C}_{12} {\bf J}_1 - {\bf J}_1 \hat{C}_{12} \; \equiv \;
  {\bf J}_2 \hat{C}_{12} -\hat{C}_{12} {\bf J}_2
  \;\;\;\;\;\; \Rightarrow  \;\;\;\;\;\;
 [I_a, \; J_b] = J_d \; X^d_{ab}  \; ,
 \ee
 which means that the elements $J_b$ form the adjoint representation
 of $\mathfrak{g}$. We also note useful generalizations of
 (\ref{yang03}) and (\ref{yang04})
  \be
 \lb{yang07}
 [{\bf I}_1, \; {\bf I}_2^k] =
 {\bf I}_2^k \, \hat{C}_{12}  - \hat{C}_{12} \, {\bf I}_2^k
 \; , \;\;\;\;\;\;
 [{\bf I}_1, \; {\bf J}_2^k] =
 {\bf J}_2^k \, \hat{C}_{12}  - \hat{C}_{12} \, {\bf J}_2^k   \; .
 \ee
 The commutator $[{\bf J}_1, \; {\bf J}_2]$ is not
 fully specified by eq. (\ref{yang01}) if we know
 the expansion of $R$-matrix (\ref{Rmatr})
 only up to the order $u^{-2}$.

 We define the Yangian $Y(\mathfrak{g})$ as
  the enveloping algebra generated by the basis
 elements $I_a|_{a=1,...,\dim \mathfrak{g}}$ of the
 Lie algebra $\mathfrak{g}$, the additional
set of elements $J_a|_{a=1,...,\dim \mathfrak{g}}$,
which form the adjoint representation (\ref{yang04}) of
$\mathfrak{g}$ and with
 a nontrivial noncommutative coproduct $\Delta$ : $Y(\mathfrak{g}) \to
Y(\mathfrak{g}) \otimes Y(\mathfrak{g})$ which is defined by
 \be
 \lb{yang05}
 \Delta(L^{i}_{\; j}(u) ) =
 L^{i}_{\; k}(u)  \otimes L^{k}_{\; j} (u) \; .
 \ee
The substitution of (\ref{yang02}) into (\ref{yang05}) gives
 \be
 \lb{yang06}
 \begin{array}{c}
 \Delta ({\bf I}^i_{\; j}) = {\bf I}^i_{\; j} \otimes  1
  + 1 \otimes {\bf I}^i_{\; j} \;\;\; \Rightarrow \;\;\;
 \Delta (I_a) = I_a \otimes  1 + 1 \otimes I_a \; , \\ [0.2cm]
 \displaystyle
  \Delta ({\bf J}^i_{\; j}) = {\bf J}^i_{\; j} \otimes  1
  + 1 \otimes {\bf J}^i_{\; j} +
  \frac{1}{2} \bigl( {\bf I}^i_{\; k} \otimes  {\bf I}^k_{\; j} -
 {\bf I}^k_{\; j} \otimes {\bf I}^i_{\; k} \bigr)
  \;\;\; \Rightarrow \;\;\; \\ [0.2cm]
 \displaystyle
 \Delta (J_a) = J_a \otimes  1 + 1 \otimes J_a
 + \frac{1}{2} X_{abc} \; I^b \otimes I^c \; ,
 \end{array}
 \ee
where $X_{abc} = X_{ab}^d \; {\sf g}_{dc}$. Finally, we remark
that the commutator $[J_a, \; J_b]$ is constrained by the requirement
that $\Delta$ be a homomorphism. Indeed, we have
 \cite{Drin85}, \cite{MaKa2}:
 $$
 \begin{array}{c}
 [J_a,\, [J_b,\, I_c]] - [I_a,\, [J_b,\, J_c]] = {a_{abc}}^{efg} \;
 \{I_e,I_f,I_g \} \; , \\ [0.2cm]
 [[J_a,\, J_b],\, [I_c,\, J_d]] + [[J_c,\, J_d],\, [I_a,\, J_b]]
  = ({a_{abh}}^{efg} \; X_{h cd} + {a_{cdh}}^{efg} \; X_{h ab})
   \{I_e,I_f,J_g \}
  \end{array}
 $$
 where ${a_{abh}}^{efg} = \frac{1}{24} X_{ia}^{e}X_{jb}^{f}X_{kh}^{g} X^{ijk}$
 and $\{x_1,x_2,x_3 \} = \sum\limits_{e\neq f \neq g} x_e,x_f,x_g$.

\subsection{\bf \em Quantum Knizhnik - Zamolodchikov equations}
\setcounter{equation}0

In Sections {\bf \ref{baxtel}} and {\bf \ref{ospbax}},
by using $R$-matrix representations for the
Hecke and Birman-Murakami-Wenzl algebras, we have found the trigonometric solutions
$R(x)$ of the Yang-Baxter equations (Baxterized $R$-matrices). In this
subsection we show that,
for every trigonometric solution $R(x)$
of the Yang-Baxter equations (\ref{3.5.2}),
one can construct the set of difference
equations which are called {\it quantum Knizhnik-Zamolodchikov equations}.
These equations are important since their solutions are related
(see e.g. \cite{JMi}, \cite{bks} and Refs. therein) to the
correlation functions in spin chain models associated with the same trigonometric matrix $R(x)$.

In this subsection we follow the presentation of the papers \cite{FreRe}, \cite{TaVa} and \cite{FAS}.

Consider a tensor function $\Psi^{1 \dots N \rangle}(z_1 , \dots , z_N)
\in V^{\otimes N}$
($z_i \in \mathbb{C}$, $i=1, \dots , N$)
which satisfies a system of difference equations
\be
\lb{qkz1}
T_{(i)} \, \Psi^{1 \dots N \rangle}(z_1 , \dots , z_N)  =
A^{(i)}_{1 \dots N}(z_1 , \dots , z_N) \,  \Psi^{1 \dots N \rangle}(z_1 , \dots , z_N) \; ,
\ee
where operator $T_{(i)}$ is defined as
\be
\lb{qkz2}
T_{(i)} \, \Psi^{1 \dots N \rangle}(z_1 , \dots , z_N)  :=
\Psi^{1 \dots N \rangle}(z_1 , \dots , z_{i-1} , p z_i , z_{i+1}, \dots , z_N) \; ,
\ee
$A^{(i)}_{1 \dots N}(z_1 , \dots , z_N) \in End(V^{\otimes N})$
is called discrete connection and
indices $1, \dots, N$ denote the numbers of the vector spaces $V$ in $V^{\otimes N}$.
A consistence condition $T_{(i)} T_{(j)} = T_{(j)} T_{(i)}$
of the system (\ref{qkz1}) requires additional constraints on
the discrete connection $A^{(j)}_{1 \dots N}(z_1 \dots z_N)$:
\be
\lb{qkz3}
T_{(i)} A^{(j)}_{1 \dots N} T_{(i)}^{-1}  A^{(i)}_{1 \dots N} =
T_{(j)} A^{(i)}_{1 \dots N} T_{(j)}^{-1}  A^{(j)}_{1 \dots N}
\;\;\;  \Rightarrow
\ee
\be
\lb{qkz33}
{\sf A}^{(j)}_{1 \dots N} \, {\sf A}^{(i)}_{1 \dots N} =
{\sf A}^{(i)}_{1 \dots N} \; {\sf A}^{(j)}_{1 \dots N}
\; , \;\;\;\;\;\; {\sf A}^{(i)} :=
T_{(i)}^{-1}  A^{(i)}_{1 \dots N}  \; .
\ee
Connections $A^{(i)}_{1...N}$ and ${\sf A}^{(j)}_{1...N}$
which satisfy (\ref{qkz3}) and (\ref{qkz33}) are called flat
(or integrable).

Now we introduce the following discrete connection \cite{FreRe}
\be
\lb{qkz4}
 \begin{array}{c}
A^{(j)}_{1 \dots N}(z_1 , \dots , z_N) =
T_{(j)} R_{j \,  j-1} \dots R_{j \,  2} \, R_{j \,  1} T_{(j)}^{-1} \, D_j \,
R^{-1}_{N\, j} R^{-1}_{N-1\, j} \dots R^{-1}_{j+1\, j} = \\ [0.3cm]
= R_{j \,  j-1}(\frac{z_j}{z_{j-1}}p) \dots R_{j \,  2}({z_j\over z_2}p)
 \, R_{j \,  1}({z_j\over z_1} \, p) \, D_j \,
R^{-1}_{N\, j}({z_N\over z_j}) R^{-1}_{N-1\, j}({z_{N-1}\over z_j})
 \dots R^{-1}_{j+1\, j}({z_{j+1}\over z_j})  \; ,
\end{array}
\ee
where $R_{i\, j} := R_{i\, j}(z_i/z_j)$ is the $R$-matrix which
acts nontrivially only in the vector spaces $V$ with numbers $i,j$
in $V^{\otimes N}$
and satisfies the Yang-Baxter eq.
(\ref{3.5.2}) in the form $(R_{ij}(x) = P_{ij} \R_{ij}(x))$
\be
\lb{qkz5}
R_{ij}(x) \, R_{ik}(xy) \, R_{jk}(y)  =  R_{jk}(y) \, R_{ik}(xy) \, R_{ij}(x) \; .
\ee
The unitarity condition $R_{ij}(x) R_{ji}(x^{-1}) = {\bf 1}$ is
also required (cf. (\ref{conuni}), (\ref{3.9.17a})) for this
 $R$-matrices. The constant matrix $D_i$ acts in
$i$-th vector space $V_i$ and obeys $R_{ij} D_i D_j = D_i D_j R_{ij}$.
Eqs. (\ref{qkz1}) with discrete connection (\ref{qkz4})
is called quantum Knizhnik - Zamolodchikov (q-KZ)
equations. It is convenient to rewrite the definition of the discrete connection (\ref{qkz4})
in the form of commutative matrices (\ref{qkz33})
 as follows:
\be
\lb{qkz44}
\begin{array}{c}
{\sf A}^{(j)}_{1 \dots N}(z_1 , \dots , z_N) =
T_{(j)}^{-1} \, A^{(j)}_{1 \dots N}(z_1 , \dots , z_N) \; = \;
 \hat{\hat{R}}_{j-1} \dots , \hat{\hat{R}}_{1} \,
\overline{\bf X} \,
\hat{\hat{R}}^{-1}_{N-1} \hat{\hat{R}}^{-1}_{N-2} \dots \hat{\hat{R}}^{-1}_{j} \; , \\ [0.2cm]
 \overline{\bf X} \; := \; T_{(1)}^{-1} \, D_1 \,
{\bf P}_{1,2} {\bf P}_{2,3} \dots {\bf P}_{N-1,N} \; ,
\end{array}
\ee
where
${\bf P}_{j,k} = P_{j,k} \cdot P_{z_j,z_{k}}$
  and $P_{z_j,z_k}$ is an operator which permutes
the spectral parameters $z_j$ and $z_{k}$:
$$
P_{z_j,z_k} \cdot f(z_1 , \dots , z_k , \dots , z_j , \dots , z_N)   =
f(z_1 , \dots , z_j , \dots , z_k , \dots , z_N) \cdot P_{z_j,z_k} \; ,
$$
such that  $\hat{\hat{R}}_j := {\bf P}_{j,j+1} R_{j,j+1}(z_j/z_{j+1})$ realize generators of
the braid group ${\cal B}_N$
(see eqs. (\ref{braidg}) in Sect. {\bf \ref{gabg1}}).
We note that operator $\overline{\bf X}$
 satisfies relations
 \be
 \lb{qkz17}
 \hat{\hat{R}}_{k+1} \, \overline{\bf X} =
 \overline{\bf X} \, \hat{\hat{R}}_{k} \;\;\;\;\;
 (k=1,...,N-2) \; , \;\;\;\;\;\;
 \hat{\hat{R}}_{1} \, \overline{\bf X}^2 =
 \overline{\bf X}^2 \, \hat{\hat{R}}_{N-1} \; ,
\ee
 and it can be considered as the image of an additional element
 which extends the group ${\cal B}_N$.
 \begin{proposition}\label{pro7kz}
(see \cite{FreRe}, \cite{TaVa}).
Discrete connection (\ref{qkz4}) is the flat
discrete connection (i.e. satisfies (\ref{qkz3})) and therefore
the system of equations (\ref{qkz1}) with
connection (\ref{qkz4}) is consistent.
\end{proposition}
{\bf Proof.}
Indeed, we have from (\ref{qkz3}) for $j >i$
$$
T_{(i)} T_{(j)} R_{j \,  j-1} \dots R_{j \,  1} T_{(j)}^{-1} D_j R_{Nj}^{-1}
\dots R^{-1}_{j+1\, j}
R_{i \, i-1} \dots R_{i \, 1} T_{(i)}^{-1} D_i R_{Ni}^{-1} \dots R^{-1}_{i+1 \, i} =
$$
$$
= T_{(i)} R_{i \, i-1} \dots R_{i \, 1} T_{(i)}^{-1} D_i T_{(j)} R_{Ni}^{-1}
\dots R^{-1}_{i+1 \, i}
R_{j \,  j-1} \dots R_{j \,  1} T_{(j)}^{-1} D_j R_{Nj}^{-1} \dots R^{-1}_{j+1\, j} \; .
$$
In the left hand side we obtain  for $j >i$
\be
\lb{qkz6}
T_{(i)} T_{(j)} R_{j \,  j-1} \dots R_{j \,  1} R_{i \, i-1} \dots R_{i \, 1} T_{(i)}^{-1}
 T_{(j)}^{-1} D_i D_j R_{N\, j}^{-1} \dots R^{-1}_{j+1\, j}
 R_{Ni}^{-1} \dots R^{-1}_{i+1 \, i} \; .
\ee
Then we use here identities for transfer-matrices
$$
 R_{j \,  j-1} \dots R_{j \,  1} (R_{i \, i-1} \dots R_{i \, 1}) =
(R_{i \, i-1} \dots R_{i \, 1}) R_{j \,  j-1}  \dots R_{j \,  i+1} R_{j \,  i-1}
\dots R_{j \,  1}  R_{ji}
$$
$$
 (R_{Nj}^{-1} \dots R^{-1}_{j+1\, j}) R_{N\, i}^{-1} \dots R^{-1}_{i+1 \, i} =
 R_{ji}^{-1}  R_{Ni}^{-1} \dots R^{-1}_{j+1 \, i} R_{j-1 \, i}^{-1} \dots R^{-1}_{i+1 \, i}
(R_{Nj}^{-1} \dots R^{-1}_{j+1\, j})
$$
and obvious relations
$[T_{(i)}^{-1} T_{(j)}^{-1} , \, R_{i\, j} ] = 0 = [D_i D_j , \,  R_{i\, j}]$. As a result we
obtain for l.h.s. (\ref{qkz6})
\be
\lb{qkz7}
\begin{array}{c}
T_{(i)} T_{(j)} (R_{i \, i-1} \dots R_{i \, 1}) R_{j \,  j-1}  \dots R_{j \,  i+1} R_{j \,  i-1}  \dots R_{j \,  1}
T_{(i)}^{-1}  T_{(j)}^{-1} D_i D_j \cdot \\ [0.2cm]
\cdot  R_{Ni}^{-1} \dots R^{-1}_{j+1 \, i} R_{j-1 \, i}^{-1} \dots R^{-1}_{i+1 \, i}
(R_{Nj}^{-1} \dots R^{-1}_{j+1\, j})
\end{array}
\ee
In the r.h.s. we use $[R_{i\, j} D_k] = 0 = [R_{i\, j} T_k]$
for $i,j \neq k$ and the identity
$$
R_{Ni}^{-1} \cdots R^{-1}_{i+1 \, i} R_{j \,  j-1} \cdots R_{j \,  1} =
 R_{j \,  j-1}  \cdots R_{j \,  i+1} R_{j \,  i-1}  \cdots R_{j \,  1}
 R_{Ni}^{-1} \cdots R^{-1}_{j+1 \, i} R_{j-1 \, i}^{-1}
 \cdots R^{-1}_{i+1 \, i} \, ,
$$
which gives for the r.h.s. just the same answer (\ref{qkz7})
as for the l.h.s. \hfill \qed

At the end of this subsection
we present the definition of the q-KZ equations
due to F.A.~Smirnov \cite{FAS}. Define the $R$-matrix and operator
$D_i$ as following
\be
\lb{qkz8}
\begin{array}{c}
\Psi^{1 \dots N \rangle}(z_1 , \dots , z_{i+1}, z_i , \dots  z_{N})  =
\R_{i \, i+1}(z_i/z_{i+1})
\Psi^{1 \dots N \rangle}(z_1 , \dots , z_{i}, z_{i+1} , \dots , z_{N})
\; , \\ [0.2cm]
\Psi^{1 \dots N \rangle}(p \, z_1 , z_{2}, \dots  , z_{N})  =
D_1 \, \Psi^{2 \dots N, 1 \rangle}(z_2 , z_3, \dots , z_{N}, z_1) \; .
\end{array}
\ee
One can explicitly show that eqs. (\ref{qkz8}) lead to the eqs.
(\ref{qkz1}), (\ref{qkz2}), (\ref{qkz4}). Indeed,
one can cyclically permute spectral parameters in
 $\Psi^{1 \dots N \rangle}$ by
means of the first eq. in (\ref{qkz8}) and
then use the second eq. in (\ref{qkz8}).
 The self-consistence of eqs. (\ref{qkz8})
can be checked directly. It also
follows from the self-consistence of the extended Zamolodchikov algebra
with generators $\{A_i(z_i), Q \}$ ($i=1, \dots ,N$):
$$
\R_{1 \, 2}(z_1/z_{2}) \, A^{1 \rangle}(z_1) \, A^{2 \rangle}(z_2) =
A^{1 \rangle}(z_2) \, A^{2 \rangle}(z_1) \; , \;\;\;
D_1 A^{1 \rangle}(z_1) \, Q  = Q \, A^{1 \rangle}(p \, z_1) \; ,
$$
and remark that eqs. (\ref{qkz8}) can be formally produced from the representation
$$
\Psi^{1 \dots N \rangle}(z_1 , \dots ,  z_{N}) =
{\rm Tr} \left( Q A^{1 \rangle}(z_1) \, A^{2 \rangle}(z_2) \dots A^{N \rangle}(z_N)  \right) \; .
$$

The semiclassical limit of the q-KZ equations (if we take the
the trigonometric $R$-matrices (\ref{3.5.5}) and (\ref{3.9.15a}) and consider their Yangian limits) gives \cite{FreRe}
the usual Knizhnik - Zamolodchikov equations.
Moreover, the
flat connections (\ref{qkz4}) (and their semiclassical limits)
can be related to Dunkl operators for
Calogero-Moser-Sutherland type models.

\vspace{0.1cm}

\noindent
{\bf Remark.} In \cite{IsKiTa}
(see also \cite{Cher3}, \cite{Kiril}) the rather general classification
for q-KZ flat connections was proposed. This classification is
based on the interpretation of q-KZ flat connections (\ref{qkz33}) as
images (in $R$-matrix representations) of commutative Jucys-Murphy
elements for affine braid groups defined by Coxeter
graphs. We discuss such braid groups below in subsection
{\bf \ref{gabg1}}.
In particular, the connections
(\ref{qkz44}) are images of the Jucys-Murphy
elements $\overline{J}_i$ for affine
braid group ${\cal B}_N(C^{(1)})$ (see Proposition 2.1 in
\cite{IsKiTa}).

\subsection{\bf \em Elliptic solutions of
the Yang-Baxter equation}
\setcounter{equation}0

In this subsection, we consider $Z_N \otimes Z_N$-symmetric
solutions of the Yang-Baxter equation
(\ref{3.5.9a}) (Ref. \cite{42}). The elements $R^{i_{1}i_{2}}_{j_{1}j_{2}}(\theta)$
of the corresponding
$R$-matrix will be expressed in terms of elliptic functions of the spectral parameter $\theta$.

We construct this solution explicitly, following the method of
Ref. \cite{42}. We consider two matrices $g$ and $h$ such that $g^{N} = h^{N} =1$:
\be
\lb{3.10.1}
g = \left(
\begin{array}{ccccc}
1 & 0 & 0 & \dots & 0 \\
0 & \omega & 0 & \dots & 0 \\
\vdots &  &  &  & \vdots \\
0 & 0 & 0 & \dots & \omega^{N-1}
\end{array}
\right) \; , \;\;
h = \left(
\begin{array}{ccccc}
0 & 1 & 0 & \dots & 0 \\
0 & 0 & 1 & \dots & 0 \\
\vdots &  &  &  &  \\
1 & 0 & 0 & \dots & 0
\end{array}
\right) \; ,
\ee
where $\omega = exp(2 \pi i/ N)$ and $h \, g = \omega \, g \, h$.
The matrices $g$ and $h$ are $Z_{N}$-graded generators of the algebra
$Mat(N)$, the graded basis for which can be chosen in the form
\be
\lb{3.10.2}
I_{\vec{\alpha}} = I_{\alpha_{1} \alpha_{2}} =
g^{\alpha_{1}} h^{\alpha_{2}} \; , \;\;
\alpha_{1,2} = 0,1, \dots , N-1 \; .
\ee
On the other hand, the matrices (\ref{3.10.2}) realize a projective representation
of the group $Z_{N} \otimes Z_{N}$:
$I_{\vec{\alpha}}I_{\vec{\beta}} = \omega^{\alpha_{2}\beta_{1}}
I_{\vec{\alpha} + \vec{\beta}}$ .
Any matrix $R_{12}(\theta) = R^{i_{1}i_{2}}_{j_{1}j_{2}}(\theta)$ can now be
written in the form
$$
R_{12}(\theta) = W_{\vec{\alpha},\vec{\beta}} (\theta) \,
I_{\vec{\alpha}} \otimes I_{\vec{\beta}} \; ,
$$
(the sum over $\alpha_i,\beta_j$ is assumed).
We consider the $Z_{N} \otimes Z_{N}$-invariant subset of such matrices:
\be
\lb{3.10.3}
R_{12}(\theta) = W_{\vec{\alpha}} (\theta)
I_{\vec{\alpha}} \otimes I^{-1}_{\vec{\alpha}} \; ,
\ee
where $I^{-1}_{\vec{\alpha}} =
h^{-\alpha_{2}} g^{- \alpha_{1}} = \omega^{\alpha_{1}\alpha_{2}}
I_{-\vec{\alpha}}$. The invariance of the matrices
(\ref{3.10.3}) is expressed by the relations
\be
\lb{3.10.4}
R_{12}(\theta) =
(I_{\vec{\gamma}} \otimes I_{\vec{\gamma}}) \, R_{12}(\theta) \,
(I_{\vec{\gamma}} \otimes I_{\vec{\gamma}})^{-1} \;\;\; \forall \vec{\gamma} \; ,
\ee
which obviously follow from the identity
$$
I_{\vec{\gamma}} \, I_{\vec{\alpha}} \, I^{-1}_{\vec{\gamma}} =
\omega^{<\alpha,\gamma>} I_{\vec{\alpha}} \; , \;\;
<\alpha, \, \gamma> = \alpha_{1} \gamma_{2} -
\alpha_{2} \gamma_{1} \; .
$$

It was noted in Ref. \cite{42} that the relations
\be
\lb{3.10.5}
\begin{array}{c}
R_{12}(\theta + 1) = g_{1}^{-1} \, R_{12}(\theta) \, g_{1} =
 g_{2} \, R_{12}(\theta) \, g_{2}^{-1} \; , \\ \\
R_{12}(\theta + \tau) =
exp(-i \pi \tau) exp(-2\pi i \theta)
h_{1}^{-1} \, R_{12}(\theta) \, h_{1} = \\ \\
exp(-i \pi \tau) exp(-2\pi i \theta)
h_{2} \, R_{12}(\theta) \, h_{2}^{-1} \; ,
\end{array}
\ee
\be
\lb{3.10.5a}
R_{12}(0) = I_{\vec{\alpha}} \otimes I^{-1}_{\vec{\alpha}}  \equiv  P_{12} \; ,
\ee
where $\tau$ is some complex parameter (period), are consistent with the
Yang-Baxter equation (\ref{3.5.9a}) and can be regarded as subsidiary conditions
to these equations (the last identity in (\ref{3.10.5a})
follows from
$(I_{\vec{\alpha}} \otimes I^{-1}_{\vec{\alpha}}) \,
I_{\vec{\beta}} \otimes I_{\vec{\gamma}} =
I_{\vec{\gamma}} \otimes I_{\vec{\beta}} \,
(I_{\vec{\alpha}} \otimes I^{-1}_{\vec{\alpha}})$).
Moreover, for the $Z_{N} \otimes Z_{N}$-invariant $R$-matrix
(\ref{3.10.3}) the conditions (\ref{3.10.5}), (\ref{3.10.5a}) determine the solution of the
Yang-Baxter equation uniquely. Indeed, substitution of (\ref{3.10.3}) in
(\ref{3.10.5}), (\ref{3.10.5a}) leads to the equations
\be
\lb{3.10.6}
\begin{array}{l}
W_{\vec{\alpha}}(\theta +1) = \omega^{\alpha_{2}} \,
W_{\vec{\alpha}} (\theta) \; , \\ \\
W_{\vec{\alpha}}(\theta +\tau) =
exp(-i \pi \tau) \, exp(-2\pi i \theta) \, \omega^{-\alpha_{1}} \,
W_{\vec{\alpha}} (\theta) \; , \;\; W_{\vec{\alpha}}(0) = 1 \; ,
\end{array}
\ee
the solution of which can be found by means of an expansion in a
Fourier series and has the form
\be
\lb{3.10.7}
W_{\vec{\alpha}} (\theta) =
\frac{\Theta_{\vec{\alpha}} (\theta + \eta)}
{\Theta_{\vec{\alpha}} (\eta)} \; , \;\;\;
\left( W_{\vec{\alpha} + \vec{\nu}} ( u ) =
W_{\vec{\alpha} + \vec{\nu} \, '} ( u ) = W_{\vec{\alpha}} ( u ) \right) \; ,
\ee
where $\vec{\nu} = (N,0)$,  $\vec{\nu} \, ' = (0,N)$,
\be
\lb{3.10.8}
\Theta_{\vec{\alpha}} ( u ) =
\sum_{m = -\infty}^{\infty} exp \left[ i \pi \tau
(m + \frac{\alpha_{2}}{N})^{2} + 2 \pi i
(m + \frac{\alpha_{2}}{N})(u + \frac{\alpha_{1}}{N} ) \right] \; ,
\ee
and we recall that $\alpha_{1,2} \in {\bf Z}_N$.
The parameter $\eta$ in (\ref{3.10.7}) is arbitrary. For $N=2$, the solution
(\ref{3.10.7}) is identical to the solution obtained by Baxter \cite{2}, \cite{Baxt1}
in connection with the investigation of the so-called eight-vertex lattice model.

Direct substitution of (\ref{3.10.3}) in the Yang-Baxter equation
(\ref{3.5.9a}) shows that the functions $W_{\vec{\alpha}}(\theta)$ must satisfy the relations
\be
\lb{3.10.9}
\sum_{\vec{\gamma}} W_{\vec{\gamma}} (\theta - \theta') \,
W_{\vec{\alpha} - \vec{\gamma}} (\theta) \, W_{\vec{\beta} + \vec{\gamma}} (\theta') \,
\left( \omega^{<\gamma , \beta >} -
\omega^{<\alpha - \gamma , \beta >} \right) = 0 \; .
\ee
As it was proved in \cite{43}, these relations hold when the functions
(\ref{3.10.7}) and (\ref{3.10.8}) are substituted.
We will see later that the identities (\ref{3.10.9}) are intimately related
to a version of the Yang-Baxter equations appeared in Interaction Round Face
models (see Sect. 5.3).

\section{GROUP ALGEBRA OF BRAID GROUP \\ AND ITS QUOTIENTS\label{gabg}}
\setcounter{equation}0
\setcounter{subsection}0

\subsection{\bf \em
Affine Braid Groups and Coxeter graphs \label{gabg1}}
\setcounter{equation}0

A braid group ${\cal B}_{M+1}$ is generated
 by elements $\sigma_i$ $(i=1, \dots M)$
subject to the relations:
\be
\lb{braidg}
\sigma_i \, \sigma_{i+1} \, \sigma_i =
\sigma_{i+1} \, \sigma_i \,  \sigma_{i+1} \; , \;\;\;
[ \sigma_{i},  \sigma_{j}] = 0 \;\;
{\rm for} \;\; |i-j| > 1 \; .
\ee
By definition the elements $\sigma_i$ are supposed to be invertible and
represented graphically as (cf. (\ref{figR}))

\unitlength=7mm
\begin{picture}(17,4)

\put(0.5,1.9){$\sigma_i \;\;  =$}

\put(3,3){$\bullet$}
\put(5,3){$\bullet$}

\put(5.7,3){$\dots$}
\put(7,3){$\bullet$}
\put(9,3){$\bullet$}
\put(10,3){$\dots$}
\put(11.5,3){$\bullet$}

\put(3,3.4){$1$}
\put(5,3.4){$2$}
\put(7,3.4){$i$}
\put(8.7,3.4){$i+1$}
\put(10.7,3.4){$M+1$}

\put(7.15,3.2){\line(1,-1){1}}
\put(8.0,2){\vector(-1,-1){0.8}}
\put(8.15,2.2){\vector(1,-1){1}}
\put(8.4,2.4){\line(1,1){0.8}}
\put(3.15,3.2){\vector(0,-1){2}}
\put(5.15,3.2){\vector(0,-1){2}}
\put(11.65,3.2){\vector(0,-1){2}}

\put(3,1){$\bullet$}
\put(5,1){$\bullet$}

\put(5.7,1){$\dots$}
\put(7,1){$\bullet$}
\put(9,1){$\bullet$}
\put(10,1){$\dots$}
\put(11.5,1){$\bullet$}
\end{picture}
\vspace{-1cm}
\be
\lb{figas}
{}
\ee


\newtheorem{def13}[def1]{Definition}
\begin{def13} \label{def13}
{\it An extension of the braid group ${\cal B}_{M}$ by one invertible generator
$\sigma_{M}$ subject to relations
\be
\lb{refl3}
\begin{array}{c}
\sigma_{M} \sigma_{M-1} \sigma_{M} = \sigma_{M-1} \sigma_{M} \sigma_{M-1} \; , \;\;\;
\sigma_{1} \sigma_{M} \sigma_{1} = \sigma_{M} \sigma_{1} \sigma_{M} \; , \\[0.2cm]
[ \sigma_M , \, \sigma_k ] = 0 \;\; (k = 2, \dots , M-2) \; ,
\end{array}
\ee
is called a periodic braid group
 $\overline{{\cal B}}_{M} \equiv {\cal B}_{M}(A^{(1)})$}.
\end{def13}


\newtheorem{def14}[def1]{Definition}
\begin{def14} \label{def14}
{\it An extension of the braid group ${\cal B}_{M+1}$
$(M\geq1)$ by one invertible generator $y_1$ which
satisfies the relations
\be
\lb{refl33}
y_1 \sigma_1 y_1 \sigma_1 = \sigma_1 y_1 \sigma_1 y_1 \; , \;\;\; [\sigma_i , \, y_1] = 0 \;\;
\forall i > 1 \; ,
\ee
is called the affine braid group
$\hat{\cal B}_{M+1} \equiv {\cal B}_{M+1}(C)$
of $C$ type. The extension of the group $\hat{\cal B}_{M+1}$
 by one more additional generator $y_{M+1}$ with
 defining relations
 \be
\lb{refl34}
y_{M+1} \sigma_M y_{M+1} \sigma_M =
\sigma_M y_{M+1} \sigma_M y_{M+1} \; , \;\;\;
[\sigma_i , \, y_{M+1}] = 0 \;\;
\forall i < M \; , \;\;\;\; [y_1 , \, y_{M+1}] = 0 \;
\ee
is the affine braid group of the $C^{(1)}$ type which is
denoted as ${\cal B}_{M+1}(C^{(1)})$.}
\end{def14}
It is clear that the affine group ${\cal B}_{M+1}(C)$
is the subgroup of the affine braid group
${\cal B}_{M+1}(C^{(1)})$, while the braid group
${\cal B}_{M+1}$
is the subgroup of ${\cal B}_{M+1}(C)$.
The defining relations (\ref{braidg}), (\ref{refl3}), (\ref{refl33})
 and (\ref{refl34}) (where we denote $y_1 = \sigma_0$ and
  $y_{M+1} = \sigma_{M+1}$) of the (affine) braid groups
can be written in the unified form as
\be
\lb{affbra}
\underbrace{\sigma_i \, \sigma_j \, \sigma_i \; \dots}_{m_{ij}
\;\; {\rm factors}} =
\underbrace{\sigma_j \, \sigma_i \, \sigma_j \; \dots}_{m_{ij}
\;\; {\rm factors}} \;\; ,
\ee
where $m_{ij} = m_{ji}$ are integers such that:
$m_{ii}=1$, $m_{ij} \geq 2$
for $i \neq j$. The set of data given by the matrix $||m_{ij}||$
is conveniently represented as the {\em Coxeter graph} with $M$
(or $M+1$, or $M+2$) nodes
associated with generators $\sigma_i$ and the nodes $i$ and $j$ are
 connected by $(m_{ij}-2)$ lines
if $m_{ij} =2,3,4$ and by 3 lines if $m_{ij}=6$.
 Thus, the Coxeter graph for the braid group relations
(\ref{braidg}) is the $A$-type graph:

\unitlength=8mm
\begin{picture}(17,2)(-2,0)
\put(1.9,1.1){\circle{0.2}}
\put(1.7,0.5){$\sigma_1$}
\put(2,1.1){\line(1,0){1}}
\put(3.1,1.1){\circle{0.2}}
\put(2.9,0.5){$\sigma_2$}
\put(3.2,1.1){\line(1,0){1}}
\put(4.3,1.1){\circle{0.2}}
\put(4.1,0.5){$\sigma_3$}
\put(4.4,1.1){\line(1,0){1}}
\put(6,1.1){. . . . . . . .}
\put(9,1.1){\line(1,0){1}}
\put(10.1,1.1){\circle{0.2}}
\put(9.7,0.5){. . .}
\put(10.2,1.1){\line(1,0){1}}
\put(11.3,1.1){\circle{0.2}}
\put(11,0.5){$\sigma_M$}
\end{picture}

\vspace{-2cm}
 \be
 \lb{typeA}
 {}
 \ee

\vspace{0.5cm}
\noindent
and, for the affine braid group relations
(\ref{refl3}), (\ref{refl33}) and (\ref{refl34}),
the Coxeter graphs are respectively

\unitlength=8mm
\begin{picture}(17,3)(-2,0)

\put(-0.5,1.5){$A^{(1)} \; = $}

\put(1.9,1.1){\circle{0.2}}
\put(1.7,0.5){$\sigma_1$}
\put(2,1.1){\line(1,0){1}}
\put(3.1,1.1){\circle{0.2}}
\put(2.9,0.5){$\sigma_2$}
\put(3.2,1.1){\line(1,0){1}}
\put(4.3,1.1){\circle{0.2}}
\put(4.1,0.5){$\sigma_3$}
\put(4.4,1.1){\line(1,0){1}}
\put(6,1.1){. . . . . . . . . .}
\put(10,1.1){\line(1,0){1}}
\put(11.1,1.1){\circle{0.2}}
\put(10.7,0.5){. . .}
\put(11.2,1.1){\line(1,0){1}}
\put(12.3,1.1){\circle{0.2}}
\put(12,0.5){$\sigma_{_{M-1}}$}

\put(2,1.1){\line(4,1){5}}
\put(7.25,2.4){\line(4,-1){5}}
\put(7.18,2.4){\circle{0.2}}
\put(7.5,2.6){$\sigma_{_M}$}
\end{picture}

\vspace{-2cm}
 \be
 \lb{A1}
 {}
 \ee

\vspace{0.7cm}

   \unitlength=5mm
\begin{picture}(17,3.5)(-6,0)

\put(-2,1.7){$C \; = $}

\put(2,2){\circle{0.4}}
\put(1.5,2.5){$\sigma_0$}
\put(4.5,2){\circle{0.4}}
\put(4.5,2.5){$\sigma_1$}
\put(2.2,1.9){\line(1,0){2.1}}
\put(2.2,2.1){\line(1,0){2.1}}
\put(4.7,2){\line(1,0){2}}
\put(7,2){$. \; . \; . \; . \; . \; . \; .$}

\put(10.5,2){\circle{0.4}}
\put(9.8,2.5){$\sigma_{_{M-1}}$}
\put(10.7,2){\line(1,0){2.1}}
\put(13,2){\circle{0.4}}
\put(12.7,2.5){$\sigma_{_M}$}

\end{picture}

\vspace{-2cm}
 \be
 \lb{CCC}
 {}
 \ee

\vspace{0.7cm}

 \unitlength=5mm
\begin{picture}(17,3.5)(-4,0)

\put(-2,1.7){$C^{(1)} \; = $}

\put(2,2){\circle{0.4}}
\put(1.5,2.5){$\sigma_0$}
\put(4.5,2){\circle{0.4}}
\put(4.5,2.5){$\sigma_1$}
\put(2.2,1.9){\line(1,0){2.1}}
\put(2.2,2.1){\line(1,0){2.1}}
\put(4.7,2){\line(1,0){2}}
\put(7,2){$. \; . \; . \; . \; . \; . \; .$}

\put(10.5,2){\circle{0.4}}
\put(9.8,2.6){$\sigma_{_{M-1}}$}
\put(10.7,2){\line(1,0){2}}
\put(13,2){\circle{0.4}}
\put(12.7,2.6){$\sigma_{_M}$}
\put(15.5,2){\circle{0.4}}
\put(15.5,2.6){$\sigma_{_{M+1}}$}
\put(13.2,1.9){\line(1,0){2.1}}
\put(13.2,2.1){\line(1,0){2.1}}

\end{picture}

\vspace{-1.5cm}

 \be
\lb{C1}
{}
\ee


 \noindent
 In the same way one can also define the affine braid groups of $B^{(1)}$
 and $D^{(1)}$ types (see, e.g., \cite{IsKiTa} and references therein).

Consider the affine braid group $\hat{\cal B}_{M+1} \equiv
{\cal B}_{M+1}(C)$.
The elements $\{ y_i\}$  ($i=1, \dots M+1$) defined by:
\be
\lb{combr}
y_1 \; , \;\;\; y_2 = \sigma_1 y_1 \sigma_1 \; , \;\;\;
y_3 = \sigma_2 \sigma_1 y_1 \sigma_1 \, \sigma_2 \; ,
\; \dots \; , \; y_{i+1} = \sigma_i y_i \sigma_i  \; ,
\ee
 are called {\it Jucys-Murphy elements} and
generate an Abelian subgroup in $\hat{\cal B}_{M+1}$. For $y_1 = 1$, the Jucys-Murphy elements (\ref{combr})
generate an Abelian subgroup in the braid group ${\cal B}_{M+1}$. Note that condition $y_n y_{n+1} = y_{n+1} y_n$
 is equivalent to the reflection equation for $y_n$:
\be
\lb{refla}
y_n \, \sigma_n \, y_n \, \sigma_n =
\sigma_n \, y_n \, \sigma_n \, y_n \; .
\ee
Then, we have $y_n \, y_{n+1} \, \sigma_n =
\sigma_n \, y_n \, y_{n+1}$, and the element $Z = y_1 y_2 \cdots y_{M+1}$
is obviously central in ${\cal B}_{M+1}$.


\begin{proposition}\label{prop8}
{\it The product of $m$ elements of $\hat{B}_{M+1}$: $y^{(m)}_{k+1}:= y_{k+1}
y_{k+2} \dots y_{k+m}$
$(k+m<M+1)$ satisfies the following relations
\be
\lb{afh1}
y^{(m)}_{k+1} = U_{(k,m)} \; y^{(m)}_{1} \; U_{(m,k)} \; ,
\ee
where (cf. (\ref{04.4a}))
\be
\lb{afh2}
U_{(k,m)} = \sigma_{(k \to m+k-1)} \dots \sigma_{(2 \to m+1)} \, \sigma_{(1 \to m)} \equiv
\sigma_{(k \leftarrow 1)} \sigma_{(k+1 \leftarrow 2)} \sigma_{(k+m-1 \leftarrow m)} \; ,
\ee
and $(k \leq n)$
$$
\sigma_{(k \to n)} = \sigma_k \sigma_{k+1} \dots \sigma_n \; , \;\;\;
\sigma_{(n \leftarrow k)} = \sigma_n \dots \sigma_{k+1} \sigma_k \; .
$$
}
\end{proposition}

\vspace{0.1cm}
\noindent
{\bf Proof.}
First of all we show that
\be
\lb{afh3}
y^{(m)}_{k+1} = \sigma_{(k \to k+m-1)} \, y^{(m)}_{k} \,
\sigma_{(k+m-1 \leftarrow k)} \; .
\ee
This identity is proved by induction. For $m=1$ we obviously have
$y_{k+1} = \sigma_k y_k \sigma_k$. Let (\ref{afh3}) be correct for some $m$.
Then, for $y^{(m+1)}_{k+1}$ we have
$$
y^{(m+1)}_{k+1} = y^{(m)}_{k+1} \, y_{k+m+1} = \sigma_{(k \to k+m-1)} \, y^{(m)}_{k} \,
\sigma_{(k+m-1 \leftarrow k)} \, \sigma_{k+m} y_{k+m} \sigma_{k+m} =
$$
$$
= \sigma_{(k \to k+m-1)}  \, y^{(m)}_{k} \,
 \sigma_{k+m} y_{k+m} \sigma_{k+m} \, \sigma_{(k+m-1 \leftarrow k)}  =
$$
$$
= \sigma_{(k \to k+m-1)}
 \sigma_{k+m} \, (y^{(m)}_{k} y_{k+m}) \, \sigma_{k+m} \sigma_{(k+m-1 \leftarrow k)} \; ,
$$
which coincides with (\ref{afh3}) for $m \rightarrow m+1$.
Applying (\ref{afh3}) several times we deduce (\ref{afh1}). \hfill \qed

One can graphically represent elements $U_{(k,m)}$ (\ref{afh2})
(by means of the rules (\ref{figas})) in the following form
(cf. (\ref{figas1b}))

\unitlength=10mm
\begin{picture}(17,4)

\put(0.5,1.9){$U_{(k,m)} \;\;  =$}

\put(3,3){$\bullet$}
\put(5.2,3){$\bullet$}

\put(4.1,3){$\bullet$}
\put(4.15,3.15){\vector(2,-1){4}}
\put(3.5,3.1){$\dots$}
\put(4.5,3.1){$\dots$}
\put(5.7,1.8){$\ddots$}
\put(6.2,2.1){$\ddots$}

\put(7.7,3){$\dots$}
\put(7,3){$\bullet$}
\put(9.2,3){$\bullet$}
\put(10.5,3){$\dots$}
\put(10,3){$\bullet$}
\put(11.55,3){$\bullet$}

\put(3,3.4){$_1$}

\put(5.3,3.4){$_k$}
\put(6.7,3.4){$_{k+1}$}
\put(8.7,3.4){$_{k+m}$}
\put(9.8,3.4){$_{k+m+1}$}
\put(11.4,3.4){$_{M+1}$}

\put(7.15,3.15){\line(-2,-1){0.7}}
\put(9.2,3.1){\line(-2,-1){1.7}}


\put(6.15,2.65){\line(-2,-1){0.4}}
\put(7.2,2.1){\line(-2,-1){0.4}}

\put(4.8,2){\vector(-2,-1){1.6}}
\put(6.13,1.6){\vector(-2,-1){0.85}}

\put(3.2,3.15){\vector(2,-1){4}}
\put(5.3,3.15){\vector(2,-1){4}}
\put(11.65,3.15){\vector(0,-1){2}}
\put(10.1,3.15){\vector(0,-1){2}}

\put(3,1){$\bullet$}
\put(5.2,1){$\bullet$}
\put(3.9,1){$\dots$}
\put(7,1){$\bullet$}
 \put(9.2,1){$\bullet$}

\put(8,1){$\bullet$}
\put(7.4,1.1){$\dots$}
\put(8.5,1.1){$\dots$}

\put(10.5,1){$\dots$}
\put(10,1){$\bullet$}
\put(11.55,1){$\bullet$}
\end{picture}
\vspace{-1.5cm}
\be
\lb{figas1}
{}
\ee

\noindent
From this representation it becomes clear that the following braid
relations hold
\be
\lb{brU}
U_{(k,m)} \, (T_m U_{(k,n)}) \, U_{(m,n)} = (T_k U_{(m,n)}) \, U_{(k,n)} \, (T_n U_{(k,m)}) \; ,
\ee
where we have introduced jump operations $T_m$: $\sigma_i \to \sigma_{i+m}$. One can check
relations
(\ref{brU}) by direct calculations.

\noindent
{\bf Remark.} The sets of commutative Jucys-Murphy elements
for all affine braid groups of the $A^{(1)}$,
 $B^{(1)}$, $C^{(1)}$ and $D^{(1)}$ types
were constructed in \cite{IsKiTa}.
These elements were realized in $R$-matrix representations
of Birman-Murakami-Wenzl algebras and then used  in \cite{IsKiTa}
for the formulation of the special $q$-KZ equations.

\subsection{\bf \em  Group Algebra of the Braid Group $B_{M+1}$
and shuffle elements \label{gabg2}}
\setcounter{equation}0

We denote the group algebra of the braid group $B_{M+1}$ over
complex numbers as $\mathbb{C}[B_{M+1}]$.
Consider the elements
$\Sigma_{m \rightarrow n} \in \mathbb{C}[B_{M+1}]$, for
 $n=m,\, m+1,..., \, M+1$, that are defined inductively
\be
\label{anti}
\Sigma_{m \to n} =  f_{m \to n} \, \Sigma_{m \to n-1} =
 f_{m \to n} \, f_{m \to n-1} \cdots  f_{m \to m+1} \,  f_{m \to m}  \; ,
\ee
where the subscript $m \to n$ denotes
the set of indices $(m,m+1, \dots , n)$,
$\Sigma_{m\to m}  = 1$ and
\be
\label{fbf}
 \begin{array}{c}
f_{k \to k} = 1 \, , \;\;\;\;\;
f_{k \to n} =
 1 + \sigma_{n-1} + \sigma_{n-2} \,
\sigma_{n-1} + \dots
+  \sigma_{k} \, \sigma_{k+ 1} \cdots \sigma_{n-1} \; = \\ [0.2cm]
= 1 + f_{k \to n-1}\, \sigma_{n-1}  \; , \;\;\;\;\;\;\;\; k < n \; .
\end{array}
\ee
 Note that,
by means of braid relations (\ref{braidg}), we derive
the mirror set of expressions
for the elements
$\Sigma_{k \to n}$ (\ref{anti})
\be
\label{danti}
\Sigma_{k \to n} =   \Sigma_{k \to n-1} \, \overline{f}_{k \to n} =
\overline{f}_{k \to k} \, \overline{f}_{k \to k+1} \cdots
\overline{f}_{k \to n-1} \, \overline{f}_{k \to n} \; ,
\ee
where
\be
\label{dfbf}
 \begin{array}{c}
 \overline{f}_{k \to k} = 1 \; , \;\;\;\;\;
\overline{f}_{k \to n} =
 1 + \sigma_{n-1} + \sigma_{n-1} \,
\sigma_{n-2} + \dots
+  \sigma_{n-1} \dots \sigma_{k+ 1} \, \sigma_{k} = \\ [0.2cm]
= 1 + \sigma_{n-1} \overline{f}_{k \to n-1} \;  , \;\;\;\;\;  k < n \; .
\end{array}
\ee

 The elements $\Sigma_{m+1 \to n}$ play the role of  symmetrizers
  in $\mathbb{C}[B_{M+1}]$ and, in view
  of the projection $\sigma_i \to 1$ for
  (\ref{anti}) -- (\ref{dfbf}), they are
  algebraic analogs of the factorials $(n-m)!$.
The important properties of the elements
 $f_{k \to n} \in \mathbb{C}[B_{M+1}]$ are
\be
\lb{ffxa}
f_{1 \to n} f_{1 \to n-1}
\dots f_{1 \to m+1} =   \Sha^{(m,n-m)}_{1\to n} \,
\Sigma_{m+1 \rightarrow n} \; , \;\;\;  (0 \leq m < n) \; ,
\ee
where we introduce elements
 $\Sha^{(m,n-m)}_{1\to n}\in \mathbb{C}[B_{M+1}]$
with initial conditions
 $\Sha^{(n,0)}_{1 \to n} =  \Sha^{(0,n)}_{1 \to n} = 1$,
$\Sha^{(m-1,1)}_{1 \to m}  = f_{1 \to m}$.
The identities (\ref{ffxa}) and definition (\ref{anti}),
written in the form
 $$
\Sigma_{1 \to n} = f_{1 \to n} f_{1 \to n-1}
\dots f_{1 \to m+1} \Sigma_{1 \to m} \; ,
 $$
 lead to the right factorization formula
\be
\lb{rfactor}
\Sigma_{1 \to n} = \Sha^{(m,n-m)}_{1 \to n}
 \, \Sigma_{1 \to m} \, \Sigma_{m +1 \rightarrow n} \; .
\ee
Thus, taking into account the interpretation of
 $\Sigma_{m \rightarrow n}$ as factorials, one can consider
 $\Sha^{(m,n-m)}_{1 \to n}$ as an algebraic analog of the
 binomial coefficient.
 \begin{proposition}\label{prosha}
The elements $\Sha^{(m,n-m)}_{1 \to n}$ are defined inductively
by using the recurrent relations \cite{IsBonn}, \cite{IsOg09}
(braid analogs of the Pascal rule)
 \be
 \lb{pasc}
 \Sha^{(m,n+1-m)}_{1 \to n+1} = \Sha^{(m,n-m)}_{1 \to n}  +
 \Sha^{(m-1,n+1-m)}_{1 \to n} \,
 \sigma_{n} \, \sigma_{n-1} \dots \sigma_{m} \; .
 \ee
 \end{proposition}
 \noindent
{\bf Proof.} We denote $f_k:=f_{1\to k}$
and consider formula (\ref{ffxa}) for $n \to (n+1)$:
 \be
 \lb{pasc1}
 \begin{array}{c}
\Sha^{(m,n+1-m)}_{1\to n+1} \,
\Sigma_{m+1 \rightarrow n+1} =
f_{n+1} f_{n} \cdots f_{m+1} =
 (1 + f_{n}\sigma_n) f_{n} f_{n-1} \cdots f_{m+1} =  \\ [0.2cm]
 = f_{n} f_{n-1} \cdots f_{m+1} +
 f_{n} \sigma_{n} (1+f_{n-1} \sigma_{n-1})
 f_{n-1}f_{n-2}\cdots f_{m+1} =  \\ [0.2cm]
 = f_{n} \cdots f_{m+1}(1 + \sigma_{n}) \; + \;
 f_{n} f_{n-1} \sigma_{n} \sigma_{n-1}
 (1+f_{n-2}\sigma_{n-2}) f_{n-2}\cdots f_{m+1} = \dots = \\ [0.2cm]
 = f_{n} \cdots f_{m+1}(1 + \sigma_{n} + \sigma_{n} \sigma_{n-1} +
 \dots + \sigma_{n} \cdots \sigma_{m+1}) +
 f_{n} \cdots f_{m} \, \sigma_{n} \cdots \sigma_{m} = \\ [0.2cm]
 =  \Sha^{(m,n-m)}_{1\to n} \, \Sigma_{m+1 \rightarrow n} \;
 \overline{f}_{m+1 \to n+1} +  \Sha^{(m-1,n-m+1)}_{1\to n} \,
 \Sigma_{m \rightarrow n}\, \sigma_{n} \sigma_{n-1}
  \cdots \sigma_{m}  = \\ [0.2cm]
 =  \Sha^{(m,n-m)}_{1\to n} \, \Sigma_{m+1 \rightarrow n+1} \;
  +  \Sha^{(m-1,n-m+1)}_{1\to n} \, \sigma_{n}
  \sigma_{n-1} \cdots \sigma_{m} \,
 \Sigma_{m+1 \rightarrow n+1} \; ,
\end{array}
 \ee
 where we used (\ref{fbf}) -- (\ref{dfbf}).
 After dividing both sides of (\ref{pasc1})
 by  $\Sigma_{m+1 \rightarrow n+1}$ from the right, we obtain (\ref{pasc}).
\hfill \qed

\vspace{0.1cm}

In particular, eq. (\ref{pasc}) is written  for $m=(n-1)$
as $\Sha^{(n-1,2)}_{1 \to n+1} = (f_{1 \to n} +
\Sha^{(n-2,2)}_{1 \to n} \, \sigma_{n} \, \sigma_{n-1})$,
$(n \geq 2)$ which gives
\be
\lb{xm}
\begin{array}{c}
 \Sha^{(m-1,2)}_{1 \to m + 1} = f_{1 \to m} +
f_{1 \to m-1} \,
(\sigma_{m} \, \sigma_{m-1}) +
f_{1 \to m-2} \, (\sigma_{m-1} \,
\sigma_{m-2}) \,
(\sigma_{m} \, \sigma_{m-1}) \\ [0.3cm]
+ \dots + f_{1 \to 2} \, (\sigma_{3} \,
\sigma_{2}) \cdots
(\sigma_{m} \, \sigma_{m-1})
+ (\sigma_{2} \, \sigma_{1}) \cdots (\sigma_{m} \,
\sigma_{m-1}) \; .
\end{array}
\ee
 The next relation is:
$\Sha^{(n-2,3)}_{1 \to n+1} =  \Sha^{(n-2,2)}_{1 \to n}   +
\Sha^{(n-3,3)}_{1 \to n} \, \sigma_{n} \,
\sigma_{n-1} \, \sigma_{n-2}$
for $(n \geq 3)$, etc.

Note that $\Sha^{(m,n-m)}_{1 \to n}$
 are sums over the braid group elements which can be considered
as quantum analogs of $(m, \, n-m)$ shuffles of
 two piles with $m$ and $(n-m)$ cards if we read all monomials
in $\Sha^{(m,n-m)}_{1 \to n}$ from right to left
(the standard shuffles are obtained by projection $\sigma_i \to s_i$,
where $s_i$ are generators of the
 symmetric group ${\cal S}_{M+1}$). As it follows from (\ref{ffxa})
the elements $f_{1 \to m} = \Sha^{(m-1,1)}_{1 \to m}$ are the sums of $(m-1,1)$ shuffles. One can
use the operators $\Sigma_{1 \to n}$ (\ref{anti})
and identities (\ref{rfactor})
for the definition of the associative products that are
analogs of the wedge products proposed by S.Woronowicz in the theory of differential calculus on
quantum groups \cite{Wor}. In view of (\ref{rfactor}) these products are related to the quantum
shuffle products
(about quantum shuffles and corresponding products see \cite{Rosso2},
\cite{IsOg09}). The associativity of these
products is provided by the identities:
\be
\lb{fi}
\Sha^{(n-m,m)}_{1 \to n} \, \Sha^{(k,m-k)}_{1 \to m} =
\Sha^{(k,n-k)}_{1 \to n} \, \Sha^{(m-k,n-m)}_{k+1 \to n} \;\;\; (k < m < n) \; ,
\ee
which are the consistence conditions for the definition
of a 3-pile shuffles $(k,m-k,n-m)$:
$$
\Sigma_{1 \to n} = \Sha^{(k,m-k,n-m)}_{1 \to n}
 \, \Sigma_{1 \to k} \, \Sigma_{k+1 \to m} \, \Sigma_{m +1 \rightarrow n} \; ,
$$
Going further one can introduce $m$-pile shuffles
$\Sha^{(n_1,n_2, \dots ,n_m)}_{1 \to n}$ of the
pack of $n$ cards $(n = n_1 + n_2 + \dots + n_m)$.
 Then, we observe that the ``symmetrizer''
$\Sigma_{1 \to n}$ (\ref{anti}) is nothing else but the
 $n$-pile shuffle $\Sha^{(1,1, \dots,1)}_{1 \to n}$.

By means of the mirror mapping (when we write generators
 $\sigma_k \in {\cal B}_{M+1}$ in all monomials in
 opposite order) we obtain from (\ref{rfactor}) a left
 factorization formula
\be
\lb{lfactor}
\Sigma_{1 \to n} =  \Sigma_{1 \to m} \, \Sigma_{m +1 \rightarrow n} \,
\overline{\Sha}^{(m,n-m)}_{1 \to n} \; ,
\ee
where the elements $\overline{\Sha}^{(m,n-m)}_{1 \to n}$ are defined by recurrence relations \cite{IsBonn}, \cite{IsOg09} (cf. (\ref{pasc}))
$$
\overline{\Sha}^{(m,n-m+1)}_{1 \to n+1} =
\overline{\Sha}^{(m,n-m)}_{1 \to n} + \sigma_m \dots \sigma_n \,
{\bf \overline{\Sha}}^{(m-1,n-m+1)}_{1 \to n}
$$
with initial conditions $\overline{\Sha}^{(n,0)}_{1 \to n}  = \overline{\Sha}^{(0,n)}_{1 \to n}  =1$,
$\overline{\Sha}^{(n-1,1)}_{1 \to n} = \overline{f}_{1 \to n}$,
and $\overline{\Sha}^{(m,n-m)}_{1 \to n}$ is a sum over
$(m,n-m)$ quantum shuffles (if we read all monomials
 from left to right). The mirror
analogs of the factorization identities (\ref{fi}) are also hold
$$
\overline{\Sha}^{(k,m-k)}_{1 \to m} \, \overline{\Sha}^{(m,n-m)}_{1 \to n} =
 \overline{\Sha}^{(m-k,n-m)}_{k+1 \to n} \,  \overline{\Sha}^{(k,n-k)}_{1 \to n} \; .
$$

\subsection{\bf \em  $A$-Type Hecke algebra
$H_{M+1}(q)$\label{AHalg}}
\setcounter{equation}0

\subsubsection{Jucys-Murphy elements, symmetrizers
and antisymmetrizers in $H_{M+1}$\label{jmel}}

$A$-Type Hecke algebra $H_{M+1}(q)$ (see e.g. \cite{Jon1} and Refs. therein) is
a quotient of the braid group algebra (\ref{braidg}) by the additional
relation
\be
\lb{ahecke}
\sigma^2_i - 1 = \lambda \,  \sigma_i \; , \;\; (i = 1, \dots , M) \; .
\ee
Here $\lambda = (q -q^{-1})$
and $q \in \mathbb{C} \backslash \{0, \pm 1 \}$ is a deformation parameter.
Note that algebras $H_{M+1}(q)$ and $H_{M+1}(-q^{-1})$
 are isomorphic to each other:
$H_{M+1}(q) \simeq H_{M+1}(-q^{-1})$.
The group algebra of ${\cal B}_{M+1}$ (\ref{braidg})
has an infinite dimension while its quotient $H_{M+1}(q)$
is finite dimensional.
It can be shown (see e.g. \cite{W1}) that $H_{M+1}(q)$
is linearly spanned by
$(M+1)!$ monomials appeared in the expansion of
$\Sigma_{1 \to M+1}$ (\ref{anti}) (or in the expansion
of (\ref{danti})).

The $A$-type Hecke algebra is a special case of
a general affine Hecke algebra. The affine Hecke algebra
is the quotient (by additional constraint (\ref{ahecke}))
of the affine braid groups with generators
 $\{ \sigma_i \}$ subject to general relations (\ref{affbra}).
 As it was shown in Sect. {\bf \ref{gabg1}},
 the Coxeter graph for the
 braid group relations (\ref{braidg}) is the $A$-type graph
 (\ref{typeA}).
That is why the Hecke algebra with defining relations (\ref{braidg}) and (\ref{ahecke})
is called $A$-type Hecke algebra. The $A$-type Hecke
algebra $H_{M+1}(q)$ is a semisimple algebra.

An essential information about a finite
dimensional semisimple associative algebras ${\cal A}$ is contained
(see, e.g., \cite{Weyl}; see also
 subsection 4.5 in \cite{IsRub2}, and references therein) in the
structure of its regular bimodule, which is decomposed into direct sums:
$$
{\cal A}=\bigoplus_{\alpha=1}^s {\cal A} \cdot e_\alpha \; , \;\;\;\;
{\cal A}=\bigoplus_{\alpha=1}^s e_\alpha \cdot {\cal A} \; ,
$$
of left and right submodules (ideals), respectively ({\it left and
right Peirce decompositions}). Here the elements
$e_\alpha \in {\cal A}$ $(\alpha=1,\dots , s)$ are mutually orthogonal
idempotents resolving the identity operator $1$:
 \be
 \lb{orind}
e_\alpha \cdot e_\beta = \delta_{\alpha \beta} e_\alpha \; , \;\;\;\;\;
1=\sum_{\alpha=1}^s e_\alpha \; .
 \ee
 Making use of  left and right Peirce
 decompositions simultaneously, we have {\it two-sided
 Peirce decomposition}
 \be\lb{resun19}
 {\cal A} = \bigoplus _{\alpha,\beta=1}^s
 e_\alpha \, {\cal A} \, e_\beta =
 \bigoplus _{\alpha,\beta=1}^s \, {\cal A}_{\alpha,\beta}
 \; , \;\;\;\;\; {\cal A}_{\alpha,\beta} =e_\alpha \,
 {\cal A} \, e_\beta \; .
 \ee
 Here the linear spaces  ${\cal A}_{\alpha,\beta}$ are, generally
 speaking, neither left nor right ideals in
 ${\cal A}$. Instead, products of elements of
 ${\cal A}_{\alpha,\beta}$ obey the relations:
 $ {\cal A}_{\alpha,\beta} \cdot {\cal A}_{\gamma,\kappa} =
 \delta_{\beta,\gamma} \; {\cal A}_{\alpha,\kappa}$, which
 resemble relations for matrix units.

The number $s$ depends on the choice of the type of idempotents
in ${\cal A}$.
There are two important types of the idempotents
in ${\cal A}$ and correspondingly two
decompositions of the identity operator: \\
(1) {\it Primitive idempotents}. An idempotent
$e_\alpha$ is primitive if it cannot be further resolved into a
sum of nontrivial mutually orthogonal idempotents. \\
(2) {\it Primitive central idempotents}. An idempotent
$e'_A$ $(A=1,...,s')$ is primitive central if it is central element
in ${\cal A}$ and primitive in
the class of central idempotents. \\
One can
expand any central idempotent $e_A$ in primitive idempotents
  $\{e_\alpha \}$: $e_A  = \sum_{\alpha \in A} e_\alpha$,
 where $A$ is a subset of indices from the set $\{1,2,\dots,s\}$;
 i.e., central orthogonal idempotents
 $\{e_A \}$ are conveniently labeled
 by non-intersecting subsets $A \subset \{1,2,\dots,s\}$,
 which cover the entire set of indices $\{1,2,\dots,s\}$.

 Let
 $A_i$ $(i=1,\dots,s')$ be non-intersecting subsets in
 $\{1,2,\dots,s\}$ which cover the entire set and define
 central idempotents $e_{A_i}$.
 Let $\alpha \in A_i$, $\beta \in A_j$ and $i \neq j$, then, in view of
 orthogonality $e_{A_i} \cdot e_{A_j} =0$,
 for any element
 $a \in {\cal A}$ we have
 $a_{\alpha,\beta} = e_\alpha \cdot a \cdot e_\beta = 0$.
This tells us that
the two-sided Peirce decomposition
(\ref{resun19}) of semisimple algebra ${\cal A}$
does not contain terms
 ${\cal A}_{\alpha,\beta}$, if
 $\alpha \in A_i$, $\beta \in A_j$ and $i \neq j$, so that we have
  \be\lb{res19p}
 {\cal A} = \bigoplus_{i=1}^{s'} e'_{A_i} \cdot {\cal A} \cdot e'_{A_i}
  = \bigoplus_{i=1}^{s'}
 \bigoplus_{\alpha,\beta \in A_i}
 e_\alpha \, {\cal A} \, e_\beta = \bigoplus_{i=1}^{s'}
 \bigoplus_{\alpha,\beta \in A_i} \, {\cal A}_{\alpha,\beta} \; ,
 \ee
 where, again, $s'$ is the number of primitive
 central idempotents $e_{A_i}$.
 Thus, the regular bimodule of the semisimple algebra ${\cal A}$
decomposes into direct sums
of irreducible sub-bimodules (two-sided ideals)
${\cal A}=\bigoplus_{i=1}^{s'} {\cal A} \cdot e'_{A_i}
=\bigoplus_{i=1}^{s'} e'_{A_i} \cdot {\cal A}$
with respect to the central idempotents $e'_A$.
For semisimple algebras ${\cal A}$ the subspaces
${\cal A}_{\alpha,\beta}$ in (\ref{res19p}) are one-dimensional
and for any $a \in {\cal A}$ we have
$e_\alpha \cdot a \cdot e_\beta = c(a) \; e_{\alpha \beta}$,
where $c(a)$ are constants and basis elements
$e_{\alpha \beta} \in {\cal A}_{\alpha,\beta}$ are normalized
such that $e_{\alpha \beta} \cdot e_{\gamma \delta} =
\delta_{\beta \gamma} \; e_{\alpha \delta}$. In view of these relations
the elements $e_{\alpha \beta} \in {\cal A}$ are called matrix units.
The diagonal matrix units coincide with the
primitive idempotents: $e_{\alpha\alpha} = e_{\alpha}$.

\vspace{0.2cm}

Now we return back to the consideration of
the Hecke algebra $H_{M+1}$
 (here and below we omit the parameter $q$ in the notation
$H_{M+1}(q)$). First of all we construct two
special primitive idempotents in the Hecke algebra $H_{M+1}$
which correspond to the symmetrizers and antisymmetrizers.
For this purpose we consider two substitutions $\sigma_i \to q \sigma_i$,
$\sigma_i \to -q^{-1} \sigma_i$ for the braid group algebra element
$\Sigma_{1 \to n}$ (\ref{anti}). As a result, for the algebra $H_{M+1}$
we obtain two sequences of operators $S_{1 \to n}$ and $A_{1 \to n}$
$(n =1, \dots M+1)$
\be
\lb{santis}
S_{1 \to n} :=
a^{-}_n \, \Sigma_{1 \to n}(q \, \sigma_i) \; , \;\;\;
A_{1 \to n} := a^{+}_n \, \Sigma_{1 \to n}(-q^{-1} \, \sigma_i) \;
\ee
$$
\left( a^{\mp}_n = \frac{q^{ \mp {n(n-1) \over 2}}}{ [n]_q !} \; , \;\;\;
[n]_q ! := [1]_q \, [2]_q \cdots [n]_q \; , \;\;\;
[n]_q = \frac{(q^n - q^{-n})}{(q - q^{-1})}  \right) ,
$$
\vspace{0.1cm}
\be
\lb{santis1}
\begin{array}{c}
\sigma_i \, S_{1 \to n} = S_{1 \to n} \, \sigma_i = q \, S_{1 \to n} \;\;
(i =1, \dots n-1) \; , \\ [0.2cm]
\sigma_i \, A_{1 \to n} = A_{1 \to n} \, \sigma_i = - \frac{1}{q} \, A_{1 \to n} \;\;
(i =1, \dots n-1) \; ,
\end{array}
\ee
which are symmetrizers and antisymmetrizers, respectively (see \cite{Gur}).
The normalization factors  $a^{\mp}_n$
have been introduced in (\ref{santis}) in order to obtain the idempotent conditions
$S_{1 \to n}^2 = S_{1 \to n}$ and $A_{1 \to n}^2 = A_{1 \to n}$.
Here we additionally suppose that $[n]_q \neq 0$, $\forall n = 1, \dots, M+1$.
The first two idempotents are
\be
\lb{3.3.111}
S_{12} =  \frac{1}{[2]_q}(q^{-1} + \sigma_1)  \; , \;\;\;
A_{12} = \frac{1}{[2]_q}(q - \sigma_1)  \; .
\ee
Note that eqs. (\ref{santis1}) immediately follow from the factorization
relations (\ref{rfactor}), (\ref{lfactor}), the form of the first
idempotents (\ref{3.3.111}) and Hecke condition (\ref{ahecke}).
The projectors $S_{1 \to n}$ and $A_{1 \to n}$ (\ref{santis}) correspond to
the Young tableaux which have only one row and one column
$$
P \left(
\begin{tabular}{|c|c|c|}
\hline
$\!\!\!$ {\small 1} $\!\!$ & $\dots$ &
$\!\!\!$ {\small n} $\!\!\!$\\
\hline
\end{tabular}
\right)
= S_{1 \to n} \; , \;\;\;
P \left(\,
\begin{tabular}{|c|}
\hline
$\!\!$ {\small 1} $\!\!$\\
\hline
$\vdots$ \\
\hline
$\!\!$ {\small n} $\!\!$\\
\hline
\end{tabular}
\, \right)
=  A_{1 \to n}  \; .
$$
It follows directly from (\ref{santis1})
that the idempotents $S_{1 \to M+1}$ and $A_{1 \to M+1}$ are central
in $H_{M+1}(q)$.

Consider now the elements $y_i$ ($i=1, \dots ,M+1$) (\ref{combr})
which generate a commutative subalgebra $Y_{M+1}$ in $H_{M+1}$.
It can be proved that $Y_{M+1}$
is a maximal commutative subalgebra in $H_{M+1}$. The elements
$y_i$ are called {\it Jucys - Murphy elements} and can be easily rewritten in
the form (by using Hecke condition (\ref{ahecke}) and braid relations (\ref{braidg}))
\be
\lb{jucmu}
\begin{array}{c}
y_1 = 1 \; , \;\;\;\;
y_i \, = \, \sigma_{i-1} \, y_{i-1} \,  \sigma_{i-1}
\, = \, \sigma_{i-1} \dots \sigma_2 \, \sigma_1^2 \, \sigma_2
\dots \sigma_{i-1} = \\ [0.2cm]
= \lambda \sigma_{i-1} \dots \sigma_2 \, \sigma_1 \, \sigma_2
\dots \sigma_{i-1} +
\sigma_{i-1} \dots \sigma_3 \, \sigma_2^2 \, \sigma_3
\dots \sigma_{i-1} = \dots = \\ [0.2cm]
= \lambda \sum\limits_{k=1}^{i-2}
\, \sigma_{i-1} \dots \sigma_{k+1} \, \sigma_k \, \sigma_{k+1}
\dots \sigma_{i-1} + \lambda \sigma_{i-1} + 1 = \\ [0.2cm]
= \lambda \sum\limits_{k=1}^{i-2}
\, \sigma_{k} \dots \sigma_{i-2} \, \sigma_{i-1} \, \sigma_{i-2}
\dots \sigma_{k} + \lambda \sigma_{i-1} + 1 \; , \;\;\;\;\;\;
i=2, \dots ,M+1 \; .
\end{array}
\ee
It is interesting that the idempotents (\ref{santis})
which correspond to the symmetrizers and antisymmetrizers (the Young
tableaux is only one row or column) can be constructed in the
different way as polynomial functions of the elements $y_n$.
\begin{proposition}\label{prop9}
 The idempotents  $S_{1 \to n}$ and $A_{1 \to n}$
$(n = 2, \dots M+1)$
(\ref{santis}) are expressed in term of the Jucys-Murphy elements as
\be
\lb{sant01}
S_{1 \to n} = \frac{(y_2 - q^{-2})}{(q^2 - q^{-2})}  \,
\frac{(y_3 - q^{-2})}{(q^4 - q^{-2})}
 \cdots \frac{(y_n - q^{-2})}{(q^{2(n-1)} - q^{-2})} \; ,
\ee
\be
\lb{sant02}
A_{1 \to n} =  \frac{(y_2 - q^{2})}{(q^{-2} - q^{2})}  \,
\frac{(y_3 - q^{2})}{(q^{-4} - q^{2})}
 \cdots \frac{(y_n - q^{2})}{(q^{2(1-n)} - q^{2})} \; .
\ee
\end{proposition}

\noindent
{\bf Proof.} We note that expressions (\ref{3.3.111}) for the
first two projectors are written as
$$
S_{12} =
\frac{(\sigma_1^2 - q^{-2})}{(q^2 - q^{-2})}  \; , \;\;\;
A_{12} =
 \frac{(\sigma_1^2 - q^{2})}{(q^{-2} - q^{2})}  \; ,
$$
and therefore eqs.  (\ref{sant01}) and (\ref{sant02}) are
valid for $n=2$.  We prove eqs. (\ref{sant01})
and (\ref{sant02}) by induction.
Let eqs. (\ref{sant01}) and (\ref{sant02}) be
correct for some $n=k$. We need to prove these
 equation for $n=k+1$ or we have to show that
\be
\lb{sant03}
S_{1 \to k+1} = S_{1 \to k} \cdot
\frac{(y_{k+1} - q^{-2})}{(q^{2k} - q^{-2})} \; , \;\;\;\;
A_{1 \to k+1} = A_{1 \to k} \cdot
\frac{(y_{k+1} - q^{2})}{(q^{-2k} - q^{2})} \; .
\ee
We prove only the first eq. in (\ref{sant03})
 (the proof of the second eq. in (\ref{sant03}) is analogous).
We substitute in (\ref{sant03})
 the last expression for Jucys-Murphy elements $y_{k+1}$
(\ref{jucmu}) and take into account (\ref{santis1}). As a result we
obtain for the first eq. in (\ref{sant03})
\be
\lb{sant04}
\begin{array}{c}
S_{1 \to k+1} = \frac{1}{(q^{2k} - q^{-2})} \, S_{1 \to k} \,
\left( \lambda (q^{k-1} \sigma_k \dots \sigma_1 +
q^{k-2} \sigma_k \dots \sigma_2 + \dots +
 \sigma_k ) + 1 - q^{-2}  \right) = \\ \\
 = \frac{q^{k}}{[k+1]_q} \, S_{1 \to k} \,
\left( q^{k} \sigma_k \dots \sigma_1 +
q^{k-1} \sigma_k \dots \sigma_2 + \dots +
q \,  \sigma_k  + 1  \right) \; ,
\end{array}
\ee
which coincides with the definition of symmetrizers
(\ref{danti}), (\ref{santis}). This ends the proof of the
induction and hence this Proposition.
\hfill \qed

\vspace{0.2cm}

Note that the idempotents  $S_{1 \to M+1}$ and $A_{1 \to M+1}$ are
central in the algebra $H_{M+1}(q)$
and represented as the polynomials
$\sim  (y_2 - t) (y_3 - t) \cdots (y_{M+1} - t)$,
where $t = q^{\mp 2}$ (see (\ref{sant01}), (\ref{sant02})),
which are symmetric functions in variables $\{ y_i \}$
$(i = 2, \dots , M+1)$. In view of this,
one can conjecture that
all symmetric functions in $y_i$
 generate the central subalgebra $Z_{M+1}$ in the Hecke algebra
$H_{M+1}(q)$. Indeed, to prove this fact we need only to
check the relations:
$[\sigma_k, \, y_n + y_{n+1}] = 0 = [\sigma_k, \, y_n  y_{n+1}]$
for all $k < n+1$.

New identities for the elements $y_i$ follow from
the representations (\ref{sant01}) and (\ref{sant02})
(if we use eqs. (\ref{santis1})):
\be
\lb{ident1}
(y_i - q^{2(i-1)}) \, S_{1 \to n} = 0 \Rightarrow
(y_i - q^{2(i-1)})  (y_2 - q^{-2}) (y_3 - q^{-2}) \cdots (y_{n} - q^{-2}) = 0 \; ,
\ee
\be
\lb{ident2}
(y_i - q^{2(1-i)}) \, A_{1 \to n} =  0 \Rightarrow
(y_i - q^{2(1-i)})  (y_2 - q^{2}) (y_3 - q^{2}) \cdots (y_{n} - q^{2}) = 0 \; ,
\ee
($i = 2, \dots , n$). Then, two new types of idempotents
(which are primitive orthogonal idempotents
for the subalgebra $H_{n} \in M_{M+1}$) are obtained from
these identities:
\be
\lb{proj1}
P \left(
\begin{tabular}{|c|c|c|}
\hline
$\!\!$ {\small 1} $\!\!$ & $\dots$ &
$\!\!\!\!$ {\small n-1} $\!\!\!\!$\\
\hline
$\!\!$ {\small n} $\!\!$ & \multicolumn{1}{c}{}\\
\cline{1-1}
\end{tabular}
\right) =
{(y_n - q^{2(n-1)} )\over (q^{-2} - q^{2(n-1)})}
\prod_{k=1}^{n-1} \, {(y_k - q^{-2}) \over (q^{2(k-1)} - q^{-2})} =
\ee
$$
= S_{1\to n-1} - S_{1 \to n} = \frac{[n-1]_q}{[n]_q} \, S_{1\to n-1} \,
\sigma_{n-1}(q)  \, S_{1\to n-1} \; ,
$$
\be
\lb{proj2}
P \left(
\begin{tabular}{|c|c|}
\hline
$\!\!$ {\small 1} $\!\!$ &
$\!\!$ {\small n} $\!\!$\\
\hline
$\vdots$  \\
\cline{1-1}
$\!\!\!\!$ {\small n-1} $\!\!\!\!$ & \multicolumn{1}{c}{}\\
\cline{1-1}
\end{tabular}
\right) =
{(y_n - q^{2(1-n)} )\over (q^{2} - q^{2(1-n)})}
\prod_{k=1}^{n-1} \, {(y_k - q^{2}) \over (q^{2(1-k)} - q^{2})} =
\ee
$$
= A_{1\to n-1} - A_{1\to n} = \frac{[n-1]_q}{[n]_q} \, A_{1\to n-1} \,
\sigma_{n-1}(q^{-1})  \, A_{1\to n-1} \; ,
$$
where \cite{Jimb1}, \cite{Jon3}
$$
\sigma_n(x) : = \lambda^{-1} \,(x^{-1} \sigma_n - x \sigma_n^{-1}) \; ,
$$
are Baxterized elements for the algebra $H_{M+1}(q)$
 (the $R$-matrix representations of these elements
are given in (\ref{3.5.5})).
We consider properties of these elements below;
see eq. (\ref{baxtH}) and further discussion.
The idempotents (\ref{proj1}), (\ref{proj2}) are not central in $H_{M+1}$
but they are the elements of the commutative subalgebra $Y_{M+1}$.

\subsubsection{Primitive orthogonal idempotents
 in $H_{M+1}$ and Young tableaux\label{idemp}}

Now we describe the general construction
(see \cite{IsOg3} and references
therein) of all primitive orthogonal idempotents
$e_\alpha \in H_{M+1}$ which are elements of $Y_{M+1}$
(i.e. functions of the elements
$y_i$). All these idempotents are common eigenidempotents of $y_i$:
$$
y_i e_\alpha = e_\alpha y_i = a_i^{(\alpha)} e_\alpha \;\;\;\;
(i=1, \dots ,M+1) \; ,
$$
 where $a_i^{(\alpha)}$ are eigenvalues.
We denote by ${\rm Spec}(y_1, \dots , y_{M+1})$ the set
 of strings of $(M+1)$ eigenvalues
$\Lambda(e_\alpha) : = (a_1^{(\alpha)}, \dots , a^{(\alpha)}_{M+1})$
$(\forall \alpha)$.
The eigenidempotents $e_\alpha$ define left (and right)
submodules $H_{M+1} \cdot e_\alpha$
(and $e_\alpha \cdot H_{M+1}$) in the regular bimodule of $H_{M+1}$.

\vspace{0.3cm}
\noindent
{\bf Lemma 1.} {\it The eigenidempotents $e$ and $e'$ with
eigenvalues $a_i = a_i'$
$(\forall i=1, \dots,M)$ and $a_{M+1} \neq a_{M+1}'$ define
different
left (right) submodules in the regular bimodule of $H_{M+1}$.}

\vspace{0.1cm}
\noindent
{\bf Proof.} We proof this Lemma only for
the left submodules  $H_{M+1} \cdot e$,
$H_{M+1} \cdot e'$. The case of right submodules can be considered analogously. Let $v$ and $v'$ be respectively elements
of submodules  $H_{M+1} \cdot e$ and
$H_{M+1} \cdot e'$.
Consider central element $Z = y_1 y_2 \cdots y_{M+1}$
(symmetric function of $y_i$). There are no elements
$X \in H_{M+1}$ such that $v' = X \, v$,
since the left action of $Z$ on elements $X \, v$ and $v'$ produces different eigenvalues. \hfill
\qed

Now we introduce the important intertwining elements \cite{Isaev1}
(in another form these elements appeared in \cite{Cher5}):
\be
\lb{impint}
U_{n+1} = \sigma_n y_n - y_n \sigma_n =
\sigma_n y_n - \sigma_n^{-1} y_{n+1}  =
y_{n+1} \sigma_n^{-1} - y_n \sigma_n =
\ee
\be
\lb{impint1}
= (y_{n+1} - y_n) \sigma_n - \lambda y_{n+1} =
 \sigma_n (y_{n} - y_{n+1}) + \lambda y_{n+1} \;\; (1 \leq n \leq M) \; ,
\ee
subject to relations\footnote{The definition (\ref{impint}) of
intertwining elements is not unique. One can multiply $U_{n+1}$ by
a function $f(y_n, y_{n+1})$: $U_{n+1}  \rightarrow U_{n+1} f(y_n, y_{n+1})$.
Then, all eqs. (\ref{importt}) -- (\ref{import}) are valid if $f$ satisfies
$f(y_n, y_{n+1}) f(y_{n+1}, y_{n}) =1$.}
\be
\lb{importt}
\begin{array}{c}
U_{n+1} y_n = y_{n+1} U_{n+1} \; , \;\;\;
U_{n+1} y_{n+1} = y_{n} U_{n+1} \; , \\ [0.2cm]
\, [U_{n+1}, \, y_k ] = 0 \;\; (k \neq n,n+1)  \; ,
\end{array}
\ee
\be
\lb{importtt}
U_n \, U_{n+1} \, U_n = U_{n+1} \, U_{n} \, U_{n+1} \; ,
\ee
\be
\lb{import}
U_{n+1}^2 = (q y_{n} - q^{-1} \, y_{n+1}) \,
(q \, y_{n+1} - q^{-1} \, y_n) \; .
\ee

\vspace{0.3cm}
\noindent
{\bf Lemma 2.} {\it The eigenidempotents $e$ and $e'$ with
eigenvalues:
\be
\lb{eigen}
\begin{array}{c}
a_i = a_i' \;\;
(\forall i=1, \dots,M-1) \; , \\ [0.2cm]
a_{M} = a_{M+1}' \; , \;\;\; a_{M+1} = a_{M}' \; , \;\;\;  a_{M} \neq q^{\pm 2}  a_{M+1} \; ,
\end{array}
\ee
belong to the same irreducible sub-bimodule in the regular
bimodule of $H_{M+1}$.}

\vspace{0.1cm}
\noindent
{\bf Proof.} Since the algebra $Y_{M+1}$ generated by $\{y_1, \dots,y_{M+1}\}$
is maximal commutative subalgebra in $H_{M+1}$ we have $e' = e''$
if $\Lambda(e') = \Lambda(e'')$.
Then, using intertwining element $U_{M+1}$ (\ref{impint}) we construct
the eigenidempotent
$$
e'' = \frac{1}{(q^2 a_M - a_{M+1})(a_{M+1} - q^{-2}a_M)}
\, U_{M+1} \; e \; U_{M+1} \; , \;\;\;\;\;\;\;
(e'')^2 = e'' \; ,
$$
which is well defined in view of the last condition in (\ref{eigen}).
The element $U_{M+1} \, e \, U_{M+1}$
is not equal to zero, since $U_{M+1}^2 \, e \, U_{M+1}^2 =
(q^2 a_M - a_{M+1})^2(a_{M+1} - q^{-2}a_M)^2 e \neq 0$. This inequality
follows from the last condition in (\ref{eigen}).
For the element $e'' \sim U_{M+1}\, e\, U_{M+1}$
we have $\Lambda(e'') = \Lambda(e')$
in view of (\ref{importt}).
Thus, $e'' = e' \Rightarrow e' \sim U_{M+1}\, e \, U_{M+1}$
and the eigenidempotents
$e$ and $e'$ belong to the same irreducible sub-bimodule
in the regular bimodule
of $H_{M+1}$. \hfill \qed

\vspace{0.2cm}

Consider a Young diagram $[\nu]_{M+1}$ with $(M+1)$ nodes.
We place the numbers $1, \dots , M+1$ into the nodes of the diagram
in such a way that these numbers are arranged along rows and columns in
ascending order in right and down directions. Such diagram
is called a standard Young tableau ${\sf T}_{[\nu]_{M+1}}$.
Then we associate a number $q^{2(n-m)}$ (the "content")
to each node of the standard Young tableau, where $(n,m)$ are coordinates of the node.
Example:

\unitlength=6mm
\begin{picture}(25,4)
\put(5,3.4){\vector(1,0){7}}
\put(5,3.4){\vector(0,-1){3.5}}
\put(12,3.5){$n$}
\put(4.2,0){$m$}

\put(5.5,1.3){$
\begin{tabular}{|c|c|c|c|}
\hline
  $\!\!\!\! \!\! ^{1}$  $_1$ & $\!\! ^2$  $_{q^2}$  &
  $\!\! ^4$  $_{q^4}$ & $\!\! ^6$  $_{q^6}$ \\[0.2cm]
\hline
  $\!\! ^3$ $\!\! _{q^{\! -2}}$  &   $\!\!\!\!\! ^5$ $_1$  & $\!\!\! ^8$ $_{q^2}$  &
  \multicolumn{1}{c}{} \\[0.2cm]
\cline{1-3}
  $\!\! ^7$ $\!\! _{q^{\! -4}}$  & \multicolumn{1}{c}{} \\[0.2cm]
\cline{1-1}
\end{tabular}
$}
\end{picture}
\vspace{-1cm}
\be
\label{111}
{}
\ee
\noindent
In general, for the tableau ${\sf T}_{[\nu]_{M+1}}$,
the $i$-th node
with coordinates $(n,m)$ looks like:
$\begin{tabular}{|c|}
\hline
  $\!\! ^{i}$  $\!\! _{q^{2(n-m)}} \!\!$  \\[0.2cm]
\hline
\end{tabular}$. Thus, to each standard Young tableau $[\nu]_n$ one can associate a
string of numbers
$\Lambda =(a_1, \dots , a_n)$ with $a_i = q^{2(n-m)}$. E.g.,
a standard Young tableau (\ref{111})
corresponds to a string
$$
\Lambda =(1,q^2,q^{-2},q^4,1,q^6,q^{-4},q^2) \; .
$$

\vspace{0.2cm}

Now we associate Young tableaux ${\sf T}_{[\nu]_{M+1}}$
(related to the primitive orthogonal
idempotents) with paths in Young-Ogievetsky graph.
By definition Young-Ogievetsky graph is a Young graph with
vertices, which are Young diagrams,
with edges, which indicate inclusions of diagrams
(or a branching of representations), and with
numbers (colours) on the edges corresponding to the eigenvalues of the
Jucys-Murphy elements\footnote{To our knowledge
O.Ogievetsky was the first
who propose to associate the eigenvalues of the
Jucys-Murphy elements to the edges of Young graph. Usually the
indices on the edges of the graphs of Young type
correspond to the multiplicity of the branching.
In this case the Young graph is called Bratelli diagram. In our case
all multiplicities are equal to $1$.}. For example,
the coloured Young-Ogievetsky graph for $H_4$ is:

\setlength{\unitlength}{1600sp}
\begingroup\makeatletter\ifx\SetFigFont\undefined%
\gdef\SetFigFont#1#2#3#4#5{%
  \reset@font\fontsize{#1}{#2pt}%
  \fontfamily{#3}\fontseries{#4}\fontshape{#5}%
  \selectfont}%
\fi\endgroup%
\begin{picture}(6466,7391)(-1768,-6894)
\thicklines
{\put(4260,-4336){\vector( 3,-4){1350}}}
{\put(4160,-4336){\vector( 0,-4){1820}}}
{\put(4126,-4411){\vector(-2,-3){1107.692}}}
{\put(2251,-4111){\vector( 1,-3){682.500}}}
{\put(2101,-4186){\vector(-1,-2){930}}
}%
{\put(5026,-2386){\vector( 2,-3){900}}
}%
{\put(4801,-2536){\vector(-1,-2){660}}
}%
{\put(3301,-2386){\vector( 1,-2){770}}}
{\put(3151,-2386){\vector(-2,-3){969.231}}
}%
{\put(4276,-961){\vector( 1,-2){570}}}
{\put(4200,-961){\vector(-3,-4){900}}
}%
{\put(4276,100){\vector( 0,-1){910}}}
{\thinlines
\put(1276,-6136){\circle*{150}}}
{\put(7351,-6811){\circle*{150}}}
{\put(7351,-6586){\circle*{150}}}
\thicklines
{\put(6051,-4486){\vector( -1,-4){380}}}
\put(2050,-5536){\makebox(0,0)[lb]{\smash{\SetFigFont{10}{12.0}{\rmdefault}{\mddefault}
{\updefault}{$q^{-2}$}%
}}}
\put(3376,-4861){\makebox(0,0)[lb]{\smash{\SetFigFont{10}{12.0}{\rmdefault}{\mddefault}
{\updefault}{$q^4$}%
}}}
\put(4951,-4936){\makebox(0,0)[lb]{\smash{\SetFigFont{10}{12.0}{\rmdefault}{\mddefault}{
\updefault}{$q^{-4}$}%
}}}
\put(5551,-2836){\makebox(0,0)[lb]{\smash{\SetFigFont{10}{12.0}{\rmdefault}{\mddefault}
{\updefault}{$q^{-4}$}%
}}}
\put(6826,-5011){\makebox(0,0)[lb]{\smash{\SetFigFont{10}{12.0}{\rmdefault}{\mddefault}
{\updefault}{$q^{-6}$}%
}}}
\put(1276,-4936){\makebox(0,0)[lb]{\smash{\SetFigFont{10}{12.0}{\rmdefault}{\mddefault}
{\updefault}{$q^6$}%
}}}
\put(2401,-2911){\makebox(0,0)[lb]{\smash{\SetFigFont{10}{12.0}{\rmdefault}{\mddefault}
{\updefault}{$q^4$}%
}}}
\put(4296,-5611){\makebox(0,0)[lb]{\smash{\SetFigFont{10}{12.0}{\rmdefault}{\mddefault}
{\updefault}{$1$}%
}}}

\put(8000,-300){\makebox(0,0)[lb]{\smash{\SetFigFont{12}{12.0}{\rmdefault}{\mddefault}
{\updefault}{$= \; y_1$}%
}}}
\put(8000,-1511){\makebox(0,0)[lb]{\smash{\SetFigFont{12}{12.0}{\rmdefault}{\mddefault}
{\updefault}{$= \; y_2$}%
}}}
\put(8000,-3011){\makebox(0,0)[lb]{\smash{\SetFigFont{12}{12.0}{\rmdefault}{\mddefault}
{\updefault}{$= \; y_3$}%
}}}
\put(8000,-5011){\makebox(0,0)[lb]{\smash{\SetFigFont{12}{12.0}{\rmdefault}{\mddefault}
{\updefault}{$= \; y_4$}%
}}}

\put(6076,-5461){\makebox(0,0)[lb]{\smash{\SetFigFont{10}{12.0}{\rmdefault}{\mddefault}
{\updefault}{$q^2$}%
}}}
\put(3126,-3511){\makebox(0,0)[lb]{\smash{\SetFigFont{10}{12.0}{\rmdefault}{\mddefault}
{\updefault}{$q^{-2}$}
}}}
\put(4276,-2911){\makebox(0,0)[lb]{\smash{\SetFigFont{10}{12.0}{\rmdefault}{\mddefault}
{\updefault}{$q^2$}%
}}}
\put(4651,-1411){\makebox(0,0)[lb]{\smash{\SetFigFont{10}{12.0}{\rmdefault}{\mddefault}
{\updefault}{$q^{-2}$}%
}}}
\put(3526,-1411){\makebox(0,0)[lb]{\smash{\SetFigFont{10}{12.0}{\rmdefault}{\mddefault}
{\updefault}{$q^2$}%
}}}
\put(4051,-286){\makebox(0,0)[lb]{\smash{\SetFigFont{10}{12.0}{\rmdefault}{\mddefault}
{\updefault}{$1$}%
}}}
{\put(6151,-4411){\vector( 3,-4){1161}}
}%
{\thinlines
\put(7351,-6361){\circle*{150}}
}%
{\put(4220,104){$\emptyset$}}
{\put(4276,-886){\circle*{150}}
}%
{\put(6001,-3886){\circle*{150}}
}%
{\put(6001,-4336){\circle*{150}}
}%
{\put(4276,-4036){\circle*{150}}
}%
{\put(4051,-4036){\circle*{150}}
}%
{\put(4051,-4261){\circle*{150}}
}%
{\put(2476,-4036){\circle*{150}}
}%
{\put(2251,-4036){\circle*{150}}
}%
{\put(2026,-4036){\circle*{150}}
}%
{\put(4876,-2236){\circle*{150}}
}%
{\put(4876,-2461){\circle*{150}}
}%
{\put(3376,-2236){\circle*{150}}
}%
{\put(3151,-2236){\circle*{150}}
}%
{\put(6001,-4111){\circle*{150}}
}%
{\put(1051,-6136){\circle*{150}}
}%
{\put(7351,-6136){\circle*{150}}
}%
{\put(5651,-6736){\circle*{150}}}
{\put(5651,-6511){\circle*{150}}}
{\put(5876,-6286){\circle*{150}}}
{\put(5651,-6286){\circle*{150}}}
{\put(4376,-6586){\circle*{150}}}
{\put(4151,-6586){\circle*{150}}}
{\put(4156,-6361){\circle*{150}}}
{\put(4371,-6361){\circle*{150}}}
{\put(3226,-6286){\circle*{150}}
}%
{\put(3001,-6286){\circle*{150}}
}%
{\put(2776,-6286){\circle*{150}}
}%
{\put(2776,-6511){\circle*{150}}
}%
{\put(1726,-6136){\circle*{150}}
}%
{\put(1501,-6136){\circle*{150}}}

\put(8401,-7036){{\bf Fig.4.1}. \it Young-Ogievetsky graph for $H_4(q)$.}
\end{picture}

\vspace{0.9cm}
\noindent
The paths (associated to Young tableaux)
start from the top vertex $\emptyset$ and finish at the
vertex labeled by the Young diagram of the same shape as the tableaux.
The dimension of the corresponding representation of the Hecke algebra
is equal to the number of standard tableaux of this shape or, as we see,
the number of paths which lead to this Young diagram from
top $\emptyset$.
For example: the path
$\{ \emptyset \stackrel{1}{\rightarrow} \bullet \stackrel{q^2}{\rightarrow} \bullet \, \bullet
\stackrel{q^{-2}}{\rightarrow}
\begin{array}{cc}
\bullet & \!\!\!\! \bullet  \\[-0.25cm]
\bullet
\end{array} \stackrel{1}{\rightarrow}
\begin{array}{cc}
\bullet & \!\!\!\! \bullet  \\[-0.25cm]
\bullet  & \!\!\!\! \bullet
\end{array} \}$ corresponds to the tableau
$\begin{tabular}{|c|c|}
\hline
$\!\!\!$
{\small 1} $\!\!\!$ & $\!\!\!$ {\small 2} $\!\!\!$ \\
\hline
$\!\!\!$ {\small 3} $\!\!\!$ & $\!\!\!$ {\small 4} $\!\!\!$ \\
\cline{1-2}
\end{tabular}$, i.e. the shape of the tableau (Young diagram) is given
by the shape of the last vertex of the path, while the numbers
in nodes of the tableau show in which sequence
the points $\bullet$ appear in the vertices along the path. The edge colours of the path (or contents of the nodes of
the standard tableau as it is explained in (\ref{111}))
are the eigenvalues of the Jucys-Murphy
elements
$y_1 =1$, $y_2 = q^2$, $y_3 = q^{-2}$, $y_4 =1$
obtained by their action on the idempotent
$P \left( \begin{tabular}{|c|c|}
\hline
$\!\!\!$
{\small 1} $\!\!\!$ & $\!\!\!$ {\small 2} $\!\!\!$ \\
\hline
$\!\!\!$ {\small 3} $\!\!\!$ & $\!\!\!$ {\small 4} $\!\!\!$ \\
\cline{1-2}
\end{tabular}\right)$ .  Then, the explicit formula
for this idempotent can be constructed by induction.
 Namely, we take the explicit form of
the previous idempotent
$
P \left(
\begin{tabular}{|c|c|}
\hline
$\!\!\!$  1 $\!\!\!$ & $\!\!\!$  2 $\!\!\!$ \\
\hline
$\!\!\!$  3 $\!\!\!$ & \multicolumn{1}{c}{}\\
\cline{1-1}
\end{tabular}
\right)
$
(related to the previous vertex of the path)
and multiply it by the factors $(y_4 - 1)$,
$(y_4 - q^4)$ and $(y_4 - q^{-4})$, which correspond to
possible colours of outgoing edges from vertex
$\begin{array}{cc}
\bullet & \!\!\!\! \bullet  \\[-0.25cm]
\bullet
\end{array}$,  to obtain characteristic identity
\be
\lb{isa01}
P \left(
\begin{tabular}{|c|c|}
\hline
$\!\!\!$  1 $\!\!\!$ & $\!\!\!$  2 $\!\!\!$ \\
\hline
$\!\!\!$  3 $\!\!\!$ & \multicolumn{1}{c}{}\\
\cline{1-1}
\end{tabular}
\right) (y_4 - 1)(y_4 - q^4)(y_4 - q^{-4}) = 0 \; .
\ee
Then, to forbid
the moving from the vertex
$\begin{array}{cc}
\bullet & \!\!\!\! \bullet  \\[-0.25cm]
\bullet
\end{array}$
along the edges with labels $q^4$ and $q^{-4}$
and move along the edge with the index $1$ to the vertex
$\begin{array}{cc}
\bullet & \!\!\!\! \bullet  \\[-0.25cm]
\bullet  & \!\!\!\! \bullet
\end{array}$ we remove from the left hand side of (\ref{isa01})
the factor $(y_4 - 1)$. As a result, we obtain
(after an obvious renormalization)
\be
\lb{import3}
P \left(
\begin{tabular}{|c|c|}
\hline
$\!\!\!$
{\small 1} $\!\!\!$ & $\!\!\!$ {\small 2} $\!\!\!$ \\
\hline
$\!\!\!$ {\small 3} $\!\!\!$ & $\!\!\!$ {\small 4} $\!\!\!$ \\
\cline{1-2}
\end{tabular}
\right) =
P \left(
\begin{tabular}{|c|c|}
\hline
$\!\!\!$  1 $\!\!\!$ & $\!\!\!$  2 $\!\!\!$ \\
\hline
$\!\!\!$  3 $\!\!\!$ & \multicolumn{1}{c}{}\\
\cline{1-1}
\end{tabular}
\right) \frac{(y_4 - q^4)(y_4 - q^{-4})}{(1 - q^4)(1 - q^{-4})} \; .
\ee
In the same way one can deduce the chain of identities
\be
\lb{import2}
P \left(
\begin{tabular}{|c|c|}
\hline
$\!\!\!$  1 $\!\!\!$ & $\!\!\!$  2 $\!\!\!$ \\
\hline
$\!\!\!$  3 $\!\!\!$ & \multicolumn{1}{c}{}\\
\cline{1-1}
\end{tabular}
\right) =
P \left(
\begin{tabular}{|c|c|}
\hline
$\!\!\!$  1 $\!\!\!$ & $\!\!\!$  2 $\!\!\!$ \\
\hline
\end{tabular}
\right) \frac{(y_3 - q^4)}{(q^{-2} - q^4)}
 =
P \left(
\begin{tabular}{|c|}
\hline
$\!\!\!$  1 $\!\!\!$ \\
\hline
\end{tabular}
\right) \frac{(y_2 - q^{-2})(y_3 - q^4)}{(q^2 - q^{-2})(q^{-2} - q^4)} \; ,
\ee
where we fix $P \left(
\begin{tabular}{|c|}
\hline
$\!\!\!$  1 $\!\!\!$ \\
\hline
\end{tabular}
\right) =1$  by definition.
Using (\ref{import2}),  the final formula for (\ref{import3})
can be written as
\be
\lb{import4}
P \left(
\begin{tabular}{|c|c|}
\hline
$\!\!\!$
{\small 1} $\!\!\!$ & $\!\!\!$ {\small 2} $\!\!\!$ \\
\hline
$\!\!\!$ {\small 3} $\!\!\!$ & $\!\!\!$ {\small 4} $\!\!\!$ \\
\cline{1-2}
\end{tabular}
\right) =  \frac{(y_2 - q^{-2})(y_3 - q^4)}{(q^2 - q^{-2})(q^{-2} - q^4)}
\frac{(y_4 - q^4)(y_4 - q^{-4})}{(1 - q^4)(1 - q^{-4})} \; .
\ee
We note that the described procedure leads automatically to the
idempotents which are orthogonal to each other.

This example has demonstrated that all information about primitive orthogonal idempotents
for the A-type Hecke algebra is encoded in the Young-Ogievetsky
(YO) graph given in Fig. 4.1.
Thus, we need to justify this graph and its edge colours.
First of all we prove the following statement.
\begin{proposition}\label{prop9b}
{\it The spectrum of the Jucys-Murphy operators $y_j$
(possible edge indices of the YO graph) for $H_{M+1}$ is such that
\be
\lb{spec}
{\rm Spec}(y_j) \subset \{ q^{2 {\bf Z}_j} \} \;\;\;\; \forall j = 1,2, \dots , M+1 \; ,
\ee
where ${\bf Z}_j$ denotes the set of integer numbers
$\{ 1-j, \dots , -2,-1,0,1,2, \dots ,j-1 \}$.}
\end{proposition}
{\bf Proof.}
We use
the important intertwining elements (\ref{impint}), (\ref{impint1})
and prove (\ref{spec}) by induction. From Hecke condition (\ref{ahecke})
we have
\be
\lb{spec3}
(y_2 - q^2)(y_2 - q^{-2}) = 0 \; .
\ee
Thus, ${\rm Spec} (y_2)$ satisfies (\ref{spec}).
Assume that the spectrum of $y_{j-1}$ satisfies (\ref{spec}) for some $j \geq 3$.
Consider a characteristic equation for $y_{j-1}$ (cf. (\ref{spec3})):
$$
f(y_{j-1}) : = \prod_\alpha (y_{j-1} - a^{(\alpha)}_{j-1}) = 0 \;\;\;
(a^{(\alpha)}_{j-1} \in {\rm Spec} (y_{j-1})) \; .
$$
Using operators $U_j$ and their properties (\ref{importt}),
(\ref{import}) we deduce
\be
\lb{spec1}
\begin{array}{c}
0= U_j f(y_{j-1}) U_j = f(y_{j}) U^2_j
= f(y_{j}) (q^2 y_{j-1} - y_{j})( y_{j} - q^{-2} y_{j-1}) \; .
\end{array}
\ee
which means that
\be
\lb{spec11}
{\rm Spec} (y_j) \subset
\left( {\rm Spec}(y_{j-1}) \cup q^{\pm 2} \cdot {\rm Spec}(y_{j-1}) \right) \; ,
\ee
and it justifies (\ref{spec}). \hfill \qed

\subsubsection{Irreducible representations of $H_{M+1}$ and
recurrence formula for primitive idempotents\label{jmel2}}

In \cite{OV}, A. Okounkov and A. Vershik developed new approach
 to the construction of the irreducible representations of
 symmetric group (we review this approach in \cite{IsRub2}, section 4.6).
 Here we generalize (following \cite{Jon1},
 \cite{Mur}, \cite{W1}, \cite{Cher5}, \cite{OgPya}, \cite{IsOg3})
  the Okounkov-Vershik approach to the case
 of the Hecke algebra.

Consider a subalgebra $\hat{H}_2^{(i)}$ in
$H_{M+1}$ with generators $y_i$, $y_{i+1}$
and $\sigma_i$ (for fixed $i \leq M$). We investigate
 (see \cite{IsOg3}) representations of
$\hat{H}_2^{(i)}$ in the case when the elements $y_i$, $y_{i+1}$ are
diagonalizable. Let $e$ be a common eigenidempotent of $y_i$, $y_{i+1}$:
$y_i e = a_i e$, $y_{i+1} e = a_{i+1} e$. Then, the left action
of $\hat{H}_2^{(i)}$ closes on elements
$v_1 = e$, $v_2 = \sigma_i e$ and is given by matrices:
\be
\lb{abmw1a}
\sigma_i =
\left(
\begin{array}{cc}
0 & 1 \\
1 & \lambda
\end{array}
\right)  , \;
y_i =
\left(
\begin{array}{cc}
a_i & - \lambda a_{i+1}  \\
0 & a_{i+1}
\end{array}
\right) , \;
y_{i+1} =
\left(
\begin{array}{cc}
a_{i+1} & \lambda a_{i+1}  \\
0 & a_i
\end{array}
\right) ,
\ee
where we have used the standard convention
$y v_i = v_j \, y_{ji}$
to produce matrix representations $||y_{ji}||$
for operators $y$.

The operators $y_i$, $y_{i+1}$ (\ref{abmw1a}) can be simultaneously
diagonalized by the transformation
$y \rightarrow V^{-1} y V$ where
$$
V =
\left(
\begin{array}{cc}
1 \! & \! \frac{\lambda \, a_{i+1}}{a_i -a_{i+1}} \\
0 \! & \! 1
\end{array}
\right)  \; , \;\;\;
V^{-1} =
\left(
\begin{array}{cc}
1 & - \frac{\lambda \, a_{i+1}}{a_i - a_{i+1}} \\
0 & 1
\end{array}
\right)  \; .
$$
As a result we obtain the following matrix representation
\be
\lb{ah1}
\sigma_i =
\left(\!\!
\begin{array}{cc}
- \frac{\lambda \, a_{i+1}}{a_i -a_{i+1}} \; & \;
\frac{(a_i - q^2 \, a_{i+1})(a_i - q^{-2} \, a_{i+1})}{
(a_i -a_{i+1})^2} \\[0.3cm]
1 \; & \; \frac{\lambda \, a_i}{a_i -a_{i+1}}
\end{array}
\!\!\right)  , \;\;
y_i =
\left(\!\!
\begin{array}{cc}
a_i & 0  \\
0 & a_{i+1}
\end{array}
\!\!\right) , \;\;
y_{i+1} =
\left(\!\!
\begin{array}{cc}
a_{i+1} & 0  \\
0 & a_i
\end{array}
\!\!\right) ,
\ee
where $a_i \neq a_{i+1}$
otherwise $y_i$, $y_{i+1}$ are not diagonalizable.
We note that the form of matrix $\sigma_i$ in (\ref{ah1})
is not unique, since one can multiply $V$ by
 any diagonal matrix $D$ from the right. For $a_{i+1} \neq q^{\pm 2} a_i$,
one can perform additional similarity transformation
of operators (\ref{ah1}) with diagonal matrix
$D = {\rm diag}(d_1,d_2)$ to make the matrix
$\sigma_i$ symmetric
$$
\sigma_i \to
D^{-1} \, \sigma_i \, D=
\left(
\begin{array}{cc}
- \frac{\lambda \, a_{i+1}}{a_i -a_{i+1}} \; & \;
\frac{d_1}{d_2} \\[0.3cm]
\frac{d_1}{d_2} \; & \; \frac{\lambda \, a_i}{a_i -a_{i+1}}
\end{array}
\right) \; , \;\;\;\; \frac{d_1}{d_2} = \pm
\frac{\sqrt{(a_i - q^2 \, a_{i+1})(a_i - q^{-2} \, a_{i+1})}}{
(a_i -a_{i+1})} \; .
$$
When $a_{i+1} = q^{\pm 2} a_i$, the
2-dimensional representation (\ref{ah1}) is reduced
into the 1-dimensional representation with $\sigma_i = \pm q^{\pm 1}$,
respectively.
We summarize the above results as following.

\begin{proposition}\label{prop10}
(q-Vershik-Okounkov \cite{OV}, \cite{IsOg3}). {\it Let
$$
\Lambda = (a_1, \dots , a_i, a_{i+1} , \dots , a_n) \in {\rm Spec}
(y_1, \dots , y_{M+1}) \; ,
$$
be a possible spectrum of the commutative set $(y_1, \dots , y_{M+1})$,
which corresponds to a primitive eigenidempotent
$e_\Lambda \in H_{M+1}$.
Then, $a_i = q^{2 m_i}$, where  $m_i \in {\bf Z}_i$ (\ref{spec}) and \\
(1) $a_i \neq a_{i+1}$ for all $i < M+1$; \\
(2) if $a_{i+1} = q^{\pm 2} a_{i}$, then
$\sigma_i \cdot e_\Lambda = \pm q^{\pm 1} e_{\Lambda}$; \\
(3) if $a_i \neq q^{\pm 2} a_{i+1}$, then
$$
\Lambda' = (a_1, \dots , a_{i+1} , a_i, \dots , a_{M+1})
\in {\rm Spec} (y_1, \dots , y_{M+1}) \; ,
$$
and the left action of the elements  $\sigma_i,y_i,y_{i+1}$
in the linear span of $v_\Lambda = e_\Lambda$ and
$$
v_{\Lambda'} = \sigma_i \, e_\Lambda +
\frac{\lambda a_{i+1}}{a_i- a_{i+1}} \, e_\Lambda \; ,
$$
is given by (\ref{ah1}).
}
\end{proposition}

From this Proposition we conclude that the only admissible subgraphs
in the Young-Ogievetsky (YO) graph
 (subgraphs which show all possible two-edges paths
with fixed initial and final vertices)  are

\unitlength=5mm
\begin{picture}(20,7.5)
\thicklines
\put(5,6.8){\vector(0,-1){2.5}}
\put(5,3.8){\vector(0,-1){2.5}}
\put(4.75,0.8){\large $\star$}
\put(4.75,6.8){\large $\star$}
\put(4.75,3.8){\large $\star$}
\put(6.5,3.9){$(\sigma_i = \pm q^{\pm 1})$}
\put(-1,2.3){$y_{i+1}=$}
\put(5.5,2.3){$a q^{\pm 2}$}
\put(-0.6,5.3){$y_i=$}
\put(5.5,5.3){$a$}

\put(7,0.4){{\bf Fig. 4.2}}

\put(13.5,2.3){$b$}
\put(13.5,5.3){$a$}
\put(16.1,5.3){$b$}
\put(16.1,2.3){$a$}
\put(18,3.9){$( b \neq q^{\pm2}a)$}
\put(13.25,3.8){\large $\star$}
\put(16.25,3.8){\large $\star$}
\put(14.8,6.9){\large $\star$}
\put(14.8,0.9){\large $\star$}
\put(15,7){\vector(-1,-2){1.4}}
\put(15,7){\vector(1,-2){1.4}}
\put(13.6,3.9){\vector(1,-2){1.3}}
\put(16.4,3.9){\vector(-1,-2){1.3}}

\put(18,0.4){{\bf Fig. 4.3}}
\end{picture}

\noindent
where stars in the vertices denote Young diagrams. These subgraphs
are related to the 1-dimensional and 2-dimensional
(it corresponds to the number of paths from the top
vertex to the bottom one in Fig.4.2 and Fig.4.3) representations of
the subalgebra generated by $\{ y_i,y_{i+1},\sigma_i \}$. In view of the
braid relations $\sigma_{i} \sigma_{i \pm 1} \sigma_i =
\sigma_{i \pm 1} \sigma_i  \sigma_{i \pm 1}$ and
possible values of $\sigma$'s presented in Fig. 4.2
for 1-dimensional (1D) representation subgraphs, we
conclude that the chains:
{\large $\star \stackrel{a}{\rightarrow} \star \stackrel{q^{\pm 2}a}{\rightarrow} \star
\stackrel{a}{\rightarrow} \star$} of the 1D representation
subgraphs in the YO graph in Fig. 4.1 are forbidden. While
admissible chains of 1D representation subgraphs are:
{\large $\star \stackrel{a}{\rightarrow} \star \stackrel{q^{\pm 2}a}{\rightarrow} \star
\stackrel{q^{\pm 4}a}{\rightarrow} \star$}.
These statements and the form
of only admissible subgraphs in Figs. 4.2,4.3
justify (for the A-type
Hecke algebra) the YO graph presented in Fig. 4.1. Indeed,
we know the top of the YO graph which consists of
3 edges with indices $1,q^2,q^{-2}$ (see Fig. 4.1).
Then, one can explicitly construct ``step by step'' moving down
the whole YO graph (with all indices on edges) by using (\ref{spec11}),
the form of only admissible subgraphs in Figs. 4.2,4.3 and rules for the chains of 1-dim.
representations. We also stress here two important properties of
the YO graph: \\
1.) to each vertex of the graph the number of
incoming edges $E_{in}$ less than the number of outgoing edges $E_{out}$ on $1$:
$E_{out}=E_{in}+1$; \\
2.) to each vertex of the graph the products of indices $a_{in}$
for incoming and $a_{out}$ for outgoing edges
are equal to each other: $\prod_{\alpha =1}^{E_{in}} a_{in}(\alpha) =
\prod_{\beta =1}^{E_{out}} a_{out}(\beta)$.

 Finally we summarize all results about the spectrum
 ${\rm Spec}(y_1, \dots,y_{M+1})$
 of the Jucys-Murphy elements $y_i \in H_{M+1}(q)$ and
 YO graph for Hecke algebra $H_{M+1}(q)$ as following.
 \begin{proposition}\label{propF}
 (\cite{OV}, \cite{IsOg3}). {\it There is a bijection between the set of
the standard Young tableaux ${\sf T}_{[\nu]_n}$  with $n$ nodes, the set
${\rm Spec} (y_1, \dots , y_n)$ and the set
of paths $X_{\vec{a}}$ in coloured YO graph which connect the
top diagram $\emptyset$ and the diagram $[\nu]_n$.}
 \end{proposition}

Since the YO graph is explicitly known, we can deduce the
expressions (in terms of the
elements $y_k$) of all orthogonal primitive idempotents
for the Hecke algebra (in the same way as it has been
done in (\ref{import3}) - (\ref{import4})).
We stress once again that the method of construction
 of explicit expressions for
such primitive orthogonal idempotents is known
and was discussed e.g. in \cite{Mur}, \cite{OgPya}, \cite{IsOg3}.
Now  we explain the operation of this method by using
an arbitrary standard Young tableau as an example.

Let ${\sf \Lambda}$ be a Young diagram with $n=n_k$ rows:
$\lambda_1 \geq \lambda_2 \geq  \dots \geq \lambda_n$ and
$|{\sf \Lambda}| := \sum_{i=1}^n \lambda_i$ be the number of its nodes.
Consider the case when
$\lambda_1 = \dots = \lambda_{n_1} = \lambda_{(1)} >
\lambda_{n_1 +1} = \lambda_{n_1 + 2} = \dots = \lambda_{n_2} = \lambda_{(2)} >
\dots > \lambda_{n_k - n_{k-1} +1} = \dots = \lambda_{|{\sf \Lambda}|} = \lambda_{(n_k)}$:

\unitlength=4mm
\begin{picture}(25,7)(-10,0)
\put(-2,3){${\sf \Lambda} \;\; = $}
\put(4,5.5){\line(1,0){5}}
\put(4,4){\line(1,0){5}}
\put(4,4){\line(0,1){1.5}}
\put(9,4){\line(0,1){1.5}}

\put(4,2.5){\line(1,0){3}}
\put(7,2.5){\line(0,1){1.5}}
\put(4,4){\line(1,0){3}}
\put(4,1){\line(0,1){3}}
\put(6,1.5){\line(0,1){1}}
\put(4,1.5){\line(1,0){2}}
\put(4.5,1){$\dots$}

\put(4,0){\line(0,1){0.5}}
\put(5,0){\line(0,1){0.5}}
\put(4,0){\line(1,0){1}}
\put(4,0.5){\line(1,0){1}}

\put(6,6.2){$_{_{\lambda_{_{(1)}}}}$}
\put(3,4.7){$_{n_{_1}}$}
\put(2,3.4){$_{_{n_2-n_1}}$}
\put(1,0.4){$_{_{n_k-n_{k-1}}}$}
\put(9.1,4){$_{n_{_1},} {_{\lambda_{_{(1)}}}}$}
\put(7.1,2.5){$_{n_{_2},} {_{\lambda_{_{(2)}}}}$}
\put(6.2,1.5){$_{n_{_3},} {_{\lambda_{_{(3)}}}}$}
\put(5.1,0){$_{n_{_k},} {_{\lambda_{_{(k)}}}}$}

\end{picture}
\vspace{-1cm}
\be
\lb{qdima01}
\ee

\noindent
Here $(n_i,\lambda_{(i)})$ are coordinates of the nodes corresponding to the inner corners of
 the diagram ${\sf \Lambda}$. Consider any
standard Young tableau ${\sf T}_{{\sf \Lambda}_M}$ of shape
(\ref{qdima01}) with $M = |{\sf \Lambda}|$ nodes.
Let $e({\sf T}_{{\sf \Lambda}_M}) \in
H_{M}$
be a primitive idempotent
corresponding to the tableau ${\sf T}_{{\sf \Lambda}_M}$.
Taking into account the branching rule implied by the
coloured Young-Ogievetsky graph
for $H_{M+1}$, we fix all possible eigenvalues
$q^{2(\lambda_{(r)} - n_{r-1})}$ $(r=1,...,k+1)$
 of the element $y_{M+1}$. Then, we
 conclude that the following identity holds
$$
e({\sf T}_{{\sf \Lambda}_M}) \prod_{r=1}^{k+1}
\left( y_{M+1} - q^{2(\lambda_{(r)} - n_{r-1})} \right) = 0 \; ,
$$
where $\lambda_{(k+1)} = n_0 = 0$.
Thus, for a new tableau ${\sf T}_{{\sf \Lambda}^j_{M+1}}$,
which is obtained by adding to the tableau
${\sf T}_{{\sf \Lambda}_{M}}$
of the shape (\ref{qdima01}) a new
node with coordinates $(n_{j-1} +1, \lambda_{(j)} +1)$,
we obtain the following primitive idempotent (after a normalization)
\be
\label{pij}
e({\sf T}_{{\sf \Lambda}^j_{M+1}}) :=
e({\sf T}_{{\sf \Lambda}_{M}})
\prod_{\stackrel{r=1}{_{r \neq j}}}^{k+1} \!
\frac{\left( y_{_{|{\sf \Lambda}|+1}}
- q^{2(\lambda_{(r)} - n_{r-1})} \right)}{\left( q^{2(\lambda_{(j)}
 - n_{j-1})} - q^{2(\lambda_{(r)} - n_{r-1})} \right)}  =
e({\sf T}_{{\sf \Lambda}_{M}}) \; \Pi_j \,  .
\ee
With the help of this formula and "initial data" $e \left(
\begin{tabular}{|c|}
\hline
$\!\!\!$  1 $\!\!\!$ \\
\hline
\end{tabular}
\right) = 1$, one can deduce step by step
explicit expressions for all primitive orthogonal idempotents for
Hecke algebras.

\vspace{0.2cm}

\noindent
{\bf Remark.} The elements found by inductive formula (\ref{pij})
give by construction a complete system of mutually orthogonal
idempotents in the Hecke algebra $H_{M}(q)$. Let
${\sf T}_{a,{\sf \Lambda}}$ $(a=1,...,f_{\sf \Lambda})$
be standard Young tableaux of the shape of the
Young diagram ${\sf \Lambda} \vdash M$ and $f_{\sf \Lambda}$
is the number of such Young tableaux of the shape ${\sf \Lambda}$.
A primitive idempotent $e({\sf T}_{a,{\sf \Lambda}}) \in H_M(q)$
corresponds  to such tableau.
 Central idempotents $e({\sf \Lambda})$ correspond to
 the Young diagrams ${\sf \Lambda} \vdash M$ and are expressed
 as the sum
 $e({\sf \Lambda}) = \sum_{a=1}^{f_{\sf \Lambda}}
 e({\sf T}_{a,{\sf \Lambda}})$. Then, the completeness of the
 primitive orthogonal idempotents $e({\sf T}_{a,{\sf \Lambda}})$
 is written as the resolution of unit operator $1$ via
 central idempotents $e({\sf \Lambda}) \in H_M(q)$:
 $$
 1 = \sum_{{\sf \Lambda} \vdash M} \sum_{a=1}^{f_{\sf \Lambda}}
 e({\sf T}_{a,{\sf \Lambda}}) =
 \sum_{{\sf \Lambda} \vdash M} e({\sf \Lambda}) \; .
 $$
 One can find explicit formula for the number $f_{\sf \Lambda}$
 in \cite{IsRub2}, Sect. 4.3.2 (see also references therein).

\subsubsection{Idempotents in
$H_{M+1}(q)$ and baxterized elements. Matrix
units in $H_{M+1}$\label{idbaxt}}

Another convenient recurrent relations for Hecke symmetrizers and antisymmetrizers
(\ref{santis}), (\ref{sant01}), (\ref{sant02})
are (see e.g. \cite{Jimb1}, \cite{Gur}, \cite{HIOPT}):
\be
\lb{santis2}
S_{1 \to n} = S_{1 \to n-1} \, \frac{\sigma_{n-1}(q^{1-n})}{ [n]_q} \, S_{1 \to n-1}
\; , \;\;\;
S_{1 \to n} = S_{2 \to n} \, \frac{\sigma_{1}(q^{1-n})}{ [n]_q} \, S_{2 \to n}
\; ,
\ee
\be
\lb{santis22}
A_{1 \to n} = A_{1 \to n-1} \, \frac{\sigma_{n-1}(q^{n-1})}{ [n]_q} \, A_{1 \to n-1}
\; , \;\;\;
A_{1 \to n} = A_{2 \to n} \, \frac{\sigma_{1}(q^{n-1})}{ [n]_q} \, A_{2 \to n}
\; ,
\ee
where $\sigma_n(x)$
are Baxterized elements \cite{Jimb1}, \cite{Jon3}
 (cf. (\ref{3.5.5}))
\be
\lb{baxtH}
\sigma_n(x) : = \lambda^{-1} \,(x^{-1} \sigma_n - x \sigma_n^{-1}) \; ,
\ee
for the algebra $H_{M+1}(q)$. We have already used these elements in
definitions of the idempotents (\ref{proj1}) and (\ref{proj2}). The elements
$\sigma_n(x)$ obey the Yang-Baxter equation (the proof of this statement
is the same as in (\ref{3.5.1}) -- (\ref{3.5.4})):
\be
\lb{ybeH}
\sigma_n(x) \, \sigma_{n-1}(xy) \, \sigma_n(y) =
\sigma_{n-1}(y) \, \sigma_n(xy) \, \sigma_{n-1}(x)  \; .
\ee
These elements are also normalized by the conditions
 $\sigma_n(\pm 1) = \pm 1$ and satisfy
$$
\sigma_i(x) = \frac{x - x^{-1}}{y-y^{-1}} \,
\sigma_i(y) + \frac{yx^{-1} - xy^{-1}}{y-y^{-1}} \; ,
\;\;\;\;\;\; \forall x,y \neq \pm 1 \; ,
$$
\be
\lb{hunit1}
\sigma_i (x) \, \sigma_i(y) =
\sigma_i (xy) + (x-x^{-1})(y-y^{-1}) \lambda^{-2} \; .
\ee
The special
case $y=x^{-1}$ of (\ref{hunit1}) gives
 the ``unitarity condition''
\be
\lb{hunit}
\sigma_i (x) \, \sigma_i(x^{-1}) = \left( 1 - \frac{(x-x^{-1})^2}{\lambda^2} \right) \; \equiv \;
\frac{(q x^{-1} - q^{-1} x)(q x-q^{-1} x^{-1})}{\lambda^2}   \; .
\ee
One can write the Baxterized elements (\ref{baxtH}) as
rational function of $\sigma_i$ (cf. (\ref{3.5.5b}))
\be
\lb{hunit2}
\sigma_i(x) = \left(\frac{a^{-1}x - a x^{-1}}{\lambda x^2} \right)
\; \frac{\sigma_i - a x^2}{\sigma_i - a x^{-2}}\;\; ,
 \;\;\;\;\;\;\;\;  \sigma_i^{(a)}(x)  \; :\, =  \;
\frac{\sigma_i - a x^2}{\sigma_i - a x^{-2}}\; ,
\ee
where $a=-q$, or $a=q^{-1}$,
and it becomes clear, for both choices of $a$,
that the normalized elements $\sigma_i^{(a)} (x)$
 obey $\sigma_i^{(a)} (x) \sigma_i^{(a)} (x^{-1})=1$.

Eqs. (\ref{hunit1}), (\ref{hunit}) and
(\ref{hunit2}) follow from the Hecke relation (\ref{ahecke}).
The equivalence of both representations for (anti)symmetrizers
given in the first and second equations of
(\ref{santis2}) and (\ref{santis22}) is demonstrated by means
of Yang-Baxter eq. (\ref{ybeH}), or by means of the obvious
mirror automorphism
$\sigma_k \rightarrow \sigma_{M-k}$ for the Hecke algebra $H_{M}(q)$.
The equivalence of (\ref{santis2}), (\ref{santis22}) and
(\ref{santis}) can be easily demonstrated if one writes
first representations of (\ref{santis2}), (\ref{santis22}) in the form
\be
\lb{santis3}
\begin{array}{c}
S_{1 \to n} = \frac{1}{ [n]_q!}  \, \sigma_{1}(q^{-1})  \, \sigma_{2}(q^{-2})
 \dots  \sigma_{n-1}(q^{1-n}) \, S_{1 \to n-1}
\; , \\ \\
A_{1 \to n} =
\frac{1}{ [n]_q!}  \, \sigma_{1}(q)  \, \sigma_{2}(q^{2})
 \dots  \sigma_{n-1}(q^{n-1}) \, A_{1 \to n-1} \; ,
\end{array}
\ee
and, then, we should
 use (\ref{santis1}) to compare (\ref{santis3}) with (\ref{santis})
and (\ref{anti}).
According to (\ref{santis2}), (\ref{santis22}) the first two
projectors are (cf. (\ref{3.3.11}), (\ref{3.3.111})):
$$
P \left(
\begin{tabular}{|c|c|}
\hline
$\!\!\!$ {\small 1} $\!\!\!$ &
$\!\!\!$ {\small 2} $\!\!\!$\\
\hline
\end{tabular}
\right)
= S_{12} =
\frac{1}{ [2]_q} \, \sigma_1(q^{-1}) \; , \;\;\;
P \left(
\begin{tabular}{|c|}
\hline
$\!\!\!$ {\small 1} $\!\!\!$\\
\hline
$\!\!\!$ {\small 2} $\!\!\!$\\
\hline
\end{tabular}
\right)
= A_{12} =
\frac{1}{ [2]_q} \, \sigma_1(q) \; ,
$$
and their orthogonality readily follows from
the Hecke condition (\ref{ahecke}), or from (\ref{hunit}).
One can also express another
types of the orthogonal idempotents (not only symmetrizers and antisymmetrizers)
in terms of the Baxterized elements:
$$
P \left(
\begin{tabular}{|c|c|}
\hline
$\!\!\!$  1 $\!\!\!$ & $\!\!\!$  2 $\!\!\!$ \\
\hline
$\!\!\!$  3 $\!\!\!$ & \multicolumn{1}{c}{}\\
\cline{1-1}
\end{tabular}
\right) =  \frac{1}{[3]_q !}
\sigma_1(q^{-1}) \sigma_2(q) \sigma_1(q^{-1})  \; , \;\;\;
P \left(
\begin{tabular}{|c|c|}
\hline
$\!\!\!$ {\small 1} $\!\!\!$ & $\!\!\!$ {\small 3} $\!\!\!$ \\
\hline
$\!\!\!$ {\small 2} $\!\!\!$ & \multicolumn{1}{c}{}\\
\cline{1-1}
\end{tabular}
\right) =  \frac{1}{[3]_q !}
\sigma_1(q) \sigma_2(q^{-1}) \sigma_1(q) \; ,
$$
$$
P \left(
\begin{tabular}{|c|c|c|}
\hline
$\!\!\!$ {\small 1} $\!\!\!$ & $\!\!\!$ {\small 2} $\!\!\!$ &
$\!\!\!$ {\small 3} $\!\!\!$\\
\hline
$\!\!\!$ {\small 4} $\!\!\!$ &  \multicolumn{1}{c}{}\\
\cline{1-1}
\end{tabular}
\right) \sim  \sigma_1(q^{-1}) \sigma_2(q^{-2}) \sigma_1(q^{-1}) \sigma_3(q)
\sigma_2(q^{-2}) \sigma_1(q^{-1}) \; ,
$$
$$
P \left(
\begin{tabular}{|c|c|c|}
\hline
$\!\!\!$ {\small 1} $\!\!\!$ & $\!\!\!$ {\small 3} $\!\!\!$ &
$\!\!\!$ {\small 4} $\!\!\!$\\
\hline
$\!\!\!$ {\small 2} $\!\!\!$ &  \multicolumn{1}{c}{}\\
\cline{1-1}
\end{tabular}
\right) \sim \sigma_1(q) \sigma_2(q^{-1}) \sigma_3(q^{-2})
\sigma_2(q^{-1}) \sigma_1(q)  \; ,
$$
$$
P \left(
\begin{tabular}{|c|c|c|}
\hline
$\!\!\!$ {\small 1} $\!\!\!$ & $\!\!\!$ {\small 2} $\!\!\!$ &
$\!\!\!$ {\small 4} $\!\!\!$\\
\hline
$\!\!\!$ {\small 3} $\!\!\!$ &  \multicolumn{1}{c}{}\\
\cline{1-1}
\end{tabular}
\right) \sim  \sigma_1(q^{-1}) \sigma_2(q) \sigma_1(q^{-1}) \sigma_3(q)
\sigma_2(q^{-2}) \sigma_1(q^{-1}) \sigma_3(q^{-3}) \; ,
$$
$$
P \left(
\begin{tabular}{|c|c|}
\hline
$\!\!\!$ {\small 1} $\!\!\!$ & $\!\!\!$ {\small 2} $\!\!\!$ \\
\hline
$\!\!\!$ {\small 3} $\!\!\!$ & $\!\!\!$ {\small 4} $\!\!\!$\\
\cline{1-2}
\end{tabular}
\right) \sim  \sigma_1(q^{-1}) \sigma_2(q) \sigma_1(q^{-1}) \sigma_3(q^3)
\sigma_2(q) \sigma_1(q^{-1}) \sigma_3(q^{-1}) \; .
$$
The method of presentation of all primitive orthogonal idempotents for Hecke algebra in terms of the
Baxterized elements was developed in \cite{IsMoOs}
(see also refs. therein) by means of the fusion procedure.

\vspace{0.1cm}

\noindent
{\bf Remark 1.} Consider the quotients of the Hecke algebra
$H_{M+1}(q)$ with respect to
the additional relations $A_{1 \to n} = 0$ ($n \leq M+1$), which are equivalent
(see (\ref{santis22})) to the identities
\be
\lb{asaa}
A_{1 \to n-1} \, \sigma_{n-1} \, A_{1 \to n-1} =
\frac{q^{n-1}}{[n-1]_q} \, A_{1 \to n-1} \; .
\ee
This is the way how the generalized Temperley-Lieb-Martin algebras \cite{Mart}
are defined. As it was mentioned in \cite{Jimb1},
the quotient of $H_{M+1}(q)$ with respect to
 the identity $A_{1 \to 3} =0$
is isomorphic to the Temperley-Lieb algebra.

\vspace{0.2cm}
\noindent
{\bf Remark 2.}
By using intertwining elements (\ref{impint}) and
Baxterized element (\ref{baxtH}) one can immediately
construct off-diagonal matrix
units\footnote{Recall that the orthogonal primitive idempotents are
diagonal matrix units.} (see Sect. {\bf \ref{jmel}}) in a
double sided Peirce decomposition of the Hecke algebra
$H_{M+1}= \bigoplus_{\alpha,\beta} e_\alpha H_{M+1} e_\beta$.
Let $P(X_{\vec{a}}):=e({\sf T}_{{\sf \Lambda}_{M+1}})$
be orthogonal primitive idempotent which corresponds
to the Young tableau ${\sf T}_{{\sf \Lambda}_{M+1}}$ or
to the path $X_{\vec{a}}$ on the coloured Young-Ogievetsky
 graph labeled by the eigenvalues
 \be
 \lb{spec5}
 \begin{array}{c}
\vec{a} = (1, a_2, \dots , a_{M+1}) \in {\rm Spec}(y_1, \dots,y_{M+1}) \; ,
\\ [0.3cm]
y_i \, P(X_{\vec{a}}) = P(X_{\vec{a}}) \, y_i = a_i \, P(X_{\vec{a}}) \;\;
(\forall i = 1, \dots , M+1) \; .
 \end{array}
 \ee
In the case $a_j \neq q^{\pm 2} a_{j+1}$
(see the proof of Lemma 2 in Subsect. {\bf \ref{idemp}})
 we introduce the element
$P(X_{s_j \cdot \vec{a}}) \in H_{M+1}$
\be
\lb{sX}
P(X_{s_j \cdot \vec{a}}) := \frac{1}{(q^2 a_j - a_{j+1})(a_{j+1} - q^{-2} a_j)} \,
U_{j+1} \, P(X_{\vec{a}}) \, U_{j+1} \;\;\;  (\forall j = 1, \dots, M)
\ee
such that
$$
P(X_{s_j \cdot \vec{a}})^2 = P(X_{s_j \cdot \vec{a}}) \; , \;\;\;
\vec{y} \, P(X_{s_j \cdot \vec{a}})  = P(X_{s_j \cdot \vec{a}}) \, \vec{y}  =
(s_j \cdot \vec{a}) P(X_{s_j \cdot \vec{a}}) \; ,
$$
$$
P(X_{s_j \cdot \vec{a}}) =
\frac{(q^2 y_j - y_{j+1})(y_{j+1} - q^{-2} y_j)}{(q^2 a_j - a_{j+1})(a_{j+1} - q^{-2} a_j)} \;
(s_{j} \! \cdot \! P)(X_{\vec{a}})  \; ,
$$
\be
\lb{PsP}
P(X_{s_j \cdot \vec{a}}) \, P(X_{\vec{a}}) =
P(X_{\vec{a}})  \, P(X_{s_j \cdot \vec{a}}) = 0 \; ,
\ee
where $s_j \cdot \vec{a} = (a_1, \dots, a_{j+1},a_j, \dots a_{M+1})
\in {\rm Spec}(y_1, \dots, y_{M+1})$
is the vector with permuted coordinates $a_j$ and $a_{j+1}$;
$(s_{j} \! \cdot \! P)(X_{\vec{a}})$ denotes the function $P(X_{\vec{a}})$
with permuted variables $y_i$ and $y_{i+1}$.
 The identity (\ref{PsP}) follows
from the fact that
 $P := P(X_{\vec{a}})  \, \, P(X_{\vec{a} \, '})  \, =0$ for all
$\vec{a} \neq \vec{a} \, '$ (i.e., $\exists \, j$: $a_j \neq a'_j$)
in view of the equations $y_j P = a_j P = a'_j P$
which follow from $y_j \, P(X_{\vec{a}}) = P(X_{\vec{a}})\, y_j$.

According to (\ref{import}) and (\ref{sX}) we obtain
\be
\lb{offdg}
\begin{array}{c}
U_{j+1} \, P(X_{\vec{a}}) = P(X_{s_j \cdot \vec{a}}) \, U_{j+1}  =:
P(X_{s_j \cdot \vec{a}} | X_{\vec{a}} ) \;\; (j=1, \dots, M) \; ,
\\ [0.3cm]
 P(X_{\vec{a}}) \, U_{j+1}  =  U_{j+1} \, P(X_{s_j \cdot \vec{a}}) =:
P(X_{\vec{a}} |  X_{s_j \cdot \vec{a}} ) \;\; (j=1, \dots, M) \; ,
\end{array}
\ee
In the case $a_j \neq q^{\pm 2} a_{j+1}$, in view of Lemma 2,
we have $s_j \cdot \vec{a} \in {\rm Spec}(y_1, \dots , y_{M+1})$ (the path
$X_{s_j \cdot \vec{a}}$ exists in Young-Ogievetsky
 graph and corresponds to the standard Young tableau).
Then, taking into account (\ref{impint}), (\ref{impint1}), we deduce
\be
\lb{offdg2}
\left.
\begin{array}{c}
\sigma_j(a_j,a_{j+1}) \; P(X_{\vec{a}}) =
- P(X_{s_j \cdot \vec{a}}) \; \sigma_j(a_{j+1},a_j)  =
P(X_{s_j \cdot \vec{a}} | X_{\vec{a}} )
 \; , \\ [0.3cm]
 - P(X_{\vec{a}}) \; \sigma_j(a_j,a_{j+1})   =
 \sigma_j(a_{j+1},a_j)  \; P(X_{s_j \cdot \vec{a}}) =
 P( X_{\vec{a}} | X_{s_j \cdot \vec{a}} ) \; , 
\end{array} \right\} \;\;\;\; \Rightarrow
\ee
\be
\lb{offdg1}
\begin{array}{c}
(\sigma_j  + \frac{\lambda a_{j+1}}{(a_j - a_{j+1})}) \, P(X_{\vec{a}}) =
P(X_{s_j \cdot \vec{a}}) \, (\sigma_j - \frac{\lambda a_j}{(a_j - a_{j+1})})  =
\frac{P(X_{s_j \cdot \vec{a}} | X_{\vec{a}} )}{(a_j - a_{j+1})}
 \; , \\ [0.3cm]
 P(X_{\vec{a}}) \, (\sigma_j  + \frac{\lambda a_{j+1}}{(a_{j} - a_{j+1})})  =
(\sigma_j  - \frac{\lambda a_{j}}{(a_{j} - a_{j+1})}) \, P(X_{s_j \cdot \vec{a}}) =
\frac{P( X_{\vec{a}} | X_{s_j \cdot \vec{a}} )}{(a_{j+1} - a_{j})} \; ,
\end{array}
\ee
 where we used the Baxterized elements (cf. (\ref{baxtH}))
 \be
 \lb{offdg3}
 \sigma_n (x,y) = x\, \sigma_n  - y \, \sigma_n^{-1} =
 (x-y)\, \sigma_n + y \, \lambda \; ,
 \ee
 subject the Yang-Baxter equation (\ref{ybeH})
 written in the form
 $$
 \sigma_n (x,y) \, \sigma_{n+1} (x,z) \, \sigma_n (y,z) =
 \sigma_{n+1} (y,z) \, \sigma_n (x,z) \, \sigma_{n+1} (x,y) \; .
 $$
 It was shown in \cite{OgPya} that
the elements $P( X_{s_j \cdot \vec{a}} | X_{\vec{a}} )$ play the role
of the off-diagonal matrix elements in the Peirce decomposition
(see Subsection {\bf \ref{jmel}}).
In the case  $a_j = q^{\pm 2} a_{j+1}$ we have
\be
\lb{pbp7}
 \begin{array}{c}
U_{j+1} \, P(X_{\vec{a}}) = 0 = P(X_{\vec{a}}) \, U_{j+1}
\;\;\; \Rightarrow \\ [0.3cm]
(\sigma_j \pm q^{\mp 1})\, P(X_{\vec{a}}) = 0 =
- P(X_{\vec{a}}) \, (\sigma_j \pm q^{\mp 1}) \; ,
\end{array}
\ee
where the second line in (\ref{pbp7})
 follow from (\ref{offdg1}) and
define the 1-dimensional
representation for the generator $\sigma_j$ corresponding
 to the Fig. 4.2 in Subsection {\bf \ref{jmel2}}.

Since the Hecke algebras $H_{M+1}$ are semisimple,
we have the following identity
 \be
 \lb{pbp}
P(X_{\vec{a}}) \, B \, P(X_{\vec{a}}) =
C_{\vec{a}}(B) \, P(X_{\vec{a}}) \; , \;\;\;\;\;\;
\forall B \in H_{M+1} \; ,
 \ee
where $P(X_{\vec{a}})$ is any primitive
 idempotent in $H_{M+1}$
 and $C_{\vec{a}}(B)$ is a constant which depends on
the element $B$ and the path $X_{\vec{a}}$
 in the coloured Young-Ogievetsky graph (i.e. it depends on
  the vector $\vec{a}$ from the spectrum (\ref{spec5})).
  To justify identity (\ref{pbp}) we check it
  for any monomial
  $B = \sigma_{i_1} \,  \sigma_{i_2} \dots \sigma_{i_r}$
  in generators $\sigma_i \in H_{M+1}$ of any order $r$.
  We require that the monomial $B$ can not be reduced to
  the polynomial of order less than $r$ by means of
  relations (\ref{braidg}) and (\ref{ahecke}). Then, we use the induction to prove (\ref{pbp}). We note that by using
  the definition (\ref{impint}) of $U_{i+1}$
  and then (\ref{impint1}) we obtain
  the base of the induction
  $$
  \begin{array}{c}
  0= P(X_{\vec{a}}) \, U_{i+1} \, P(X_{\vec{a}}) =
  P(X_{\vec{a}}) \,
  ((y_{i+1} - y_i) \sigma_i - \lambda y_{i+1}) \, P(X_{\vec{a}})
    =  \\ [0.2cm]
  = \displaystyle (a_{i+1} - a_i)P(X_{\vec{a}}) \,
  \sigma_i \, P(X_{\vec{a}}) -
  \lambda a_{i+1} \, P(X_{\vec{a}}) \;\;\; \Rightarrow \;\;\;
  P(X_{\vec{a}}) \,  \sigma_i \, P(X_{\vec{a}}) =
  \frac{\lambda a_{i+1}}{(a_{i+1} - a_i)} P(X_{\vec{a}}) \; .
  \end{array}
  $$
  Let the identity (\ref{pbp}) be correct for all
  monomials $B = \sigma_{i_1} \cdots \sigma_{i_k}$ when $k \leq r$.
  We need to prove (\ref{pbp}) for
  monomials $B = \sigma_{i_1} \cdots \sigma_{i_{r+1}}$ of order
  $(r+1)$. Consider the element
  $P(X_{\vec{a}}) \, U_{i_1+1} \cdots U_{i_{r+1}+1}
  \, P(X_{\vec{a}})$ and start to commute left idempotent
  $P(X_{\vec{a}})$ to the right with the help
  of (\ref{offdg}). We have two possibilities.

  \noindent
  {\bf 1.} The first one is
  \be
  \lb{pbp4}
  \begin{array}{c}
  P(X_{\vec{a}}) \, U_{i_1+1} \cdots U_{i_{r+1}+1}
  \, P(X_{\vec{a}}) = U_{i_1+1} \cdots
  U_{i_k+1}  P(X_{\vec{a}^{(k)}}) U_{i_{k+1}+1} \cdots\,
   P(X_{\vec{a}}) = \\ [0.3cm]
  = \dots = U_{i_1+1} \cdots \cdots U_{i_{r+1}+1} \,
   P(X_{\vec{a}^{(r+1)}}) \; P(X_{\vec{a}}) = 0 \; ,
      \end{array}
  \ee
  $$
 \vec{a}^{(k)} :=s_{i_{k}} \cdot \vec{a}^{(k-1)} \, , \;\;\;
  \vec{a}^{(0)} := \vec{a}  \, , \;\;\;
  s_{i} \cdot (v_1,...,v_{i},v_{i+1},...) =
  (v_1,...,v_{i+1},v_{i},....) \, ,
 $$
  where
  $(\vec{a}^{(k)})_{i_{k+1}} \neq
  q^{\pm 2}(\vec{a}^{(k)})_{i_{k+1}+1}$ ($\forall k=0,...,r$) and
  we used  the orthogonality $P(X_{\vec{a}'}) \,
  \cdot P(X_{\vec{a}}) =0$ for $\vec{a}\,' \neq \vec{a}$
  in the last equality in (\ref{pbp4}). In this case
  by using (\ref{offdg2}) we deduce
 $$
  \begin{array}{c}
  0= P(X_{\vec{a}}) \, U_{i_1+1} \cdots U_{i_{r+1}+1}
  \, P(X_{\vec{a}})
   = U_{i_1+1} \cdots U_{i_{_r}+1} \,
   P(X_{\vec{a}^{(r)}}) \;
   U_{i_{_{r+1}}+1} \; P(X_{\vec{a}})
    =  \\ [0.2cm]
  = (-1) \; U_{i_1+1} \cdots U_{i_{r-1}+1} \,
   P(X_{\vec{a}^{(r-1)}}) \;
   U_{i_{r}} \sigma_{i_{r+1}}
   ( a^{(r)}_{i_{r+1}},a^{(r)}_{i_{_{r+1}}+1}) \; P(X_{\vec{a}})  =
   \dots =
  \end{array}
  $$
  \be
  \lb{pbp2}
  = (-1)^{r+1} P(X_{\vec{a}})
  \; \sigma_{i_{1}}( a^{(0)}_{i_{1}},a^{(0)}_{i_{_{1}}+1}) \cdots
   \sigma_{i_{r}}( a^{(r-1)}_{i_{r}},a^{(r-1)}_{i_{_{r}}+1}) \;
    \sigma_{i_{r+1}}( a^{(r)}_{i_{r+1}},a^{(r)}_{i_{_{r+1}}+1}) \;
    \; P(X_{\vec{a}}) \; ,
  \ee
  where $\sigma_i(x,y)$ are Baxterized
  elements (\ref{offdg3}). The substitution
  of the r.h.s. of (\ref{offdg3}) gives
 \be
 \lb{pbp3}
P(X_{\vec{a}}) \; \sigma_{i_{1}} \cdots \sigma_{i_{r+1}}
\; P(X_{\vec{a}}) = \frac{1}{
(a^{(0)}_{i_{_{1}}+1} - a^{(0)}_{i_{1}})
\cdots (a^{(r)}_{i_{_{r+1}}+1}- a^{(r)}_{i_{r+1}})}
P(X_{\vec{a}}) \; \overline{B} \; P(X_{\vec{a}})
 \ee
  where $\overline{B}$ is a polynomial in $\sigma_i \in H_{M+1}$
  of order less than $(r+1)$ and therefore in view
  of the induction conjecture we obtain (\ref{pbp}).

    \noindent
  {\bf 2.} The second possibility occurs if at some step $k$
  in (\ref{pbp4}) the condition $(\vec{a}^{(k)})_{i_{k+1}} =
  q^{\pm 2}(\vec{a}^{(k)})_{i_{k+1}+1}$ arises. In this case,
  in view of (\ref{pbp7}), we obtain for any element
  $A \in H_{M+1}$
   \be
  \lb{pbp5}
  \begin{array}{c}
  P(X_{\vec{a}}) \, U_{i_1+1} \cdots U_{i_{k+1}+1} \, A
  \, P(X_{\vec{a}}) = U_{i_1+1} \cdots
  U_{i_k+1}  P(X_{\vec{a}^{(k)}}) U_{i_{k+1}+1} \, A \,
   P(X_{\vec{a}}) = 0 \; ,
      \end{array}
  \ee
  We take here $A = \sigma_{i_{k+1}+1} \cdots \sigma_{i_{r+1}}$
  and write (\ref{pbp5}) with the help of (\ref{offdg2}) and (\ref{pbp7}) (in the same way as in (\ref{pbp2})) in the form
  $$
   \begin{array}{c}
  0=P(X_{\vec{a}}) \, U_{i_1+1} \cdots U_{i_{k+1}+1} \, \sigma_{i_{k+1}+1} \cdots \sigma_{i_{r+1}} \, P(X_{\vec{a}}) =
  \\ [0.2cm]
 = (-1)^{k-1}P(X_{\vec{a}})\,\sigma_{i_{1}}( a^{(0)}_{i_{1}},a^{(0)}_{i_{_{1}}+1}) \cdots
   \sigma_{i_{k-1}}( a^{(k-2)}_{i_{k-1}},a^{(k-2)}_{i_{_{k-1}}+1})
   (\sigma_{i_{k}} \pm q^{\mp1})\,
 \sigma_{i_{k+1}+1} \cdots \sigma_{i_{r+1}} \,
 P(X_{\vec{a}}) \;\;\; \Rightarrow
   \end{array}
  $$
 \be
 \lb{pbp6}
 P(X_{\vec{a}}) \; \sigma_{i_{1}} \cdots \sigma_{i_{r+1}}
\; P(X_{\vec{a}}) = \frac{1}{
(a^{(0)}_{i_{_{1}}+1} - a^{(0)}_{i_{1}})
\cdots (a^{(k-1)}_{i_{_{k}}+1}- a^{(k-1)}_{i_{k}})}
P(X_{\vec{a}}) \; \overline{B}_{\pm} \; P(X_{\vec{a}}) \; ,
 \ee
  where $\overline{B}_{\pm}$ are polynomials in
  $\sigma_i \in H_{M+1}$ of order less than $(r+1)$
  and thus, in view
  of the induction conjecture, we again prove (\ref{pbp}).

  Finally we note that equations
  (\ref{pbp2}), (\ref{pbp3}) and (\ref{pbp6})
  give us the possibility to
  calculate explicitly the constant $C_{\vec{a}}(B)$ in
  (\ref{pbp}).

\subsubsection{Affine Hecke algebras and reflection equation\label{ahare}}

{\bf 1.} In this subsection we consider
the  infinite dimensional Hecke algebra,
which corresponds to the
affine $A^{(1)}$-type Coxeter graph (\ref{A1}),
with generators $\sigma_i$ $(i=1,...,M)$ subject relations (\ref{ahecke}), (\ref{affbra}).
Thus, this algebra is the quotient of the algebra
$\mathbb{C}[{\cal B}_{M+1}(A^{(1)})]$ with respect to additional
Hecke relations (\ref{ahecke}).
We call this algebra a periodic $A$-type Hecke algebra\footnote{As we will see
below, in Sec. 5.1,  this algebra appears in a formulation of
the integrable periodic spin chain models.}
and denote it as $AH_{M+1}$. For the algebra $AH_{M+1}$ one can construct
the set of $(M-1)$ commuting elements
\be
\lb{integr}
I_k = \sum_{i=1}^M \, \sigma_i \sigma_{i+1} \dots \sigma_{i+k}
\;\;\;\;\;\;\;\;\;
(k = 0, \dots , M-2) \; ,
\ee
where we have identified $\sigma_{M+i} = \sigma_i$. The first
element in (\ref{integr}) is
$I_0 = \sum_{i=1}^M \sigma_i$ and, in the $R$-matrix representation,
this element gives a Hamiltonian for periodic spin chain (see
(\ref{hami4}) in subsection {\bf \ref{perspc}}).

Let $\{\sigma_1 , \dots , \sigma_{M-1} \}$ be generators of
the braid group ${\cal B}_{M}$.
We extend the group ${\cal B}_{M}$ by the element $X$ such that
\be
\lb{refl1}
X \, \sigma_k \, X^{-1} = \sigma_{k-1}
\;\;\; (\forall k =2,\dots,M-1) \; ,
\ee
\be
\lb{refl1a}
X \, \sigma_1 \, X^{-1} = X^{-1} \, \sigma_{M-1} \, X =: \sigma_M \; .
\ee
It is not hard to check that the new element $\sigma_M$ satisfies
eqs. (\ref{refl3}) and, therefore,
the elements $\{ \sigma_1 , \dots , \sigma_M \}$
(where $\sigma_M$ has been defined in (\ref{refl1a}))
generate the periodic braid group $\overline{{\cal B}}_{M}=
{\cal B}_{M}(A^{(1)})$.

Note that $X$ (\ref{refl1}), (\ref{refl1a}) can be realized as
the inner element of ${\cal B}_{M} \subset {\cal B}_{M}(A^{(1)})$.
 Indeed, the operator $X$
which solves equations (\ref{refl1}), (\ref{refl1a}) can be taken in the form
$X = \sigma_{M-1 \leftarrow 1} \in {\cal B}_M$, where the notation
$\sigma_{M-1 \leftarrow m} := \sigma_{M-1} \dots \sigma_{m+1} \, \sigma_m$
has been used.
Then, we define the additional generator $\sigma_M$ (\ref{refl1a}) as:
\be
\lb{sigM}
\sigma_M := X \, \sigma_1 \, X^{-1}
= \sigma_{M-1 \leftarrow 1} \, \sigma_1 \, \sigma_{M-1 \leftarrow 1}^{-1} =
\sigma_{M-1 \leftarrow 1} \, \sigma_{M-1 \leftarrow 2}^{-1}
\ee
and its graphical representation is:

\unitlength=7.5mm
\begin{picture}(17,4)

\put(0.5,1.9){$\sigma_M \;\;  =$}

\put(3,1){$\bullet$}
\put(4.7,1){$\bullet$}

\put(5.7,1){$\dots$}
\put(6.5,1){$\bullet$}
\put(9.5,1){$\bullet$}
\put(8,1){$\dots$}
\put(11.5,1){$\bullet$}

\put(3,3){$\bullet$}
\put(4.7,3){$\bullet$}

\put(5.7,3){$\dots$}
\put(6.5,3){$\bullet$}
\put(9.5,3){$\bullet$}
\put(8,3){$\dots$}
\put(11.5,3){$\bullet$}

\put(3,3.4){$_1$}
\put(5,3.4){$_2$}
\put(6.8,3.4){$_i$}
\put(9,3.4){$_{M-1}$}
\put(11.7,3.4){$_{M}$}

\put(3.15,3.2){\line(4,-1){1.5}}
 \put(4.99,2.74){\line(4,-1){1.5}}
\put(6.75,2.3){\line(4,-1){2.8}}
\put(9.85,1.55){\vector(4,-1){1.7}}

\put(11.65,3.2){\line(-4,-1){1.8}}
\put(9.45,2.65){\line(-4,-1){1.8}}
\put(7.17,2.08){\line(-4,-1){0.5}}
\put(6.45,1.9){\line(-4,-1){1.4}}
\put(4.7,1.5){\vector(-4,-1){1.45}}

 \put(4.85,3.15){\vector(0,-1){1.9}}

 \put(6.6,3.15){\vector(0,-1){1.9}}

 \put(9.65,3.15){\vector(0,-1){1.9}}

\end{picture}
\vspace{-1cm}
\be
\lb{figasp}
{}
\ee

\noindent
It is evident that $\sigma_M$ satisfies (\ref{refl3}) in view of its
graphical representation (\ref{figasp}). According to
eqs. (\ref{braidg}) and (\ref{refl3})
the elements $\{ \sigma_1 , \dots , \sigma_M \}$ of the group ${\cal B}_{M}$
(where $\sigma_M$ has been defined in (\ref{sigM}))
generate the periodic braid group $\overline{{\cal B}}_{M}$
and, therefore, eq. (\ref{sigM}) defines the homomorphism:
$\overline{{\cal B}}_{M} \rightarrow {\cal B}_{M}$.


\vspace{0.1cm}

\noindent
{\bf 2.} Another infinite dimensional Hecke algebra is
an affine Hecke algebra $\hat{H}_{M+1}(q)$.
We recall that the affine
Hecke algebra $\hat{H}_{M+1}(q)$ is defined\footnote{This algebra
is isomorphic to the quotient of the braid group algebra
 $\mathbb{C}[{\cal B}_{M+1}(C)]$, where the generators
  $\sigma_i \in {\cal B}_{M+1}(C)$ ($i=1,...,M$, $i \neq 0$) are constrained
  by additional Hecke conditions (\ref{ahecke}). The definition
  of ${\cal B}_{M+1}(C)$ is given in Sect. {\bf \ref{gabg1}}
  and is related to the Coxeter graph (\ref{CCC}).}
(see e.g. \cite{ChPr}, Chapter 12.3) as algebra generated by
elements $\sigma_i$ $(i=1, \dots , M)$
of $H_{M+1}(q)$ and additional generators $y_k$ $(k=1, \dots , M+1)$
subject relations (cf. (\ref{combr}), (\ref{refla})):
\be
\lb{afheck}
y_{k+1} = \sigma_k \, y_k \, \sigma_k \; , \;\;\; y_k \, y_j = y_j \, y_k \; , \;\;\;
 y_j \, \sigma_i  =  \sigma_i \, y_j \;\; (j \neq i,i+1) \; ,
\ee
(the generators $\{ y_k \}$ form a commutative subalgebra in $\hat{H}_{M+1}(q)$). We note that, in view of the first relation
in (\ref{afheck}), the minimal set of
generators of $\hat{H}_{M+1}(q)$ is $\{\sigma_1,...,\sigma_M,y_1 \}$.
Symmetric functions of the elements $y_k$ generate the center of the algebra
$\hat{H}_{M+1}(q)$. Below, to be short, we omit the parameter $q$
in the notation $H_{M+1}(q)$ and $\hat{H}_{M+1}(q)$.
 The interesting property of the algebra $\hat{H}_{M+1}$
is the existence of the important intertwining elements \cite{Isaev1}
(cf. (\ref{impint}) and elements $\phi_i$ in \cite{Cher5}, Prop. 3.1):
$$
U_{i+1} = (\sigma_i y_i - y_i \sigma_i) \; f(y_i,y_{i+1})
 \; , \;\;\;\;\; (1 \leq i \leq M) \; ,
$$
 where function $f$ satisfies:
$f(y_i,y_{i+1}) f(y_{i+1},y_{i})=1$. The elements $U_{i}$ obey
the  same relations (\ref{importt})--(\ref{import})
as in the case of the
non-affine A-type Hecke algebra $H_{M+1}$.


 Now we describe the procedure how one can construct $(M+1)$-dimensional
representation for the Hecke
 subalgebra $H_{M+1} \subset \hat{H}_{M+1}$.
 Let $v$ be a vector
in the space of 1-dimensional representation of
 $H_{M+1}$ such that
$\sigma_i \, v = q \, v$ $(\forall i =1, \dots, M)$.
 Consider the induced
$(M+1)$-dimensional space with the basis
 $\{v_1,v_2, \dots, v_{M+1} \}$,
  where $v_k := y_k v$. Then, according to (\ref{afheck})
and Hecke condition (\ref{ahecke})
we obtain $(M+1)$-dimensional representation for generators $\sigma_i$:
$$
\sigma_i \, v_k = q \, v_k  \;\; (k \neq i,i+1) \; , \;\;\;
\sigma_i \, v_i = q^{-1} \, v_{i+1} \; , \;\;\;
\sigma_i \, v_{i+1} = \lambda \, v_{i+1} + q \, v_{i} \; ,
$$
which is called {\it Burau representation} of $H_{M+1}$.
The matrix form of this representation is
 \be
 \lb{bura}
\sigma_i = {\rm diag} \Bigl( \underbrace{q,\dots,q,}_{i-1}
\left(\!\!
\begin{array}{cc}
0 & q \\
1/q & \lambda
\end{array}
\!\!\right),\underbrace{q,\dots,q}_{M-i} \Bigr) \; .
 \ee
  One can start from another possible one-dimensional representation
  $\sigma_i \, v = - q^{-1} \, v$ $(\forall i)$ of
  $H_{M+1}$, which leads to a new Burau representation
  resulting from (\ref{bura}) by replacing $q \to (-1/q)$.

The affine Hecke algebra $\hat{H}_{M+1}$
 makes it possible to formulate the universal
 Baxterized solution of the reflection equation
 (see (\ref{rea14}), (\ref{4.8}) below)
\be
\lb{reflH}
\sigma_n(x \, z^{-1}) \, K_n(x) \, \sigma_{n}(x \, z) \, K_n(z) =
K_{n}(z) \, \sigma_n(x \, z) \, K_n(x) \, \sigma_{n}(x \, z^{-1})  \; ,
\ee
where $x,z$ are spectral parameters and
Baxterized elements $\sigma_n(x) \in H_{M+1}$ are defined in (\ref{baxtH}).
 The reflection equation (\ref{reflH})
 appears, e.g., in the theory of integrable
 spin chains with boundaries \cite{48a} and in
 $2D$ quantum integrable field theories \cite{48'}
 (see also below Sec. {\bf \ref{factsc}}).
Taking the reflection operator $K_n(x)$ in the form
\be
\lb{reflH''}
K_n(x) = \frac{y_n - \xi \, x^{2}}{y_n - \xi \, x^{-2}} \; ,
\ee
where $\xi$ is any constant,
we find \cite{IsaO} that this $K_n(x)$ is a solution of (\ref{reflH})
if $y_n$ are the affine generators of $\hat{H}_{M+1}$.
In particular one can easily reduce (\ref{reflH''}) to the solution
\be
\lb{reflH'}
K_n(x) = y_n + \frac{\beta_0/\xi + \xi  + \beta_1 x^2}{x^2 - x^{-2}}
\ee
of the reflection equation (\ref{reflH})
if in addition we require that $y_n$ satisfies a quadratic
characteristic equation
$y_n^2 + \beta_1 y_n + \beta_0 = 0$
($\forall \beta_0,\beta_1 \in \mathbb{C}\backslash 0$).
The solution (\ref{reflH''})
is obviously regular: $K_n(1)=1$, and obeys a "unitary condition":
$$
K_n(x) K_n(x^{-1}) = 1 \; .
$$
We stress that the simplest
solution (\ref{reflH'}) of the reflection equation (\ref{reflH})
was considered in \cite{LevMar} (another special solutions
where found in \cite{KulMud}).

If one has a solution of eq. (\ref{reflH}) for $n=m$, then, a solution for
$n = m+1$ can be constructed by means of the formula:
$$
K_{m+1}(x) = (\lambda x)^2 \, \sigma_{m}(x) K_{m}(x) \sigma_{m}(x) \; .
$$
In particular one can take $K_{n-1}(x)=1$ and,
using (\ref{ybeH}) and (\ref{hunit1}),  directly check that (cf. (\ref{reflH'}))
$$
\frac{K_n(x)}{x^2} = \lambda^2 \, \sigma_{n-1}^2(x) =
\lambda^2 \left( \sigma_{n-1}(x^2) + \frac{(x - x^{-1})^2}{\lambda^2} \right) =
\sigma_{n-1}^2 + \frac{2-(2+\lambda^2) \, x^2}{x^2 - x^{-2}} \; ,
$$
solves eq. (\ref{reflH}).

\vspace{0.1cm}

\noindent
{\bf Remark 1.} Consider the following inclusions of the subalgebras
$\hat{H}_{1} \subset \hat{H}_{2} \subset \dots \subset \hat{H}_{M+1}$:
$$
\{y_1; \sigma_1, \dots ,\sigma_{n-1}\} = \hat{H}_{n} \subset \hat{H}_{n+1} =
\{y_1; \sigma_1, \dots ,\sigma_{n-1},\sigma_n \} \; .
$$
Then, following \cite{IsaO},
\cite{Isa7} we equip the algebra $\hat{H}_{M+1}$
by linear mappings
$$
Tr_{D(n+1)}: \;\; \hat{H}_{n+1} \to \hat{H}_{n} \; ,
$$
from the algebras $\hat{H}_{n+1}$ to its subalgebras $\hat{H}_{n}$,
such that for all $X,X' \in \hat{H}_{n}$ and $Y \in \hat{H}_{n+1}$ we have
\be
\label{mapp}
\begin{array}{c}
Tr_{D(n+1)} ( X ) = Z^{(0)} \, X  \, , \;\;
Tr_{D(n+1)} ( X \, Y \, X' ) = X \, Tr_{D(n+1)}(Y) \, X' \;\;  \, , \\[0.1cm]
Tr_{D(n+1)} ( \sigma_n^{\pm 1} X \sigma_n^{\mp 1}) =
Tr_{D(n)} (X)  \; , \;\;\;
Tr_{D(n+1)} (X \sigma_n X') =  X \, X' \; , \\[0.1cm]
Tr_{D(1)} (y_1^k)= Z^{(k)} \; , \;\;
Tr_{D(n)}  Tr_{D(n+1)} ( \sigma_n Y ) =
Tr_{D(n)}  Tr_{D(n+1)} ( Y \sigma_n) \; ,
\end{array}
\ee
where $Z^{(k)} \in {\bf C}\backslash \{0\}$ $(k \in \mathbb{Z})$ are constants.
We stress that $Z^{(k)}$ could be considered as additional generators
of an abelian subalgebra $\hat{H}_{0}$ which extends $\hat{H}_{M+1}$
and be central in $\hat{H}_{M+1}$, but for us it is enough to put $Z^{(k)}$ to constants.

Using the maps $Tr_{D(n+1)}$ one can
construct elements
\be
\lb{tau1}
\tau_n(x) = Tr_{_{D(n+1)}} \bigl( \sigma_{n}(x) \cdots
\sigma_1(x) \, K_1(x) \,
\sigma_1(x)  \cdots \sigma_{n}(x) \bigr)
\; \in \; \hat{H}_{n}  \; ,
\ee
where $\sigma_i(x)$ are Baxterized elements (\ref{baxtH}) and
$K_1(x)$ is a solution (\ref{reflH''})
of the reflection equation (\ref{reflH})
for $n=1$.
The elements (\ref{tau1}) are generating functions
 for a commutative family of elements in $\hat{H}_{n}$
 since we have (see \cite{IsaO}, \cite{Isa7})
$$
[\tau_n(x), \, \tau_n(z)]=0 \;\;\;\;\; (\forall x,z) \; .
$$
Moreover the elements (\ref{tau1}) are analogs of
Sklyanin's transfer-matrices \cite{48a} and, making use of the
elements $\tau_n(x)$, one
 can formulate \cite{Isa7} the integrable open Hecke chain models
  with nontrivial boundary conditions. These models generalize the
  quantum integrable spin models of the Heisenberg type.
The local Hamiltonian of the open Hecke chain is
\be
\lb{ham2}
{\cal H}_n =   \sum_{m=1}^{n-1} \sigma_m -
\frac{\lambda}{2} \, y'_1(1)  \; .
\ee
This Hamiltonian (up to a normalization factor and additional constant)
can be obtained by differentiating $\tau_{n}(x)$
with respect to spectral
parameter $x$ at the point $x=1$. The Hamiltonian (\ref{ham2})
describes the open chain model with nontrivial boundary condition on the first site
(given by the second term in (\ref{ham2}))
and free boundary condition on the last site of the chain.
In \cite{Isa7}, we show that the transfer matrix
elements $\tau_{n}(x)$ satisfy functional relations generalizing functional relations ($T-Q$ relations) for transfer matrices in
solvable open spin chain models (see, e.g.,
\cite{Zhou}, \cite{Nepo} and the references therein).

\vspace{0.1cm}

\noindent
{\bf Remark 2.}
Interrelations of periodic
$AH_{M}$ (see point {\bf 1.} above)
and affine $\hat{H}_{M}$ (see point {\bf 2.} above)
Hecke algebras has been discussed in
\cite{OgPya}.
Here we present more explicit construction \cite{Levy} of these interrelations
which is valid even for the braid group case (when the Hecke condition (\ref{ahecke})
is relaxed).

Consider the affine braid group
$\hat{\cal B}_M = {\cal B}_M(C)$ (see Definition
{\bf \em \ref{def14}}
in Section {\bf \ref{gabg1}}) with generators
 $\{\sigma_1 , \dots , \sigma_{M-1}, y_{1}\}$.
The generator $y_{1}$ satisfies reflection equation and locality
conditions
$$
\sigma_{1} \, y_{1} \, \sigma_{1} \, y_{1} =
y_{1} \, \sigma_{1} \, y_{1} \,  \sigma_{1} \; , \;\;\; [y_1 , \, \sigma_k]=0 \;\;
(k = 2, \dots ,M-1) \; .
$$
Then, the operator
$$
X = \sigma_{M-1 \leftarrow 1} \, y_{1} \in \hat{\cal B}_M \; ,
$$
solves eqs. (\ref{refl1}), (\ref{refl1a}) and one can introduce new
generator $\sigma_M \in \hat{\cal B}_M$ according to (\ref{refl1a}):
\be
\lb{refl4}
\sigma_M = \sigma_{M-1 \leftarrow 1} \, y_{1} \,
\sigma_1 \, y_1^{-1} \, \sigma^{-1}_{M-1 \leftarrow 1} \; ,
\ee
which satisfies (\ref{refl3}). Thus, eq. (\ref{refl4})
defines the homomorphism $\overline{\cal B}_M \rightarrow
\hat{\cal B}_M$.

This homomorphism of affine braid groups
is readily carried over to the Hecke algebra case.
Indeed, the definition (\ref{refl1a}) of the additional generator $\sigma_M$
(needed to close the set of the generators
 $\sigma_k \in H_{M}(q)$ to the periodic chain)
looks like the similarity transformation of $\sigma_1$. Thus, the characteristic
Hecke identity (\ref{ahecke}) for the elements $\sigma_1$ and $\sigma_M$ coincides.

\subsubsection{q-Dimensions of idempotents in $H_{M}(q)$
and knot/link polynomials\label{qdlink}}

Here we follow the approach presented in \cite{IsOg3}.
Consider a linear map $Tr_{_{D(n+1)}}$: $H_{n+1}(q) \to H_n(q)$
from the Hecke algebra $H_{n+1}(q)$ to its subalgebra $H_n(q)$
which is defined by formulas (\ref{mapp}),
where we take $y_1=1$ (it means that $Z^{(k)} = Z^{(0)}$, $\forall k$)
and fix the constant
\be
\lb{z000}
Z^{(0)} = \frac{1 - q^{-2d}}{q - q^{-1}} \; ,
\;\;\;\;\; Z^{(0)} \equiv  Tr_{_{D(n)}}({\bf 1}) \; ,
\ee
 for later convenience.
Then one can define an Ocneanu's trace
${\cal T}r^{(M)}$: $H_{M}(q) \to \mathbb{C}$
as a sequence of maps
 \be
 \lb{ocntr0}
{\cal T}r^{(M)} := Tr_{D(1)} Tr_{D(2)} \cdots Tr_{D(M)} \; .
\ee

\begin{proposition}\label{ocntr1}
 The Jucys-Murphy elements $y_k \in H_{M+1}$ satisfy
the following identity \cite{IsOg3}
 \be
 \lb{ident5}
1 + \lambda \, Tr_{_{D (M+1)}}\left( \frac{y_{M+1}}{t-y_{M+1}} \right)
 = \frac{(t - q^{-2d})}{(t - 1)}
\prod_{k=1}^M \, \frac{(t-y_{k})^2}{(t-q^2 y_{k}) (t-q^{-2} y_{k})} \; ,
 \ee
 where $\lambda=q-q^{-1}$ and $t$ is a parameter.
\end{proposition}
{\bf Proof.}
Taking into account the definition
(\ref{jucmu}) of the generators $y_M$ we have equations
\be
\lb{qdim3}
\frac{1}{(t-y_{M+1})} \sigma^{-1}_M = \sigma^{-1}_M \frac{1}{(t-y_M)} +
\frac{\lambda \, y_M}{(t-y_M)} \frac{1}{(t-y_{M+1})} \ee
\be
\lb{qdim4}
 \frac{1}{(t-y_{M+1})} \sigma_M = \sigma^{-1}_M  \frac{1}{(t-y_M)} +
\frac{ \lambda \, t}{(t-y_M)} \frac{1}{(t-y_{M+1})} \; .
\ee
We multiply the both sides of eq. (\ref{qdim3}) from
the right by $\sigma_M$. Then, in the r.h.s.
of the result, we substitute eq. (\ref{qdim4}) and
apply the map  $Tr_{D (M+1)}$ (\ref{mapp}). Finally
we obtain a recurrence relation
\be
\lb{qdim7}
\frac{(t-q^2 y_{M}) (t-q^{-2} y_{M})}{(t-y_{M})^2}  Z_{M+1}  = Z_M +
\frac{\lambda \, y_M \, (1 - \lambda \, Z^{(0)})}{(t-y_{M})^2} \; ,
\;\;\;\;\;\;
Z_M : = Tr_{_{D (M)}}
\Bigl( \frac{1}{t-y_M} \Bigr)  \; ,
\ee
where the parameter $Z^{(0)}$ is introduced in (\ref{mapp}),
 (\ref{z000}).
Eq. (\ref{qdim7}) is simplified by the substitution
$Z_M = \tilde{Z}_M -
( 1 - \lambda \, Z^{(0)} )/(\lambda \, t)$ and we have
$$
\frac{(t-q^2 y_{M}) (t-q^{-2} y_{M})}{(t-y_{M})^2} \,  \tilde{Z}_{M+1}  = \tilde{Z}_M  \;\; , \;\;\;\;\;\;\;
\tilde{Z}_1 = \frac{1}{\lambda \, t}
\Bigl( 1 + \frac{\lambda \, Z^{(0)}}{(t-1)}\Bigr)\; .
$$
This equation can be easily solved and finally
 we obtain formula
$$
Z_{M+1} = \frac{1}{\lambda \, t}
\left( 1 + \frac{\lambda \, Z^{(0)}}{(t-1)} \right)\prod_{k=1}^M \,
 \frac{(t-y_{k})^2}{(t-q^2 y_{k}) (t-q^{-2} y_{k})}
 - \frac{1}{\lambda \, t} \left[ 1 - \lambda \,
 Z^{(0)} \right] \; ,
 $$
 which is equivalent to (\ref{ident5}).
\hfill \qed

\vspace{0.2cm}

\noindent
We note that the r.h.s. of (\ref{ident5}) is the symmetric function
in $y_k$ $(k=1,...,M)$. It means that the element
(\ref{ident5}) belongs to the center
 of the Hecke algebra $H_M \subset H_{M+1}$.

\begin{proposition}\label{ocntr}
Ocneanu's traces of idempotents $e(T_{{\sf \Lambda}})$
 and $e(T_{{\sf \Lambda}}^{\, \prime})$, corresponding to
 different Young tableaux
 $T_{{\sf \Lambda}}$ and $T_{{\sf \Lambda}}^{\, \prime}$
of the same shape ${\sf \Lambda} \vdash M$,
 coincide. Thus, characteristics
 \be
 \lb{shape}
{\rm qdim}({\sf \Lambda}) := {\cal T}r^{(M)}
e(T_{\sf \Lambda}) =
{\cal T}r^{(M)} e(T_{\sf \Lambda}^{\, \prime})  \; ,
 \ee
depends only on the Young diagram ${\sf \Lambda}$, are
called $q$-dimension of ${\sf \Lambda} \vdash M$ and
we have \cite{W1}
\be
\lb{qdimW}
{\rm qdim}({\sf \Lambda}) = q^{-\,M\, d} \,
 \prod_{n,m \in {\sf \Lambda}} \frac{[d + m-n]_q}{[h_{n,m}]_q}
 \;\; , \;\;\;\;\;\;\;\;
  [h]_q := \frac{q^h - q^{-h}}{q-q^{-1}}\; .
\ee
Here $h_{n,m}$ are hook lengths
of nodes $(n,m)$ of the diagrams ${\sf \Lambda}$, the product runs
 over all nodes of ${\sf \Lambda}$ and the constant
 $d$ is defined in (\ref{z000}).
\end{proposition}
{\bf Proof.}  We follow the proof presented in \cite{IsOg3}.
Idempotents $e(T_{\sf \Lambda})$ and
$e(T_{\sf \Lambda}^{\, \prime})$
 corresponding to two different tableaux
 $T_{\sf \Lambda}$ and
 $T_{\sf \Lambda}^{\, \prime}$ having the same shape
  ${\sf \Lambda}$ are related by several similarity
  transformations with operators $U_j$
  (see l.h.s of (\ref{offdg})). This implies (\ref{shape}).

To calculate the characteristic
"qdim" (\ref{shape}) for the diagram (\ref{qdima01})
 with $M$ nodes and $n$ rows,
we find the right action to the both sides of
(\ref{ident5}) by the idempotent
 $e(T_{\sf \Lambda})$,
where $T_{\sf \Lambda}$ is any Young tableau of
the shape of Young diagram (\ref{qdima01}).
We take the "row-standard" tableau $T_{\sf \Lambda}$
corresponding to the eigenvalues of
 $y_k$ arranged along the rows from left to right
and from the top to bottom:
\be
\lb{qdim25}
\begin{array}{l}
y_1 = 1, \; y_2 = q^2, \; y_3 = q^4, \; \dots , \; y_{\lambda_1-1} = q^{2(\lambda_1 -2)},
\; y_{\lambda_1} = q^{2(\lambda_1 -1)}, \\
y_{\lambda_1 +1} = q^{-2}, \; y_{\lambda_1 +1} = 1, \;
\dots , \; y_{\lambda_1 +\lambda_2} = q^{2(\lambda_2 -2)}, \\
\dots\dots\dots\dots\dots\dots , \\
y_{_{M- \lambda_n +1}} = q^{-2(n-1)}, \; \dots , \;
y_{_M} = q^{2(\lambda_n -n)} \; ,
\end{array}
\ee
where $n$ is the number of rows in ${\sf \Lambda}$.
After substitution of the eigenvalues (\ref{qdim25})
 into the r.h.s. of (\ref{ident5}),
 which is the product over
  all $M$ nods of the Young diagram (\ref{qdima01}),
  and cancelation of many factors we obtain the
 result $(n_k =n,\;\; n_0 := 0)$
\be
\lb{qdim11}
Tr_{\!\!_{D(M+1)}} \! \left( \! \sum_j  P_{j}
\frac{(q-q^{-1}) \, \mu_j }{
t- \mu_j} \!\! \right) =
e(T_{\sf \Lambda}) \!
\left( \! \frac{t-  q^{-2d}}{t- q^{-2n}} \,
\prod_{r=1}^k \frac{t- q^{2(\lambda_{(r)}-n_r)}}{
t-  q^{2(\lambda_{(r)} -n_{r-1})}} - 1
\!\! \right) .
\ee
We inserted into the l.h.s. of
 (\ref{ident5}) the spectral decomposition of
the idempotent $e(T_{\sf \Lambda})$ (see (\ref{pij})):
$$
e(T_{\sf \Lambda}) =
 e(T_{\sf \Lambda}) \sum_j  \Pi_j = \sum_j P_{j} \, ,
 \;\;\;\;\;\; P_{j} \, y_{_{M+1}} =
P_{j} \, \mu_j \, , \;\;\;\;\;\;
 \mu_j := q^{2(\lambda_{(j)} -n_{j-1})} \; .
$$
The idempotent
$P_j = e(T_{{\sf \Lambda}^{(j)}}) \in H_{M+1}$ projects
 $y_{M+1}$ on its eigenvalue $\mu_j$ which also appeared
 in the denominator in the r.h.s. of (\ref{qdim11}) for $r=j$.

 Let us discuss in more detail how one can
 deduce the expression in the r.h.s. of
 (\ref{qdim11}). It is obtained if we evaluate
 the action of the idempotent $e(T_{\sf \lambda})$
 on the element in r.h.s. of (\ref{ident5})
  for each rectangular block  in
  the diagram ${\sf \Lambda}$ (\ref{qdima01}) with
  all rows having the same length
  $\lambda_{(m)}$ and  the number of rows equals to
  $(n_m-n_{m-1})$. The result of such
  evaluation,
  given in the r.h.s. of (\ref{qdim11}), is the
  product of the factor
  $\frac{t-  q^{-2d}}{t- q^{-2n}}$ and all factors
  which visualized as figure in (\ref{block}) and correspond
  to all rectangular blocks in the diagram
  (\ref{qdima01})).

\unitlength=9.5mm
\begin{picture}(25,3.3)(1,0)

 \put(2,2.5){\line(1,0){5}}
 \put(2,1){\line(1,0){5}}
 \put(2,1){\line(0,1){1.5}}

 \put(7,1){\line(0,1){1.5}}

\put(4,3){$_{_{\lambda_{_{(m)}}}}$}

\put(1.2,2.8){$_{n_{m-1}}$}
\put(1.4,1.2){$_{n_{_m} }$}

\put(7.5,1.9){$_{(n_{_m-1} +1\; ,} {_{\; \lambda_{_{(m)}} +1)}}$}
\put(7.1,0.9){$_{(n_{_m}\; ,} {_{\; \lambda_{_{(m)}}) }}$}

 \put(2,0.5){\line(0,1){0.5}}
 \put(2,0.5){\line(1,0){0.5}}
 \put(2.5,0.5){\line(0,1){0.5}}
 \put(2,1){\line(1,0){0.5}}
 \put(2.05,0.8){$_{_{- 1}}$}

 \put(2,2){\line(0,1){0.5}}
 \put(2,2){\line(1,0){0.5}}
 \put(2.5,2){\line(0,1){0.5}}
 \put(2,2.5){\line(1,0){0.5}}
 \put(2.05,2.3){$_{_{+ 1}}$}

  \put(7,2){\line(0,1){0.5}}
 \put(7,2){\line(1,0){0.5}}
 \put(7.5,2){\line(0,1){0.5}}
 \put(7,2.5){\line(1,0){0.5}}
 \put(7.05,2.3){$_{_{- 1}}$}

  \put(6.5,1){\line(0,1){0.5}}
 \put(6.5,1){\line(1,0){0.5}}
 \put(7,1){\line(0,1){0.5}}
 \put(6.5,1.5){\line(1,0){0.5}}
 \put(6.55,1.3){$_{_{+ 1}}$}

 \put(11,2.2){$= \;\;\; \prod\limits_{i=1}^4 \,
 (t - \mu_i')^{\alpha_i} \;\; = $}

 \put(10.5,1){$= \;\;\;
 \frac{(t - q^{-2n_{m-1}})(t - q^{2(\lambda_{(m)}-n_m)})}{
 (t - q^{-2n_{m}})(t - q^{2(\lambda_{(m)}-n_{m-1})})}$}

\end{picture}

\vspace{-2.5cm}

\be
\lb{block}
{}
\ee

\vspace{1cm}

\noindent
Each rectangular block contributes to
the r.h.s. of (\ref{qdim11}) four factors
which correspond to four cells
indicated in (\ref{block})
 by indices $\alpha_i = \pm 1$
and having contents
$\mu_i'$  $(i=1,...,4)$. Indices $\alpha_i = \pm 1$
 are powers of the factors $(t - \mu_i')^{\alpha_i}$, in the r.h.s. of (\ref{block}). Two factors
 which corresponds sells with contents
 $(-n_{m-1})$ and $(-n_{m})$ are canceled in r.h.s. of (\ref{qdim11}) by neighboring
 blocks from top and bottom if any. The other cells
 of the block (\ref{block})
 have the powers equal to zero, the corresponding
 factors are canceled and
 do not contribute to the r.h.s. of (\ref{qdim11}).

Now we compare the residues at $t=\mu_j$
in both sides of eq. (\ref{qdim11}) and deduce
$$
Tr_{\!_{D(M+1)}} \bigl( e(T_{{\sf \Lambda}^{(j)}}) \bigr)
= \frac{\mu_j^{-1}}{(q-q^{-1})} e(T_{\sf \Lambda}) \,
 \left. 
 \frac{(t-  q^{-2d})
 \prod\limits_{r=1}^k(t- q^{2(\lambda_{(r)}-n_r)})}{(t-q^{-2n})
 \prod\limits_{\stackrel{r\neq j}{r=1}}^k
 (t- q^{2(\lambda_{(r)} -n_{r-1})})}
 \right|_{t = \mu_j\equiv q^{2(\lambda_{(j)}-n_{j-1})}}
 \!\!\!  =
$$
\be
\lb{qdim13}
=e(T_{\sf \Lambda}) \cdot q^{-d} \, [q^{(\lambda_{(j)}-n_{j-1}+d)}]_q \,
 \frac{ \prod_{n,m \in {\sf \Lambda}}  \; [h_{n,m}]_q}{
 \prod_{n,m \in {\sf \Lambda}^{(j)}} \; [h_{n,m}']_q} \; ,
\ee
where $h_{n,m}$ and $h_{n,m}'$ are hook
lengths\footnote{The hook length
of the node $(n,m)$ of the diagram
${\sf \Lambda}=[\lambda_1,\lambda_2,...]$ is defined as
$h_{n,m}=(\lambda_n + \lambda_m^{\vee}-n-m+1)$,
where ${\sf \Lambda}^{\vee} =[\lambda_1^{\vee},\lambda_2^{\vee},...]$ is the transpose partition of the partition ${\sf \Lambda}$.}
of nodes $(n,m)$ of the diagrams
${\sf \Lambda} \vdash M$ and
 ${\sf \Lambda}^{(j)} \vdash (M+1)$. The diagram
 ${\sf \Lambda}^{(j)}$ is obtained by adding to
 ${\sf \Lambda}$ (shown in (\ref{qdima01}))
 a new node with coordinates
$(n_{j-1} +1, \lambda_{(j)} +1)$ as it
is shown in the picture

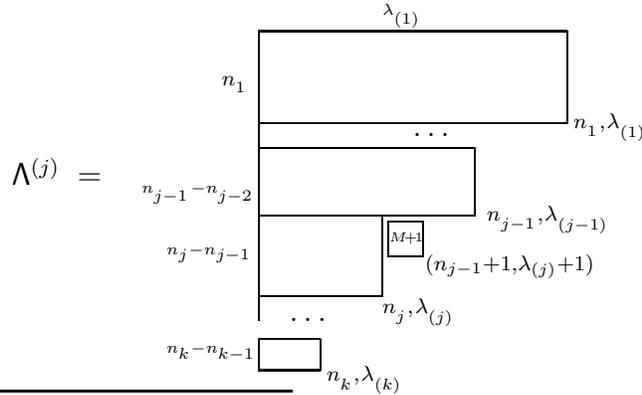
\unitlength=8.2mm
\begin{picture}(25,6.1)(-2,0)
\put(1,3){${\sf \Lambda}^{(j)}\;  = $}

 \put(5,5.5){\line(1,0){5}}
 \put(5,4){\line(1,0){5}}
 \put(5,4){\line(0,1){1.5}}

 \put(10,4){\line(0,1){1.5}}

 \put(5,3.6){\line(1,0){3.5}}
 \put(5,2.5){\line(1,0){3.5}}

 \put(7.5,3.8){$\dots$}
 \put(8.5,2.5){\line(0,1){1.1}}

 \put(5,0.8){\line(0,1){3.2}}
 \put(7,1.2){\line(0,1){1.3}}
 \put(5,1.2){\line(1,0){2}}
\put(5.5,0.8){$\dots$}

 \put(7.1,1.85){\framebox(0.55,0.55)}
 \put(7.11,2.2){\tiny $_{_{M\! +\! 1}}$}
 \put(7.7,1.7){$_{(n_{j-1} +1, \lambda_{(j)} +1)}$}

 \put(5,0){\line(0,1){0.5}}
 \put(6,0){\line(0,1){0.5}}
\put(5,0){\line(1,0){1}}
\put(5,0.5){\line(1,0){1}}

\put(7,5.9){$_{_{\lambda_{_{(1)}}}}$}
\put(4.4,4.7){$_{n_{_1}}$}
\put(3.1,3){$_{_{n_{_{j-1}}-n_{_{j-2}}}}$}
 \put(3.5,2){$_{_{n_{j}-n_{j-1}}}$}
\put(3.5,0.4){$_{_{n_k-n_{k-1}}}$}
\put(10.1,4){$_{n_{_1},} {_{\lambda_{_{(1)}}}}$}
\put(8.7,2.5){$_{n_{_{j-1}},} {_{\lambda_{_{(j-1)}}}}$}
\put(7,1){$_{n_{_j},} {_{\lambda_{_{(j)}}}}$}
\put(6.1,-0.1){$_{n_{_k},} {_{\lambda_{_{(k)}}}}$}

\end{picture}

\vspace{-2cm}
\be
\lb{lambj}
{}
\ee

\vspace{1cm}

\noindent
To deduce the last formula in (\ref{qdim13}) we need to check
the identity
\be
\lb{qdim10}
q^{(n_{j-1}-\lambda_{(j)})} \,
 \frac{(q^{2(\lambda_{(j)} -n_{j-1})} - q^{2(\lambda_{(j)}-n_j)})}{(q^{2(\lambda_{(j)}-n_{j-1})}
 - q^{-2n_k})} \,
\prod\limits_{\stackrel{r=1}{r \neq j}}^k \frac{(q^{2(\lambda_{(j)} -n_{j-1})} - q^{2(\lambda_{(r)}-n_r))}}{
(q^{2(\lambda_{(j)} -n_{j-1})} -
q^{2(\lambda_{(r)} - n_{r-1})})}   =
\frac{\prod\limits_{(r,m) \in {\sf \Lambda}}
[h_{r,m}]_q}{
 \prod\limits_{(r,m) \in {\sf \Lambda}^{(j)}}
  [h'_{r,m}]_q} \; ,
\ee
 where $n_0 \equiv 0$, $n_k \equiv n$,
 while $h_{r,m}$ and $h'_{r,m}$ are hook lengths for cells
  with coordinates $(r,m)$ in the diagrams
  ${\sf \Lambda}$ and
 ${\sf \Lambda}^{(j)}$, respectively.
  The products in the r.h.s. of (\ref{qdim10})
  run over all cells of ${\sf \Lambda}$ and
 ${\sf \Lambda}^{(j)}$. To prove (\ref{qdim10})
 we note that the lengths of hooks
 $h_{r,m}$ and $h'_{r,m}$ for diagrams ${\sf \Lambda}$ and ${\sf \Lambda}^{(j)}$
 differ only for cells, for which $r = n_{j-1}+1$,
  or $m=\lambda_{(j)}+1$, i.e., for
   cells located in the same row,
   or in the same column with additional cell
 $(n_{j-1}+1,\lambda_{(j)}+1)$. Thus, we have
  \be
\lb{qdim14}
 \frac{\prod\limits_{r,m \in \lambda_n}
 [h_{r,m}]_q}{
 \prod\limits_{r,m \in \lambda_{n+1}}
  [h'_{r,m}]_q} =
  \prod\limits_{r=1}^{n_{j-1}}  \frac{ [h_{r,\lambda_{(j)}+1}]_q}{
  [h'_{r,\lambda_{(j)}+1}]_q}
  \prod\limits_{m =1}^{\lambda_{(j)}}  \frac{ [h_{n_{j-1}+1,m}]_q}{ [h'_{n_{j-1}+1,m}]_q}
  \; .
 \ee
 Further, if rows with numbers $r$ and $r+1$
  in the diagram ${\sf \Lambda}$ have the same length,
 then we obviously have
 $h_{r,\lambda_{(j)}+1} = h_{r+1,\lambda_{(j)}+1}'$.
 And analogously, if the columns with numbers
 $m$ and $m+1$ in the diagram ${\sf \Lambda}$ have the
 same height, then we have
 $h_{n_{j-1}+1,m} = h_{n_{j-1}+1,m+1}'$. That is
  why a lot of factors in the
 r.h.s. of (\ref{qdim14}) are canceled, according
 to the block form of diagram (\ref{qdima01}), and we
 obtain contributions only from the first and
 the last row in blocks (located above the additional
 cell in ${\sf \Lambda}^{(j)}$)
 of the diagram ${\sf \Lambda}$
  \be
  \lb{qdim15a}
 \prod\limits_{r=1}^{n_{j-1}}
 \frac{[h_{r,\lambda_{(j)}+1}]_q}{ [h'_{r,\lambda_{(j)}+1}]_q} =
 \prod\limits_{p=1}^{j-1}
 \frac{[h_{n_p,\lambda_{(j)}+1}]_q}{
 [h'_{n_{p-1}+1,\lambda_{(j)}+1}]_q} =
 q^{n_{j-1}} \prod\limits_{p=1}^{j-1}
 \frac{(q^{2(\lambda_{(p)}-n_{p})}
 -q^{2(\lambda_{(j)} - n_{j-1})})}{
 (q^{2(\lambda_{(p)}-n_{p-1})}
 -q^{2(\lambda_{(j)} - n_{j-1})})} \; ,
 \ee
 and contributions only from the first and the last
 column in the blocks (located to the left of
 the additional cell) of the diagram $\lambda_n$
 \be
  \lb{qdim15b}
 \prod\limits_{m =1}^{\lambda_{(j)}}
 \frac{ [h_{n_{j-1}+1,m}]_q}{ [h'_{n_{j-1}+1,m}]_q} =
 \prod\limits_{p =j}^{k}
 \frac{ [h_{n_{j-1}+1,\lambda_{(p)}}]_q}{ [h'_{n_{j-1}+1,\lambda_{(p+1)}+1}]_q} =
 q^{-\lambda_{(j)}} \prod\limits_{p =j}^{k}
 \frac{(q^{2(\lambda_{(j)}-n_{j-1})}
 - q^{2(\lambda_{(p)}-n_{p})})}{(q^{2(\lambda_{(j)}-n_{j-1})}
 - q^{2(\lambda_{(p+1)}-n_{p})}} \; .
 \ee
 The substitution of (\ref{qdim15a}) and
  (\ref{qdim15b}) into (\ref{qdim14})
  gives (\ref{qdim10}).
Finally we apply the Ocneanu's trace ${\cal T}r^{(M)}$
to both sides of eq. (\ref{qdim13}) and find the
recurrence relation:
$$
{\rm qdim}({\sf \Lambda}^{(j)}) = {\rm qdim}({\sf \Lambda}) \;
q^{-d} \; [\lambda_{(j)}-n_{j-1}+d]_q \;
 \frac{ \prod\limits_{n,m \in {\sf \Lambda}}\;  [h_{n,m}]_q}{
 \prod\limits_{n,m \in {\sf \Lambda}^{(j)}}  [h_{n,m}']_q} \; ,
$$
which is uniquely solved (up to a constant
 multiplier\footnote{We fix this multiplier by the
  condition ${\rm qdim}(\Box)= q^{-d}[d]_q = Z^{(0)}$;
  see (\ref{z000}).}) as in (\ref{qdimW}). \hfill \qed

\vspace{0.2cm}

The $R$-matrix representations (see \cite{Jimb1}, \cite{10})
of the Hecke algebra $H_{M+1}(q)$
were discussed
 in Section {\bf \ref{qgrsl}} in the
 context of the quantum group $GL_q(N)$ and
 in Section {\bf \ref{qsuper}} in the
 context of the quantum supergroup $GL_q(N|K)$.
For the $R$-matrices (\ref{3.6.1aa}) related
to the quantum supergroup $GL_q(N|K)$, the parameter $d$ is
equal to $(N-K)$. This fact follows from eqs.
(\ref{mapp}) and (\ref{z000})
in the limit $q \to 1$. It also justifies our
choice of the parametrization of $Z_0$ in (\ref{z000}).

The statement (\ref{shape})
in Proposition {\bf \em \ref{ocntr}} can be generalized. Let $T$ be a quantum matrix satisfying
\be
\lb{qgr}
\hat{R}_{12} \, T_1 \, T_2 =  T_1 \, T_2 \, \hat{R}_{12} \; ,
\ee
where $\hat{R}_{12} = \rho(\sigma_1)$ is the $R$-matrix representation
of the Hecke algebra.


\begin{proposition}\label{ocntr2}
 The quantum traces
(for the definition of the quantum traces
see Section {\bf \ref{qtrace}}) of the matrices
$[T_1 \cdots T_{_{M}} \,
\rho (e(T_{\sf \Lambda}))]$ and
$[T_1 \cdots T_{_{M}} \,
\rho (e(T_{\sf \Lambda}^{\, \prime}))]$,
where different tableaux $T_{\sf \Lambda}$ and
$T_{\sf \Lambda}^{\, \prime}$ are of the same shape
 ${\sf \Lambda} \vdash M$, coincide
\be
\lb{qgr2}
\chi_{\sf \Lambda}(T) :=
Tr_{{\cal D}(1...M)} \left( T_1 \cdots
T_{M} \,
\rho (e(T_{\sf \Lambda})) \right) =
Tr_{{\cal D}(1...M)} \left( T_1 \cdots
T_{M} \,
\rho (e(T_{\sf \Lambda}^{\, \prime})) \right) \; .
 \ee
Thus, the element $\chi_{\sf \Lambda}(T)$ depends only on the
 shape of the diagram ${\sf \Lambda}$.
\end{proposition}
According to Proposition {\bf \em \ref{comset}}
(see Subsect. {\bf \ref{cencom}}) the elements
$\chi_{\sf \Lambda}(T)$ for
all Young diagrams ${\sf \Lambda} \vdash M$ $(M=1,2,3,...)$
 generate the commutative subalgebra in
$RTT$ algebra (\ref{qgr}).

\vspace{0.2cm}

Consider the $GL_q(N)$ quantum group (\ref{qgr}) with the standard $GL_q(N)$ Drinfeld-Jimbo $\hat{R}_{12}$ matrix  (\ref{3.3.6b}).
As we mentioned above (see Sect.
{\bf \ref{qgrsl}} and \cite{Jimb1}, \cite{10})
this standard $GL_q(N)$ matrix $\hat{R}_{12}$
 ($R$-matrix in the defining representation)
 gives the representation of the
Hecke algebra. We note that the $GL_q(N)$ quantum matrix $T$ can be realized by arbitrary
numerical diagonal $(N \times N)$ matrix
$Y = {\rm diag}(x_1,...,x_N)$. Then
$\chi_{\sf \Lambda}(Y)$ is a numerical
function of the deformation parameter $q$ and the
 entries $\{ x_i \}$ of $Y$. In the classical
limit $q \to 1$ the operator
$\rho (e(T_{\sf \Lambda}))$ tends to the Young
projector, and the function $\chi_{\sf \Lambda}(Y)$
coincides with a character of the
element $Y \in GL(N)$ in the representation corresponding to the diagram ${\sf \Lambda}$; i.e.
$\chi_{\sf \Lambda}(Y)|_{q \to 1}$ coincides with the Schur
polynomial $S_{\sf \Lambda}(x_1,...x_N)$.

\vspace{0.1cm}

\noindent
{\bf Remark 1.} The hook formula (\ref{qdimW})
for the $q$-dimension of ${\sf \Lambda} \vdash M$
is written in the remarkable form
(which is more convenient for calculations)
\be
\lb{best01}
{\rm qdim}({\sf \Lambda}) = q^{-\,M\, d} \,
\prod_{i=1}^k \frac{[d+i-1]_q! }{[d-\lambda^{\vee}_i +i-1]_q![\lambda^{\vee}_i+k-i]_q!} \prod_{i<j} [\lambda^{\vee}_i-\lambda^{\vee}_j+(j-i)]_q \; ,
\ee
where ${\sf \Lambda}^\vee = (\lambda^{\vee}_1,\lambda^{\vee}_2,...,\lambda^{\vee}_k)$
is the  transpose partition of ${\sf \Lambda}$
and $[h]_q := \frac{q^h - q^{-h}}{q-q^{-1}}$.

\vspace{0.1cm}

\noindent
{\bf Remark 2.}
At the end of this Subsection, we derive a universal analogue
(in terms of Hecke algebra generators) of the
 formula (\ref{qxinv}) for knot/link invariants.
 Let $B_{1 \to M}$ be monomial written as a product of
 generators $\sigma_i \in H_{M}$. It is clear that
 $B_{1 \to M}$ is visualized as a braid with
  $M$ strands. Then, by means of the Ocneanu's trace
  (\ref{ocntr0}) we construct the Hecke algebraic
  analog of (\ref{qxinv}) in the form
  \be
  \lb{qxinv2}
  {\sf Q}(B_{1 \to M}) :=
 {\cal T}r^{(M)} \bigl( B_{1 \to M} \bigr)  \; .
  \ee
Insert in the right hand side of (\ref{qxinv2}) the resolution
of the unit operator (see second equation in (\ref{orind}))
 \be
 \lb{orind2}
 1 = \sum_{{\sf \Lambda} \vdash M}
 \sum_{{\sf T}^{^{\vec{a}}}_{\sf \Lambda}}
 e({\sf T}^{^{\vec{a}}}_{\sf \Lambda}) =
 \sum_{\vec{a}} P(X_{\vec{a}}) \; ,
 \ee
 where $e({\sf T}^{^{\vec{a}}}_{\sf \Lambda}) = P(X_{\vec{a}})$
 are mutually orthogonal idempotents
 related to the standard Young tableau
 ${\sf T}^{^{\vec{a}}}_{\sf \Lambda}$
 (with a content $\vec{a}= (a_1,...,a_{M})$)
 having the shape of the Young diagram
 ${\sf \Lambda} \vdash M$ (or related to the path $X_{\vec{a}}$
 in the coloured Young graph for $H_{M}$).
 The sum in the r.h.s. of (\ref{orind2}) is going over all
 standard tableaux with $M$ nodes,
 or equivalently over all their contents $\vec{a} \in {\rm Spec}(y_1,...,y_M)$.
 As a result we obtain for knot/link invariants
 (\ref{qxinv2}) the expressions
 \be
  \lb{qxinv3}
  \begin{array}{c}
  {\sf Q}(B_{1 \to M}) = \sum\limits_{\vec{a}}
 {\cal T}r^{(M)} \bigl(B_{1 \to M}\, P(X_{\vec{a}})  \bigr)
  = \sum\limits_{\vec{a}}
 {\cal T}r^{(M)} \bigl(P(X_{\vec{a}}) \, B_{1 \to M}\, P(X_{\vec{a}})  \bigr) = \\ [0.3cm]
 = \sum\limits_{\vec{a}} C_{\vec{a}}(B_{1 \to M})
 {\cal T}r^{(M)} \bigl(P(X_{\vec{a}}))  \bigr) =
 \sum\limits_{\sf \Lambda \vdash M}
 {\rm qdim}({\sf \Lambda})
 \sum\limits_{\vec{a}({\sf \Lambda})}
 C_{\vec{a}({\sf \Lambda})}(B_{1 \to M})  \; ,
 \end{array}
  \ee
  where we used relation $P(X_{\vec{a}})^2=P(X_{\vec{a}})$,
  cyclic property of ${\cal T}r^{(M)}$, the identity (\ref{pbp})
 for the diagonal matrix units of $H_{M}$ and substitute
   ${\cal T}r^{(M)} \bigl(P(X_{\vec{a}({\sf \Lambda})}))
   ={\rm qdim}({\sf \Lambda})$ for the Young diagrams
  ${\sf \Lambda} \vdash M$.
  In the last equality of (\ref{qxinv3}) we split the
  sum over all contents $\vec{a}$ of the standard
  Young tableaux with $M$ nodes
  into the sum over Young diagrams
  ${\sf \Lambda} \vdash M$ and sum
  over all contents $\vec{a}({\sf \Lambda})$
  of the Young tableaux ${\sf T}_{\sf \Lambda}$
  having the fixed shape
  ${\sf \Lambda} \vdash M$. We note that in (\ref{qxinv3})
  the q-dimensions
  ${\rm qdim}({\sf \Lambda})$ (given by formula (\ref{qdimW}))
 are independent of the braid $B_{1 \to M}$ and all dependence
  on $B_{1 \to M}$ contains in the coefficients
  $C_{\vec{a}({\sf \Lambda})}(B_{1 \to M})$. As it was
  indicated in Subsection {\bf \ref{idbaxt}} the coefficients
  $C_{\vec{a}({\sf \Lambda})}(B_{1 \to M})$ can be
  explicitly calculated with the help of equations (\ref{pbp2}),
 (\ref{pbp3}) and (\ref{pbp6}).  The $R$-matrix
 version \cite{RoJo} of the formula (\ref{qxinv3})
 is extensively used in \cite{MMM} (see also references therein)
 for calculations of HOMFLY ($GL_q(N)$, $N=d$) knot/link polynomials.

\subsection{\bf \em Birman-Murakami-Wenzl algebras
$B\!M\!W_{M+1}(q,\nu)$\label{bmwalg}}
\setcounter{equation}0

\subsubsection{Definition}

The Birman-Murakami-Wenzl algebra
$B\!M\!W_{\!\!_{M+1}}(q,\nu)$
is generated by the elements $\kappa_i$ $(i=1,\dots,M)$
and invertible elements $\sigma_i$ (\ref{braidg}) which
satisfy the following relations \cite{BW1}, \cite{W2}, \cite{Mur1}
\be
\lb{bmw1}
\kappa_i \, \sigma_i =  \sigma_i \, \kappa_i =
\nu \, \kappa_i \; ,
\ee
\be
\lb{bmw2}
\kappa_i \, \sigma_{i-1}^{\pm 1} \, \kappa_i = \nu^{\mp 1} \, \kappa_i \; ,
\ee
\be
\lb{bmw3}
\sigma_i - \sigma_i^{-1} = \lambda \, ( 1 - \kappa_i) \; ,
\ee
where
$\nu \in \mathbb{C}\backslash \{0, \pm q^{\pm 1} \}$
is an additional parameter of the algebra and
$\lambda = q -q^{-1}$. The following relations
can be derived from (\ref{braidg}), (\ref{bmw1}) - (\ref{bmw3})
\begin{eqnarray}
\lb{bmw4}
& & \kappa_i \, \kappa_{i} =  \mu \, \kappa_i
 \; , \\ \vspace{0.2cm}
\lb{bmw4a}
& &
( \; \mu = (\lambda + \nu^{-1} - \nu)/\lambda = - (\nu +q^{-1})(\nu-q) \, (\lambda \nu)^{-1} \; )
 \; , \\ \vspace{0.2cm}
\lb{bmw5}
& & \kappa_i \, \sigma_{i \pm 1} \, \sigma_i =
\sigma_{i \pm 1} \, \sigma_i \, \kappa_{i \pm 1} \; , \\ \vspace{0.2cm}
\lb{bmw6}
& & \kappa_i \, \sigma_{i \pm 1} \, \sigma_{i} =
\kappa_i \, \kappa_{i \pm 1} \; , \\ \vspace{0.2cm}
\lb{bmw7}
& & \kappa_i \, \sigma^{-1}_{i \pm 1} \, \sigma^{-1}_{i} =
\kappa_i \, \kappa_{i \pm 1} \; , \\  \vspace{0.2cm}
\lb{bmw8a}
& & \sigma_{i \pm 1} \, \kappa_i  \, \sigma_{i \pm 1} =
\sigma^{-1}_i \, \kappa_{i \pm 1} \, \sigma^{-1}_i \; , \\  \vspace{0.2cm}
\lb{bmw8}
& & \kappa_i \, \kappa_{i \pm 1} \, \kappa_i =  \kappa_{i} \; , \\ \vspace{0.2cm}
\lb{bmw8b}
& & \kappa_{i \pm 1} \, \kappa_i  \, (\sigma_{i \pm 1} - \lambda) =
\kappa_{i \pm 1} \, (\sigma_i - \lambda) \; , \\  \vspace{0.2cm}
\lb{bmw9}
& & (\sigma_{i} - \lambda) \, \kappa_{i-1} \,  (\sigma_i - \lambda) =
(\sigma_{i- 1} -\lambda) \, \kappa_{i} \, (\sigma_{i-1} - \lambda) \; .
\end{eqnarray}
Eq. (\ref{bmw4}) is deduced by the action of the element $\kappa_i$
on (\ref{bmw3}) and using (\ref{bmw1}).
Relations (\ref{bmw5}) follow from (\ref{braidg}) and (\ref{bmw3}).
Relations (\ref{bmw6}) and (\ref{bmw7}) with lower signs are obtained
by multiplying (\ref{bmw2}) with
$\sigma_{i-1}^{\mp 1} \sigma_i^{\mp 1}$ from the right and using
(\ref{bmw1}) and (\ref{bmw5}). Eq. (\ref{bmw8a})
follows from (\ref{bmw6}), (\ref{bmw7}).
Combining the pair of relations (\ref{bmw2}) in the form:
$\kappa_i \, (\sigma_{i+1} - \sigma_{i+1}^{-1}) \, \kappa_i=
(\nu^{-1} - \nu) \, \kappa_i$ and using (\ref{bmw3}) and (\ref{bmw4})
we derive (\ref{bmw8}). Eq. (\ref{bmw8b}) is proved as following
$$
\kappa_{i \pm 1} \, \kappa_i  \, (\sigma_{i \pm 1} - \lambda) =
\kappa_{i \pm 1} \, \kappa_i  \, (\sigma^{-1}_{i \pm 1} - \lambda \kappa_{i \pm 1}) =
\kappa_{i \pm 1} \, (\sigma_i - \lambda) \; ,
$$
where we have used (\ref{bmw3}), (\ref{bmw6}) and (\ref{bmw8}). Eq. (\ref{bmw9}) is
deduced by means of eq. (\ref{bmw8b}), its mirror counterpart and eq. (\ref{bmw8}).
The pairs of eqs. in (\ref{bmw5}) - (\ref{bmw8}) (with upper and lower
signs) are related to each other by the similarity transformations
$$
\sigma_{i+1} = V_i \sigma_{i-1} V_i^{-1} \; , \;\;\;
\sigma_{i} = V_i \sigma_{i} V_i^{-1}
$$
where  $V_i = \sigma_{i-1} \sigma_i \sigma_{i+1} \sigma_i \sigma_{i-1} \sigma_i$
(the only braid relations (\ref{braidg}) should be used). We also present the relations
\begin{eqnarray}
\lb{bmw8bb}
& & \kappa_{i \pm 1} \, \kappa_i  \, (\sigma^{-1}_{i \pm 1} + \lambda) =
\kappa_{i \pm 1} \, (\sigma^{-1}_i + \lambda) \; , \\  \vspace{0.2cm}
\lb{bmw9bb}
& & (\sigma^{-1}_{i} + \lambda) \, \kappa_{i-1} \,  (\sigma^{-1}_i + \lambda) =
(\sigma^{-1}_{i- 1} +\lambda) \, \kappa_{i} \, (\sigma^{-1}_{i-1} + \lambda) \; ,
\end{eqnarray}
which are related to (\ref{bmw8b}), (\ref{bmw9}) via the obvious
isomorphism $(\sigma_i, q,\nu) \leftrightarrow
(\sigma_i^{-1}, q^{-1},\nu^{-1})$ of the algebras
$B\!M\!W_{\!\!_{M+1}}(q,\nu) \simeq
B\!M\!W_{\!\!_{M+1}}(q^{-1},\nu^{-1})$, which can be
checked by the substitution $\sigma_i \to \sigma_i^{-1}$
in (\ref{braidg}), (\ref{bmw1}) -- (\ref{bmw3}).

In fact the pair of relations (\ref{bmw2}) (in the definition
of the Birman-Murakami-Wenzl algebra) are not independent for
the case $\nu \neq \lambda$ \cite{27aa}. Indeed, using
$\kappa_i \, \sigma_{i-1} \, \kappa_i = \nu^{- 1} \, \kappa_i$ and (\ref{braidg})
one can deduce
$\sigma_{i-1} \sigma_i \kappa_{i-1} \kappa_i =
\nu \kappa_i \sigma_{i-1} \kappa_i = \kappa_i$,
which is written in the form
$\sigma^{-1}_{i-1} \kappa_i = \sigma_i \kappa_{i-1} \kappa_i$.
Acting to this relation by $\lambda \kappa_i$ from the left, we deduce
$$
\lambda \kappa_i \sigma^{-1}_{i-1} \kappa_i = \lambda \nu \kappa_i \kappa_{i-1} \kappa_i
 = \nu \kappa_i ( \sigma^{-1}_{i-1} - \sigma_{i-1} + \lambda) \kappa_i
 = \nu \kappa_i \sigma^{-1}_{i-1} \kappa_i
+ \nu (\lambda \mu - \nu^{-1}) \kappa_i \; ,
$$
which is equivalent to
$(\lambda - \nu) (\kappa_i \sigma^{-1}_{i-1} \kappa_i - \nu \kappa_i) = 0$ and,
thus, to the above statement.

\subsubsection{Symmetrizers, antisymmetrizers
and Baxterized elements in $B\!M\!W_{M+1}$\label{abmw0}}

Below, for brevity,  we often omit in the notation $B\!M\!W_{\!\!_{M+1}}(q,\nu)$
the dependence on the parameters $q,\nu$.
One can construct the analogs of the symmetrizers and antisymmetrizers for
the algebra $B\!M\!W_{\!\!_{M+1}}$ using the inductive relations
similar to that we have considered in the Hecke case (\ref{santis3}):
\be
\lb{bmws}
S_{1 \to n} =f_{1 \to n}^{(-)}  \, S_{1 \to n-1}  =
 S_{1 \to n-1} \, \overline{f}_{1 \to n}^{(-)} \; ,
\ee
\be
\lb{bmwa}
A_{1 \to n} = f_{1 \to n}^{(+)}  \, A_{1 \to n-1} =
 A_{1 \to n-1} \, \overline{f}_{1 \to n}^{(+)}
\; ,
\ee
where 1-shuffles are
$$
f_{1 \to n}^{(\pm)} = \frac{1}{ [n]_q!} \, \sigma^{(\pm)}_{1}(q^{\pm 1})  \cdots
 \sigma^{(\pm)}_{n-2}(q^{\pm(n-2)})  \,  \sigma^{(\pm)}_{n-1}(q^{\pm(n-1)}) \; ,
$$
$$
\overline{f}_{1 \to n}^{(\pm)} = \frac{1}{ [n]_q!} \,
 \sigma^{(\pm)}_{n-1}(q^{\pm(n-1)}) \,  \sigma^{(\pm)}_{n-2}(q^{\pm(n-2)})  \cdots
 \sigma^{(\pm)}_{1}(q^{\pm 1})  \; ,
$$
and $\sigma^{(\pm)}_i(x)$ are Baxterized elements (cf. (\ref{3.9.15a}), (\ref{3.9.15d}), (\ref{3.9.16}))
for the algebra $B\!M\!W_{\!\!_{M+1}}(q,\nu)$
(see \cite{Jon3}, \cite{Jon2}, \cite{Mur2}, \cite{CGX},
\cite{IsFirst}):
\be
\lb{bmwbax}
\sigma^{(\pm)}_i(x) = \frac{1}{\lambda} (x^{-1} \, \sigma_i - x \, \sigma^{-1}_i) +
\frac{(\nu \pm q^{\pm 1})}{(\nu x \pm q^{\pm 1} x^{-1})} \, \kappa_i =
\ee
\be
\lb{bmwbax5}
 = \left. \frac{x^{-1} -x}{\lambda}
  \left( \sigma_i + \frac{\lambda}{(x^{-2}-1)}   \, {\bf 1} +
\frac{\nu \lambda}{(\nu  - a x^{-2})} \, \kappa_i \right)
\right|_{a = \mp q^{\pm 1}} =
 \ee
$$
 = x \left( {\bf 1} + \frac{1}{\lambda} (x^{-2}-1)  \, \sigma_i \right)
\left( 1 + \frac{\lambda (x^{-2}-1)}{
(\lambda - \nu (1-x^{-2}))(1  \pm q^{\pm 1} \nu^{-1} x^{-2})} \, \kappa_i \right) \; .
$$
These elements are normalized by the conditions
$\sigma^{(\pm)}(\pm 1) = \pm {\bf 1}$, satisfy the Yang-Baxter equations
\be
\lb{ybebmw}
\sigma^{(\pm)}_n(x) \, \sigma^{(\pm)}_{n-1}(xy) \, \sigma^{(\pm)}_n(y) =
\sigma^{(\pm)}_{n-1}(y) \, \sigma^{(\pm)}_n(xy) \, \sigma^{(\pm)}_{n-1}(x) \; .
\ee
and obey
\be
\lb{bmwunit}
\sigma^{(\pm)}_i (x) \, \sigma^{(\pm)}_i(x^{-1}) = \left( 1 - \lambda^{-2} \,
(x-x^{-1})^2 \right) \, {\bf 1} \; .
\ee
Let parameter $a$ be $(-q)$, or $q^{-1}$ (see
 (\ref{bmwbax5})), and we respectively denote
$\sigma^{(a)}_i (x)=\sigma^{(\pm)}_i (x)$. Then, the Baxterized elements (\ref{bmwbax5})
 are written after an additional
  normalization in the form (cf. (\ref{hunit2}))
  \be
\lb{bmwbax3}
\sigma^{(a)\prime}_i(x) =
 \frac{\lambda x^{2}}{(a^{-1}x - a x^{-1})}
 \sigma^{(a)}_i(x) =
  \frac{(\sigma_i - a \, x^{2})}{(\sigma_i - a \, x^{-2})}  \; .
 \ee
New normalized elements (\ref{bmwbax3}) obviously
satisfy ``unitarity conditions'':
$\sigma^{(a)\prime}_i(x) \; \sigma^{(a)\prime}_i(x^{-1}) = {\bf 1}$
and $\sigma^{(a)\prime}(\pm 1) = {\bf 1}$.
Identities (\ref{bmwbax}) -- (\ref{bmwbax3})
are checked with the help of relations
(\ref{bmw1}) -- (\ref{bmw9}).

Note that the elements $\sigma^{(+)}_i (x)$ and
$\sigma^{(-)}_i (x)$ (\ref{bmwbax}) are related to each other by the
transformation $q \leftrightarrow -q^{-1}$, which corresponds to the
isomorphism of algebras $B\!M\!W_{\!\!_{M+1}}(q,\nu)$ $\simeq$
$B\!M\!W_{\!\!_{M+1}}(-q^{-1},\nu)$ and we also have
$$
\sigma^{(+)}_i (x) - \sigma^{(-)}_i (x) =
\frac{\nu (q+q^{-1})(x-x^{-1})}{(x\nu + q x^{-1})(x\nu -q^{-1} x^{-1})}
\, \kappa_i \; .
$$
We also stress that the both inequivalent sets $(\pm)$ of the Baxterized elements
(\ref{bmwbax}) are important for explicit constructions of
 (anti)symmetrizers (\ref{bmws}), (\ref{bmwa}).
To my knowledge these both sets
(\ref{bmwbax}) were firstly presented in paper \cite{CGX}
(see also the very first version \cite{IsFirst} of these lectures).
The only one of these sets was presented in \cite{Jon3},
\cite{Jon2} and in \cite{Mur2}.

It follows from eqs. (\ref{bmw1}) - (\ref{bmw3}) that the
algebra $B\!M\!W_{\!\!_{M+1}}(q,\nu)$ ($\nu \neq \lambda$) is a quotient
of the braid group algebra (\ref{braidg}) if the additional
relations on $\sigma_i$ are imposed
\be
\lb{bmwchar}
(\sigma_i - q)(\sigma_i +q^{-1})(\sigma_i - \nu) = 0 \; ,
\ee
$$
(\sigma_i^{-1} + \lambda - \sigma_i) \left( \sigma^{\pm 1}_{i+1} \,
(\sigma_i^{-1} + \lambda - \sigma_i) - \lambda \nu^{\mp 1} \right) = 0 \; .
$$
This quotient is finite dimensional and the dimension of
$B\!M\!W_{\!\!_{M+1}}(q,\nu)$ is $(2M+1)!! = 1 \cdot 3 \cdots (2M+1)$
(this dimension evidently follows from the graphical representation
(\ref{figa}) of the $B\!M\!W_{\!\!_{M+1}}(q,\nu)$ elements).
The whole set of basis elements for the algebra  $B\!M\!W_{\!\!_{M+1}}(q,\nu)$ appear
in the expansion of the symmetrizer $S_{M+1}$  (\ref{bmws})
(or antisymmetrizer
$A_{M+1}$ (\ref{bmwa})). Note that
the quotient of the Birman-Murakami-Wenzl algebra $B\!M\!W_{\!\!_{M+1}}(q,\nu)$
(\ref{bmw1}) - (\ref{bmw3}) by an ideal generated
by $\kappa_i$ is isomorphic to the $A$-type
Hecke algebra $H_{M+1}(q)$.

The first symmetrizer and antisymmetrizer for
the algebra $B\!M\!W_{\!\!_{M+1}}(q,\nu)$ are (cf. eqs. (\ref{3.7.2}), (\ref{3.9.17}))
\be
\lb{first1}
\begin{array}{c}
S_{1 \to 2} = \frac{1}{[2]_q} \sigma_1^{(-)}(q^{-1})
 = \frac{(\sigma^2_1 - q^{-2})(\sigma^2_1 - \nu^2)}{(q^2 - q^{-2})(q^2-\nu^2)} = \\ \\
=  \frac{1}{[2]_q}
(q^{-1} + \sigma_1 + \frac{\lambda}{1- q \nu^{-1}} \kappa_1)
=  \frac{1}{q^2 - q^{-2}}
(\sigma_1^2 - q^{-2})(1 - \mu^{-1} \, \kappa_1) \; ,
\end{array}
\ee
\be
\lb{first2}
\begin{array}{c}
A_{1 \to 2} = \frac{1}{[2]_q} \sigma_1^{(+)}(q)
 = \frac{(\sigma^2_1 - q^2)(\sigma^2_1 - \nu^2)}{(q^{-2}-q^2)(q^{-2}-\nu^2)} = \\ \\
 = \frac{1}{[2]_q}
(q - \sigma_1 - \frac{ \lambda}{1 + q^{-1} \nu^{-1}} \kappa_1)
=  \frac{1}{q^{-2} - q^{2}}
(\sigma_1^2 - q^{2})(1 - \mu^{-1} \, \kappa_1) \; .
\end{array}
\ee
They are obviously orthogonal to each other and to the element $\kappa_1$
in view of the characteristic equation (\ref{bmwchar}).
The following eqs. also hold
$$
\sigma^{(-)}_1(q) \, S_{1 \to 2}  = 0 = \kappa_1 \, S_{1 \to 2}  \; , \;\;\;
\sigma^{(+)}_1(q^{-1}) \, A_{1 \to 2}   = 0 = \kappa_1 \, A_{1 \to 2}  \; ,
$$
which can be deduced from the "unitarity conditions" (\ref{bmwunit})
 and first equalities in (\ref{first1}), (\ref{first2}).
In fact, these eqs. are special cases
of the more general relations (for $i=1, \dots, n-1$)
$$
\begin{array}{l}
\sigma^{(-)}_i(q) \, S_{1 \to n}  =  S_{1 \to n} \, \sigma^{(-)}_i(q) = 0 \; , \\ \\
\sigma^{(+)}_i(q^{-1}) \, A_{1 \to n} =  A_{1 \to n} \, \sigma^{(+)}_i(q^{-1})  =0 \; ,
\end{array}
$$
which equivalent to the equations $(i = 1, \dots n-1)$:
\be
\lb{anul}
\begin{array}{l}
(\sigma_i -q) \, S_{1 \to n}  = 0 =  S_{1 \to n} \, (\sigma_i -q) \; , \;\;\;
\kappa_i \, S_{1 \to n}  = 0 =  S_{1 \to n} \, \kappa_i
\; , \\ \\
(\sigma_i +q^{-1}) \, A_{1 \to n} =  0 = A_{1 \to n} \, (\sigma_i +q^{-1}) \; , \;\;\;
\kappa_i \, A_{1 \to n} =  0 = A_{1 \to n} \, \kappa_i  \; ,
\end{array}
\ee
and demonstrate that $S_{1 \to M+1}$, $A_{1 \to M+1}$ are central idempotents.
Eqs. (\ref{anul}) can be
readily proved by means of the analogs of the factorization relations
(\ref{rfactor}), (\ref{lfactor}) or by the induction with using of (\ref{bmws}), (\ref{bmwa})
and the Yang-Baxter equations (\ref{ybebmw}).

We note that the idempotents (\ref{bmws}), (\ref{bmwa}) can be easily written in the form
(cf. (\ref{santis2}), (\ref{santis22})):
\be
\lb{pyat1}
\begin{array}{c}
S_{1 \to n} = S_{1 \to n-1} \,
\frac{\sigma^{(-)}_{n-1}(q^{-(n-1)})}{ [n]_q} \, S_{1 \to n-1}
\; ,
\end{array}
\ee
\be
\lb{pyat2}
\begin{array}{c}
A_{1 \to n} = A_{1 \to n-1} \,
\frac{\sigma^{(+)}_{n-1}(q^{n-1})}{ [n]_q} \,
A_{1 \to n-1} \ .
\end{array}
\ee
This inductive definition of the idempotents (\ref{bmws}), (\ref{bmwa}) was also used in
\cite{OgPyaSO} and in \cite{WT} (see Lemma 7.6).
 Note that, in view of the definitions (\ref{bmwbax})
of baxterized elements $\sigma^{(\pm)}_k(x)$ , expressions (\ref{pyat1}) and
(\ref{pyat2}) have singularities for $q^{2k}=1$, $\nu = q^{2k-3}$ and
$q^{2k}=1$, $\nu = -q^{-2k+3}$ ($k=2,\dots,n$), respectively.
It means that the representation theory of the BMW algebras
has to be modified
for $q^{2k}=1$ and $\nu = \pm q^{\pm 2k-3}$.

Using the representations (\ref{pyat1}), (\ref{pyat2}) we prove the analog of Proposition {\bf \ref{prop9}} 
about symmetrizers and antisymmetrizers for the case of the
Birman-Murakami-Wenzl algebra.

\begin{proposition}\label{prop11}
{\it The idempotents  $S_{1 \to n}$ and $A_{1 \to n}$
$(n = 2, \dots M+1)$  (\ref{pyat1}), (\ref{pyat2}) for the Birman-Murakami-Wenzl algebra
are expressed in term of the Jucys-Murphy elements $y_k$  $(k = 2, \dots ,M)$:
\be
\lb{bmwy}
y_1 = 1 \; , \;\;\;\;\; y_{k+1} = \sigma_{k} y_{k} \sigma_k
\; , \;\;\;\;\; [y_k, \, y_m]=0 \; ,
\ee
 as follows
\be
\lb{bmws2}
S_{1 \to n} = \prod_{i=2}^{n} \left( \frac{(y_i - q^{-2})}{(q^{2(i-1)} - q^{-2})}  \,
\frac{(y_i - \nu^2 q^{-2(i-2)})}{(q^{2(i-1)} - \nu^2 q^{-2(i-2)})} \right)  \; ,
\ee
\be
\lb{bmwa2}
A_{1 \to n} =  \prod_{i=2}^{n} \left( \frac{(y_i - q^{2})}{(q^{-2(i-1)} - q^{2})}  \,
\frac{(y_i - \nu^2 q^{2(i-2)})}{(q^{-2(i-1)} - \nu^2 q^{2(i-2)})} \right)  \; .
\ee
}
\end{proposition}

\noindent
{\bf Proof.} 
To prove identity (\ref{bmwa2}) we show that it
 is equivalent to (\ref{pyat2}).
The identity (\ref{bmws2}) for the symmetrizers (\ref{pyat1}) can be justified analogously.
The equations (\ref{first2}) demonstrate that (\ref{bmwa2})
coincides with (\ref{pyat2}) for $n=2$. Then, we use the
induction. Let (\ref{bmwa2})
coincides with the formula (\ref{pyat2}) for
$A_{1 \to n}$  for some fixed
$n \geq 2$ and, thus,
it is the element which satisfies (\ref{anul}). We prove that
the formulas (\ref{pyat2}) and  (\ref{bmwa2}) are equivalent for
$A_{1 \to  n+1}$.
In view of the induction
conjecture and obvious properties
$[A_{1 \to n}, \, y_{n+1}]=0$ (since
$A_{1 \to n}$ is a function of $y_i$) we obtain from
(\ref{bmwa2}):
\begin{equation}
\label{bmwa22}
A_{1 \to n+1} = A_{1 \to n} \, \frac{(y_{n+1} - q^{2})}{(q^{-2 n} - q^{2})}  \,
\frac{(y_{n+1} - \nu^2 q^{2(n-1)})}{(q^{-2 n} - \nu^2 q^{2(n-1)})} \,
A_{1  \to n} \; .
\end{equation}
We need the identities
$$
\sigma_n \dots \sigma_2 \sigma_1^2 \sigma_2
\dots  \sigma_n \, \kappa_n =
 \nu^2 \, \sigma^{-1}_{n-1} \dots \sigma^{-1}_2 \sigma^{-2}_1 \sigma^{-1}_2
\dots \sigma^{-1}_{n-1} \, \kappa_n  \; \Rightarrow
$$
\begin{equation}
\label{bmwa33}
\begin{array}{c}
y_n \, \sigma_n \, y_n \, \kappa_n = \nu \, \kappa_n  \;\; \Rightarrow \;\;
y_{n+1} \, y_n \, \kappa_n = \nu^2 \, \kappa_n\; ,
\end{array}
\end{equation}
which follow from eq. $\sigma_k \kappa_{k+1} = \sigma_{k+1}^{-1}
\kappa_k \kappa_{k+1}$. We also deduce the
analogs of the  identities (\ref{jucmu}) for the
Birman-Murakami-Wenzl algebra case:
\begin{equation}
\label{bmwa5}
\begin{array}{c}
y_{n+1} =
\sigma_n  \dots \sigma_2 \sigma_1^2 \sigma_2 \dots  \sigma_{n} =
1 + \lambda (\sum_{i=1}^{n-1}
\sigma_i \dots \sigma_{n-1} \sigma_n \sigma_{n-1} \dots \sigma_i
 + \sigma_n) - \\ \\
 - \lambda \nu \, ( \sum_{i=1}^{n-1}
\, \sigma^{-1}_i \dots \sigma_{n-1}^{-1} \kappa_n
 \sigma^{-1}_{n-1} \dots \sigma^{-1}_{i} + \kappa_n) \; .
\end{array}
\end{equation}
Using eqs. (\ref{bmwa33}), (\ref{bmwa5}) and $A_{1 \to n} y_n = q^{2 (1-n)}  A_{1 \to n}$ (see eqs. (\ref{anul}) for $A_{1 \to n}$) we obtain
\begin{equation}
\label{bmwa35}
 A_{1 \to n} \, y_{n+1} \, A_{1  \to n}
 = A_{1 \to n} \left(1 +  q \, (1-q^{-2 n}) \, \sigma_n
 + \frac{\nu}{q} \, (1-q^{2 n}) \kappa_n  \right) A_{1 \to n} \; ,
\end{equation}
$$
 A_{1 \to n} \, y_{n+1}^2 \, A_{1  \to n} =
 A_{1 \to n} y_{n+1} \left(1 +  q \, (1-q^{-2 n}) \, \sigma_n
 + \frac{\nu}{q} \, (1-q^{2 n}) \kappa_n  \right)  A_{1 \to n} =
$$
\begin{equation}
\label{bmwa36}
\begin{array}{c}
 = A_{1 \to n} [(1 +  \lambda q  (1-q^{-2 n}))  +q  (1-q^{-2 n}) (q^2 + q^{-2n})
 \sigma_n - \\ \\
  + \frac{\nu}{q} \, (1-q^{2 n}) (q^2 - q^{-2(n-1)} + q^{-2n} + \nu  (\lambda
  + \nu   q^{2(n-1)})  \kappa_n )  ]  A_{1 \to n} \; .
\end{array}
\end{equation}
Then, we substitute (\ref{bmwa35}) and (\ref{bmwa36}) into (\ref{bmwa22}) and
 finally deduce
\begin{equation}
\label{bmwa22'}
A_{1 \to n+1} = \frac{q^{-1}\lambda}{(1-q^{-2(n+1)})} \, A_{1 \to n} \,   \,
\left( 1 + \frac{(q^{-2n}-1)\sigma_n}{\lambda}
+ \frac{\nu (q^{-2n}-1)\kappa_n}{(q^{-2 n+1} + \nu)} \right)
A_{1  \to n} \; ,
\end{equation}
which coincides with (\ref{pyat2}).
  \hfill \qed

  \vspace{0.2cm}

One can  prove directly the identities
(\ref{anul}) for elements (\ref{bmws2}),
(\ref{bmwa2}). We again use the induction.
Let (\ref{anul}) valid for (\ref{bmwa2})
for some fixed $n \geq 2$ (it is obviously correct for $n=2$).
Then, we have to prove
the identities (\ref{anul}) only for $n \to n+1$ and $i=n$.
One can deduce
\begin{equation}
\label{bmwa3}
\begin{array}{c}
A_{1 \to n} \,
(y_{n+1} - \nu^2 q^{2(n-1)}) \kappa_n =
A_{1 \to n} \, (y_{n+1} - \nu^2 \, y_n^{-1}) \kappa_n = 0   \; ,
\end{array}
\end{equation}
where we have applied identities (\ref{bmwa33})
and $A_{1 \to n} y_n = A_{1 \to n} q^{-2(n-1)}$. Using eq.
(\ref{bmwa3}) and the relation
$[ A_{1 \to n}, \, y_{n+1}]=0$ we prove that
$A_{1 \to n+1} \kappa_{n}=0$ for (\ref{bmwa22}).
Now consider the following chain of relations
\begin{equation}
\label{bmwa4}
\begin{array}{c}
A_{1 \to n} \, (y_{n+1} - q^{2}) (\sigma_n + q^{-1}) =
A_{1 \to n} \,
(\sigma_n y_{n} \sigma_n - q^{2}) (\sigma_n + q^{-1}) = \\ \\
= A_{1 \to n} \, ( q \, \sigma_n y_n \sigma_n + \sigma_n y_n -
\lambda \nu \sigma_n y_n \kappa_n - q^2 \sigma_n - q)  =  \\ \\
= A_{1 \to n} \, (- \lambda \nu) \, \kappa_n \,
(\sum_{i=1}^{n-1} (-1)^{n-i} \,
q^{i+1 -n} \sigma^{-1}_{n-1} \dots \sigma_i^{-1} + q + \\ \\
+ \sum_{i=1}^{n-2} (-1)^{n-i} \,
q^{i -n} \kappa_{n-1} \sigma^{-1}_{n-2} \dots \sigma_i^{-1} + q \, \kappa_{n-1}
+ \nu q^{2(1-n)} ) \; ,
\end{array}
\end{equation}
where we have used eqs. (\ref{anul}), (\ref{bmwa33}) and
(\ref{bmwa5}).
Multiplying eq. (\ref{bmwa4})
by the factor $(y_{n+1} - \nu^2 q^{2(n-1)})$
from the left and taking into account (\ref{bmwa3}) we obtain
$A_{1 \to n+1} (\sigma_{n}+ q^{-1}) = 0$.


\vspace{0.2cm}
\noindent
{\bf Remark 1.}
The idempotents
$S_{1 \to n}$ and $A_{1 \to n}$ for the Birman-Murakami-Wenzl algebra
have been also constructed in another form in \cite{HS}. The authors of \cite{HS}
(as well as the authors of \cite{WT})
have not used the Baxterized or Jucys-Murphy elements and, thus, their expressions for $S_{1 \to n}$ and $A_{1 \to n}$
look rather cumbersome. The construction of the primitive
 idempotents $S_{1 \to n}$ and $A_{1 \to n}$ in
 terms of the Baxterized elements
(\ref{bmwbax}) has been proposed by P.N.Pyatov in fall of 2001
 and used, e.g., in \cite{OgPyaSO}. After the substitution of
(\ref{bmwbax}) to
(\ref{bmws}), (\ref{bmwa}) and direct calculations one
can derive the formulas for $S_{1 \to n}$ and
$A_{1 \to n}$ presented in \cite{HS}.

\vspace{0.2cm}
\noindent
{\bf Remark 2.} Assume that the projectors $A_{1 \to n+1}$ (or $S_{1 \to n+1}$) are equal to zero
for some $n$ while
$A_{1 \to n}\neq 0 \neq S_{1 \to n}$. It leads to certain constraints on the
parameter $\nu$. Indeed, from
conditions $\kappa_{n+1} A_{1 \to n+1} \kappa_{n+1} =0$ and $\kappa_{n+1} S_{1 \to n+1}
\kappa_{n+1} =0$ we obtain constraints
$\kappa_{n+1} \sigma^{(+)}_{n}(q^{n}) \kappa_{n+1} = 0$ and
$\kappa_{n+1} \sigma^{(-)}_{n}(q^{- n}) \kappa_{n+1} = 0$, respectively. These constraints are
equivalent to equations $(n>0)$
$$
\begin{array}{c}
(\nu + q^{- (2n+1)})(\nu - q^{- (n-1)})(\nu + q^{- (n-1)}) = 0 \; , \\[0.2cm]
(\nu - q^{(2n+1)})(\nu - q^{(n-1)})(\nu + q^{(n-1)}) = 0 \; .
\end{array}
$$
It means that for $k>n$ all antisymmetrizers $A_{1 \to k}$ could be equal to zero only if
$\nu$ takes one of the values
$\nu = -q^{-(2n+1)}, \pm q^{1-n}$, and, respectively, for $k>n$ all symmetrizers  $S_{1 \to k}$
could be equal to zero only if
$\nu = q^{(2n+1)}, \pm q^{n-1}$. Recall (see Subsection 3.9) that $\nu = q^{1-n}$
and $\nu = -q^{-1-2n}$ specify Birman-Murakami-Wenzl $R$-matrices for $SO_q(n)$ and
$Sp_q(2(n+1))$ groups, respectively. The parameter $\nu = q^{n-1}$ could be related to the
$Osp_q(2(m+1)-n|2m)$ $R$-matrix (\ref{Rosp}) with the choice (\ref{er1}), (\ref{er2}).


\subsubsection[Affine algebras $\alpha BMW_{M+1}$
 and their central elements. Baxterized solution of RE]{Affine algebras $\alpha BMW_{M+1}$
 and their central elements. Baxterized solution of reflection equation\label{abmw1}}

In Subsections {\bf \ref{abmw1}} and {\bf \ref{abmw2}}
 we follow the presentation of the paper \cite{IsOgi}.

Affine Birman-Murakami-Wenzl algebras
$\alpha B\! M \! W_{M+1}(q,\nu)$ are extensions of the algebras
$BMW_{M+1}(q,\nu)$. The algebras  $\alpha B\! M \! W_{M+1}$ are
generated by the elements
$\{\sigma_i,\kappa_i\}$ ($i=1,...,M$)
with relations (\ref{braidg}), (\ref{bmw1}) -- (\ref{bmw3})
and the affine element $y_1$ which satisfies
\be
\lb{bmw02}
\begin{array}{c}
\sigma_1 \, y_1 \, \sigma_1 \, y_1 = y_1 \, \sigma_1 \, y_1 \, \sigma_1 \; , \;\;\;
[\sigma_k, \, y_1]=0 \;\;\; \mathrm{for} \;\; k > 1 \; ,
\\ [0.2cm]
\kappa_1 \; y_1 \, \sigma_1 \, y_1 \, \sigma_1 = c \, \kappa_1 =
\sigma_1 \, y_1 \, \sigma_1 \,  y_{1} \; \kappa_1  \; ,
\\ [0.2cm]
\kappa_1 \; y_1^n \; \kappa_1 = \hat{z}^{(n)} \kappa_1\ \ ,
\;\;\; n=1,2,3,\dots  \; .
\end{array}
\ee
where $c$, $\hat{z}^{(n)}$ are central elements.
Initially, for the Brauer algebras, the affine
version was introduced by M. Nazarov \cite{Nazar}.
Below we use the set of affine elements
\be
\lb{abmw01}
y_1 \; , \;\;\;\;\; y_{k+1} = \sigma_k \,
y_k \, \sigma_k \; \in \;
\alpha B\! M \! W_{M+1} \; , \;\;\; k=1,2,\dots ,M \; .
\ee
These elements generate a commutative subalgebra $Y_{M+1}$ in
$\alpha B\! M \! W_{M+1}$.

We need some information about the center of $\alpha B\!M\!W$.
\begin{proposition}\label{centr1}
{\it The elements
\be
\lb{abmw02}
\hat{\cal Z} = y_1 \cdot y_2 \cdots y_{M} \; ,
\;\;\; \hat{\cal Z}^{(n)}_{M} =
\sum_{k=1}^{M} \left( y^n_{k} - c^n \,
y^{-n}_{k} \right) \;\; ,\;\; n \in {\mathbb N}\ \ ,
\ee
are central in the $\alpha B\!M\!W_{M}$ algebra.}
\end{proposition}
{\bf Proof.} One can directly
 check the centrality of (\ref{abmw02}) by making use
of the relations (\ref{bmw02}) and (\ref{abmw01}). \hfill \qed

\vspace{0.1cm}

\noindent
{\bf Remark 1.} The set of central
"power sums" $\hat{\cal Z}^{(n)} = \sum_k (y_{k}^n - c^n \, y_{k}^{-n})$ is produced by the generating function
$$
{\cal Z}(t) = \sum_{n=1} \hat{\cal Z}^{(n)} t^{n-1}
  = \frac{d}{d t}  \log \left( \prod_{k=1} \frac{y_{k}- c \, t}{1-y_{k}t} \right) \ .
$$

Consider an ascending chain of subalgebras
$$
\alpha B\!M\!W_0\subset \alpha B\!M\!W_1\subset
\alpha B\!M\!W_2\subset\dots\subset
\alpha B\!M\!W_{M}\subset \alpha B\!M\!W_{M+1} \; ,
$$
where $\alpha B\!M\!W_{0},\alpha B\!M\!W_{1}$ and $\alpha B\!M\!W_{j}$ $(j>1)$ are respectively generated by $\{
c, \hat{z}^{(n)}\}$, $\{ c, \hat{z}^{(n)}, y_1\}$ and
$\{ c, \hat{z}^{(n)}, y_1,
 \sigma_1,\sigma_2, \dots , \sigma_{j-1} \}$.
For the corresponding commutative subalgebras we have $Y_1\subset Y_2\subset\dots\subset
Y_{M}\subset Y_{M+1}$.

\begin{proposition}\label{centr2}
{\it Let $\hat{Z}_{k}^{(n)}$ be central elements in the algebra $\alpha B\!M\!W_{k}$, $\alpha B\!M\!W_{k} \subset \alpha B\!M\!W_{k+2}$, defined by the relations
 \be
 \lb{gfun3}
\kappa_{k+1} y_{k+1}^{n} \kappa_{k+1} = \hat{Z}_{k}^{(n)} \, \kappa_{k+1} \in \alpha B\!M\!W_{k+2}
\;\;\;\;\;\; (\hat{Z}_{0}^{(n)} \equiv \hat{z}^{(n)} \; ,
\;\;\; \hat{Z}_{k}^{(0)} \equiv \hat{z}^{(0)} = \mu) \; .
 \ee
Then, the generating function for the elements $\hat{Z}_{k}^{(n)}$
is
\be
\lb{gfun}
\begin{array}{l}
{\displaystyle  \sum_{n=0}^\infty \hat{Z}_{k}^{(n)} t^{n} = - \frac{\nu}{(q-q^{-1})} +
\frac{1}{(1 - c \, t^2)}  +
\left( \sum_{n=0}^{\infty} t^{n} \, \hat{z}^{(n)}  + \frac{\nu}{(q-q^{-1})} - \frac{1}{(1
-c \, t^2)} \right)  }
\\
\hspace{1.2cm}
{\displaystyle \cdot \prod_{r=1}^{k}
\frac{(1-y_{r} t)^2(q^2  - c \, y^{-1}_r t)
(q^{-2}  -c \, y^{-1}_r t)}{
(1-c \, y^{-1}_r t)^2(q^2  - y_r t) (q^{-2}  - y_r t)} } \; .
\end{array}
\ee
}
\end{proposition}
{\bf Proof.} We define the following function of central
elements in $\alpha B\!M\!W_{k}$
$$
Q_k(t) = \sum_{n=0}^\infty \hat{Z}_{k}^{(n)} t^{n} + \frac{\nu}{(q-q^{-1})} - \frac{1}{(1 - c \, t^2)} \; .
$$
Then one can deduce (see the method in \cite{BB})
the recursive formula
 \be
 \lb{gfunk}
Q_k(t) = \frac{(1-y_{k} t)^2(q^2  - c \, y^{-1}_k t)
(q^{-2}  -c \, y^{-1}_k t)}{
(1-c \, y^{-1}_k t)^2(q^2  - y_k t) (q^{-2}  - y_k t)} \;
Q_{k-1}(t) \; ,
 \ee
 where $Q_{k-1}(t) \in \alpha B\!M\!W_{k-1} \subset
 \alpha B\!M\!W_{k}$. From (\ref{gfunk}) we immediately
 obtain (\ref{gfun}). \hfill \qed

 \vspace{0.2cm}

\noindent
 {\bf Remark 2.}
The evaluation map $\alpha B\!M\!W_{M}\rightarrow B\!M\!W_{M}$ is defined by
 \be
 \lb{evmap}
  y_1\mapsto 1 \;\; \Rightarrow \;\; c \mapsto \nu^2 \; , \;\;\;
  \hat{z}^{(n)} \mapsto
  1+\frac{\nu^{-1} -\nu}{q-q^{-1}} \equiv \mu \; .
 \ee
Under this map the function (\ref{gfun}) transforms into the generating function presented in
\cite{BB}, where it is used
for a proof of the Wenzl formula for the quantum dimensions of the
$B\!M\!W_{M}$ primitive idempotents.

\vspace{0.2cm}
\noindent
{\bf Remark 3.} The homomorphisms of the periodic $\overline{B\!M\!W}_{\!\!_{M+1}}$ algebra
to the algebra $B\!M\!W_{\!\!_{M+1}}$ and to the affine algebra
$\alpha B\!M\!W_{\!\!_{M+1}}$
are defined by the same eqs.  (\ref{sigM}) and (\ref{refl4})
as in the case of the group algebra of the braid group.
Indeed, for the periodic $\overline{B\!M\!W}_{\!\!_{M+1}}$
algebra, the characteristic identity for $\sigma_M$ is the same as for $\sigma_1$, while relations
$$
\kappa_1 \sigma^{\pm 1}_M \kappa_1 = \nu^{\mp 1} \kappa_1 \; , \;\;\;
\kappa_{M-1} \sigma^{\pm 1}_M \kappa_{M-1} = \nu^{\mp 1} \kappa_{M-1} \; ,
$$
can be checked directly.

\vspace{0.2cm}
\noindent
{\bf Remark 4.} We redefine the Baxterized
elements in (\ref{bmwbax}), (\ref{bmwbax3})
as follows
 \be
 \lb{bmwba}
 \sigma^{(a)}_i(x) =
 (\sigma_i - x \, \sigma^{-1}_i) +
\frac{\lambda (\nu - a)}{(\nu - a\, x^{-1})}
\, \kappa_i = (a^{-1} - a x^{-1})
\frac{\sigma_i - a \, x}{\sigma_i
- a \, x^{-1}} \; ,
\ee
where we change the spectral parameter
$x^2 \to x$ and denote by
$a$ the solution of the
equation $a^{-1} - a=\lambda \equiv q-q^{-1}$.
It was discovered in \cite{IsaO} that the element of the affine
BMW algebra
\be\lb{rea15}
y_j(u)=f(u) \,
\frac{y_{j} - \xi_a \, u}{y_{j}
- \xi_a \, u^{-1}} \; ,
\ee
(here $\xi_a^2 := a\, c/\nu$ and
 $f(u)$ is any numerical function)
solves the reflection equation
(cf. (\ref{reflH}), (\ref{4.8}))
\be\lb{rea14}
y_j(u)\, \sigma_j\bigl(u\, v\bigr)\,y_j(v) \, \sigma_j\bigl(v \, u^{-1}\bigr)=
\sigma_j\bigl(v \, u^{-1}\bigr)\, y_j(v)\,
\sigma_j\bigl(u\, v\bigr)\, y_j(u) \; .
\ee
This fact
is important in the study of the evaluation homomorphisms for the quantum universal enveloping algebras, see \cite{IsOgM1}
for the classical counterpart.
The main ingredients of the fusion procedure
 \cite{IsOgM2} -- the elements
 $$
 {\cal Y}_j(u_1,\dots ,u_{j-1},u)
 :=\sigma_{j-1}(u \, u_{j-1})\,
{\cal Y}_{j-1}(u_1,\dots ,u_{j-2},u)\, \sigma_{j-1}(u\, u_{j-1}^{-1}) \; ,
 \;\;\;\;\;\;\; {\cal Y}_1(u) : = y_1(u) \; ,
 $$
 $(j=1,\dots,n-1)$, also satisfy the reflection equation
\be\lb{rea14b}
\begin{array}{l}
{\cal Y}_j(u_1,\dots ,u_{j-1},u)\,
\sigma_j \bigl(u\, v\bigr)
\,{\cal Y}_j(u_1,\dots ,u_{j-1},v)\, \sigma_j\bigl(v\, u^{-1}\bigr) =
\\[1em] \hspace{1cm}=
\sigma_j\bigl(v\, u^{-1}\bigr)\, {\cal Y}_j(u_1,\dots ,u_{j-1},v)\,
\sigma_j \bigl(u\, v\bigr) \, {\cal Y}_j(u_1,\dots ,u_{j-1},u)\ .
\end{array}
 \ee
This is shown by induction on $j$.

\subsubsection{Intertwining operators in
$\alpha B\!M\!W_{M+1}$ algebra\label{abmw2}}

Introduce the {\it intertwining} elements
 $U_{k+1} \in  \alpha B\!M\!W_{M+1}$ $(k=1,\dots,M)$
 (cf. (\ref{impint}))
\be
\lb{intw}
U_{k+1} = [\sigma_k , \, y_{k} - c \, y_{k+1}^{-1}]  \; .
\ee

\begin{proposition}\label{intbmw}
{\it The elements $U_k$ satisfy
(cf. (\ref{importt})--(\ref{import}))
$$
U_{k+1} y_k = y_{k+1} \, U_{k+1} \; , \;\;\; U_{k+1} y_{k+1} = y_{k} \, U_{k+1} \; ,
\;\;\;
U_{k+1} y_i = y_{i} \, U_{k+1} \;\;\; {\mathrm{for}}\ \  i \neq k,k+1 \; ,
$$
\be
\lb{bmw05}
U_{k+1} \, [\sigma_k , \, y_{k}] = (q y_k - q^{-1} y_{k+1})(q y_{k+1} - q^{-1} y_{k})
\Bigl(1 - \frac{c}{y_k \, y_{k+1}} \Bigr) \; ,
\ee
$$
U_{k+1} \, U_k \, U_{k+1} = U_k \, U_{k+1} \, U_k  \; ,
$$
$$
\kappa_k \, U_{k+1} = U_{k+1} \, \kappa_k = 0 \; .
$$}
\end{proposition}
The elements $U_{k}$ provide an important information about the spectrum of the affine elements
$\{ y_j \}$ defined in (\ref{abmw01}).

\noindent
{\bf Lemma 3.} (cf. Proposition
{\bf \em \ref{prop9b}}). {\it The spectrum of the elements $y_j \in  \alpha B\!M\!W_{M+1}$ satisfies
\begin{equation}
\label{spec1d}
{\rm Spec} (y_j) \subset \{ q^{2 {\mathbb Z}} \cdot
 {\rm Spec} (y_1) \, , \;\;\; c \, q^{2 {\mathbb Z}}
\cdot {\rm Spec} (y_1^{-1}) \} \; ,
\end{equation}
where ${\mathbb Z}$ is the set of integer numbers.}

\noindent
{\bf Proof.} We prove it by induction in $j$.
 Eq. (\ref{spec1d}) obviously holds for $y_1$. Assume that
$$
{\rm Spec} (y_{j-1}) \subset \{ q^{2 {\mathbb Z}} \cdot {\rm Spec} (y_1)\, , \;\; c \, q^{2 {\mathbb Z}}
\cdot {\rm Spec} (y_1^{-1}) \} \; , \;\; j >1 \; .
$$
Let $f$ be the characteristic polynomial of $y_{j-1}$, $f(y_{j-1}) = 0$. Then
$$
\begin{array}{c}
0=U_j f(y_{j-1}) [\sigma_{j-1} , y_{j-1} ] = f(y_{j}) U_j [\sigma_{j-1} , y_{j-1} ] \\[.5em]
= f(y_{j}) (q^2 y_{j-1} - y_{j})( y_{j} - q^{-2} y_{j-1})
\left( y_{j} - c \, y^{-1}_{j-1}  \right)  y_{j}^{-1} \; .
\end{array}
$$
Here we used (\ref{bmw05}). Thus
${\rm Spec} (y_j) \subset {\rm Spec}(y_{j-1}) \cup q^{\pm 2} \cdot {\rm Spec}(y_{j-1}) \cup
c \cdot {\rm Spec}(y_{j-1}^{-1})$. \hfill \qed

\vspace{0.2cm}

  We denote the image of $w \in \alpha B\!M\!W_{M}$ under the evaluation map (\ref{evmap}) by
$\tilde{w}$, e.g., $y_j \mapsto \tilde{y}_j$.
The Jucys-Murphy (JM) elements $\tilde{y}_j$ $(j=2,\dots,M)$
defined in (\ref{bmwy}) are the images of $y_j$:
$$
\tilde{y}_{j}=\sigma_{j-1}\dots\sigma_2\,\sigma_1^2\,
\sigma_2\dots\sigma_{j-1}\;\in\; B\!M\!W_{M}\; .
$$
Lemma 3 provides the information about the spectrum of JM elements
$\tilde{y}$'s.

\noindent
{\bf Corollary.} Since $\tilde{y}_1=1$ and $\tilde{c}=\nu^2$, it follows from (\ref{spec1d}) that
 \be
\lb{spec2}
{\rm Spec} (\tilde{y}_j) \subset \{ q^{2 {\mathbb Z}} , \; \nu^2 q^{2 {\mathbb Z}} \} \; .
 \ee

\subsection{\bf \em Representation theory
of $B\!M\!W_{M+1}$ algebras\label{repbmw}}

The representation theory for the Birman-Murakami-Wenzl algebra was constructed in \cite{Mur1}
(see also \cite{RaWe}, \cite{LedRam}).
The approach considered in this Section (the colored Young
graph, the analog of Proposition {\bf \ref{prop10}}, 
the explicit formulas for all primitive idempotents in terms
of the Jucys-Murphy elements, intertwiner
 operators $U_k$ (\ref{importt}) -- (\ref{import}), etc.)
 similar to that presented for the Hecke algebra case in Subsection
 {\bf \ref{AHalg}} was
developed in \cite{IsOgi}.


\subsubsection{Representations of affine algebra
$\alpha B\!M\!W_2$}

{\bf A. $\alpha B\!M\!W_2$ algebra and its modules $V_D$}

The elements $\{y_i, y_{i+1}, \sigma_i, \kappa_i \}
\in \alpha B\!M\!W_M$ (for fixed $i < M$) satisfy
\be
\lb{bmw06}
(q-q^{-1}) \kappa_i = \sigma_i^{-1} - \sigma_i + (q-q^{-1}) \; ,
\ee
\be
\lb{bmw08}
y_{i+1} = \sigma_i y_i \sigma_i \; , \;\;\; y_i y_{i+1} = y_{i+1} y_i \; ,
\;\;\;
 \kappa_i y_i^n \kappa_i = \hat{Z}_{i-1}^{(n)} \kappa_i \; ,
\ee
\be
\lb{bmw07}
y_i y_{i+1} \kappa_i = c \, \kappa_i
= \kappa_i y_{i+1} y_i \; .
\ee
The elements $c$ and $\hat{Z}_{i-1}^{(n)}$ commute with
$\{y_i, y_{i+1}, \sigma_i, \kappa_i \}$. The elements $\{y_i, y_{i+1},
\sigma_i, \kappa_i \} \in \alpha B\!M\!W_M$ generate a subalgebra isomorphic to $\alpha B\!M\!W_2 $.

Below we investigate representations $\rho$ of $\alpha B\!M\!W_2$ for which the generators
$\rho(y_i)$ and
$\rho(y_{i+1})$ are diagonalizable and $\rho(c)=\nu^2 \cdot {\text{Id}}$. Let $\psi$
be a common eigenvector of $\rho(y_i)$ and
$\rho(y_{i+1})$ with some eigenvalues $a$ and $b$:
$$
\rho(y_i) \, \psi = a \, \psi \; , \;\;\; \rho(y_{i+1}) \, \psi = b
\, \psi \; .
$$
The element $\hat{z}=y_i y_{i+1}$ is central in $\alpha \! B\!M\!W_2$. There are two possibilities:
\be
\lb{bmw09}
\begin{array}{l}
{\mathbf 1.} \;\;\;\; \rho(\kappa_i) \neq 0 \;\;\; \stackrel{\rm eq. \;
(\ref{bmw07})}{\Longrightarrow}
\;\;\; \rho(y_i y_{i+1})  = \nu^2 \cdot {\text{Id}} \;\; \Rightarrow \;\; \underline{a \, b = \nu^2}; \\ [0.2cm]
{\mathbf 2.} \;\;\;\; \rho(\kappa_i) = 0 \;\; , \;\;\underline{{\mathrm{the\ product}}\ \ a \, b
\;\;  {\rm is} \;
{\rm not} \; {\rm fixed}}.
\end{array}
\ee

Further for brevity we often omit the symbol $\rho$ and denote
 the operator
$\rho(x)$ for $x \in \alpha B\!M\!W$ by the same letter $x$; this should not lead to a confusion.

Applying the operators from $\alpha B\!M\!W_2$ to
the vector $\psi$ we produce, in general infinite
dimensional, $\alpha B\!M\!W_2$-module $V_{\infty}$ spanned by
$$
\begin{array}{ll}
  & \;\;\; e_2 = \psi \; ,\\
e_1 = \kappa_i \psi \; , & \;\;\; e_3 = \sigma_i \psi \; , \\
e_4 = y_i \kappa_i \psi \; , & \;\;\; e_5 = \sigma_i y_i \kappa_i  \psi \; , \\
e_6 = y^2_i \kappa_i \psi \; , & \;\;\; e_7 = \sigma_i y^2_i \kappa_i  \psi \; ,  \\
\;\;\; \dots\dots\dots \; ,  & \;\;\;\;\;\;  \dots\dots\dots  \; , \\
e_{2k +2} =  y^k_i \kappa_i \psi \; , & \;\;\; e_{2k +3} =  \sigma_i y^k_i \kappa_i  \psi \;\;
(k \geq 1) \; , \dots\dots \; .
\end{array}
$$
Using relations (\ref{bmw06}) -- (\ref{bmw07}) for $\alpha B\!M\!W_2$ one can write down the left
action of elements
$\{y_i, y_{i+1}$, $\sigma_i, \kappa_i \}$ on $V_{\infty}$.
Our aim is to understand when the sequence $e_j$ can terminate giving therefore rise to a finite
dimensional module $V_D$ (of dimension  $D$) of $\alpha B\!M\!W_2$ and investigate the
(ir)reducibility of $V_D$.

We distinguish 3 cases for the module $V_D$:
 \begin{itemize}
  \item[\bf (i)]
$\kappa_i V_D=0$ (i.e., $\kappa_i \, e = 0 \;\;
\forall \; e \in V_D$) and in particular $\kappa_i \, \psi = 0$. Therefore, $e_j = 0$
for all $j\neq 2,3$ and  $V_{\infty}$ reduces to a 2-dim module with the basis $\{ e_2,e_3\}$. In
view of (\ref{bmw09}) the product $a \, b$ is not fixed and the irreps coincide with the irreps of
the affine Hecke algebra
$\alpha H_2$ considered in Sect. {\bf \ref{jmel2}}
 and in \cite{IsOg3}.
 \item[\bf (ii)]
$\kappa_i V_D \neq 0$ (i.e., $\exists \; e \in V_D$:
$\kappa_i e \neq 0$). The module $V_D$ is extracted
from $V_{\infty}$ by constraints
 \be
 \lb{bmw10}
e_{2k+4} = \sum_{m=1}^{2k+3} \alpha_m \, e_{m} \;\;\; (k \geq -1) \; , \;\; ab = \nu^2 \; ,
\ee
with some parameters $\alpha_m$. The independent basis vectors are
$(e_1,e_2,\dots,e_{2k+3})$. The module $V_D$ has  \underline{odd dimension}.
 \item[\bf (iii)]
$\kappa_i V_D \neq 0$ and additional constraints are
 \be
 \lb{bmw11}
e_{2k+3} = \sum_{m=1}^{2k+2} \alpha_m \, e_{m} \;\;\; (k \geq 0) \; , \;\; ab = \nu^2 \; .
 \ee
The independent basis vectors are $(e_1,e_2,\dots,e_{2k+2})$. The module
$V_D$ has \underline{even dimension}.
 \end{itemize}
Below we consider a version $\alpha' BMW_2$ of the affine BMW algebra. The additional requirement
for this algebra concerns the spectrum of $y_i,y_{i+1} \in \alpha' BMW_2$:
$$
{\rm Spec} (y_j) \subset \{ q^{2 {\mathbb Z}} , \; \nu^2 q^{2 {\mathbb Z}} \} \; .
  $$
  The evaluation map (\ref{evmap}) descends to the algebra $\alpha' BMW$ (cf. Corollary
  after Lemma 3).
  In particular for the cases {\bf (ii)} and {\bf (iii)} we have
$$
a = \nu^2 q^{2 z}  \; , \;\;\; b = q^{- 2 z}  \;\;\; {\rm or} \;\;\;
 a = q^{2 z} \; ,
\;\;\; b = \nu^2  q^{- 2 z}  \;
$$
for some $z \in {\mathbb Z}$.

\vspace{0.2cm}

\noindent
{\bf B. The case $\kappa_i V_D = 0$:
Hecke algebra case \cite{IsOg3}
(see also Sect. {\bf \ref{jmel2}}).}

{} Representations of $\alpha B\!M\!W_2$ with
$\kappa_i V_D =0$ reduce to representations of the
affine Hecke algebra $\alpha H_2$. In the basis of two vectors
$(e_2,e_3) = (\psi,\sigma_i \psi)$ the matrices
of the generators are (cf. (\ref{abmw1a}))
 \be
 \lb{bmw12}
\sigma_i =
\left(
\begin{array}{cc}
0 & 1 \\
1 & q-q^{-1}
\end{array}
\right)  , \;
y_i =
\left(
\begin{array}{cc}
a & - (q-q^{-1}) b  \\
0 & b
\end{array}
\right) , \;
y_{i+1} =
\left(
\begin{array}{cc}
b & (q-q^{-1}) b  \\
0 & a
\end{array}
\right) ,
 \ee
where $a \neq b$ (otherwise $y_i,y_{i+1}$ are not diagonalizable). By Lemma 3, we have for
$y_i,y_{i+1} \in \alpha' B\!M\!W_2$ the eigenvalues
$a,b \in \{ q^{2 \mathbb Z}, \nu^2 q^{2 \mathbb Z} \}$. The 2-dimensional
representation (\ref{bmw12}) contains a 1-dimensional subrepresentation iff $a = q^{\pm 2} b$.
Graphically these 1- and 2-dimensional irreps of $\alpha' B\!M\!W_2$ are visualized by the same pictures as in Fig.4.2 and Fig.4.3
in Subsection {\bf \ref{jmel2}}.
Different paths going from the upper vertex to the lower vertex correspond to different
eigenvectors of $y_i,y_{i+1}$. The indices on the edges are eigenvalues of
$y_i,y_{i+1}$.

\vspace{0.2cm}

\noindent
{\bf C. $\kappa_i V_D \neq 0$: odd dimensional representations for $\alpha' B\!M\!W_2$}

Using condition (\ref{bmw10}) for the reduction $V_\infty$ to $V_{2m+1}$, one can describe odd
dimensional representations of
$\alpha' B\!M\!W_2$, determine matrices for the
action of $y_i$, $y_{i+1}$ on $V_{2m+1}$ and calculate
\be
\lb{detev}
\det (y_i) = \prod_{r=1}^{2m+1} y^{(r)}_i = \nu^{2m} \; , \;\;\;
\det(y_{i+1}) = \prod_{r=1}^{2m+1} y^{(r)}_{i+1} = \nu^{2m+2} \; .
\ee
Here for eigenvalues $y^{(r)}_i$, $y^{(r)}_{i+1}$ $(r = 1,2, \dots, 2m+1)$ of
$y_i$ and $y_{i+1}$ we have constraints
$$
y^{(r)}_i \, y^{(r)}_{i+1} = \nu^2 \; , \;\;\; r=1, \dots, 2m+1 \;
$$
and (see eq. (\ref{spec2}))
 $$
y^{(r)}_i \in \{ q^{2 {\mathbb Z}} , \; \nu^2 q^{2 {\mathbb Z}} \} \; , \;\; r=1, \dots, 2m+1 \; .
  $$

These odd-dimensional irreps are visualized as graphs:

\vspace{0.2cm}

\setlength{\unitlength}{0.00016in}%

\begin{picture}(18225,15105)(1426,-15322)

\put(14176,-6736){\tiny $\nu^2 \! q^{2z_2}$}
\put(13351,-5461){\tiny $\nu^2 \! q^{2z_1}$}
\put(17126,-6661){\tiny $\nu^2 \! q^{2z_m}$}
\put(21176,-6736){\tiny $q^{2z_{m+1}}$}
\put(26051,-6811){\tiny $q^{2z_{2m}}$}
\put(28376,-6586){\tiny $q^{2z_{2m+1}}$}
\put(27101,-10336){\tiny $\nu^2 \! q^{\! -2z_{2m+1}}$}
\put(21776,-11086){\tiny $\nu^2 \! q^{\! -2z_{2m}}$}
\put(21251,-9511){\tiny $\nu^2 \! q^{\! -2z_{m+1}}$}
\put(19026,-9511){\tiny $q^{\! -2z_{m}}$}
\put(15926,-8936){\tiny $q^{\! -2z_2}$}
\put(13401,-10611){\tiny $q^{\! -2z_1}$}

\thicklines
\put(20571,-606){$\star$}
\put(20651,-8011){$\star$}
\put(29051,-8011){$\star$}
\put(26351,-8011){$\star$}
\put(15526,-7711){. . . . . . .}
\put(21526,-7711){. . . . . . . . . .}
\put(18326,-8011){$\star$}
\put(14426,-8011){$\star$}
\put(12026,-8011){$\star$}
\put(20651,-15486){$\star$}
\put(20651,-736){\line(-6,-5){8225.410}}
\put(20726,-736){\line(-5,-6){5784.836}}
\put(20951,-736){\line( 0,-1){6825}}
\put(21026,-736){\line( 6,-5){8181.148}}
\put(20951,-736){\line( 4,-5){5495.122}}
\put(20801,-811){\line(-1,-3){2212.500}}
\put(12376,-7861){\line( 6,-5){8424.590}}
\put(14801,-7786){\line( 5,-6){6061.475}}
\put(18501,-7711){\line( 1,-3){2407.500}}
\put(20951,-7786){\line( 0,-1){7200}}
\put(26726,-7786){\line(-4,-5){5765.854}}
\put(29551,-7711){\line(-6,-5){8535.246}}

\put(5751,-656){$\star$}
\put(5751,-7786){$\star$}
\put(5751,-15036){$\star$}
\put(6051,-656){\line( 0,-1){14200}}

\put(2501,-5161){$y_i=1$}
\put(1801,-11461){$y_{i+1}=\nu^2$}

\put(28401,-14286){\bf Fig. 4.4}

\end{picture}

\noindent
where $z_r \in {\mathbb Z}$ and $\sum_{r=1}^{2m+1} z_r =0$ as it follows from (\ref{detev}).
Different paths going from the top vertex to the bottom vertex correspond to different
 common eigenvectors of
$y_i,y_{i+1}$. Indices on upper and lower edges of these paths are the eigenvalues of $y_i$ and
$y_{i+1}$, respectively.

\noindent
{\bf Remark.} In view of the braid relations $\sigma_{i} \sigma_{i \pm 1} \sigma_i =
\sigma_{i \pm 1} \sigma_i  \sigma_{i \pm 1}$ and
possible eigenvalues of $\sigma$'s for 1-dimensional representations  (described in Subsections 3.2
and 3.3) we conclude that the following chains of 1-dimensional representations are forbidden

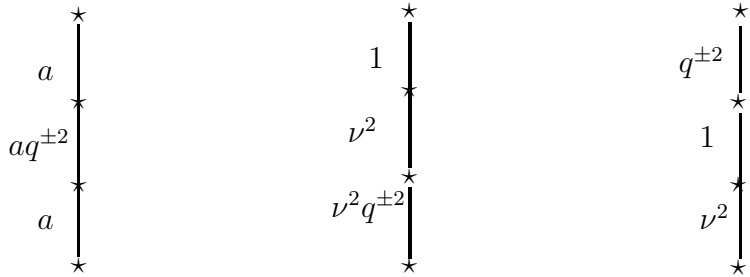
\unitlength=5.5mm
\begin{picture}(20,7.5)
\thicklines
\put(5,1.15){\line(0,1){5.6}}

\put(4.8,0.75){$\star$}
\put(4.8,6.8){$\star$}
\put(4.8,2.7){$\star$}
\put(4.8,4.7){$\star$}

\put(4,1.8){$a$}
\put(3.3,3.6){$a q^{\pm 2}$}
\put(4,5.4){$a$}

\put(13,1.1){\line(0,1){1.7}}
\put(13,3.3){\line(0,1){3.5}}
\put(12.8,0.75){$\star$}
\put(12.8,6.9){$\star$}
\put(12.8,2.9){$\star$}
\put(12.8,5){$\star$}
\put(11.1,2.1){$\nu^2 q^{\pm 2}$}
\put(11.5,3.9){$\nu^{2}$}
\put(12,5.7){$1$}

\put(21,1.1){\line(0,1){3.5}}
\put(21,5.1){\line(0,1){1.6}}
\put(20.75,0.7){$\star$}
\put(20.8,6.9){$\star$}
\put(20.75,2.7){$\star$}
\put(20.75,4.7){$\star$}
\put(20,1.8){$\nu^2 $}
\put(20,3.8){$1$}
\put(19.5,5.6){$q^{\pm 2}$}

\end{picture}

\vspace{0.5cm}
\noindent
where $a= q^{2z}$ or $a = \nu^2 q^{2z}$ $(z \in {\mathbb Z})$.


{\bf D. $\kappa_i V_D \neq 0$: even dimensional representations of $\alpha' B\!M\!W_2$}

With the help of conditions (\ref{bmw11}) we reduce $V_\infty$ to $V_{2m}$, then explicitly
construct
$(2m) \times (2m)$ matrices for the operators $y_i,y_{i+1}$ and calculate their determinants
\be
\lb{bmw14}
\det(y_i) =
\prod_{r=1}^{2m} y_i^{(r)} =  \epsilon q^{\epsilon} \, \nu^{2m-1}  \; , \;\;
\det(y_{i+1}) =
\prod_{r=1}^{2m} y_{i+1}^{(r)} = - \epsilon q^{\epsilon} \, \nu^{2m+1}
\; ,
\ee
where $y_i^{(r)}$, $y_{i+1}^{(r)}$ are eigenvalues of $y_i$, $y_{i+1}$ (we have two possibilities:
$\epsilon=\pm 1$). We see from
(\ref{bmw14}) that all
$(2m)$ eigenvalues of $y_i$, $y_{i+1}$ cannot belong to the spectrum  (\ref{spec2}).
More precisely there is at least one eigenvalue $y_i^{(r)}$ of $y_i$
(and the eigenvalue $y_{i+1}^{(r)}$ of $y_{i+1}$) such that
$$
y_i^{(r)},y_{i+1}^{(r)} \notin \{ q^{2 {\mathbb Z}} , \; \nu^2 q^{2 {\mathbb Z}} \} \; .
$$
Thus, even dimensional irreps of
$\alpha B\!M\!W_2$ subject to the conditions (\ref{bmw11}),
are not admissible for $\alpha' B\!M\!W_2$.

 \subsubsection{Spec$(y_1,\dots,y_{n})$ and rules for
 strings of eigenvalues}

Now we reconstruct the representation theory of $B\!M\!W$ algebras using an approach which
generalizes the approach of Okounkov -- Vershik \cite{OV} for symmetric groups.

The  JM elements {$\{
\tilde{y}_1,\dots,\tilde{y}_{n}\}$} generate a commutative
subalgebra in $B\!M\!W_n$. The basis in the space of an irrep of
$B\!M\!W_n$ can be chosen to be the common eigenbasis of all
$\tilde{y}_i$. Each common  eigenvector $v$ of $\tilde{y}_i$,
$$
\tilde{y}_i \,  v = a_i \, v \; , \;\; i=1,\dots,n \; ,
$$
defines a string $(a_1, \dots,a_n) \in {\mathbb C}^n$. Denote by ${\rm
Spec}(\tilde{y}_1,\dots,\tilde{y}_{n})$ the set of such strings.

We summarize our results about representations of $\alpha' B\!M\!W_2$ and the spectrum of the JM
elements $\tilde{y}_i$ in the following Proposition.

\begin{proposition}\label{spbmw1}
{\it Consider the string
$$
\alpha = (a_1,...,a_i,a_{i+1},...,a_n) \in  {\rm
Spec}(\tilde{y}_1,...,\tilde{y}_i,\tilde{y}_{i+1},...,\tilde{y}_{n}) \; .
$$
Let $v_\alpha$ be the corresponding eigenvector of $\tilde{y}_i$:
$\tilde{y}_i \, v_\alpha = a_i \, v_\alpha$. Then
$$
\!\!\!\!\!\!\! \begin{array}{l}
 {\bf (1)} \;\;\;\; a_i \in \{ q^{2 {\mathbb Z}} ,
\; \nu^2 q^{2 {\mathbb Z}} \} \, ; \\
 {\bf (2)} \;\;\;\; a_i \neq a_{i+1} \; , \;\; i =1,\dots,n-1 \; ; \\
 {\bf (3a)} \;\; a_i a_{i+1} \neq \nu^2 \; ,
 \;\;\; a_{i+1} = q^{\pm 2} \, a_{i} \Rightarrow
\sigma_i \cdot v_\alpha = \pm q^{\pm 1} v_\alpha \, , \;\;
\kappa_i \cdot v_\alpha = 0 \, ; \\ [0.2cm]
 {\bf (3b)} \;\; a_i a_{i+1} \neq \nu^2 \, ,
  \;\;\; a_{i+1} \neq q^{\pm 2} \, a_{i} \Rightarrow \\
\;\;\;\;\;\;\;\;\alpha\, '= (a_1,..., a_{i+1},a_i,...,a_n) \in  {\rm
Spec}(\tilde{y}_1,...,\tilde{y}_i,\tilde{y}_{i+1},...,\tilde{y}_{n})
\, , \; \kappa_i \cdot v_\alpha = 0 \, , \; \kappa_i \cdot v_{\alpha\, '} = 0 \, ; \\ [0.25cm]
 {\bf (4)} \;\;\;\; a_i
a_{i+1} \! = \! \nu^2 \;\; \Rightarrow \;\; \exists \;\; {\rm odd} \; {\rm number} \; {\rm of} \;
{\rm strings} \;\;
\alpha^{(k)} \;\; (k=1,2,\dots, 2m+1): \\ [0.2cm]
\;\;\;\;\;\;\;\;\;
 \alpha^{(k)}=(a_1,...,a_{i-1},a^{(k)}_{i}\!\!
,a^{(k)}_{i+1},a_{i+2},...,a_n) \in {\rm Spec}(\tilde{y}_1,...,\tilde{y}_n) \;\; \forall k \; ,
\\ [0.3cm]
\;\;\;\;\;\;\;\;
\;\; \alpha \in \{\alpha^{(k)}\} \, , \;\; a^{(k)}_{i} a^{(k)}_{i+1}
\! = \! \nu^2 \, , \;\; \prod\limits_{k=1}^{2m+1} a^{(k)}_{i} =
\nu^{2m} \, , \;\; \prod\limits_{k=1}^{2m+1} a^{(k)}_{i+1} =
\nu^{2m+2} \; .
\end{array}
$$}
\end{proposition}

The necessary and sufficient conditions for a string to belong to the common spectrum of
$\tilde{y}_i$ are formulated in the following way.

\begin{proposition}\label{spbmw2}
{\it The string $\alpha =(a_1,a_2, \dots , a_n)$, where $a_i \in
(q^{2\mathbb Z}, \nu^2 q^{2\mathbb Z})$, belongs to the set
${\rm Spec}(\tilde{y}_1,\tilde{y}_2,...,\tilde{y}_{n})$
 iff $\alpha$ satisfies the following conditions ($z \in {\mathbb Z}$)
$$
\begin{array}{l}
{\bf (1)} \;\; a_1 = 1  \; ; \;\;\;\;\;\; {\bf (2)} \;\;\;\; a_i = \nu^2 q^{-2 z}
\;\;
\Rightarrow q^{2z} \in \{a_1, ... , a_{i-1} \} \; ; \\ [0.2cm]
{\bf (3)} \;\; a_i = q^{2z} \Rightarrow
\{ a_i q^2, a_i q^{-2} \} \cap \{a_1, ... , a_{i-1} \} \neq \O
 \; , \;\;  z \neq 0 ; \\ [0.2cm]
\!\!\!\!\!\! {\bf (4a)} \;\; a_i \! = \! a_j \! = \! q^{2z} \; (i < j) \! \Rightarrow
\!\! \left\{ \!
\begin{array}{l}
\!\! {\rm either} \; \{ q^{^{2(z + 1)}},q^{^{2(z - 1)}} \} \! \subset \!
\{a_{_{i+1}}, ... , a_{_{j-1}} \} \; , \\ [0.2cm]
\!\! {\rm or} \; \nu^2 q^{-2z}
\in \{a_{_{i+1}}, ... , a_{_{j-1}} \} \; ;
\end{array}
\right.
\\ [0.6cm]
\!\!\!\!\!\! {\bf (4b)} \; a_i \! = \! a_j \! = \! \nu^2 q^{2z} \, (i < j) \Rightarrow
\!\! \left\{ \!
\begin{array}{l}
\!\! {\rm either} \; \{ \nu^2 q^{^{2(z + 1)}} \!\! , \nu^2 q^{^{2(z - 1)}} \} \! \subset \!
\{a_{_{i+1}}, ... , a_{_{j-1}} \} \; , \\ [0.2cm]
\!\! {\rm or} \;  q^{-2z}
\in \{a_{_{i+1}}, ... , a_{_{j-1}} \} \; ;
\end{array}
\right.
\\ [0.5cm]
\!\!\!\!\!\! {\bf (5a)} \;\; a_i = \nu^2 q^{-2z} \; , \;\; a_j = q^{2 z'} \; (i < j)
\Rightarrow
q^{2z} \; {\rm or} \; \nu^2 q^{-2 z'}\in \{a_{_{i+1}}, ... , a_{_{j-1}} \}  \; ; \\ [0.5cm]
\!\!\!\!\!\! {\bf (5b)} \;\; a_i = q^{2z} \; , \;\; a_j = \nu^2 q^{-2 z'} \; (i < j)
\Rightarrow
\nu^2 q^{-2z} \; {\rm or} \;  q^{2 z'} \in \{a_{_{i+1}}, ... , a_{_{j-1}} \}  \; .
\end{array}
$$
where in {\bf (5a)} and {\bf (5b)} we set $z' = z \pm 1$}.
\end{proposition}

 \subsubsection{Coloured Young graph for $B\!M\!W$ algebras}

We illustrate the above considerations on the example of the colored (in the sense of
\cite{IsOg3}) Young graph for the algebra
$BMW_5$. This graph contains the whole information about the irreps of $BMW_5$ and the branching rules
$BMW_5 \downarrow BMW_4$.

\setlength{\unitlength}{0.00023495in}%
\begingroup\makeatletter\ifx\SetFigFont\undefined%

\fi\endgroup%
\begin{picture}(15266,24019)(4668,-23168)

 \put(23801,339){$\emptyset$}
\put(24151,-436){\tiny $1$}
\put(23476,-3061){\circle*{300}}
\put(23101,-3061){\circle*{300}}
 \put(22951,-2236){\tiny $q^2$}
  \put(22126,-5311){\tiny $q^4$}
    \put(21901,-10636){\tiny $q^{\!-2}$}
 \put(20626,-20361){\tiny $1$}
 \put(22501,-6361){\circle*{300}}
 \put(22126,-6361){\circle*{300}} \put(22876,-6361){\circle*{300}}
\put(24001,-1111){\circle*{300}}

\put(21001,-22186){\circle*{300}} \put(20626,-22186){\circle*{300}}
\put(20251,-22186){\circle*{300}} \put(20251,-22561){\circle*{300}}
\put(20626,-22561){\circle*{300}}

 \put(22501,-11836){\circle*{300}}
 \put(22501,-11461){\circle*{300}}
  \put(22876,-11461){\circle*{300}}
\put(23251,-11461){\circle*{300}}

 \put(22501,-6661){\line(0,-1){4500}}
  \put(24001,239){\line( 0,-1){1050}}
 \put(23251,-2761){\line( 2, 5){600}}

 \put(23176,-3436){\line(-1,-4){675}}

 \put(23026,-11986){\line(-1,-4){2452.941}}

  \put(17026,-7836){\tiny $q^{-2}$}
 \put(15301,-9361){\tiny $q^2$}

 \put(21226,-10261){\tiny $q^6$}
 \put(25876,-4636){\tiny $q^{-4}$}
 \put(24901,-4936){\tiny $q^2$}
   \put(23476,-5911){\tiny $q^{-2}$}

 \put(27301,-8386){\tiny $q^{-6}$}
 \put(24526,-9961){\tiny $1$}
 \put(25201,-8761){\tiny $q^{-4}$}
 \put(26476,-9361){\tiny $q^2$}
 \put(23476,-9436){\tiny $q^4$}

\put(24901,-2236){\tiny $q^{-2}$}
\put(16751,-4561){\tiny $1$}
\put(18301,-20461){\tiny $q^{\!-2}$}
\put(17026,-21136){\tiny $q^8$}
\put(15101,-20861){\tiny $q^{\!-4}$}
\put(13176,-16861){\tiny $q^{\!-2}$}
\put(23551,-20986){\tiny $q^4$}
\put(28076,-16411){\tiny $q^{\!-8}$}
\put(27201,-17461){\tiny $q^2$}
\put(26506,-14886){\tiny $q^{\!-6}$}

\put(21451,-20536){\tiny $q^4$}
\put(22676,-18586){\tiny $q^{\!-4}$}
\put(19651,-20461){\tiny $q^6$}
\put(25801,-17011){\tiny $1$}
\put(24651,-18151){\tiny $q^{\!\!-4}$}
\put(15201,-14986){\tiny $q^2$}
\put(10106,-20556){\tiny $q^{4}$}
\put(10201,-16036){\tiny $1$}

 \thicklines
 \put(24451,-12286){\line(-2,-5){3791.379}}
   \put(17251,-11536){\line(-1,-2){4860}}
 \put(15151,-11311){\line(-1,-2){5085}}

 \put(26176,-12736){\line(-1,-3){3000}}
  \put(24901,-3586){\line(-1,-6){450}}

 \put(24451,-6961){\line(2,-5){1800}}
 \put(26326,-6961){\line( 0,-1){4425}}
 \put(24301,-6961){\line( 0,-1){4425}}
 \put(24151,-6961){\line(-1,-3){1350}}

 \put(22276,-6586){\line(-1,-3){1425}}
 \put(16576,-6586){\line( 1,-6){675}}

 \put(16351,-6586){\line(-1,-4){1050}}
 \put(25051,-3511){\line( 1,-2){1200}}

 \put(23326,-3286){\line(1,-3){975}}
   \put(26326,-12661){\line(-1,-6){1485.811}}

 \put(24601,-12286){\line(0,-1){9150}}
 \put(26401,-12361){\line( 0,-1){9000}}
 \put(27901,-12136){\line(-1,-6){1526.351}}
 \put(27976,-12136){\line(0,-1){9075}}
 \put(16501,-3361){\line( 0,-1){2700}}

   \put(24151,-1261){\line( 2,-5){600}}

 \put(23176,-11911){\line( 0,-1){9750}}

 \put(22801,-11911){\line(-2,-5){3987.931}}
 \put(20626,-11461){\line(-1,-5){2082.692}}
 \put(20176,-11461){\line(-1,-3){3442.500}}
 \put(17326,-11611){\line(-1,-4){2497.059}}

 \put(15451,-11461){\line(-1,-3){3292.500}}
 \put(12451,-11461){\line(-2,-5){3946.552}}
 \put(26626,-6811){\line(1,-3){1200}}

 \put(12401,-11301){$\emptyset$}
 \put(16301,-3161){$\emptyset$}

  \put(21226,-11236){\circle*{300}}
 \put(24301,-11686){\circle*{300}}
\put(24301,-12061){\circle*{300}} \put(24676,-11686){\circle*{300}}
 \put(24676,-12061){\circle*{300}}
 \put(26101,-11761){\circle*{300}}
\put(26101,-12136){\circle*{300}} \put(26101,-12511){\circle*{300}}
\put(26476,-11761){\circle*{300}} \put(27901,-10636){\circle*{300}}
\put(27901,-11011){\circle*{300}} \put(27901,-11386){\circle*{300}}
\put(27901,-11761){\circle*{300}} \put(9826,-21586){\circle*{300}}
\put(12676,-21661){\circle*{300}} \put(28276,-21511){\circle*{300}}
\put(28276,-21886){\circle*{300}} \put(20851,-11236){\circle*{300}}

   \put(24901,-2911){\circle*{300}}
\put(24901,-3286){\circle*{300}}

\put(24301,-6361){\circle*{300}} \put(24676,-6361){\circle*{300}}
\put(24301,-6736){\circle*{300}} \put(26401,-6361){\circle*{300}}
\put(26401,-5986){\circle*{300}} \put(26401,-6736){\circle*{300}}
\put(16501,-6361){\circle*{300}} \put(15751,-11161){\circle*{300}}
\put(15376,-11161){\circle*{300}} \put(17401,-11011){\circle*{300}}
\put(17401,-11386){\circle*{300}} \put(20101,-11236){\circle*{300}}
\put(20476,-11236){\circle*{300}} \put(28276,-22261){\circle*{300}}
\put(19126,-22186){\circle*{300}} \put(18001,-22186){\circle*{300}}
\put(18001,-22561){\circle*{300}} \put(17026,-22036){\circle*{300}}
\put(16651,-22036){\circle*{300}} \put(16276,-22036){\circle*{300}}
\put(15901,-22036){\circle*{300}} \put(15526,-22036){\circle*{300}}
\put(14476,-21811){\circle*{300}} \put(14476,-22186){\circle*{300}}
\put(14476,-22561){\circle*{300}} \put(12301,-21661){\circle*{300}}
\put(12301,-22036){\circle*{300}} \put(10201,-21586){\circle*{300}}
\put(10576,-21586){\circle*{300}} \put(8326,-21586){\circle*{300}}
\put(18376,-22186){\circle*{300}} \put(28276,-22636){\circle*{300}}
\put(28276,-23011){\circle*{300}} \put(26326,-21661){\circle*{300}}
\put(24601,-21736){\circle*{300}} \put(24601,-22111){\circle*{300}}
\put(24226,-21736){\circle*{300}} \put(24226,-22111){\circle*{300}}
\put(24226,-22486){\circle*{300}} \put(23101,-21961){\circle*{300}}
 \put(22726,-21961){\circle*{300}} \put(22351,-21961){\circle*{300}}
 \put(22351,-22336){\circle*{300}} \put(22351,-22711){\circle*{300}}

 \put(26326,-22036){\circle*{300}}
\put(26326,-22411){\circle*{300}} \put(26326,-22786){\circle*{300}}
\put(26701,-21661){\circle*{300}} \put(18751,-22186){\circle*{300}}

 \put(24676,-3586){\line(-3,-1){7852.500}}
  \put(24076,-6811){\line(-5,-3){6529.412}}
   \put(26326,-12661){\line(-5,-4){11485.811}}
 \put(24076,-12136){\line(-6,-5){11161.475}}
 \put(22276,-11836){\line(-1,-1){9525}}
  \put(26032,-12820){\line(-3,-2){13309.615}}
 \put(20251,-11461){\line(-1,-1){9780}}
  \put(17776,-10861){\line( 2,1){8400}}
 \put(12676,-10861){\line( 5, 6){3750}}
 \put(15676,-10861){\line( 2, 1){8400}}
 \put(15451,-10861){\line( 3,2){6525}}
  \put(16651,-5986){\line( 5, 2){6375}}
   \put(22201,-11836){\line(-6,-5){11419.672}}
    \put(16976,-11461){\line(-5,-6){8250}}
 \put(15076,-11386){\line(-2,-3){6600}}
  \put(27976,-12136){\line( -4,-3){13075}}
  \put(16801,-2911){\line( 4,1){6900}}

\put(19051,-4411){\tiny $\nu^2q^{-2}$}
\put(15726,-21536){\tiny $\nu^2\!q^{6}$}
\put(16001,-20361){\tiny $\nu^2\!q^{\!-2}$}
\put(13701,-19636){\tiny $\nu^2\!q^{\!\!-4}$}
 \put(18976,-10536){\tiny $\nu^2 q^4$}
 \put(17576,-10356){\tiny $\nu^2\! q^{\!-2}$}
 \put(20246,-8236){\tiny $\nu^2q^2$}
 \put(18376,-8236){\tiny $\nu^2q^{\!-4}$}
 \put(14101,-8311){\tiny $\nu^2$}
 \put(19651,-5611){\tiny $\nu^2q^2$}
 \put(19876,-1786){\tiny $\nu^2$}
 \put(13351,-21486){\tiny $\nu^2q^4$}
\put(11326,-21136){\tiny $\nu^2q^2$}
\put(11651,-19300){\tiny $\nu^2\!q^{\!-6}$}
\put(9026,-21211){\tiny $\nu^2\!\!q^2$}
\put(11401,-15661){\tiny $\nu^2\!q^{\!-2}$}
\put(16576,-18661){\tiny $\nu^2$}

\put(4551,-20786){{\bf Fig.4.5}}
\end{picture}

\noindent
A vertex $\{\lambda;5\}$ on the lowest level of this graph is labeled by some Young diagram
$\lambda$; this vertex corresponds to the irrep $W_{\{\lambda;5\}}$ of
$BMW_5$ (the notation $\{\lambda;5\}$ is designed to encode the diagram $\lambda$ and the level
on which this diagram is located; the levels are counted starting from 0). Paths going down from
the top vertex $\emptyset$ to the lowest level
(that is, paths of length 5) correspond to common eigenvectors of the JM elements $\tilde{y}_1,
\dots , \tilde{y}_5$. Paths ending at $\{\lambda;5\}$ label the basis in $W_{\{\lambda;5\}}$.
In particular, the number of different paths going down from the top $\emptyset$ to
$\{\lambda;5\}$ is equal to the dimension
of the irrep $W_{\{\lambda;5\}}$.

Note that the colored Young graph in Fig.4.5 contains subgraphs presented in Fig.4.2, Fig.4.3 and Fig.4.4.
For example, in Fig.4.5 one recognizes rhombic subgraphs (the vertices on the subgraphs are obtained
from one another by a rotation)

\vspace{0.3cm}

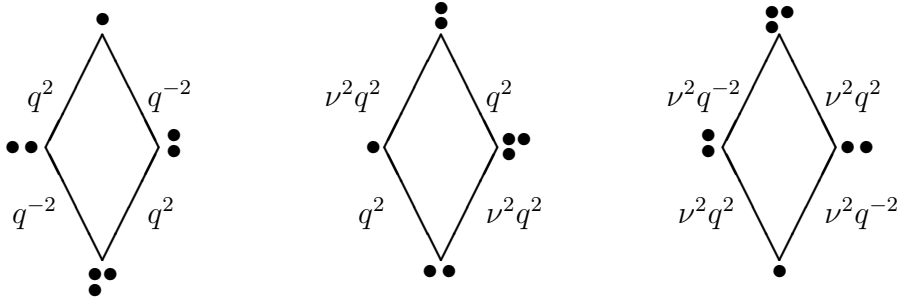
\unitlength=5mm
\begin{picture}(20,7.5)

\thicklines

\put(1.6,2){$q^{-2}$}
 \put(2,5.1){$q^{2}$}
 \put(5.2,5.1){$q^{-2}$}
 \put(5.2,2){$q^{2}$}

\put(5.7,4.1){$\bullet$}
\put(5.7,3.7){$\bullet$}

\put(1.4,3.8){$\bullet$}
\put(1.9,3.8){$\bullet$}

\put(3.6,0.4){$\bullet$}
\put(3.6,0.0){$\bullet$}
\put(4,0.4){$\bullet$}

\put(4,7){\line(-1,-2){1.5}}
\put(4,7){\line(1,-2){1.5}}
\put(4,1){\line(-1,2){1.5}}
\put(4,1){\line(1,2){1.5}}
\put(3.8,7.2){$\bullet$}


\put(10.8,2){$q^{2}$}
 \put(9.9,5.1){$\nu^2 q^{2}$}
 \put(14.2,5.1){$q^{2}$}
 \put(14.2,2){$\nu^2 q^{2}$}
 \put(11,3.8){$\bullet$}

\put(14.6,4){$\bullet$}
\put(14.6,3.6){$\bullet$}
\put(15,4){$\bullet$}

\put(12.5,0.5){$\bullet$}
\put(13,0.5){$\bullet$}

\put(13,7){\line(-1,-2){1.5}}
\put(13,7){\line(1,-2){1.5}}
\put(13,1){\line(-1,2){1.5}}
\put(13,1){\line(1,2){1.5}}

\put(12.8,7.5){$\bullet$}
\put(12.8,7.1){$\bullet$}


\put(19.3,2){$\nu^2 q^{2}$}
 \put(19,5.1){$\nu^2 q^{-2}$}
 \put(23.2,5.1){$\nu^2 q^{2}$}
 \put(23.2,2){$\nu^2 q^{-2}$}

\put(19.9,4.1){$\bullet$}
\put(19.9,3.7){$\bullet$}

\put(23.6,3.8){$\bullet$}
\put(24.1,3.8){$\bullet$}

\put(21.6,7.4){$\bullet$}
\put(21.6,7.0){$\bullet$}
\put(22,7.4){$\bullet$}

\put(22,7){\line(-1,-2){1.5}}
\put(22,7){\line(1,-2){1.5}}
\put(22,1){\line(-1,2){1.5}}
\put(22,1){\line(1,2){1.5}}
\put(21.8,0.5){$\bullet$}

\end{picture}

\vspace{0.3cm}

\noindent
of the type presented in Fig.4.3.

Let $(s,t)$ be coordinates of a node in the Young diagram $\lambda$. To the node $(s,t)$ of the
diagram $\lambda$ we associate a number
$q^{2(s-t)}$ which is called "content":

\vspace{0.5cm}

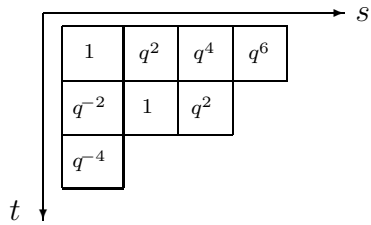
\unitlength=5mm
\begin{picture}(25,4)
\put(7,4){\vector(1,0){8}}
\put(7,4){\vector(0,-1){5.5}}
\put(15.3,3.8){$s$}
\put(6.1,-1.5){$t$}

\put(7.5,1.3){$
\begin{tabular}{|c|c|c|c|}
\hline
  $\!\!\! $  $_1$ & $\!\! $  $_{q^2}$  &
  $\!\! $  $_{q^4}$ & $\!\! $  $_{q^6}$ \\[0.2cm]
\hline
  $\!\!\! $ $ _{q^{\! -2}}$  &   $\!\!\!$ $_1$  & $\!\!\! $ $_{q^2}$  &
  \multicolumn{1}{c}{} \\[0.2cm]
\cline{1-3}
  $\!\!\! $ $ _{q^{\! -4}}$  & \multicolumn{1}{c}{} \\[0.2cm]
\cline{1-1}
\end{tabular}
$}
\end{picture}

\vspace{1cm}

\noindent
Then according to the colored Young graph on Fig.4.5, at each step down along the path one can add or
remove one node (therefore this graph is called the "oscillating" Young graph) and the eigenvalue
of the corresponding JM element is determined by the content of the node:

\vspace{0.2cm}

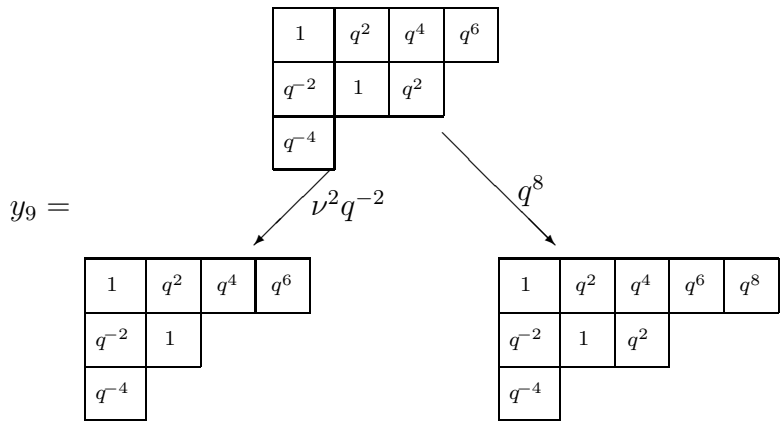
\unitlength=5mm
\begin{picture}(25,4)

\put(7.5,1.3){$
\begin{tabular}{|c|c|c|c|}
\hline
  $\!\!\! $  $_1$ & $\!\! $  $_{q^2}$  &
  $\!\! $  $_{q^4}$ & $\!\! $  $_{q^6}$ \\[0.2cm]
\hline
  $\!\!\! $ $ _{q^{\! -2}}$  &   $\!\!\!$ $_1$  & $\!\!\! $ $_{q^2}$  &
  \multicolumn{1}{c}{} \\[0.2cm]
\cline{1-3}
  $\!\!\! $ $ _{q^{\! -4}}$  & \multicolumn{1}{c}{} \\[0.2cm]
\cline{1-1}
\end{tabular}
$}
\end{picture}

\vspace{1.3cm}

\unitlength=5mm
\begin{picture}(25,4)

\put(9,6){\vector(-1,-1){2}}
\put(12,7){\vector(1,-1){3}}

\put(0.5,4.8){$y_9 =$}
\put(8.5,4.8){$\nu^2 q^{-2}$}
\put(14,5.2){$q^8$}

\put(2.5,1.3){$
\begin{tabular}{|c|c|c|c|}
\hline
  $\!\!\! $  $_1$ & $\!\! $  $_{q^2}$  &
  $\!\! $  $_{q^4}$ & $\!\! $  $_{q^6}$ \\[0.2cm]
\hline
  $\!\!\! $ $ _{q^{\! -2}}$  &   $\!\!\!$ $_1$  &
  \multicolumn{1}{c}{} \\[0.2cm]
\cline{1-2}
  $\!\!\! $ $ _{q^{\! -4}}$  & \multicolumn{1}{c}{} \\[0.2cm]
\cline{1-1}
\end{tabular}
$}

\put(13.5,1.3){$
\begin{tabular}{|c|c|c|c|c|}
\hline
  $\!\!\! $  $_1$ & $\!\! $  $_{q^2}$  &
  $\!\! $  $_{q^4}$ & $\!\! $  $_{q^6}$ & $\!\! $  $_{q^8}$ \\[0.2cm]
\hline
  $\!\!\! $ $ _{q^{\! -2}}$  &   $\!\!\!$ $_1$  & $\!\!\! $ $_{q^2}$  &
  \multicolumn{1}{c}{} \\[0.2cm]
\cline{1-3}
  $\!\!\! $ $ _{q^{\! -4}}$  & \multicolumn{1}{c}{} \\[0.2cm]
\cline{1-1}
\end{tabular}
$}
\end{picture}

\vspace{1cm}

\noindent
The eigenvalue corresponding to the addition or removal of the $(s,t)$ node is $q^{2(s-t)}$ or
$\nu^2 q^{-2(s-t)}$, respectively.

 Let $X(n)$ be the set of paths of length $n$ starting from the top vertex $\emptyset$ and going down in
the Young graph of oscillating Young diagrams. Now we formulate the following Proposition.

 \begin{proposition}\label{spbmw3}
 {\it There is a bijection between the set
$\, {\mathrm Spec} (\tilde{y}_1, \dots , \tilde{y}_n)$
 and the set $X(n)$.}
\end{proposition}

 \subsubsection{Primitive idempotents}

The colored Young graph (as in Fig. 4.5) gives also the rule of construction of a complete set of
 orthogonal primitive idempotents for the
$BMW$ algebra. The completeness of the set of orthogonal primitive idempotents is equivalent to the
maximality of the commutative set of JM elements.
Let $\{{\sf \Lambda};n\}$ be a vertex in the Young
graph with

{\unitlength=4mm
\begin{picture}(25,6.5)(-2,0)
\put(-0.5,3){${\sf \Lambda} \;\; = $} \put(4,5.5){\line(1,0){5}}
\put(4,4){\line(1,0){5}} \put(4,4){\line(0,1){1.5}}
\put(9,4){\line(0,1){1.5}}

\put(4,2.5){\line(1,0){3}} \put(7,2.5){\line(0,1){1.5}}
\put(4,4){\line(1,0){3}} \put(4,1){\line(0,1){3}}
\put(6,1.5){\line(0,1){1}} \put(4,1.5){\line(1,0){2}}
\put(4.5,1){$\dots$}

\put(4,0){\line(0,1){0.5}} \put(5,0){\line(0,1){0.5}}
\put(4,0){\line(1,0){1}} \put(4,0.5){\line(1,0){1}}

\put(6,6.2){$_{_{\lambda_{_{(1)}}}}$} \put(3,4.7){$_{n_{_1}}$}
\put(9.1,4){$_{n_{_1},} {_{\lambda_{_{(1)}}}}$}
\put(7.1,2.5){$_{n_{_2},} {_{\lambda_{_{(2)}}}}$}
\put(6.2,1.5){$_{n_{_3},} {_{\lambda_{_{(3)}}}}$}
\put(5.1,0){$_{n_{_k},} {_{\lambda_{_{(k)}}}}$}

 \put(13,3.4){$(n_i,\lambda_{(i)})$ are coordinates of the nodes}
\put(13.5,2.2){which are in the corners of the diagram}
\put(13.5,0.8){${\sf \Lambda}=[\lambda_{(1)}^{n_1},
\lambda_{(2)}^{n_2-n_1},\dots,\lambda_{(k)}^{n_k-n_{k-1}}]$ .}

\end{picture}}

\vspace{-1.5cm}
\be
\lb{lambjj}
{}
\ee

\vspace{0.7cm}

\noindent Consider any path  $T_{\{{\sf \Lambda};n\}}$ going down from the top $\emptyset$ to this vertex. Let
$E_{T_{\{{\sf \Lambda};n\}}} \in B\! M\! W_{n}$ be the primitive
idempotent corresponding to $T_{\{{\sf \Lambda};n\}}$. Using the branching rule implied by the Young
graph for
 $B\! M\! W_{n+1}$ we know all possible eigenvalues of the
 element $\tilde{y}_{n+1}$ and, therefore, obtain the identity
$$
E_{T_{\{{\sf \Lambda};n\}}} \cdot \prod_{r=1}^{k+1} \left(
\tilde{y}_{n+1} - q^{2(\lambda_{(r)} - n_{r-1})} \right)
\prod_{r=1}^{k} \left( \tilde{y}_{n+1} - \nu^2 q^{2(n_{r} -
\lambda_{(r)})} \right)= 0 \; ,
$$
where  $\lambda_{(k+1)} = n_0 = 0$. So, for a new diagram
${\sf \Lambda}'$ obtained by adding to
${\sf \Lambda}$ a new node with coordinates $(n_{j-1} +1,
\lambda_{(j)}+1)$ the corresponding primitive idempotent (after an appropriate normalization) reads
$$
E_{T_{\{{\sf \Lambda}';n+1\}}} =  E_{T_{\{{\sf \Lambda};n\}}} \cdot \prod_{\stackrel{r=1}{_{r \neq j}}}^{k+1}
\! \frac{\left( \tilde{y}_{_{n+1}} - q^{2(\lambda_{(r)} -
n_{r-1})} \right)}{\left( q^{2(\lambda_{(j)} - n_{j-1})} -
q^{2(\lambda_{(r)} - n_{r-1})} \right)}
 \prod_{r=1}^{k} \! \frac{\left(
\tilde{y}_{_{n+1}} - \nu^2 q^{2(n_{r} - \lambda_{(r)} )}
\right)}{\left( q^{2(\lambda_{(j)} - n_{j-1})} - \nu^2 q^{2(n_{r} -
\lambda_{(r)} )} \right)}\,  .
 $$
For a new diagram ${\sf \Lambda}''$ which is obtained from
${\sf \Lambda}$ by removing a node with coordinates
$(n_{j},\lambda_{(j)})$ we construct the primitive idempotent
$$
E_{T_{\{{\sf \Lambda}'';n+1\}}} =  E_{T_{\{{\sf \Lambda};n\}}} \cdot \prod_{r=1}^{k+1}
\! \frac{\left( \tilde{y}_{_{n+1}} - q^{2(\lambda_{(r)} -
n_{r-1})} \right)}{\left( \nu^2 q^{2(n_{j}-\lambda_{(j)})} - q^{2(\lambda_{(r)} - n_{r-1})}
\right)}
 \prod_{\stackrel{r=1}{_{r \neq j}}}^{k} \! \frac{\left(
\tilde{y}_{_{n+1}} - \nu^2 q^{2(n_{r} - \lambda_{(r)} )}
\right)}{\left( \nu^2 q^{2(n_{j}-\lambda_{(j)})} - \nu^2 q^{2(n_{r} -
\lambda_{(r)} )} \right)}\,  .
 $$
 Using these
formulas and the "initial data" $E_{T_{\{\emptyset;0\}}} =1$, one can deduce step by step explicit
expressions for the primitive orthogonal idempotents related to the paths in the BMW Young graph.

\noindent
{\bf Remark.}
In this subsection we reconstructed the representation
 theory of the tower of the $BMW$ algebras using
the properties of the commutative subalgebras
 generated by the Jucys--Murphy elements in the
$BMW$ algebras. This representation theory is of use in the representation theory of the quantum
groups
$U_q(osp(N|K))$ due to the Brauer - Schur - Weyl duality, but
 also finds applications in physical models.
Recently \cite{IsaO}, we have formulated integrable chain models with nontrivial boundary
conditions in terms of the affine Hecke algebras $\hat{H}_n$ and the affine BMW algebras $\alpha BMW_n$. The Hamiltonians for these models are special elements of the algebras $\hat{H}_n$ and
$\alpha BMW_n$. E.g., for the $\alpha BMW_n$ algebra we deduced
\cite{IsaO} the Hamiltonians
\be
\lb{ham5}
{\cal H} = \sum_{m=1}^{n-1} \left( \sigma_m + \frac{(q-q^{-1}) \nu}{\nu +a} \kappa_m \right) +
\frac{(q-q^{-1}) \xi}{y_1 - \xi}  \; ,
\ee
where $\xi^2 = -a \, c/\nu$ and the parameter $a$ can take one of two values $a=\pm q^{\pm 1}$. Now
different representations $\rho$ of the
algebra $\alpha BMW_n$ give different integrable spin
chain models with Hamiltonians $\rho({\cal H})$ which in particular possess $U_q(osp(N|K))$
symmetries for some $N$ and $K$. So, the representations
$\rho$ of the algebra $\alpha BMW_n$ are related to the spin chain models of $osp$ - type with $n$ sites
and nontrivial boundary conditions. The BMW chains (chains based on the BMW algebras in the $R$-matrix representations) describe in a
unified way spin chains with $U_q(osp(N|K))$ symmetries.

The Hamiltonians for the Hecke chain models are obtained
from the Hamiltonians for the BMW chain models by
taking the quotient $\kappa_j =0$. These models were considered in \cite{IsOg2}, \cite{Isa7}. The
Hecke chains (chain models based on the Hecke algebras) describe in a unified way spin chains
with
$U_q(sl(N|K))$ symmetries. In \cite{IsOg2}, \cite{Isa7},
 \cite{BIsKrN}, we investigated the integrable open
chain models formulated in terms of the generators
of the Hecke algebra (non-affine case, $y_1=1$).
  For the open Hecke chains of a finite size,
the spectrum of the Hamiltonians with free boundary conditions is determined \cite{IsOg2} for
special (corner type) irreducible representations of the Hecke algebra. In
\cite{Isa7}, we investigated the functional equations for the transfer matrix type elements of the Hecke
algebra appearing in the theory of Hecke chains.


\subsubsection{q-Dimensions of idempotents
in the $BMW$ algebra}

Consider the following inclusions of the subalgebras
$\alpha\! BMW_{1} \subset \alpha\! BMW_{2} \subset \dots \subset \alpha\! BMW_{n+1}$:
$$
\{y_1; \sigma_1, \dots ,\sigma_{k-1}\} \in
 \alpha\! BMW_{k} \subset \alpha\! BMW_{k+1} \ni
\{y_1; \sigma_1, \dots ,\sigma_{k-1},\sigma_k \} \; .
$$
For the subalgebras $\alpha\! BMW_{k+1}$
we introduce linear mapping (quantum trace)
$$
{\rm Tr}_{(k+1)}: \;\; \alpha\! BMW_{k+1} \to \alpha\! BMW_{k} \; ,
\;\;\;\; (k=1,2,\dots,n) \; ,
$$
which is defined by the formula (cf. (\ref{3.7.29}))
 \be
 \lb{qtr}
\kappa_{k+1} \, X_{k+1} \, \kappa_{k+1} = \frac{1}{\nu} \, Tr_{({k+1})}(X_{k+1}) \,
\kappa_{k+1} \; , \;\;\; \forall X_{k+1} \in \alpha\! BMW_{k+1} \; .
 \ee
  \begin{proposition}\label{spbmw4}
  {\it (see \cite{IsKiTa})
 For the map $Tr_{(k+1)}$: $\alpha\! BMW_{k+1} \to \alpha\! BMW_{k}$
 we have the following properties $(\forall X_{k},X_k' \; \in \; \alpha\! BMW_{k} \; ,
 \;\; \forall Y_{k+1} \; \in \; \alpha\! BMW_{k+1})$
 \be
 \lb{qtr1}
 \begin{array}{c}
 Tr_{(k+1)}(\sigma_{k})=1  \; , \;\;
 Tr_{(k+1)}(\sigma_{k}^{-1})=\nu^2 \; , \;\;
   {\rm Tr}_{(k+1)} ( X_k ) = \nu \, \mu \, X_k  \; , \\[0.1cm]
   \displaystyle
     Tr_{(k+1)}(\kappa_k)= \nu  \; ,
   \;\; {\rm Tr}_{(1)} (y_1^k)= \nu \, \hat{z}^{(k)} \; ,
   \;\; {\rm Tr}_{(1)} (1)= \nu \, \mu =
   (1+ \lambda \nu -\nu^2) \lambda^{-1} \; ,
 \end{array}
 \ee
 \be
 \lb{qtr2}
 Tr_{(k+1)}(\sigma_{k} \, X_{k} \, \sigma_{k}^{-1}) =
 Tr_{(k)}(X_{k})  =
  Tr_{(k+1)}(\sigma_{k}^{-1} \, X_{k} \, \sigma_{k})   \; ,
 \ee
 \be
 \lb{qtr3}
 Tr_{(k+1)}(\sigma_{k} \, X_{k} \, \kappa_{k}) =
  Tr_{(k+1)}(\kappa_{k} \, X_{k} \, \sigma_{k}) \; ,
 \ee
 \be
\lb{mapn}
\begin{array}{c}
{\rm Tr}_{(k+1)} ( X_k \cdot Y_{k+1} \cdot X'_k )
= X_k \cdot {\rm Tr}_{(k+1)}(Y_{k+1}) \cdot X_k' \;\;  \, , \\[0.1cm]
{\rm Tr}_{(k)}  {\rm Tr}_{(k+1)} ( \sigma_k \cdot Y_{k+1} ) = {\rm Tr}_{(k)}  {\rm Tr}_{(k+1)}
( Y_{k+1} \cdot \sigma_k)  \; .
\end{array}
\ee}
\end{proposition}

By using the mapping (\ref{qtr}), definitions (\ref{gfun3}),
(\ref{bmw02}) and evaluation (\ref{evmap}),
we write relation
(\ref{gfun}) in the form (cf. (\ref{ident5}))
\be
\lb{gfun5}
\begin{array}{c}
{\displaystyle
\lambda  {\rm Tr}_{_{(M+1)}}
 \Bigl( \frac{y_{M+1}}{t - y_{M+1}} \Bigr)
+ 1 - \frac{\lambda \nu^3}{(t^2 - \nu^2)}  = }\\ [0.4cm]
= {\displaystyle
\frac{(t-\nu^2 )(t-q^{-1} \nu)(t+q \nu)}{(t-1)(t-\nu)(t+\nu)}
\cdot \prod_{r=1}^{M}
\frac{(t-y_{r} )^2(q^2 t - \nu^2 \, y^{-1}_r )
(q^{-2}t  -\nu^2 \, y^{-1}_r )}{
(t-\nu^2 \, y^{-1}_r )^2(q^2 t - y_r ) (q^{-2} t - y_r )} } \; ,
\end{array}
\ee
where we change variable $t \to t^{-1}$, index
 $k \to M$ and for simplicity denote
 $y_r = \tilde{y}_r$. Then, we act to both sides of
 (\ref{gfun5}) by the idempotent
 $E_{T_{\{{\sf \Lambda},M\}}}$, where
 $T_{\{{\sf \Lambda},M\}}$ is the path of length $M$ in
 the coloured Young graph (of the type Fig.4.5) with the
 final vertex labeled by the Young diagram
 ${\sf \Lambda}$ (\ref{lambjj}). According to
 the branching rule, which is implied by the
 coloured Young graph for $B\! M\! W_{M+1}$, we use
 the expansion
 $$
 E_{T_{\{{\sf \Lambda},M\}}} =
 \sum_{j=1}^{k+1}  E_{T_{\{{\sf \Lambda}_j',M+1\}}}
 + \sum_{j=1}^k E_{T_{\{{\sf \Lambda}_j'',M+1\}}} \; ,
 $$
 where the Young diagram ${\sf \Lambda}'_j$ is obtained
 by adding a node to the outer
 corners $(n_{j-1}+1, \lambda_{(j)}+1)|_{j=1,...,k+1}$
 of the diagram ${\sf \Lambda}$ and
 the  diagram ${\sf \Lambda}_j''$ is obtained
 by removing a node from inner corners
 $(n_{j}, \lambda_{(j)})|_{j=1,...,k}$
 of ${\sf \Lambda}$.
  As a result we obtain
 \be
\lb{gfun6}
\begin{array}{c}
{\displaystyle
  {\rm Tr}_{_{(M+1)}}
 \Bigl( \lambda \frac{y_{M+1}}{t - y_{M+1}}
 \Bigl(\sum_{j=1}^{k+1}  E_{T_{\{{\sf \Lambda}_j',M+1\}}}
 + \sum_{j=1}^k E_{T_{\{{\sf \Lambda}_j'',M+1\}}} \Bigr)\Bigr)
 = }\\ [0.4cm]
{\displaystyle \Bigl(
\frac{(t-\nu^2 )(t-q^{-1} \nu)(t+q \nu)}{(t-1)(t-\nu)(t+\nu)}
\cdot \prod_{r=1}^{M}
\frac{(t-y_{r} )^2 (t - q^{-2}\nu^2 \, y^{-1}_r )
(t  - q^2 \nu^2 \, y^{-1}_r )}{
( t - q^{-2}  y_r ) (t - q^2 y_r )(t-\nu^2 \, y^{-1}_r )^2}
- c(t) \Bigr) E_{T_{\{{\sf \Lambda},M\}}}  } \; ,
\end{array}
\ee
where $c(t) := 1 - \frac{\lambda \nu^3}{(t^2 - \nu^2)}$.
Now we note that, if a cell with
content $q^{2a}$ was added on the step $i$
in the path $T_{\{{\sf \Lambda},M\}}$
and further this cell was removed on the step $j>i$
(it means that $y_i = q^{2a}$ and $y_j = \nu^2 q^{-2a}$),
then the factors with $r=i$ and $r=j$ are canceled
 in the product in the r.h.s. of (\ref{gfun6}). Thus,
the only factors contribute
  in this product that correspond to adding cells
 to form the diagram ${\sf \Lambda}$ . In this case
 we can substitute the eigenvalues (\ref{qdim25}) of $y_r$
 and consider $M$ as a number of cells in the diagram
 ${\sf \Lambda} \vdash M$. After a cancelation of many factors
 in the r.h.s. of (\ref{gfun6})
 (see the derivation of (\ref{qdim11})) we write (\ref{gfun6})
 in the form
 \be
\lb{gfun7}
\begin{array}{c}
{\displaystyle
  \sum_{j=1}^{k+1}
 \Bigl(\frac{\lambda  q^{2(\lambda_{(j)}-n_{j-1})}}{
 t - q^{2(\lambda_{(j)}-n_{j-1})}} \Bigr)
  {\rm Tr}_{_{(M+1)}} E_{T_{\{{\sf \Lambda}_j',M+1\}}} +
 \sum_{j=1}^k
 \Bigl(\frac{\lambda  \nu^2 q^{2(n_{j} - \lambda_{(j)})}}{
 t - \nu^2 q^{2(n_{j} - \lambda_{(j)})}} \Bigr)
 {\rm Tr}_{_{(M+1)}}  E_{T_{\{{\sf \Lambda}_j'',M+1\}}}
 = }\\ [0.4cm]
= {\displaystyle \Bigl(
\frac{(t-\nu^2 q^{2n})(t-q^{-1} \nu)(t+q \nu)}{
(t-q^{-2n})(t-\nu)(t+\nu)}
\cdot \prod_{r=1}^{k}
\frac{(t - q^{2(\lambda_{(r)}-n_r)})
(t  - \nu^2 \, q^{2(n_{r-1}-\lambda_{(r)})})}{
( t - q^{2(\lambda_{(r)}-n_{r-1})} )
(t-\nu^2 \, q^{2(n_{r}-\lambda_{(r)})})}
- c(t) \Bigr) E_{T_{\{{\sf \Lambda},M\}}}  } \; ,
\end{array}
\ee
where $k$ is a number of blocks in the
diagram ${\sf \Lambda}$ (\ref{lambjj}) and
$n_0=0$.  In the l.h.s. of (\ref{gfun7}) we take into account
 that $y_{M+1} = q^{2(\lambda_{(j)}-n_{j-1})}$,
if we add a new cell in the outer corner
$(n_{j-1}+1, \lambda_{(j)}+1)|_{j=1,...,k+1}$ of ${\sf \Lambda}$,
and
$y_{M+1} = \nu^2 q^{2(n_{j} - \lambda_{(j)})}$, if  we
remove the cell in the inner corner
$(n_{j}, \lambda_{(j)})|_{j=1,...,k}$ of ${\sf \Lambda}$.
Now we compare the residues at
$t=q^{2(\lambda_{(j)}-n_{j-1})} = :\mu_j$ and
$t= \nu^2 q^{2(n_{j} - \lambda_{(j)})} = : \nu^2 \bar{\mu}_j$
in both sides of eq. (\ref{gfun7}) and deduce
 \be
\lb{gfun8}
Tr_{\!_{D(M+1)}} \bigl( E_{T_{\{{\sf \Lambda}_j',M+1\}}} \bigr)
= E_{T_{\{{\sf \Lambda},M\}}} \,
\frac{1-(\mu_j\bar{\mu}_j)^{-1}}{\lambda}
f(\mu_j,q,\nu)
 \prod\limits_{\stackrel{r\neq j}{r=1}}^k
 \frac{\mu_j-\bar{\mu}_r^{-1}}{\mu_j-\mu_r}
 \;\; \prod\limits_{r=1}^k
 \frac{\mu_j- \nu^2 \mu_r^{-1}}{\mu_j- \nu^2 \bar{\mu}_r} \; ,
\ee
 \be
\lb{gfun9}
\!\!
Tr_{\!_{D(M+1)}} \bigl( E_{T_{\{{\sf \Lambda}_j'',M+1\}}} \bigr)
= E_{T_{\{{\sf \Lambda},M\}}} \,
\frac{1-(\mu_j\bar{\mu}_j)^{-1}}{\lambda}
f(\nu^2 \bar{\mu}_j,q,\nu)
\prod\limits_{r=1}^k
 \frac{\nu^2 \bar{\mu}_j- \bar{\mu}_r^{-1}}{
 \nu^2 \bar{\mu}_j- \mu_r}
 \prod\limits_{\stackrel{r\neq j}{r=1}}^k
 \frac{\bar{\mu}_j-\mu_r^{-1}}{\bar{\mu}_j -\bar{\mu}_r} \; ,
\ee
where $f(t,q,\nu):=\frac{(t-\nu^2 q^{2n})(t-q^{-1} \nu)
(t+q \nu)}{(t-q^{-2n})(t-\nu)(t+\nu)}$.
We apply the Ocneanu's (Markov) trace
${\rm Tr}_{_{(1)}} \cdots{\rm Tr}_{_{(M)}}$
to both sides of eq. (\ref{gfun8}) and find the
recurrence relation:
 \be
\lb{gfun10}
{\rm qdim}({\sf \Lambda}_j') =
{\rm qdim}({\sf \Lambda}) \,
\frac{(1-(\mu_j\bar{\mu}_j)^{-1})}{\lambda}
f(\mu_j,q,\nu)
 \prod\limits_{\stackrel{r\neq j}{r=1}}^k
 \frac{\mu_j-\bar{\mu}_r^{-1}}{\mu_j-\mu_r}
 \;\; \prod\limits_{r=1}^k
 \frac{\mu_j- \nu^2 \mu_r^{-1}}{\mu_j- \nu^2 \bar{\mu}_r} \; ,
 \ee
where the diagram ${\sf \Lambda}_j'$ is obtained
by adding one cell
in the outer corner $(n_{j-1}+1, \lambda_{(j)}+1)$
of the diagram ${\sf \Lambda}$. Note that
applying the Ocneanu's (Markov) trace
to both sides of the second eq. (\ref{gfun9}) we deduce the
recurrence relation that is equivalent to the
relation (\ref{gfun10}).
It was
shown in \cite{BB} that the solution of the
recurrence relation (\ref{gfun10})
is given (up to some factor) by the Wenzl's formula
\cite{W2}, \cite{BB}:
 \be
\lb{gfun11}
{\rm qdim}({\sf \Lambda}) =
\prod_{(i,j)\in {\sf \Lambda}}
\frac{q^{\frac{1}{2}d_{\sf \Lambda}(i,j)} -
\nu q^{-\frac{1}{2}d_{\sf \Lambda}(i,j)}}{
q^{\frac{1}{2}h_{i,j}}-q^{-\frac{1}{2}h_{i,j}}}
\prod_{(i,j)\in {\sf \Lambda}}
\frac{\nu^{-1} q^{\frac{1}{2}d'_{\sf \Lambda}(i,j)} +
q^{-\frac{1}{2}d'_{\sf \Lambda}(i,j)}}{
q^{\frac{1}{2}h_{i,j}}+q^{-\frac{1}{2}h_{i,j}}} \;\; ,
\ee
where $h_{i,j}=(\lambda_i + \lambda_j^{\vee}-i-j+1)$
is a hook length (here ${\sf \Lambda}^{\vee} =[\lambda_1^{\vee},\lambda_2^{\vee},...]$ is the transpose
partition of the partition ${\sf \Lambda}$) and
$$
d_{\sf \Lambda}(i,j) =
\left\{
\begin{array}{l}
f_{i,j} \;\; {\rm if} \;\; i \leq j \\
f^{\vee}_{i,j}\;\; {\rm if} \;\; i > j
\end{array}
\right. \; , \;\;\;\;\;
d_{\sf \Lambda}'(i,j) =
\left\{
\begin{array}{l}
f_{i,j} \;\; {\rm if} \;\; i < j \\
f^{\vee}_{i,j}\;\; {\rm if} \;\; i \geq j
\end{array}
\right.  \; ,
$$
where $f_{i,j} =\lambda_i + \lambda_j-i-j+1$
and $f^{\vee}_{i,j} = -\lambda_i^{\vee} - \lambda_j^{\vee}+i+j-1$.

\vspace{0.1cm}

\noindent
{\bf Remark 1.} In the $R$-matrix representation
of the the $BMW_M(\nu)$ algebras, the generators $\sigma_i$
are given by the $SO_q(N)$, $Sp_q(2n)$
$R$-matrices (\ref{frtosph})
(by the $R$-matrices (\ref{Rosp}) in the $Osp_q(N,2m)$
case). For these representations
the parameter $\nu$ is fixed (see (\ref{3.7.1}) and
 Remark {\bf 2} in Subsection {\bf \ref{abmw0}})
and we have $\nu = \epsilon q^{\epsilon-N}$,
where $\epsilon = +1$ and $\epsilon = -1$ correspond to
$SO_q$ and $Sp_q$ cases, respectively. In particular,
 the formula (\ref{gfun11}) is written, for
the $SO_q(N)$ $R$-matrix representation
of $BMW_M$,
in more explicit form
(cf. (\ref{best01}))
\be
\lb{best02}
\begin{array}{c}
{\rm qdim}({\sf \Lambda}) =
 \prod\limits_{i=1}^k  \frac{[N+2(i-1)]_q!}{
[\lambda^{\vee}_i+k-i]_q![N-\lambda^{\vee}_i+k-2+i]_q!} \cdot
 \prod\limits_{i<j} [\lambda^{\vee}_i-\lambda^{\vee}_j + j-i]_{_q}
[N-\lambda^{\vee}_i-\lambda^{\vee}_j + i+j-2]_{_q} \; ,
\end{array}
\ee
where ${\sf \Lambda}^{\vee} = [\lambda^{\vee}_1,\lambda^{\vee}_2,...,\lambda^{\vee}_k]
\vdash M$
is the  transpose partition of ${\sf \Lambda}$
and $[h]_q := \frac{q^h - q^{-h}}{q-q^{-1}}$.

\vspace{0.1cm}

\noindent
{\bf Remark 2.} The analogs of the statements (\ref{shape})
and (\ref{qxinv3}) for the Hecke algebras
 are fairly easy to reformulate and prove
for the case of the $BMW_M$ algebras.

\section{APPLICATIONS AND CONCLUSIONS\label{appli}}
\setcounter{equation}0
\setcounter{subsection}0

In the previous sections of the review, we have presented the fundamentals of the
theory of quantum groups. We have also considered how to
obtain trigonometric and rational (Yangian) solutions of the Yang-Baxter
equation on the basis of the theory of quantum Lie groups. Unfortunately,
in the previous sections it was not possible for us to discuss in detail
the numerous applications of the theory of quantum groups and the Yang-Baxter
equation in both theoretical and mathematical physics. In this final section,
we shall merely give a brief list of such applications that, in the author's
opinion, have some interest.

Before we do this, we recall that in the physics of condensed matter,
exactly solvable two-dimensional models are used to describe various
layered structures, contact surfaces in electronics, surfaces of
superconducting liquids like He II, etc. Two-dimensional integrable field
theories are used to describe dynamical effects in one-dimensional
spatial systems (such as light tubes, nerve fibers, etc.). In addition,
such field theories (and also integrable systems on one-dimensional chains)
can also arise on reductions of multidimensional field theories
(see, for example, Ref. \cite{44}). Quite recently it has been argued
that the one-loop dilatation operator
(anomalous dimension operator) of the ${\cal N}=4$ Super
Yang-Mills theory may be identified, in some restricted cases, with the hamiltonians of various
integrable quantum (super) spin chains \cite{44a}.
Similar spin chain models (related to the
noncompact Lie groups) have previously appeared in the QCD context \cite{44b}, \cite{44bbb}.

\subsection{\bf \em Quantum periodic spin chains\label{perspc}}
\setcounter{equation}0

 We have already mentioned that the quantum inverse scattering method \cite{1}
(an introduction to this method that
can be readily understood by a wide range
of readers can be found in Refs. \cite{45a}, \cite{45}) is designed as a constructive procedure
for solving quantum two-dimensional integrable systems. In addition, the quantum
inverse scattering method makes it possible to construct quantum integrable systems
on one-dimensional chains (see, for example,
Refs. \cite{28}, \cite{36}, and \cite{46}).
Here we discuss the case of periodic chains.
The generalization to the case of
open chains will be mentioned in the next subsection 5.2.
The initial point
is the relation (\ref{3.5.10}) for the $L$ operators, which can be written
in the form\footnote{Usually, this equation is written in the form in which the matrix
$R_{12}(\theta)$ is substituted by $R_{21}(\theta)$. This
is not important, since $R_{21}(\theta)=R_{12}^{-1}(-\theta)$
 satisfies, up to the change of
spectral parameters, the same Yang-Baxter equation
(\ref{3.5.9a}) as $R_{12}(\theta)$ and all formulas below can be easily adapted to the standard case.}
\be
\lb{4.1}
R_{12}(\theta - \theta') \, L_{K2}(\theta) \, L_{K1}(\theta') =
L_{K1}(\theta') \, L_{K2}(\theta) \, R_{12}(\theta -\theta') \; .
\ee
Here $L_{Ki}(\theta)$ are the $N \times N$ matrices
in the auxiliary vector space
$V_i$, with the matrix coefficients that
are the operators in the space of states of the $K$-th site of a chain consisting of $M$ sites:
\be
\lb{4.2}
L_{Ki}(\theta) = I^{\hat{\otimes} (K-1)} \hat{\otimes} L_{i}(\theta) \hat{\otimes}
I^{\hat{\otimes} (M - K)} \; \rightarrow \; [L_{Ki} , \, L_{K'i'}] = 0 \;\;
(K \neq K') \; .
\ee
In (\ref{4.2}), the symbol $\hat{\otimes}$ denotes
a direct product of operator spaces.

To construct an integrable system, we
introduce the monodromy matrix
\be
\lb{4.4}
T_{i}(\theta) = D^{(1)}_{i} \, L_{1i}(\theta) \,
D^{(2)}_{i} \, L_{2i} \dots
D^{(M)}_{i} \, L_{Mi}(\theta)  \; .
\ee
If the matrices $D^{(K)}$ $(1 \leq K \leq M)$ satisfy the relations
\be
\lb{4.5a}
R_{ij} (\theta) \, D^{(K)}_{j} \, D^{(K)}_{i} =
D^{(K)}_{j} \, D^{(K)}_{i} \, R_{ij} (\theta) \; ,
\ee
$$
[ D^{(K)}_i , \, D^{(J)}_j ] = [ D^{(J)}_i , \, L_{Kj} ] = 0 \; ,
$$
then it follows from (\ref{4.1}) that
\be
\lb{4.5}
\begin{array}{c}
R_{ij}(\theta -\theta') \, T_{j}(\theta) \, T_{i}(\theta') =
T_{i}(\theta') \, T_{j}(\theta) R_{ij}(\theta -\theta') \; .
\end{array}
\ee

The trace of the monodromy matrix (\ref{4.4}) over the auxiliary space $i$ forms
the transfer matrix
\be
\lb{tranm1}
t(\theta) = Tr_{(i)} \left(  T_{i}(\theta) \right)
\ee
which gives a commuting family of operators:
$ [ t(\theta), \, t(\theta') ] = 0$.
The commutativity of the transfer matrices
follows directly from eq. (\ref{4.5}) if we multiply it by the
matrix $(R_{ij}(\theta -\theta'))^{-1}$ from the right and take the
trace $Tr_{(i,j)} \left(  \dots \right)$. Using the
family of commuting operators $t(\theta)$
a certain local operator $H$ can be constructed, which is interpreted as the
Hamiltonian of the system. The locality
of the Hamiltonian is a natural physical
requirement and means that $H$ describes the interaction of only
nearest-neighbor sites of the chain. The remaining operators
in the commuting set $t(\theta)$ give an infinite set of integrals of
motion indicating the integrability of the constructed system.
In many well-known cases, the commuting set is associated with the
coefficients in the expansions of $t(\theta)$
over the spectral parameter $\theta$.
For example, one can consider logarithmic derivatives of the
transfer matrix:
\be
\lb{hami1}
{\cal I}_n = \frac{d^n}{d \theta^n } \,
\ln \left( t(\theta) \, t(0)^{-1} \right) \, |_{\theta = 0} \; .
\ee
and identify the local Hamiltonian with the first
logarithmic derivative of the transfer matrix
\be
\lb{hami}
H \equiv {\cal I}_1 = \frac{d}{d \theta } \,
\ln \left( t(\theta) \, t(0)^{-1} \right) \, |_{\theta = 0} \; ,
\ee
where the matrix $t(0)^{-1}$ is introduced in order to obtain the
local charges ${\cal I}_n$ \cite{Lush}.

Now we consider explicit examples of integrable periodic spin chains.
It is clear that from the Yang-Baxter equation (\ref{3.5.9a}) there always follow
representations for the $L$ operators (\ref{4.1}) in the form of $R$-matrices:
\be
\lb{4.3}
\rho_{V_k} \left( L_{Ki}(\theta) \right) =  R_{ki}(\theta) \; , \;\;
\overline{\rho}_{V_k} \left( L_{Ki}(\theta) \right) = (R_{ik}(\theta))^{-1}  \; .
\ee
In this case, the representations of
$L_{Ki}(\theta)$ act nontrivially in the space
$V_{k} \otimes V_{i}$.
We choose for $L$ operators
the first representation in (\ref{4.3}) and obtain for $T_i(\theta)$
(\ref{4.4}):
$$
T_i(\theta) = D^{(1)}_{i} \, R_{1i}(\theta) \,
D^{(2)}_{i} \, R_{2i}(\theta) \dots
D^{(M)}_{i} \, R_{Mi}(\theta) =
$$
$$
=  \R'_{i1}(\theta) \,  \R'_{12}(\theta) \, \R'_{23}(\theta) \dots
\R'_{M-1 M}(\theta) P_{M-1 \, M} \dots P_{23}P_{12} P_{1i} \; ,
$$
where $P_{i j}$ are the permutation matrices and
$\R'_{ij}(\theta) = D^{(j)}_{i} \, \R_{ij}(\theta)$.
Taking the trace $Tr_{(i)}$, we deduce
\be
\lb{tranm}
t(\theta) = Tr_{(i)} \left( \R'_{i1}(\theta) \, \R'_{12}(\theta) \,
\R'_{23}(\theta) \dots \R'_{M-1 M} (\theta) P_{Mi} \right)
P_{M-1 M} \dots P_{23}P_{12} \; .
\ee
We consider a rather general case of $R$-matrices
which can be normalized so that
(see, e.g., (\ref{3.5.10a}), (\ref{3.5.11}), (\ref{3.9.15b}), (\ref{3.9.19}))
\be
\lb{hami22}
\R_{ij}(\theta) =
I + \theta \, h_{ij} + \theta^2 \dots \; .
\ee
These $R$-matrices are called {\it regular} \cite{32}. For the regular $R$-matrices, using
(\ref{tranm}) we obtain
$$
t(\theta) \, t(0)^{-1} = I + \theta
\left(\sum_{k=1}^{M} \, h'_{k \, k+1} \right)
+ \theta^2 \dots \;\; , \;\;\;\;
h'_{k \, k+1} := D^{(k+1)}_k \, h_{k \, k+1}  \,
\left( D^{(k+1)}_k \right)^{-1} \; ,
$$
where $D^{(M+1)}_M := D^{(1)}_M$,
$h_{M \, M+1} := h_{M \, 1}$
and the local Hamiltonian (\ref{hami}) is
\be
\lb{hami2}
H = \sum_{k=1}^{M} \, h'_{k \, k+1}  \; .
\ee
If we choose the $R$-matrix in (\ref{tranm})
in the form of the trigonometric solution (\ref{3.9.15b}), then we obtain
\be
\lb{hami3}
h_{j \, j+1} = {1 \over 2} \left( \R_{j \, j+1} + \R^{-1}_{j \, j+1} -
\lambda \, \beta_{\pm} \, \K_{j \, j+1} \right) \; ,
\ee
where $\beta_{\pm} = \frac{\alpha_{\pm} -1}{\alpha_{\pm} +1}$,
$\alpha_{\pm}= \pm q^{\pm 1}\, \nu^{-1}$ and the parameter $\nu$ is
fixed for different quantum (super)groups in (\ref{nu}).
We note that the Hamiltonians (\ref{hami2}) with the densities (\ref{hami3})
(and $D^{(k)} = 1$) are the $R$-matrix images of the operators:
$$
H_{\pm} = {1 \over 2} \, \sum_{j=1}^{M} \,
\left( \sigma_{j} + \sigma^{-1}_{j} +
\lambda \, \frac{\nu \mp q^{\pm 1}}{\nu \pm q^{\pm 1}} \,
\kappa_{j} \right)  \; ,
$$
where $\sigma_i$, $\kappa_i$ $(i=1, \dots , M)$
obey (\ref{braidg}), (\ref{bmw1}) -- (\ref{bmw3})
with periodic identifications: $\sigma_{M+i} = \sigma_i$,
$\kappa_{M+i} = \kappa_i$. It is natural to call the algebra with
such generators as the periodic Birman-Murakami-Wenzl algebra. The case $\kappa_i = 0$
corresponds to the periodic system with the Hamiltonian:
\be
\lb{hami4}
H =  \sum_{j=1}^{M} \, \sigma_{j} -
{\lambda M \over 2}  \; ,
\ee
where
$\sigma_i$ are the generators of the periodic $A$-type Hecke algebra
$AH_{M+1}$ (see Sec. 4.2).
In the $R$-matrix representation: $\sigma_i \rightarrow \R_i$,
$\sigma_M \rightarrow \R_{M \, 1}$, where $\R$ is the $GL_q(2)$
matrix (\ref{3.3.6a}),
this Hamiltonian describes the periodic $XXZ$ Heisenberg model.

For the Yangian $R$-matrices (\ref{3.9.19})
we obtain $SO(N)$ $(\epsilon = +1)$ and
$Sp(N)$ $(N=2n, \; \epsilon = -1)$ invariant spin chain models
with local Hamiltonian densities (see, e.g. \cite{28}):
$$
h_{l\, l+1} =
\left( P_{l\, l+1} + \frac{2}{2\epsilon-N} K^{(0)}_{l\, l+1} \right)
$$
where, as usual, $P_{l\, l+1}$ are the
transposition matrices, the matrices $K^{(0)}_{l\, l+1}$
were defined in (\ref{3.7.6a}), and
for closed chains we imply $O_{M\, M+1} = O_{M\, 1}$.
The $Osp(N|2m)$ invariant spin
chain model corresponds to the densities
$h_{l,l+1} =
\left( {\cal P}_{l\, l+1} +
\frac{2}{2 + 2 m - N} \, {\cal K}^{(0)}_{l\, l+1} \right)$
which are deduced from (\ref{3.9.19s}).
These Yangian models are generalizations of the XXX Heisenberg models of magnets.
We recall that the XXX model
can be obtained if we take the special limit $q \rightarrow 1$
in the XXZ model or choose the
$gl(2)$ Yangian $R$-matrix (\ref{3.5.11}) as
a representation of $L$ operators in (\ref{4.3}).

By using (in formulas (\ref{4.3}) and (\ref{tranm}))
 the elliptic solution (\ref{3.10.3}),
(\ref{3.10.7})
of the Yang-Baxter equation, we recover for $N=2$
the XYZ spin chain model \cite{Baxt1} while for $N >2$
we obtain its integrable generalizations.

At the end of this subsection we stress that
using the transfer matrix (\ref{tranm1}), one can construct an
integrable 2-dimensional statistical model on the $(M \times L)$ lattice
with periodic boundary conditions. Namely, one should
consider the partition function
$$
Z = Tr_{(1 \dots M)} ( \; \underbrace{t(\theta_0) \cdots t(\theta_0)}_L \; ) =
Tr_{(1 \dots M)} \left( \; \prod_{i=1}^L Tr_{(i)} ( D_i^{(1)} L_{1i}(\theta_0)
\dots D_i^{(M)} L_{Mi}(\theta_0) ) \; \right) \; ,
$$
where the combination $D_i^{(K)} L_{Ki}(\theta_0)$ (for a special value
of the spectral parameter
$\theta = \theta_0$) defines the weight of the statistical
system in the site
$(K,i)$ and $Tr_{(1 \dots M)}$ are the traces over
the operator spaces.

\subsection{\bf \em Factorizable scattering: $S$-matrix and
boundary $K$ matrix\label{factsc}}
\setcounter{equation}0

 The Yang-Baxter equation (\ref{3.5.9}):
\be
\lb{4.6}
S_{23}(\theta -\theta') \, S_{13}(\theta) \, S_{12}(\theta') =
S_{12}(\theta') \, S_{13}(\theta) \, S_{23}(\theta -\theta')
\ee
together with the subsidiary relations of unitarity and crossing symmetry
\be
\lb{4.7}
S_{12}(\theta) \, S_{21}(-\theta) = I_{12} \; , \;\;
S_{12}(\theta) = \left( S_{21}(i \pi -\theta) \right)^{t_{1}} \; ,
\ee
uniquely determine factorizable $S$ matrices (with a minimal set of poles) describing
the scattering of particle-like excitations in $(1 + 1)$-dimensional
integrable relativistic models \cite{3}.
Equations (\ref{4.7}) guarantee that the $S$ matrix $S_{12}(\theta)$
is invertible and skew-invertible (see (\ref{skewa})). The matrix
$S^{i_{1}i_{2}}_{j_{1}j_{2}}(\theta)$ is interpreted
as the $S$ matrix for the scattering of two neutral particles
with isotopic spins $i_{1}$ and $i_{2}$ to two particles
with spins $j_{1}$ and $j_{2}$, and the spectral parameter $\theta$ is none other than
the difference of the rapidities of these particles.
For charged particles, the crossing symmetry relation (\ref{4.7})
should be written in the form
$S_{12}(\theta) = ( S_{2\overline{1}}(i \pi -\theta))^{t_{1}}$,
where the $S$-matrix $S_{2\overline{1}}= S^{i_2\bar{i}_1}_{k_2\bar{k}_1}$
describes particle-antiparticle scattering.
The many-particle
$S$ matrices decompose into products of two-particle matrices (factorization).
In this sense, the Yang-Baxter equation (\ref{4.6})
is the condition for the uniqueness
of the determination of the many-particle $S$ matrices.

The reflection equation
\cite{47}
-- \cite{48K},
which depends on the spectral parameters,
\be
\lb{4.8}
S_{12}(\theta - \theta') \, K_{1}(\theta) \,
S_{21}(\theta + \theta') \, K_{2}(\theta') =
K_{2}(\theta') \, S_{12}(\theta + \theta') \,
K_{1}(\theta) \, S_{21}(\theta - \theta')   \;
\ee
determines, together with the unitarity condition
\be
\lb{unitK}
K^i_j(\theta) \, K^j_m(-\theta) = \delta^i_m \; ,
\ee
and relations (\ref{4.6}) and (\ref{4.7}), the factorizable
scattering of particles (solitons) on a half-line
(see, e.g., \cite{47}, \cite{48'}, \cite{48McK}, \cite{48M}). In this case, the
operator matrix $K^i_{j}(\theta)$ describes
reflection of a particle with rapidity $\theta$ at a boundary point of the half-line.
Graphically, relation (\ref{4.8}) can be represented in the form \\
\unitlength=9mm
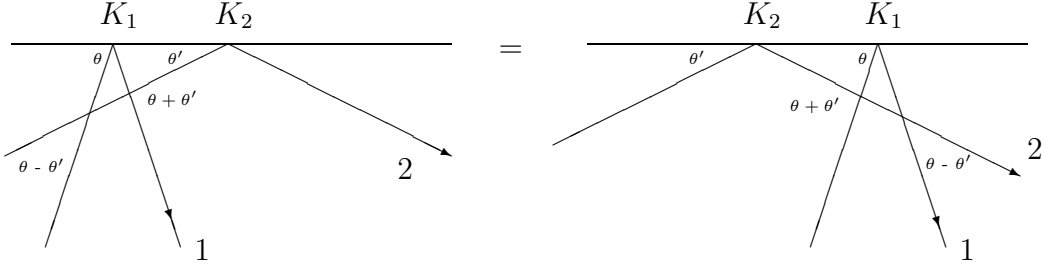
\begin{picture}(16,4.5)
\put(-0.5,3.7){\line(1,0){6.5}}
\put(1,3.7){\line(1,-3){1}}
\put(1,3.7){\line(-1,-3){1}}
\put(1.76,1.4){\vector(1,-3){0.1}}
\put(0.8,4){$K_{1}$}
\put(2.2,0.5){$1$}
\put(1.8,3.4){\tiny $\theta'$}
\put(0.7,3.4){\tiny $\theta$}
\put(-0.4,1.8){\tiny $\theta$ - $\theta'$}
\put(1.5,2.8){\tiny $\theta+\theta'$}
\put(2.7,3.7){\vector(2,-1){3.3}}
\put(2.7,3.7){\line(-2,-1){3.3}}
\put(2.5,4){$K_{2}$}
\put(5.2,1.7){$2$}
\put(6.7,3.5){$=$}

\put(8,3.7){\line(1,0){6.5}}
\put(12.3,3.7){\line(1,-3){1}}
\put(12.3,3.7){\line(-1,-3){1}}
\put(13.1,1.3){\vector(1,-3){0.1}}
\put(12.1,4){$K_{1}$}
\put(13.5,0.5){$1$}
\put(9.5,3.4){\tiny $\theta'$}
\put(12,3.4){\tiny $\theta$}
\put(13,1.85){\tiny $\theta$ - $\theta'$}
\put(11,2.7){\tiny $\theta+\theta'$}
\put(10.5,3.7){\vector(2,-1){3.9}}
\put(10.5,3.7){\line(-2,-1){3}}
\put(10.3,4){$K_{2}$}
\put(14.5,2){$2$}
\end{picture}

We recall \cite{48'} that factorizable scattering
on a half-line can be described by a
Zamolodchikov algebra with generators $\{ A^i(\theta) \}$ $(i=1,\dots,N)$ and boundary operator $B$
that satisfy the defining relations
\be
\lb{zamal}
\begin{array}{c}
A^i(\theta) \, A^j(\theta') = S^{ij}_{kl}(\theta - \theta') \,
A^l(\theta') \, A^k(\theta) \; , \;\;\;
A^i(\theta) \, B = K^i_j(\theta) \, A^j(-\theta) \, B \;\; \Rightarrow \\ \\
A_{1 \rangle}(\theta) \, A_{2\rangle}(\theta') = S_{12}(\theta - \theta') \, A_{2 \rangle}(\theta')
\, A_{1 \rangle}(\theta) \; ,
\;\;\; A_{1 \rangle}(\theta) \, B = K_1(\theta) \, A_{1 \rangle}(-\theta) \, B
\; .
\end{array}
\ee
The consistency conditions for this algebra give rise to the Yang-Baxter equation (\ref{4.6}),
the unitarity conditions (\ref{4.7}), (\ref{unitK}) and the reflection equation (\ref{4.8}) for the
matrices $S$ and $K$.

The reflection equation (\ref{4.8}) can be used \cite{48a} -- \cite{48K}, \cite{48''} for
construction of quantum group invariant integrable
spin systems (see, e.g., \cite{5}) on the
chains with nonperiodic boundary conditions. Indeed, let $T(\theta)$ be a solution of (\ref{4.5})
(for $R_{ij}(\theta) \equiv S_{ji}(\theta)$) and $K(\theta)$ satisfies (\ref{4.8}). Then the matrix
\be
\lb{sklmm}
{\cal T}(\theta) = T(\theta) \, K(\theta) \, [T(-\theta)]^{-1} \; ,
\ee
is also a solution of
(\ref{4.8}). It can be checked directly, but it also follows from the symmetry transformation
$$
A(\theta) \; \rightarrow \; [T(\theta)]^{-1} \, A(\theta) \; , \;\;\;
B \; \rightarrow \; B \; ,
$$
for the algebra (\ref{zamal}) (relations
(\ref{zamal}) are obviously invariant under this transformation if we simultaneously substitute
$K(\theta) \to T(\theta) K(\theta) [T(-\theta)]^{-1}$).

The matrix ${\cal T}(\theta)$ (\ref{sklmm})
is called {\it the Sklyanin monodromy matrix}. By means of
this matrix one can construct a partition function for the integrable lattice model with nontrivial
boundary conditions defined by the reflection matrix $K(\theta)$. The set of commuting integrals
(including the Hamiltonian of the model) is given by the transfer matrix $t(\theta)$ which is
constructed as a special trace of ${\cal T}(\theta)$:
\be
\lb{transre}
t(\theta) = Tr \left(  {\cal T}(\theta) \, \overline{K}(\theta) \right) =
Tr \left(  T(\theta) \, K(\theta) \, [T(-\theta)]^{-1} \,
\overline{K}(\theta) \, \right) \; ,
\ee
where the matrix $\overline{K}(\theta)$ is any solution
of the conjugated reflection equation
\cite{48a} -- \cite{48K}, \cite{48''}:
\be
\lb{conjre}
S^t_{12}(\theta - \theta') \, \overline{K}^t_{2}(\theta') \,
\Psi^t_{12}(\theta + \theta') \, \overline{K}^t_{1}(\theta) =
\overline{K}^t_{1}(\theta) \, \Psi^t_{21}(\theta + \theta') \,
\overline{K}^t_{2}(\theta')  \, S^t_{21}(\theta - \theta') \; .
\ee
Here we require that $\overline{K}(\theta)$ has commutative entries: $[\overline{K}^i_j(\theta), \,
\overline{K}^m_n(\theta')] =0$,
$$
\,
[\overline{K}^i_j(\theta), \, K^m_n(\theta')] = 0 = [\overline{K}^i_j(\theta), \, T^m_n(\theta')]
\; \Rightarrow \; [\overline{K}^i_j(\theta), \, {\cal T}^m_n(\theta')] = 0 \; .
$$
In (\ref{conjre}), we have used the notation
$S^t_{12} := S_{12}^{t_1 t_2}$, and the matrix
$\Psi_{12}$ is the skew-inverse matrix for
$S_{12}$ (cf. (\ref{skewa})):
\be
\lb{psis}
\Psi^{t_1}_{12}(\theta) \, S^{t_1}_{12}(\theta) = I_{12} =
S^{t_1}_{12}(\theta) \, \Psi^{t_1}_{12}(\theta) \; , \;\;\;
\Psi_{12}(\theta) = (S^{t_1}_{12}(\theta)^{-1})^{t_1} \; .
\ee
We also assume (see, e.g., \cite{48'}) that the matrix $S_{12}(\theta)$ satisfies the
cross-unitarity condition
(cf. (\ref{cross01}),  (\ref{cross02}))
\be
\lb{cross04}
S_{12}^{t_1}(\theta) \left( D_1^{-1} S_{21}(b-\theta) D_1 \right)^{t_1} =
\eta(\theta,b) I_{12} \; ,
\ee
where $\eta(\theta,b)$ is the scalar function, $b$ is the special parameter which depends on the form
of the matrix $S_{12}$ and
$D$ is the constant matrix such that:
$[D_1 D_2, \, S_{12}(\theta)]=0$. Comparing eqs. (\ref{psis}) and (\ref{cross04}),
one can identify
$$
\Psi_{12}(\theta) = \frac{1}{\eta(\theta,b)} \, D_1^{-1} S_{21}(b-\theta) D_1 \; ,
$$
and, then, rewrite the conjugated reflection equation (\ref{conjre}) in the form
\be
\lb{conjre1}
S^t_{12}(\theta - \theta') \, \widetilde{K}^t_{2}(\theta') \,
S^t_{21}(b-\theta -\theta') \, \widetilde{K}^t_{1}(\theta) =
\widetilde{K}^t_{1}(\theta) \, S^t_{12}(b-\theta -\theta') \,
\widetilde{K}^t_{2}(\theta')  \, S^t_{21}(\theta - \theta') \; ,
\ee
where $\widetilde{K}(\theta) = D^{-1} \overline{K}(\theta)$.
Note that Eq. (\ref{conjre1}) is also
one of the consistence conditions but for the "left-boundary" Zamolodchikov algebra (cf.
(\ref{zamal}))
\be
\lb{leftZ}
\widetilde{B} \, \widetilde{A}_{1 \rangle}(\theta) =
\widetilde{B} \, \widetilde{K}_1(\theta) \, \widetilde{A}_{1 \rangle}(b-\theta) \; , \;\;\;
\widetilde{A}_{1 \rangle}(\theta) \, \widetilde{A}_{2 \rangle}(\theta') = S_{21}(\theta - \theta')
\widetilde{A}_{2 \rangle}(\theta') \widetilde{A}_{1 \rangle}(\theta) \; ,
\ee
with generators $\widetilde{A}^i$ $(i=1,\dots,N)$ and
left boundary operator $\widetilde{B}$ (we
need to consider the condition for the unique
reordering of the third order monomial
$\widetilde{B} \, \widetilde{A}_{1 \rangle}(\theta) \, \widetilde{A}_{2 \rangle}(\theta')$).

The proof of the identity $[t(\theta), \, t(\theta')]=0$ for the transfer matrix $t(\theta)$
(\ref{transre}) is straightforward \cite{48a}
 ($\theta^{\pm} = \theta \pm \theta'$)
$$
t(\theta') \, t(\theta) = Tr_{12} \left( {\cal T}_2(\theta') \, {\cal T}^t_1(\theta) \,
\overline{K}_2(\theta') \, \overline{K}^t_1(\theta)   \right) =
$$
$$
= Tr_{12} \left( {\cal T}_2(\theta') \, {\cal T}^t_1(\theta) \,
 S^{t_1}_{12}(\theta^+) \, \Psi^{t_1}_{12}(\theta^+) \,
\overline{K}_2(\theta') \, \overline{K}^t_1(\theta) \right) =
$$
$$
= Tr_{12} \left(
 \left( {\cal T}_2(\theta') \, S_{12}(\theta^+) \,
 {\cal T}_1(\theta) \right)^{t_1} \,
\left( \overline{K}^t_2(\theta') \, \Psi^{t}_{12}(\theta^+) \,
\overline{K}^t_1(\theta)  \right)^{t_2}  \right) =
$$
using eq. (\ref{4.8}) for $K(\theta) \to {\cal T}(\theta)$ we deduce
$$
= Tr_{12} \left( {\cal T}_{1}(\theta) \,
S_{21}(\theta^+) \, {\cal T}_{2}(\theta') \, S^{-1}_{21}(\theta^-)
\left( \overline{K}^t_2(\theta') \, \Psi^{t}_{12}(\theta^+) \,
\overline{K}^t_1(\theta)  \right)^{t} \,  S_{12}(\theta^-)
 \right) =
$$
and applying here the conjugated reflection equation (\ref{conjre}) and transpositions
we finally obtain
$$
= Tr_{12} \left( \left( {\cal T}_{1}(\theta) \,
S_{21}(\theta^+) \, {\cal T}_{2}(\theta') \right)^{t_2} \,
\left( \overline{K}^t_1(\theta) \, \Psi^{t}_{21}(\theta^+) \,
\overline{K}^t_2(\theta')  \right)^{t_1}
 \right) =
$$
$$
= Tr_{12} \left( {\cal T}_{1}(\theta) \,  {\cal T}^t_{2}(\theta') \,
 S^{t_2}_{21}(\theta^+) \, \Psi^{t_2}_{21}(\theta^+) \, \overline{K}_1(\theta) \,
 \overline{K}^t_2(\theta')  \right) = t(\theta) \, t(\theta') \; .
$$

Now we take in (\ref{4.8}) the limit $\theta, \, \theta' \rightarrow \pm \infty$
in such a way that $\theta - \theta' \rightarrow \pm \infty$,
and at the same time we set
$$
K(\theta) |_{\theta \rightarrow \infty} = L \; , \;\;
S_{12}(\theta) |_{\theta \rightarrow \infty} = R_{12} \; .
$$
$$
K(\theta) |_{\theta \rightarrow - \infty} = L^{-1} \; , \;\;
S_{12}(\theta) |_{\theta \rightarrow - \infty} = (R_{21})^{-1} \; .
$$
Then (\ref{4.8}) goes over into (\ref{3.1.23a}), and this is the reason why all
algebras with defining relations of type
(\ref{3.1.23a}) are called the reflection equation algebras \cite{48} --
\cite{48K}.

Note that each solution of the Yang-Baxter equation (\ref{4.6})
with the conditions (\ref{4.7}) determines an equivalence class of
quantum integrable systems with the given factorizable $S$ matrix.
Thus, each classification of solutions to the Yang-Baxter equation is,
to some extent, a classification of integrable systems with the
properties indicated above.

The $3D$ analog of the Yang-Baxter (triangle) equation (\ref{zzz}), (\ref{4.6}) is called the {\em
tetrahedron equation} \cite{4}, \cite{51} and
defines the consistence condition for $3D$ factorizable
scattering of strings.The 3-dimensional model of
such factorizable scattering was first
proposed by A.Zamo\-lodchikov in \cite{4}. Then, this $3D$ model was generalized in
\cite{BaBax}. New solutions of the tetrahedron equation were also considered in
\cite{BazhSer}.  A 3-dimensional version of the
$2D$ reflection equation (\ref{4.8}) (the tetrahedron reflection equation)
was proposed in \cite{49}. Combinatorial and algebraic aspects of the
$3D$ reflection equation were considered in \cite{Kuni},
\cite{Kuni2}. Special solutions of the tetrahedron reflection equation
were found in \cite{Yone}.

From a mathematical point of view,
higher dimensional generalizations of the Yang-Baxter
equations are related to the Manin-Schechtman higher braid groups \cite{MaSch}, $n$-categories
\cite{ncat}, \cite{Baez1}, \cite{Baez2},
 and also appeared in the theory of quasi-triangular
Hopf algebras (see Remark at the end of Subsection 2.6).

\subsection{\bf \em  Yang-Baxter equations and calculations of multiloop Feynman diagrams\label{ybfd}}
\setcounter{equation}0

 We mention the application of the Yang-Baxter equation
in multiloop calculations in quantum field theory. There is a form of
the Yang-Baxter equation (see Refs. \cite{2}, \cite{50}, and \cite{51}) that
can also be represented in the form of the triangle equation
(\ref{zzz}), but the indices $x,x_i,y_i$
are ascribed not to the "lines"
but to the "faces": \\
\unitlength=9mm
\begin{picture}(15,4.3)
\put(4.5,0.1){$y_{1}$}
\put(4.5,1.9){$x_{3}$}
\put(4.5,3.8){$y_{2}$}
\put(6.2,1){$x_{2}$}
\put(5.4,1){\tiny $\theta$}
\put(7.5,1.9){$y_{3}$}
\put(6.7,1.9){\tiny $\theta'$}
\put(5.6,1.9){$x$}
\put(5.2,4){\vector(0,-1){4}}
\put(4.3,3.4){\vector(3,-2){3.2}}
\put(6.2,3){$x_{1}$}
\put(4.3,0.6){\vector(3,2){3.2}}
\put(8.3,1.9){$=$}

\put(11.6,4){\vector(0,-1){4}}
\put(9.2,2.7){\vector(3,-2){3}}
\put(9.2,1.3){\vector(3,2){3}}
\put(11.1,1.9){$x$}
\put(10.5,1.9){\tiny $\theta'$}
\put(9,1.9){$x_{3}$}
\put(12,1.9){$y_{3}$}
\put(11.8,0.1){$x_{2}$}
\put(11.8,3.6){$x_{1}$}
\put(11.8,2.7){\tiny $\theta$}
\put(10.1,1){$y_{1}$}
\put(10.1,3){$y_{2}$}
\put(15.5,1.9){(5.3.0)}
\end{picture}

\noindent
where $\theta,\theta'$ are angles (spectral parameters),
summation is over the index $x$, and

\unitlength=8mm
\begin{picture}(17,2.5)
\put(3.5,0.8){$R^{x y}_{u z}(\theta) \; =$}
\put(7.6,2){$y$}
\put(8.5,1){$z$}
\put(6.8,1){$x$}
\put(7.6,0.3){$u$}
\put(7.7,1.2){\tiny $\theta$}
\put(7,0.2){\vector(1,1){1.5}}
\put(8.5,0.2){\vector(-1,1){1.5}}
\end{picture}

\noindent
The analytical form of (5.3.0) is
\be
\lb{facem}
\sum_x \, R^{x_1 \, x}_{y_2 \, x_3}(\theta - \theta') \,
R^{x \, x_2}_{x_3 \, y_1}(\theta) \,
R^{x_1 y_3}_{x \, x_2}(\theta') =
\sum_x \, R^{y_2 \, x}_{x_3 \, y_1}(\theta')  \, R^{x_1 \, y_3}_{y_2 \, x}(\theta)
 \, R^{y_3 \, x_2}_{x \, y_1} (\theta - \theta') \; .
\ee
We have already considered the solution of this equation in Sec. 3.13.
Indeed, one can show that eq. (\ref{3.10.9}) is equivalent to eq. (\ref{facem}) if we put
(for the notation see Sec. 3.13):
\be
\lb{facem1}
R^{x_1 \, x}_{x_2 \, x_3}(\theta) = \omega^{{1 \over 2}
<x-x_2 \, , \, x_1 +x_3 > + <x_2, x>}
\, W_{x+x_2 - x_1 - x_3}(\theta) \; ,
\ee
where the indices $x,x_i$ are 2-dimensional vectors, e.g. $x = (\alpha_1,\alpha_2) \in {\bf Z}_N^2$.
 Thus, (\ref{facem1}), (\ref{3.10.7}) and
(\ref{3.10.8}) solve the face-type Yang-Baxter equation (\ref{facem}).

There is a transformation from the vertex-type Yang-Baxter equation
(\ref{3.5.9a}) to the face type (\ref{facem}) using intertwining
vectors $\psi_i^{x_1 x_2}$ (see e.g. \cite{48K} and Refs. therein),
where $i$ is a vertex index while $x_1,x_2$ are face indices. The vectors
$\psi_i^{x_1 x_2}$ satisfy the intertwining relations
\be
\lb{inttw}
 \psi_{\langle 2}^{x_2 x_1}(\theta - \theta') \,
\psi_{\langle 1}^{x_3 x_2}(\theta) \, R_{12}(\theta') =
\sum_x  \, R^{x_1 \, x}_{x_2 x_3} (\theta')
\, \psi_{\langle 1}^{x x_1} (\theta) \,
\psi_{\langle 2}^{x_3 x}(\theta - \theta')
\ee
which are represented graphically in the form (here the angles are the same as
in (5.3.0)) \\
\unitlength=6mm
\begin{picture}(15,4.3)
\put(4.3,0.1){$x_{3}$}
\put(4.3,1.9){$x_{2}$}
\put(4.3,3.8){$x_{1}$}
\put(7.3,0.6){$2$}
\put(5.2,4){\vector(0,-1){4}}
\put(5.3,4){\vector(0,-1){4}}
\put(4.3,3.4){\vector(3,-2){3.2}}
\put(7.3,3){$1$}
\put(4.3,0.6){\vector(3,2){3.2}}
\put(8.5,1.9){$=$}

\put(12.5,4){\vector(0,-1){4}}
\put(12.6,4){\vector(0,-1){4}}
\put(10.2,2.7){\vector(3,-2){3}}
\put(10.2,1.3){\vector(3,2){3}}
\put(12,1.9){$x$}
\put(10,1.9){$x_{2}$}
\put(13.3,0.3){$2$}
\put(13.3,3.4){$1$}
\put(11.1,1){$x_{3}$}
\put(11.1,3){$x_{1}$}
\put(15,1.9){,}

\put(16.2,1.8){$\psi_{i}^{x_1 x_2}(\theta)  :=$}
\put(21.3,0.8){$i$}
\put(20.5,1.5){\tiny $\theta$}
\put(20,2.4){$x_1$}
\put(21.6,2.4){$x_2$}
\put(21,2.8){\vector(0,-1){2}}
\put(22,1.9){\vector(-1,0){2}}
\put(22,2){\vector(-1,0){2}}
\put(23,1.9){.}
\end{picture}

\noindent
Then, the face type Yang-Baxter equation (\ref{facem}) is obtained from
the vertex equation (\ref{3.5.9a}) if we act on it by
$(\psi_{\langle 3}^{y_2 x_1} \, \psi_{\langle 2}^{x_3 y_2} \,
\psi_{\langle 1}^{y_1 x_3})$
from the left.

Relations (5.3.0) and (\ref{facem}),
like (\ref{zzz}), give the conditions of integrability
of two-dimensional lattice statistical systems
(interaction round face models) with weights determined by
the $R$-matrices $R^{xy}_{uz}(\theta)$. In this case, the transfer matrix
has the form
$$
t^{y_1 y_2 \dots y_M}_{x_1 x_2 \dots x_M}(\theta) =
R^{y_1 y_2}_{x_1 x_2}(\theta) \, R^{y_2 y_3}_{x_2 x_3}(\theta) \,
R^{y_3 y_4}_{x_3 x_4}(\theta) \dots R^{y_M y_1}_{x_M x_1}(\theta) \; ,
$$
and its graphical representation is

\unitlength=6mm
\begin{picture}(15,4.3)

\put(1,2){\vector(3,0){11}}
\put(1.3,3){$y_{1}$}
\put(2.8,3){$y_2$}
\put(4.8,3){$y_3$}
\put(1.3,1){$x_{1}$}
\put(2.8,1){$x_2$}
\put(4.8,1){$x_3$}
\put(2.3,0.6){\vector(0,1){3.2}}
\put(4.3,0.6){\vector(0,1){3.2}}
\put(6.3,0.6){\vector(0,1){3.2}}
\put(8.5,0.6){\vector(0,1){3.2}}
\put(10.5,0.6){\vector(0,1){3.2}}
\put(6.8,1){$\dots$}
\put(6.8,3){$\dots$}
\put(8.9,1){$x_{M}$}
\put(8.9,3){$y_{M}$}
\put(6.5,1){\oval(13.5,2)}
\put(6,0){\vector(-1,0){2}}

\put(2.5,2.2){\tiny $\theta$}
\put(4.5,2.2){\tiny $\theta$}
\put(6.5,2.2){\tiny $\theta$}
\put(8.8,2.2){\tiny $\theta$}
\put(10.8,2.2){\tiny $\theta$}
\put(14.7,1.9){,}

\end{picture}

\noindent
while the partition function for the periodic system on
the $(M \times K)$ lattice is given
by the standard formula: $Z = Tr_{1 \dots M} (t(\theta))^K$.

We now note that the Yang-Baxter equation (5.3.0), (\ref{facem})
has a solution in the form
$R^{xy}_{uz}(\theta) = G^{y}_{u}(\theta) \, \bar{G}^{x}_{z}(\pi - \theta)$, where the
matrices $G^{y}_{u}$, $\bar{G}^{x}_{z} = \bar{G}^{z}_{x}$ satisfy the star triangle relation
(see, for example, Refs. \cite{2} and \cite{50})
\be
\lb{4.9}
f(\theta,\theta') \bar{G}^{x_{1}}_{x_{3}}(\pi - \theta + \theta') \,
G^{x_{2}}_{x_{3}}(\theta) \,
\bar{G}^{x_{1}}_{x_{2}}(\pi - \theta') =
\sum_{x} G^{x}_{x_3}(\theta') \,
\bar{G}^{x_1}_{x}(\pi - \theta) \,
G^{x_{2}}_{x}(\theta - \theta') \; ,
\ee
and $f(.,.)$ is an arbitrary function such that
$f(\theta,\theta') = f(\theta, \theta - \theta')$.
The relations (\ref{4.9}) for $f=1$ can
be represented graphically in the form \\
\unitlength=9mm
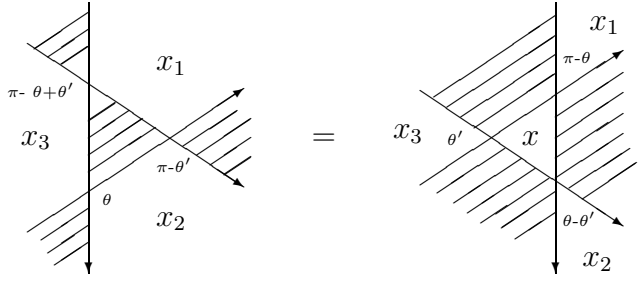
\begin{picture}(15,4.3)
\put(4.2,1.9){$x_{3}$}
\put(6.2,1.5){\tiny $\pi$-$\theta'$}
\put(4,2.6){\tiny $\pi$- $\theta$+$\theta'$}
\put(6.2,0.7){$x_{2}$}
\put(5.4,1){\tiny $\theta$}

\put(5.2,4){\vector(0,-1){4}}
\put(5.2,0.95){\line(-3,-2){0.7}}
\put(5.2,0.7){\line(-3,-2){0.6}}
\put(5.2,0.45){\line(-3,-2){0.5}}

\put(4.3,3.4){\vector(3,-2){3.2}}
\put(5.2,3.35){\line(-3,-2){0.4}}
\put(5.2,3.6){\line(-3,-2){0.6}}
\put(5.2,3.85){\line(-3,-2){0.8}}

\put(6.2,3){$x_{1}$}
\put(4.3,0.6){\vector(3,2){3.2}}
\put(5.2,2.25){\line(3,2){0.4}}
\put(5.2,2){\line(3,2){0.6}}
\put(5.2,1.75){\line(3,2){0.8}}
\put(5.2,1.5){\line(3,2){0.99}}

\put(6.6,1.87){\line(3,2){0.8}}
\put(6.8,1.75){\line(3,2){0.8}}
\put(7,1.6){\line(3,2){0.6}}
\put(7.2,1.45){\line(3,2){0.4}}

\put(8.5,1.9){$=$}

\put(12.1,4){\vector(0,-1){4}}
\put(12,1.4){\line(-3,-2){1}}
\put(11.8,1.55){\line(-3,-2){1}}
\put(11.6,1.7){\line(-3,-2){1}}
\put(11.35,1.85){\line(-3,-2){1}}
\put(12.1,1.2){\line(-3,-2){0.8}}
\put(12.1,0.9){\line(-3,-2){0.6}}

\put(10.1,2.7){\vector(3,-2){3}}
\put(12.5,1.1){\line(3,2){0.8}}
\put(12.3,1.25){\line(3,2){0.9}}
\put(12.1,1.4){\line(3,2){1}}
\put(12.1,1.7){\line(3,2){1}}
\put(12.1,2){\line(3,2){1}}
\put(12.1,2.3){\line(3,2){1}}

\put(10.1,1.3){\vector(3,2){3}}
\put(10.9,2.15){\line(3,2){1.2}}
\put(10.7,2.3){\line(3,2){1.4}}
\put(10.5,2.45){\line(3,2){1.6}}
\put(10.3,2.6){\line(3,2){1.8}}

\put(11.6,1.9){$x$}
\put(9.7,2){$x_{3}$}
\put(10.5,1.9){\tiny $\theta'$}
\put(12.5,0.1){$x_{2}$}
\put(12.2,0.7){\tiny $\theta$-$\theta'$}
\put(12.2,3.1){\tiny $\pi$-$\theta$}
\put(12.6,3.6){$x_{1}$}
\end{picture}

The Feynman diagrams, which will be considered
here, are graphs with vertices
connected by lines labeled by numbers (indices). With each vertex we
associate the point in the D-dimensional space $R^D$ while the lines of the graph
(with index $\alpha$) are associated with the massless Feynman propagator

\unitlength=1cm
\begin{picture}(25,1.5)

\put(1,0.5){\line(1,0){1.8}}
\put(1.7,0.7){$\alpha$}

\put(0.5,0.5){$x$}
\put(3.2,0.5){$x'$}
\put(5,0.5){$=$}
\put(6,0.5){$\frac{\Gamma(\alpha)}{(x-x')^{2\alpha}}$}
\end{picture}
\noindent
(which is a function of two points $x,x'$ in $D$ dimensional space-time):
\be
\lb{4.10}
G_{D}(x-x'|\alpha) =
\frac{\Gamma(\alpha)}{(x-x')^{2\alpha}} =
\frac{\Gamma(\alpha)}{(\sum_\mu (x-x')_{\mu}(x-x')^{\mu})^{\alpha}} \; ,
\ee
where $\Gamma(\alpha)$ is the Euler gamma-function,
$D = 4-2\epsilon$ is the dimension of space-time, $(x)_{\mu}$
$(\mu = 1,2, \dots , D)$ are its
coordinates, $\alpha = D/2 -1+\eta$, and $\epsilon$ and $\eta$ are, respectively, the parameters of
the dimensional and analytic regularizations. The boldface vertices $\bullet$ denote that the
corresponding points $x$ are integrated over $R^D$:
$\frac{1}{\pi^{D/2}} \, \int d^{D} x$.
These diagrams are called the Feynman diagrams in the configuration space.

The propagator (\ref{4.10}) satisfies the relation
\be
\lb{4.11}
\begin{array}{c}
 \int  \frac{d^{D} x}{\pi^{D/2}} \,
\prod\limits_{i=1}^{3} G_{D}(x-x_{i}|\alpha_{i})
\;\; \stackrel{\sum \alpha_{i} = D}{=} 
G_{D}(x_{1}-x_{2}|\alpha'_{3}) \,
G_{D}(x_{2}-x_{3}|\alpha'_{1}) \,
G_{D}(x_{3}-x_{1}|\alpha'_{2}) \; ,
\end{array}
\ee
which is represented as
the star-triangle identity for the Feynman diagrams:

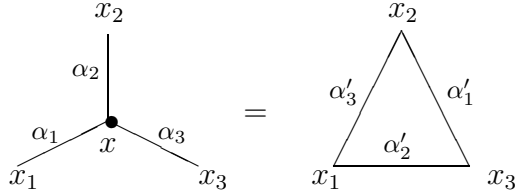
\unitlength=6mm
\begin{picture}(20,4.9)(-6,0)

\put(8,1){\line(1,0){3}}

\put(9.2,4.3){$x_2$}
\put(7.5,0.6){$x_1$}
\put(11.4,0.6){$x_3$}

\put(8,1){\line(1,2){1.5}}
\put(7.9,2.5){\footnotesize $\alpha'_3$}
\put(9.1,1.3){\footnotesize $\alpha'_2$}
\put(10.5,2.5){\footnotesize $\alpha'_1$}

\put(11,1){\line(-1,2){1.5}}

\put(6,2){$=$}

\put(1,1){\line(2,1){2}}
\put(3,2){\line(2,-1){2}}
\put(3,1.9){\line(0,1){2}}
\put(2.9,1.8){$\bullet$}

\put(4.1,1.6){\footnotesize $\alpha_3$}
\put(1.3,1.6){\footnotesize $\alpha_1$}
\put(2.2,3){\footnotesize $\alpha_2$}
\put(2.8,1.3){$x$}
\put(2.7,4.3){$x_2$}
\put(0.8,0.6){$x_1$}
\put(5,0.6){$x_3$}

\end{picture}

\noindent
where  $\alpha_1 + \alpha_2 + \alpha_3 = D$ and
$\alpha_i' := D/2 - \alpha_i$.
Equation (\ref{4.11}) can be readily derived if
we put $(x_{3})^\mu =0$ $\forall \mu$
and make in the right-hand side of
(\ref{4.11}) a simultaneous inversion
transformation of the variables of integration,
$ (x)_{\mu} \rightarrow (x)_{\mu}/x^{2}$
and the coordinates  $(x_{1,2})^\mu$.
Relations (\ref{4.9}) and (\ref{4.11}) are equivalent if we set
\be
\lb{4.12}
G^{x}_{x'} (\theta) = \bar{G}^{x}_{x'} (\theta) = G_{D} ( x-x'|
\frac{D}{2} ( 1 - \frac{\theta}{\pi} ) ) \; , \;\;
f(\theta,\, \theta') =1 \; , \;\; \sum_x = \int \frac{d^{D} x}{\pi^{D/2}}
\; .
\ee
Thus, the analytically and dimensionally regularized massless propagator
(\ref{4.10}) satisfies the infinite-dimensional star-triangle relation (\ref{4.9}) and
accordingly, on the basis of (\ref{4.10}) and (\ref{4.12}),
we can construct solutions of the
Yang-Baxter equation (5.3.0), (\ref{facem}). This remark was made
in Ref. \cite{50}, in which calculations
were carried out of vacuum diagrams with an infinite number
(in the thermodynamical limit) of vertices corresponding to a
planar square lattice ($\phi^{4}$ theory, $D=4$), a planar
triangular lattice ($\phi^{6}$ theory, $D=3$),
and a honeycomb lattice ($\phi^{3}$ theory, $D=6$). The
star-triangle relation (\ref{4.11}) (known
also as {\it the uniqueness relation}) was used in addition for analytic calculation of the
diagrams that contribute to the 5-loop $\beta$- function
of the $\phi^{4}_{D=4}$ theory \cite{52} and of
massless ladder diagrams \cite{52a}, \cite{52''}. By means of identity (\ref{4.11})
the symmetry groups of
dimensionally and analytically regularized
massless diagrams were investigated \cite{52''}, \cite{52'},
\cite{Broad}\footnote{Here the symmetry of diagrams
means the symmetry of the corresponding perturbative integrals.}.
We emphasize that an extremely interesting problem is that of massive
deformation of the propagator function (\ref{4.10}) and the corresponding
deformation of the star-triangle relation (\ref{4.11}).

There is an elegant operator interpretation \cite{52''}, \cite{IsDia} of the star-triangle identity
(\ref{4.11}). Indeed, consider the $D$-dimensional
Heisenberg algebra ${\cal H}_D$
as the algebra of functions of the generators $\hat{q}_\mu =
\hat{q}_\mu^\dagger$ and $\hat{p}_\mu = \hat{p}_\mu^\dagger$
($\mu = 1, \dots , D$) subject to the defining relations
\begin{equation}
\label{gr0}
[\hat{q}_\mu , \, \hat{p}_\nu ] = {\bf i} \, \delta_{\mu \nu}  \; , \;\;\;
(\mu,\nu = 1,2, \dots , D) \; ,
\end{equation}
where $\hat{q}_\mu$ and $\hat{p}_\mu$ are the operators
of the coordinate and momentum, respectively.
Consider a representation of the algebra (\ref{gr0}) in the linear vector space of complex
functions $\psi(x) := \psi(x_\mu)$ on $\mathbb{R}^D$:
$$
\hat{q}_\mu \, \psi(x) = x_\mu \, \psi(x) \; , \;\;\;
\hat{p}_\mu \, \psi(x) = - {\bf i} \, \partial_\mu \, \psi(x) \; .
$$
It is convenient to realize the action of elements $\hat{A} \in {\cal H}_D$ as the action of
integral operators:
$\hat{A} \, \psi(x) = \int d^D y \, \langle x | \hat{A} | y \rangle \, \psi(y)$.
The integral kernels $\langle x | \hat{A} | y \rangle$ can be considered as matrix elements of
$\hat{A}$ for the states $|x \rangle := |\{ x_\mu \}\rangle$ and
$\langle y | = |y \rangle^{\dagger}$
such that
\begin{equation}
\label{gr1}
\langle y |x\rangle  = \delta^D(y-x) \; , \;\;\;
\hat{q}_\mu |x\rangle  = x_\mu \, |x\rangle \; , \;\;\;
\int d^D x \, |x\rangle \, \langle x| = \hat{1}  \; .
\end{equation}
We extend the algebra ${\cal H}_D$ by the elements
$\hat{q}^{2\alpha} := (\hat{q}^\mu \hat{q}_\mu)^\alpha$ and pseudo-differential operators
$\hat{p}^{-2\beta} :=(\hat{p}^\mu \hat{p}_\mu)^{-\beta}$ $(\forall \alpha,\beta \in \mathbb{C})$.
The corresponding integral kernels are:
\begin{equation}
\label{gr4}
\langle x| \hat{q}^{2 \alpha} |y\rangle = x^{2\alpha} \, \delta^D(x-y) \; , \;\;\;
\langle x| \frac{1}{\hat{p}^{2 \beta}} |y\rangle = a(\beta) \,
\frac{1}{(x-y)^{2 \beta'}} \; ,
\end{equation}
where $a(\beta) =
\frac{\Gamma(\beta') }{ \pi^{D/2} \, 2^{2\beta} \, \Gamma(\beta)}$,
$\beta' = D/2 - \beta$
and $\beta' \neq 0, -1, -2, \dots$.

For the extended Heisenberg algebra one can prove
 \cite{52''}, \cite{IsDia} that the operators
$H_{\alpha} := \hat{p}^{2\alpha} \, \hat{q}^{2\alpha}$
$(\forall \alpha \in \mathbb{C})$ form a commutative family.
The commutativity condition
$[H_{\alpha} , \, H_{-\beta} ] = 0$
is represented in the form
\begin{equation}
\label{uniq}
\hat{p}^{2 \alpha} \hat{q}^{2 \gamma} \hat{p}^{2 \beta} =
\hat{q}^{2 \beta} \hat{p}^{2 \gamma} \hat{q}^{2 \alpha} \; , \;\;\;
(\gamma = \alpha + \beta) \; .
\end{equation}
Then it is not hard to see that this identity, written for integral kernels
by means of (\ref{gr4}), is equivalent to the
star-triangle relation (\ref{4.11}). One should act on (\ref{uniq})
by vectors $\langle x_1-x_3|$ and $|x_2-x_3 \rangle$ from the left and right,
respectively, and insert, in the l.h.s. of (\ref{uniq}),
the unit $\hat{1}$ (\ref{gr1}):
$$
\begin{array}{c}
\langle x_1-x_3|\hat{p}^{2 \alpha} \left( \int d^D x \, |x\rangle \, \langle x| \right)
\hat{q}^{2 \gamma} \hat{p}^{2 \beta}|x_2-x_3 \rangle =
\langle x_1-x_3| \hat{q}^{2 \beta} \hat{p}^{2 \gamma} \hat{q}^{2 \alpha}|x_2-x_3 \rangle \;\;
\Rightarrow  \\ [0.3cm]
\int d^D x \,
\langle x_1-x_3|\hat{p}^{2 \alpha} |x\rangle \,
x^{2 \gamma} \, \langle x| \hat{p}^{2 \beta}|x_2-x_3 \rangle =
(x_1-x_3)^{2 \beta} \langle x_1-x_3| \hat{p}^{2 \gamma}|x_2-x_3 \rangle
(x_2-x_3)^{2 \alpha} \; .
 \end{array}
$$
Applying here the second eq. in (\ref{gr4}), we obtain (\ref{4.11})
for $\alpha = -\alpha_1'$,
$\beta = -\alpha_2'$ and $\gamma = -\alpha_3$.

Consider the set of Heisenberg algebras ${\cal H}_D$ with
the generators $\{ \hat{q}_{(a)}^\mu ,
\hat{p}_{(b)}^\nu\}$  $(a,b=1,2,\dots,N)$ such that:
$[\hat{q}_{(a)}^\mu ,  \, \hat{p}_{(b)}^\nu ] =
{\rm i} \, \delta^{\mu \nu} \, \delta_{ab}$.
Then, the star-triangle identity (\ref{uniq}) is obviously generalized as
\be
\lb{lip1}
(\hat{q}_{(ab)})^{2\alpha} (\hat{p}_{(b)})^{2(\alpha +\beta)} (\hat{q}_{(ab)})^{2\beta} =
(\hat{p}_{(b)})^{2\beta} (\hat{q}_{(ab)})^{2(\alpha+\beta)} (\hat{p}_{(b)})^{2\alpha} \; ,
\ee
where $\hat{q}_{(ab)}^\mu = \hat{q}_{(a)}^\mu-\hat{q}_{(b)}^\mu$.
 Taking into account (\ref{lip1}),
one can directly check that for an
arbitrary parameter $\xi$ the operator
\be
\lb{Rlip}
\begin{array}{c}
R_{ab}(\alpha;\xi):=
(\hat{q}_{(ab)})^{2(\alpha+\xi)} (\hat{p}_{(a)})^{2\alpha} (\hat{p}_{(b)})^{2\alpha}
(\hat{q}_{(ab)})^{2(\alpha-\xi)} =  1 + \alpha \, h_{(ab)}(\xi) + \alpha^2 \dots \; ,
\end{array}
\ee
 is a regular (see (\ref{hami22})) solution of the Yang-Baxter equation:
\be
\lb{lip2}
R_{ab}(\alpha;\xi) \, R_{bc}(\alpha + \beta;\xi) \, R_{ab}(\beta;\xi) = R_{bc}(\beta;\xi)
\, R_{ab}(\alpha + \beta;\xi) \, R_{bc}(\alpha;\xi) \; .
\ee
 The solution (\ref{Rlip}) for arbitrary $D$ and $\xi=1$
 was found in \cite{IsDia} and for any $\xi$ in \cite{DICh}.
The factorized form of the solution (\ref{Rlip}) (for $D=1$)
reminds the factorization of R-matrices observed
in \cite{Derk1}, \cite{Derk}.

 Using the standard procedure (see eqs. (\ref{hami22}), (\ref{hami2})),
 one can construct an
integrable system with a Hamiltonian that is related to the $R$-matrix (\ref{Rlip})
\be
\lb{Hlip}
H(\xi) = \sum_{a=1}^{N-1} \, h_{(a,a+1)}(\xi) \; ,
\ee
where the Hamiltonian densities $h_{(ab)}(x)$ are derived from (\ref{Rlip})
\be
\lb{lip3}
\begin{array}{c}
h_{(ab)}(\xi) =  2 \, \ln(\hat{q}_{(ab)})^{2}  + (\hat{q}_{(ab)})^{2\xi} \,
\ln (\hat{p}_{(a)}^{\; 2} \, \hat{p}_{(b)}^{\; 2} ) \, (\hat{q}_{(ab)})^{-2\xi}  = \\ [0.3cm]
=  \hat{p}_{(a)}^{-2\xi} \, \ln(\hat{q}_{(ab)})^{2}  \, \hat{p}_{(a)}^{\; 2\xi}
+ \hat{p}_{(b)}^{-2\xi} \, \ln(\hat{q}_{(ab)})^{2}  \, \hat{p}_{(b)}^{\; 2\xi} +
\ln (\hat{p}_{(a)}^{\; 2} \, \hat{p}_{(b)}^{\; 2} ) \; .
\end{array}
\ee
For $D=1$ and $\xi=1/2$ the Hamiltonian (\ref{Hlip}) reproduces the Hamiltonian for the Lipatov
integrable model \cite{44b}.

A remarkable fact is that for the algebra with the generators
$\{ \hat{q}_{(a)}^\mu , \hat{p}_{(b)}^\nu\}$ one can define a trace.
In particular, we need to define correctly the D-dimensional integral:
\be
\lb{trace}
 \int d^{^D} \!\! x  \, \langle x |
\hat{q}^{2 \alpha_1} \hat{p}^{2 \beta_1} \hat{q}^{2 \alpha_2}
 \hat{p}^{2 \beta_2} \! \dots
\hat{q}^{2 \alpha_n} \hat{p}^{2 \beta_n}
| x \rangle
= {\displaystyle c(\alpha_i,\beta_j)\!
\int \! \frac{d^{^D}\! x}{x^{2\gamma}} }
 \equiv {\rm Tr}(\hat{q}^{2 \alpha_1} \hat{p}^{2 \beta_1} \! \dots
\hat{q}^{2 \alpha_n} \hat{p}^{2 \beta_n} ) ,
\ee
where $\gamma = D/2 + \sum_{i} (\beta_i - \alpha_i)$
and $c(\alpha_i,\beta_j)$ is the coefficient
function.
Recall that the dimension regularization scheme requires the identity \cite{Hooft}:
\be
\lb{gi11}
\int \frac{d^{^D} x}{x^{^{2(D/2 + \alpha)}}} = 0 \;\;\;
\forall \alpha \neq 0 \; ,
\ee
and the integral (\ref{trace}) looks meaningless.
However, we can extend the definition for
(\ref{gi11}) at the point $\alpha = 0$ and, thus,
define the formal expression (\ref{trace}).
The definition is \cite{52'}:
\be
\lb{gi}
\int \frac{d^{^D} x}{x^{^{2(D/2 + \alpha)}}} =
\pi \Omega_{_D} \delta(|\alpha|) \; ,
\ee
where $\Omega_{_D} = \frac{2\pi^{^{D/2}}}{\Gamma(D/2)}$
is the area of the unit hypersphere in
$\mathbb{R}^D$, $\alpha = |\alpha| e^{^{i \arg(\alpha)}}$
and $\delta(.)$ is the one-dimensional
delta-function. The cyclic property $Tr(AB) = Tr(BA)$
for the trace (\ref{trace}) can be checked directly.
The trace operation (\ref{trace}) permits one
 to reduce \cite{52'} the evaluation of propagator-type
perturbative integrals
(and searching for their symmetries) to the evaluation
of vacuum perturbative integrals.
Further, breaking any of the propagators in the vacuum diagram,
one can obtain many remarkable nontrivial
relations between the propagator type integrals.
Sometimes these relations are called
"glue-and-cut" symmetry (for details see \cite{52'}, \cite{Chet})).

One can deduce another star-triangle relation \cite{Kash}
($x_i \in \mathbb{R}^D$)
$$
\left(\frac{ 2 \bar{\alpha}_1 \bar{\alpha}_3 }{ \alpha_2 }
\right)^{D/2}
\, W(x_3 - x_1 | \alpha_1 ) \, W(x_1 - x_2 | \alpha_2 ) \, W(x_2 - x_3 | \alpha_3 ) =
$$
\be
\lb{gr16}
= \int \, \frac{d^D x}{\pi^{D/2}} \, W(x_1 -x | \bar{\alpha}_3) \,
W(x_3 - x | \bar{\alpha}_2 ) \, W(x_2 -x | \bar{\alpha}_1 ) \; ,
\ee
where
$W(x|\alpha) = \exp \left( - x^2 / (2 \alpha ) \right)$ and the map
\be
\lb{map}
\bar{\alpha}_i = \frac{\alpha_1 \alpha_2 \alpha_3}{
 (\alpha_1 + \alpha_2 + \alpha_3)} \, \frac{1}{\alpha_i} \; , \;\;\;
\alpha_i = \frac{ \bar{\alpha}_1\bar{\alpha}_3 + \bar{\alpha}_2\bar{\alpha}_3
+\bar{\alpha}_1\bar{\alpha}_2 }{ \bar{\alpha}_i } \; ,
\ee
is the well known star-triangle transformation for resistances in electric networks.
The identity (\ref{gr16}) is related to the
local Yang-Baxter equation \cite{Maill}
and is also rewritten in the operator form \cite{Kash}
\be
\lb{gr17}
W(\hat{q} \, | \, \alpha_1 ) \,
W(\hat{p} \, | \, \alpha_2^{-1} ) \, W(\hat{q} \, | \, \alpha_3 ) =
W(\hat{p} \, | \, \bar{\alpha}_3^{-1} ) \,
W(\hat{q} \, | \, \bar{\alpha}_2 ) \, W(\hat{p} \, | \, \bar{\alpha}_1^{-1} )  \; .
\ee
To obtain (\ref{gr16}) from (\ref{gr17}),
 we have used the representations
$$
\langle x | e^{ \frac{1}{\alpha}\hat{q}^2 }  |y \rangle =
e^{ \frac{1}{\alpha}(x)^2 } \, \delta^D(x-y) \; , \;\;\;
\langle x | e^{- \frac{1}{2} \, \alpha \hat{p}^2 }  |y \rangle =
(2\pi \alpha)^{-D/2}  \,
e^{ -\frac{1}{2 \alpha} (x-y)^2} \; .
$$
It is tempting to apply identities
(\ref{gr16}) -- (\ref{gr17}) for investigation of
symmetries and analytical calculations of massive
perturbative multi-loop integrals written in the
$\alpha$- representation. Besides, we hope that the
local star-triangle relations (\ref{gr16}),
(\ref{gr17}) will help in constructing a massive deformation of the star-triangle relations
(\ref{4.11}), (\ref{uniq}). The generalizations
 of the star-triangle relations
(\ref{4.11}), (\ref{uniq}) for spinorial and tensor
particles were considered in \cite{52'',FP,DICh2,Mitra}.

\vspace{0.3cm}
$\bullet$ Note that we have not considered at all the numerous applications of
quantum Lie groups and algebras with deformation parameters $q$ satisfying the
conditions $q^{N} =1$, i.e., when the parameters $q$ are equal to the
roots of unity. These applications (see, for example, Ref. \cite{53})
appear mostly in the context of the topological and 2D conformal field theories and are
associated with the specific theory of representations of such quantum
groups that, generally speaking, can no longer be regarded as the
deformation of the classical Lie groups and algebras.


\section{REFERENCES}

\noindent
{\large \bf  Selected books} \\

{\small

{\bf 1.} E. Sweedler, Hopf Algebras, Benjamin (1969).

{\bf 2.} E. Abe, Hopf Algebras, Cambridge University Press (1977).

{\bf 3.} R.J. Baxter, Exactly Solved Models in Statistical Mechanics,
Academic Press, New York (1982).

{\bf 4.} L.D. Faddeev and L. Takhtajan, Hamiltonian methods in the theory
of solitons, Springer (1987).

{\bf 5.} P. Martin, Potts Models and Related Problems in
Statistical Mechanics, World Scientific, Singapore (1991).

{\bf 6.} Yu.I. Manin,  Topics in noncommutative geometry,
M.B. Porter Lecture Series, Princeton University Press, New Jersey (1991).

{\bf 7.} J. Madore, An introduction to noncommutative differential geometry and its physical applications (Vol. 257),
Cambridge University Press (1999).

{\bf 8.} G. Landi, An introduction to noncommutative spaces and their geometries (Vol. 51), Springer Science \& Business Media
(2003).

{\bf 9.} G. Lusztig, Introduction to quantum groups, Progress in Math. Birkhauser Boston (1993).

{\bf 10.} V.E. Korepin, A.G. Izergin, and N.M. Bogoliubov, Quantum Inverse Scattering
Method, Correlation Functions and Algebraic Bethe Ansatz, Cambridge Univ. Press (1993).

{\bf 11.} Zhong-Qi Ma, Yang-Baxter Equation and quantum enveloping algebras, World Scientific (1993).

{\bf 12.} A. Connes, Noncommutative geometry, Academic Press (1994).

{\bf 13.} V. Chari and A. Pressley, A guide to quantum groups,
Cambridge Univ. Press (1994).

{\bf 14.} C. Kassel, Quantum groups, Springer-Verlag GTM 155 (1994).

{\bf 15.} Sh. Majid, Foundations of quantum group theory, Cambridge Univ. Press (1995).

{\bf 16.} C. Gomez, M. Ruiz-Altaba and G. Sierra, Quantum Groups in Two-Dimensional Physics, Cambridge
University Press (1996).

{\bf 17.} C. Kassel, M. Rosso and V. Turaev, Quantum groups and knot invariants, Panoramas et
Synth\`{e}ses 5, Soci\'{e}t\'{e} Math\'{e}matique de France (1997).

{\bf 18.} A. Klimyk and K. Schm\"{u}dgen, Quantum groups and their representations, Springer (1997).

{\bf 19.} P. Etingof and O. Schiffmann, Lectures on Quantum Groups, Intern. Press (1998).

{\bf 20.} A. Molev, Yangians and Classical Lie Algebras,
Math. S. and M., 143, AMS (2007).

{\bf 21.} A. Molev, Sugawara Operators for Classical Lie Algebras,
Math. S. and M., 229, AMS (2018).
 }



\end{document}